%% file: RH_DOP_4.tex
\documentclass{article}
\usepackage{fullpage,amssymb,amsmath,amsthm}
\usepackage{epsfig,color}
\usepackage{url}
\usepackage{supertabular}

\begin{document}
\urldef{\miller}\url{millerpd@umich.edu}
\urldef{\baik}\url{baik@umich.edu}
\urldef{\kriecherbauer}\url{Thomas.Kriecherbauer@ruhr-uni-bochum.de}
\urldef{\mclaughlin}\url{mcl@amath.unc.edu}
\providecommand{\del}{{\Delta}} \providecommand{\nab}{{\nabla}}
\newtheorem{theorem}{Theorem}
\newtheorem{rhp}[theorem]{Riemann-Hilbert Problem}
\newtheorem{ip}[theorem]{Interpolation Problem}
\newtheorem{prop}[theorem]{Proposition} \newcommand{\cx}{{\mathbb C}}
\newtheorem{corollary}[theorem]{Corollary}
\newtheorem{Lemma}[theorem]{Lemma} \newtheorem{definition}[theorem]{Definition}

\numberwithin{theorem}{section}

\newenvironment{remark}{$\triangleleft$ {\bf
Remark:}}{$\triangleright$\medskip} \newcommand{\nat}{{\mathbb N}}
\newcommand{\Z}{{\mathbb Z}} \newcommand{\mat}[1]{{\bf #1}}

\author{J. Baik\thanks{ Department of Mathematics, University of Michigan.
Email: \baik.  }\and T.
Kriecherbauer\thanks{Fakult\"at f\"ur Mathematik,
Ruhr-Universit\"at Bochum.  Email: \kriecherbauer.  } \and K. T.-R.
McLaughlin\thanks{Department of Mathematics, University of North
Carolina at Chapel Hill.  Email: \mclaughlin. } \and P. D.
Miller\thanks{Department of Mathematics, University of Michigan.
Email: \miller. }}
\title{Uniform Asymptotics for Polynomials Orthogonal With Respect
to a General Class of Discrete Weights and Universality Results
for Associated Ensembles}
\date{\today} \maketitle

\begin{abstract}
  A general framework is developed for the asymptotic analysis of
  systems of polynomials orthogonal with respect to measures supported
  on finite sets of nodes.  Starting from a purely discrete
  interpolation problem for rational matrices whose solution encodes
  the polynomials, we show how the poles can be removed in favor of
  discontinuities along certain contours, turning the problem into an
  equivalent Riemann-Hilbert problem that we analyze with the help of
  an appropriate equilibrium measure related to weighted logarithmic
  potential theory.  For a large class of general weights and general
  distributions of nodes (not necessarily uniform), we calculate
  leading-order asymptotic formulae for the polynomials, with error
  bound inversely proportional to the number of nodes.  We obtain a
  number of asymptotic formulae that are valid in different
  overlapping regions whose union is the entire complex plane.  We
  prove exponential convergence of zeros to the nodes of
  orthogonalization in saturated regions where the equilibrium measure
  achieves a certain upper constraint.  Two of the asymptotic formulae
  for the polynomials display features distinctive of discrete
  weights: one formula uniformly valid near the endpoints of the
  interval of accumulation of nodes where the upper constraint is
  active is written in terms of the Euler gamma function, and another
  formula uniformly valid near generic band edges where the upper
  constraint becomes active are written in terms of both Airy
  functions $Ai(z)$ and $Bi(z)$ (by contrast $Ai(z)$ appears alone at
  band edges where the lower constraint becomes active, as with
  continuous weights).
  We illustrate our methods with the Krawtchouk polynomials and two
  families of polynomials belonging to the Hahn class.  We calculate
  the equilibrium measure for the Hahn weight.

  Universality of a number of statistics derived from so-called
  discrete orthogonal polynomial ensembles (discrete analogues of
  random matrix ensembles) is established using asymptotics for the
  discrete orthogonal polynomials.  In particular, we establish the
  universal nature of the discrete sine and Airy kernels as models for
  the correlation functions in certain regimes, and we prove
  convergence of distributions of extreme particles near band edges to
  the Tracy-Widom law.  We apply these results to the problem of
  computing asymptotics of statistics for random rhombus tilings of a
  large hexagon.  This problem is described in terms of discrete
  orthogonal polynomial ensembles corresponding to Hahn-type
  polynomials.  Therefore, combining the universality theory with our
  specific calculations of the equilibrium measure for the Hahn
  weights yields new error estimates and edge fluctuation phenomena
  for this statistical model.
\end{abstract}

\tableofcontents

\section{Introduction}
Our basic aim is to deduce asymptotic properties of polynomials that are
orthogonal with respect to pure point measures supported on finite
sets.  Let \label{symbol:number} $N\in\nat$, and consider $N$ distinct real
nodes $x_{N,0}<x_{N,1}<\dots <x_{N,N-1}$ to be given; together the
nodes make up the support of the pure point measures we consider.  We
use the notation
\begin{equation}
X_N:=
\{x_{N,n}\}_{n=0}^{N-1}\,,\hspace{0.2 in}
\text{where $x_{N,j}<x_{N,k}$ whenever $j<k$,}
\end{equation}
\label{symbol:nodeset}
\label{symbol:node}
for the support set.  Along with nodes we are given positive weights
\label{symbol:weight} $w_{N,0},w_{N,1},\dots, w_{N,N-1}$, which are
the magnitudes of the point masses located at the corresponding nodes.
We will occasionally use the alternate notation $w(x)$, $x\in X_N$ 
\label{symbol:weightfunc} for
a weight on the set of nodes $X_N$; thus
\begin{equation}
  w(x_{N,n})= w_{N,n}\,,\hspace{0.2 in}n=0,1,2,\dots,N-1\,.
\label{eq:weightrewrite}
\end{equation}
One should not infer from this notation that $w(x)$ has any meaning
for any $x$, complex or real, other than those $x\in X_N$; even if $w$
has a convenient functional form, we will only ever evaluate $w(x)$
when $x\in X_N$.  The {\em discrete orthogonal polynomials} associated
with this data are polynomials
\label{symbol:DOP} $\{p_{N,k}(z)\}_{k=0}^{N-1}$ where $p_{N,k}(z)$ is
of degree exactly $k$ with a positive leading
coefficient and where
\begin{equation}
\sum_{n=0}^{N-1}p_{N,k}(x_{N,n})p_{N,l}(x_{N,n})w_{N,n}=\delta_{kl}\,.
\label{eq:ortho}
\end{equation}
Writing \label{symbol:coeffs} $p_{N,k}(z)=c_{N,k}^{(k)}z^k +
\dots + c_{N,k}^{(0)}$, we introduce distinguished notation for the
positive leading coefficient:
\begin{equation}
\gamma_{N,k}:=c_{N,k}^{(k)}\,,
\end{equation}
\label{symbol:leadingcoeff}
and we denote by $\pi_{N,k}(z)$ the associated monic polynomial:
\begin{equation}
\pi_{N,k}(z):=\frac{1}{\gamma_{N,k}}p_{N,k}(z)\,.
\end{equation}
\label{symbol:MDOP}
The discrete orthogonal polynomials exist and are uniquely
determined by the orthogonality conditions because the inner product
associated with (\ref{eq:ortho}) is positive definite on ${\rm
span}(1,z,z^2,\dots,z^{N-1})$ but is degenerate on larger spaces of
polynomials.  The polynomials $p_{N,k}(z)$ may be constructed from the
monomials by a Gram-Schmidt process.  A general reference for
properties of orthogonal polynomials specific to the discrete case is
the book of Nikiforov, Suslov, and Uvarov
\cite{NikiforovSU91}.

One well-known elementary property of the discrete orthogonal
polynomials is an exclusion principle for the zeros that forbids more
than one zero from lying between adjacent nodes:
\begin{prop}
Each discrete orthogonal polynomial $p_{N,k}(z)$ has $k$ simple real
zeros.  All zeros lie in the range $x_{N,0} < z < x_{N,N-1}$ and no more
than one zero lies in the closed interval $[x_{N,n},x_{N,n+1}]$
between any two consecutive nodes.
\label{prop:confine}
\end{prop}

\begin{proof}
From the Gram-Schmidt process it follows that the coefficients of
$p_{N,k}(z)$ are all real.  Suppose that $p_{N,k}(z)$ were to vanish
to $n$th order for some nonreal $z_0$.  Then it follows that
$p_{N,k}(z)$ also vanishes to the same order at $z_0^*$, and thus that
$p_{N,k}(z)/[(z-z_0)^n(z-z_0^*)^n]$ is a polynomial of lower degree,
$k-2n\ge 0$.  By orthogonality, we must have on the one hand
\begin{equation}
\sum_{n=0}^{N-1}p_{N,k}(x_{N,n})\cdot\frac{p_{N,k}(x_{N,n})}
{|x_{N,n}-z_0|^{2n}}
\cdot w_{N,n}=0\,.
\end{equation}
On the other hand, the left-hand side is strictly positive because $k<N$
so $p_{N,k}(z)$ cannot vanish at all of the nodes.  So we have a contradiction
and the roots must be real.

The necessarily real roots are simple for a similar reason.  If $z_0$
is a real root of $p_{N,k}(z)$ of order greater than one, the quotient
$p_{N,k}(z)/(z-z_0)^2$ is a polynomial of degree $k-2\ge 0$, which must
be orthogonal to $p_{N,k}(z)$ itself:
\begin{equation}
\sum_{n=0}^{N-1}p_{N,k}(x_{N,n})\cdot\frac{p_{N,k}(x_{N,n})}{(x_{N,n}-z_0)^2}
\cdot w_{N,n}=0\,,
\end{equation}
but the left-hand side is manifestly positive, which gives the desired
contradiction.

If a simple real zero $z_0$ of $p_{N,k}(z)$ satisfies either
$z_0\le x_{N,0}$ or $z_0\ge x_{N,N-1}$, then we repeat the above argument
considering the polynomial $p_{N,k}(z)/(z-z_0)$ of degree $k-1\ge 0$ to
which $p_{N,k}(z)$ must be orthogonal, but for which the inner product
is strictly of one sign.

Finally if more than one zero of $p_{N,k}(z)$ were to lie between the
consecutive nodes $x_{N,n}$ and $x_{N,n+1}$, then we can certainly select
two of them, say $z_0$ and $z_1$, and construct the polynomial
$p_{N,k}(z)/[(z-z_0)(z-z_1)]$ of degree $k-2\ge 0$.  Again, this polynomial
must be orthogonal to $p_{N,k}(z)$, but the corresponding inner product
is of one definite sign, leading to a contradiction.
\end{proof}

Our goal is to establish the asymptotic behavior of the polynomials
$p_{N,k}(z)$ or their monic counterparts $\pi_{N,k}(z)$ in the limit
of large degree, assuming certain asymptotic properties of the nodes
and the weights.  In particular, the number of nodes must necessarily
increase to admit polynomials with arbitrarily large degree, and the
weights we consider involve an exponential factor with exponent
proportional to the number of nodes (such weights are sometimes called
{\em varying weights}).  We will obtain pointwise asymptotics with
precise error bound uniformly valid in the whole complex plane.  Our
assumptions on the nodes and weights include as special cases all
relevant classical discrete orthogonal polynomials, but are
significantly more general; in particular, we will consider nodes that
are not necessarily equally spaced.  While the number of nodes will
become large, it is important for our methods that this number is
finite; general weights supported on an
infinite discrete set of nodes require modifications of the methods we
will describe and will be considered in a subsequent paper.

\subsection{Basic assumptions.}
\label{sec:basicassumptions}
We will establish rigorous asymptotics for the discrete orthogonal polynomials
subject to the following fundamental assumptions.
\subsubsection{The nodes.}
\label{sec:C1}
We suppose the existence of a {\em node density function} 
\label{symbol:rho0}
$\rho^0(x)$ that is real-analytic in a complex neighborhood
of a closed interval \label{symbol:nodeinterval} $[a,b]$, and that satisfies:
\begin{equation}
\int_a^b\rho^0(x)\,dx=1\,,
\label{eq:rho0prob}
\end{equation}
and
\begin{equation}
\text{$\rho^0(x)>0$ strictly, for all $x\in[a,b]$\,.}
\label{eq:rho0nonzero}
\end{equation}
The nodes are then defined precisely in terms of the density function
$\rho^0(x)$ by the quantization rule
\begin{equation}
\int_a^{x_{N,n}}\rho^0(x)\,dx = \frac{2n+1}{2N}
\label{eq:BS}
\end{equation}
for $N\in\nat$ and $n=0,1,2,\dots,N-1$.
Thus, the nodes lie in a bounded open interval $(a,b)$ and are distributed
with density $\rho^0(x)$.
\subsubsection{The weights.}
\label{sec:C2}
Without loss of generality, we write the weights in the form
\begin{equation}
w_{N,n}=(-1)^{N-1-n}e^{-NV_N(x_{N,n})}\mathop{\prod_{m=0}^{N-1}}_{m\neq n}
(x_{N,n}-x_{N,m})^{-1}=e^{-NV_N(x_{N,n})}\mathop{\prod_{m=0}^{N-1}}_{m\neq n}
|x_{N,n}-x_{N,m}|^{-1}\,.
\label{eq:weightform}
\end{equation}
\label{symbol:VN}
No generality has been sacrificed with this representation because the
family of functions $\{V_N(x)\}$ is {\em apriori} specified only at
the nodes; in other words, given positive weights $\{w_{N,n}\}$ one may
solve \eqref{eq:weightform} uniquely for the $N$ quantities $\{V_{N}(x_{N,n})\}$.  However,
we now assume that for each sufficiently large $N$, $V_N(x)$ may be
taken to be a real-analytic function defined in a complex neighborhood $G$ of
the closed interval $[a,b]$, and that
\begin{equation}
V_N(x)=V(x)+\frac{\eta(x)}{N}
\label{eq:VNexpand}
\end{equation}
where $V(x)$ \label{symbol:V} is a fixed real-analytic function
defined in $G$, and
\begin{equation}
\limsup_{N\rightarrow\infty}\,\sup_{z\in G}|\eta(z)|<\infty\,.
\label{eq:etaNcontrol}
\end{equation}
\label{symbol:eta}
Note that in general the correction $\eta(z)$ may depend on $N$,
although $V(x)$ may not.  In some cases
({\em e.g.} Krawtchouk polynomials; see
\S~\ref{sec:Krawtchouk}) it is possible to take $V_N(x)\equiv V(x)$
for all $N$, in which case $\eta(x)\equiv 0$.  However, the freedom of
assuming $\eta(x)\not\equiv 0$ is useful to handle other cases ({\em
e.g.} the Hahn and associated Hahn polynomials; see
\S~\ref{sec:Hahn}).
While \eqref{eq:weightform} may be written
for any system of positive weights, the condition that
\eqref{eq:VNexpand} should hold restricts attention to systems of
weights that have analytic continuum limits in a certain precise
sense.

\begin{remark}
The familiar examples of classical discrete orthogonal polynomials
correspond to nodes that are equally spaced, say on $(a,b)=(0,1)$
(in which case we have $\rho^0(x)\equiv 1$).  In this special case,
the product factor on the right-hand side of (\ref{eq:weightform})
becomes simply
\begin{equation}
\mathop{\prod_{m=0}^{N-1}}_{m\neq n}|x_{N,n}-x_{N,m}|^{-1}=
\frac{N^{N-1}}{n!(N-n-1)!}\,.
\end{equation}
Using Stirling's formula to take the continuum limit of this factor
(that is, considering $N\rightarrow\infty$ with $n/N\rightarrow x$)
shows that in these cases the leading term in the formula
\eqref{eq:weightform}
is a continuous weight on $(0,1)$:
\begin{equation}
w_{N,n}\sim w(x):=C\left(\frac{e^{-V(x)}}{x^{x}(1-x)^{1-x}}\right)^{N}
\label{eq:clweight}
\end{equation}
as $N\rightarrow\infty$ and $n/N\rightarrow x\in (0,1)$, where $C$ is
independent of $x$.  However, taking the continuum limit of the weight
first to arrive at a formula like \eqref{eq:clweight}, and then
obtaining asymptotics of the polynomials of degree proportional to $N$
as $N\rightarrow\infty$ is not equivalent to the double-scaling limit
process we will consider here.  Our results display new phenomena
because we simultaneously take the continuum limit as the degree of
the polynomials grows.
\end{remark}

Our choice of the form (\ref{eq:weightform}) for the weights is
motivated by several specific examples of classical discrete
orthogonal polynomials.  The form (\ref{eq:weightform}) is
sufficiently general for us to carry out useful calculations related
to proofs of universality conjectures arising in certain types of
random tiling problems, random growth models, and last-passage
percolation problems.

\subsubsection{The degree.}
We assume that the degree $k$ of the polynomial of interest is tied to
the number $N$ of nodes by a relation of the form
\begin{equation}
k=cN + \kappa
\end{equation}
\label{symbol:k}
where $c\in (0,1)$ \label{symbol:c} is a fixed parameter, and $\kappa$
\label{symbol:kappa} remains bounded as $N\rightarrow\infty$.

\subsection{Simplifying assumptions of genericity.}
In order to keep our exposition as simple as possible, we make
further assumptions that exclude certain nongeneric triples $(\rho^0(x),V(x),c)$.  These assumptions depend on the functions 
$\rho^0(x)$ and $V(x)$, and on the parameter $c$, in an implicit manner 
that is
easier to describe once some auxiliary quantities have been introduced.  
They will be given in \S~\ref{sec:C3}.  

In regard to these particular assumptions, we want to stress two
points.  First, the excluded triples are nongeneric in the sense that
any perturbation of, say, the parameter $c$ will immediately return us
to the class of triples for which all of our results are valid.
The discussion at the beginning of \S~\ref{sec:Airy} provides some
insight into the generic nature of our assumptions.
Second, the discrete orthogonal polynomials corresponding to
nongeneric triples can be analyzed by the same basic method that we
use here, with many of the same results.  To do this, the proofs we
present will require modifications to include additional local
analysis near certain isolated points in the complex $z$-plane.  Some
such modifications have already been described in detail in the
context of asymptotics for polynomials orthogonal with respect to
continuous weights in \S~5 of \cite{DeiftKMVZ99}. The remaining
modifications have to do with nongeneric behavior near the endpoints
of the interval $[a,b]$, and while the corresponding local analysis
has not been done before, it can be expected to be of a similar
character.

\subsection{The goal.}
Given an interval $[a,b]$, appropriate fixed functions $\rho^0(x)$ and
$V(x)$, appropriate sequences $\eta(x)=\eta_N(x)$ and
$\kappa=\kappa_N$, and a constant $c\in (0,1)$, we wish to find
accurate asymptotic formulae, valid in the limit $N\rightarrow\infty$
with rigorous error bounds, for the polynomial $\pi_{N,k}(z)$.  These
formulae should be uniformly valid in overlapping regions of the
complex $z$-plane.  We will also require asymptotic formulae for
related quantities, like the zeros of $\pi_{N,k}(z)$, the three-term
recurrence coefficients, and the reproducing kernels $K_{N,k}(x,y)$.

\subsection{Motivation.}  Our work on this subject is connected with
three different themes of current research.  First of all, in the
context of approximation theory, there has been recent activity
\cite{DKMVZstrong,DeiftKMVZ99} in the study of polynomials orthogonal
on the real axis with respect to general continuous varying weights
and the corresponding large degree pointwise asymptotics.  The setting
for this work is the characterization of the orthogonal polynomials in
terms of the solution of a certain matrix-valued Riemann-Hilbert
problem \cite{FIK}.  These methods are not at all particular to any
special classical formulae for weights; they are completely general.
Thus, a natural question to ask is whether it is possible to
generalize the methods further to handle discrete weights.  Discrete
weights are of such a fundamentally different character than their
continuous counterparts that this would require the development of new
analytical tools.  In particular, each point mass added to the weight
amounts to a pole in the matrix solution of the Riemann-Hilbert
problem, so the problem is to analyze the asymptotics of an
accumulation of poles.

Secondly, there has been some recent progress
\cite{KamvissisMM03,Miller02} in the integrable systems literature
concerning the problem of computing asymptotics for solutions of
integrable nonlinear partial differential equations ({\em e.g.} the
nonlinear Schr\"odinger equation) in the limit where the spectral data
associated with the solution via the inverse-scattering transform is
made up of a large number of discrete eigenvalues.  Significantly,
inverse-scattering theory also exploits much of the theory of matrix
Riemann-Hilbert problems, and it turns out that the discrete
eigenvalues appear as poles in the corresponding matrix-valued
unknown.  So, the methods recently developed in the context of
inverse-scattering actually suggest a general scheme by means of which
an accumulation of poles in the matrix unknown can be analyzed.

Finally, a number of problems in probability theory have recently been
identified that are in some sense solved in terms of discrete
orthogonal polynomials, and certain statistical questions can be
translated into corresponding questions about the asymptotic behavior
of the polynomials.  The particular problems we have in mind are
related to statistics of random tilings of various shapes, and also to
certain natural measures on sets of partitions.  The joint probability
distributions in these problems are examples of so-called discrete
orthogonal polynomial ensembles \cite{Johansson02}.  Roughly speaking,
the analogy is that the relationship between universal asymptotic
properties of discrete orthogonal polynomials and universal statistics
of so-called discrete orthogonal polynomial ensembles is the same as
the relationship between universal asymptotic properties of
polynomials orthogonal with respect to continuous weights and
universal eigenvalue statistics of certain random matrix ensembles.  We will
give many more details later on, but for now the point is that the techniques required for computing
asymptotics of discrete orthogonal polynomials with general weights have become available at
just the time when questions that can be answered with these tools are
appearing in the applied literature.

\subsection{Methodology.}
%We study the discrete orthogonal polynomials subject to a number of
%assumptions on the distribution of nodes (see \S~\ref{sec:C1}), the
%form of the weights (see \S~\ref{sec:C2}), and the behavior of the
%so-called equilibrium measure (see \S~\ref{sec:C3}).  We also restrict
%attention to the limit of $N\rightarrow\infty$ with $k=Nc+\kappa$,
%where $N$ is the number of nodes $k$ is the degree of the polynomial,
%$c\in (0,1)$ is held fixed and $\kappa$ is bounded.  The assumptions
%concerning the equilibrium measure are technical and can be removed at
%the cost of some extra work.  The assumptions we make of the nodes and
%weights are more rigid, although we want to emphasize that there is
%enough generality to handle all the classical cases of discrete
%orthogonal polynomials and more.
\subsubsection{The basic interpolation problem.}
Given a natural number $N$, a set $X_N$ of nodes, and a set of
corresponding weights $\{w_{N,n}\}$, consider the possibility of
finding the matrix $\mat{P}(z;N,k)$ \label{symbol:P} solving the
following problem, where $k$ is an integer.
\begin{ip}
Find a $2\times 2$ matrix
$\mat{P}(z;N,k)$ with the following
properties:
\begin{enumerate}
\item
{\bf Analyticity}: $\mat{P}(z;N,k)$ is an analytic function of $z$ for
$z\in\cx\setminus X_N$.
\item
{\bf Normalization}: As $z\rightarrow\infty$,
\begin{equation}
\mat{P}(z;N,k)\left(\begin{array}{cc}z^{-k} & 0 \\\\ 0 & z^k\end{array}
\right)={\mathbb I} + O\left(\frac{1}{z}\right)\,.
\label{eq:norm}
\end{equation}
\item
{\bf Singularities}: At each node $x_{N,n}\in X_N$, the first column of
$\mat{P}(z;N,k)$ is analytic and the second column of 
$\mat{P}(z;N,k)$ has a simple
pole, where the residue satisfies the condition
\begin{equation}
\mathop{\rm Res}_{z=x_{N,n}}\mat{P}(z;N,k)=\lim_{z\rightarrow x_{N,n}}
\mat{P}(z;N,k)\left(\begin{array}{cc}0 & w_{N,n}\\ \\0 & 0\end{array}\right)
=\left(\begin{array}{cc}0 & w_{N,n}P_{11}(x_{N,n};N,k)\\\\
0 & w_{N,n}P_{21}(x_{N,n},N,k)\end{array}\right)
\label{eq:poles}
\end{equation}
for $n=0,\dots,N-1$.
\end{enumerate}
\label{rhp:DOP}
\end{ip}
This problem is a
discrete version of the Riemann-Hilbert problem appropriate for orthogonal polynomials with continuous weights that was first
used in \cite{FIK} (see also \cite{DKMVZstrong,DeiftKMVZ99}).  The solution
of this problem encodes all quantities of relevance to a study of the discrete orthogonal polynomials, as we will now see. 
\begin{prop}
Interpolation Problem~\ref{rhp:DOP} has a unique solution when $0\le
k\le N-1$.  In this case,
\begin{equation}
\mat{P}(z;N,k)=
\left(\begin{array}{cc}
\pi_{N,k}(z) & \displaystyle
\sum_{n=0}^{N-1}\frac{w_{N,n}\pi_{N,k}(x_{N,n})}{z-x_{N,n}}\\\\
\gamma_{N,k-1}p_{N,k-1}(z) &\displaystyle \sum_{n=0}^{N-1}
\frac{w_{N,n}\gamma_{N,k-1}p_{N,k-1}(x_{N,n})}{z-x_{N,n}}
\end{array}\right)
\label{eq:usoln}
\end{equation}
if $k>0$ and
\begin{equation}
\mat{P}(z;N,0)=
\left(\begin{array}{cc}
1 &\displaystyle \sum_{n=0}^{N-1}\frac{w_{N,n}}{z-x_{N,n}}\\\\
0 & 1\end{array}\right)\,.
\label{eq:usolnzero}
\end{equation}
\label{prop:solnrhp}
\end{prop}
\begin{proof}
Consider the first row of $\mat{P}(z;N,k)$.  According to
(\ref{eq:poles}), the function $P_{11}(z;N,k)$ is an entire function
of $z$.  Because $k\ge 0$ it follows from the normalization condition
(\ref{eq:norm}) that in fact $P_{11}(z;N,k)$ is a monic polynomial of
degree exactly $k$.  Similarly, from the characterization
(\ref{eq:poles}) of the simple poles of $P_{12}(z;N,k)$, we see that
$P_{12}(z;N,k)$ is necessarily of the form
\begin{equation}
P_{12}(z;N,k)=e_1(z)+\sum_{n=0}^{N-1}\frac{w_{N,n}P_{11}
(x_{N,n};N,k)}{z-x_{N,n}}
\end{equation}
where $e_1(z)$ is an entire function.  The normalization condition
(\ref{eq:norm}) for $k\ge 0$ immediately requires, via Liouville's
Theorem, that $e_1(z)\equiv 0$, and then when $|z|>\max_n |x_{N,n}|$ we
have by geometric series expansion that
\begin{equation}
P_{12}(z;N,k)=\sum_{m=0}^{\infty}\left(\sum_{n=0}^{N-1}P_{11}(x_{N,n};N,k)x_{N,n}^m w_{N,n}\right)\frac{1}{z^{m+1}}\,.
\end{equation}
According to the normalization condition (\ref{eq:norm}),
$P_{12}(z;N,k)=o(z^{-k})$ as $z\rightarrow \infty$; therefore it follows that the monic polynomial $P_{11}(z;N,k)$ of degree exactly $k$ must satisfy \begin{equation}
\sum_{n=0}^{N-1}P_{11}(x_{N,n};N,k)x_{N,n}^m w_{N,n}=0\hspace{0.2 in}
\mbox{for}\hspace{0.2 in}m=0,1,2,\dots,k-1\,.
\end{equation}
As long as $k\le N-1$, these conditions uniquely identify $P_{11}(z;N,k)$
with the monic discrete orthogonal polynomial $\pi_{N,k}(z)$.  The
existence and uniqueness of $\pi_{N,k}(z)$ for such $k$ is guaranteed
given distinct orthogonalization nodes and positive weights
(definiteness of the inner product).

The second row of $\mat{P}(z;N,k)$ is studied similarly.  The function
$P_{21}(z;N,k)$ is seen from (\ref{eq:poles}) to be an entire function
of $z$, that according to the normalization condition (\ref{eq:norm})
must be a polynomial of degree at most $k-1$ (for the special case of
$k=0$ these conditions immediately imply that $P_{21}(z;N,0)\equiv
0$).  The characterization (\ref{eq:poles}) implies that
$P_{22}(z;N,k)$ can be expressed in the form
\begin{equation}
P_{22}(z;N,k)=e_2(z)+\sum_{n=0}^{N-1}\frac{w_{N,n}P_{21}(x_{N,n};N,k)}
{z-x_{N,n}}
\end{equation}
where $e_2(z)$ is an entire function.  If $k=0$, then
$P_{22}(z;N,0)=e_2(z)$ and then according to the normalization
condition (\ref{eq:norm}) we must take $e_2(z)\equiv 1$.  On the other
hand, if $k>0$, then (\ref{eq:norm}) implies that $P_{22}(z;N,k)$
decays for large $z$ and therefore we must take $e_2(z)\equiv 0$ in
this case.  Expanding the denominators in geometric series for
$|z|>\max_n |x_{N,n}|$, we find
\begin{equation}
P_{22}(z;N,k)=\sum_{m=0}^{\infty}\left(\sum_{n=0}^{N-1}P_{21}(x_{N,n};N,k)
x_{N,n}^mw_{N,n}\right)\frac{1}{z^{m+1}}\,.
\end{equation}
Imposing the normalization conditions (\ref{eq:norm}) we now insist that
$P_{22}(z;N,k)=z^{-k} + O(z^{-k-1})$ as $z\rightarrow\infty$; therefore
\begin{equation}
\sum_{n=0}^{N-1}P_{21}(x_{N,n};N,k)
x_{N,n}^mw_{N,n}=0\,,\hspace{0.2 in}\mbox{for}\hspace{0.2 in}
m=0,1,2,\dots,k-2\,,
\label{eq:secondroworthogonality}
\end{equation}
and
\begin{equation}
\sum_{n=0}^{N-1}P_{21}(x_{N,n};N,k)x_{N,n}^{k-1}w_{N,n}=1\,.
\label{eq:secondrownorm}
\end{equation}
Using (\ref{eq:secondroworthogonality}), the condition (\ref{eq:secondrownorm})
can be replaced by
\begin{equation}
\sum_{n=0}^{N-1}P_{21}(x_{N,n};N,k)\pi_{N,k-1}(x_{N,n})w_{N,n}=1\hspace{0.2 in}
\mbox{or}\hspace{0.2 in}
\sum_{n=0}^{N-1}\left[\frac{1}{\gamma_{N,k-1}}P_{21}(x_{N,n};N,k)
\right]p_{N,k-1}(x_{N,n})w_{N,n}=1
\,.
\end{equation}
These conditions therefore uniquely identify
$P_{21}(z;N,k)/\gamma_{N,k-1}$ with the orthogonal polynomial
$p_{N,k-1}(z)$.

The interpolation problem is thus solved uniquely by the the matrix
explicitly given by (\ref{eq:usoln}) for $k>0$ and by
(\ref{eq:usolnzero}) for $k=0$.
\end{proof}

An important feature of all systems of orthogonal polynomials, that is
present whether the weights are discrete or continuous, is the
well-known three-term recurrence relation.  See \cite{S91} for
details.  There are constants $a_{N,0},a_{N,1},\dots,a_{N,N-2}$
\label{symbol:aNk} and positive constants
$b_{N,0},b_{N,1},\dots,b_{N,N-2}$ \label{symbol:bNk} such that
\begin{equation}
zp_{N,k}(z)=b_{N,k}p_{N,k+1}(z)+a_{N,k}p_{N,k}(z)+b_{N,k-1}p_{N,k-1}(z)
\end{equation}
holds for $k=1,\dots,N-2$, while for $k=0$ one has
\begin{equation}
zp_{N,0}(z)=b_{N,0}p_{N,1}(z) + a_{N,0}p_{N,0}(z)\,.
\end{equation}
Necessarily, one has $b_{N,k}=\gamma_{N,k}/\gamma_{N,k+1}$.  The
constants in the three-term recurrence relation are also encoded in
the solution of Interpolation Problem~\ref{rhp:DOP}.  By expansion for
large $z$ of the explicit solution given in
Proposition~\ref{prop:solnrhp} we have the following:
\begin{corollary}
Let $A_{N,k}$, $B_{N,k}$, $C_{N,k}$, and $D_{N,k}$ denote certain
terms in the large $z$ expansion of the matrix elements of
$\mat{P}(z;N,k)$:
\begin{equation}
\begin{array}{rcl}
z^kP_{12}(z;N,k)&=&\displaystyle \frac{A_{N,k}}{z} + \frac{B_{N,k}}{z^2} + O\left(\frac{1}{z^3}\right)\\\\
\displaystyle\frac{1}{z^k}P_{11}(z;N,k)&=&\displaystyle 1 +\frac{C_{N,k}}{z} +O\left(\frac{1}{z^2}\right)\\\\
\displaystyle\frac{1}{z^k}P_{21}(z;N,k)&=&\displaystyle\frac{D_{N,k}}{z} + O\left(\frac{1}{z^2}\right)
\end{array}
\end{equation}
as $z\rightarrow\infty$.  Then
\begin{equation}
\begin{array}{rclrcl}
\gamma_{N,k} &= &\displaystyle \frac{1}{\sqrt{A_{N,k}}}\,,\hspace{0.2 in}
&\gamma_{N,k-1}&=&\sqrt{D_{N,k}}\,,\\\\
a_{N,k}&=&\displaystyle C_{N,k}+\frac{B_{N,k}}{A_{N,k}}\,,\hspace{0.2 in}
&b_{N,k-1}&=&\sqrt{A_{N,k}D_{N,k}}\,.
\end{array}
\label{eq:3termcoeffs}
\end{equation}
\label{cor:three-term}
\end{corollary}

\subsubsection{Triangularity of residue matrices and dual polynomials.}
%In this section, we present and then modify for our later convenience
%a mathematical problem, involving a meromorphic matrix-valued unknown, that 
%is equivalent to the orthogonality
%conditions \eqref{eq:ortho} in that its solution is built from the
%same set of orthogonal polynomials.  On the other hand, we will soon
%see that one gains from
%this alternate characterization
%a considerable analytical advantage in asymptotic calculations over other 
%approaches based more directly on the orthogonality property.
\label{sec:dual}
Interpolation Problem~\ref{rhp:DOP} involves residue matrices that are
upper-triangular.  An essential aspect of our methodology will be to
modify the matrix $\mat{P}(z;N,k)$ in order to selectively reverse the
triangularity of the residue matrices near certain individual nodes
$x_{N,n}$.  
Let $\del\subset {\mathbb Z}_N$ \label{symbol:del}\label{symbol:ZN} where
\begin{equation}
{\mathbb Z}_N:=\{0,1,2,\dots,N-1\}\,,
\end{equation}
and denote the number of elements in $\del$ by \label{symbol:Numdel}
$\#\del$.  We will reverse the triangularity for those nodes $x_{N,n}$
for which $n\in\del$.  Consider the matrix $\mat{Q}(z;N,k)$
\label{symbol:Q} related to the solution $\mat{P}(z;N,k)$ of
Interpolation Problem~\ref{rhp:DOP} as follows:
\begin{equation}
\mat{Q}(z;N,k):=\mat{P}(z;N,k)\left[
\prod_{n\in\del} (z-x_{N,n})\right]^{-\sigma_3}=
\mat{P}(z;N,k)\left(\begin{array}{cc} \displaystyle \prod_{n\in\del}
(z-x_{N,n})^{-1} & 0 \\\\ 0 &\displaystyle \prod_{n\in\del}(z-x_{N,n})
\end{array}\right)\,.
\label{eq:PtoQ}
\end{equation}
Here $\sigma_3$ \label{symbol:sigma3} is a Pauli matrix:
\begin{equation}
\sigma_3:=\left(\begin{array}{cc} 1 & 0 \\\\ 0 & -1\end{array}\right)\,.
\end{equation}
It is easy to check that
the matrix $\mat{Q}(z;N,k)$ so defined is an analytic function of $z$ for
$z\in\cx\setminus X_N$ that satisfies the normalization condition
\begin{equation}
\mat{Q}(z;N,k)\left(\begin{array}{cc}z^{\#\del-k} & 0 \\\\ 0 &
z^{k-\#\del}\end{array} \right)={\mathbb I} + O\left(\frac{1}{z}\right)\,,\hspace{0.2 in}\mbox{as $z\rightarrow\infty$.}
\label{eq:norm-S}
\end{equation}
Furthermore, at each node $x_{N,n}$, the matrix $\mat{Q}(z;N,k)$ has
a simple pole.  If $n$ belongs to the complementary set
\begin{equation}
\nab:={\mathbb Z}_N\setminus\del\,,
\end{equation} 
\label{symbol:nab}
then the first column
is analytic at $x_{N,n}$ and the pole is in the second column such
that the residue satisfies the condition
\begin{equation}
\mathop{\rm Res}_{z=x_{N,n}}\mat{Q}(z;N,k)=\lim_{z\rightarrow x_{N,n}}
\mat{Q}(z;N,k)\left(\begin{array}{cc}0 & \displaystyle w_{N,n}
\prod_{m\in\del}(x_{N,n}-x_{N,m})^2\\\\ 0 & 0\end{array}\right)
\label{eq:polesinS}
\end{equation}
for $n\in\nab$.
If $n\in\del$, then the second column is analytic at $x_{N,n}$ and the
pole is in the first column such that the residue satisfies the condition
\begin{equation}
\mathop{\rm Res}_{z=x_{N,n}}\mat{Q}(z;N,k)=\lim_{z\rightarrow x_{N,n}}
\mat{Q}(z;N,k)\left(\begin{array}{cc}0 & 0 \\\\
\displaystyle
\frac{1}{w_{N,n}}\mathop{\prod_{m\in\del}}_{m\neq n}(x_{N,n}-x_{N,m})^{-2}
& 0\end{array}\right)
\label{eq:polesnotinS}
\end{equation}
for $n\in\del$.
Thus, the triangularity of the residue matrices has been reversed for
nodes in $\del\subset X_N$.

The relation between the solution $\mat{P}(z;N,k)$ of Interpolation
Problem~\ref{rhp:DOP} and the matrix $\mat{Q}(z;N,k)$ obtained
therefrom by selective reversal of residue triangularity gives rise in
a special case to a remarkable duality between pairs of weights
$\{w_{N,n}\}$ defined on the same set of nodes and their corresponding
families of discrete orthogonal polynomials that comes up in
applications.  Given nodes $X_N$ and weights $\{w_{N,n}\}$, the dual
polynomials arise by taking $\del={\mathbb Z}_N$ in the change of variables
(\ref{eq:PtoQ}), and then defining \label{symbol:Pbar}\label{symbol:bark}
\begin{equation}
\overline{\mat{P}}(z;N,\overline{k}):=
\sigma_1 \mat{Q}(z;N,k) \sigma_1\,,\hspace{0.2 in}\mbox{where}
\hspace{0.2 in}\overline{k}:=N-k\,.
\end{equation}
Here $\sigma_1$ \label{symbol:sigma1} is another Pauli matrix:
\begin{equation}
\sigma_1:=\left(\begin{array}{cc} 0 & 1 \\\\ 1 & 0\end{array}\right)\,.
\end{equation}
Thus, we are reversing the triangularity at all of the nodes, and
swapping rows and columns of the resulting matrix.  It is easy to check
that $\overline{\mat{P}}(z;N,\overline{k})$ satisfies
\begin{equation}
\overline{\mat{P}}(z;N,\overline{k})\left(\begin{array}{cc}
z^{-\overline{k}} & 0 \\ 0 & z^{\overline{k}}\end{array}\right) =
{\mathbb I} + O\left(\frac{1}{z}\right)\hspace{0.2 in}\mbox{as}
\hspace{0.2 in}z\rightarrow\infty
\end{equation}
and is a matrix with simple poles in the second column at all nodes,
such that
\begin{equation}
\mathop{\rm Res}_{z=x_{N,n}}\overline{\mat{P}}(z;N,\overline{k})=
\lim_{z\rightarrow x_{N,n}}
\overline{\mat{P}}(z;N,\overline{k})
\left(\begin{array}{cc}0 & \overline{w}_{N,n}\\ \\0 & 0\end{array}\right)
\end{equation}
holds for $n\in {\mathbb Z}_N$, where the ``dual weights'' 
\label{symbol:dualweights}
$\{\overline{w}_{N,n}\}$ are defined by the identity
\begin{equation}
w_{N,n}\overline{w}_{N,n}\mathop{\prod_{m=0}^{N-1}}_{m\neq n}(x_{N,n}-x_{N,m})^2 = 1\,.
\label{eq:dualweightsdefine}
\end{equation}
Comparing with Interpolation Problem~\ref{rhp:DOP} we see that
$\overline{P}_{11}(z;N,\overline{k})$ is the monic orthogonal
polynomial $\overline{\pi}_{N,\overline{k}}(z)$ \label{symbol:pidual}
of degree $\overline{k}$ associated with the dual weights
$\{\overline{w}_{N,j}\}$ (and the same set of nodes $X_N$).  In this
sense, families of discrete orthogonal polynomials always come in dual
pairs.  An explicit relation between the dual polynomials comes from
the representation of $\mat{P}(z;N,k)$ given by
Proposition~\ref{prop:solnrhp}:
\begin{equation}
\begin{array}{rcl}
\displaystyle
\overline{\pi}_{N,\overline{k}}(z)&=&\displaystyle
\overline{P}_{11}(z;N,\overline{k})\\\\
&=&\displaystyle
P_{22}(z;N,k)\prod_{n=0}^{N-1}(z-x_{N,n})\\\\
&=&\displaystyle
\sum_{m=0}^{N-1}w_{N,m}\gamma_{N,k-1}^2  \pi_{N,k-1}(x_{N,m})
\mathop{\prod_{n=0}^{N-1}}_{n\neq m}(z-x_{N,n})\,.
\end{array}
\label{eq:explicitdual}
\end{equation}
Since the left-hand side is a monic polynomial of degree
$\overline{k}=N-k$ and the right-hand side is apparently a polynomial
of degree $N-1$, equation (\ref{eq:explicitdual}) furnishes $k$
relations among the weights and the normalization constants
$\gamma_{N,k}$.

In particular, if we evaluate (\ref{eq:explicitdual}) for $z=x_{N,l}$
for some $l\in {\mathbb Z}_N$, then only one term from the sum on the
right-hand side survives and we find
\begin{equation}
\overline{\pi}_{N,\overline{k}}(x_{N,l})=\gamma_{N,k-1}^2
w_{N,l}\mathop{\prod_{n=0}^{N-1}}_{n\neq l}
(x_{N,l}-x_{N,n})\cdot\pi_{N,k-1}(x_{N,l})\,,
\label{eq:Borodin}
\end{equation}
an identity relating values of each discrete orthogonal polynomial and
a corresponding dual polynomial at any given node.  The identity
\eqref{eq:Borodin} has also been derived by Borodin \cite{Borodin01}.

Furthermore, by using (\ref{eq:explicitdual}) twice, along with the
fact that $\overline{\overline{\pi}}_k(z)\equiv \pi_k(z)$ ({\em i.e.}
duality is an involution), we can obtain some additional identities
involving the discrete orthogonal polynomials and their duals.  By
involution, (\ref{eq:explicitdual}) implies that
\begin{equation}
\pi_{N,k}(z)=\overline{\gamma}_{N,\overline{k}-1}^2\gamma_{N,k}^2
\sum_{m=0}^{N-1}\pi_{N,k}(x_{N,m})\mathop{\prod_{n=0}^{N-1}}_{n\neq
m}\frac{z-x_{N,n}}{x_{N,m}-x_{N,n}}\,.
\label{eq:dualofdual}
\end{equation}
\label{symbol:gammadual}
The sum on the right-hand side is the Lagrange interpolating
polynomial of degree $N-1$ (at most) that agrees with $\pi_{N,k}(z)$
at all $N$ nodes.  Of course this identifies the sum with
$\pi_{N,k}(z)$ itself, upon which we find the relation
\begin{equation}
\overline{\gamma}_{N,\overline{k}-1}=\frac{1}{\gamma_{N,k}}
\label{eq:gammagammabar}
\end{equation}
between the normalization constants for the discrete orthogonal
polynomials and their duals.

\begin{remark}
We want to point out that the notion of duality described here is
different from that explained in \cite{NikiforovSU91}.  The latter
generally involves relationships between families of discrete
orthogonal polynomials with two different sets of nodes of
orthogonalization.  For example, the Hahn polynomials are orthogonal
on a lattice of equally spaced points, and the polynomials dual to
the Hahn polynomials by the scheme of \cite{NikiforovSU91} are
orthogonal on a quadratic lattice for which $x_{N,n}-x_{N,n-1}$ is
proportional to $n$.  However, the polynomials dual to the Hahn
polynomials under the scheme described above are the associated Hahn
polynomials, which are orthogonal on the same equally-spaced nodes
as are the Hahn polynomials themselves.  The notion of duality described above coincides with that described in
\cite{Borodin01} and is also equivalent to the ``hole/particle
transformation'' considered by Johansson \cite{Johansson02}.
\end{remark}

\subsubsection{Outline of approach.}
The characterization of the discrete orthogonal polynomials in terms
of Interpolation Problem~\ref{rhp:DOP} is the starting point for our
asymptotic analysis.
%
%Our starting point will be a representation of the discrete orthogonal
%polynomials in terms of the solution of a matrix interpolation problem %rather than a Riemann-Hilbert
%problem that is of the traditional sort; instead of the unknown
%being piecewise analytic with given jumps along certain contours, it
%is a meromorphic function whose only singularities are simple poles.
%The residues of the poles are encoded in terms of certain
%upper-triangular matrices.  

Our rigorous analysis of $\mat{P}(z;N,k)$ consists of three steps:
\begin{enumerate}
\item 
We introduce a change of variables, transforming $\mat{P}(z;N,k)$ into
$\mat{X}(z)$, another matrix function of $z$.  The transformation mediating
between $\mat{P}(z;N,k)$ and $\mat{X}(z)$ is explicit and exact.  The matrix $\mat{X}(z)$ is shown to satisfy a matrix Riemann-Hilbert problem that is equivalent to Interpolation Problem~\ref{rhp:DOP}.  
\item
We construct an explicit model $\hat{\mat{X}}(z)$ for $\mat{X}(z)$ on the basis of formal asymptotics.  We call $\hat{\mat{X}}(z)$ a {\em global parametrix} for $\mat{X}(z)$.  
\item
We compare $\mat{X}(z)$ to the global parametrix $\hat{\mat{X}}(z)$ by
considering the error $\mat{E}(z):=\mat{X}(z)\hat{\mat{X}}(z)^{-1}$, which should be close to the identity matrix if the formally obtained global parametrix $\hat{\mat{X}}(z)$ is indeed a good approximation of $\mat{X}(z)$.  We rigorously analyze $\mat{E}(z)$ by viewing its definition
in terms of $\mat{X}(z)$ as another change of variables, since $\hat{\mat{X}}(z)$ is known explicitly from step 2.  
This means that we may pose an equivalent Riemann-Hilbert problem for $\mat{E}(z)$.  We prove that this
Riemann-Hilbert problem may be solved by a convergent Neumann series if $N$ is sufficiently large.  The series for $\mat{E}(z)$ is also an asymptotic series whose first term is the identity matrix, such that $\mat{E}(z)-\mathbb{I}$ is of order $1/N$ in a suitable precise sense.  This gives an asymptotic formula for the unknown matrix $\mat{X}(z) = \mat{E}(z)\hat{\mat{X}}(z)$.  Inverting the explicit change of variables
from step 1 linking $\mat{X}(z)$ with $\mat{P}(z;N,k)$, we finally arrive at an asymptotic formula for $\mat{P}(z;N,k)$.
\end{enumerate}
The first step in this process is the most crucial, since the explicit transformation from $\mat{P}(z;N,k)$ and $\mat{X}(z)$ has to result in
a problem that has been properly prepared for asymptotic analysis.  The transformation is best presented as a composition of several subsequent transformations:
\begin{itemize}
\item[1(a).] A transformation \eqref{eq:PtoQ} is introduced from $\mat{P}(z;N,k)$ to a
new unknown matrix $\mat{Q}(z;N,k)$ having the effect of moving poles
at some of the nodes in $X_N$ from the second column of
$\mat{P}(z;N,k)$ to the first column of $\mat{Q}(z;N,k)$.  This
transformation turns out to be necessary in our approach to take into
account subintervals of $[a,b]$ that are saturated with zeros of
$\pi_{N,k}(z)$ in the sense that there is a zero between each pair of
neighboring nodes (recall Proposition~\ref{prop:confine}).  The
saturated regions are not known in advance, but are detected by the
equilibrium measure (see step 1(c) below).
\item[1(b).] The matrix $\mat{Q}(z;N,k)$ is transformed into $\mat{R}(z)$,
a matrix that has, instead of polar singularities, a jump discontinuity across a contour in the complex $z$-plane along which $\mat{R}(z)$ takes continuous boundary values.  To see how a pole may be removed at the cost of a jump across a contour, consider a point $x_0$ at which a matrix function 
$\mat{M}(z)$ is meromorphic, 
having a simple pole in the second column such that for some given constant $w_0$:
\begin{equation}
\mathop{\rm Res}_{z=x_0}\mat{M}(z)=\lim_{z\rightarrow x_0}\mat{M}(z)\left(\begin{array}{cc} 0 & w_0 \\\\ 0 & 0\end{array}\right)\,.
\label{eq:rescond}
\end{equation}
If $f(z)$ is a scalar 
function analytic in the region $0<|z-x_0|<\epsilon$ for some $\epsilon>0$
having a simple pole at $x_0$ with
residue $w_0$ (obviously there are many such functions and consequently significant freedom in making a choice), then we may try to define a new matrix function $\mat{N}(z)$
by choosing some positive $\delta<\epsilon$ sufficiently small and setting
\begin{equation}
\mat{N}(z)=\left\{\begin{array}{ll}\mat{M}(z)\,,&\text{for $|z-x_0|>\delta$\,,}\\\\
\displaystyle\mat{M}(z)\left(\begin{array}{cc} 1 & -f(z) \\\\0 & 1
\end{array}\right)\,, &\text{for $|z-x_0|<\delta$\,.}
\end{array}\right.
\label{eq:poleremoval}
\end{equation}
It follows that the singularity of $\mat{N}(z)$ at $z=x_0$ is removable.
Therefore $\mat{N}(z)$ may be considered to be 
analytic in the region $|z-x_0|<\delta$, and also at each point of the region $|z-x_0|>\delta$ where additionally $\mat{M}(z)$ is known to be analytic.  In place of the residue condition
\eqref{eq:rescond}, we now have a known jump discontinuity across the circle $|z-x_0|=\delta$ along which $\mat{N}(z)$ takes continuous boundary values from the inside (denoted $\mat{N}_+(z)$) and the outside 
(denoted $\mat{N}_-(z)$):
\begin{equation}
\mat{N}_+(z)=\mat{N}_-(z)\left(\begin{array}{cc}1 & -f(z)
\\\\ 0 & 1\end{array}\right)\,,\hspace{0.2 in}\text{for $|z-x_0|=\delta$}\,.
\label{eq:jumpcondition}
\end{equation}
Obviously, the disc $|z-x_0|<\delta$ can be replaced by another domain
$D$ containing $x_0$.  This technique of removing poles was first
introduced in \cite{DeiftKKZ96}.  

The problem at hand is more complicated because the number of poles
grows in the limit of interest, and in this limit the poles accumulate
on a fixed set, and thus it is not feasible to encircle each with its
own circle of fixed size.  In \cite{KamvissisMM03} a generalization of
the technique described above was developed precisely to allow for the
simultaneous removal of a large number of poles in a way that is
asymptotically advantageous as the number of poles increases.  This
generalization employs a single function $f(z)$ with simple poles at
$x_{N,n}$ for $n=0,\dots,N-1$ having corresponding residues $w_{N,n}$,
and makes the change of variables \eqref{eq:poleremoval} in a common
domain $D$ containing all of the points $x_{N,0},\dots,x_{N,N-1}$.
The essential asymptotic analysis is then related to the nature of the
jump condition that generalizes \eqref{eq:jumpcondition} for $z$ on the
boundary of $D$.  This jump condition can have different asymptotic properties
in the limit $N\rightarrow\infty$ according to the placement of the boundary
of $D$ in the complex plane.

Further difficulties arise here because it turns out that the correct
location for the boundary of $D$ needed to facilitate the asymptotic
analysis in the limit $N\rightarrow\infty$ coincides in part with the
interval $[a,b]$ that contains the poles, and in the context of the
method of \cite{KamvissisMM03} this leads to singularities
in both the boundary values of the matrix unknown and also in the jump
matrix relating the boundary values.  These singularites are an
obstruction to further analysis.  Therefore, the transformation we
will introduce from $\mat{Q}(z;N,k)$ to $\mat{R}(z)$ uses a further
variation of the pole removal technique developed in \cite{Miller02}
in which two different residue interpolating functions $f_1(z)$ and
$f_2(z)$ are used in respective disjoint domains $D_1$ and $D_2$ such
that all of the poles $x_{N,n}$ are common boundary points of both
domains.  This version of the pole removal technique ultimately
enables subsequent detailed analysis in the neighborhood of the
interval $[a,b]$ in which $\mat{Q}(z;N,k)$ has poles.
\item[1(c).] $\mat{R}(z)$ is transformed into $\mat{S}(z)$ by a change of variables that is written explicitly in terms of the equilibrium measure.
The equilibrium measure is the solution of a variational problem of
logarithmic potential theory that is posed in terms of the functions
$\rho^0(x)$ and $V(x)$ given on the interval $[a,b]$ and the constant
$c\in (0,1)$.  The fundamental properties of the equilibrium measure
are well-known in general, and for particular cases of $\rho^0(x)$,
$V(x)$, and $c$, it is not difficult to calculate the equilibrium
measure explicitly.  The purpose of introducing the equilibrium
measure is that the variational problem it satisfies entails some
constraints that impose strict inequalities on variational
derivatives.  These variational derivatives ultimately appear in the
problem with a factor of $N$ in certain exponents, and the
inequalities lead to desirable exponential decay as
$N\rightarrow\infty$.

The technique of preparing a matrix Riemann-Hilbert problem for
subsequent asymptotic analysis with the introduction of an appropriate
equilibrium measure first appeared in the paper \cite{DeiftVZ97}, and
was subsequently applied to the computation of asymptotics for
orthogonal polynomials with continuous weights in
\cite{DKMVZstrong,DeiftKMVZ99}.  The key quantity in all of these
papers is the complex logarithmic potential of the equilibrium
measure, the so-called ``$g$-function''.  In order to apply these
methods in the discrete weights context, we need to modify the
relationship between the $g$-function and the equilibrium measure (see
\eqref{eq:rhodef} and \eqref{eq:gdef} below) to reflect the local
reversal of triangularity described in 1(a) above.  This amounts to a
further generalization of the technique introduced in
\cite{DeiftVZ97}.
\item[1(d).] The final transformation explicitly relates $\mat{S}(z)$ to a matrix $\mat{X}(z)$.  The matrix $\mat{S}(z)$ is apparently difficult to analyze in the neighborhood of subintervals of $[a,b]$ where constraints
in the variational problem are not active and consequently exponential decay is not obvious.  A model for this kind of situation is a matrix $\mat{M}(z)$
that takes continuous boundary values on an interval $I$ of the real axis from above (denoted $\mat{M}_+(z)$) and below (denoted $\mat{M}_-(z)$) that satisfy a jump relation of the form
\begin{equation}
\mat{M}_+(z)=\mat{M}_-(z)\left(\begin{array}{cc}e^{iN\theta(z)} & 1 \\\\
0 & e^{-iN\theta(z)}\end{array}\right)\,,
\end{equation}
where $\theta(z)$ is a real-analytic function that is strictly increasing
for $z\in I$.  This is therefore a rapidly oscillatory jump relation that has no obvious limit as $N\rightarrow\infty$  However, noting the algebraic factorization
\begin{equation}
\left(\begin{array}{cc} e^{iN\theta(z)} & 1 \\\\ 0 & e^{-iN\theta(z)}\end{array}\right)=
\left(\begin{array}{cc} 1 & 0 \\\\ e^{-iN\theta(z)}& 1\end{array}\right)\left(\begin{array}{cc} 0 & 1
\\\\ -1 & 0\end{array}\right)\left(\begin{array}{cc} 1 & 0\\\\
e^{iN\theta(z)} & 1\end{array}\right)\,,
\end{equation}
and using the analyticity of $\theta(z)$, we may choose some sufficiently small $\epsilon>0$ and define a new unknown by setting
\begin{equation}
\mat{N}(z):=\left\{\begin{array}{ll}
\displaystyle \mat{M}(z)\left(\begin{array}{cc} 1 & 0 \\\\
-e^{iN\theta(z)} & 1\end{array}\right)\,,
&\text{for $\Re(z)\in I$ and $0<\Im(z)<\epsilon$\,,}\\\\
\displaystyle \mat{M}(z)\left(\begin{array}{cc} 1 & 0 \\\\
e^{-iN\theta(z)} & 1 \end{array}\right)\,,&\text{for $\Re(z)\in I$ and
$-\epsilon<\Im(z)<0$\,,}\\\\
\mat{M}(z)\,, & \text{otherwise\,.}
\end{array}\right.
\end{equation}
The matrix $\mat{N}(z)$ has jump discontinuities along the three parallel contours $I$, $I+i\epsilon$, and $I-i\epsilon$.  If on any of these we indicate the boundary value taken by $\mat{N}(z)$ from above as $\mat{N}_+(z)$ and from below as $\mat{N}_-(z)$, then the oscillatory jump condition for $\mat{M}(z)$ in $I$ is replaced by the three different formulae:
\begin{equation}
\mat{N}_+(z)=\mat{N}_-(z)\left(\begin{array}{cc} 1 & 0 \\\\ e^{iN\theta(z)} & 1\end{array}\right)\,,\hspace{0.2 in}z\in I+i\epsilon\,,
\label{eq:N1}
\end{equation}
\begin{equation}
\mat{N}_+(z)=\mat{N}_-(z)\left(\begin{array}{cc} 1 & 0 \\\\ e^{-iN\theta(z)} & 1\end{array}\right)\,,\hspace{0.2 in}z\in I-i\epsilon\,,
\end{equation}
\begin{equation}
\mat{N}_+(z)=\mat{N}_-(z)\left(\begin{array}{cc} 0 & 1 \\\\ -1 & 0\end{array}\right)\,,\hspace{0.2 in}z\in I\,.
\label{eq:N3}
\end{equation}
The Cauchy-Riemann equations satisfied by $\theta(z)$ in $I$ imply that
$\Re(i\theta(z))$ is negative for $\Im(z)=\epsilon$ and positive for $\Im(z)=-\epsilon$.  Thus the jump conditions \eqref{eq:N1}--\eqref{eq:N3}
all have obvious asymptotics as $N\rightarrow\infty$.  

The replacement of an oscillatory jump matrix by an exponentially
decaying one on the basis of an algebraic factorization is the essence
of the steepest descent method for Riemann-Hilbert problems first
proposed in \cite{DeiftZ93}.  Our transformation from $\mat{S}(z)$ to
$\mat{X}(z)$ is based on this key idea, but involves more complicated
factorizations of both upper and lower triangular matrices.  
\end{itemize}

These three steps of our analysis of the matrix $\mat{P}(z;N,k)$
solving Interpolation Problem~\ref{rhp:DOP} will be carried out in
\S~\ref{sec:preparation} and \S~\ref{sec:asymptotics}.  

But first, we will present the results of our analysis.  The detailed
asymptotic behavior in the limit $N\rightarrow\infty$ of the discrete
orthogonal polynomials in overlapping sets that cover the entire
complex plane will be discussed in \S~\ref{sec:polytheorems}. 
After some important definitions and notation are
established in \S~\ref{sec:equilibrium} and
\S~\ref{sec:hyperelliptic}, the results themselves will be given in
in \S~\ref{sec:actualtheorems}.  In \S~\ref{sec:examplepolys} we show
how the general theory applies in some classical cases, specifically
the Krawtchouk polynomials and two types of polynomials in the Hahn
family.  The equilibrium measures for the Hahn polynomials are
also described here in Theorem~\ref{theorem:hahn}.  

Further results of our analysis in the context of statistical
ensembles associated with families of discrete orthogonal polynomials
are discussed in \S~\ref{sec:DOPensembles}.  First, we introduce the
notion of a discrete orthogonal polynomial ensemble in
\S~\ref{sec:DOPensemblesintro}, and describe the ensembles associated
with dual polynomials in \S~\ref{sec:DOPensemblesholes}.  In
\S~\ref{sec:tiling}, we discuss rhombus tilings of a hexagon as a
specific application of discrete orthogonal polynomial ensembles and
their duals.  Our general results on universality of various
statistics in the limit $N\rightarrow\infty$ are explained in
\S~\ref{sec:universalityactualtheorems}.  The specific results implied
by the general ones in the context of the hexagon tiling problem are
described in \S~\ref{sec:Hextheorems}.

As mentioned above, \S~\ref{sec:preparation} and
\S~\ref{sec:asymptotics} contain the complete asymptotic analysis of
the matrix $\mat{P}(z;N,k)$ in the limit $N\rightarrow\infty$.  This
analysis is then used in \S~\ref{sec:asymptoticspi} to establish the
results presented in \S~\ref{sec:actualtheorems}, and then is used again in
\S~\ref{sec:universality} to establish the results presented in
\S~\ref{sec:universalityactualtheorems}.

For those asymptotic results given here that correspond to theorems
already stated in our announcement \cite{BaikKMM03}, we generally obtain
significantly sharper error estimates.  Since we published our
announcement, we have learned how to circumvent certain technical
difficulties related to the continuum limit of the discrete
orthogonality measures and the possibility of transition points where
triangularity of residue matrices changes abruptly.  In our opinion,
these technical innovations do more than make the error estimates sharper; 
they also make the proofs more elegant.

\subsection{Acknowledgements.}
J. Baik was supported in part by the National Science Foundation under
grant DMS-0208577.  T. Kriecherbauer was supported in part by the
Deutsche Forschungsgemeinschaft under grant SFB/TR 12.
K. T.-R. McLaughlin was supported in part by the National Science
Foundation under grant numbers DMS-9970328 and DMS-0200749.
K. T.-R. McLaughlin wishes to thank T. Paul, F. Golse, and the staff
of the \'Ecole Normal Superieur, Paris, for their kind hospitality.
P. D. Miller was supported in part by the National Science Foundation
under grant DMS-0103909, and by a grant from the Alfred P. Sloan
Foundation.  

We wish to thank several individuals for their comments: A. Borodin,
P. Deift, M. Ismail, K. Johansson, and A. Kuijlaars.  We also thank
J. Propp for providing us with Figure~\ref{fig:Propp}.

\section{Asymptotics of General Discrete Orthogonal Polynomials in the Complex Plane}
\label{sec:polytheorems}
In this section we state our results on the asymptotics of the
discrete orthogonal polynomials subject to the conditions described in
the introduction, in the limit $N\rightarrow\infty$.  The asymptotic
formulae we will present characterize the polynomials in terms of an
equilibrium measure (described below in \S~\ref{sec:equilibrium}) and
also the function theory of a hyperelliptic Riemann surface associated
with the equilibrium measure (described below in
\S~\ref{sec:hyperelliptic}).  The results themselves will follow in 
\S~\ref{sec:actualtheorems}.
\subsection{The equilibrium energy problem.}
\label{sec:equilibrium}

\subsubsection{The equilibrium measure.}
It has been recognized for some time (see
\cite{Rakhmanov96,DragnevS97} as well as the review article \cite{KuijlaarsR98}) that the asymptotic
behavior of discrete orthogonal polynomials in the limit $N\rightarrow\infty$ with $k/N\rightarrow c\in (0,1)$, and in
particular the distribution of zeros in $(a,b)$, is related to a
constrained equilibrium problem for logarithmic potentials in a field
$\varphi(x)$ given by the formula \label{symbol:varphi}
\begin{equation}
\varphi(x):=V(x)+\int_a^b\log|x-y|\rho^0(y)\,dy
\label{eq:fielddef}
\end{equation}
for $x\in (a,b)$.  Under our assumptions on the weights, we can also
view $\varphi(x)$ as being defined via a continuum limit:
\begin{equation}
\varphi(x):=-\lim_{N\rightarrow\infty}\frac{\log(w_{N,n})}{N}
\label{eq:fieldlim}
\end{equation}
where $w_{N,n}$ is expressed in terms of $x_{N,n}$ which in turn is
identified with $x$.  Eliminating $\varphi(x)$ between
\eqref{eq:fielddef} and \eqref{eq:fieldlim} amounts to a more general
version of the limiting statement \eqref{eq:clweight}.  The external
field $\varphi(x)$ we need here is analogous to the continuum limit of
that usually encountered in the potential theory of orthogonal
polynomials
\cite{S91}.

Here, the field $\varphi(x)$ is a
real-analytic function in the open interval $(a,b)$ because $V(x)$ and
$\rho^0(x)$ are (by assumption) real-analytic functions in a neighborhood of $[a,b]$.
Unlike $V(x)$ and $\rho^0(x)$, however, the field $\varphi(x)$ does
not extend analytically beyond the endpoints of $(a,b)$ due to the
condition \eqref{eq:rho0nonzero}.

Given $c\in (0,1)$ and $\varphi(x)$ as above, consider the quadratic
functional \label{symbol:Ec}
\begin{equation}
E_c[\mu]:=c\int_a^b\int_a^b\log\frac{1}{|x-y|}
\, d\mu(x)\,d\mu(y)+ \int_a^b \varphi(x)\,d\mu(x)
\label{eq:energy}
\end{equation}
of Borel measures $\mu$ on $[a,b]$.  The subscript denotes the dependence
of the energy functional on the parameter $c$.
Let $\mu^c_{\rm min}$ \label{symbol:minimizer} be the measure that minimizes $E_c[\mu]$ over the
class of measures satisfying the upper and lower constraints
\begin{equation}
0\le \int_{x\in{\mathcal B}}d\mu(x)\le\frac{1}{c}\int_{x\in{\mathcal
B}}\rho^0(x)\,dx
\label{eq:constraints}
\end{equation}
for all Borel sets ${\mathcal B}\subset [a,b]$, and the normalization
condition
\begin{equation}
\int_a^bd\mu(x)=1\,.
\label{eq:Lagrange}
\end{equation}
The superscript on the minimizer indicates the value of the parameter
$c$ for which the energy functional \eqref{eq:energy} is minimized.
The existence of a unique minimizer under the conditions enumerated in
\S~\ref{sec:C1} and \S~\ref{sec:C2} follows from the Gauss-Frostman
Theorem; see \cite{ST} and \cite{DragnevS97} for details.  We will
refer to the minimizer as the {\em equilibrium measure}.

That a variational problem plays a central role in asymptotic behavior
is a familiar theme in the theory of orthogonal polynomials.  The key
new feature contributed by discreteness is the appearance of the upper
constraint on the equilibrium measure ({\em i.e.} the upper bound in
\eqref{eq:constraints}).  Since the equilibrium measure gives the distribution of zeros of $\pi_{N,k}(z)$ in $[a,b]$, the upper constraint 
can be traced to the
exclusion principle for zeros described in Proposition~\ref{prop:confine}.

The theory of the ``doubly constrained'' variational problem we are
considering is well-established.  In particular, the analytic
properties we assume of $V(x)$ and $\rho^0(x)$ turn out to be unnecessary for
the mere existence of the equilibrium measure.
However, analyticity of $V(x)$ and $\rho^{0}(x)$ provides additional
regularity that we wish to exploit.  In particular, we have
the following result from the paper \cite{Kuijlaars00}.
\begin{prop}[Kuijlaars]
Let $V(x)$ and $\rho^0(x)$ be functions analytic in a complex neighborhood
of $[a,b]$ with $\rho^{0}(x)>0$ in $[a,b]$.  Then, the equilibrium
measure $\mu_{\rm min}^c$ is continuously
differentiable with respect to $x\in (a,b)$.  Moreover, the derivative
$d\mu_{\rm min}^c/dx$ is piecewise analytic, with a finite number of
points of nonanalyticity that may not occur at any $x$ where both
(strict) inequalities $d\mu_{\rm min}^c/dx(x)>0$ and $d\mu_{\rm
min}^c/dx(x)<\rho^0(x)/c$ hold.
\end{prop}

At a formal level, finding $\mu$ minimizing $E_c[\mu]$ subject to the constraint \eqref{eq:Lagrange} may be viewed as seeking a critical point
of the modified functional
\begin{equation}
F_c[\mu]:=E_c[\mu] -\ell_c\int_a^bd\mu(x)
\end{equation}
\label{symbol:Fc}
where $\ell_c$ \label{symbol:ellc} is a Lagrange multiplier.  When $\mu=\mu_{\rm min}^c$
and $\ell_c$ is an appropriate associated real constant, variations of $F_c[\mu]$
vanish in subsets of $[a,b]$ where neither the upper nor the lower constraints are active.  The Lagrange multiplier $\ell_c$ (the subscript
indicates the dependence on the parameter $c$) is known
in potential theory as the {\em Robin constant}.

\subsubsection{Simplifying assumptions on the equilibrium measure.}
\label{sec:C3}
For simplicity of exposition we want to exclude certain nongeneric
phenomena that may occur even under conditions of analyticity of
$V(x)$ and $\rho^0(x)$.  Let $\underline{\mathcal F}\subset [a,b]$ \label{symbol:underF} denote the closed set of
$x$-values where
$d\mu_{\rm min}^c/dx(x)=0$, and let $\overline{\mathcal F}\subset [a,b]$ \label{symbol:overF} denote
the closed set of $x$-values where $d\mu_{\rm min}^c/dx(x)=\rho^0(x)/c$.
We will make the following assumptions:
\begin{enumerate}
\item
Each connected component of $\underline{\mathcal F}$ and
$\overline{\mathcal F}$ has a nonempty interior.  Therefore
$\underline{\mathcal F}$ and $\overline{\mathcal F}$ are both finite unions
of closed intervals with each interval containing more than one point.
Note that this does not exclude the possibility of either $\underline{\mathcal F}$ or $\overline{\mathcal F}$ being empty.
\item
For each open subinterval $U$ of $(a,b)\setminus (\underline{\mathcal F}\cup
\overline{\mathcal F})$ and each limit point $z_0\in\underline{\mathcal F}$
of $U$, we have
\begin{equation}
\lim_{x\rightarrow z_0, x\in U}\frac{1}{\sqrt{|x-z_0|}}
\frac{d\mu_{\rm min}^c}{dx}(x) = K\hspace{0.2 in}\mbox{with $0<K<\infty$}
\label{eq:rootlower}
\end{equation}
and for each limit point $z_0\in\overline{\mathcal F}$ of $U$, we have
\begin{equation}
\lim_{x\rightarrow z_0, x\in U}\frac{1}{\sqrt{|x-z_0|}}
\left[\frac{1}{c}\rho^0(x)-\frac{d\mu_{\rm min}^c}{dx}(x)\right] =
K\hspace{0.2 in}\mbox{with $0<K<\infty$}\,.
\label{eq:rootupper}
\end{equation}
Therefore the derivative of the equilibrium measure meets each constraint
exactly like a square root.
\item
A constraint is active at each endpoint: $\{a,b\}\subset
\underline{\mathcal F}\cup\overline{\mathcal F}$.
\end{enumerate}
It is difficult to translate these conditions on $\mu_{\rm min}^c$
into sufficient conditions on $c$, $V(x)$, and $\rho^0(x)$.  However,
there is a sense in which the conditions above are satisfied
generically.  By genericity, we mean that given $V(x)$ and
$\rho^{0}(x)$, the set of values of $c$ for which the conditions fail
is discrete.  For the analogous problem in the continuous weights
case, the conditions 1 and 2 above are proved to be generic in
\cite{KuijlaarsM00} and we expect the same of the discrete weights.
For further arguments supporting the claim of the generic nature of
these two conditions, see the discussion at the beginning of
\S~\ref{sec:Airy}.  However, the condition 3 above has no counterpart
in continuous weight cases.  Nevertheless, for the classical discrete
weights of the Krawtchouk and Hahn classes this condition indeed holds
for all but a finite number of $c\in (0,1)$.

\begin{remark}
Relaxing the condition that a constraint should be active at each
endpoint requires specific local analysis near these two points.  We
expect that a constraint being active at each endpoint is a generic
phenomenon in the sense that the opposite situation occurs only for
isolated values of $c$.  We know this statement to be true in all
relevant classical cases.  For the Krawtchouk polynomials only the
values $c=p$ or $c=q=1-p$ correspond to an equilibrium measure that is
not constrained at both endpoints (see
\cite{DragnevS00} and \S~\ref{sec:Krawtchouk}).  The situation is similar
for the Hahn polynomials, where only the values $c=c_A$ and $c=c_B$ defined by (\ref{eq:cAcB}) are exceptional (see \S~\ref{sec:Hahn}).  
While the values of $c$ for which no constraint is active at an endpoint
of $(a,b)$ are exceptional, the behavior of the equilibrium measure 
near that endpoint at the exceptional values of $c$ is again, in a sense,
generic.  In particular, we have the following result.
\begin{prop}
Suppose that $\rho^0(x)$ and $V(x)$ are analytic functions for $x\in
(a,b)$ having, along with $V'(x)$, continuous extensions to $[a,b]$.
Suppose also that  one of the endpoints $a$
or $b$ is not contained in $\underline{\cal F}\cup\overline{\cal F}$.  Then if $\rho^0(x)$ and $d\mu_{\rm
min}^c/dx$ are H\"older continuous at this endpoint with exponent
$\nu>0$, then as $x\in (a,b)$ tends toward the endpoint,
\begin{equation}
\rho^0(x)-2c\frac{d\mu_{\rm min}^c}{dx}(x)\rightarrow 0\,.
\end{equation}
Thus, if neither constraint is active at an endpoint of $(a,b)$, then
at that endpoint the equilibrium measure takes on the average value of the upper
and lower constraints.
\label{prop:average}
\end{prop}
\begin{proof}
Since for some neighborhood $U$ of the endpoint no constraint is active 
in $(a,b)\cap U$, an Euler-Lagrange derivative of $E_c[\mu]$ satisfies
a variational equilibrium condition (see \eqref{eq:equilibrium}) at each
point of $(a,b)\cap U$.
Differentiating this condition with
respect to $x$ yields
\begin{equation}
\mbox{P. V.}\int_a^b \frac{f(y)\,dy}{x-y} + V'(x)\equiv 0\hspace{0.2 in}\mbox{where}
\hspace{0.2 in}f(x):=\rho^0(x)-2c\frac{d\mu_{\rm min}^c}{dx}(x)
\end{equation}
which holds for all $x$ in the interior of the subinterval.  Removing the
singularity from the principal value integral gives the identity
\begin{equation}
f(x)\log\left(\frac{x-a}{b-x}\right) + \int_a^b\frac{f(y)-f(x)}{x-y}\,dy
+ V'(x)\equiv 0\,,
\label{eq:regularized}
\end{equation}
where the integral is nonsingular by virtue of the H\"older condition,
and is uniformly bounded.  If $K>0$ is the H\"older constant for $f$,
then we have
\begin{equation}
\begin{array}{rcl}
\displaystyle \left|\int_a^b\frac{f(y)-f(x)}{x-y}\,dy\right| & \le &
\displaystyle \int_a^b\left|\frac{f(y)-f(x)}{x-y}\right|\,dy\\\\
&\le &\displaystyle
K\int_a^b |x-y|^{\nu-1}\,dy\\\\
&=&\displaystyle \frac{K}{\nu}\left[
(x-a)^\nu+(b-x)^\nu\right]\\\\
&\le &\displaystyle \frac{K(b-a)^\nu}{2^{\nu-1}\nu}\,.
\end{array}
\end{equation}
Letting $x$ tend toward the endpoint of interest in
(\ref{eq:regularized}), all terms but the first one on the left-hand
side therefore remain bounded.  Hence it is necessary that the first
term involving the logarithm be bounded as well; by the H\"older
condition satisfied by $f(x)$ at the endpoint it follows that $f(x)$
tends to zero as $x$ tends toward the endpoint.
\end{proof}
This fact provides the key to modifications of the analysis we will
present in \S~\ref{sec:asymptotics} that are
necessary to handle the exceptional values of $c\in (0,1)$.  We will
describe these modifications in a subsequent paper.
%In any case, the nongeneric situation can
%be handled as well, and in 
%\S~\ref{sec:conclusion} we will describe the
%behavior of the minimizer near the endpoints in this case
%that will be necessary for the local analysis at endpoints needed to
%generalize our analysis.  For the exceptional values of $c$ one expects to %obtain a local asymptotic description of the polynomials, uniformly valid %near the unconstrained endpoints of $[a,b]$, in terms of novel special %functions.
\end{remark}

\subsubsection{Voids, bands, and saturated regions.}
Under the conditions enumerated in \S~\ref{sec:C1}, \S~\ref{sec:C2}, and
\S~\ref{sec:C3}, the equilibrium measure $\mu^c_{\rm min}$ partitions
$(a,b)$ into three kinds of subintervals, a finite number of each, and
each having a nonempty interior.  The three types are defined as follows:
\begin{definition}[Voids]
A void $\Gamma$ is an open subinterval of $[a,b]$ of maximal length in
which $\mu^c_{\rm min}(x)\equiv 0$, and thus the equilibrium measure realizes
the lower constraint.  
\end{definition}
\begin{definition}[Bands]
A band $I$ is an open subinterval of $[a,b]$ of maximal length where
$\mu^c_{\rm min}(x)$ is a measure with a real-analytic density
satisfying $0<d\mu^c_{\rm min}/dx<\rho^0(x)/c$.
\end{definition}
\begin{definition}[Saturated regions]
A saturated region $\Gamma$ is an open subinterval of $[a,b]$ of
maximal length in which $d\mu^c_{\rm min}/dx\equiv \rho^0(x)/c$, and
thus the equilibrium measure realizes the upper constraint.  
\end{definition}
Voids and saturated regions will also be called {\em gaps} when it is
not necessary to distinguish between these two types of intervals.
The closure of the union of all subintervals of the three types
defined above is the interval $[a,b]$.  From condition 1 in
\S~\ref{sec:C3} above, bands cannot be adjacent to each other, and from condition 3 in \S~\ref{sec:C3} a band may not be adjacent to an endpoint
of $[a,b]$.  Thus subject to our assumptions, a band always has on each side either
a void or a saturated region, and the equilibrium measure thus determines a set of numbers in $(a,b)$ \label{symbol:G} 
\label{symbol:alphas} \label{symbol:betas}
\begin{equation}
a<\alpha_0<\beta_0<\alpha_1<\beta_1<\cdots <\alpha_G<\beta_G<b
\end{equation}
that are the endpoints of the bands.  Thus the bands are open intervals
of the form \label{symbol:bands}
\begin{equation}
I_j:=(\alpha_j,\beta_j)\,,\hspace{0.2 in}\text{for $j=0,\dots,G$\,.}
\end{equation}
The corresponding gaps are
the intervals $(a,\alpha_0)$, $(\beta_G,b)$, which we refer to as the {\em exterior gaps}, and \label{symbol:gaps}
\begin{equation}
\Gamma_j:=(\beta_{j-1},\alpha_j)\,,\hspace{0.2 in}\text{for $j=1,\dots,G$\,,}
\end{equation}
which we refer to as the {\em interior gaps}.

\subsubsection{Quantities derived from the equilibrium measure.}

The variational derivative of $E_c[\mu]$ evaluated on the equilibrium measure $\mu=\mu_{\rm min}^c$ is given by \label{symbol:varderiv}
\begin{equation}
\frac{\delta E_c}{\delta\mu}(x):=-2c\int_a^b\log|x-y|\,d\mu_{\rm min}^c(y) +
\varphi(x)\,,
\label{eq:variationalderivative}
\end{equation}
and we may define an analytic logarithmic potential of the equilibrium measure by the formula \label{symbol:Lc(z)}
\begin{equation}
L_c(z):=c\int_a^b\log(z-x)\,d\mu_{\rm min}^c(x)\,,\hspace{0.2 in}\text{for $z\in\mathbb{C}\setminus (-\infty,b]$\,.}
\label{eq:Lcdef}
\end{equation}
From any gap $\Gamma$ we may introduce a function $\overline{L}_c^\Gamma(z)$ \label{symbol:LcbarGamma} analytic for $z$ with $\Re(z)\in\Gamma$ and $|\Im(z)|$ sufficiently small that satisfies
\begin{equation}
\overline{L}_c^\Gamma(z)=c\int_a^b\log|z-x|\,d\mu_{\rm min}^c(x)\,,\hspace{0.2 in}\text{for $z\in\Gamma$\,.}
\label{eq:LcbarGammadef}
\end{equation}
And from any band $I$ we may introduce a function $\overline{L}_c^I(z)$
\label{symbol:LcbarI} analytic for $z$ with $\Re(z)\in I$ and $\|\Im(z)|$
sufficiently small that satisfies
\begin{equation}
\overline{L}_c^I(z)=c\int_a^b\log|z-x|\,d\mu_{\rm min}^c(x)\,,
\hspace{0.2 in}\text{for $z\in I$\,.}
\label{eq:LcbarI}
\end{equation}

Recall the Lagrange multiplier $\ell_c$.  If $\Gamma$ is a void,
then admissible variations of $\mu_{\rm min}^c$ are positive, and
a simple variational calculation shows that and for $x\in\Gamma$ we have the strict inequality
\begin{equation}
\frac{\delta E_c}{\delta\mu}(x)>\ell_c\,.
\label{eq:voidinequality}
\end{equation}
Thus, for each void $\Gamma$ we may introduce a positive function
having an analytic extension from the interior: \label{symbol:xiGamma}
\begin{equation}
\xi_\Gamma(x):=\frac{\delta E_c}{\delta\mu}(x)-\ell_c\,.
\end{equation}

In a band $I$, variations of the equilibrium measure in $I$ are free (may be of either
sign).  Thus, for $x\in I$ we have the equilibrium condition
\begin{equation}
\frac{\delta E_c}{\delta\mu}(x)\equiv\ell_c\,.
\label{eq:equilibrium}
\end{equation}
For each band $I$ we may introduce two positive functions having analytic extensions
from the interior: \label{symbol:psiI}\label{symbol:psiIbar}
\begin{equation}
\psi_I(x):=\frac{d\mu_{\rm min}^c}{dx}(x)\hspace{0.2 in}\mbox{and}
\hspace{0.2 in}
\overline{\psi}_I(x):=\frac{1}{c}\rho^0(x)-\frac{d\mu_{\rm min}^c}{dx}(x)\,.
\label{eq:psiIpsiIbar}
\end{equation}

If $\Gamma$ is a saturated region, then variations of the equilibrium measure in $\Gamma$ are strictly negative, and for $x\in\Gamma$ we
have the strict variational inequality
\begin{equation}
\frac{\delta E_c}{\delta\mu}(x)<\ell_c\,.
\label{eq:saturatedregioninequality}
\end{equation}
It follows that for each saturated region $\Gamma$ we may introduce a
positive function having an analytic extension from the interior:
\begin{equation}
\xi_\Gamma(x):=\ell_c-\frac{\delta E_c}{\delta\mu}(x)\,.
\end{equation}

In addition to the functions $\overline{L}_c^\Gamma(x)$ and
$\xi_\Gamma(x)$ that extend analytically from each gap $\Gamma$, and
the functions $\overline{L}_c^I(x)$, $\psi_I(x)$, and
$\overline{\psi}_I(x)$ that extend analytically from each band $I$, we
may define for each band endpoint a function analytic in a
neighborhood of this point.  If $z=\alpha$ is a left band edge
separating a void $\Gamma$ (for real $z<\alpha$) from a band $I$ (for
real $z>\alpha$), then according to the generic assumption
\eqref{eq:rootlower} in \S~\ref{sec:C3}, the function defined by
\begin{equation}
\tau_\Gamma^{\nab,L}(z):=\left(2\pi Nc\int_\alpha^z \psi_I(x)\,dx\right)^{2/3}\,,\hspace{0.2 in}\text{$\tau_\Gamma^{\nab,L}(z)>0$ for $z>\alpha$\,,}
\label{eq:zetadefAL}
\end{equation}
\label{symbol:taunabL}
extends to a neighborhood of $z=\alpha$ as an invertible conformal
mapping.  If $z=\beta$ is a right band edge separating a void $\Gamma$
(for real $z>\beta$) from a band $I$ (for real $z<\beta$), then
\eqref{eq:rootlower} implies that the function defined by 
\begin{equation}
\tau_\Gamma^{\nab,R}(z):=\left(-2\pi Nc\int_\beta^z\psi_I(x)\,dx\right)^{2/3}
\,,\hspace{0.2 in}\text{$\tau_\Gamma^{\nab,R}(z)>0$ for $z<\beta$\,,}
\label{eq:zetadefAR}
\end{equation}
\label{symbol:taunabR}
extends to a neighborhood of $z=\beta$ as an invertible conformal mapping.
If $z=\alpha$ is a left band edge separating a saturated region $\Gamma$
(for real $z<\alpha$) from a band $I$ (for real $z>\alpha$), then according
to the generic assumption \eqref{eq:rootupper} in \S~\ref{sec:C3}, the function
defined by
\begin{equation}
\tau_\Gamma^{\del,L}(z):=\left(2\pi Nc\int_\alpha^z \overline{\psi}_I(x)\,dx
\right)^{2/3}\,,\hspace{0.2 in}\text{$\tau_\Gamma^{\del,L}(z)>0$ for $z>\alpha$\,,}
\label{eq:zetadefBL}
\end{equation}
\label{symbol:taudelL}
extends to a neighborhood of $z=\alpha$ as an invertible conformal mapping.
If $z=\beta$ is a right band edge separating a saturated region $\Gamma$ (for real $z>\beta$) from a band $I$ (for real $z<\beta$), then \eqref{eq:rootupper}
implies that the function defined by
\begin{equation}
\tau_\Gamma^{\del,R}(z):=\left(-2\pi Nc\int_\beta^z\overline{\psi}_I(x)\,dx\right)^{2/3}\,,\hspace{0.2 in}\text{$\tau_\Gamma^{\del,R}(z)>0$ for $z<\beta$\,,}
\label{eq:zetadefBR}
\end{equation}
\label{symbol:taudelR}
extends to a neighborhood of $z=\beta$ as an invertible conformal mapping.

For later use it will be useful to define a real constant
$\theta_\Gamma$ corresponding to each gap.  For each void
$\Gamma=\Gamma_j$ surrounded by bands $I_{j-1}$ and $I_j$, we define a
constant by
\label{symbol:thetaconst}
\begin{equation}
\theta_{\Gamma_j}:=-2\pi c\int_{\alpha_j}^b \frac{d\mu_{\rm min}^c}{dx}(x)\,dx\,.
\label{eq:thetaGammajvoid}
\end{equation}
Similarly, for each saturated region $\Gamma=\Gamma_j$ surrounded by
bands $I_{j-1}$ and $I_j$, we define a constant by
\begin{equation}
\theta_{\Gamma_j}:=2\pi c\int_{\alpha_j}^b\left[\frac{\rho^0(x)}{c}-
\frac{d\mu_{\rm min}^c}{dx}(x)\right]\,dx\,.
\label{eq:thetaGammajsat}
\end{equation}
There is no difficulty using the same symbol $\theta_\Gamma$ on the
left-hand side of these two definitions, because a given gap
$\Gamma=\Gamma_j$ is either a void or a saturated region, but cannot
be both. For the gaps $\Gamma=(a,\alpha_0)$ and $\Gamma=(\beta_G,b)$
that exist according to the genericity assumptions stated in
\S~\ref{sec:C3}, it will also be useful to define associated constants
$\theta_\Gamma$.  Whether each of these gaps is a void or a saturated
region, we define
\begin{equation}
\theta_{(a,\alpha_0)}:=-2\pi c\,,\hspace{0.2 in}\text{and}\hspace{0.2 in}
\theta_{(\beta_G,b)}:=0\,.
\label{eq:thetaGammaleftright}
\end{equation}
\label{symbol:thetaGammaleftright}

\subsubsection{Dual equilibrium measures.}
According to (\ref{eq:dualweightsdefine}), if $V_N(x)$ is
associated with the weights $\{w_{N,n}\}$ and if $\overline{V}_N(x)$
\label{symbol:VbarN}
is associated with the dual weights $\{\overline{w}_{N,n}\}$, then we
have the simple identity $\overline{V}_N(x)=-V_N(x)$.  
This leads to the useful fact that knowing the equilibrium measure for one
family of discrete orthogonal polynomials is equivalent to knowing the
equilibrium measure for the dual discrete orthogonal polynomials.
We have the following specific result.
\begin{prop}
Let $E_c[\mu;V,\rho^0]$ denote the energy functional
(\ref{eq:energy}) with external field $\varphi(x)$ given in terms
of analytic functions $V(x)$ and $\rho^0(x)$ by
(\ref{eq:fielddef}), and let $\ell_c[V,\rho^0]$ denote the corresponding Lagrange multiplier.  Let $P(c,V,\rho^0)$ denote the problem of
finding the measure $\mu$ on $(a,b)$ minimizing
$E_c[\mu;V,\rho^0]$, subject to the conditions
(\ref{eq:constraints}) and (\ref{eq:Lagrange}).  If for all $c\in
(0,1)$, $\mu_{\rm min}^c$ is the solution of the problem
$P(c,V,\rho^0)$, then the measure $\bar{\mu}_{\rm min}^{1-c}$ with
density
\begin{equation}
\frac{d\bar{\mu}_{\rm min}^{1-c}}{dx}(x):=
\frac{1}{1-c}\left(\rho^0(x)-c\frac{d\mu_{\rm
min}^{c}}{dx}(x)\right)
\label{eq:taudef}
\end{equation}
is the solution of the problem $P(1-c,-V,\rho^0)$. Also
\begin{equation}\label{eq:dualVari}
  \frac{\delta E_{1-c}[\bar{\mu};-V,\rho^0]}{\delta \bar{\mu}} 
\biggl|_{\bar{\mu}_{\rm min}^{1-c}}
  = - \frac{\delta E_{c}[\mu;V,\rho^0]}{\delta \mu} \biggl|_{\mu_{\rm
  min}^{c}},
  \qquad \ell_{1-c}[-V,\rho^0] = -\ell_c[V,\rho^0]\,.
\end{equation}
\label{prop:flip}
\end{prop}
\label{symbol:mubar}
\begin{proof}
Clearly, the measure with density (\ref{eq:taudef}) satisfies both
conditions (\ref{eq:constraints}) and (\ref{eq:Lagrange}).  A
direct calculation then shows that $E_{1-c}[\bar{\mu};-V,\rho^0]$, when
considered as a functional of $\mu$ by the relation
\begin{equation}
\frac{d\bar{\mu}}{dx}(x)=\frac{1}{1-c}\left(\rho^0(x)-c\frac{d\mu}{dx}(x)\right)
\end{equation}
is linearly related to the functional $E_{c}[\mu;V,\rho^0]$:
\begin{equation}
E_{1-c}[\bar{\mu};-V,\rho^0]=\frac{c}{1-c}E_{c}[\mu;V,\rho^0]-\frac{1}{1-c}\int_a^b
V(x)\rho^0(x)\,dx\,.
\end{equation}
Since $c$ and $1-c$ are both positive, we obtain
\eqref{eq:taudef}. The proof of \eqref{eq:dualVari} is 
similar.
\end{proof}

\subsection{Elements of hyperelliptic function theory.}
\label{sec:hyperelliptic}
Let the analytic function $R(z)$ be defined for $z\in\cx\setminus\cup_jI_j$
to satisfy \label{symbol:sqrt}
\begin{equation}
R(z)^2 = \prod_{j=0}^G(z-\alpha_j)(z-\beta_j)\,,\hspace{0.2 in}\text{and}
\hspace{0.2 in} \text{$R(z)\sim z^{G+1}$ as $z\rightarrow\infty$\,,}
\label{eq:Rofzdefine}
\end{equation}
and for $z$ in the same domain define
\begin{equation}
h'(z):=\frac{1}{2\pi iR(z)}\int_{\cup_jI_j}\frac{\eta'(x)R_+(x)}{x-z}\,dx
+ \frac{1}{R(z)}\left[\kappa z^G + \sum_{m=0}^{G-1}f_mz^m\right]\,,
\label{eq:hprime}
\end{equation}
where $R_+(x)$ denotes the boundary value taken by $R(z)$ from the upper half-plane, and where the constants $f_m$, $m=0,\dots,G-1$, are chosen (uniquely, see Appendix~\ref{sec:thetasolve}) so that 
\begin{equation}
\int_{\Gamma_j}h'(z)\,dz = 0\,,\hspace{0.2 in}\text{for $j=1,\dots,G$\,.}
\label{eq:intszero}
\end{equation}
Then, we define a function for $z\in\cx\setminus (-\infty,\beta_G]$ by
the integral \label{symbol:h(z)}
\begin{equation}
h(z):=\kappa\log(z) + \int_z^\infty\left[\frac{\kappa}{s}-h'(s)\right]\,ds
\label{eq:hformula}
\end{equation}
where the path of integration lies in $\cx\setminus (-\infty,\beta_G]$.
Furthermore, we define a constant $\gamma$ by \label{symbol:gamma}
\begin{equation}
\gamma:=\eta(\beta_G)-2h(\beta_G)\,.
\label{eq:gammadefine}
\end{equation}
The combination $N\ell_c+\gamma$ plays an important role in what
follows.  Since $\gamma$ remains bounded as $N\rightarrow\infty$, we
may interpret $\gamma$ as a higher-order correction to the scaled
Robin constant $N\ell_c$.

It may be verified that for $z\in\Gamma_j$, the difference of boundary values taken by $h(z)$ depends on $j$ but is independent of $z$.  Thus there
are constants $c_j$, $j=1,\dots,G$, such that
\begin{equation}
h_+(z)-h_-(z):=\lim_{\epsilon\downarrow 0} h(z+i\epsilon)-h(z-i\epsilon) = ic_j\,,\hspace{0.2 in}\text{for $z\in\Gamma_j$\,.}
\label{eq:csdef}
\end{equation}
Moreover, it can be checked directly that the constants $c_j$ are real and linear in $\kappa$ so that we may write
\begin{equation}
c_j=c_j^{(0)}+\omega_j\kappa
\end{equation}
for some other real constants $c_j^{(0)}$ \label{symbol:cjzero} and $\omega_j$ \label{symbol:omegaj} that are independent of
$\kappa$.  We define a vector $\mat{r}$ with components \label{symbol:vectorr}
\begin{equation}
r_j:=N\theta_{\Gamma_j}-c_j^{(0)}
\label{eq:vectorr}
\end{equation}
for $j=1,\dots,G$, 
and a vector $\mat{\Omega}$ \label{symbol:vectorOmega} with components $\omega_j$ for $j=1,\dots,G$.

The function $iR_+(z)$ may be analytically continued from any band $I$ to
the complex plane with the semi-infinite intervals $(-\infty,\alpha_0]$ and
$[\beta_G,\infty)$ and the closures of the gaps $\Gamma_1,\dots,\Gamma_G$
deleted.  We call this analytic continuation $y(z)$, \label{symbol:y(z)} and for $z$ in this domain of definition, we introduce a vector function $\mat{m}(z)$ \label{symbol:vecm(z)} having components
$m_p(z):=z^{p-1}/y(z)$ for $p=1,2,\dots,G$.  Next, a constant $G\times G$ matrix $\mat{A}=(\mat{a}^{(1)},\mat{a}^{(2)},\dots,\mat{a}^{(G)})$ \label{symbol:Amatrix} is defined by insisting that the linear equations
\begin{equation}
\mat{A}\int_{\beta_{j-1}}^{\alpha_j}\mat{m}_-(z)\,dz = \pi i\mat{e}^{(j)}
\,,\hspace{0.2 in}\text{for $j=1,\dots,G$\,,}
\label{eq:Adetermine}
\end{equation}
are satisfied where $\mat{m}_-(z)$ denotes the boundary value taken on the real axis from the lower half-plane, and where $\mat{e}^{(j)}$ is column $j$ of the $G\times G$ identity matrix.  This determines vectors $\mat{b}^{(j)}$ by
the definition
\begin{equation}
\mat{b}^{(j)}:=-2\mat{A}\sum_{m=1}^j\int_{\alpha_{m-1}}^{\beta_{m-1}}
\mat{m}(z)\,dz\,,
\label{eq:Bdetermine}
\end{equation}
and we obtain a second $G\times G$ constant matrix from these column vectors
by setting $\mat{B}:=(\mat{b}^{(1)},\mat{b}^{(2)},\dots,\mat{b}^{(G)})$.  \label{symbol:Bmatrix}  
A vector $\mat{k}$ may now be defined by the formula \label{symbol:Riemannk}
\begin{equation}
\mat{k}:=\left\{
\begin{array}{ll}
\displaystyle \pi i\sum_{\text{$j$ odd}}\mat{e}^{(j)} + 
\frac{1}{2}\sum_{j=1}^G\mat{b}^{(j)}\,, & \text{$G$ odd}\,,\\\\
\displaystyle \pi i\sum_{\text{$j$ even}}\mat{e}^{(j)} +
\frac{1}{2}\sum_{j=1}^G\mat{b}^{(j)}\,,&\text{$G$ even}\,.
\end{array}
\right.
\label{eq:RiemannK}
\end{equation}
The matrix $\mat{B}$ is real, symmetric, and negative definite, so we may use it to define a Riemann theta function for $\mat{w}\in\cx^G$ by the Fourier series \label{symbol:thetafunction}
\begin{equation}
\Theta(\mat{w}):=\sum_{\mat{n}\in\mathbb{Z}^G}t_{\mat n}e^{\mat{n}^T\mat{w}}\,,\hspace{0.2 in}\text{with Fourier coefficients
$\displaystyle t_\mat{n}:=\exp\left(\frac{1}{2}\mat{n}^T\mat{B}\mat{n}
\right)\,.$}
\label{eq:thetafunction}
\end{equation}
Next, for $z\in\mathbb{C}\setminus\mathbb{R}$, set \label{symbol:Abelmap}
\begin{equation}
\mat{w}(z):=\int_{\alpha_0}^z\mat{A}\mat{m}(s)\,ds
\label{eq:Abelmap}
\end{equation}
where the path of integration lies in the half-plane $\Im(s)=\Im(z)$ but is otherwise arbitrary.  As special cases we set \label{symbol:wpminfty}
\begin{equation}
\mat{w}_+(\infty):=\mathop{\lim_{z\rightarrow\infty}}_{\Im(z)>0}\mat{w}(z)
\hspace{0.2 in}\text{and}\hspace{0.2 in}
\mat{w}_-(\infty):=\mathop{\lim_{z\rightarrow\infty}}_{\Im(z)<0}\mat{w}(z)\,.
\label{eq:wpminfty}
\end{equation}

Let $\lambda(z)$ be defined in the same domain as $y(z)$ by 
\begin{equation}
\lambda(z)^4 = \prod_{j=0}^G\frac{z-\alpha_j}{z-\beta_j}\,,\hspace{0.2 in}
\text{and $\lambda(z)\rightarrow 1$ as $z\rightarrow\infty$ with $\Im(z)>0$\,.}
\label{eq:lambda}
\end{equation}
In terms of $\lambda(z)$ we define two functions in the same domain by
setting \label{symbol:u(z)} \label{symbol:v(z)}
\begin{equation}
u(z):=\frac{1}{2}\left[\lambda(z)+\frac{1}{\lambda(z)}\right]
\hspace{0.2 in}\text{and}\hspace{0.2 in}
v(z):=\frac{1}{2i}\left[\lambda(z)-\frac{1}{\lambda(z)}\right]\,.
\label{eq:uv}
\end{equation}

The polynomial equation
\begin{equation}
\prod_{j=0}^G (x-\alpha_j) - \prod_{j=0}^G (x-\beta_j) = 0
\label{eq:polynomialequation}
\end{equation}
of degree $G$ has exactly one root $x=x_j$ in each gap $\Gamma_j$ for $j=1,\dots,G$.
Denoting the boundary values of $\mat{w}(z)$ taken on the real axis from the upper and lower half-planes by $\mat{w}_+(z)$ and $\mat{w}_-(z)$ respectively, we define two vectors by setting \label{symbol:qu}
\label{symbol:qv}
\begin{equation}
\mat{q}_u:=\sum_{j=1}^G\mat{w}_-(x_j)+\mat{k}\hspace{0.2 in}
\text{and}\hspace{0.2 in}
\mat{q}_v:=\sum_{j=1}^G\mat{w}_+(x_j)+\mat{k}\,.
\label{eq:quv}
\end{equation}

In terms of these ingredients we may now define two functions that turn out to extend analytically to the domain $\cx\setminus (-\infty,\beta_G]$:
\label{symbol:W(z)}
\begin{equation}
W(z):=\left\{\begin{array}{ll}
u(z)e^{h(z)}\frac{\displaystyle \Theta(\mat{w}_+(\infty)-\mat{q}_u)}
{\displaystyle \Theta(\mat{w}_+(\infty)-\mat{q}_u-i\mat{r}+i\kappa\mat{\Omega})}\,
\frac{\displaystyle \Theta(\mat{w}(z)-\mat{q}_u-i\mat{r}+i\kappa\mat{\Omega})}
{\displaystyle \Theta(\mat{w}(z)-\mat{q}_u)}\,, &\Im(z)>0\,,\\\\
-v(z)e^{h(z)}\frac{\displaystyle \Theta(\mat{w}_-(\infty)-\mat{q}_v)}
{\displaystyle\Theta(\mat{w}_-(\infty)-\mat{q}_v+i\mat{r}-i\kappa\mat{\Omega})}\,
\frac{\displaystyle \Theta(\mat{w}(z)-\mat{q}_v+i\mat{r}-i\kappa\mat{\Omega})}
{\displaystyle \Theta(\mat{w}(z)-\mat{q}_v)}\,,&\Im(z)<0\,,
\end{array}\right.
\label{eq:Wdefine}
\end{equation}
and
\label{symbol:Z(z)}
\begin{equation}
Z(z):=\left\{\begin{array}{ll}
iv(z)e^{-h(z)}\frac{\displaystyle \Theta(\mat{w}_-(\infty)-\mat{q}_v)}
{\displaystyle \Theta(\mat{w}_-(\infty)-\mat{q}_v+i\mat{r}-i\kappa\mat{\Omega})}\,
\frac{\displaystyle \Theta(\mat{w}(z)-\mat{q}_v+i\mat{r}-i\kappa\mat{\Omega})}
{\displaystyle \Theta(\mat{w}(z)-\mat{q}_v)}\,,&\Im(z)>0\,,\\\\
iu(z)e^{-h(z)}\frac{\displaystyle \Theta(\mat{w}_+(\infty)-\mat{q}_u)}
{\displaystyle \Theta(\mat{w}_+(\infty)-\mat{q}_u-i\mat{r}+i\kappa\mat{\Omega})}\,
\frac{\displaystyle \Theta(\mat{w}(z)-\mat{q}_u-i\mat{r}+i\kappa\mat{\Omega})}
{\displaystyle \Theta(\mat{w}(z)-\mat{q}_u)}\,, &\Im(z)<0\,.
\end{array}\right.
\label{eq:Zdefine}
\end{equation}
Finally, for any gap $\Gamma$ we may define two corresponding functions
$H^\pm_\Gamma(z)$ in terms of $W(z)$ and $Z(z)$:
\begin{equation}
H^\pm_\Gamma(z):=\frac{W(z)}{\sqrt{2}}
e^{(\gamma-\eta(z)-iN{\rm sgn}(\Im(z))\theta_\Gamma)/2} \pm
\frac{Z(z)}{\sqrt{2}}
e^{-(\gamma-\eta(z)-iN{\rm sgn}(\Im(z))\theta_\Gamma)/2}\,.
\label{eq:HGamma}
\end{equation}
\label{symbol:HGammapm}

\subsection{Results on asymptotics of discrete orthogonal polynomials.}
\label{sec:actualtheorems}
Subject to the basic assumptions described in
\S~\ref{sec:basicassumptions} and the simplifying assumptions
described in \S~\ref{sec:C3}, we have the following results, the
proofs of which will be given in \S~\ref{sec:asymptoticspi}.
\begin{theorem}[Outer asymptotics of $\pi_{N,k}(z)$]
Let $K$ be a closed set with $K\cap [a,b]=\emptyset$.  Then there exists
a constant $C_K>0$ such that
\begin{equation}
\pi_{N,k}(z)=e^{NL_c(z)}\left[W(z)+\varepsilon_N(z)\right]
\end{equation}
where the estimate
\begin{equation}
\sup_{z\in K}|\varepsilon_N(z)|\le\frac{C_K}{N}
\end{equation}
holds for sufficiently large $N$, and
$W(z)$ defined by \eqref{eq:Wdefine} is a function that is nonvanishing and uniformly bounded in $K$ independently of $N$.  Furthermore the product
$e^{NL_c(z)}W(z)$ is analytic for $z\in\cx\setminus[a,b]$.
\label{theorem:outside}
\end{theorem}

\begin{theorem}[Asymptotics of leading coefficients and recurrence coefficients]
If there is only a single band of unconstrained support of the equilibrium measure $\mu_{\rm min}^c$ in $[a,b]$, with endpoints $\alpha_0<\beta_0$, then
\begin{equation}
\gamma_{N,k}^2=\frac{4}{\beta_0-\alpha_0}e^{N\ell_c+\gamma}\left(1+\varepsilon_N^{(1)}\right)\,,
\end{equation}
\begin{equation}
\gamma_{N,k-1}^2 = \frac{\beta_0-\alpha_0}{4}e^{N\ell_c+\gamma}\left(1+\varepsilon_N^{(2)}\right)\,,
\end{equation}
\begin{equation}
b_{N,k-1}=\frac{\beta_0-\alpha_0}{4}\left(1+\varepsilon_N^{(3)}\right)\,,
\end{equation}
and
\begin{equation}
a_{N,k}=\frac{\beta_0+\alpha_0}{2}+\varepsilon_N^{(4)}\,,
\end{equation}
where there is a constant $C>0$ such that the estimates $|\varepsilon_N^{(m)}|\le C/N$, $m=1,2,3,4$, all hold for sufficiently large $N$.
More generally, if for some $G>0$ there are $G+1$ disjoint bands with endpoints $\alpha_0<\beta_0<\alpha_1<\beta_1<\dots<\alpha_G<\beta_G$, then
\begin{equation}
\gamma_{N,k}^2 = \frac{4e^{N\ell_c+\gamma}}{\displaystyle
\sum_{j=0}^G (\beta_j-\alpha_j)}
\frac{\Theta(\mat{w}_-(\infty)-\mat{q}_v+i\mat{r}-i\kappa\mat{\Omega})
\Theta(\mat{w}_+(\infty)-\mat{q}_v)}
{\Theta(\mat{w}_+(\infty)-\mat{q}_v+i\mat{r}-i\kappa\mat{\Omega})
\Theta(\mat{w}_-(\infty)-\mat{q}_v)}\left(1+\varepsilon_N^{(1)}\right)\,,
\end{equation}
\begin{equation}
\gamma_{N,k-1}^2=\frac{e^{N\ell_c+\gamma}}{4}\left[\sum_{j=0}^G
(\beta_j-\alpha_j)\right]\frac{\Theta(\mat{w}_+(\infty)-\mat{q}_v-i\mat{r}+i\kappa\mat{\Omega})\Theta(\mat{w}_-(\infty)-\mat{q}_v)}
{\Theta(\mat{w}_-(\infty)-\mat{q}_v-i\mat{r}+i\kappa\mat{\Omega})
\Theta(\mat{w}_+(\infty)-\mat{q}_v)}\left(1+\varepsilon_N^{(2)}\right)\,,
\end{equation}
\begin{equation}
\begin{array}{rcl}
b_{N,k-1}&=&\displaystyle\frac{1}{4}\left[\sum_{j=0}^G(\beta_j-\alpha_j)\right]
\frac{\Theta(\mat{w}_-(\infty)-\mat{q}_v)}{\Theta(\mat{w}_+(\infty)-\mat{q}_v)}\\\\
&&\displaystyle\hspace{0.2 in}\times\,\,\,\sqrt{\frac{\Theta(\mat{w}_+(\infty)-\mat{q}_v-i\mat{r}+i\kappa\mat{\Omega})\Theta(\mat{w}_+(\infty)-\mat{q}_v+i\mat{r}-i\kappa\mat{\Omega})}
{\Theta(\mat{w}_-(\infty)-\mat{q}_v-i\mat{r}+i\kappa\mat{\Omega})
\Theta(\mat{w}_-(\infty)-\mat{q}_v+i\mat{r}-i\kappa\mat{\Omega})}}
\left(1+\varepsilon^{(3)}_N\right)\,,
\end{array}
\end{equation}
and
\begin{equation}
\begin{array}{rcl}\displaystyle
a_{N,k}&=&\displaystyle
\frac{i\mat{a}^{(G)}\cdot\nabla\Theta(\mat{w}_+(\infty)+\mat{q}_v-i\mat{r}+i\kappa\mat{\Omega})}
{\Theta(\mat{w}_+(\infty)+\mat{q}_v-i\mat{r}+i\kappa\mat{\Omega})}-
\frac{i\mat{a}^{(G)}\cdot\nabla\Theta(\mat{w}_+(\infty)+\mat{q}_v)}
{\Theta(\mat{w}_+(\infty)+\mat{q}_v)}\\\\
&&\displaystyle\,\,\,+\,\,\,
\frac{i\mat{a}^{(G)}\cdot\nabla\Theta(\mat{w}_+(\infty)-\mat{q}_v+i\mat{r}-i\kappa\mat{\Omega})}
{\Theta(\mat{w}_+(\infty)-\mat{q}_v+i\mat{r}-i\kappa\mat{\Omega})}-
\frac{i\mat{a}^{(G)}\cdot\nabla\Theta(\mat{w}_+(\infty)-\mat{q}_v)}
{\Theta(\mat{w}_+(\infty)-\mat{q}_v)}\\\\
&&\displaystyle
\,\,\,+\,\,\,\frac{1}{2}\frac{\displaystyle\sum_{j=0}^G(\beta_j^2-\alpha_j^2)}{\displaystyle\sum_{j=0}^G(\beta_j-\alpha_j)}
+\varepsilon_N^{(4)}
\end{array}
\end{equation}
where $\mat{a}^{(G)}\cdot\nabla\Theta$ denotes the directional derivative of $\Theta(\cdot)$ in the direction of $\mat{a}^{(G)}$ in $\cx^G$, and where there is a constant $C>0$ such that the estimates $|\epsilon_N^{(m)}|\le C/N$, $m=1,2,3,4$, all hold for sufficiently large $N$.
\label{theorem:cs}
\end{theorem}

For an interval $J\subset[a,b]$ and any $\delta>0$ define the compact set
\begin{equation}
K_J^\delta:=\bigcup_{w\in J}\{\text{$z\in\mathbb{C}$ such that $|z-w|\le\delta $}\}\,.
\label{eq:fleshedJ}
\end{equation}
\label{symbol:KJdelta}

\begin{theorem}[Asymptotics of $\pi_{N,k}(z)$ in voids]
Let $J\subset[a,b]$ be a closed interval, and let $\Gamma$ be a void.
If $\Gamma=(a,\alpha_0)$, then assume that $J\subset [a,\alpha_0)$.
If $\Gamma=(\beta_G,b)$, then assume that $J\subset (\beta_G,b]$.
Finally, if $\Gamma=\Gamma_j=(\beta_{j-1},\alpha_j)$ for some $j=1,\dots,G$, then assume
that $J\subset\Gamma$.   There is a positive $\delta$
and a constant $C_J^\delta>0$ such that for $z\in K_J^\delta$ defined by
\eqref{eq:fleshedJ} we have
\begin{equation}
\pi_{N,k}(z)=e^{N\overline{L}^\Gamma_c(z)}\left[A^\nab_\Gamma(z)+\varepsilon_N(z)\right]
\end{equation}
where the estimate
\begin{equation}
\sup_{z\in K_J^\delta}|\varepsilon_N(z)|\le\frac{C_J^\delta}{N}
\label{eq:voidtheoremestimate}
\end{equation}
holds for sufficiently large $N$, and
\begin{equation}
A^\nab_\Gamma(z):=e^{N(L_c(z)-\overline{L}_c^\Gamma(z))}W(z)
\end{equation}
with $W(z)$ given by \eqref{eq:Wdefine} is a function that is real-analytic and uniformly bounded in $K_J^\delta$ independently of $N$.  If $\Gamma$ is adjacent to either endpoint, $z=a$ or $z=b$, then $A^\nab_\Gamma(z)$ does not vanish in $K_J^\delta$.  Otherwise,
$A^\nab_\Gamma(z)$ has at most one (real) zero in $K_J^\delta$.
\label{theorem:lower}
\end{theorem}

The possible lone zero of $A^\nab_\Gamma(z)$ in the void
$\Gamma$ is analogous to a {\em spurious zero} of
approximation theory.  The motion of a spurious zero through an
interior gap $\Gamma$ as parameters (like the degree $k$) are varied corresponds to the spontaneous emission of a zero from one band and
its subsequent capture by an adjacent band separated by a void.  At
most one zero can be in transit in $\Gamma$ for each choice of
parameters.

For $z$ in the domain of analyticity of $\rho^0(z)$, we set
\label{symbol:theta0}
\begin{equation}
\theta^0(z):=2\pi\int_z^b\rho^0(s)\,ds\,.
\label{eq:theta0def}
\end{equation}

\begin{theorem}[Asymptotics of $\pi_{N,k}(z)$ in saturated regions]
Let $J\subset\Gamma\subset[a,b]$ be a closed interval, and let $\Gamma$
be a saturated region. 
There is a positive $\delta$ and there are constants
$C_J^\delta>0$, $D_J^\delta>0$, and $E_J^\delta>0$ such that
for $z\in K_J^\delta$ defined by \eqref{eq:fleshedJ} we have
\begin{equation}
\pi_{N,k}(z)=e^{N\overline{L}^\Gamma_c(z)}\left[\left(A^\del_\Gamma(z)+\varepsilon_N(z)\right)\cos\left(\frac{N\theta^0(z)}{2}\right)+\delta_N(z)\right]
\label{eq:piNkupperasymp}
\end{equation}
where the estimates
\begin{equation}
\sup_{z\in K_J^\delta}|\varepsilon_N(z)|\le\frac{C_J^\delta}{N}
\hspace{0.2 in}\mbox{and}\hspace{0.2 in}
\sup_{z\in K_J^\delta}|\delta_N(z)|\le D_J^\delta e^{-NE_J^\delta}
\label{eq:sattheoremestimates}
\end{equation}
hold for sufficiently large $N$, and
\begin{equation}
A^\del_\Gamma(z):=2e^{N(L_c(z)-\overline{L}_c^\Gamma(z))}e^{-iN{\rm
sgn}(\Im(z))\theta^0(z)/2}W(z)
%\dot{X}_{11}(z)e^{\kappa g(z)}
\end{equation}
with $W(z)$ given by \eqref{eq:Wdefine} is a function that is
real-analytic and uniformly bounded in $K_J^\delta$ independently of
$N$.  If $\Gamma$ is adjacent to either endpoint, $z=a$ or $z=b$, then
$A^\del_\Gamma(z)$ does not vanish in $K_J^\delta$.  Otherwise,
$A^\del_\Gamma(z)$ has at most one (real) zero in $K_J^\delta$.
\label{theorem:upper}
\end{theorem}

When a saturated region meets an endpoint of $[a,b]$, we say that
there is a {\em hard edge} at that endpoint.  This terminology is
borrowed from random matrix theory, where it refers to an ensemble of
matrices all of which share a certain common bound on their spectra.
For example, random Wishart matrices of the form
$\mat{W}=\mat{X}^T\mat{X}$ for some real matrix $\mat{X}$ necessarily have
nonnegative spectra, and for certain types of matrices $\mat{X}$ the
asymptotic density of eigenvalues $z$ of $\mat{W}$ can have a
jump discontinuity at $z=0$, being identically zero for $z<0$ and
strictly positive for $z>0$ however small.  Thus $z=0$ is a hard edge
for the spectrum in Wishart random matrix ensembles.  We will see in
\S~\ref{sec:DOPensembles} that the density of the scaled equilibrium measure 
$\mu_{\rm min}^c/c$ plays the same role in certain discrete random
processes as the asymptotic density of eigenvalues plays in random
matrix theory.  Since the upper constraint always corresponds to a
strictly positive density, and since the support of the equilibrium
measure is a subset of $[a,b]$, an active upper constraint at either
$z=a$ or $z=b$ implies a jump discontinuity in the density of the
equilibrium measure at the corresponding endpoint, which explains our
terminology.

\begin{theorem}[Asymptotics of $\pi_{N,k}(z)$ near hard edges]
Suppose either that $\Gamma=(a,\alpha_0)$ is a saturated region and
$J=[a,t]$ for some $t\in \Gamma$, or that $\Gamma=(\beta_G,b)$ is a
saturated region and $J=[t,b]$ for some $t\in\Gamma$.  Set
\begin{equation}
\zeta:=\left\{\begin{array}{ll}
\displaystyle N\int_a^z\rho^0(s)\,ds\,, &\hspace{0.2 in}\mbox{if $a\in J$}\\\\
\displaystyle N\int_z^b\rho^0(s)\,ds\,, &\hspace{0.2 in}\mbox{if $b\in J$.}
\end{array}\right.
\end{equation}
There is a positive $\delta$ and there are constants
$C_J^\delta>0$, $D_J^\delta>0$, and $E_J^\delta>0$ such that for
$z\in K_J^\delta$ defined by \eqref{eq:fleshedJ} we have
\begin{equation}
\pi_{N,k}(z)=e^{N\overline{L}_c^\Gamma(z)}\zeta^{-\zeta}\left[
\left(\tilde{A}_\Gamma^\del(z) +\tilde{\varepsilon}_N(z)\right)
\frac{\Gamma(1/2+\zeta)\cos(\pi\zeta)}{\sqrt{2\pi}e^{-\zeta}} + \zeta^\zeta\delta_N(z)\right]\,,
\label{eq:piNkhardedge}
\end{equation}
where $\tilde{\varepsilon}_N(z)$ and $\zeta^\zeta\delta_N(z)$ extend from $\zeta>0$ as functions analytic in $K_J^\delta$ such that the estimates
\begin{equation}
\sup_{z\in K_J^\delta}|\tilde{\varepsilon}_N(z)|\le\frac{C_J^\delta}{N}
\hspace{0.2 in}\mbox{and}\hspace{0.2 in}
\sup_{z\in K_J^\delta}|\delta_N(z)|\le D_J^\delta e^{-NE_J^\delta}
\end{equation}
hold for sufficiently large $N$, and where
\begin{equation}
\tilde{A}_\Gamma^\del(z):=2e^{N(L_c(z)-\overline{L}_c^\Gamma(z))}
e^{-iN\pi {\rm sgn}(\Im(\zeta))\zeta}W(z)
%\dot{X}_{11}(z)e^{\kappa g(z)}
\end{equation}
with $W(z)$ given by \eqref{eq:Wdefine} is a function that is
real-analytic, nonvanishing, and uniformly bounded in $K_J^\delta$
indpendently of $N$.  Finally, note that
$e^{N\overline{L}^\Gamma_c(z)}\zeta^{-\zeta}$ extends from $\zeta>0$
as an analytic function in $K_J^\delta$, and that $\delta_N(z)$
represents exactly the same function as in
Theorem~\ref{theorem:upper}.
\label{theorem:hardedges}
\end{theorem}

\begin{remark}
The fact that the asymptotic formulae presented in
Theorem~\ref{theorem:hardedges} are in terms of the Euler gamma
function is directly related to the discrete nature of the weights.
In a sense, the poles of the functions $\Gamma(1/2+\zeta)$ are
``shadows'' of the poles of the matrix $\mat{P}(z;N,k)$ solving
Interpolation Problem~\ref{rhp:DOP}.
\end{remark}

\begin{theorem}[Exponential confinement of zeros in saturated regions]
Let $J\subset[a,b]$ be a closed interval, and let $\Gamma$ be a
saturated region.  
If $\Gamma=(a,\alpha_0)$, then assume that $J\subset[a,\alpha_0)$.
If $\Gamma=(\beta_G,b)$, then assume that $J\subset(\beta_G,b]$.
Finally, if $\Gamma=\Gamma_j=(\beta_{j-1},\alpha_j)$ for some $j=1,\dots,G$,
then assume that $J\subset\Gamma$.  There are positive constants $D_J$, 
$E_J$, and $N_0$ such that 
for every node $x_{N,n}\in X_N\cap J$ there exists a zero $z_0$ of the
monic discrete orthogonal polynomial $\pi_{N,k}(z)$ with
\begin{equation}
|z_0-x_{N,n}|\le D_Je^{-NE_J}\,,\hspace{0.2 in}\text{whenever $N>N_0$}\,.
\label{eq:twinning}
\end{equation}
Moreover:
\begin{enumerate}
\item
If $\Gamma=(a,\alpha_0)$, so that $\min J\ge a$, 
then for each node $x_{N,n}\in J\cap X_N$ there is a  
zero $z_0$ of
$\pi_{N,k}(z)$ such that
\begin{equation}
x_{N,n}<z_0<x_{N,n}+D_Je^{-NE_J}\,,\hspace{0.2 in}\text{whenever $N>N_0$}\,,
\label{eq:twinningright}
\end{equation}
and for each zero $z_0\in J$ of $\pi_{N,k}(z)$, there is a node $x_{N,n}\in X_N$ such that \eqref{eq:twinningright} holds.
%lies to the right of the node to which it is
%attracted, and there are no other zeros of $\pi_{N,k}(z)$ in $J$.
\item
If $\Gamma=(\beta_G,b)$, so that $\max J\le b$, 
then for each node $x_{N,n}\in J\cap X_N$ there is a zero $z_0$ of
$\pi_{N,k}(z)$ such that
\begin{equation}
x_{N,n}-D_Je^{-NE_J}<z_0<x_{N,n}\,,\hspace{0.2 in}
\text{whenever $N>N_0$}\,,
\label{eq:twinningleft}
\end{equation}
and for each zero $z_0\in J$ of $\pi_{N,k}(z)$, there is a node $x_{N,n}\in X_N$ such that \eqref{eq:twinningleft} holds.
%then each zero of
%$\pi_{N,k}(z)$ in $J$ lies to the left of the node to which it is
%attracted, and there are no other zeros of $\pi_{N,k}(z)$ in $J$.
\item
If $\Gamma=\Gamma_j=(\beta_{j-1},\alpha_j)$ for some $j=1,\dots,G$,
then exactly one of the following two mutually exclusive possibilities
holds:
\begin{enumerate}
\item
There is a node $x_{N,m}\in\Gamma\cap X_N$ such that
$\pi_{N,k}(x_{N,m})=0$.  
For
each node $x_{N,n}\in J\cap X_N$ with $x_{N,n}>x_{N,m}$
there is a zero $z_0$ of $\pi_{N,k}(z)$ such that
\eqref{eq:twinningright} holds, and for each zero $z_0\in J$ of
$\pi_{N,k}(z)$ with $z_0>x_{N,m}$ there is a
node $x_{N,n}\in X_N$ such that \eqref{eq:twinningright} holds.
For each node $x_{N,n}\in J\cap X_N$ with 
$x_{N,n} < x_{N,m}$ there is a zero $z_0$ of $\pi_{N,k}(z)$
such that \eqref{eq:twinningleft} holds, and for each zero $z_0\in J$ of
$\pi_{N,k}(z)$ with $z_0 < x_{N,m}$ there is a
node $x_{N,n}\in X_N$ such that \eqref{eq:twinningleft} holds.  
%
%There is a node $x_{N,n}$ (not necessarily contained in $J$) such that
%$\pi_{N,k}(x_{N,n})=0$, while in the interval $[\min J,x_{N,n})$ there
%is exactly one zero of $\pi_{N,k}(z)$ lying exponentially near and to
%the left of each node in this interval and no other zeros, and in the
%interval $(x_{N,n},\max J]$ there is exactly one zero of
%$\pi_{N,k}(z)$ lying exponentially near and to the right of each node
%in this interval and no other zeros.
\item
There is a consecutive pair of nodes $x_{N,m}\in\Gamma\cap X_N$ and
$x_{N,m+1}\in\Gamma\cap X_N$ such that
\begin{itemize}
\item
For each node $x_{N,n}\in J\cap X_N$ with $x_{N,n}\ge x_{N,m+1}$ there is a 
zero $z_0$ of $\pi_{N,k}(z)$ such that \eqref{eq:twinningright} holds,
and for each zero $z_0\in J$ of $\pi_{N,k}(z)$ with $z_0\ge x_{N,m+1}$ there
is a node $x_{N,n}\in X_N$ such that \eqref{eq:twinningright} holds.
%In the closed interval $[\min J,x_{N,n}]$ there is a zero of
%$\pi_{N,k}(z)$ exponentially near and to the left of each node in
%$X_N\cap[\min J,x_{N,n})$ and no other zeros.
\item
For each node $x_{N,n}\in J\cap X_N$ with $x_{N,n}\le x_{N,m}$ there is a
zero $z_0$ of $\pi_{N,k}(z)$ such that \eqref{eq:twinningleft} holds,
and for each zero $z_0\in J$ of $\pi_{N,k}(z)$ with $z_0\le x_{N,m}$ there
is a node $x_{N,n}\in X_N$ such that \eqref{eq:twinningleft} holds.
%In the closed interval $[x_{N,n+1},\max J]$ there is a zero of
%$\pi_{N,k}(z)$ exponentially near and to the right of each node in
%$X_N\cap[x_{N,n+1},\max J]$ and no other zeros.
\item
There is at most one zero $z_0$ of $\pi_{N,k}(z)$ in the closed interval
$[x_{N,m},x_{N,m+1}]$, and if it exists then $z_0\in (x_{N,m},x_{N,m+1})$.
%There may or may not be a single zero
%of $\pi_{N,k}(z)$ in the closed interval $[x_{N,n},x_{N,n+1}]$, and if
%it exists it may lie anywhere in the interval including the endpoints.
%That is, it is not necessarily exponentially attracted to either node.
\end{itemize}
\end{enumerate}
\end{enumerate}
\label{theorem:exponential}
\end{theorem}
Note that in case 3(b), if there is a zero $z_0$ of $\pi_{N,k}(z)$ with
$x_{N,m}<z_0<x_{N,m+1}$ there need not be any node $x_{N,n}\in X_N$ such
that \eqref{eq:twinning} holds.  This particular zero, and only this one,
is not necessarily exponentially close to any node.

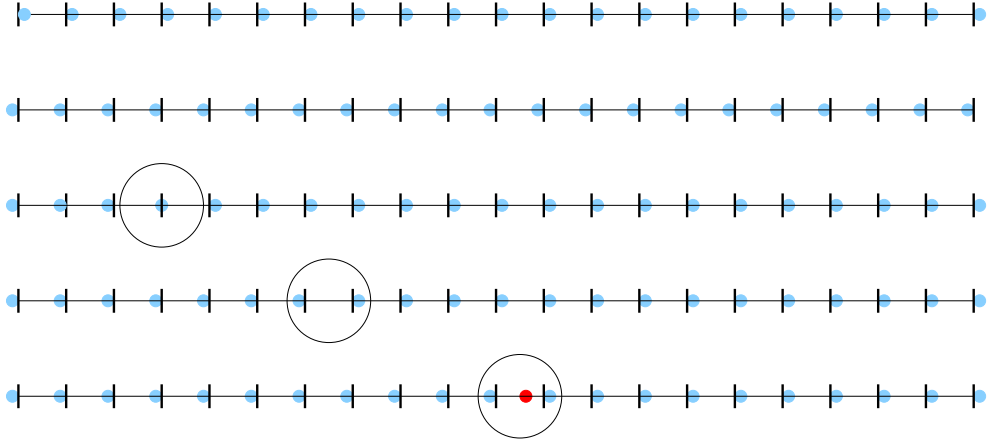
\begin{figure}[h]
\begin{center}
\input{defects.pstex_t}
\end{center}
\caption{\em First line:  the pattern of zeros of $\pi_{N,k}(z)$ 
(blue dots) and nodes (vertical segments) in a saturated region
adjacent to the left endpoint $z=a$.  Second line: same as above, but
for a saturated region adjacent to the right endpoint $z=b$.  Third
line: a pattern of zeros and nodes in a saturated region between two
bands; there is a single defect (circled) corresponding to a zero of
$\pi_{N,k}(z)$ occurring exactly at a node.  Fourth line: a pattern as
above in which the defect (circled) is an interval
$[x_{N,m},x_{N,m+1}]$ that does not carry any spurious zero.  Fifth
line: same as above, but a case in which the defect (circled) carries
a spurious zero (distinguished by red shading).}
\label{fig:exponentialzeros}
\end{figure}

We refer to the node $x_{N,m}$ in 3(a) and also to the interval
$[x_{N,m},x_{N,m+1}]$ in 3(b), both of which serve to separate the two
directions of perturbation of the zeros of $\pi_{N,k}(z)$ from the
nodes, as {\em defects}, and to the zero possibly carried by the
defect in 3(b) as a {\em spurious zero}.  The remaining zeros
correspond in a one-to-one fashion with the nodes; we refer to them as
{\em Hurwitz zeros} by analogy with the approximation theory
literature.  See Figure~\ref{fig:exponentialzeros}.

\begin{remark}
It should perhaps be stressed that there is nothing in principle that
prevents a zero of $\pi_{N,k}(z)$ from coinciding {\em exactly} with
one of the nodes $x_{N,n}\in X_N$.  Indeed, this is the case in 3(a)
above.  However, Theorem~\ref{theorem:exponential} shows that in
saturated regions adjacent to endpoints $z=a$ or $z=b$ all zeros
become asymptotically distinct from (yet paradoxically converge
rapidly to) nodes as $N\rightarrow\infty$.  In saturated regions lying
between two bands it is asymptotically only possible for a single zero
to coincide exactly with a node.
\end{remark}

The precise location of a defect within a saturated region depends on
all the parameters of the problem, and in some circumstances it may
make sense for a parameter, say appearing in the function $V(x)$
defined in (\ref{eq:VNexpand}), to be continuously varied.  This is
interesting because it can imply corresponding dynamics of the defects
and any spurious zeros they may carry.  If continuous deformation of a
parameter leads to that of the phase vector
$\mat{r}-\kappa\mat{\Omega}$ then the defect will move continuously
through the saturated region $\Gamma$ as well.

\begin{figure}[h]
\begin{center}
\input{defectmotion.pstex_t}
\end{center}
\caption{\em Left: the motion of a defect (circled) 
that does not carry any spurious zero.  Right: the motion of a defect
(circled) carrying a spurious zero, which exchanges its identity with
a Hurwitz zero in each step.}
\label{fig:defectmotion}
\end{figure}
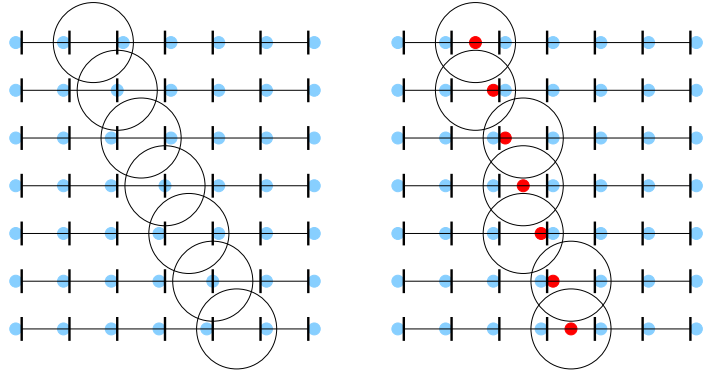

How does a defect move?  If there is no spurious zero, then a defect
$[x_{N,m},x_{N,m+1}]$ can move to the right to become a defect
$[x_{N,m+1},x_{N,m+2}]$ as the zero of $\pi_{N,k}(z)$ just to the
right of the node $x_{N,m+1}$ moves continuously to the left through
the node.  Then the same process then occurs near the node $x_{N,m+2}$
and so on.  Thus a defect without a spurious zero moves to the right
by a process in which Hurwitz zeros move to the left an exponentially
small amount, passing through the corresponding nodes, one after the
other.  During the continuous motion of a defect without a spurious
zero the situation described in 3(a) above occurs only at isolated
values of the deformation parameter on which the phase vector
$\mat{r}-\kappa\mat{\Omega}$ continuously depends.  See the left
diagram in Figure~\ref{fig:defectmotion}.

If the defect $[x_{N,m},x_{N,m+1}]$ contains a spurious zero, then the
motion of the defect to the right occurs by a change-of-identity
process in which the spurious zero moves to the right through the
defect toward $x_{N,m+1}$, and when it is exponentially close to
$x_{N,m+1}$ it becomes a Hurwitz zero and the previously Hurwitz zero
just to the right of $x_{N,m+1}$ becomes a spurious zero belonging to
the new defect $[x_{N,m+1},x_{N,m+2}]$.  Thus a defect carrying a
spurious zero moves to the right by a process in which zeros move to
the right by an amount proportional to $1/N$ one after the other.  
During the continuous motion of a defect containing a spurious zero
the situation described in 3(a) above never occurs at all.  See
the right diagram in Figure~\ref{fig:defectmotion}.

In fact, a defect carrying a spurious zero that reaches an endpoint of
$\Gamma$ under deformation will generally be reflected back into
$\Gamma$ as a defect without a spurious zero.  Thus as parameters are
deformed, a defect may oscillate back and forth within a saturated
region acting like a conveyor belt, carrying a spurious zero from one
band to the next, and returning empty to pick up the next zero.

\begin{theorem}[Asymptotics of $\pi_{N,k}(z)$ in bands]
Let $J$ be a closed interval and let $I=(\alpha_j,\beta_j)$ be a band.
Assume that $J\subset I$.  There is a positive $\delta$ and a constant
$C_J^\epsilon>0$ such that for $z\in K_J^\delta$ defined by
\eqref{eq:fleshedJ} we have
\begin{equation}
\pi_{N,k}(z)=e^{N\overline{L}^I_c(z)}\left[A_I(z)\cos\left(\Phi_I(z)+ N\pi c
\mu_{\rm min}^c([x,b])-N\pi c\int_x^z\psi_I(s)\,ds\right) + \varepsilon_I(z)\right]
\end{equation}
where $x$ is any point (or endpoint) of $I$, and $\varepsilon_I(z)$
satisfies the estimate
\begin{equation}
\sup_{z\in K_J^\delta}|\varepsilon_I(z)|\le\frac{C_J^\epsilon}{N}\,.
\label{eq:bandestimate}
\end{equation}
Here, $A_I(z)$ and $\Phi_I(z)$ are real-analytic functions of $z$ that
are uniformly bounded independently of $N$ in $K_J^\delta$.
Moreover $A_I(z)$ is strictly positive for real $z$.  For real $z$,
the identity
\begin{equation}
W_+(z) = \frac{1}{2}A_I(z)e^{i\Phi_I(z)}
\end{equation}
holds, where $W(z)$ is defined by \eqref{eq:Wdefine} and $W_+(z)$
indicates the boundary value taken from the upper half-plane.  This
relation serves as a definition of the analytic functions $A_I(z)$ and
$\Phi_I(z)$.
\label{theorem:band}
\end{theorem}

\begin{theorem}[Asymptotic description of zeros of $\pi_{N,k}(z)$ in bands]
Let $J$ be a closed interval and let $I=(\alpha_j,\beta_j)$ be a band.
Assume that $J\subset I$.  Then, the zeros of $\pi_{N,k}(z)$ in $J$
correspond in a one-to-one fashion with those of the model function
\begin{equation}
C_I(z):=\cos\left(\Phi_I(z)+ N\pi c
\mu_{\rm min}^c([x,b])-N\pi c\int_x^z\psi_I(s)\,ds\right)\,,
\end{equation}
where $\Phi_I(z)$ is as in the statement of Theorem~\ref{theorem:band}
above, and $x$ is any point (or endpoint) of $I$.  Moreover, there
exists a constant $D_J^\epsilon>0$ such that if $N$ is sufficiently
large, each pair of corresponding zeros $z_0$ of $\pi_{N,k}(z)$ and
$\tilde{z}_0$ of $C_I(z)$ in $J$ satisfies the estimate
\begin{equation}
|z_0-\tilde{z}_0|\le\frac{D_J^\epsilon}{N^2}\,.
\end{equation}
\label{theorem:bandzeros}
\end{theorem}

Before we state the next results, we point out that the function
$\overline{L}^I_c(z)$ defined for each band $I=(\alpha,\beta)$ (see
(\ref{eq:LcbarI})) may be considered to be analytic in a complex
neighborhood of the closed interval $[\alpha,\beta]$. In particular,
$\overline{L}_c^I(z)$ is analytic in a neighborhood of each endpoint
of the band.  The analytic continuation to a neighborhood $U_\alpha$
of $z=\alpha$ is accomplished
by the identity
\begin{equation}
\overline{L}^I_c(z)=\overline{L}^\Gamma_c(z)+
\frac{1}{2N}\left(-\tau_\Gamma^{\nab,L}(z)\right)^{3/2}\,,
\hspace{0.2 in}
\text{for $z\in U_\alpha$ with $\Im(z)\neq 0$\,,}
\end{equation}
if the adjacent gap $\Gamma$ is a void, and by the identity
\begin{equation}
\overline{L}_c^I(z)=\overline{L}_c^\Gamma(z)-\frac{1}{2N}
\left(-\tau_\Gamma^{\del,L}(z)\right)^{3/2}\,,\hspace{0.2 in}
\text{for $z\in U_\alpha$ with $\Im(z)\neq 0$\,,}
\end{equation}
if the adjacent gap $\Gamma$ is a saturated region.
Similarly, the analytic continuation to a neighborhood $U_\beta$ of
$z=\beta$ is accomplished by
the identity
\begin{equation}
\overline{L}_c^I(z)=\overline{L}_c^\Gamma(z)+\frac{1}{2N}
\left(-\tau_\Gamma^{\nab,R}(z)\right)^{3/2}\,,
\hspace{0.2 in}\text{for $z\in U_\beta$ with $\Im(z)\neq 0$\,,}
\end{equation}
if the adjacent gap $\Gamma$ is a void, and by the identity
\begin{equation}
\overline{L}_c^I(z)=\overline{L}_c^\Gamma(z)-\frac{1}{2N}
\left(-\tau_\Gamma^{\del,R}(z)\right)^{3/2}\,,
\hspace{0.2 in}\text{for $z\in U_\beta$ with $\Im(z)\neq 0$\,,}
\end{equation}
if the adjacent gap $\Gamma$ is a saturated region.

\begin{theorem}[Asymptotics of $\pi_{N,k}(z)$ near band/void edges] 
Let $z=\alpha$ be the left endpoint of a band $I$, and suppose that a
void $\Gamma$ lies immediately to the left of $z=\alpha$.  There exist
constants $r>0$ and $C>0$ such that when $|z-\alpha|\le r$,
\begin{equation}
\begin{array}{rcl}
\displaystyle \pi_{N,k}(z)&=&\displaystyle
e^{N\overline{L}_c^I(z)}\Bigg[N^{1/6}
\left(A_\Gamma^{\nab,L}(z)+\varepsilon_A(z)\right)
Ai\left(-\left(\frac{3}{4}\right)^{2/3}\tau_{\Gamma}^{\nab,L}(z)
\right) \\\\
&&\displaystyle\,\,\,+\,\,\, 
N^{-1/6}\left(B_\Gamma^{\nab,L}(z)+\varepsilon_B(z)
\right)Ai'\left(-\left(\frac{3}{4}\right)^{2/3}
\tau_{\Gamma}^{\nab,L}(z)\right)\Bigg]\,,
\end{array}
\label{eq:AirynabLuniform}
\end{equation}
where the estimates
\begin{equation}
\sup_{|z-\alpha|\le r}|\varepsilon_A(z)|\le\frac{C}{N}\hspace{0.2 in}
\text{and}\hspace{0.2 in}
\sup_{|z-\alpha|\le r}|\varepsilon_B(z)|\le\frac{C}{N}
\end{equation}
both hold for all $N$ sufficiently large, and where the leading
coefficient functions defined by
\begin{equation}
\begin{array}{rcl}
\displaystyle A_\Gamma^{\nab,L}(z)&:=&\displaystyle
\left(\frac{3}{4}\right)^{1/6}\sqrt{2\pi}e^{(\eta(z)-\gamma)/2}H_\Gamma^-(z)
\cdot N^{-1/6}\left(-\tau_\Gamma^{\nab,L}(z)\right)^{1/4}\,,\\\\
\displaystyle B_\Gamma^{\nab,L}(z)&:=&\displaystyle
-\left(\frac{3}{4}\right)^{-1/6}\sqrt{2\pi}e^{(\eta(z)-\gamma)/2}
H_\Gamma^+(z)\cdot N^{1/6}\left(-\tau_\Gamma^{\nab,L}(z)\right)^{-1/4}\,,
\end{array}
\end{equation}
both are real-analytic functions for $|z-\alpha|\le r$ that remain uniformly
bounded in this disc as $N\rightarrow\infty$.  
Furthermore, we may also write
\begin{equation}
\pi_{N,k}(z)=e^{N\overline{L}_c^I(z)}\left[N^{1/6}A_\Gamma^{\nab,L}(\alpha)
Ai\left(-\left(\frac{3}{4}\right)^{2/3}\tau_{\Gamma}^{\nab,L}(z)\right)
+\delta(z)\right]
\label{eq:AirynabLinner}
\end{equation}
where the estimate
\begin{equation}
\sup_{|z-\alpha|\le rN^{-2/3}}|\delta(z)|\le\frac{C}{N^{1/6}}
\end{equation}
holds for all sufficiently large $N$.

Let $z=\beta$ be the right endpoint of a band $I$, and suppose that a
void $\Gamma$ lies immediately to the right of $z=\beta$.  There exist
constants $r>0$ and $C>0$ such that when $|z-\beta|\le r$,
\begin{equation}
\begin{array}{rcl}
\pi_{N,k}(z)&=&\displaystyle e^{N\overline{L}_c^I(z)}\Bigg[
N^{1/6}\left(A_\Gamma^{\nab,R}(z)+\varepsilon_A(z)\right)
Ai\left(-\left(\frac{3}{4}\right)^{2/3}\tau_\Gamma^{\nab,R}(z)\right)\\\\
&&\displaystyle\,\,\,+\,\,\,
N^{-1/6}\left(B_\Gamma^{\nab,R}(z)+\varepsilon_B(z)\right)
Ai'\left(-\left(\frac{3}{4}\right)^{2/3}\tau_\Gamma^{\nab,R}(z)\right)
\Bigg]\,,
\end{array}
\label{eq:AirynabRuniform}
\end{equation}
where the estimates
\begin{equation}
\sup_{|z-\beta|\le r} |\varepsilon_A(z)|\le\frac{C}{N}\hspace{0.2 in}\text{and}
\hspace{0.2 in}
\sup_{|z-\beta|\le r} |\varepsilon_B(z)|\le\frac{C}{N}\hspace{0.2 in}
\end{equation}
both hold for all $N$ sufficiently large, and where the leading coefficient
functions defined by
\begin{equation}
\begin{array}{rcl}
\displaystyle A_\Gamma^{\nab,R}(z)&:=&\displaystyle
\left(\frac{3}{4}\right)^{1/6}\sqrt{2\pi}e^{(\eta(z)-\gamma)/2}H_\Gamma^+(z)
\cdot N^{-1/6}\left(-\tau_\Gamma^{\nab,R}(z)\right)^{1/4}\,,\\\\
\displaystyle B_\Gamma^{\nab,R}(z)&:=&\displaystyle
-\left(\frac{3}{4}\right)^{-1/6}\sqrt{2\pi}e^{(\eta(z)-\gamma)/2}H_\Gamma^-(z)
\cdot N^{1/6}\left(-\tau_\Gamma^{\nab,R}(z)\right)^{-1/4}\,,
\end{array}
\end{equation}
both are real-analytic functions for $|z-\beta|\le r$ that remain uniformly
bounded in this disc as $N\rightarrow\infty$.  Furthermore, we may also write
\begin{equation}
\pi_{N,k}(z)=e^{N\overline{L}_c^I(z)}\left[N^{1/6}A_\Gamma^{\nab,R}(\beta)
Ai\left(-\left(\frac{3}{4}\right)^{2/3}\tau_\Gamma^{\nab,R}(z)\right)
+\delta(z)\right]
\label{eq:AirynabRinner}
\end{equation}
where the estimate
\begin{equation}
\sup_{|z-\beta|\le rN^{-2/3}}|\delta(z)|\le\frac{C}{N^{1/6}}
\end{equation}
holds for all sufficiently large $N$.
\label{theorem:Airynab}
\end{theorem}

Note that these asymptotic formulae are similar in nature to the
corresponding asymptotic formulae found in \cite{DeiftKMVZ99} for
polynomials orthogonal with respect to analytic weights on the whole
real line.  On the other hand, there is no analogue of a saturated
region for continuous weights. The asymptotic band edge behavior
between a band $I$ and a saturated region $\Gamma$ involves the Airy
function $Bi(\cdot)$ as well as $Ai(\cdot)$ and is the subject of the
next theorem.

\begin{theorem}[Asymptotics of $\pi_{N,k}(z)$ near band/saturated region edges]
Let $z=\alpha$ be the left endpoint of a band $I$, and suppose that a
saturated region $\Gamma$ lies immediately to the left of $z=\alpha$.
There exist constants $r>0$ and $C>0$ such that when $|z-\alpha|\le r$,
\begin{equation}
\pi_{N,k}(z)=e^{N\overline{L}_c^I(z)}\left[N^{1/6}\left(A_\Gamma^{\del,L}(z)+
\varepsilon_A(z)\right)F_A^L(z) +
N^{-1/6}\left(B_\Gamma^{\del,L}(z)+\varepsilon_B(z)\right)F_B^L(z)\right]
\label{eq:AirydelLuniform}
\end{equation}
with
\begin{equation}
\begin{array}{rcl}
\displaystyle F_A^L(z)&:=&\displaystyle\cos\left(\frac{N\theta^0(z)}{2}\right)
Bi\left(-\left(\frac{3}{4}\right)^{2/3}\tau_\Gamma^{\del,L}(z)\right)-
\sin\left(\frac{N\theta^0(z)}{2}\right)Ai\left(-\left(\frac{3}{4}\right)^{2/3}
\tau_\Gamma^{\del,L}(z)\right)\,,\\\\
\displaystyle F_B^L(z)&:=&\displaystyle\cos\left(\frac{N\theta^0(z)}{2}\right)
Bi'\left(-\left(\frac{3}{4}\right)^{2/3}\tau_\Gamma^{\del,L}(z)\right)-
\sin\left(\frac{N\theta^0(z)}{2}\right)Ai'\left(-\left(\frac{3}{4}\right)^{2/3}
\tau_\Gamma^{\del,L}(z)\right)\,,
\end{array}
\label{eq:FALBL}
\end{equation}
where the estimates
\begin{equation}
\sup_{|z-\alpha|\le r}|\varepsilon_A(z)|\le\frac{C}{N}\hspace{0.2 in}\text{and}
\hspace{0.2 in}
\sup_{|z-\alpha|\le r}|\varepsilon_B(z)|\le\frac{C}{N}
\end{equation}
both hold for all $N$ sufficiently large, and where the leading
coefficient functions defined by
\begin{equation}
\begin{array}{rcl}
\displaystyle A_\Gamma^{\del,L}(z)&:=& \displaystyle
\left(\frac{3}{4}\right)^{1/6}\sqrt{2\pi}e^{(\eta(z)-\gamma)/2}H_\Gamma^-(z)
\cdot N^{-1/6}\left(-\tau_{\Gamma}^{\del,L}(z)\right)^{1/4}\,,\\\\
\displaystyle B_\Gamma^{\del,L}(z)&:=& \displaystyle
\left(\frac{3}{4}\right)^{-1/6}\sqrt{2\pi}e^{(\eta(z)-\gamma)/2} H_\Gamma^+(z)
\cdot N^{1/6}\left(-\tau_\Gamma^{\del,L}(z)\right)^{-1/4}\,,
\end{array}
\end{equation}
both are real-analytic functions for $|z-\alpha|\le r$ that remain uniformly
bounded in this disc as $N\rightarrow\infty$.  Furthermore, we may also write
\begin{equation}
\pi_{N,k}(z)=e^{N\overline{L}_c^I(z)}\left[N^{1/6}A_\Gamma^{\del,L}(\alpha)
F_A^L(z)
+\delta(z)\right]
\label{eq:AirydelLinner}
\end{equation}
where the estimate
\begin{equation}
\sup_{|z-\alpha|\le rN^{-2/3}} |\delta(z)|\le\frac{C}{N^{1/6}}
\end{equation}
holds for all sufficiently large $N$.

Let $z=\beta$ be the right endpoint of a band $I$, and suppose that a
saturated region $\Gamma$ lies immediately to the right of $z=\beta$.
There exist constants $r>0$ and $C>0$ such that when $|z-\beta|\le r$,
\begin{equation}
\pi_{N,k}(z)=e^{N\overline{L}_c^I(z)}\left[N^{1/6}\left(A_\Gamma^{\del,R}(z)
+\varepsilon_A(z)\right)F_A^R(z)+N^{-1/6}\left(B_\Gamma^{\del,R}(z)+
\varepsilon_B(z)\right)F_B^R(z)\right]
\label{eq:AirydelRuniform}
\end{equation}
with
\begin{equation}
\begin{array}{rcl}
\displaystyle F_A^R(z)&:=&\displaystyle
\cos\left(\frac{N\theta^0(z)}{2}\right)
Bi\left(-\left(\frac{3}{4}\right)^{2/3}\tau_\Gamma^{\del,R}(z)\right)+
\sin\left(\frac{N\theta^0(z)}{2}\right)
Ai\left(-\left(\frac{3}{4}\right)^{2/3}\tau_\Gamma^{\del,R}(z)\right)\,,\\\\
\displaystyle F_B^R(z)&:=&\displaystyle
\cos\left(\frac{N\theta^0(z)}{2}\right)
Bi'\left(-\left(\frac{3}{4}\right)^{2/3}\tau_\Gamma^{\del,R}(z)\right)+
\sin\left(\frac{N\theta^0(z)}{2}\right)
Ai'\left(-\left(\frac{3}{4}\right)^{2/3}\tau_\Gamma^{\del,R}(z)\right)\,,
\end{array}
\label{eq:FARBR}
\end{equation}
where the estimates
\begin{equation}
\sup_{|z-\beta|\le r}|\varepsilon_A(z)|\le\frac{C}{N}\hspace{0.2 in}\text{and}
\hspace{0.2 in}
\sup_{|z-\beta|\le r}|\varepsilon_B(z)|\le\frac{C}{N}
\end{equation}
both hold for all $N$ sufficiently large, and where the leading coefficient
functions defined by
\begin{equation}
\begin{array}{rcl}
\displaystyle A_\Gamma^{\del,R}(z)&:=&\displaystyle
\left(\frac{3}{4}\right)^{1/6}\sqrt{2\pi}e^{(\eta(z)-\gamma)/2}
H_\Gamma^-(z)\cdot N^{-1/6}\left(-\tau_\Gamma^{\del,R}(z)\right)^{1/4}\,,\\\\
\displaystyle B_\Gamma^{\del,R}(z)&:=&\displaystyle
-\left(\frac{3}{4}\right)^{-1/6}\sqrt{2\pi}e^{(\eta(z)-\gamma)/2}
H_\Gamma^+(z)\cdot N^{1/6}\left(-\tau_\Gamma^{\del,R}(z)\right)^{-1/4}\,,
\end{array}
\end{equation}
both are real-analytic functions for $|z-\beta|\le r$ that remain uniformly
bounded in this disc as $N\rightarrow\infty$.  Furthermore, we may also
write
\begin{equation}
\pi_{N,k}(z)=e^{N\overline{L}_c^I(z)}\left[N^{1/6}A_\Gamma^{\del,R}(\beta)
F_A^R(z)+\delta(z)\right]
\label{eq:AirydelRinner}
\end{equation}
where the estimate
\begin{equation}
\sup_{|z-\beta|\le rN^{-2/3}}|\delta(z)|\le\frac{C}{N^{1/6}}
\end{equation}
holds for all sufficiently large $N$.
\label{theorem:Airydel}
\end{theorem}

With the proper choice of the closed set $K$ in
Theorem~\ref{theorem:outside}, we see that the whole complex $z$-plane
has been covered with overlapping closed sets, in each of which there
is an associated asymptotic formula for $\pi_{N,k}(z)$ with rigorous error
bounds.  

\subsection{Equilibrium measures for classical discrete orthogonal 
polynomials.}
\label{sec:examplepolys}
%In this section, we illustrate briefly how two families of classical discrete
%orthogonal polynomials fit into the general framework we have developed
%so far in this paper.
Since the asymptotic behavior of the discrete orthogonal polynomials
is determined by the equilibrium measure $\mu_{\rm min}^c$
corresponding to the functions $\rho^0(x)$, $V(x)$, the interval
$[a,b]$, and the constant $c$, it will be useful to demonstrate that
the results stated in \S~\ref{sec:actualtheorems} can be made
effective by a concrete calculation of the equilibrium measure.  We
consider below two classical cases.  The equilibrium measure for the
Krawtchouk polynomials was obtained by Dragnev and Saff in
\cite{DragnevS00}.  The equilibrium measure for the Hahn polynomials
has not appeared in the literature before (to our knowledge) and we
present it below as well.
\subsubsection{The Krawtchouk polynomials.}
\label{sec:Krawtchouk}
The Krawtchouk polynomials \cite{AbramowitzS65} are orthogonal on
a finite set of equally
spaced nodes in the interval $(0,1)$:
\begin{equation}
x_{N,n}:=\frac{2n+1}{2N}\hspace{0.2 in}\text{for}\hspace{0.2 in}
n=0,1,2,\dots,N-1\,.
\end{equation}
The analytic probability density on $(0,1)$ is then given simply by
$\rho^0(x)\equiv 1$.  The corresponding weights are given by
\begin{equation}
w^{\rm Kraw}_{N,n}(p,q):=
\frac{N^{N-1}\sqrt{pq}}{q^N\Gamma(N)}
\binom{N-1}{n}
p^nq^{N-1-n}
\end{equation}
\label{symbol:wKraw}
where $p$ and $q$ are positive parameters.  The first factor that
depends only on $N$, $p$, and $q$, is not present in the classical
formula \cite{AbramowitzS65} for the weights; we include it for
convenience since the lattice spacing for our nodes is $1/N$ rather
than being fixed.  In any case, since
\begin{equation}
\mathop{\prod_{m=0}^{N-1}}_{m\neq n} (x_{N,n}-x_{N,m})=\frac{(-1)^{N-1-n}}{N^{N-1}}n!(N-1-n)!\,,
\end{equation}
the weights may also be written in the form
\begin{equation}
w^{\rm Kraw}_{N,n}(p,q)=
e^{-NV^{\rm
Kraw}_N(x_{N,n};l)}\mathop{\prod_{m=0}^{N-1}}_{m\neq n}
|x_{N,n}-x_{N,m}|^{-1}
\end{equation}
where
\begin{equation}
V_N^{\rm Kraw}(x;l):=l x\hspace{0.2 in}\text{with}\hspace{0.2 in}
l:=\log\frac{q}{p}\,.
\end{equation}
\label{symbol:VKraw}
Note that in this case the function $V_N^{\rm Kraw}(x;l)$ is
coincidentally independent of $N$, so that $V^{\rm Kraw}(x;l) = l x$ and $\eta^{\rm Kraw}(x;l)\equiv 0$.  These weights are therefore of the
required form ({\em cf.}  (\ref{eq:weightform})) for our analysis.
Since for the dual family of discrete orthogonal polynomials we should
simply take the opposite sign of the function $V_N(x)$, we see that the
polynomials dual to the Krawtchouk polynomials with parameter $l$
are again Krawtchouk polynomials with parameter $-l$.  A number
of different involutions of the primitive parameters $p$ and $q$
correspond to changing the sign of $l$.  For example, one could
have $p\leftrightarrow 1/p$ and $q\leftrightarrow 1/q$, or simply
$p\leftrightarrow q$.  The latter involution is consistent with the
typical assumption that $0\le p\le 1$ and $p+q=1$.

For the typical case when $0\le p\le 1$ and $p+q=1$, the above
self-duality of the Krawtchouk polynomials implies that it is
sufficient in fact to consider $0 \le p\le 1/2$.  This fact was used in the paper 
\cite{DragnevS00}, where the equilibrium measure was explicitly
constructed for all $p$ in this range, and for all $c\in (0,1)$.  To
summarize the results, it has been shown that there is a single band
$I\subset [0,1]$, with endpoints $\alpha=\alpha(p,c)<\beta(p,c)=\beta$ for which there are
explicit formulae.  The behavior of the equilibrium measure in
$(0,1)\setminus I$ depends on the relationship between $c$ and $p$ in
the following way:
\begin{itemize}
\item If $0\le c < p$:
The intervals $(0,\alpha)$ and $(\beta,1)$ are both voids.
\item If $p < c<q$:
The interval $(0,\alpha)$ is a saturated region and the interval $(\beta,1)$
is a void.
\item  If $q < c\le 1$:
The intervals $(0,\alpha)$ and $(\beta,1)$ are both saturated regions.
\end{itemize}
This information supports our argument that the situation of having a
constraint active at both endpoints of the interval is generic with
respect to small perturbations of $c$.  The borderline cases of $c=p$
and $c=q$ are interesting also.  In the paper \cite{DragnevS00}
it is shown that 
\begin{equation}
\text{$\alpha\rightarrow 0$ as $c\rightarrow p$}\hspace{0.2 in}
\text{and}\hspace{0.2 in}
\text{$\beta\rightarrow 1$ as $c\rightarrow q$\,,}
\end{equation}
%\begin{equation}
%\begin{array}{rclcrcl}
%\alpha&\rightarrow &0&\hspace{0.2 in}\text{as}\hspace{0.2 in}&
%c\rightarrow p\,,\\
%\beta&\rightarrow &1&\hspace{0.2 in}\text{as}\hspace{0.2 in}&
%c\rightarrow q\,,
%\end{array}
%\end{equation}
and for $c=p$ the density $d\mu_{\rm min}^c/dx$ of the equilibrium
measure is equal to the average of the constraints at $x=0$, while for
$c=q$ it is equal to the average of the constraints at $x=1$.  
These are thus both special cases of the general result
stated in Proposition~\ref{prop:average}.

There exists an integral representation for the Krawtchouk
polynomials, and an exhaustive asymptotic analysis of the polynomials
has been carried out using this formula and the classical method of
steepest descent; see \cite{IsmailS98}.  Our formulae for the
leading-order terms agree with those of \cite{IsmailS98} in the
interior of all bands, voids, and saturated regions (the steepest
descent analysis is carried out with $z$ held fixed away from all band
edges and from the endpoints of the interval of accumulation of the
nodes).  The relative error obtained in \cite{IsmailS98} is typically
of the order $O(N^{-1/2})$, although it is stated that under some
circumstances this can be improved to $O(N^{-1})$ in some voids and
saturated regions.  The relative error estimates associated with the
asymptotic formulae presented in \S~\ref{sec:actualtheorems} thus
generally sharpen those of \cite{IsmailS98} in regions where the
$O(N^{-1/2})$ relative error bound is obtained.
%Steepest descent analysis of an
%integral is no doubt a more direct task than what we have been
%presenting here.  
It should be noted that while integral representations like that
analyzed in \cite{IsmailS98} are not available for more general
(nonclassical) discrete orthogonal polynomials, the methods to be
developed in \S~\ref{sec:preparation} and \S~\ref{sec:asymptotics} and
that lead to the general theorems stated in
\S~\ref{sec:actualtheorems} apply in absence of any such
representation.

\subsubsection{The Hahn and associated Hahn polynomials.}
\label{sec:Hahn}
Now we consider a semi-infinite lattice of equally-spaced nodes
\begin{equation}
x_{N,n}:=\frac{2n+1}{2N}\hspace{0.2 in}\text{for}\hspace{0.2 in}
n=0,1,2,\dots\,,
\label{eq:rawHahnnodes}
\end{equation}
and consider a corresponding three-parameter family of
weights \cite{AbramowitzS65}
\begin{equation}
w_{N,n}(b,c,d):=
\frac{N^{N-1}}{\Gamma(N)}\cdot\frac{\Gamma(b)\Gamma(c+n)\Gamma(d+n)}
{\Gamma(n+1)\Gamma(b+n)\Gamma(c)\Gamma(d)}\,,
\label{eq:Hahnraw}
\end{equation}
\label{symbol:Hahnraw}
where $b$, $c$, and $d$ are real parameters.  The prefactor depending
only on $N$ is included as a convenient
normalization factor that takes into account the fact that the lattice
spacing in \eqref{eq:rawHahnnodes} is $1/N$.

Although the measure corresponding to the weight function
\eqref{eq:Hahnraw} is supported on an infinite set, there are always
only a finite number of orthogonal polynomials.  For example, if one takes the
parameters $b$, $c$, and $d$ to be positive, then Stirling's formula
shows that the weight only decays for large $n$ if a certain
inequality is satisfied among $b$, $c$, and $d$, and then it decays
only algebraically, like $n^{-p}$ with the power $p$ depending on $b$,
$c$, and $d$.  Therefore for positive parameters the weight function
\eqref{eq:Hahnraw} has only a finite number of finite moments, and consequently only a finite
number of powers of $n$ may be orthogonalized.

We consider here a different way of arriving at a finite family of orthogonal polynomials starting from \eqref{eq:Hahnraw}.  If one takes a limit in the
parameter space, letting the parameter $c$ in \eqref{eq:Hahnraw} tend toward the negative integer $1-N$, then one
finds
\begin{equation}
w_{N,n}(b,1-N,d):=\lim_{c\rightarrow 1-N}w_{N,n}(b,c,d)= \left\{
\begin{array}{ll}
\displaystyle \frac{N^{N-1}}{\Gamma(N)}\binom{N-1}{n}\cdot
(-1)^n\cdot
\frac{\Gamma(b)\Gamma(d+n)}
{\Gamma(b+n)\Gamma(d)}
&\mbox{if $n\in{\mathbb Z}_N$,}\\\\
0 & \mbox{if $n\ge N$.}
\end{array}\right.
\end{equation}
The limiting weights are thus supported on ${\mathbb Z}_N$ rather than
on an infinite lattice and according to \eqref{eq:rawHahnnodes} the $N$ nodes $x_{N,0}<\dots<x_{N,N-1}$ are equally spaced with spacing 
$1/N$
and $x_{N,0}=1/(2N)$.  Therefore the node density function
is $\rho^0(x)\equiv 1$.
Note that the weights $w_{N,n}(b,1-N,d)$ are not positive for all $n\in\mathbb{Z}_N$ unless further conditions are placed on the remaining real parameters $b$ and $d$.  Insisting that $w_{N,n}(b,1-N,d)>0$ for
all $n\in\mathbb{Z}_N$ identifies two disjoint regions in the $(b,d)$-plane. 

One of these regions is delineated by the inequalities $d>0$ and $b<2-N$.  In this case, we refer to $w_{N,n}(b,1-N,d)$ as the {\em Hahn weight} and we call the corresponding polynomials the Hahn polynomials.
Let $P$ and $Q$ be positive parameters.  Setting $d=P$
and $b=2-N-Q$ in the limiting formula for $w_{N,n}(b,1-N,d)$, we
arrive at a simple formula for the Hahn weights:
\begin{equation}
w_{N,n}(b,1-N,d)=w_{N,n}^{\rm Hahn}(P,Q):=\frac{N^{N-1}}{\Gamma(N)}\cdot
\frac{\displaystyle\binom{n+P-1}{n}\binom{N+Q-2-n}{N-1-n}}{\displaystyle\binom{N+Q-2}{Q-1}}\,,\hspace{0.2 in}
\text{for $n\in{\mathbb Z}_N$\,.}
\label{eq:wHahn}
\end{equation}
\label{symbol:wHahn}
Note that by taking $P=Q=1$, the Hahn weights become independent
of $n$, so in this special case the Hahn polynomials are up to a factor the
(discrete) Tchebychev polynomials; this same family of polynomials arises
as a special case of the Krawtchouk polynomials with $p=q=1/2$.

The other region of the $(b,d)$-plane for which the weights $w_{N,n}(b,1-N,d)$ are positive for all $n\in \mathbb{Z}_N$ is delineated by the
inequalities $b>0$ and $d<2-N$.  In this case, we refer to $w_{N,n}(b,1-N,d)$ as the {\em associated Hahn weight} and we call the corresponding
polynomials the associated Hahn polynomials.  Again, let $P$ and 
$Q$ be positive parameters.  Setting $d=2-N-Q$ and $b=P$
in the limiting formula for $w_{N,n}(b,1-N,d)$, the associated Hahn weights  are
\begin{equation}
w_{N,n}(b,1-N,d)=w^{\rm Assoc}_{N,n}(P,Q):=\frac{N^{N-1}}{\Gamma(N)}
\cdot\frac{\Gamma(N)\Gamma(N+Q-1)\Gamma(P)}
{\Gamma(n+1)\Gamma(P+n)\Gamma(N-n)\Gamma(N+Q-1-n)}\,,\hspace{0.2 in}
\text{for $n\in{\mathbb Z}_N$\,.}
\label{eq:wAssoc}
\end{equation}
\label{symbol:wAssoc}

Note that
\begin{equation}
w^{\rm Hahn}_{N,n}(P,Q)w^{\rm Assoc}_{N,n}(P,Q)
\prod_{m\neq n}(x_{N,m}-x_{N,n})^2 = 1
\end{equation}
for all $n\in {\mathbb Z}_N$, and all $P>0$ and $Q>0$.  This means that the associated Hahn
polynomials are dual to the Hahn polynomials ({\em cf.} the general
definition (\ref{eq:dualweightsdefine}) of dual weights in
\S~\ref{sec:dual}).

%In fact, if one examines the current literature, one finds that the
%general formula \eqref{eq:Hahnraw} is hardly used at all, and more
%commonly the name ``Hahn polynomials'' is reserved for the two special
%cases described above.  

%\subsubsection{Hahn polynomials.}
%\label{sec:finiteHahn}

Writing the Hahn weights \eqref{eq:wHahn} in the form (\ref{eq:weightform}), we have
\begin{equation}
V_N^{\rm
Hahn}(x_{N,n};P,Q)=\frac{1}{N}\log\left(\frac{\Gamma(P)
\Gamma(N+Q-1)}{\Gamma(Nx_{N,n}+P-1/2)\Gamma(N(1-x_{N,n})+Q-1/2)}
\right)\,.
\end{equation}
\label{symbol:VNHahn}
The interesting case is when $P$ and $Q$ are large.  We
therefore set $P=NA+1$ and $Q=NB+1$ for $A$ and $B$ fixed
positive parameters, and from Stirling's formula, we then have
\begin{equation}
V_N^{\rm Hahn}(x;NA+1,NB+1) = V^{\rm Hahn}(x;A,B) +
\frac{\eta^{\rm Hahn}(x;A,B)}{N}
\end{equation}
where
\begin{equation}
V^{\rm Hahn}(x;A,B):=
A\log(A)+(B+1)\log(B+1)-(A+x)\log(A+x)-(B+1-x)\log(B+1-x)
\end{equation}
\label{symbol:VHahn}
and
\begin{equation}
\eta^{\rm Hahn}(x;A,B):=
\frac{1}{2}\log\left(\frac{A}{B+1}\right) + O\left(\frac{1}{N}\right)\,.
\end{equation}
\label{symbol:etaHahn}
The convergence is uniform for $x$ in compact subsets of ${\mathbb
C}\setminus ((-\infty,-A)\cup(B+1,+\infty))$.

\begin{remark}
The fact that the leading term in $\eta^{\rm Hahn}(x;A,B)$ is 
independent of $x$
can be traced back to the particular choice of the order one
terms in $P$ and $Q$ that we have made.  Other choices
consistent with the same leading-order scaling (say, simply taking
$P=NA$ and $Q=NB$) introduce genuine analytic $x$ dependence
into the leading term of the correction 
$\eta^{\rm Hahn}(x;A,B)$.
\end{remark}

For the associated Hahn weight \eqref{eq:wAssoc}, the case of $P=NA+1$ and
$Q=NB+1$ is also of interest.  By duality,
\begin{equation}
V_N^{\rm Assoc}(x;AN+1,BN+1)=-V_N^{\rm Hahn}(x;AN+1,BN+1)\,,
\end{equation}
\label{symbol:VNAssoc}
and therefore we also have $V^{\rm Assoc}(x;A,B)=-V^{\rm Hahn}(z;A,B)$
\label{symbol:VAssoc} at the level of the leading term as
$N\rightarrow\infty$.  According to Proposition~\ref{prop:flip}, if
the equilibrium measure corresponding to the function $V^{\rm
Hahn}(x;A,B)$ and the node density function $\rho^0(z)$ is known for
all values of the parameter $c$, then that corresponding to the
function $V^{\rm Assoc}(x;A,B)$ (and the same node density function)
is also known for all values of the parameter $c$, essentially by
means of the involution $c\leftrightarrow 1-c$.

We have computed the equilibrium measure corresponding to $V^{\rm
Hahn}(x;A,B)$ and $\rho^0(x)\equiv 1$ for $x\in (0,1)$.  To describe
it, we first define two positive constants $c_A$ and $c_B$ \label{symbol:cAcB}
by
\begin{equation}
\begin{array}{rcl}
  c_A &:=& \displaystyle \frac{-(A+B)+\sqrt{(A+B)^2+4A}}{2}\,, \\\\
  c_B &:=& \displaystyle \frac{-(A+B)+\sqrt{(A+B)^2+4B}}{2}\,. 
\end{array}
\label{eq:cAcB}
\end{equation}
It is direct to check that $0<c_A, c_B<1$ for $A,B> 0$, and $c_A<c_B$
if $0<A<B$.  Now let us assume that $A\le B$ (see the remark
below). Then $c_A\le c_B$, and we consider the three distinct
possibilities: $c\in(0,c_A)$, $c\in(c_A,c_B)$, or $c\in(c_B,1)$. It
turns out that in each of these cases, there is one band interval,
denoted by $(\alpha, \beta)\subset (0,1)$, on both sides of which are either
saturated regions or voids.
\begin{itemize}
\item
For $c\in (0,c_A)$, the interval $(\alpha,\beta)$ is the band, and the
intervals $(0,\alpha)$ and $(\beta,1)$ are voids.  We refer to this configuration
as {\em void-band-void}.
\item
For $c\in (c_A, c_B)$, the interval $(\alpha,\beta)$ is the band, $(0,\alpha)$
is a saturated region, and $(\beta,1)$ is a void.  We refer to this configuration
as {\em saturated-band-void}.
\item
For $c\in (c_B, 1)$, the interval $(\alpha,\beta)$ is the band, and the
intervals $(0,\alpha)$ and $(\beta,1)$ are both saturated regions.  We refer to
this configuration as {\em saturated-band-saturated}.
\end{itemize}
As in the Krawtchouk case, the critical values of $c=c_A$ or $c=c_B$ are
somewhat special because either $\alpha= 0$ or $\beta= 1$.

\begin{remark}
For the case when $A\ge B$, we have $c_B\le c_A$ and the mid-regime
for $c$ becomes the interval $(c_B, c_A)$. For $c\in (c_B, c_A)$, the
interval $(\alpha,\beta)$ is the band, $(0,\alpha)$ is a void, and $(\beta,1)$ is a
saturated region, and we refer to this configuration as {\em
void-band-saturated}.  However, there is a symmetry in this problem:
if one swaps $A\leftrightarrow B$ and $x\leftrightarrow (1-x)$, then the
field $\varphi(x)$ is changed only by a constant which can be absorbed
into the Lagrange multiplier $\ell_c$.  Therefore it is sufficient to
consider $A\le B$.
\end{remark}

For all values of $c$, the band edge points $\alpha$ and $\beta$ are the two
(real) solutions of the following quadratic equation in $X$:
\begin{equation}\label{eq;qeqLR}
  X^2-2\frac{A(A+B)+(A+B)(B-A+2)c+(B-A+2)c^2}{(A+B+2c)^2} X
  + \left(\frac{c^2+(A+B)c-A}{(A+B+2c)^2} \right)^2=0\,.
\end{equation}
It is straightforward to check that $\alpha, \beta\in [0,1]$ for all $0<c<1$ and 
$A,B>0$, and the formulae for $\alpha$ and $\beta$ are
\begin{equation}\label{eq;L}
  \alpha := \frac{(B-A+2)c^2+(A+B)(B-A+2)c+A(A+B)
  - 2\sqrt{D}}{(A+B+2c)^2}
\end{equation}
and
\begin{equation}\label{eq;R}
  \beta := \frac{(B-A+2)c^2+(A+B)(B-A+2)c+A(A+B)
  + 2\sqrt{D}}{(A+B+2c)^2}\,,
\end{equation}
where the discriminant is given by
\begin{equation}
D:=c(1-c)(A+c)(B+c)(A+B+c)(A+B+c+1)\,.
\end{equation}
Now we describe the density of the equilibrium measure
$d\mu^c_{\rm min}/dx$, assuming (without loss of generality according to
the remark above) that $A\le B$.
It is useful to introduce the following notation.  Let the positive
function $T(x)$ be defined by
\begin{equation}
T(x):=\sqrt{\frac{\beta-x}{x-\alpha}}\hspace{0.2 in}\text{for}\hspace{0.2 in}
\alpha<x<\beta\,,
\end{equation}
and define four positive constants by
\begin{equation}
\begin{array}{rclrcl}
  k_1 &:=& \displaystyle \sqrt{\frac{1+B-\alpha}{1+B-\beta}}\,, & 
  k_2 &:=& \displaystyle \sqrt{\frac{1-\alpha}{1-\beta}}\,, \\\\
  k_3 &:=& \displaystyle \sqrt{\frac{A+\alpha}{A+\beta}}\,, &
  k_4 &:=& \displaystyle \sqrt{\frac{\alpha}{\beta}}\,.
\end{array}
\end{equation}
\begin{theorem}
For the functions $V(x)\equiv V^{\rm Hahn}(x)$ and $\rho^0(x)\equiv 1$
on $x\in [0,1]$ the solution of the variational problem of
\S~\ref{sec:equilibrium} is given by the following formulae when the
parameters satisfy $A\le B$.  Let the constants $c_A$ and $c_B$ be
given by \eqref{eq:cAcB} and let $\alpha$ and $\beta$ be defined by
\eqref{eq;L} and \eqref{eq;R}.  If $0<c<c_A$ (void-band-void) then
for $\alpha\le x\le \beta$,
\begin{equation}
\frac{d\mu^c_{\rm min}}{dx}(x):=
\frac{1}{\pi c}\left[\arctan\left(k_2T(x)\right)+\arctan\left(k_3T(x)\right)
-\arctan\left(k_1T(x)\right)-\arctan\left(k_4T(x)\right)\right]\,,
\label{eq:Hahnequilibrium1}
\end{equation}
and $d\mu_{\rm min}^c/dx\equiv 0$ if $0\le x\le \alpha$ or $\beta\le x\le 1$.
The corresponding Lagrange multiplier is given by
\begin{equation}
\ell_c:=(\beta-\alpha)\left\{\left[\log(\beta-\alpha)-1\right]K^{(1)}_{\rm VBV}-
2\log(2)K^{(2)}_{\rm VBV} + 2K^{(3)}_{\rm VBV}\right\} + \varphi(\beta)
\end{equation}
where
\begin{equation}
\begin{array}{rcl}
K^{(1)}_{\rm VBV}&:=&\displaystyle \frac{k_1}{1+k_1}-\frac{k_2}{1+k_2}
-\frac{k_3}{1+k_3}+\frac{k_4}{1+k_4}\,,\\\\
K^{(2)}_{\rm VBV}&:=&\displaystyle\frac{k_1}{1-k_1^2}-\frac{k_2}{1-k_2^2}
-\frac{k_3}{1-k_3^2}+\frac{k_4}{1-k_4^2}\,,\\\\
K^{(3)}_{\rm VBV}&:=&\displaystyle\frac{\log(1+k_1)}{1-k_1^2}-
\frac{\log(1+k_2)}{1-k_2^2}
-\frac{\log(1+k_3)}{1-k_3^2}+\frac{\log(1+k_4)}{1-k_4^2}\,,
\end{array}
\end{equation}
and $\varphi(\cdot)$ is the external field given in terms of $V(\cdot)=V^{\rm Hahn}(\cdot;A,B)$ and
$\rho^0(\cdot)\equiv 1$ by \eqref{eq:fielddef}.  If $c_A<c<c_B$ (saturated-band-void),
then for $\alpha\le x\le \beta$,
\begin{equation}
\frac{d\mu^c_{\rm min}}{dx}(x):=
\frac{1}{\pi c}\left[\arctan\left(k_2T(x)\right)+\arctan\left(k_3T(x)\right)
-\arctan\left(k_1T(x)\right)+\arctan\left(k_4T(x)\right)\right]\,,
\label{eq:Hahnequilibrium2}
\end{equation}
and $d\mu_{\rm min}^c/dx\equiv 1/c$ if $0\le x\le \alpha$ and
$d\mu_{\rm min}^c/dx\equiv 0$ if $\beta\le x\le 1$.
The corresponding Lagrange multiplier is given by
\begin{equation}
\begin{array}{rcl}
\ell_c &:= &\displaystyle
(\beta-\alpha)\left\{\left[\log(\beta-\alpha)-1\right]K^{(1)}_{\rm SBV}-
2\log(2)K^{(2)}_{\rm SBV} + 2K^{(3)}_{\rm SBV}\right\} \\
\\
&&\displaystyle\hspace{0.3 in}+\,\,\, \varphi(\beta)
+2(\beta-\alpha)\log(\beta-\alpha) + 2\alpha-2\beta\log(\beta)\,,
\end{array}
\end{equation}
where
\begin{equation}
\begin{array}{rcl}
K^{(1)}_{\rm SBV}&:=&\displaystyle \frac{k_1}{1+k_1}-\frac{k_2}{1+k_2}
-\frac{k_3}{1+k_3}-\frac{k_4}{1+k_4}\,,\\\\
K^{(2)}_{\rm SBV}&:=&\displaystyle\frac{k_1}{1-k_1^2}-\frac{k_2}{1-k_2^2}
-\frac{k_3}{1-k_3^2}-\frac{k_4}{1-k_4^2}\,,\\\\
K^{(3)}_{\rm SBV}&:=&\displaystyle\frac{\log(1+k_1)}{1-k_1^2}-
\frac{\log(1+k_2)}{1-k_2^2}
-\frac{\log(1+k_3)}{1-k_3^2}-\frac{\log(1+k_4)}{1-k_4^2}\,.
\end{array}
\end{equation}
Finally, if $c_B<c<1$ (saturated-band-saturated), then for $\alpha\le x\le \beta$,
\begin{equation}
\frac{d\mu_{\rm min}^c}{dx}(x):=\frac{1}{c}+\frac{1}{\pi c}
\left[\arctan\left(k_3T(x)\right)+\arctan\left(k_4T(x)\right)-
\arctan\left(k_1T(x)\right)-\arctan\left(k_2T(x)\right)\right]\,,
\label{eq:Hahnequilibrium3}
\end{equation}
and $d\mu^c_{\rm min}/dx\equiv 1/c$ if $0\le x\le \alpha$ and $\beta\le x\le 1$.
The corresponding Lagrange multiplier is given by
\begin{equation}
\begin{array}{rcl}
\ell_c&:=&\displaystyle (\beta-\alpha)\left\{
\left[\log(\beta-\alpha)-1\right]K^{(1)}_{\rm SBS} -2\log(2)K^{(2)}_{\rm SBS}
+2K^{(3)}_{\rm SBS}\right\}\\\\
&&\displaystyle\hspace{0.3 in}+\,\,\,
\varphi(\beta)+
2-2\beta\log(\beta)-2(1-\beta)\log(1-\beta)\,,
\end{array}
\end{equation}
where
\begin{equation}
\begin{array}{rcl}
K^{(1)}_{\rm SBS}&:=&\displaystyle \frac{k_1}{1+k_1}+\frac{k_2}{1+k_2}
-\frac{k_3}{1+k_3}-\frac{k_4}{1+k_4}\,,\\\\
K^{(2)}_{\rm SBS}&:=&\displaystyle\frac{k_1}{1-k_1^2}+\frac{k_2}{1-k_2^2}
-\frac{k_3}{1-k_3^2}-\frac{k_4}{1-k_4^2}\,,\\\\
K^{(3)}_{\rm SBS}&:=&\displaystyle\frac{\log(1+k_1)}{1-k_1^2}+
\frac{\log(1+k_2)}{1-k_2^2}
-\frac{\log(1+k_3)}{1-k_3^2}-\frac{\log(1+k_4)}{1-k_4^2}\,.
\end{array}
\end{equation}
\label{theorem:hahn}
\end{theorem}

\begin{figure}[h]
\begin{center}
\mbox{\psfig{file=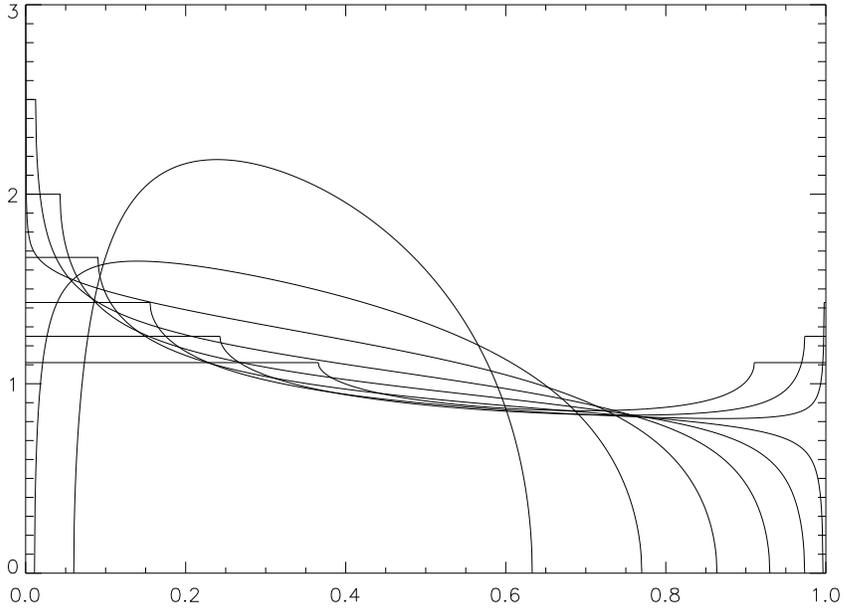,width=5 in}}
\end{center}
\caption{\em The density of the equilibrium measure for the Hahn polynomials
for parameter values $A=3$ and $B=7$.  Pictured are the measures for the
values $c=0.1,0.2,\dots,0.9$.}
\label{fig:hahn}
\end{figure}

The shapes of the equilibrium measures for the Hahn weights are
illustrated in Figure~\ref{fig:hahn}, which shows the way the measures
change as $c$ is varied for fixed $A<B$.
The proof of Theorem~\ref{theorem:hahn} is simply to check directly that the
following essential conditions are indeed verified:
\begin{itemize}
\item
The variational inequality \eqref{eq:voidinequality} holds in all voids,
the variational inequality \eqref{eq:saturatedregioninequality} holds
in all saturated regions, and the equilibrium condition \eqref{eq:equilibrium}
holds for $\alpha<x<\beta$.
\item
The measure satisfies the normalization condition \eqref{eq:Lagrange}.
\item
For $\alpha<x<\beta$, the measure has a density lying strictly between
the upper and lower constraints
\eqref{eq:constraints}.
\end{itemize}
So, rather than checking these conditions, we indicate some of the
techniques we used to deduce the formulae.  The equilibrium measure
may be computed either via an integral formula \cite{KuijlaarsV99} 
relating it to the asymptotics
of the recursion coefficients (which are known for the Hahn 
polynomials), or by directly solving the variational problem.
In Appendix~\ref{sec:hahnderivation} we follow similar
reasoning as in \cite{DeiftKM98} to derive the relevant formulae
recorded in Theorem~\ref{theorem:hahn}.

\begin{remark}
The Hahn polynomials have not been studied in the literature to the
same extent as the Krawtchouk polynomials.  There exists an integral
representation of the Hahn polynomials, but it is apparently more
difficult to analyze carefully than, for example, the corresponding
integral formula for the Krawtchouk polynomials studied 
in \cite{IsmailS98}.  We believe that the formulae
for the Hahn equilibrium measure presented in
Theorem~\ref{theorem:hahn} and the corresponding Plancherel-Rotach
type asymptotics that are formulated in \S~\ref{sec:actualtheorems},
are new in the literature.
\end{remark}

%\subsubsection{Associated Hahn polynomials.}
%\label{sec:assocHahn}

\section{Universal Asymptotic Properties of Discrete Orthogonal Polynomial Ensembles}
\label{sec:DOPensembles}
\subsection{Discrete orthogonal polynomial ensembles and particle statistics.}
\label{sec:DOPensemblesintro}
Consider the joint probability distribution of finding $k$ particles
at positions $x_1,\dots,x_k$ in $X_N$ to be given
by the following expression:
\begin{equation}
\begin{array}{rcl}
\displaystyle
\mathbb{P}(\mbox{there are particles at each of the nodes $x_1,\dots,x_k$})&=&
p^{(N,k)}(x_1,\dots,x_k)\\\\
&:=&\displaystyle
\frac{1}{Z_{N,k}}\prod_{1\le i<j\le k}(x_i-x_j)^2\cdot
\prod_{j=1}^k w(x_j)\,,
\end{array}
\label{eq:prob}
\end{equation}
\label{symbol:jointprobparticles}\label{symbol:Prob}
(we are using the symbol $\mathbb{P}(\mbox{event})$ to denote the
probability of an event) where $Z_{N,k}$ \label{symbol:partitionfunction}
is a normalization constant
(or partition function) chosen so that
\begin{equation}
\mathop{\sum_{x_1<\dots <x_k}}_{x_j\in X_N}
p^{(N,k)}(x_1,\dots,x_k)=1\,.
\label{eq:pnorm}
\end{equation}
Note that the particles are all indistinguishable from each other.
%$\{y_{N,j}\}_{j=1}^{k}$ that lie on $X_N$ whose density function is
%give by
%\begin{equation}
%  p^{(N,k)}(y_{N,1}, \cdots, y_{N,k})
%= \frac1{Z_{N,k}} \prod_{1\le i<j\le k} (y_{N,i}-y_{N,j})^2
%\prod_{j=1}^k w_N(y_{N,j})
%\end{equation}
%where $Z_{N,k}$ is the normalization constant.
The statistical ensemble associated with the density function
\eqref{eq:prob} is called a {\em discrete orthogonal polynomial ensemble}.

Discrete orthogonal polynomial ensembles arise in a number of specific
contexts (see, for example,
\cite{BorodinO01,Johansson00,Johansson01,Johansson02}), with particular choices
of the weight function $w(\cdot)$ related (in cases we are aware of)
to classical discrete orthogonal polynomials.  For instance,
\begin{itemize}
\item
The Meixner weight
\begin{equation}
w(x)=\binom{x+M-N}{x}q^x
\end{equation}
for $x=0,1,2,\dots$ arises in the directed last passage site percolation
model in the two-dimensional finite lattice $\mathbb{Z}_M\times\mathbb{Z}_N$
with independent geometric random variables as passage times for each site
\cite{Johansson00}.  The rightmost node occupied by a particle in the ensemble, $x_{\rm max}:={\rm max}_j x_j$, is a random variable having
the same distribution as the last passage time to travel from the site
$(0,0)$ to the site $(M-1,N-1)$.
\item
The Charlier weight
\begin{equation}
w(x)=\frac{t^x}{x!}
\end{equation}
for $x=0,1,2,\dots$ arises in the longest random word problem
\cite{Johansson01}.
\item
The Krawtchouk weight
\begin{equation}
w(x)=\binom{K}{x}p^xq^{K-x}
\end{equation}
for $x=0,1,\dots,K$ arises in the random domino tiling of the Aztec diamond
\cite{Johansson02}.
\item
The Hahn weight
\begin{equation}
w(x)=\binom{x+\alpha}{x}\binom{N+\beta-x}{N-x}
\end{equation}
for $x=1,2,3,\dots,N$ arises in the random rhombus tiling of a hexagon
\cite{Johansson02}.  See also \S~\ref{sec:tiling} for more details.
\end{itemize}
The first two cases (Meixner and Charlier) are examples of the
so-called {\em Schur measure} \cite{BorodinO01,Okounkov01} on the set
of partitions.  On the other hand, in special limiting cases the
Meixner and Charlier ensembles both become the {\em Plancherel
measure}, which describes the longest increasing subsequence of a
random permutation\footnote{Strictly speaking, this is not a discrete
orthogonal polynomial ensemble in the sense we have described because
as a consequence of the limiting process involved in the definition
the number of particles $k$ is not fixed in advance, but is itself
a random variable.}
\cite{BaikDJ99, BorodinOO00, Johansson01}.
Clearly it would be of some
theoretical interest to determine properties of the ensembles that are
more or less independent of the particular choice of weight function,
at least within some class.  Such properties are said to support the
conjecture of {\em universality} within the class of weight functions
under consideration.  Note that since the Meixner and Charlier weights
involve a semi-infinite lattice of nodes, a study of the corresponding
ensembles requires a generalization of the asymptotic methods we will describe in \S~\ref{sec:preparation} and \S~\ref{sec:asymptotics}.  Consequently universality results for such
ensembles will not be discussed here but will be developed in a
subsequent paper.

Some common properties of discrete orthogonal polynomial ensembles can
be read off immediately from the formula \eqref{eq:prob}.  For
example, the presence of the Vandermonde factor means that the
probability of finding two particles at the same site in $X_N$ is
zero.  Thus a discrete orthogonal polynomial ensemble always describes
an exclusion process.  This phenomenon is the discrete analogue of the
familiar {\em level repulsion} phenomenon in random matrix
theory. Moreover, due to the discreteness of the underlying space, the
particles are separated at least by the distance between consecutive
nodes. This strong exclusion due to the discreteness of the space
imposes the condition that the density of the states of the particles
has an upper bound, the limiting density of the nodes. This is the new
feature in the discrete orthogonal polynomial ensembles that is not
present in the orthogonal polynomial ensembles associated with
continuous weights ({\em i.e.} random matrix theory). Also, since the
weights are associated with nodes, the interpretation is that
configurations where particles are concentrated in sets of nodes where
the weight is larger are more likely.

Our goal will be to establish asymptotic formulae for
various statistics associated with the ensemble \eqref{eq:prob} for a
general class of weights in the continuum limit $N\rightarrow\infty$
subject to the basic assumptions enumerated in \S~\ref{sec:basicassumptions}
and the generic simplifying assumptions described in \S~\ref{sec:C3}.
%with the number of particles $k$ chosen so that for some fixed $c\in
%(0,1)$, we have $k=cN+\kappa$ where $\kappa$ remains bounded as
%$N\rightarrow\infty$.  We use the same assumptions on the nodes and
%weights as in the rest of the paper (see \S~\ref{sec:C1},
%\S~\ref{sec:C2}, and \S~\ref{sec:C3}).

Of basic interest is the \emph{$m$-point correlation function},
\label{symbol:RmNk} defined for $m\le k$ by
\begin{equation}\label{eq:correldef}
\begin{array}{rcl}
\displaystyle  R_m^{(N,k)}(x_1,\dots, x_m)
&:=&\displaystyle
\mathbb{P}(\text{there are particles at each of the nodes $x_1,\dots,x_m$})\\\\&=&\displaystyle
\mathop{\sum_{x_{m+1}<\dots < x_k}}_{x_j\in X_N}
p^{(N,k)}(x_1, \dots, x_k)\,,
\end{array}
\end{equation}

\begin{remark}
  In random matrix theory \cite{Mehta91,TracyW98} the correlation
  functions $R_m^{(N,k)}$ are usually introduced with a prefactor of
  $k!/(k-m)!$ which mediates between a density function for which
  particles (eigenvalues) are considered to be distinguishable
  (unordered) and statistics for which order is irrelevant.  Since we
  introduced $p^{(N,k)}(x_1,\dots,x_k)$ from the start with the
  interpretation that the particles are indistinguishable, this factor
  is not present in (\ref{eq:correldef}).
\end{remark}

In particular, the one-point function $R^{(N,k)}_1(x)$ \label{symbol:onepoint}
denotes the \emph{density of the states}, which is the probability
that there is a particle at $x$. One can also verify the following
interpretations: for any set $B\subset X_N$,
\begin{equation}
  \sum_{x\in B} R^{(N,k)}_1(x)
= \mathbb{E}( \text{number of particles in $B$})\,,
\end{equation}
and
\begin{equation}
  \mathop{\sum_{x<y}}_{x,y\in B} R^{(N,k)}_2(x,y)
= \mathbb{E}( \text{number of pairs of particles in
$B$}),
\end{equation}
where $\mathbb{E}$ denotes the expected value. \label{symbol:Expected}

The fundamental calculation of random matrix theory in the case
of so-called $\beta=2$ ensembles, due to Gaudin and Mehta (see, for example, \cite{Mehta91} or \cite{TracyW98}), shows that the correlation functions may equivalently be represented in the form
\begin{equation}
  R_m^{(N,k)}(x_1,\dots, x_m)
= \det \bigl( K_{N,k} (x_i, x_j) \bigr)_{1\le i,j\le m}\,,
\label{eq:determinantalformulacorrelations}
\end{equation}
where the so-called \emph{reproducing kernel} (Christoffel-Darboux
kernel) \label{symbol:reproducing} is defined for nodes $x$ and $y$ by
\begin{equation}\label{eq:reproducingdef}
K_{N,k}(x,y):=\sqrt{w(x)w(y)}\sum_{n=0}^{k-1}p_{N,n}(x)p_{N,n}(y)\,.
\end{equation}
Using the Christoffel-Darboux formula \cite{S91}, which holds for all
orthogonal polynomials, even in the discrete case, the sum on the
right telescopes.  Thus for distinct nodes $x\neq y$,
\begin{equation}
\begin{array}{rcl}
\displaystyle K_{N,k}(x,y)&=&\displaystyle \sqrt{w(x)w(y)}
\frac{\gamma_{N,k-1}}{\gamma_{N,k}}\cdot
\frac{p_{N,k}(x)p_{N,k-1}(y)-p_{N,k-1}(x)p_{N,k}(y)}{x-y} \\\\
&=&\displaystyle \sqrt{w(x)w(y)} \frac{\pi_{N,k}(x)\cdot
\gamma_{N,k-1}p_{N,k-1}(y)-
\gamma_{N,k-1}p_{N,k-1}(x)\cdot\pi_{N,k}(y)}{x-y}\\\\
&=&\displaystyle \sqrt{w(x)w(y)}
\frac{P_{11}(x;N,k)P_{21}(y;N,k)-P_{21}(x;N,k)P_{11}(y;N,k)}
{x-y}\,,\end{array}
\end{equation}
where the last line follows from Proposition~\ref{prop:solnrhp}.
Similarly, for  any node $x$,
\begin{equation}
K_{N,k}(x,x)=w(x)\left[P_{11}'(x;N,k)P_{21}(x;N,k)-P_{21}'(x;N,k)P_{11}(x;N,k)\right]\,.
\end{equation}
Note that the resulting formulae are expressed in terms of the
first column of the solution $\mat{P}(x;N,k)$ of
Interpolation Problem~\ref{rhp:DOP} for a single value of
$k$ and that
\begin{equation}\label{eq;kernelasPmatrix21}
\begin{array}{rcl}
\displaystyle  P_{11}(x;N,k)P_{21}(y;N,k)-P_{21}(x;N,k)P_{11}(y;N,k)
&  = &\displaystyle\bigl( \mat{P}(x;N,k)^{-1}\mat{P}(y;N,k)\bigr)_{21}\,,\\\\
\displaystyle P_{11}'(x;N,k)P_{21}(x;N,k)-P_{21}'(x;N,k)P_{11}(x;N,k)&=&
\displaystyle -\bigl(\mat{P}(x;N,k)^{-1}\mat{P}'(x;N,k)\bigr)_{21}\,.
\end{array}
\end{equation}
Therefore, the correlation functions are written explicitly in terms
of the discrete orthogonal polynomials associated with the nodes
$X_N$ and the weights $w_{N,n}=w(x_{N,n})$, and consequently
these formulae can be analyzed rigorously in an appropriate continuum limit
by using the methods we will present in detail in
\S~\ref{sec:preparation} and \S~\ref{sec:asymptotics}.

Consider a set $B\subset X_N$ and an integer $m$ with $0\le m\le{\rm
min}(\#B,k)$.  Another interesting statistic of a discrete orthogonal
polynomial ensemble is then \label{symbol:AmNkB}
\begin{equation}\label{eq;Am1}
A_m^{(N,k)}(B) := \mathbb{P} ( \text{there are precisely $m$
particles in the set $B$})\,,
\end{equation}
which vanishes automatically if $m>\#B$ by exclusion.
This statistic is also well-known to be expressible by the exact formula
\begin{equation}\label{eq;Am2}
  A_m^{(N,k)}(B) =
\frac1{m!} \biggl(-\frac{d}{dt}\biggr)^m\biggl|_{t=1} \det \bigl(
1-tK_{N,k}\bigl|_B \bigr)\,,
\end{equation}
where $K_{N,k}$ is the operator (in this case a finite matrix,
since $B$ is contained in the finite set $X_N$) acting in
$\ell^2(X_N)$ given by the kernel $K_{N,k}(x,y)$, and
$K_{N,k}\bigl|_B$ denotes the restriction of $K_{N,k}$ to
$\ell^2(B)$.

This is by no means an exhaustive list of statistics
that can be directly expressed in terms of the orthogonal
polynomials associated with the (discrete) weight $w(\cdot)$.
For example, one may consider the fluctuations and in particular
the variance of the number of particles in an interval $B\subset
X_N$.  The continuum limit asymptotics for this statistic were
computed in \cite{Johansson02} for the Krawtchouk ensemble (see
Proposition 2.5 of that paper) with the result that the
fluctuations are Gaussian; it would be of some interest to
determine whether this is special property of the Krawtchouk
ensemble, or a universal property of a large class of ensembles.
Also, there are convenient formulae for statistics associated with
the spacings between particles; the reader can find such formulae
in section 5.6 of the book \cite{Deift99}.

\subsection{Dual ensembles and hole statistics.}
\label{sec:DOPensemblesholes}
Since the nodes $X_N$ are finite in number, the distribution of the positions
$x_1,\dots,x_k$ of the particles naturally induces a
distribution of the positions $y_1,\dots,y_{\bar{k}}$ of the
\emph{holes} (that is, the nodes not occupied by particles).
Here $\bar{k}=N-k$, and $\{x_1,\dots, x_k\}\cup\{y_1,\dots,
y_{\bar{k}}\}=X_N$. It is interesting to determine the explicit
formula of the hole distribution. We will show that when the particle
locations $x_j$ are distributed according to the probability density function
$p^{(N,k)}(x_1,\dots,x_k)$ as in \eqref{eq:prob}, the density function
of the hole locations $y_j$ is always of the same form with only a different
choice of weight function.

Let us define \label{symbol:jointprobholes}
\begin{equation}
  \overline{p}^{(N,\bar{k})}(y_1,\cdots,y_{\bar{k}}) := \mathbb{P} (
  \text{there are holes at each of the nodes $y_1,\dots,
  y_{\bar{k}}$})\,.
\label{eq:pbardef}
\end{equation}
Given two complementary sets of nodes $\{x_1, \dots, x_k\} \cup \{y_1,
\dots, y_{\bar{k}}\}=X_N$, from
the definition (\ref{eq:pbardef}),
\begin{equation}\label{eq:pPpH1}
\begin{array}{rcl}
\displaystyle \overline{p}^{(N,\bar{k})}(y_1, \dots, y_{\bar{k}})&=&
 \displaystyle p^{(N,k)}(x_1,\cdots, x_k)\\\\
 & =& \displaystyle
\frac{1}{Z_{N,k}}\prod_{1\le i<j\le k}(x_i-x_j)^2\cdot\prod_{j=1}^k
w(x_j)\,.
\end{array}
\end{equation}
As
\begin{equation}
  \prod_{j=1}^k w(x_j) = C_N \prod_{j=1}^{\bar{k}}
  \frac1{w(y_j)},
  \qquad C_N:= \prod_{j=0}^{N-1}w_{N,j}\,,
\end{equation}
we find that
\begin{equation}
\begin{array}{rcl}
 \displaystyle \overline{p}^{(N,\bar{k})}(y_1,\dots,y_{\bar{k}})&=&
\displaystyle
\frac{C_N}{Z_{N,k}}\prod_{1\le i<j\le k} (x_i-x_j)^2
\cdot \prod_{j=1}^{\bar{k}} \frac1{w(y_j)} \\\\
& = &\displaystyle
\frac{C_N}{Z_{N,k}} \prod_{1\le i<j\le k}(x_i-x_j)^2
\cdot\prod_{j=1}^{\bar{k}} \prod_{\substack{n=0 \\ y_j\neq
x_{N,n}}}^{N-1} (y_j-x_{N,n})^2 \\\\
&&\displaystyle \,\,\,\cdot \,\,\,\prod_{j=1}^{\bar{k}} \biggl[
\frac1{w(y_j)} \prod_{\substack{n=0 \\ y_j\neq x_{N,n}}}^{N-1}
\frac1{(y_j-x_{N,n})^2} \biggr] \, .
\end{array}
\end{equation}
A little algebra shows that ({\em cf.} (9.42) of \cite{Baik99} or
Lemma 2.2 of \cite{Johansson01})
\begin{equation}
\begin{split}
  \prod_{1\le i<j\le k} |x_i-x_j|
\cdot \prod_{j=1}^{\bar{k}} \prod_{\substack{n=0 \\ y_j\neq
x_{N,n}}}^{N-1} |y_j-x_{N,n}|
  & = \prod_{1\le i<j\le k} |x_i-x_j| \cdot \prod_{j=1}^{\bar{k}}
  \prod_{\substack{i=1\\ i\neq j}}^{\bar{k}} |y_j-y_i| \cdot
  \prod_{j=1}^{\bar{k}}  \prod_{i=1}^{k} |y_j-x_i| \\
  & = \prod_{1\le i<j\le k} |x_i-x_j| \cdot \prod_{1\le i<j \le \bar{k}}
  |y_j-y_i|^2\cdot
  \prod_{j=1}^{\bar{k}}  \prod_{i=1}^{k} |y_j-x_i| \\
  & = D_N \prod_{1\le i<j \le\bar{k}}
  |y_j-y_i|
\end{split}
\end{equation}
where $D_N$ is the Vandermonde determinant of the nodes
\begin{equation}
  D_N := \prod_{0\le i<j\le N-1} |x_{N,i}-x_{N,j}|
\end{equation}
and the identity
\begin{equation}
D_N  = \prod_{1\le i<j\le k} |x_i-x_j| \cdot \prod_{1\le i<j\le \bar{k}}
  |y_i-y_j| \cdot \prod_{i=1}^k\prod_{j=1}^{\bar{k}} |x_i-y_j|
\end{equation}
is used in the last line. Therefore, the density of the holes is
given by
\begin{equation}
  \overline{p}^{(N,\bar{k})}(y_1, \cdots, y_{\bar{k}})
  = \displaystyle
\frac{1}{\overline{Z}_{N,\bar{k}}} \prod_{1\le i<j\le \bar{k}}
(y_i-y_j)^2 \cdot \prod_{j=1}^{\bar{k}} \overline{w}(y_j) \, ,
\end{equation}
where the normalization constant \label{symbol:Zbar} is
\begin{equation}
  \overline{Z}_{N,\bar{k}} = \frac{Z_{N,k}}{C_ND_N^2} = Z_{N,k} \prod_{j=0}^{N-1}
  \frac1{w_{N,j}} \cdot \prod_{0\le i<j\le N-1}
  \frac1{|x_{N,i}-x_{N,j}|^2},
\end{equation}
and the weight function is
\begin{equation}
  \overline{w}(y_j)
  = \frac1{w(y_j)} \prod_{\substack{n=0 \\ y_j\neq x_{N,n}}}^{N-1}
\frac1{(y_j-x_{N,n})^2} .
\end{equation}
Note that this new weight function is precisely the dual weight
defined in \eqref{eq:dualweightsdefine} of \S~\ref{sec:dual}. Hence
when the particles are distributed according a discrete orthogonal
polynomial ensemble, the holes are distributed according to the
discrete orthogonal polynomial ensemble corresponding to the dual
weights.  We will say that the ensembles governed by the density
functions $p^{(N,k)}(x_1,\dots,x_k)$ and
$\overline{p}^{(N,\bar{k})}(y_1,\dots,y_{\bar{k}})$ are dual to each
other.  Since dual ensembles correspond to weights of similar form,
but with the involutions $c\leftrightarrow 1-c$ and
$V(x)\leftrightarrow -V(x)$, their statistics are analyzed in exactly
the same way.  Therefore, the universality properties of the particle
distribution that we will establish below will automatically imply
corresponding universality properties of the of the hole distribution.

\subsection{Random rhombus tilings of a hexagon.  Relation to the Hahn and associated Hahn ensembles.}
\label{sec:tiling}
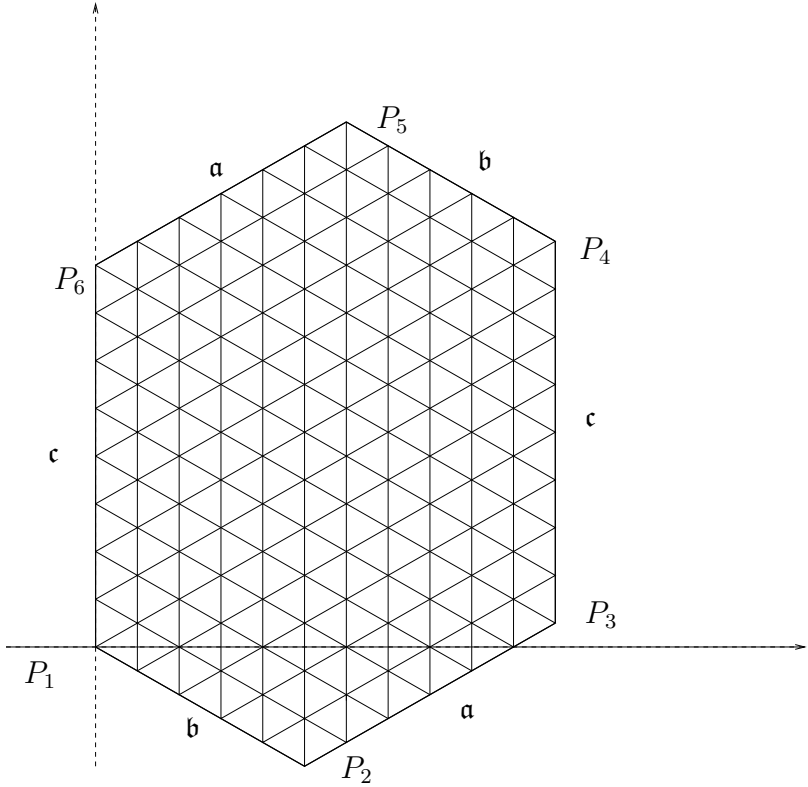
\begin{figure}[ht]
\begin{center}
\input{tilingbase.pstex_t}
\end{center}
% \centerline{\epsfig{file=tilingbase.eps, width=11cm}}
 \caption{\em
The $\mathfrak{abc}$-hexagon with vertices $P_1,\cdots, P_6$, and the lattice
$\mathcal{L}$} \label{fig-tilingbase}
\end{figure}

We briefly digress to describe a concrete occurrance in probability
theory of discrete orthogonal polynomial ensembles, in particular those
corresponding to the Hahn and associated Hahn weights.
Let $\mathfrak{a}$, $\mathfrak{b}$, and $\mathfrak{c}$ be positive
integers, \label{symbol:Hexabc} 
and consider the hexagon (see Figure \ref{fig-tilingbase})
with the following vertices written as points in the complex plane:
\begin{equation}
\begin{array}{rclrclrcl}
 P_1 &=& 0\,, & P_2&=&\mathfrak{b}e^{-i\pi/6}\,,& P_3&=&P_2+ 
\mathfrak{a}e^{i\pi/6}\,, \\\\
P_4&=& P_3+i\mathfrak{c}\,,& P_5&=&P_4+ \mathfrak{b}e^{5\pi i/6}\,,& P_6&=& i\mathfrak{c}\,.
\end{array}
\end{equation}
\label{symbol:vertices}
All interior angles of this hexagon are equal and measure $2\pi/3$
radians, and the lengths of the sides are, starting with the side
$(P_1,P_2)$ and proceeding in counter-clockwise order,
$\mathfrak{b},\mathfrak{a},\mathfrak{c},\mathfrak{b},\mathfrak{a},\mathfrak{c}$.
We call this the {\em $\mathfrak{abc}$-hexagon}.  Denote by $\mathcal{L}$ the
part of the set of lattice points (see Figure \ref{fig-tilingbase})
\begin{equation}
  \left\{ ke^{i\pi/6} + je^{-i\pi/6} \right\}_{k,j\in\mathbb{Z}}
   = \left\{ \frac{\sqrt 3}2n +  \frac{i}{2}n'  \right\}_{n,n' \in\mathbb{Z}}
\,.
\end{equation}
\label{symbol:hexlattice}
that lies within the hexagon, including the sides $(P_6,P_1)$,
$(P_1,P_2)$, $(P_2,P_3)$, and $(P_3,P_4)$, but excluding the sides
$(P_4,P_5)$ and $(P_5,P_6)$.  See Figure~\ref{fig-tilingbase}.

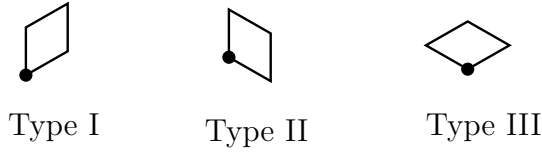
\begin{figure}[ht]
\begin{center}
\input{tilingrhombi.pstex_t}
\end{center}
% \centerline{\epsfig{file=tilingrhombi.eps, width=7cm}}
 \caption{\em The three types of rhombi; the position of each rhombus is
indicated with a dot.} \label{fig-tilingrhombi}
\end{figure}

Consider tiling the $\mathfrak{abc}$-hexagon with rhombi having sides of unit
length. Such rhombi come in three different types (orientations)
that we refer to as type I, type II, and type III; see
Figure~\ref{fig-tilingrhombi}.
Rhombi of types I and II are sometimes collectively called
\emph{horizontal rhombi}, while rhombi of type III are sometimes
called \emph{vertical rhombi}. The ``position" of each rhombus
tile in the hexagon is a specific lattice point in $\mathcal{L}$
defined as indicated in Figure \ref{fig-tilingrhombi}. 

MacMahon's formula \cite{MacMahon60} gives the total number of all
possible rhombus tilings of the $\mathfrak{abc}$-hexagon as the expression
\begin{equation}
  \prod_{i=1}^\mathfrak{a} \prod_{j=1}^\mathfrak{b} \prod_{k=1}^\mathfrak{c} 
\frac{i+j+k-1}{i+j+k-2}\,.
\end{equation}
Consider the set of all rhombus tilings equipped with uniform
probability. Hence we choose a tiling of the $\mathfrak{abc}$-hexagon at
random. It is of some current interest to determine the behavior
of various corresponding statistics of this ensemble in the limit
as $\mathfrak{a},\mathfrak{b},\mathfrak{c}\to\infty$.

\begin{figure}[ht]
\begin{center}
\epsfig{file=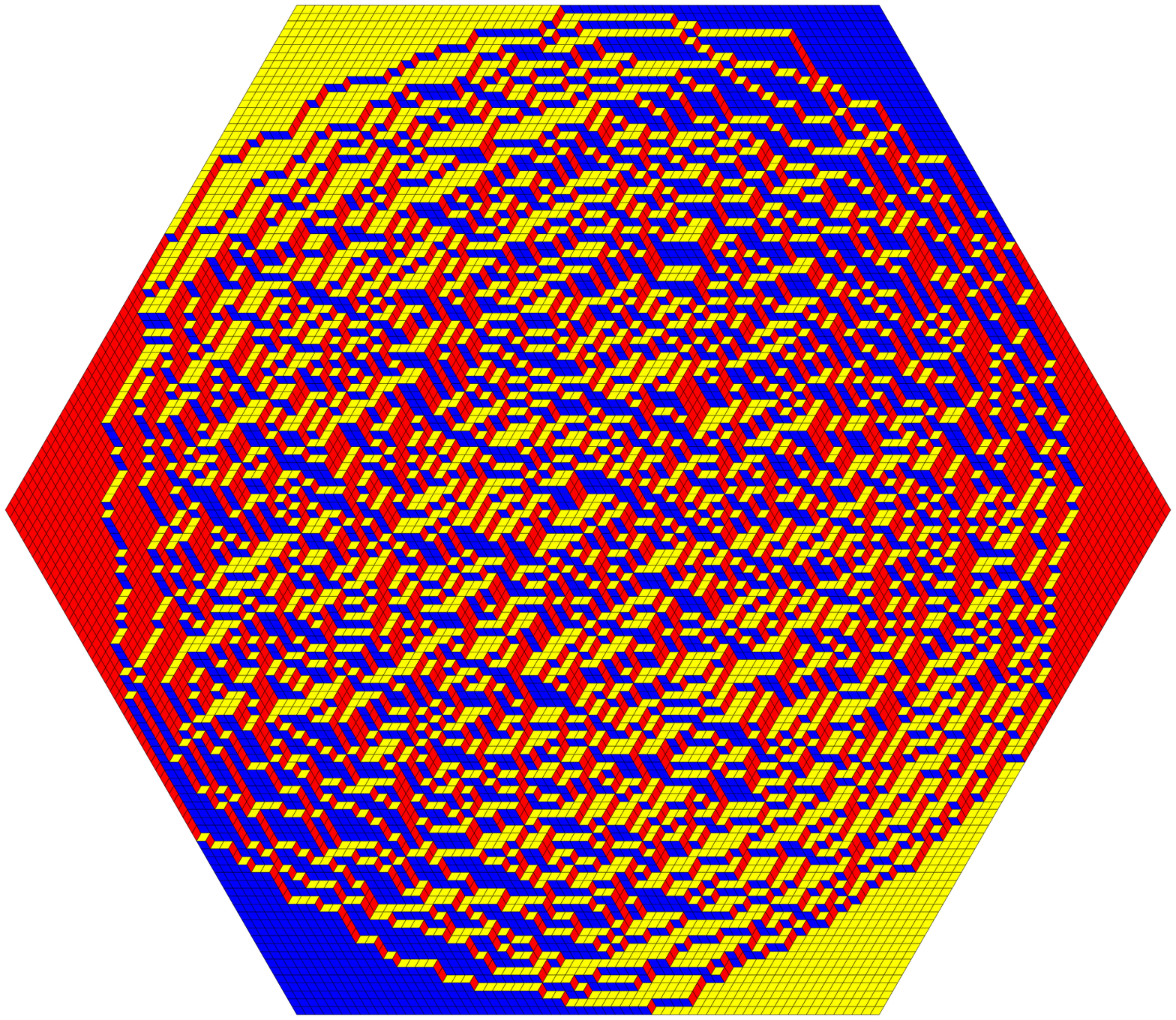,width=3 in,angle=30}
\end{center}
\caption{\em A rhombus tiling of a large $\mathfrak{abc}$-hexagon with 
$\mathfrak{a}=\mathfrak{b}=\mathfrak{c}=64$. Type I tiles are yellow, type II tiles are red, and type III tiles are blue.  Image provided by J. Propp.}
\label{fig:Propp}
\end{figure}

In the scaling limit of $n\to\infty$ where
\begin{equation}\label{eq;largehexagon}
  \mathfrak{a}= \mathfrak{A} n\,, \qquad \mathfrak{b}= \mathfrak{B} n\,, \qquad 
\mathfrak{c}=\mathfrak{C} n\,,
\end{equation}
\label{symbol:alphabetagamma}
with fixed $\mathfrak{A}, \mathfrak{B}, \mathfrak{C} >0$, the regions near the six
corners are ``frozen'' or ``polar'' zones ({\em i.e.}, regions in
which only one type of tile is present), while toward the center of
the hexagon is a ``temperate'' zone ({\em i.e.}, a region containing
all three types of tiles).  The random tiling shown in
Figure~\ref{fig:Propp} dramatically illustrates the two types of
regions, and the asymptotically sharp nature of the boundary between them.
Cohn, Larsen and Propp
\cite{CohnLP98} showed that in such a limit, upon scaling by $1/n$,
the expected shape of the boundary separating the polar zones from the
temperate zone is given by the inscribed ellipse. Moreover, the same
authors also computed the expected number of vertical rhombi in an
arbitrary set $U\in \mathbb{R}^2$. However, this calculation was
provided without specific error bounds.  Subsequently Johansson
\cite{Johansson00} proved a large deviation result for the boundary
shape, and also proved weak convergence of the marginal probability of
finding, say, a vertical tile near a given location in the temperate
zone.  The same paper also investigates a different tiling model,
namely the Aztec diamond tiling model for which a finer result is
proved.  It is proved that the fluctuation of the boundary between the
polar zones and the temperate zone in the Aztec tiling model is
governed (in a proper scaling limit) by the so-called Tracy-Widom law
of random matrix theory \cite{TracyW94}. One of the results implied by our
analysis of general discrete orthogonal polynomial ensembles (see Theorem~\ref{theorem:Hexagonboundaryfluctuations}) is that the same 
Tracy-Widom law holds for rhombus
tilings of the $\mathfrak{abc}$-hexagon.

The method of \cite{Johansson00} is to express the induced probability
for certain configurations of rhombi in the $\mathfrak{abc}$-hexagon or 
of rectangles in an Aztec diamond in terms of particular
discrete orthogonal polynomial ensembles. The weights corresponding to
the Aztec diamond are Krawtchouk weights, and those corresponding to
the $\mathfrak{abc}$-hexagon are Hahn or associated Hahn weights.
Johansson applied the classical steepest-descent method to the
integral representation of the Krawtchouk polynomials in order to
obtain various asymptotic results for the Krawtchouk ensemble.
However, even though the Hahn polynomials are also classical
polynomials, their integral representation does not seem to be so
straightforward to analyze asymptotically using the classical
steepest-descent method.  Hence questions of asymptotics for Hahn and
associated Hahn ensembles have not been adequately answered to date.
But as the Hahn and associated Hahn weights are special cases of
the general 
weights under study (see \S~\ref{sec:Hahn} for the relevant equilibrium measures), the universality results to be described below in \S~\ref{sec:universalityactualtheorems} apply to the Hahn and associated Hahn
ensembles as special cases, and hence we will obtain new results for the
random rhombus tiling of the $\mathfrak{abc}$-hexagon (see \S~\ref{sec:Hextheorems} below).

%In \cite{Johansson00}, Johansson expresses the induced probability
%for a given configuration of vertical or horizontal rhombi on a
%given sublattice in terms of discrete orthogonal polynomial
%ensembles with Hahn or associated Hahn weights. Even though the
%Hahn weight is a classical weight, the relevant asymptotics for
%Hahn polynomials have not been previously established.  However,
%the asymptotics of the previous sections may now be applied to the
%special case of the Hahn polynomials, and this yields new results
%for the asymptotic properties of the hexagon tiling problem (see
%Theorems~\ref{theorem:hexagonstrong} and \ref{theorem:hexagonsine}
%below).

%\textcolor{blue}{\textsf{We may like to remove the subsequent
%discussion up to Proposition 2, and directly go to our results.
%How do you think ?}}

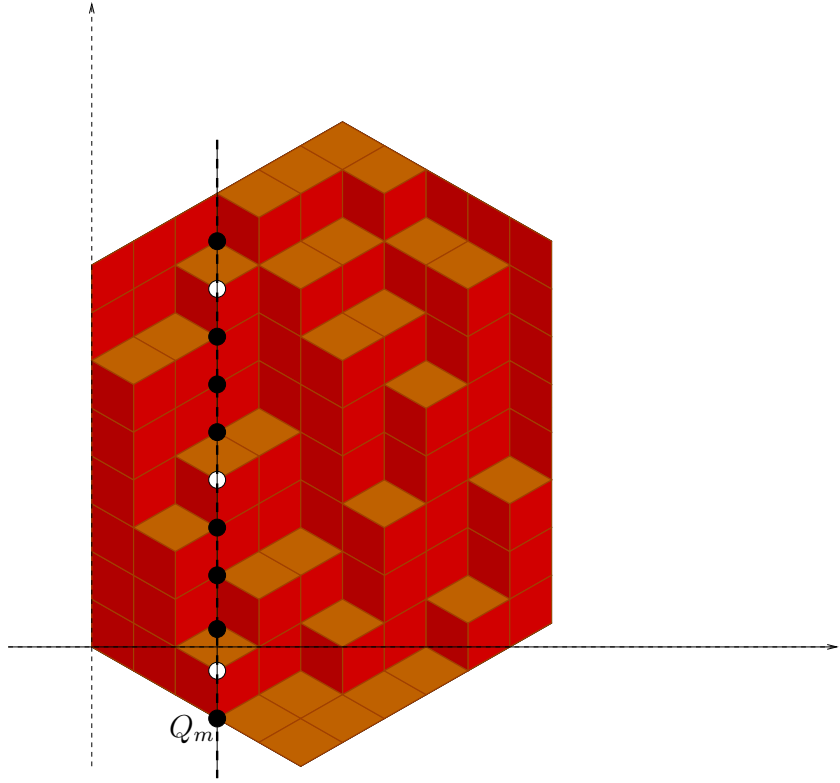
\begin{figure}[ht]
\begin{center}
\input{tilingandLm.pstex_t}
\end{center}
\caption{\em A rhombus tiling of the $\mathfrak{abc}$-hexagon, and
the lattice $\mathcal{L}_m$ when $m=3$; holes are represented by
white dots and particles are represented by black dots.}
\label{fig-tilingandLm}
\end{figure}

We first state the result of \cite{Johansson00} providing expressions
for probability density functions related to rhombus tilings of the
$\mathfrak{abc}$-hexagon in terms of discrete orthogonal polynomial
ensembles.  We will assume without loss of generality that
$\mathfrak{a}\ge \mathfrak{b}$ (by the symmetry of the hexagon, the
case when $\mathfrak{a}< \mathfrak{b}$ is completely analogous).
Consider the $m^{\rm th}$ vertical line of the lattice $\mathcal{L}$
counted from the left. We denote by $\mathcal{L}_m$
\label{symbol:hexlatticem} the intersection of this line and the
lattice $\mathcal{L}$. The number of points in $\mathcal{L}_m$ 
\label{symbol:cardLm} is
\begin{equation}
N=N(\mathfrak{a},\mathfrak{b},\mathfrak{c},m):=\mathfrak{c}+
\frac{\mathfrak{a}-\mathfrak{a}_m}{2}+
\frac{\mathfrak{b}-\mathfrak{b}_m}{2}\,,
\end{equation}
where \label{symbol:ambm}
\begin{equation}
\mathfrak{a}_m:=|m-\mathfrak{a}|\hspace{0.2 in}\text{and}\hspace{0.2 in}
\mathfrak{b}_m:=|m-\mathfrak{b}|\,.
\end{equation}
In a given tiling, the $N$ points in $\mathcal{L}_m$ correspond to
positions (in the sense defined above) of a number of rhombi of types
I, II, and III.  We call the positions of horizontal rhombi (types I
and II) the \emph{particles}, and the positions of vertical rhombi
(type III) the \emph{holes}.  See Figure~\ref{fig-tilingandLm} for an
example of $\mathcal{L}_m$ when $m=3$, illustrating the corresponding
particles and holes.

%The uniform probability distribution on the ensemble of tilings
%induces the probability distribution for finding particles and
%holes at particular locations in the one-dimensional finite
%lattice $\mathcal{L}_m$. A surprising result is due to Johansson
%\cite{Johansson00} which states that the induced probability
%distribution functions for holes and particles are both discrete
%orthogonal polynomial ensembles with Hahn and associated Hahn
%weight functions respectively (see \eqref{eq;Hahn1} and
%\eqref{eq;aHahn1}).

Let $Q_m$ \label{symbol:Qm} be the lowest point in the sublattice
$\mathcal{L}_m$.  On the sublattice $\mathcal{L}_m$, there are always
exactly $\mathfrak{c}$ particles, and $L_m:=N-\mathfrak{c}$
\label{symbol:holesinLm} holes.  Now, let $x_1<\cdots <
x_\mathfrak{c}$, where $x_j\in\{ 0,1,2,\dots, N-1\}$, denote the
(ordered) distances of the particles in $\mathcal{L}_m$ from $Q_m$,
and let $\xi_1< \cdots < \xi_{L_m}$, where $\xi_j \in\{ 0,1,2,\dots,
N-1\}$, denote the distances of the holes in $\mathcal{L}_m$ from
$Q_m$. In particular, we then have $\{ x_1,
\dots, x_\mathfrak{c}\}\cup \{ \xi_1, \dots, \xi_{L_m}\} =
\{0,1,2,\dots, N-1\}$. The uniform probability distribution on the
ensemble of tilings induces probability distributions for finding
particles and holes at particular locations in the one-dimensional
finite lattice $\mathcal{L}_m$. Let $\tilde{P}_m(x_1,\dots,
x_\mathfrak{c})$ \label{symbol:particleprob} denote the probability of
finding the particle configuration $x_1, \cdots, x_\mathfrak{c}$, and
let $P_m(\xi_1,\dots,
\xi_{L_m})$ \label{symbol:holeprob} 
denote the probability of finding the hole configuration
$\xi_1, \cdots,
\xi_{L_m}$.

\begin{prop}[Theorem 4.1 of \cite{Johansson00}]
  Let $\mathfrak{a},\mathfrak{b},\mathfrak{c}\ge 1$ be given integers
  with $\mathfrak{a}\ge \mathfrak{b}$. 
 Then
\begin{equation}
  \tilde{P}_m(x_1,\dots, x_\mathfrak{c})
 = \frac1{\tilde{Z}_m} \prod_{1\le j<k\le \mathfrak{c}} (x_j-x_k)^2
  \prod_{j=1}^{c} \tilde{w}(x_j)\,,
\end{equation}
where $\tilde{Z}_m$ is the normalization constant (partition
function), and where the weight function is the associated Hahn
weight (see \eqref{eq:wAssoc})
\begin{equation}
\tilde{w}(x):=
 w^{\rm Assoc}_{N,x}(\mathfrak{a}_m+1, b_\mathfrak{m}+1)
= \frac{\tilde{C}}{x!(\mathfrak{a}_m+x)!(N-x-1)!(N-x-1+\mathfrak{b}_m)!}\,,
\end{equation}
for a certain constant $\tilde{C}$. Also,
\begin{equation}
  P_m(\xi_1,\dots, \xi_{L_m}) =
  \frac1{Z_m} \prod_{1\le j<k\le L_m} (\xi_j-\xi_k)^2
  \prod_{j=1}^{L_m} w(\xi_j)\,,
\end{equation}
where $Z_m$ is the normalization constant, and where the weight
function is the Hahn weight (see \eqref{eq:wHahn})
\begin{equation}
  w(\xi)
:= w^{\rm Hahn}_{N,\xi}(\mathfrak{a}_m+1, \mathfrak{b}_m+1) =
C\frac{(\xi+\mathfrak{a}_m)!(N-\xi-1+\mathfrak{b}_m)!}{\xi!(N-\xi-1)!}\,,
\end{equation}
for a certain constant $C$.
\end{prop}

%\textcolor{blue}{\textsf{If we remove some parts in this section,
%it's probably up to here.}}

Together with the scaling \eqref{eq;largehexagon}, we set
\begin{equation}
  m=\tau n\,,
\end{equation}
for some fixed \label{symbol:tauparameter} $\tau>0$. 
The mean density of particles in $\mathcal{L}_m$ is then
\begin{equation}
\overline{c}:=\frac{\mathfrak{c}}{N}=\frac{2\mathfrak{C}}{2\mathfrak{C}
+\mathfrak{A}+\mathfrak{B}-
|\tau-\mathfrak{A}|-|\tau-\mathfrak{B}|}\,,
\label{eq:cbarval}
\end{equation}
and the mean density of holes in $\mathcal{L}_m$ is
\begin{equation}
c:=\frac{N-\mathfrak{c}}{N}=\frac{\mathfrak{A}+\mathfrak{B}-|\tau-\mathfrak{A}|-|\tau-\mathfrak{B}|}{2\mathfrak{C}+\mathfrak{A}+\mathfrak{B}-|\tau-\mathfrak{A}|-|\tau-\mathfrak{B}|}\,.
\label{eq:cval}
\end{equation}

\subsection{Results on asymptotic universality for general weights.}
\label{sec:universalityactualtheorems}
The following theorems all describe the asymptotic behavior as $N\rightarrow\infty$ of various
statistical quantities connected with the discrete orthogonal polynomial ensemble corresponding to nodes $X_N\subset [a,b]$ characterized
by the function $\rho^0(\cdot)$ and weights characterized by the 
function $V(\cdot)$.  These quantities, and the parameter $c$ (asymptotic value of $k/N$ where $k$ is the number of particles in the ensemble) are presumed to satisfy the same basic assumptions set forth in \S~\ref{sec:C1} and the simplifying assumptions set forth in \S~\ref{sec:C3}.  The theorems stated in this section will be proved below in \S~\ref{sec:universality}.

Let $\xi_N$ and $\eta_N$ be elements of a discrete subset $D_N$
of $\mathbb{R}$, such that $\max D_N-\min D_N$ remains bounded and the distance between neighboring points of $D_N$ converges to a constant as $N\rightarrow\infty$.  The expression \label{symbol:DSK}
\begin{equation}
S(\xi_N,\eta_N):=\frac{\sin(\pi(\xi_N-\eta_N))}{\pi(\xi_N-\eta_N)}
\end{equation}
is called the {\em discrete sine kernel} (``discrete'' reminds us that $\xi_N$ and $\eta_N$ lie in a discrete set $D_N$).  We extend the definition of the discrete sine kernel to the diagonal by setting
\begin{equation}
S(\xi_N,\xi_N):= 1\,.
\end{equation}
\begin{theorem}[Universality of the discrete sine kernel in bands]
\label{theorem;bulk}
Suppose that $x_1,\dots,x_l$ and $x_{l+1},\dots,x_m$ are disjoint sets of nodes in a fixed closed interval $F$ in the interior of any band $I$, and
denote by $\delta_N$ the distance between the two sets,
\begin{equation}
\delta_N:=\mathop{\min_{1\le i\le l}}_{l+1\le j\le m} |x_i-x_j|\,.
\end{equation}
Then,
\begin{equation}
R_m^{(N,k)}(x_1,\dots,x_m)=R_l^{(N,k)}(x_1,\dots,x_l)R_{m-l}^{(N,k)}(x_{l+1},\dots,x_m)+O\left(\frac{1}{N\delta_N}\right)\,.
\end{equation}
Fix $x$ in the interior of any band $I$, let
\begin{equation}
\delta(x):=\left[c\frac{d\mu_{\rm min}^c}{dx}(x)\right]^{-1}\,,
\label{eq:scaledparticlespacing}
\end{equation}
and for some integer $n\ge 1$
consider $\xi^{(1)}_N,\dots,\xi^{(n)}_N$ all to lie in a fixed bounded set $D\subset\mathbb{R}$ such that
\begin{equation}
x_j:=x+\xi^{(j)}_N\frac{\delta(x)}{N}\,,\hspace{0.2 in}j=1,\dots,n
\end{equation}
all satisfy $x_j\in X_N$ and $x_j\rightarrow x$ as $N\rightarrow\infty$.
Then there is a constant $C_{D,n}>0$ such that
for all $N$ sufficiently large,
\begin{equation}
\max_{\xi_N^{(1)},\dots,\xi_N^{(n)}\in D}
\left|R_n^{(N,k)}(x_1,\dots,x_n)-\left[\frac{c}{\rho^0(x)}\frac{d\mu_{\rm min}^c}{dx}(x)\right]^n\det(S(\xi_N^{(i)},\xi_N^{(j)}))_{1\le i,j\le n}\right|\le\frac{C_{D,n}}{N}\,.
\end{equation}
Thus particles separated by distances large compared to $1/N$ are asymptotically statistically independent, and the asymptotically nontrivial correlations among particles separated by distances comparable to $1/N$ are determined by the discrete sine kernel and the value of the one-point function.
\end{theorem}

Let the operator $\mathcal{S}(x)$ \label{symbol:DSKindex} act on
$\ell^2(\mathbb{Z})$ with the kernel (see, e.g. \cite{BorodinOO00})
\begin{equation}
\begin{array}{rcl}\displaystyle
  \mathcal{S}_{ij}(x)
  &:= &\displaystyle\frac{c}{\rho^0(x)}\frac{d\mu_{\rm min}^c}{dx}(x) S\left(\frac{c}{\rho^0(x)}\frac{d\mu_{\rm min}^c}{dx}(x)\cdot i,\frac{c}{\rho^0(x)}\frac{d\mu_{\rm min}^c}{dx}(x)\cdot j\right)\\\\
  &=&\displaystyle
  \frac{\displaystyle\sin \left(\frac{c}{\rho^0(x)}\frac{d\mu_{\rm min}^c}{dx}(x)\cdot\pi (i-j)\right)}{\pi (i-j)}\,,
  \end{array}
  \label{eq:sinekerneloperator}
\end{equation}
where $i,j\in\mathbb{Z}$.
\begin{theorem}[Asymptotics of local occupation probabilities in bands]
\label{theorem;bulk2}
Let $B_N\subset X_N$ be a set of $M$ nodes of the form
\begin{equation}
B_N=\{x_{N,j},x_{N,j+k_1},x_{N,j+k_2},\dots,x_{N,j+k_{M-1}}\}
\end{equation}
where $\#B_N=M$ is
held fixed as $N\rightarrow\infty$, and where
\begin{equation}
0<k_1<k_2<\cdots <k_{M-1}
\end{equation}
are fixed integers.  Set
$\mathbb{B}:=\{0,k_1,k_2,\dots,k_{M-1}\}\subset \mathbb{Z}$.
Suppose also that as $N\rightarrow\infty$, $x_{N,j}=\min
B_N\rightarrow x$ with $x$ lying in a band (and hence the same holds
for $x_{N,j+k_{M-1}}=\max B_N$).  Then, as $N\rightarrow\infty$,
\begin{equation}\label{eq;limittoS2}
  A_m^{(N,k)}(B_N)
  = \frac1{m!}\left(-\frac{d}{dt}\right)^m\biggl|_{t=1}
  \det\left( 1 -t \mathcal{S}(x)\bigl|_{\mathbb{B}} \right) +
  O\left(\frac{1}{N}\right)\,.
\end{equation}
\end{theorem}

\begin{theorem}[Uniform exponential bounds for the correlation functions in voids]
Let $F$ be a fixed closed interval in a void $\Gamma$ that is bounded away from all bands.  Then there is a constant $C_{F,m}>0$ such that for all $N$ sufficiently large,
\begin{equation}
  \max_{x_1,\dots,x_m\in X_N\cap F}\left|R^{(N,k)}_m(x_1,\cdots,x_m)\right|\le C_{F,m}\frac{e^{-mK_FN}}{N^m}\,,
\end{equation}
where the constant $K_F$ is defined by \begin{equation}\label{eq:voidKFconst}
  K_F:= \min_{z\in F} \biggl[ \frac{\delta E_c}{\delta \mu}(z) - \ell_c
  \biggr]\,.
\end{equation}
Note that $K_F>0$ because $F$ is closed and disjoint from the support of
the equilibrium measure $\mu_{\rm min}^c$.
 \label{theorem:onepointvoid}
\end{theorem}

For any $x\in (a,b)$, any $H>0$ and any $N>0$, let \label{symbol:Eint}
\begin{equation}
E_{\rm int}([A,B];x,H,N):=\mathbb{E}\left(
\text{number of particles at nodes $z$ of the form $\displaystyle z=
x+\frac{\xi_N}{H\sqrt{N}}$ with $A\le\xi_N\le B$}\right)\,.
\end{equation}
\begin{theorem}[Normal particle number distribution near interior local minima of $\delta E_c/\delta\mu$ in voids]
\label{theorem:normalvoid}
There is a finite set $Q$ such that for each
point $x$ in the interior of a void $\Gamma$ and with $x\not\in Q$, where
\begin{equation}
\frac{\delta E_c}{\delta\mu}(z)-\ell_c = W+H^2 \cdot (z-x)^2 + O\left((z-x)^3\right)
\label{eq:variationalTaylor}
\end{equation}
holds with $\mu=\mu_{\rm min}^c$ for some $H>0$ as $z\rightarrow x$, there
is a subsequence of integers $N$ tending to infinity for which we have
\begin{equation}
\frac{E_{\rm int}([C,D]\subset[A,B];x,H,N)}{E_{\rm int}([A,B];x,H,N)}=\frac{\displaystyle\int_C^De^{-\xi^2}\,d\xi}{\displaystyle\int_A^Be^{-\xi^2}\,d\xi} + O\left(\frac{1}{\sqrt{N}}\right)\,.
\label{eq:normalvoid}
\end{equation}
That is, the expected number of particles in a certain interval of size $1/\sqrt{N}$ near $x$ is given by a normal distribution.
\end{theorem}

\begin{remark}
Whether in the interior of a given void $\Gamma$ there may exist a local minimum of $\delta E_c/\delta\mu-\ell_c$ depends on the parameter $c$ and the nature of the functions $V(x)$ and $\rho^0(x)$ characterizing the equilibrium measure.
\end{remark}

The higher (multipoint) correlation functions for particles in a neighborhood of size $1/\sqrt{N}$ of the interior local minimum $x$ are smaller
in magnitude by a factor proportional to $1/\sqrt{N}$ than the one-point function.  This implies that although the one-point function is Gaussian, the statistics of distinct particles near $x$ are far from independent.

\begin{remark}
Another interesting possibility would be a local minimum of $\delta E_c/\delta\mu-\ell_c$ occurring at either endpoint $a$ or $b$ or the interval of accumulation of nodes, if this endpoint lies in a void.  But a direct
calculation gives, for $x$ in a void $\Gamma$,
\begin{equation}
\frac{d}{dx}\left[\frac{\delta E_c}{\delta\mu}(x)-\ell_c\right]=
{\rm P. V.}\int_a^b\frac{\rho^0(y)\,dy}{x-y}-2c\int_a^b\frac{\mu'(y)\,dy}{x-y} + V'(x)\,.
\end{equation}
Here $\mu=\mu_{\rm min}^c$.  The second integral is nonsingular because $x$ lies outside the support of the equilibrium measure.  As $x$ tends to an endpoint of $[a,b]$ in a void $\Gamma$, the latter two terms remain finite and the first term tends to $-\infty$ as $x\downarrow a$ and to $+\infty$ as $x\uparrow b$ (under our assumptions on $V(x)$ and $\rho^0(x)$).  Thus, neither endpoint can be a local minimum.
\end{remark}

The analogue of Theorem~\ref{theorem:onepointvoid} for saturated regions is the following.
\begin{theorem}[Uniform exponential bounds for the correlation functions in saturated regions]
\label{theorem:onepointsaturatedregion}
Let $F$ be a fixed closed interval in a saturated region $\Gamma$ that is bounded away from all bands.  Then there is a constant $C_{F,m}>0$ such that for all $N$ sufficiently large,
\begin{equation}
\max_{x_1,\dots,x_m\in X_N\cap F}\left|R_m^{(N,k)}(x_1,\dots,x_m)-1\right|\le C_{F,m}\frac{e^{-L_FN}}{N}\,,
\end{equation}
where the constant $L_F$ is defined by
\begin{equation}
L_F:=-\max_{z\in F}\left[\frac{\delta E_c}{\delta\mu}(z)-\ell_c\right]\,.
\label{eq:satKFconst}
\end{equation}
Note that $L_F>0$ because $F$ is a closed subinterval of an interval in which the the variational inequality \eqref{eq:saturatedregioninequality}
holds.
\end{theorem}

For $x\in (a,b)$, any $H>0$ and any $N>0$, let \label{symbol:Mint}
\begin{equation}
M_{\rm int}([A,B];x,H,N):=\#\left\{\text{nodes $z$ of the form $\displaystyle z=
x+\frac{\xi_N}{H\sqrt{N}}$
with $A<\xi_N<B$}\right\}
\end{equation}
which is asymptotically proportional to $\sqrt{N}$ for fixed $H$ and fixed $A<B$.
Then, the analogue of Theorem~\ref{theorem:normalvoid} for saturated regions is the following.
\begin{theorem}[Normal particle number deviations near interior local maxima of $\delta E_c/\delta\mu$ in saturated regions]
\label{theorem:normalsaturatedregion}
There is a finite set $Q$ such that for each
point $x$ in the interior of a saturated region $\Gamma$ with
$x\not\in Q$, where
\begin{equation}
\frac{\delta E_c}{\delta\mu}(z)-\ell_c = -W -H^2\cdot(z-x)^2 + O\left((z-x)^3\right)
\label{eq:interiorlocalmax}
\end{equation}
 holds with $\mu=\mu_{\rm min}^c$ for some $H>0$ as $z\rightarrow x$, there is
a subsequence of integers $N$ tending to infinity for which we have
\begin{equation}
\frac{M_{\rm int}([C,D]\subset[A,B];x,H,N)-E_{\rm int}([C,D]\subset[A,B];x,H,N)}
{M_{\rm int}([A,B];x,H,N)-E_{\rm int}([A,B];x,H,N)} =
\frac{\displaystyle\int_C^De^{-\xi^2}\,d\xi}{\displaystyle\int_A^Be^{-\xi^2}\,d\xi}+O\left(\frac{1}{\sqrt{N}}\right)\,.
\end{equation}
That is, the deviation of the expected number of particles from the number of available nodes in a certain interval of size $1/\sqrt{N}$ near $x$ is given by a normal distribution.
\end{theorem}

\begin{remark}
It is not possible for a local maximum to occur at an endpoint of $[a,b]$ lying in a saturated region, since for $x$ in a saturated region $\Gamma$
\begin{equation}
\frac{d}{dx}\left[\frac{\delta E_c}{\delta\mu}(x)-\ell_c\right]=
-c\,{\rm P. V.}\int_a^b\frac{\mu'(y)\,dy}{x-y}+\int_a^b\frac{\rho^0(y)-c\mu'(y)}{x-y}\,dy + V'(x)
\end{equation}
where $\mu=\mu_{\rm min}^c$ and the second term is nonsingular  because the upper constraint is satisfied by the equilibrium measure in saturated regions.  The latter two terms remain finite as $x$ tends to an endpoint of $[a,b]$, but the first term tends to $+\infty$ as $x\downarrow a$ and to $-\infty$ as $x\uparrow b$.  This shows that a local maximum may not occur at either endpoint in saturated regions.
\end{remark}

The expression \label{symbol:AiryK}
\begin{equation}
A(\xi_N,\eta_N):=\frac{Ai(\xi_N)Ai'(\eta_N)-Ai'(\xi_N)Ai(\eta_N)}{\xi_N-\eta_N}
\label{eq:Airykernel}
\end{equation}
is called the {\em Airy kernel}.

\begin{theorem}[Universality of the Airy kernel near band edges adjacent to voids]
\label{theorem:Airyunivcorrvoid}
For each fixed $M>0$, each left band edge $\alpha$ separating the band from a void, and each positive integer $m$, there is a constant $G^m_\alpha(M)>0$ such that for sufficiently large $N$,
\begin{equation}
\mathop{\max_{x_1,\dots,x_m\in X_N}}_{\alpha-MN^{-1/2}<x_j<\alpha+MN^{-2/3}\,,\forall j}
\left|R_m^{(N,k)}(x_1,\dots,x_m)-\left[\frac{\left(\pi c B^L_\alpha\right)^{2/3}}{N^{1/3}\rho^0(\alpha)}\right]^m\det\left(A(\xi_N^{(i)},\xi_N^{(j)})\right)_{1\le i,j\le m}\right|\le\frac{G^m_\alpha(M)}{N^{(m+1)/3}}\,,
\end{equation}
where 
\begin{equation}
B^L_\alpha:=\lim_{x\downarrow\alpha}\frac{1}{\sqrt{x-\alpha}}\frac{d\mu_{\rm min}^c}{dx}(x) >0\,,
\label{eq:Bleft}
\end{equation}
and $\xi_N^{(j)}=-\left(N\pi c B^L_\alpha\right)^{2/3}(x_j-\alpha)$.  Similarly, for each fixed $M>0$, each right band edge $\beta$ separating the band from a void, and each positive integer $m$, there is a constant $G^m_\beta(M)>0$ such that for sufficiently large $N$,
\begin{equation}
\mathop{\max_{x_1,\dots,x_m\in X_N}}_{\beta-MN^{-2/3}<x_j<\beta+MN^{-1/2}\,,\forall j}
\left|R_m^{(N,k)}(x_1,\dots,x_m)-\left[\frac{\left(\pi c B^R_\beta\right)^{2/3}}{N^{1/3}\rho^0(\beta)}\right]^m\det\left(A(\xi_N^{(i)},\xi_N^{(j)})\right)_{1\le i,j\le m}\right|\le\frac{G^m_\beta(M)}{N^{(m+1)/3}}\,,
\end{equation}
where 
\begin{equation}
B^R_\beta:=\lim_{x\uparrow\beta}\frac{1}{\sqrt{\beta-x}}\frac{d\mu_{\rm min}^c}{dx}(x)>0\,,
\label{eq:Bright}
\end{equation}
and $\xi_N^{(j)}=\left(N\pi c B^R_\beta\right)^{2/3}(x_j-\beta)$.
\end{theorem}

\begin{theorem}[Universality of the Airy kernel near band edges adjacent to saturated regions]
\label{theorem:Airyunivcorrsat}
For each fixed $M>0$, each left band edge $\alpha$ separating the band from a saturated region, and each positive integer $m$, there is a constant
$H_\alpha^m(M)>0$ such that for sufficiently large $N$,
\begin{equation}
\mathop{\max_{x_1,\dots,x_m\in X_N}}_{\alpha-MN^{-1/2}<x_j<\alpha+MN^{-2/3}\,,\forall j}
\left|R_m^{(N,k)}(x_1,\dots,x_m)-1+\frac{
\left(\pi \bar{c} \bar{B}^L_\alpha\right)^{2/3}}{N^{1/3}\rho^0(\alpha)}\sum_{j=1}^m
A(\xi_N^{(j)},\xi_N^{(j)})\right|\le\frac{H_\alpha^m(M)}{N^{2/3}}\,,
\end{equation}
where $\bar{c}:=1-c$,
\begin{equation}
\bar{B}^L_\alpha:=\lim_{x\downarrow\alpha}\frac{1}{\sqrt{x-\alpha}}\frac{c}{\bar{c}}\left[\frac{1}{c}\rho^0(x)-\frac{d\mu_{\rm min}^c}{dx}(x)\right]>0\,,
\label{eq:barBleft}
\end{equation}
and $\xi_N^{(j)}=-(N\pi \bar{c}\bar{B}^L_\alpha)^{2/3}(x_j-\alpha)$.  Similarly,
for each fixed $M>0$, each right band edge $\beta$ separating the band from a saturated region, and each positive integer $m$, there is a constant
$H_\beta^m(M)>0$ such that for sufficiently large $N$,
\begin{equation}
\mathop{\max_{x_1,\dots,x_m\in X_N}}_{\beta-MN^{-2/3}<x_j<\beta+MN^{-1/2}\,,\forall j}
\left|R_m^{(N,k)}(x_1,\dots,x_m)-1+\frac{
\left(\pi\bar{c}\bar{B}^R_\beta\right)^{2/3}}{N^{1/3}\rho^0(\beta)}\sum_{j=1}^mA(\xi_N^{(j)},\xi_N^{(j)})\right|\le\frac{H_\beta^m(M)}{N^{2/3}}\,,
\end{equation}
where again $\bar{c}=1-c$,
\begin{equation}
\bar{B}^R_\beta:=\lim_{x\uparrow\beta}\frac{1}{\sqrt{\beta-x}}\frac{c}{\bar{c}}\left[\frac{1}{c}\rho^0(x)-\frac{d\mu_{\rm min}^c}{dx}(x)\right]>0\,,
\label{eq:barBright}
\end{equation}
and $\xi_N^{(j)}=(N\pi\bar{c}\bar{B}^R_\beta)^{2/3}(x_j-\beta)$.
\end{theorem}

A statistic more interesting than the correlation functions near a
band edge is the limiting distribution of the location of the leftmost
or rightmost particle or hole.  It is well-known that the distribution
of the largest eigenvalue of a random matrix from the Gaussian unitary
ensemble converges, after proper centering and scaling, to a certain
one-parameter family of Fredholm determinants constructed from the
Airy kernel.  The dependence of the determinant on the parameter can
also be expressed in terms of a particular solution to the Painlev\'e
II equation \cite{TracyW94}. This universal distribution function is
known as the {\em Tracy-Widom distribution}. We claim that the
distribution of the location of the leftmost or rightmost particle or
hole has the same limit for general discrete orthogonal polynomial
ensembles of the type corresponding to the assumptions on the nodes,
weights, and equilibrium measures described in \S~\ref{sec:basicassumptions} and \S~\ref{sec:C3}.

Let $x_{\rm min}\in X_N$ and $x_{\rm max}\in X_N$
\label{symbol:xminmax} be the nodes occupied by the leftmost particle
and the rightmost particle respectively. Also denote by ${\cal
A}|_{[s,\infty)}$ \label{symbol:calA} the (trace class) integral
operator acting on $L^2[s,\infty)$ with the Airy kernel
\eqref{eq:Airykernel}. Recall the generic assumption that the equilibrium measure of
the $k$-particle ensemble has either a void or a saturated region
adjacent to each endpoint of the interval $[a,b]$ in which the
nodes accumulate.  Then we have the following result.

\begin{theorem}[Tracy-Widom distribution of the leftmost and rightmost particles]\label{thm:Airydet}
If the left endpoint $a$ is adjacent to a void $(a,\alpha)$, then for
each fixed $s\in\mathbb{R}$,
\begin{equation}\label{eq:Airydetlimit}
  \lim_{N\to\infty}
  \mathbb{P} \left( (x_{\rm min}-\alpha)\cdot(\pi NcB^L_\alpha)^{2/3} \ge -s
  \right)
  = \det(1-{\cal A}|_{[s,\infty)})\,,
\end{equation}
where $B^L_\alpha$ is defined by \eqref{eq:Bleft}.  If the right endpoint $b$ is adjacent to a void $(\beta,b)$,
then for each fixed $s\in\mathbb{R}$,
\begin{equation}\label{eq:Airydetlimit2}
  \lim_{N\to\infty}
  \mathbb{P} \left( (x_{\rm max}-\beta)\cdot(\pi NcB^R_\beta)^{2/3} \le s
  \right)
  = \det(1-{\cal A}|_{[s,\infty)})\,,
\end{equation}
where $B^R_\beta$ is defined by \eqref{eq:Bright}.
\end{theorem}

We also obtain a similar result for the
leftmost and the rightmost holes. Let $h_{\min}$ and $h_{\max}$ 
\label{symbol:hminmax} be
the nodes occupied by the leftmost and the rightmost hole
respectively.

\begin{theorem}[Tracy-Widom distribution of the locations of the leftmost and rightmost holes]\label{thm:Airydetforhole}
If the left endpoint $a$ is adjacent to a saturated region
$(a,\alpha)$, then for each fixed $s\in\mathbb{R}$,
\begin{equation}\label{eq:Airydetlimit3}
  \lim_{N\to\infty}
  \mathbb{P} \left( (h_{\rm min}-\alpha)\cdot(\pi N\bar{c}\bar{B}^L_\alpha)^{2/3} \ge -s
  \right)
  = \det(1-{\cal A}|_{[s,\infty)})\,,
\end{equation}
where $\bar{B}^L_\alpha$ is defined by \eqref{eq:barBleft} and $\bar{c}=1-c$.
If the right endpoint $b$ is adjacent to a saturated region
$(\beta,b)$, then for each fixed $s\in\mathbb{R}$,
\begin{equation}\label{eq:Airydetlimit4}
  \lim_{N\to\infty}
  \mathbb{P} \left( (h_{\rm max}-\beta)\cdot(\pi N\bar{c}\bar{B}^R_\beta)^{2/3} \le s
  \right)
  = \det(1-{\cal A}|_{[s,\infty)})\,,
\end{equation}
where $\bar{B}^R_\beta$ is defined by \eqref{eq:barBright} and $\bar{c}=1-c$.
\end{theorem}

\subsection{Random rhombus tilings of a hexagon.  Statistical asymptotics.}
\label{sec:Hextheorems}
The general asymptotic results stated in \S~\ref{sec:universalityactualtheorems} combined with
the specific calculations of the equilibrium measure for the Hahn weight in \S~\ref{sec:Hahn} imply several facts in the random tiling of the $\mathfrak{abc}$-hexagon.
Firstly, Theorems \ref{theorem;bulk},
\ref{theorem:onepointvoid} and \ref{theorem:onepointsaturatedregion}
predict the asymptotic behavior of the one-point correlation function, implying that as $n\to\infty$, the
one-dimensional lattice $\mathcal{L}_m$, rescaled to a
finite size independent of $n$, consists of three disjoint intervals:
one band, surrounded by two gaps (either saturated regions or voids,
depending on the parameters $\alpha$, $\beta$, $\gamma$, and $\tau$). The saturated regions and voids correspond
to the polar zones, while the central band is a section of the
temperate zone. Hence in particular, the endpoints of the band (see equations \eqref{eq;qeqLR}, \eqref{eq;L}, and \eqref{eq;R}, where $A:=\mathfrak{a}_m/N$ and $B:=\mathfrak{b}_m/N$ are functions of $\mathfrak{A}$, $\mathfrak{B}$, $\mathfrak{C}$, and $\tau$ only) when
considered as functions of $\tau$ for fixed $\mathfrak{A}$, $\mathfrak{B}$, and $\mathfrak{C}$ determine the typical shape of the
boundary between the polar and temperate zones of the rescaled
$\mathfrak{abc}$-hexagon. It may be checked that this curve, as calculated directly from the quadratic equation \eqref{eq;qeqLR}, coincides with the inscribed ellipse first shown to be the expected shape of the boundary by Cohn, Larsen and Propp \cite{CohnLP98}.

Moreover, we find that the one-point functions for particles and holes converge
\emph{pointwise} except at the band edges to
the equilibrium measures respectively for the associated Hahn weight corresponding to the value of $\overline{c}$ given in \eqref{eq:cbarval} and for the Hahn weight corresponding to the value of $c$ given in \eqref{eq:cval}, and we obtain a precise error bound. This
result thus improves upon those obtained in \cite{CohnLP98} and
\cite{Johansson00}. We expect that with additional analysis of the
same formulae it should be possible to show that the error is
locally uniform with respect to $\tau$, in which case the same
bounds should hold for more general regions $U\in \mathbb{R}^2$.
We state our result in this direction as follows.

\begin{theorem}[Strong asymptotics of the one-point function
  in the $\mathfrak{abc}$-hexagon] Consider holes on the line
  $\mathcal{L}_m$ of length $N$, where $m=\tau n$ and $\tau$ is fixed as
  $n\to\infty$.  The corresponding one-point function $R_1^{(N,cN)}(\xi)$ satisfies
  \begin{equation}
  R_1^{(N,cN)}(\xi)\rightarrow c\frac{d\mu_{\rm min}^c}{dx}(x)\hspace{0.2 in}
  \text{where}\hspace{0.2 in} x=\frac{\xi}{N}
  \end{equation}
  as $n\rightarrow\infty$ with $\mathfrak{a}=\mathfrak{A} n$,
  $\mathfrak{b}=\mathfrak{B} n$, and $\mathfrak{c}=\mathfrak{C} n$, and
  $\mathfrak{A}$, $\mathfrak{B}$, and $\mathfrak{C}$ are held fixed.  Here, the equilibrium measure is
  that corresponding to the Hahn weight with parameters $A=\mathfrak{a}_m/N$ and $B=\mathfrak{b}_m/N$ (see
  \eqref{eq:Hahnequilibrium1}, \eqref{eq:Hahnequilibrium2}, and
  \eqref{eq:Hahnequilibrium3} in \S~\ref{sec:Hahn}).  The convergence
  is uniform for $\xi=0,1,2,\dots,N-1$.  Note that the limit function
  $c d\mu_{\rm min}^c/dx (x)$ is identically equal to one in the polar
  zones near the vertices $P_2$ and $P_5$ and is identically equal to
  zero in the polar zones near the vertices $P_1$, $P_3$, $P_4$, and
  $P_6$.  The rate of convergence is uniformly exponentially fast (the
  error is of the order $O(e^{-Kn})$ for some $K>0$) for $\xi$ in any
  polar zone such that $x=\xi/N$ is uniformly bounded away from the
  temperate zone as $n\rightarrow\infty$.  For $\xi$ in the
  temperate zone such that $x=\xi/N$ is uniformly bounded away from
  all polar zones as $n\rightarrow\infty$ the rate of convergence is
  such that the error is uniformly of the order $O(1/n)$.
  \label{theorem:hexagonstrong}
\end{theorem}

In the temperate zone, in addition to the one-point function, we
can control all the multipoint correlation functions under proper
scaling (see Theorem \ref{theorem;bulk}). One consequence of this is
the following theorem concerning the scaling limit for the locations of
the holes (see Theorem \ref{theorem;bulk2}) in the line $\mathcal{L}_m$.

\begin{theorem}[Local occupation probabilities in the temperate zone of the
$\mathfrak{abc}$-hexagon]
  Consider a vertical line $\mathcal{L}_m$ of length $N$ in the
  $\mathfrak{abc}$-hexagon with $\mathfrak{a}=\mathfrak{A} n$,
  $\mathfrak{b}=\mathfrak{B} n$, $\mathfrak{c}=\mathfrak{C} n$ and $m=\tau n$ for
  fixed positive $\mathfrak{A}$, $\mathfrak{B}$, $\mathfrak{C}$, and $\tau$.  Let $x>0$
  be fixed such that $Nx\in \mathbb{Z}_N$ and such that the location
  $\xi=Nx$ units above $Q_m$ in $\mathcal{L}_m$ lies in the temperate zone
  bounded away from the expected boundary between the polar and
  temperate zones by a distance proportional to $n$.  Let $B= \{Nx,
  Nx+j_1, Nx+j_2, \dots, Nx+j_M\}$, where $\mathbb{B}= \{0, j_1, j_2,
  \cdots, j_M\}\subset \mathbb{Z}_N$ is a fixed set of integers.  Then
\begin{equation}
  \lim_{n\to\infty}
\mathbb{P}(\text{there are precisely $p$ holes in the set $B$})
= \frac1{p!} \biggl( -\frac{d}{dt}\biggr)^p\biggl|_{t=1} \det
\left( 1-t \mathcal{S}(x)|_{\mathbb{B}} \right),
\end{equation}
where $S(x)$ acts on $\ell^2(\mathbb{Z})$ with the kernel
\begin{equation}
  \mathcal{S}_{ij}(x) = \frac{\sin ( \pi q(x)(i-j))}{\pi (i-j)}\,,\hspace{0.2 in}\text{for $i,j\in \mathbb{Z}$\,,}
\end{equation}
where $q(x)=c d\mu_{\rm min}^c/dx(x)$ is the limiting one-point
function, or the density of states at $x$. \label{theorem:hexagonsine}
\end{theorem}

Finally we obtain the limiting distribution of the fluctuation of the
boundary separating the polar and temperate zones.  From
\ref{thm:Airydet}, we have the following result which was conjectured
in \cite{Johansson00}.  Recall that the Fredholm determinant
$\det\bigl( 1- {\cal A}|_{[x,\infty)}\bigr)$ (see \eqref{eq:HexTW}
below) has an alternative expression in terms of a particular solution
of the Painlev\'e II equation in the independent variable $x$, which
is referred to in random matrix theory as the Tracy-Widom law.

\begin{theorem}[Tracy-Widom distribution of extreme particles and holes in the $\mathfrak{abc}$-hexagon]
  Consider a vertical line $\mathcal{L}_m$ of length $N$ in the
  $\mathfrak{abc}$-hexagon with $\mathfrak{a}=\mathfrak{A} n$,
  $\mathfrak{b}=\mathfrak{B} n$, $\mathfrak{c}=\mathfrak{C} n$ and $m=\tau n$ for
  fixed positive $\mathfrak{A}$, $\mathfrak{B}$, $\mathfrak{C}$, and $\tau$.  Suppose
  further that $\tau$ is sufficiently small or sufficiently large that
  the polar zone at the top of $\mathcal{L}_m$ is a void for holes
  (equivalently, is saturated with particles).  Denote by $\xi_{*}$
  the height above the point $Q_m$ of the topmost hole in
  $\mathcal{L}_m$, and recall that for $\beta$
  defined by \eqref{eq;R} in \S~\ref{sec:Hahn} with $A=\mathfrak{a}_m/N$ and $B=\mathfrak{b}_m/N$, the limiting expected
  height above $Q_m$ of the boundary between the temperate and polar
  zones is $N\beta$.
  Then, for some constant $t>0$,
\begin{equation}
  \lim_{n\to\infty} 
\mathbb{P}\biggl( \frac{\xi_{*}-N\beta}{(tn)^{1/3}}
\le x \biggr)
= \det\bigl( 1- {\cal A}|_{[x,\infty)}\bigr)
\label{eq:HexTW}
\end{equation}
for each $x\in\mathbb{R}$, where ${\cal A}|_{[x,\infty)}$ is the Airy
operator acting on s $L^2[x,\infty)$ with the Airy kernel
\eqref{eq:Airykernel}.

The above result applies to the boundary between the polar zones near
the vertices $P_4$ and $P_6$ and the temperate zone.  The analogous
results hold for the boundary near $P_1$ and $P_3$ with the use of the
other endpoint $\alpha$ (see \eqref{eq;L} in \S~\ref{sec:Hahn}) in place of
$\beta$, a change of sign in the inequality, the interpretation of
$\xi_{*}$ as the location of the bottommost hole in $\mathcal{L}_m$
and a proper adjustment of the constant $t$.  Similarly, for the
boundary near $P_2$ and $P_5$ where the polar zones are voids for
particles (or packed with holes) the analogous results hold with the
interpretation of $\xi_{*}$ as the height above $Q_m$ of the
bottommost or topmost particle.
\label{theorem:Hexagonboundaryfluctuations}
\end{theorem}

%\begin{remark}
%All of the results we have written down for holes have analogous
%statements in terms of particles using the duality relation
%between the Hahn and associated Hahn weights that has been
%explained in \S~\ref{sec:universality}.
%\end{remark}

\begin{remark}
Similar results for domino tilings of the Aztec diamond are obtained
in \cite{Johansson01}.  In \cite{OkounkovR01}, 
a $q$-version or grand canonical
ensemble version of the uniform probability measure on the set of rhombus tilings of the
$\mathfrak{abc}$-hexagon is considered; thus the size of
the hexagon also becomes a random variable.  These authors computed the
correlation functions of holes in the temperate region that do not
necessarily lie along the same line, in a proper limit that
corresponds to the limit $N\rightarrow\infty$.  The result of this
calculation is a kernel built from the incomplete beta function,
referred to as the ``discrete incomplete beta kernel''. This kernel
reduces to the discrete sine kernel when the holes all lie along the
same line.  We expect that the same kernel should appear in the Hahn
ensemble if one computes the asymptotic correlation function for holes
lying in a two-dimensional region.  The Airy limit of the boundary of
the polar zones for this model was obtained by Ferrari and Spohn
\cite{FerrariS02}.
\end{remark}

\section{An Equivalent Riemann-Hilbert Problem}
\label{sec:preparation}
In this section we introduce a sequence of exact transformations
relating the matrix $\mat{P}(z;N,k)$
 to a matrix $\mat{X}(z)$ satisfying an equivalent
Riemann-Hilbert problem.  The Riemann-Hilbert problem characterizing
the matrix $\mat{X}(z)$ will be amenable to asymptotic analysis in the
joint limit of large degree $k$ and large parameter $N$.  This
asymptotic analysis will be carried out in \S~\ref{sec:asymptotics}.

\subsection{Choice of $\del$:  the transformation $\mat{P}(z;N,k)\to
\mat{Q}(z;N,k)$.}
\label{sec:choiceofdel}
It turns out that Interpolation Problem~\ref{rhp:DOP} will only be
amenable to analysis without any modification of the triangularity of
some of the residue matrices if the equilibrium measure never realizes
its upper constraint.  This is because the variational inequality
\eqref{eq:saturatedregioninequality} associated with this constraint
leads to exponential growth as $N\rightarrow\infty$ in each situation
that we wish to exploit the inequality \eqref{eq:voidinequality} to
obtain exponential decay.  This difficulty was recognized, for
example, in \cite{BorodinO02}, where for a specific weight it was
circumvented using representations of the corresponding polynomials in
terms of hypergeometric functions.  We need to handle the problem of
the upper constraint in full generality, and we will do so by using an
explicit transformation of the form \eqref{eq:PtoQ} to reverse the
triangularity of the residue matrices near only those poles where the
upper constraint is active, and leaving the triangularity of the
remaining residues unchanged.  The result of the change of variables
\eqref{eq:PtoQ} is a matrix $\mat{Q}(z;N,k)$ that depends on the choice
of a subset $\del\subset\mathbb{Z}_N$.  Our immediate goal is to
describe how the set $\del$ must be chosen to prepare for the
subsequent asymptotic analysis to be described in
\S~\ref{sec:asymptotics} in the limit $N\rightarrow\infty$.

The continuity of $d\mu_{\rm min}^c/dx$ (which follows from our basic
assumptions outlined in \S~\ref{sec:basicassumptions}, 
see also \S~\ref{sec:equilibrium})
along with the assumption
\eqref{eq:rho0nonzero}
implies that voids and saturated regions cannot be adjacent to each
other, but must always be separated by bands.  A band that lies
between a void and a saturated region (rather than between two voids
or between two saturated regions) will be called a {\em transition
band}.  In each transition band, we select arbitrarily a fixed point
$y_k$.  There are a finite number, say $M$, of transition bands, and
we label the points we select one from each in increasing order:
\label{symbol:yk} $y_1,\dots,y_M$.

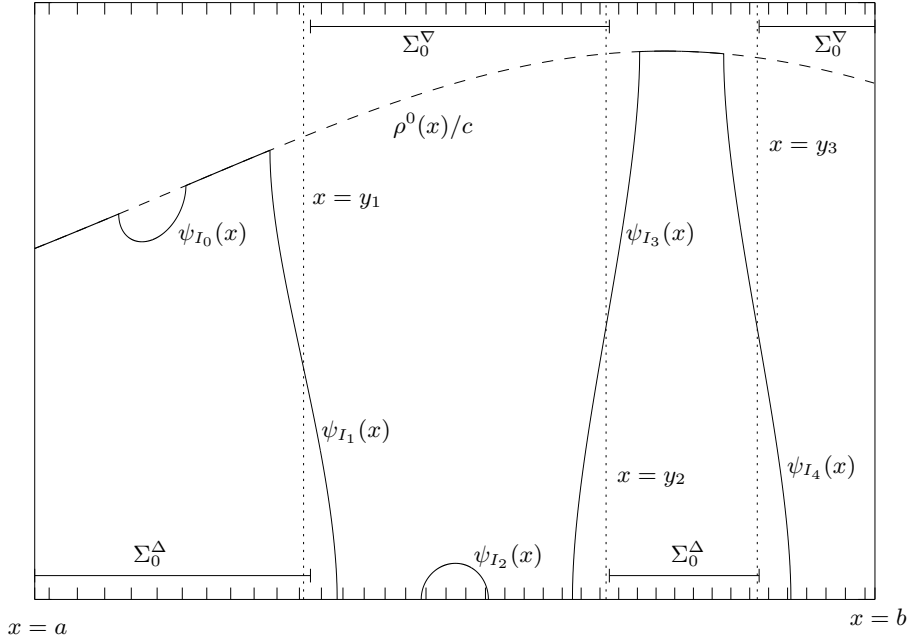
\begin{figure}[h]
\begin{center}
\input{minimizer.pstex_t}
\end{center}
\caption{\em A diagram showing the relation of a hypothetical equilibrium measure
$\mu^c_{\rm min}$ to the interval systems $\Sigma_0^\nab$ and
$\Sigma_0^\del$.  The nodes $X_N$ are indicated at the top and bottom
of the figure with tick marks; their density is proportional to the
upper constraint.  The endpoints of subintervals of $\Sigma_0^\nab$
and $\Sigma_0^\del$ are the transition points $Y_N$ that converge as
$N\rightarrow\infty$ to the fixed points $x=y_k$ whose positions
within each transition band are indicated with vertical dotted lines.
The analytic unconstrained components $\psi_I(x)$ of the density
$d\mu^c_{\rm min}/dx(x)$ are also indicated.}
\label{fig:intervals}
\end{figure}

With each $y_k$ we associate a sequence $\{y_{k,N}\}_{N=0}^\infty$
that converges to $y_k$ as $N\rightarrow\infty$.  Each element of the
sequence is defined by the quantization rule: \label{symbol:ykN}
\begin{equation}
N\int_a^{y_{k,N}}\rho^0(x)\,dx = \left\lceil N\int_a^{y_k}\rho^0(x)\,dx
\right\rceil
\label{eq:yquantize}
\end{equation}
where $\lceil u\rceil$ denotes the least integer greater than or equal
to $u$.  We call the points $y_{k,N}$ {\em transition points}, and use
the notation $Y_N$ for the set $\{y_{k,N}\}_{k=1}^M$, and $Y_\infty$
for the set $\{y_k\}_{k=1}^M$.  Since $\rho^0(x)$ is analytic and
nonzero in $(a,b)$, we have $y_{k,N}=y_k + O(1/N)$ as
$N\rightarrow\infty$.  Also, comparing with the condition
(\ref{eq:BS}) that defines the nodes $X_N$, we see that each of the
transition points $y_{k,N}$ asymptotically lies halfway between two
adjacent nodes.  Note that if only one constraint is active in
$[a,b]$, then there are no transition bands at all, and therefore no
transition points, so $Y_N=\emptyset$.  For all sufficiently large
fixed $N$, the transition points $y_{k,N}$ are ordered in the same way
as the points $y_k$.  For each $N$, we take the transition points in
$Y_N$ to be the common endpoints of two complementary systems
$\Sigma_0^\nab$ and $\Sigma_0^\del$ of open subintervals of $(a,b)$:
\label{symbol:Sigmanaughts}
\begin{definition}[The systems of subintervals $\Sigma_0^\nab$ and
$\Sigma_0^\del$]
The set $\Sigma_0^\nab$ is the union of those open subintervals
$(y_{k,N},y_{k+1,N})$ or $(a,y_{1,N})$ or $(y_{M,N},b)$ (or $(a,b)$ if
there are no transition points) that contain no saturated regions.
The set $\Sigma_0^\del$ is the union of those open subintervals
$(y_{k,N},y_{k+1,N})$ or $(a,y_{1,N})$ or $(y_{M,N},b)$ (or $(a,b)$ if
there are no transition points) that contain no voids.
\end{definition}
See Figure~\ref{fig:intervals}.
The sets $\Sigma_0^\nab$ and $\Sigma_0^\del$ depend on $N$ in a very
mild way, but they depend more crucially on the fixed parameter $c$ and
on the analytic functions $V(x)$ and $\rho^0(x)$.

With this notation, we now describe how we will choose the set $\del$
involved in the change of variables (\ref{eq:PtoQ}) from
$\mat{P}(z;N,k)$ to $\mat{Q}(z;N,k)$.  The set $\del$ will be taken to
contain precisely those indices $n$ corresponding to nodes $x_{N,n}$
contained in $\Sigma_0^\del$:
\begin{equation}
\del:=\{n\in {\mathbb Z}_N\hspace{0.1 in}\mbox{such that}
\hspace{0.1 in}x_{N,n}\in \Sigma_0^\del\}\,.
\label{eq:deldefine}
\end{equation}
In particular, this choice has the effect of reversing the
triangularity of the residue matrices at those nodes $x_{N,n}$ where
the upper constraint is active.  Note that $\#\del$ is roughly
proportional to $N$; we will define a rational constant $d_N$ by
writing \label{symbol:dN}
\begin{equation}
d_N:=\frac{\#\del}{N}\,.
\end{equation}
Note that $d_N$ has a limiting value $d$ as $N\rightarrow\infty$; for
technical reasons (see (\ref{eq:rhodef}) below) we will assume without
loss of generality (because we have considerable freedom in choosing
the points in $Y_\infty$) that $d\neq c$.

\subsection{Removal of poles in favor of jumps on contours:  the transformation $\mat{Q}(z;N,k)\to\mat{R}(z)$.}
The transformation in this section is based on an idea first used in
\cite{KamvissisMM03}. In that monograph, an analytic function was used to
simultaneously interpolate the residues of many poles at the pole
locations.  A generalization of this procedure involving two distinct
analytic interpolants was introduced in \cite{Miller02}.  The approach
we take in this section will also use two interpolants.

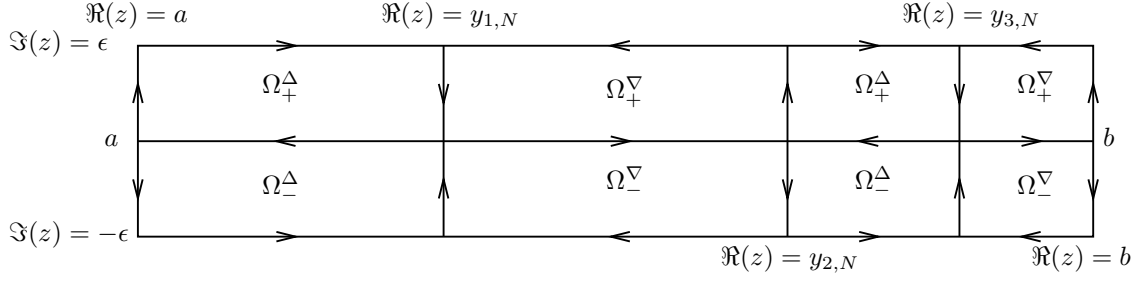
\begin{figure}[h]
\begin{center}
\input{Sigma.pstex_t}
\end{center}
\caption{\em The oriented contour $\Sigma$ and regions $\Omega_\pm^\nab$ and $\Omega_\pm^\del$.}\label{fig:Sigma}
\end{figure}

Note that by definition of the nodes $x_{N,j}\in X_N$ (see
\S~\ref{sec:C1}), and using
\eqref{eq:theta0def}, we have
\begin{equation}
ie^{-iN\theta^0(x_{N,n})/2}=-ie^{iN\theta^0(x_{N,n})/2}=(-1)^{N-1-n}\,,\hspace{0.2 in}
\text{for $N\in\nat$ and $n\in\mathbb{Z}_N$.}
\label{eq:dualinterp}
\end{equation}
Let $\epsilon>0$ \label{symbol:epsilon} be a fixed parameter
(independent of $N$) and consider the contour $\Sigma$
\label{symbol:Sigma} \label{symbol:Omegapmnabdel} illustrated in Figure~\ref{fig:Sigma}.
The figure is drawn to correspond to the hypothetical equilibrium
measure illustrated in Figure~\ref{fig:intervals}.  The contour
$\Sigma$ consists of the subintervals $\Sigma_0^\nab$ and
$\Sigma_0^\del$ and additional horizontal segments with
$|\Im(z)|=\epsilon$ and vertical line segments with $\Re(z)=a$,
$\Re(z)=b$, and $\Re(z)\in Y_N$.  We take the parameter $\epsilon$ to
be sufficiently small so that the contour $\Sigma$ lies entirely in
the region of analyticity of $V(x)$ and $\rho^0(x)$.  Further
restrictions will be placed on $\epsilon$ later on.

From the solution of Interpolation Problem~\ref{rhp:DOP} transformed
into the matrix $\mat{Q}(z;N,k)$ via (\ref{eq:PtoQ}) using the choice
of $\del$ given in (\ref{eq:deldefine}), we define a new matrix
$\mat{R}(z)$ \label{symbol:matrixR} as follows.  Set
\begin{equation}
\mat{R}(z):=\mat{Q}(z;N,k)\left(\begin{array}{cc}
1 & \displaystyle \mp ie^{\mp iN\theta^0(z)/2}e^{-NV_N(z)}\frac{\displaystyle
\prod_{n\in\del}(z-x_{N,n})}{\displaystyle\prod_{n\in\nab}
(z-x_{N,n})}\\\\
0 & 1\end{array}\right)\hspace{0.2 in}\mbox{for}\hspace{0.2 in}
z\in \Omega_\pm^\nab\,,
\label{eq:QtoR1}
\end{equation}
\begin{equation}
\mat{R}(z):=\mat{Q}(z;N,k)\left(\begin{array}{cc}
1 & 0\\\\
\displaystyle \mp ie^{\mp iN\theta^0(z)/2}e^{NV_N(z)}\frac{\displaystyle
\prod_{n\in\nab}(z-x_{N,n})}
{\displaystyle\prod_{n\in\del}(z-x_{N,n})} & 1
\end{array}\right)\hspace{0.2 in}\mbox{for}\hspace{0.2 in}
z\in \Omega_\pm^\del\,,
\label{eq:QtoR2}
\end{equation}
and for all other $z\in\mathbb{C}\setminus\Sigma$ set $\mat{R}(z):=\mat{Q}(z;N,k)$.

The significance of this explicit change of variables is that all poles have completely
disappeared from the problem.  Using the residue conditions (\ref{eq:polesinS}) and (\ref{eq:polesnotinS}) in conjunction with the ``interpolation'' identity (\ref{eq:dualinterp}), it is easy to check that
$\mat{R}(z)$ is an analytic function for $z\in \mathbb{C}\setminus\Sigma$
that takes continuous, and in fact analytic, boundary values on $\Sigma$.  In fact, $\mat{R}(z)$ can easily be seen to be the solution
of a Riemann-Hilbert problem relative to the contour $\Sigma$.  This problem
is sufficiently
similar to that introduced in \cite{FIK} for the continuous weight case
that it may, in principle, be analyzed by methods like those used in
\cite{DKMVZstrong,DeiftKMVZ99}.  We now proceed to describe the steps required for the corresponding analysis in the discrete case.

\subsection{Use of the
equilibrium measure:  the transformation $\mat{R}(z)\to\mat{S}(z)$.}
\label{sec:useofequilibrium}
\subsubsection{The complex potential $g(z)$ and the matrix $\mat{S}(z)$.}
The parameter $c$ and the analytic functions
$V(z)$ and $\rho^0(z)$ all influence the large $N$ behavior of the orthogonal polynomials. Thus, we recall the equilibrium measure $\mu^c_{\rm min}$ obtained in terms of these quantities in \S~\ref{sec:equilibrium}, and
for $x\in \Sigma_0^\nab\cup\Sigma_0^\del\subset[a,b]$ we define the piecewise real-analytic
density function \label{symbol:rho} as follows:
\begin{equation}
\rho(x) := \left\{\begin{array}{ll}
\displaystyle \frac{c}{c-d_N}\frac{d\mu^c_{\rm min}}{dx}(x)\,,&x\in\Sigma_0^\nab\\\\
\displaystyle \frac{c}{c-d_N}\left(\frac{d\mu^c_{\rm min}}{dx}(x)-
\frac{1}{c}\rho^0(x)
\right)\,, &x\in\Sigma_0^\del\,.\end{array}\right.
\label{eq:rhodef}
\end{equation}
We extend the domain of $\rho$ to the whole interval $[a,b]$, say by defining the function at its jump discontinuities to be the average of
its left and right limits.  Noting the
denominators in \eqref{eq:rhodef} we recall that we have assumed
without any loss of generality that $\lim_{N\rightarrow\infty}d_N\neq c$.
Since $\mu^c_{\rm min}$ is a probability measure, and since we may equivalently express $d_N$ in the form
\begin{equation}
d_N=\int_{\Sigma_0^\del}\rho^0(x)\,dx\,,
\end{equation}
we see that
\begin{equation}
\int_a^b\rho(x)\,dx = 1\,.
\label{eq:rhonorm}
\end{equation}
We also introduce the associated complex logarithmic
potential \label{symbol:gdef}
\begin{equation}
g(z):=\int_a^b\log(z-x)\rho(x)\,dx\,.
\label{eq:gdef}
\end{equation}
The logarithm in (\ref{eq:gdef}) is the principal branch; thus this function is analytic for $z\in\mathbb{C}\setminus (-\infty,b]$.  As a consequence of (\ref{eq:rhonorm}), we have $g(z)\sim \log(z)$ as $z\rightarrow\infty$.  The function $g(z)$ takes boundary values on $(-\infty,b]$ that are H\"older continuous with any exponent $\alpha<1$.

Recall the constant $\gamma$ defined in \eqref{eq:gammadefine}.  This
constant remains bounded as $N\rightarrow\infty$.  Consider the
transformation \label{symbol:matrixS}
\begin{equation}
\mat{S}(z):=e^{(N\ell_c+\gamma)\sigma_3/2}\mat{R}(z)
e^{(\#\del-k)g(z)\sigma_3}e^{-(N\ell_c+\gamma)\sigma_3/2}\,.
\label{eq:SfromR}
\end{equation}
Now, the identity (\ref{eq:rhonorm}) implies that the exponential $e^{(\#\del-k)g(z)}$ is analytic for $z\in\mathbb{C}\setminus [a,b]$ (in fact, since we are assuming a constraint to be active at both ends of the interval, the support of $\rho(x)$ is a closed subinterval of $(a,b)$ and we may replace $[a,b]$ by ${\rm supp}(\rho(x))$ in this statement).  Thus, like $\mat{R}(z)$, the matrix $\mat{S}(z)$ is
also analytic for $z\in\mathbb{C}\setminus\Sigma$, and the boundary values
taken on $\Sigma$ are continuous.  However, since $\mat{R}(z)z^{(\#\del-k)\sigma_3}\rightarrow\mathbb{I}$ as $z\rightarrow\infty$, we see that $\mat{S}(z)$ satisfies the normalization
condition
\begin{equation}
\mat{S}(z)=\mathbb{I} + O\left(\frac{1}{z}\right)\hspace{0.2 in}
\mbox{as $z\rightarrow\infty$.}
\end{equation}

\subsubsection{The jump of $\mat{S}(z)$ on the real axis.}
The point of introducing the equilibrium measure in this way is that
the matrix $\mat{S}(z)$ satisfies jump conditions across the voids,
bands, and saturated regions of $[a,b]$ that are analytically
tractable as a consequence of the variational inequalities that
$\mu^c_{\rm min}$ imposes on $\delta E_c/\delta\mu$ in the gaps.  To
describe these jump conditions, we first introduce for $z\in[a,b]$ the
functions \label{symbol:thetaofz} \label{symbol:phiofz}
\begin{equation}
\theta(z):=2\pi (d_N-c)\int_z^b\rho(s)\,ds\hspace{0.2 in}\mbox{and}
\hspace{0.2 in}
\phi(z):=-2\pi \kappa\int_z^b\rho(s)\,ds\,.
\label{eq:theta}
\end{equation}
Recalling the upper and lower constraints on the
equilibrium measure, the definition \eqref{eq:rhodef} implies that the function
$\theta(z)$ is real and nondecreasing for $z\in\Sigma_0^\nab$ and real
and nonincreasing for $z\in\Sigma_0^\del$.
Next, for $z\in\Sigma_0^\nab$ we define the function
\label{symbol:Tnab}
\begin{equation}
T_\nab(z):=2\cos\left(\frac{N\theta^0(z)}{2}\right)
\frac{\displaystyle\prod_{n\in\del}(z-x_{N,n})}{\displaystyle
\prod_{n\in\nab}(z-x_{N,n})}
\exp\left(N\left[\int_{\Sigma_0^\nab}\log|z-x|\rho^0(x)\,dx-
\int_{\Sigma_0^\del}\log|z-x|\rho^0(x)\,dx\right]\right)\,,
\label{eq:TAfuncdef}
\end{equation}
and for $z\in\Sigma_0^\del$ we define the function \label{symbol:Tdel}
\begin{equation}
T_\del(z):=2\cos\left(\frac{N\theta^0(z)}{2}\right)
\frac{\displaystyle\prod_{n\in\nab}(z-x_{N,n})}{\displaystyle
\prod_{n\in\del}(z-x_{N,n})}
\exp\left(-N\left[\int_{\Sigma_0^\nab}\log|z-x|\rho^0(x)\,dx-
\int_{\Sigma_0^\del}\log|z-x|\rho^0(x)\,dx\right]\right)\,.
\label{eq:TBfuncdef}
\end{equation}
Note that both $T_\nab(z)$ and $T_\del(z)$ are positive real-analytic
functions throughout their respective intervals of definition (the
cosine function cancels the poles contributed by the denominator in
each case).

Now, denoting the boundary value taken by $\mat{S}(z)$ on $\Sigma$ from the left by $\mat{S}_+(z)$ and that taken from the right by $\mat{S}_-(z)$, we can easily derive the relation
\begin{equation}
\mat{S}_+(z)=\mat{S}_-(z)
\left(\begin{array}{cc} e^{iN\theta(z)}e^{i\phi(z)} &
\displaystyle -iT_\nab(z)e^{\gamma-\eta(z) + \kappa(g_+(z)+g_-(z))}
\exp\left(N\left[\ell_c-\frac{\delta E_c}{\delta\mu}(z)\right]\right)\\\\
0 & e^{-iN\theta(z)}e^{-i\phi(z)}\end{array}\right)
\label{eq:Sjumpnab}
\end{equation}
holding for $z$ in any subinterval of $\Sigma_0^\nab$.  Similarly,
if $z$ is in any subinterval of $\Sigma_0^\del$, then
\begin{equation}
\mat{S}_+(z)=\mat{S}_-(z)\left(\begin{array}{cc}
e^{-iN\theta(z)}e^{-i\phi(z)} & 0 \\\\\displaystyle
iT_\del(z)e^{\eta(z)-\gamma-\kappa(g_+(z)+g_-(z))}\exp\left(N\left[
\frac{\delta E_c}{\delta\mu}(z)-\ell_c\right]\right) & e^{iN\theta(z)}e^{i\phi(z)}
\end{array}\right)\,.
\label{eq:Sjumpdel}
\end{equation}
Here, $g_+(z)+g_-(z)$ is the sum of the upper and lower boundary
values taken by the complex potential $g(z)$ on the real axis, and the
variational derivative is evaluated on the equilibrium measure
$\mu_{\rm min}^c$.

As $z$ varies within a gap $\Gamma$, the definition \eqref{eq:rhodef}
implies that the functions $\theta(z)$ and $\phi(z)$ remain constant.
In particular, to each gap $\Gamma$ we may assign a constant 
\label{symbol:phiGamma}
\begin{equation}
\phi_\Gamma:=\phi(z)\,,\hspace{0.2 in}\text{for $z\in\Gamma$\,.}
\label{eq:phiGamma}
\end{equation}
The constant values of $\theta(z)$ in the gaps have essentially already
been defined.  Recalling the definitions \eqref{eq:thetaGammajvoid},
\eqref{eq:thetaGammajsat}, and \eqref{eq:thetaGammaleftright} depending
on whether $\Gamma$ is (respectively) a void between two bands, a
saturated region between two bands, or one of the intervals
$(a,\alpha_0)$ or $(\beta_G,b)$, we see from \eqref{eq:yquantize}
that 
\begin{equation}
\theta(z)\equiv \theta_{\Gamma} \,\,\,\left(\displaystyle{\rm mod}\,\,\,
\frac{2\pi}{N}\right)\,,
\end{equation}
for $z$ in any gap $\Gamma$.  Note that the constants $\theta_\Gamma$
are by definition independent of the transition points $Y_N$.  Note
also that
\begin{equation}
e^{\pm iN\theta_\Gamma}e^{\pm i\phi_\Gamma} = 1\,,\hspace{0.2 in}
\text{when $\Gamma=(a,\alpha_0)$ or $\Gamma=(\beta_G,b)$\,,}
\end{equation}
because $\#\del$ and $k$ are both integers.

Now, for $z$ in a void $\Gamma$, the strict variational inequality
(\ref{eq:voidinequality}) holds.  Subject to the claim that
$T_\nab(z)$ remains bounded as $N\rightarrow\infty$ (this claim is
established in Proposition~\ref{prop:tasymp} below), we therefore see
that the jump matrix relating the boundary values in
(\ref{eq:Sjumpnab}) is exponentially close to the constant matrix
$e^{iN\theta_\Gamma\sigma_3}e^{i\phi_\Gamma\sigma_3}$ as
$N\rightarrow\infty$.  Similarly, for $z$ in a saturated region
$\Gamma$, the strict variational inequality
(\ref{eq:saturatedregioninequality}) holds, which shows that the jump
matrix relating the boundary values in (\ref{eq:Sjumpdel}) is
exponentially close to the constant matrix
$e^{-iN\theta_\Gamma\sigma_3}e^{-i\phi_\Gamma\sigma_3}$ in the limit
$N\rightarrow\infty$.

A band interval $I$ can be contained in $\Sigma_0^\nab$, in
$\Sigma_0^\del$, or (if it is a transition band) partly in
$\Sigma_0^\nab$ and partly in $\Sigma_0^\del$.  Throughout $I$, the
equilibrium condition (\ref{eq:equilibrium}) holds identically.  Thus,
for $z\in I\cap\Sigma_0^\nab$, we have a factorization of the jump
condition:
\begin{equation}
\mat{S}_+(z)=\mat{S}_-(z)
\left(\begin{array}{cc} e^{iN\theta(z)}e^{i\phi(z)} &
\displaystyle -iT_\nab(z)e^{\gamma-\eta(z)+\kappa(g_+(z)+g_-(z))}
\\ \\0 & e^{-iN\theta(z)}e^{-i\phi(z)}\end{array}\right) =
\mat{S}_-(z)\mat{L}_-(z)\mat{J}(z)\mat{L}_+(z)\,,
\label{eq:bandjumpnab}
\end{equation}
where for $z\in\Sigma_0^\nab$, \label{symbol:Lpm}
\begin{equation}\mat{L}_\pm(z):=
\left(\begin{array}{cc}T_\nab(z)^{\mp 1/2} & 0\\\\\displaystyle
iT_\nab(z)^{-1/2}e^{\eta(z)-\gamma-2\kappa g_\pm(z)}e^{\pm iN\theta(z)}
& T_\nab(z)^{\pm 1/2}\end{array}\right)\,,
\end{equation}
and \label{symbol:Jmatrix}
\begin{equation}
\mat{J}(z):=
\left(\begin{array}{cc}0 &\displaystyle -ie^{\gamma-\eta(z)+\kappa(g_+(z)+g_-(z))}\\\\
\displaystyle
-ie^{\eta(z)-\gamma-\kappa(g_+(z)+g_-(z))} & 0\end{array}\right)\,.
\label{eq:Jdef}
\end{equation}
As noted earlier, the function $T_\nab(z)$ is a strictly positive
analytic function throughout $I\cap\Sigma_0^\nab$, and we take
$T_\nab(z)^{\pm 1/2}$ to also be positive.  Similarly, for $z\in
I\cap\Sigma_0^\del$, (\ref{eq:Sjumpdel}) becomes
\begin{equation}
\displaystyle \mat{S}_+(z)=\mat{S}_-(z)
\left(\begin{array}{cc}e^{-iN\theta(z)}e^{-i\phi(z)} & 0\\\\
\displaystyle iT_\del(z)e^{\eta(z)-\gamma-\kappa(g_+(z)+g_-(z))} &
e^{iN\theta(z)}e^{i\phi(z)}
\end{array}\right) =
\mat{S}_-(z)\mat{U}_-(z)\mat{J}(z)^{-1}\mat{U}_+(z)\,,
\label{eq:bandjumpdel}
\end{equation}
where, for $z\in\Sigma_0^\del$, \label{symbol:Upm}
\begin{equation}
\mat{U}_\pm(z):=\left(\begin{array}{cc}
T_\del(z)^{\pm 1/2} & -iT_\del(z)^{-1/2}e^{\gamma-\eta(z)+2\kappa
g_\pm(z)}e^{\pm iN\theta(z)}\\\\ 0 & T_\del(z)^{\mp
1/2}\end{array}\right)\,,
\end{equation}
and $\mat{J}(z)$ is defined as in (\ref{eq:Jdef}).
Note that since $T_\del(z)$ is strictly positive for $z\in I\subset \Sigma_0^\del$, we are choosing the square roots $T_\del(z)^{\pm 1/2}$ to
also be positive.

\subsubsection{Important properties of the functions $T_\nab(z)$ and 
$T_\del(z)$.}  Here we establish for later use several properties of
$T_\nab(z)$ and $T_\del(z)$.  We first introduce the related function
$Y(z)$ defined by \label{symbol:Y(z)}
\begin{equation}
Y(z):=\frac{\displaystyle\prod_{n\in\del}(z-x_{N,n})}{\displaystyle
\prod_{n\in\nab}(z-x_{N,n})}
\exp\left(N\left[\int_{\Sigma_0^\nab}\log(z-x)\rho^0(x)\,dx-
\int_{\Sigma_0^\del}\log(z-x)\rho^0(x)\,dx\right]\right)\,,
\label{eq:Yfuncdef}
\end{equation}
for all $z$ in the domain of analyticity of $\rho^0(z)$ with
$\Im(z)\neq 0$.

We begin by explicitly relating $T_\nab(z)$, $T_\del(z)$, and $Y(z)$.
\begin{prop}[Analytic Properties of $T_\nab(z)$, $T_\del(z)$, and $Y(z)$]
There exists an open complex neighborhood $G$ of the closed interval
$[a,b]$ such that the following statements are true.
\begin{enumerate}
\item
$T_\nab(z)$ admits analytic continuation
to the domain $D_\nab:=(\mathbb{C}\setminus (\Sigma_0^\del\cup (-\infty,a] \cup
[b,+\infty)))\cap G$.
\item
$T_\del(z)$ admits analytic
continuation to the domain $D_\del:=(\mathbb{C}\setminus (\Sigma_0^\nab\cup
(-\infty,a]\cup[b,+\infty)))\cap G$.
\item
$Y(z)$ admits analytic
continuation to the domain $G\setminus [a,b]$.  
\item
The function $T_\nab(z)$ is real and positive for
$z\in\Sigma_0^\nab\subset D_\nab$ and the function $T_\del(z)$ is real
and positive for $z\in\Sigma_0^\del\subset D_\del$, and the
continuations of $T_\nab(z)$ and $T_\del(z)$ map the open domains
$D_\nab$ and $D_\del$ respectively into the cut plane
$\mathbb{C}\setminus (-\infty,0]$.
\item
The square roots $T_\nab(z)^{1/2}$ and $T_\del(z)^{1/2}$ exist as
analytic functions defined in the open domains $D_\nab$ and $D_\del$
respectively that are real and positive for $z\in \Sigma_0^\nab\subset
D_\nab$ and $z\in \Sigma_0^\del\subset D_\del$ respectively.
\item
We have the
identities
\begin{equation}
T_\nab(z)^{1/2}T_\del(z)^{1/2} = T_\nab(z)Y(z)^{-1}=T_\del(z)Y(z)=1+e^{-iN\theta^0(z)}
\,,\hspace{0.2 in}\mbox{for $z\in G$ with $\Im(z)>0$,}
\label{eq:Tsqrtidentity_up}
\end{equation}
and
\begin{equation}
T_\nab(z)^{1/2}T_\del(z)^{1/2}=T_\nab(z)Y(z)^{-1}=T_\del(z)Y(z)=1+e^{iN\theta^0(z)}\,,\hspace{0.2 in}\mbox{for $z\in G$ with $\Im(z)<0$.}
\label{eq:Tsqrtidentity_down}
\end{equation}
These formulae hold also on the real axis in the sense of boundary
values taken from the upper and lower half-planes.
\end{enumerate}
\label{prop:tanal}
\end{prop}
\begin{proof}
We take the domain $G$ to be contained in the domain of analyticity of
$\rho^0(z)$.  Let $z\in\Sigma_0^\nab$.  We then have
\begin{equation}
\begin{array}{l}\displaystyle
\int_{\Sigma_0^\nab}\log|z-x|\rho^0(x)\,dx-\int_{\Sigma_0^\del}\log|z-x|\rho^0(x)\,dx\,\,\, = \\\\
\displaystyle\hspace{0.3 in}\lim_{\epsilon\downarrow 0}\left[
\int_{\Sigma_0^\nab}\log(z\pm i\epsilon-x)\rho^0(x)\,dx-\int_{\Sigma_0^\del}\log(z \pm i\epsilon-x)\rho^0(x)\,dx\right]
\mp \frac{i\theta^0(z)}{2} \pm 2\pi iM\,,
\end{array}
\end{equation}
where
\begin{equation}
M = \int_{z<x\in\Sigma_0^\del}\rho^0(x)\,dx\,.
\end{equation}
The integral $M$ is a constant since $z\in\Sigma_0^\nab$, and by
virtue of the quantization condition (\ref{eq:yquantize}) it is an
integer.  This proves that $T_\nab(z)$ may be analytically continued
from any subinterval of $\Sigma_0^\nab$ to all of the open domain
$D_\nab$, and that the continuation does not depend on the particular
subinterval of $\Sigma_0^\nab$ from which the continuation is
performed.  The analytic continuation of $T_\del(z)$ to the open
domain $D_\del$ is obtained in a similar way.  The function $Y(z)$
clearly admits analytic continution to $z>b$, and for $z<a$ we have
\begin{equation}
\lim_{\epsilon\downarrow 0} \frac{Y(z+i\epsilon)}{Y(z-i\epsilon)}=
\exp\left(2\pi iN\left[\int_{\Sigma_0^\nab}\rho^0(x)\,dx 
-\int_{\Sigma_0^\del}\rho^0(x)\,dx\right]\right) = 1\,,
\end{equation}
where the last equality follows from the quantization condition
(\ref{eq:yquantize}) that determines the endpoints of the subintervals
of $\Sigma_0^\nab$ and $\Sigma_0^\del$.  This proves statements 1, 2, and 3.

These arguments immediately establish several of the identities
claimed in statement 6, namely that
$T_\nab(z)Y(z)^{-1}=T_\del(z)Y(z)=1+e^{-iN\theta^0(z)}$ for $\Im(z)>0$
and that $T_\nab(z)Y(z)^{-1}=T_\del(z)Y(z)=1+e^{iN\theta^0(z)}$ for
$\Im(z)<0$.  Combining these, one easily obtains the identities
\begin{equation}
T_\nab(z)T_\del(z)=(1+e^{-iN\theta^0(z)})^2\,,\hspace{0.2 in}
\text{for $\Im(z)>0$,}
\label{eq:Tnabdelup}
\end{equation}
and
\begin{equation}
T_\nab(z)T_\del(z)=(1+e^{iN\theta^0(z)})^2\,,\hspace{0.2 in}
\text{for $\Im(z)<0$.}
\label{eq:Tnabdeldown}
\end{equation}
Let $G_+$ and $G_-$ denote the intersections of the neighborhood $G$
with the upper and lower open half-planes respectively.  By choosing
$G$ to be sufficiently small but independent of $N$, we may ensure
(because the analytic function $\rho^0(z)$ is strictly positive for
$z\in [a,b]$) that for all $N>0$ the function $w=1+e^{\mp
iN\theta^0(z)}$ maps the open set $G_\pm$ into the open disk
$|w-1|<1$.  It follows that the image of $G_\pm$ under the map
$(1+e^{\mp iN\theta^0(z)})^2$ is an open set disjoint from the
negative real axis.  In particular, from \eqref{eq:Tnabdelup} we see
that the analytic functions $T_\nab(z)$ and $T_\del(z)$ have no zeros
in the open set $G_+$, and similarly from \eqref{eq:Tnabdeldown} we
see that neither function has any zeros in the open set $G_-$.  Now,
the strict positivity of $T_\nab(z)$ for $z\in\Sigma_0^\nab$ is a
simple consequence of the definition
\eqref{eq:TAfuncdef}, and that of $T_\del(z)$ for $z\in\Sigma_0^\del$
is a simple consequence of the definition \eqref{eq:TBfuncdef}.  So
while $T_\nab(z)$ has no zeros in $G$ away from the real axis or in
$\Sigma_0^\nab$, \eqref{eq:Tnabdelup} and \eqref{eq:Tnabdeldown} show
that the boundary values taken by $T_\nab(z)$ on any subinterval of
$\Sigma_0^\del$ from above or below have many double zeros.
Similarly, the boundary values taken by $T_\del(z)$ on any subinterval
of $\Sigma_0^\nab$ have many double zeros.  However, it is clear from
the preceding statements and from \eqref{eq:Tnabdelup} that if $C$ is
a contour homotopic to a subinterval of $\Sigma_0^\del$ that lies
(with the exception of its endpoints) in the open upper half-plane,
and if $C$ is close enough to the real axis, then $T_\nab(z)$ maps $C$
into the cut plane $\mathbb{C}\setminus (-\infty,0]$.  If instead $C$
lies in the lower half-plane, then \eqref{eq:Tnabdeldown} shows that
it is again mapped by $T_\nab(z)$ into the cut plane
$\mathbb{C}\setminus (-\infty,0]$ if it lies close enough to the real
axis.  Similar arguments show that contours in $D_\del$ homotopic to
subintervals of $\Sigma_0^\nab$ and close enough to the real axis are mapped
by $T_\del(z)$ into the cut plane $\mathbb{C}\setminus (-\infty,0]$.
This is sufficient to establish statement 4.

Statement 5 follows from statement 4 with an appropriate choice of the
square root.  The remaining identities in statement 6 are then
obtained by taking the square root of \eqref{eq:Tnabdelup} and
\eqref{eq:Tnabdeldown} and choosing the sign to be consistent with
taking the limit $z\rightarrow y_{k,N}$ in which the left-hand side is
positive.
\end{proof}

In a suitable precise sense the functions $T_\nab(z)$, $T_\del(z)$,
and $Y(z)$ may all be regarded as being approximately equal to one
when $N$ is large.  This is the content of the following proposition:
\begin{prop}[Asymptotic Properties of $T_\nab(z)$, $T_\del(z)$, and $Y(z)$]  
There exists an open complex neighborhood $G$ of the closed interval $[a,b]$
such that the following statements are true.
\begin{enumerate}
\item
(Asymptotics away from the boundary.)  For any fixed compact subset
$K\subset D_\nab := (\mathbb{C}\setminus (\Sigma_0^\del\cup
(-\infty,a]\cup[b,+\infty)))\cap G$, there exists a constant
$C^\nab_K>0$ for which the estimate
\begin{equation}
\sup_{z\in K}|T_\nab(z)-1|\le\frac{C^\nab_K}{N}\,,
\label{eq:Tnabout}
\end{equation}
holds for all sufficiently large $N$.  Similarly, for any fixed
compact subset $K\subset D_\del:=(\mathbb{C}\setminus
(\Sigma_0^\nab\cup (-\infty,a]\cup[b,+\infty)))\cap G$,  there
exists a constant $C^\del_K>0$ for which the estimate
\begin{equation}
\sup_{z\in K}|T_\del(z)-1|\le\frac{C^\del_K}{N}\,,
\label{eq:Tdelout}
\end{equation}
holds for all sufficiently large $N$.  Finally, for any fixed compact
subset $K\subset G\setminus[a,b]$, there is a constant $C_K>0$ for which
the estimate
\begin{equation}
\sup_{z\in K}|Y(z)-1|\le\frac{C_K}{N}\,,
\label{eq:Sout}
\end{equation}
holds for all sufficiently large $N$.
\item
(Asymptotics near $z=a$ and $z=b$.)
If $K\subset G$ is a compact neighborhood of $z=a$ and $\Sigma_0^\del$ is bounded away from $K$, then there is a constant $C^{\nab,a}_K>0$ and for each $\delta>0$ there is a constant $C^{\nab,a}_{K,\delta}$ such that for sufficiently large $N$,
\begin{equation}
\begin{array}{c}\displaystyle
\sup_{z\in K,|\arg(z-a)|<\pi}\left|T_\nab(z)-\frac{\sqrt{2\pi}e^{-\zeta_a}\zeta_a^{\zeta_a}}{\Gamma(\zeta_a+1/2)}\right|\le\frac{C^{\nab,a}_K}{N}\,,
\\\\\displaystyle
\sup_{z\in K,\delta\le |\arg(z-a)|\le\pi}
\left|Y(z)-\frac{\Gamma(1/2-\zeta_a)}{\sqrt{2\pi}e^{\zeta_a}(-\zeta_a)^{-\zeta_a}}
\right|\le\frac{C^{\nab,a}_{K,\delta}}{N}\,,
\end{array}
\label{eq:Tnaba}
\end{equation}
where
\begin{equation}
\zeta_a:=N\int_a^z\rho^0(s)\,ds\,.
\end{equation}
If instead it is $\Sigma_0^\nab$ that is bounded away from $K$, then there is a constant $C^{\del,a}_K>0$ and for each $\delta>0$ there is a constant
$C^{\del,a}_{K,\delta}>0$ such that for sufficiently large $N$,
\begin{equation}
\begin{array}{c}\displaystyle
\sup_{z\in K,|\arg(z-a)|<\pi}\left|T_\del(z)-\frac{\sqrt{2\pi}e^{-\zeta_a}\zeta_a^{\zeta_a}}{\Gamma(\zeta_a+1/2)}\right|\le\frac{C^{\del,a}_K}{N}\,,
\\\\\displaystyle
\sup_{z\in K,\delta\le|\arg(z-a)|\le\pi}\left|Y(z)^{-1}-
\frac{\Gamma(1/2-\zeta_a)}{\sqrt{2\pi}e^{\zeta_a}(-\zeta_a)^{-\zeta_a}}
\right|\le\frac{C^{\del,a}_{K,\delta}}{N}\,.
\end{array}
\label{eq:Tdela}
\end{equation}
Similarly,
if $K\subset G$ is a compact neighborhood of $z=b$ and $\Sigma_0^\del$ is bounded away from $K$, then there is a constant $C^{\nab,b}_K>0$ and for each $\delta>0$ there is a constant $C^{\nab,b}_{K,\delta}>0$ such that for
sufficiently large $N$,
\begin{equation}
\begin{array}{c}\displaystyle
\sup_{z\in K,|\arg(b-z)|<\pi}\left|T_\nab(z)-\frac{\sqrt{2\pi}e^{-\zeta_b}\zeta_b^{\zeta_b}}{\Gamma(\zeta_b+1/2)}\right|\le
\frac{C^{\nab,b}_K}{N}\,,
\\\\\displaystyle
\sup_{z\in K,\delta\le|\arg(b-z)|\le\pi}\left|Y(z)-
\frac{\Gamma(1/2-\zeta_b)}{\sqrt{2\pi}e^{\zeta_b}(-\zeta_b)^{-\zeta_b}}
\right|\le\frac{C^{\nab,b}_{K,\delta}}{N}\,,
\end{array}
\label{eq:Tnabb}
\end{equation}
where
\begin{equation}
\zeta_b:=N\int_z^b\rho^0(s)\,ds\,.
\end{equation}
If instead it is $\Sigma_0^\nab$ that is bounded away from $K$, then there
is a constant $C^{\del,b}_K>0$ and for each $\delta>0$ there is a constant $C^{\del,b}_{K,\delta}>0$ such that for sufficiently large $N$,
\begin{equation}
\begin{array}{c}\displaystyle
\sup_{z\in K,|\arg(b-z)|<\pi}\left|T_\del(z)-\frac{\sqrt{2\pi}e^{-\zeta_b}\zeta_b^{\zeta_b}}{\Gamma(\zeta_b+1/2)}\right|\le
\frac{C^{\del,b}_K}{N}\,,
\\\\\displaystyle
\sup_{z\in K,\delta\le|\arg(b-z)|\le\pi}\left|Y(z)^{-1}-
\frac{\Gamma(1/2-\zeta_b)}{\sqrt{2\pi}e^{\zeta_b}(-\zeta_b)^{-\zeta_b}}
\right|\le\frac{C^{\del,b}_{K,\delta}}{N}\,.
\end{array}
\label{eq:Tdelb}
\end{equation}
\end{enumerate}
\label{prop:tasymp}
\end{prop}

\begin{proof}
For $z,x\in G$, let \label{symbol:differencequotient}
\begin{equation}
D(z,x):=\frac{1}{z-x}\int_x^z\rho^0(s)\,ds\,.
\end{equation}
This function is analytic in both variables, and since $D(z,x)$ is strictly positive for $z$ and $x$ both in $[a,b]$ we may choose $G$ to be a sufficiently small neighborhood of $[a,b]$ to ensure that $\Re(D(z,x))$ is strictly positive for all $z$ and $x$ in $G$.  In particular, $D(z,x)$ is nonzero.  It follows that $\log(D(z,x))$ is well-defined as an analytic function for $z$ and $x$ in $G$.  Next, we define an analytic function for $x\in G$ by the integral
\label{symbol:m(x)}
\begin{equation}
m(x):=\int_a^x\rho^0(s)\,ds\,.
\end{equation}
since $\rho^0(s)$ is strictly positive in $[a,b]$, there is a unique analytic inverse function which we denote by $x(m)$ which is defined for $m\in m(G)$, where $m(G)$ is an open complex neighborhood of $[0,1]$.  It follows that
\begin{equation}
\frac{
\partial^2}{\partial m^2}\log(D(z,x(m)))\hspace{0.2 in}
\mbox{is uniformly bounded for $z\in G$ and $m\in m(G)$.}
\label{eq:fact1}
\end{equation}
Using this fact, we see that there are constants $C^\nab>0$ and $C^\del>0$ such that for sufficiently large $N$,
\begin{equation}
\sup_{z\in G}\left|\int_{m(\Sigma_0^\nab)}\log(D(z,x(s)))\,ds - \sum_{n\in\nab}
\log(D(z,x(s_{N,n})))\frac{1}{N}\right|\le \frac{C^\nab}{N^2}\,,
\label{eq:Dnabestimate}
\end{equation}
and
\begin{equation}
\sup_{z\in G}\left|\int_{m(\Sigma_0^\del)}\log(D(z,x(s)))\,ds - \sum_{n\in\del}
\log(D(z,x(s_{N,n})))\frac{1}{N}\right|\le \frac{C^\del}{N^2}\,,
\label{eq:Ddelestimate}
\end{equation}
where \label{symbol:straightened}
\begin{equation}
s_{N,n}:=\frac{2n+1}{2N}\,.
\end{equation}
Indeed, these estimates follow from (\ref{eq:fact1}) because the sums
are Riemann sum estimates of the corresponding integrals with the
midpoints $s_{N,n}$ of the subintervals $(n/N,(n+1)/N)$ chosen as
sample points.  The midpoint rule is second-order accurate if the
second derivative of the integrand is uniformly bounded.  The
constants $C^\nab$ and $C^\del$ depend on $\max |(\partial^2/\partial
m^2)\log(D(z,x(m)))|$ for $z\in G$ and $m\in m(G)$.  

For $z\in\Sigma_0^\nab$ we define \label{symbol:TildeTnab}
\begin{equation}
\begin{array}{rcl}
\tilde{T}_\nab(z)&:=&\displaystyle 2\cos(\pi N-\pi N m(z))
\frac{\displaystyle\prod_{n\in\del}
(m(z)-s_{N,n})}{\displaystyle\prod_{n\in\nab}(m(z)-s_{N,n})}\\\\
&&\,\,\,\times\,\,\,\displaystyle
\exp\left(N\left[\int_{m(\Sigma_0^\nab)}\log|m(z)-s|\,ds-
\int_{m(\Sigma_0^\del)}\log|m(z)-s|\,ds\right]\right)\,,
\end{array}
\end{equation}
which is extended by analytic continuation to
$z\in D_\nab$. The estimates
(\ref{eq:Dnabestimate}) and (\ref{eq:Ddelestimate}) imply that
uniformly for all $z\in D_\nab$,
\begin{equation}
T_\nab(z)=\tilde{T}_\nab(z)\left(1+O\left(\frac{1}{N}\right)\right)
\hspace{0.2 in}\text{as $N\rightarrow\infty$\,.}
\label{eq:TnabTtildenab}
\end{equation}
Indeed, some straightforward calculations show that
\begin{equation}
\begin{array}{rcl}
\displaystyle\log\left(\frac{T_\nab(z)}{\tilde{T}_\nab(z)}\right) &=&
\displaystyle\left(\sum_{n\in\nab}\log(D(z,x_{N,n}))-N\int_{m(\Sigma_0^\nab)}
\log(D(z,x(s)))\,ds\right)\\\\
&&\displaystyle\,\,\,-\,\,\,
\left(\sum_{n\in\del}\log(D(z,x_{N,n}))-N\int_{m(\Sigma_0^\del)}
\log(D(z,x(s)))\,ds\right)\,,
\end{array}
\end{equation}
from which \eqref{eq:TnabTtildenab} follows.
Similarly for $z\in\Sigma_0^\del$ we define \label{symbol:TildeTdel}
\begin{equation}
\begin{array}{rcl}
\tilde{T}_\del(z)&:=&\displaystyle 2\cos(\pi N-\pi N m(z))
\frac{\displaystyle\prod_{n\in\nab}
(m(z)-s_{N,n})}{\displaystyle\prod_{n\in\del}(m(z)-s_{N,n})}\\\\
&&\,\,\,\times\,\,\,\displaystyle
\exp\left(N\left[\int_{m(\Sigma_0^\del)}\log|m(z)-s|\,ds-
\int_{m(\Sigma_0^\nab)}\log|m(z)-s|\,ds\right]\right)
\end{array}
\end{equation}
which is extended to $z\in D_\del$ by analytic continuation, and we
have 
\begin{equation}
T_\del(z)=\tilde{T}_\del(z)\left(1+O\left(\frac{1}{N}\right)\right)
\hspace{0.2 in}\text{as $N\rightarrow\infty$\,,}
\label{eq:TdelTtildedel}
\end{equation}
holding uniformly for $z\in D_\del$.  The uniform asymptotic relations
\eqref{eq:TnabTtildenab} and \eqref{eq:TdelTtildedel} effectively reduce the 
asymptotic analysis of the functions $T_\nab(z)$ and $T_\del(z)$ to
that of the functions $\tilde{T}_\nab(z)$ and $\tilde{T}_\del(z)$.
This is advantageous because the discrete points $s_{N,n}$ are equally
spaced while the nodes $x_{N,n}$ are not necessarily so.

Thus it remains to study $\tilde{T}_\nab(z)$ and $\tilde{T}_\del(z)$.
We will consider $\tilde{T}_\nab(z)$, since the analysis of
$\tilde{T}_\del(z)$ is similar.  Assume that $K$ is a compact subset
of the open set $D_\nab$.  Let
\begin{equation}
\delta_K:=\frac{1}{2}\inf_{z\in K\cap\Sigma_0^\nab,w\in \Sigma_0^\del\cup\{a,b\}}|z-w|>0
\end{equation}
be half the minimum distance of $K\cap\Sigma_0^\nab$ from the boundary of
$\Sigma_0^\nab$.  Also, define the open covering $U$ by
\begin{equation}
U := \bigcup_{z\in K\cap\Sigma_0^\nab}(z-\delta_K,z+\delta_K)
\end{equation}
and let $F=\overline{U}$ be the closure.
Finally, set
\begin{equation}
\epsilon_K := \inf_{z\in K,\Re(z)\notin F} |\Im(z)|>0\,.
\end{equation}
This is strictly positive because $K$ is compact and can only touch the real axis in the interior of subintervals of $\Sigma_0^\nab$.
Thus, each $z\in K$ satisfies either $|\Im(z)|\ge\epsilon_K>0$ (because $\Re(z)\notin F$) or
\begin{equation}
\inf_{w\in\Sigma_0^\del\cup\{a,b\}}|w-\Re(z)|\ge \delta_K>0
\end{equation}
(because $\Re(z)\in F$).

We may extend $\tilde{T}_\nab(z)$ into the complex plane from $\Sigma_0^\nab$ by the following formula:
\begin{equation}
\begin{array}{rcl}
\tilde{T}_\nab(z)&=&\displaystyle \left(1+e^{2\pi i{\rm sgn}(\Im(z))N m(z)}\right)\\\\
&&\displaystyle
\,\,\,\times\,\,\,\exp\left(N\left[\int_{m(\Sigma_0^\nab)}\log(m(z)-s)\,ds-
\sum_{n\in\nab}\log(m(z)-s_{N,n})\frac{1}{N}\right]\right)\\\\
&&\displaystyle
\,\,\,\times\,\,\,\exp\left(N\left[\sum_{n\in\del}\log(m(z)-s_{N,n})\frac{1}{N}-\int_{m(\Sigma_0^\del)}\log(m(z)-s)\,ds\right]\right)\,.
\end{array}
\label{eq:tildetnabupper}
\end{equation}
Suppose that for some $\epsilon>0$, we have $|\Im(z)|\ge\epsilon$, a
condition that also bounds $\Im(m(z))$ away from zero. Therefore
$\log(m(z)-s)$ has a second derivative with respect to $s$ that is
uniformly bounded for all $s\in [a,b]$.  The bound on the second
derivative will depend on $\epsilon$ and the function $\rho^0(s)$ used
to define the function $m(z)$. In any case, an argument involving
midpoint-rule Riemann sums shows that the second and third lines of
(\ref{eq:tildetnabupper}) are each uniformly of the form $1+O(1/N)$ as
$N\rightarrow\infty$ for $|\Im(z)|\ge\epsilon$.  Furthermore, a
Cauchy-Riemann argument shows that the first line of
(\ref{eq:tildetnabupper}) is exponentially close to one as
$N\rightarrow\infty$ for $z\in G$ with $|\Im(z)|\ge\epsilon$.  Thus we
have shown that there is a constant $C_\epsilon>0$ such that for
sufficiently large $N$,
\begin{equation}
\sup_{z\in G, |\Im(z)|\ge\epsilon>0}|\tilde{T}_\nab(z)-1|\le\frac{C_\epsilon}{N}\,.
\end{equation}

Next suppose that $\Re(z)\in\Sigma_0^\nab$, bounded away from
$\Sigma_0^\del\cup\{a,b\}$ by a distance $\delta>0$.  Let $J$ denote
the maximal component interval of $\Sigma_0^\nab$ that contains
$\Re(z)$, and suppose that the corresponding index subset of $\nab$
consists of the contiguous list of integers $A,A+1,\dots,B-1,B$. Then,
from the representation (\ref{eq:tildetnabupper}), one sees once again
by a midpoint-rule Riemann sum argument that the factor on the third
line of (\ref{eq:tildetnabupper}) is of the form $1+O(1/N)$ as
$N\rightarrow\infty$ with a constant on the $O(1/N)$ term that depends
on $\delta$.  A similar argument applies to the factor on the second
line of (\ref{eq:tildetnabupper}) with the exception of the
contribution of the integral over $J$ and the corresponding discrete
sum.  Thus, uniformly for $\Re(z)$ as above, we have
\begin{equation}
\tilde{T}_\nab(z)=\tilde{T}^J_\nab(z)\left(1+O\left(\frac{1}{N}\right)\right)
\hspace{0.2 in}\text{as $N\rightarrow\infty$\,,}
\label{eq:TnabTJnab}
\end{equation}
where \label{symbol:TildeTnabJ}
\begin{equation}
\tilde{T}^J_\nab(z):=\frac{2N^{B+1-A}e^{i\pi {\rm sgn}(\Im(z))Nm(z)}\cos(\pi Nm(z))}
{\displaystyle\prod_{n=A}^{B}\left(Nm(z)-n-\frac{1}{2}\right)}
\exp\left(N\int_{\frac{A}{N}}^{\frac{B+1}{N}}\log(m(z)-s)\,ds\right)\,.
\label{eq:TJnab}
\end{equation}
Evaluating the integral exactly and rewriting the product in terms of the Euler gamma function, this becomes
\begin{equation}
\begin{array}{rcl}
\tilde{T}^J_\nab(z)&=&\displaystyle
2(-1)^{B+1}e^{-(B+1-A)}\frac{\displaystyle \Gamma\left(Nm(z)-B-\frac{1}{2}\right)}{\displaystyle\Gamma\left(Nm(z)-A+\frac{1}{2}\right)}\cos(\pi Nm(z))\\\\
&&\,\,\,\times\,\,\,\displaystyle e^{(Nm(z)-A)\log(Nm(z)-A)}
e^{(B+1-Nm(z))\log(B+1-Nm(z))}\,,
\end{array}
\end{equation}
and with the use of the reflection identity
$\Gamma(1/2+z)\Gamma(1/2-z)=\pi\sec(\pi z)$, we get
\begin{equation}
\tilde{T}^J_\nab(z)=\frac{2\pi e^{-(B+1-A)} e^{(Nm(z)-A)\log(Nm(z)-A)}
e^{(B+1-Nm(z))\log(B+1-Nm(z))}}
{\displaystyle \Gamma\left(Nm(z)-A+\frac{1}{2}\right)
\Gamma\left(B+1-Nm(z)+\frac{1}{2}\right)}\,.
\end{equation}
Now the condition that $\Re(z)$ be bounded away from the endpoints of
$J$ by at least $\delta>0$ fixed implies that $Nm(z)-A+1/2$ and
$B+1-Nm(z)+1/2$ are both quantities in the right half-plane that scale
like $N$; an application of Stirling's formula then gives, uniformly
for such $z$,
\begin{equation}
\tilde{T}_\nab^J(z)=1+O\left(\frac{1}{N}\right)\hspace{0.2 in}
\text{as $N\rightarrow\infty$\,.}
\end{equation}
Taking $\delta=\delta_K$ and $\epsilon=\epsilon_K$ then completes the
proof of (\ref{eq:Tnabout}), that $T_\nab(z)-1$ is uniformly of order
$1/N$ as $N\rightarrow\infty$ for $z\in K$, where $K$ is bounded away
from $(-\infty,a)\cup\Sigma_0^\del\cup(b,+\infty)$.  Analogous
arguments establish the corresponding result (\ref{eq:Tdelout}) for
$T_\del(z)$.  Using (\ref{eq:Tsqrtidentity_up}) and
(\ref{eq:Tsqrtidentity_down}) then proves (\ref{eq:Sout}).  Thus
statement 1 is established.

If $K$ is a compact set containing the left endpoint $z=a$ and bounded
away from $\Sigma_0^\del$ (so that the lower constraint is active at
the left endpoint), then one may follow nearly identical arguments to
arrive at the asymptotic relation \eqref{eq:TnabTJnab} now holding
uniformly for $z\in K$, where $J$ is the leftmost subinterval of
$\Sigma_0^\nab$ and $\tilde{T}^J_\nab(z)$ is defined by
\eqref{eq:TJnab}.  In this case we have $A=0$, so we only expand the
gamma function involving $B$.  This proves the first line of
(\ref{eq:Tnaba}); the second line follows upon using
(\ref{eq:Tsqrtidentity_up}) and (\ref{eq:Tsqrtidentity_down}).  On the
other hand, if $K$ contains $z=b$ where the lower constraint is
active, then again we have \eqref{eq:TnabTJnab} holding uniformly for
$z\in K$ where now $J$ is the rightmost subinterval of
$\Sigma_0^\nab$.  Thus, $B=N-1$ and we only expand the gamma function
involving $A$ to prove (\ref{eq:Tnabb}).  The analogous statements
(\ref{eq:Tdela}) and (\ref{eq:Tdelb}) are proved similarly.  Thus
statement 2 is established.
\end{proof}

\subsection{Steepest descent:  the transformation $\mat{S}(z)\to\mat{X}(z)$.}
Now, from any band interval $I\cap\Sigma_0^\nab$, the matrix
$\mat{L}_+(z)$ admits an analytic continuation into the upper
half-plane, and the matrix $\mat{L}_-(z)$ admits an analytic
continuation into the lower half-plane.  Since the function
$\theta(z)$ is real and increasing in $I\cap\Sigma_0^\nab$, its
analytic continuation from $I$, which we denote by $\theta^\nab_I(z)$,
\label{symbol:thetanabI}
will have a positive imaginary part near the real axis in the upper
half-plane, and a negative imaginary part near the real axis in the
lower half-plane, as a simple Cauchy-Riemann argument shows.  Thus,
the factors $e^{\pm iN\theta^\nab_I(z)}$ present in $\mat{L}_\pm(z)$
continued into their respective half-planes become exponentially small
as $N\rightarrow\infty$.  Subject to the claim that the analytic
function $T_\nab(z)-1$ remains uniformly small upon analytic
continuation, we see that the analytic continuation of $\mat{L}_+(z)$
and $\mat{L}_-(z)$ into the upper and lower half-planes respectively
become small perturbations of the identity matrix.

Similarly, from a band interval $I\cap\Sigma_0^\del$, the analytic
continuation of the matrix $\mat{U}_-(z)$ into the upper half-plane
and that of $\mat{U}_+(z)$ into the lower half-plane will be small
perturbations of the identity matrix in the limit
$N\rightarrow\infty$, because the real function $\theta(z)$ is
strictly decreasing.  This implies that the analytic continuation of
$\theta(z)$, which in this case we refer to as $\theta^\del_I(z)$, 
\label{symbol:thetadelI} has
an imaginary part that is positive in the lower half-plane and
negative in the upper half-plane.

Therefore, if the factors $\mat{U}_+(z)$ and $\mat{L}_+(z)$ can be
deformed into the upper half-plane, and at the same time if the
factors $\mat{U}_-(z)$ and $\mat{L}_-(z)$ can be deformed into the
lower half-plane, then the rapidly oscillatory jump matrix for
$\mat{S}(z)$ in the bands will be resolved into near-identity factors
and a central slowly-varying factor.  This idea is the essence of the
steepest descent method for matrix Riemann-Hilbert problems developed
by Deift and Zhou.

To carry out the deformation, it will be convenient to introduce some
explicit formulae for the analytic continuations $\theta^\nab_I(z)$
and $\theta^\del_I(z)$.  If $I\subset \Sigma_0^\nab$ is a band
containing a point (or endpoint) $x$, then we have
\begin{equation}
\theta^\nab_I(z):=\theta(x)+2\pi c\int_x^z\psi_I(s)\,ds\,.
\label{eq:thetanabI}
\end{equation}
If $I\subset\Sigma_0^\del$ is a band containing a point (or endpoint) $x$, then we have
\begin{equation}
\theta^\del_I(z):=\theta(x)-2\pi c\int_x^z\overline{\psi}_I(s)\,ds\,.
\label{eq:thetadelI}
\end{equation}
If $I$ is a transition band, then it is divided into two halves,
$I\cap\Sigma_0^\nab$ and $I\cap\Sigma_0^\del$, by the transition point
$y_{k,N}$ therein.  From $I\cap\Sigma_0^\nab$ we obtain a continuation
$\theta^\nab_I(z)$ of $\theta(z)$ using the formula
(\ref{eq:thetanabI}) for $x\in I\cap\Sigma_0^\nab$, and from
$I\cap\Sigma_0^\del$ we obtain a continuation $\theta^\del_I(z)$ of
$\theta(z)$ using the formula (\ref{eq:thetadelI}) for $x\in
I\cap\Sigma_0^\del$.

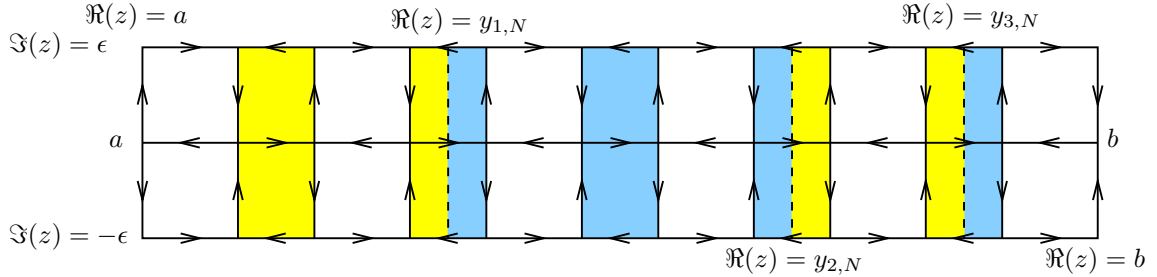
\begin{figure}[h]
\begin{center}
\input{SigmaSD.pstex_t}
\end{center}
\caption{\em The oriented contour $\Sigma_{\rm SD}$ consists of the 
interval $[a,b]$, corresponding horizontal segments $\Im(z)=\pm\epsilon$, and
vertical segments aligned at the edges of all band intervals.  The
dashed vertical lines separating the yellow and blue regions are not
part of $\Sigma_{\rm SD}$.}
\label{fig:SigmaSD}
\end{figure}

Based on the factorizations \eqref{eq:bandjumpnab} and
\eqref{eq:bandjumpdel}, we now carry out the steepest descent deformation, 
introducing a final change of variables defining a new unknown
$\mat{X}(z)$ \label{symbol:matrixX} 
in terms of $\mat{S}(z)$ with the aim of obtaining a jump
condition for $\mat{X}(z)$ in the bands involving only the matrix
$\mat{J}(z)$.  Let $\Sigma_{\rm SD}$ \label{symbol:SigmaSD}
be the oriented contour
illustrated in Figure~\ref{fig:SigmaSD}.
For each band interval $I\subset (a,b)$, we make the following definitions.  If $z$ lies in the
open rectangle $I\cap\Sigma_0^\nab + i(0,\epsilon)$ (these are the blue-shaded rectangles lying in the upper half-plane in Figure~\ref{fig:SigmaSD}) we set
\begin{equation}
\mat{X}(z):=\mat{S}(z)\left(\begin{array}{cc}
T_\nab(z)^{1/2} & 0 \\\\
-iT_\nab(z)^{-1/2}e^{\eta(z)-\gamma-2\kappa g(z)}e^{iN\theta^\nab_I(z)} &
T_\nab(z)^{-1/2}\end{array}\right)\,.
\label{eq:X-nab-upper-def}
\end{equation}
If $z$ lies in the open rectangle $I\cap\Sigma_0^\nab -i (0,\epsilon)$
(these are the blue-shaded rectangles in the lower half-plane) we set
\begin{equation}
\mat{X}(z):=\mat{S}(z)\left(\begin{array}{cc}
T_\nab(z)^{1/2} & 0 \\\\
iT_\nab(z)^{-1/2}e^{\eta(z)-\gamma-2\kappa g(z)}e^{-iN\theta^\nab_I(z)} &
T_\nab(z)^{-1/2}\end{array}\right)\,.
\end{equation}
Next, if $z$ lies in the open rectangle $I\cap\Sigma_0^\del+i(0,\epsilon)$
(these are the yellow-shaded rectangles in the upper half-plane) we set
\begin{equation}
\mat{X}(z):=\mat{S}(z)\left(\begin{array}{cc}
T_\del(z)^{-1/2} & -iT_\del(z)^{-1/2}e^{\gamma-\eta(z)+2\kappa g(z)}e^{-iN\theta^\del_I(z)} \\\\
0 & T_\del(z)^{1/2}\end{array}\right)\,.
\label{eq:X-del-lower-def}
\end{equation}
And if $z$ lies in the open rectangle $I\cap\Sigma_0^\del-i(0,\epsilon)$
(the yellow-shaded rectangles in the lower half-plane) we set
\begin{equation}
\mat{X}(z):=\mat{S}(z)\left(\begin{array}{cc}
T_\del(z)^{-1/2} & iT_\del(z)^{-1/2}e^{\gamma-\eta(z)+2\kappa g(z)}e^{iN\theta^\del_I(z)} \\\\
0 &
T_\del(z)^{1/2}\end{array}\right)\,.
\end{equation}
Finally for all remaining $z\in\mathbb{C}\setminus\Sigma_{\rm SD}$
we set
\begin{equation}\label{eq:XfromSlast}
  \mat{X}(z):=\mat{S}(z).
\end{equation}

\subsection{Properties of $\mat{X}(z)$.}
This change of variables is the last of a sequence of exact and
explicit transformations relating $\mat{P}(z;N,k)$ to
$\mat{Q}(z;N,k)$, $\mat{Q}(z;N,k)$ to $\mat{R}(z)$, $\mat{R}(z)$ to
$\mat{S}(z)$, and finally $\mat{S}(z)$ to $\mat{X}(z)$.  For future
reference it will be useful to summarize this sequence of
transformations by presenting the explicit formulae directly giving
$\mat{X}(z)$ in terms of $\mat{P}(z;N,k)$, the solution of Interpolation
Problem~\ref{rhp:DOP}.  In general, the transformation may be written
as
\begin{equation}
\mat{X}(z)=e^{(N\ell_c+\gamma)\sigma_3/2}
\mat{P}(z;N,k)\mat{D}(z)e^{(N(d_N-c)-\kappa)g(z)\sigma_3}
e^{-(N\ell_c+\gamma)\sigma_3/2}\,,
\label{eq:XD}
\end{equation}
where the matrix $\mat{D}(z)$ \label{symbol:matrixD} 
takes different forms in different
regions of the complex plane as follows.  For $z$ in the unbounded
component of $\mathbb{C}\setminus\Sigma_{\rm SD}$, we have
\begin{equation}
\mat{D}(z):=\left(\begin{array}{cc}
\displaystyle \prod_{n\in\del}(z-x_{N,n})^{-1} & 0 \\\\ 
0 & \displaystyle \prod_{n\in\del}(z-x_{N,n})
\end{array}\right)\,.
\label{eq:Dfirst}
\end{equation}
For $z$ in the regions $\Omega_\pm^\nab$ such that $\Re(z)$ lies in a
void of $[a,b]$, we have
\begin{equation}
\mat{D}(z):=\left(\begin{array}{cc}
\displaystyle \prod_{n\in\del}(z-x_{N,n})^{-1} &
\displaystyle \mp ie^{\mp iN\theta^0(z)/2}e^{-NV_N(z)}
\prod_{n\in\nab}(z-x_{N,n})^{-1}\\\\
0 & \displaystyle\prod_{n\in\del}(z-x_{N,n})\end{array}\right)\,.
\label{eq:Dnexttovoid}
\end{equation}
For $z$ in the regions $\Omega_\pm^\del$ such that $\Re(z)$ lies in
a saturated region of $[a,b]$, we have
\begin{equation}
\mat{D}(z):=\left(\begin{array}{cc}
\displaystyle\prod_{n\in\del}(z-x_{N,n})^{-1} & 0 \\\\
\displaystyle \mp ie^{\mp iN\theta^0(z)/2}e^{NV_N(z)}
\prod_{n\in\nab}(z-x_{N,n}) &
\displaystyle\prod_{n\in\del}(z-x_{N,n})
\end{array}\right)\,.
\label{eq:Dnexttosat}
\end{equation} 
For $z$ in the regions $\Omega_+^\nab$ such that
$\Re(z)$ lies in a band $I$ of $[a,b]$ (the blue regions in the upper
half-plane in Figure~\ref{fig:SigmaSD}), we have
\begin{equation}
\begin{array}{rcl}
D_{11}(z)&:=&
\displaystyle T_\nab(z)^{1/2}\prod_{n\in\del}(z-x_{N,n})^{-1} \\\\
&&\displaystyle \,\,\,-\,\,\, T_\nab(z)^{-1/2}
e^{N(\ell_c-2(d_N-c)g(z)-V(z)+i\theta_I^\nab(z)-i\theta^0(z)/2)}
\prod_{n\in\nab}(z-x_{N,n})^{-1}\,, \\\\
D_{12}(z)&:=&
\displaystyle -iT_\nab(z)^{-1/2}e^{-\eta(z)}e^{-N(V(z)+i\theta^0(z)/2)}
\prod_{n\in\nab}(z-x_{N,n})^{-1}\,, \\\\
D_{21}(z)&:=&
\displaystyle -iT_\nab(z)^{-1/2}e^{\eta(z)}
e^{N(\ell_c-2(d_N-c)g(z)+i\theta_I^\nab(z))}\prod_{n\in\del}(z-x_{N,n})
\,,
\\\\
D_{22}(z)&:=&
\displaystyle T_\nab(z)^{-1/2}\prod_{n\in\del}(z-x_{N,n})\,.
\end{array}
\label{eq:Xblueupper}
\end{equation}
For $z$ in the regions $\Omega_-^\nab$ such that $\Re(z)$ lies in a
band $I$ of $[a,b]$ (the blue regions in the lower half-plane in
Figure~\ref{fig:SigmaSD}), we have
\begin{equation}
\begin{array}{rcl}
D_{11}(z)&:=&\displaystyle
T_\nab(z)^{1/2}\prod_{n\in\del}(z-x_{N,n})^{-1}\\\\
&&\displaystyle\,\,\,-\,\,\, T_\nab(z)^{-1/2}
e^{N(\ell_c-2(d_N-c)g(z)-V(z)-i\theta_I^\nab(z)+i\theta^0(z)/2)}
\prod_{n\in\nab}(z-x_{N,n})^{-1}\,,\\\\
D_{12}(z)&:=&\displaystyle
iT_\nab(z)^{-1/2}e^{-\eta(z)}e^{-N(V(z)-i\theta^0(z)/2)}\prod_{n\in\nab}(z-x_{N,n})^{-1}\,,\\\\
D_{21}(z)&:=&\displaystyle iT_\nab(z)^{-1/2}e^{\eta(z)}
e^{N(\ell_c-2(d_N-c)g(z)-i\theta_I^\nab(z))}\prod_{n\in\del}(z-x_{N,n})\,,\\\\
D_{22}(z)&:=& \displaystyle
T_\nab(z)^{-1/2}\prod_{n\in\del}(z-x_{N,n})\,.
\end{array}
\label{eq:Xbluelower}
\end{equation}
For $z$ in the regions $\Omega_+^\del$ such that $\Re(z)$ lies in a
band $I$ of $[a,b]$ (the yellow regions in the upper half-plane in
Figure~\ref{fig:SigmaSD}), we have
\begin{equation}
\begin{array}{rcl}
D_{11}(z)&:=&\displaystyle
T_\del(z)^{-1/2}\prod_{n\in\del}(z-x_{N,n})^{-1}\,,
\\\\
D_{12}(z)&:=&\displaystyle -iT_\del(z)^{-1/2}e^{-\eta(z)}
e^{-N(\ell_c-2(d_N-c)g(z)+i\theta_I^\del(z))}
\prod_{n\in\del}(z-x_{N,n})^{-1}\,,
\\\\
D_{21}(z)&:=&\displaystyle -iT_\del(z)^{-1/2}e^{\eta(z)}
e^{N(V(z)-i\theta^0(z)/2)}\prod_{n\in\nab}(z-x_{N,n})\,,\\\\
D_{22}(z)&:=&\displaystyle
T_\del(z)^{1/2}\prod_{n\in\del}(z-x_{N,n})\\\\
&&\displaystyle\,\,\,-\,\,\, T_\del(z)^{-1/2}
e^{-N(\ell_c-2(d_N-c)g(z)-V(z)+i\theta_I^\del(z)+i\theta^0(z)/2)}
\prod_{n\in\nab}(z-x_{N,n})\,.
\end{array}
\label{eq:Xyellowupper}
\end{equation}
Finally, for $z$ in the regions $\Omega_-^\del$ such that $\Re(z)$
lies in a band $I$ of $[a,b]$ (the yellow regions in the lower
half-plane in Figure~\ref{fig:SigmaSD}), we have
\begin{equation}
\begin{array}{rcl}
D_{11}(z)&:=&\displaystyle
T_\del(z)^{-1/2}\prod_{n\in\del}(z-x_{N,n})^{-1}\,,
\\\\
D_{12}(z)&:=&\displaystyle iT_\del(z)^{-1/2}e^{-\eta(z)}
e^{-N(\ell_c-2(d_N-c)g(z)-i\theta_I^\del(z))}
\prod_{n\in\del}(z-x_{N,n})^{-1}\,,
\\\\
D_{21}(z)&:=&\displaystyle iT_\del(z)^{-1/2}e^{\eta(z)}
e^{N(V(z)+i\theta^0(z)/2)}\prod_{n\in\nab}(z-x_{N,n})\,,\\\\
D_{22}(z)&:=&\displaystyle
T_\del(z)^{1/2}\prod_{n\in\del}(z-x_{N,n})\\\\
&&\displaystyle\,\,\,-\,\,\, T_\del(z)^{-1/2}
e^{-N(\ell_c-2(d_N-c)g(z)-V(z)-i\theta_I^\del(z)-i\theta^0(z)/2)}
\prod_{n\in\nab}(z-x_{N,n})\,.
\end{array}
\label{eq:Xyellowlower}
\end{equation}

Unlike the contour $\Sigma$, the new contour $\Sigma_{\rm SD}$ does
not contain the vertical segments $Y_N\pm i(0,\epsilon)$ that form the
common boundary of the yellow and blue rectangles and that are
illustrated with dashed lines in Figure~\ref{fig:SigmaSD}.  Since the
matrix $\mat{X}(z)$ is defined by different formulae in the yellow and
blue regions, one should suspect that $\mat{X}(z)$ cannot be defined
on the common boundary so as to make $\mat{X}(z)$ continuous there.
In other words, it would seem that there should be a jump
discontinuity of $\mat{X}(z)$ on these vertical segments.  On the
contrary, we have the following result.
\begin{prop}
The matrix $\mat{X}(z)$ defined from
(\ref{eq:X-nab-upper-def})--(\ref{eq:XfromSlast}) extends to a
function analytic in $\mathbb{C}\setminus\Sigma_{\rm SD}$.  In
particular, $\mat{X}(z)$ is continuous and analytic on the vertical
segments $Y_N\pm i(0,\epsilon)$.  Moreover, on each subset of
$\Sigma_{\rm SD}$ that contains no self-intersection points, the
matrix-valued ratio of boundary values taken by $\mat{X}(z)$ is an
analytic function of $z$.
\label{prop:disappear}
\end{prop}
\begin{proof}
Let $\mat{X}^{\nab,+}(z)$ denote the matrix $\mat{X}(z)$ defined by
\eqref{eq:XD} with $\mat{D}(z)$ given by \eqref{eq:Xblueupper},
and let $\mat{X}^{\del,+}(z)$ denote the matrix $\mat{X}(z)$ defined
by
\eqref{eq:XD} with $\mat{D}(z)$ given by \eqref{eq:Xyellowupper}.
We will show that $\mat{X}^{\nab,+}(z)$ and $\mat{X}^{\del,+}(z)$ are
the same analytic function in the common region $0<\Im(z)<\epsilon$
and $\Re(z)\in I$, where $I$ is a transition band.  By direct
calculation, we obtain
\begin{equation}
e^{(N(d_N-c)-\kappa)g(z)\sigma_3}e^{-(N\ell_c+\gamma)\sigma_3/2}
\mat{X}^{\nab,+}(z)^{-1}\mat{X}^{\del,+}(z)e^{(N\ell_c+\gamma)\sigma_3/2}
e^{-(N(d_N-c)-\kappa)g(z)\sigma_3} =\mat{A}^+(z)\,,
\end{equation}
where
\begin{equation}
\begin{array}{rcl}
A^+_{11}(z)&:=&\displaystyle
\frac{1+e^{-iN\theta^0(z)}}{T_\nab(z)^{1/2}T_\del(z)^{1/2}}\,,\\\\
A^+_{12}(z)&:=&\displaystyle
ie^{-\eta(z)}e^{-N(\ell_c-2(d_N-c)g(z)+i\theta_I^\del(z))}
\left[F^+_\del(z)^{-1}-
\frac{1+e^{-iN\theta^0(z)}}{T_\nab(z)^{1/2}T_\del(z)^{1/2}}\right]\,,\\\\
A^+_{21}(z)&:=&\displaystyle
ie^{\eta(z)}e^{N(\ell_c-2(d_N-c)g(z)+i\theta_I^\nab(z))}
\left[\frac{1+e^{-iN\theta^0(z)}}{T_\nab(z)^{1/2}T_\del(z)^{1/2}}-
F^+_\nab(z)^{-1}\right]\,,\\\\ A^+_{22}(z)&:=&\displaystyle
\frac{1+e^{-iN\theta^0(z)}}{T_\nab(z)^{1/2}T_\del(z)^{1/2}}
e^{iN(\theta_I^\nab(z)-\theta_I^\del(z))} +
T_\nab(z)^{1/2}T_\del(z)^{1/2} -e^{-iN\theta^0(z)}\left[F^+_\nab(z)+
F^+_\del(z)\right]\,,
\end{array}
\label{eq:A+}
\end{equation}
and where
\begin{equation}
\begin{array}{rcl}
F^+_\nab(z)&:=&\displaystyle
\frac{\displaystyle T_\del(z)^{1/2}\prod_{n\in\del}(z-x_{N,n})}
{\displaystyle T_\nab(z)^{1/2}\prod_{n\in\nab}(z-x_{N,n})}
e^{N(\ell_c-2(d_N-c)g(z)+i\theta_I^\nab(z)-V(z)+i\theta^0(z)/2)}\,,
\\\\
F^+_\del(z)&:=&\displaystyle
\frac{\displaystyle T_\nab(z)^{1/2}\prod_{n\in\nab}(z-x_{N,n})}
{\displaystyle T_\del(z)^{1/2}\prod_{n\in\del}(z-x_{N,n})}
e^{-N(\ell_c-2(d_N-c)g(z)+i\theta_I^\del(z)-V(z)-i\theta^0(z)/2)}\,.
\end{array}
\label{eq:Fplus}
\end{equation}
Now, taking the base points $x$ in the formulae \eqref{eq:thetanabI}
and \eqref{eq:thetadelI} to both coincide with the transition point
$y_{k,N}$ in the transition band $I$, then recalling
\eqref{eq:psiIpsiIbar} and using the quantization condition
\eqref{eq:yquantize}, and finally comparing with the definition 
\eqref{eq:theta0def} of $\theta^0(z)$, we obtain the identity
\begin{equation}
e^{iN(\theta_I^\nab(z)-\theta_I^\del(z))} = e^{-iN\theta^0(z)}
\,,\hspace{0.2 in}
\text{for $|\Im(z)|<\epsilon$, $\Re(z)\in I$, and $N\in\mathbb{Z}$\,.}
\label{eq:thetaInabdel}
\end{equation}
Taking this identity into account, along with the identity
\eqref{eq:Tsqrtidentity_up} from Proposition~\ref{prop:tanal} 
valid for $\Im(z)>0$, we therefore see that the matrix elements
\eqref{eq:A+} simplify:
\begin{equation}
\begin{array}{rcl}
A_{11}^+(z)&=&1\,,\\\\ A_{12}^+(z)&=&\displaystyle
ie^{-\eta(z)}e^{-N(\ell_c-2(d_N-c)g(z)+i\theta_I^\del(z))}
\left[F^+_\del(z)^{-1}-1\right]\,,\\\\
A_{21}^+(z)&=&\displaystyle
ie^{\eta(z)}e^{N(\ell_c-2(d_N-c)g(z)+i\theta_I^\nab(z))}
\left[1-F_\nab^+(z)^{-1}\right]\,,\\\\
A_{22}^+(z)&=&\displaystyle
1+e^{-iN\theta^0(z)}\left[2-F_\nab^+(z)-F_\del^+(z)\right]\,.
\end{array}
\end{equation}
Thus, to show that $\mat{X}^{\nab,+}(z)\equiv\mat{X}^{\del,+}(z)$, it
suffices to show that $F^+_\nab(z)\equiv 1$ and $F^+_\del(z)\equiv 1$.

Let us calculate the boundary value taken by the function
$F_\nab^+(z)$ on the real interval $I\cap\Sigma_0^\nab$ from the upper
half-plane.  For such $z$, we have three facts at our disposal, namely
the identity $\theta_I^\nab(z)\equiv\theta(z)$, the identity
\eqref{eq:Tsqrtidentity_up} from Proposition~\ref{prop:tanal}, and the
formula
\eqref{eq:TAfuncdef}.  Applying these, and in particular first using 
the latter to eliminate the ratio of products in the definition
\eqref{eq:Fplus} of $F_\nab^+(z)$, we obtain simply
\begin{equation}
F_\nab^+(z)=\exp\left(N\left[\int_{\Sigma_0^\del}\log|z-x|\rho^0(x)\,dx-
\int_{\Sigma_0^\nab}\log|z-x|\rho^0(x)\,dx+
\ell_c-2(d_N-c)g_+(z)+i\theta(z)-V(z)\right]
\right)\,,
\end{equation}
where $g_+(z)$ indicates a boundary value taken from the upper
half-plane.  Using \eqref{eq:rhodef}, \eqref{eq:gdef}, and
\eqref{eq:theta}, we see that for real $z\in [a,b]$,
\begin{equation}
2(d_N-c)g_+(z)-i\theta(z)=-2c\int_a^b\log|z-x|\,d\mu_{\rm min}^c(x) +
2\int_{\Sigma_0^\del}\log|z-x|\rho^0(x)\,dx\,.
\end{equation}
Therefore, recalling
\eqref{eq:fielddef} and
\eqref{eq:variationalderivative}, we have simply
\begin{equation}
F_\nab^+(z)=\exp\left(N\left[\ell_c-\frac{\delta
E_c}{\delta\mu}(z)\right]
\right)\,,
\end{equation}
where the variational derivative is evaluated on the equilibrium
measure $\mu_{\rm min}^c$.  It follows that $F_\nab^+(z)\equiv 1$ as a
consequence of \eqref{eq:equilibrium}, since $z$ is in a band $I$.  By
analytic continuation this identity holds in the whole region
$0<\Im(z)<\epsilon$ with $\Re(z)\in I$.

We may also compute a boundary value of the function $F_\del^+(z)$,
letting $z$ tend toward the real interval $I\cap\Sigma_0^\del$ from
the upper half-plane.  In this case, instead of
\eqref{eq:TAfuncdef}, we use the identity
\eqref{eq:TBfuncdef} to eliminate the ratio of products, and we may write
$\theta_I^\del(z)\equiv \theta(z)$. The rest of the argument is
exactly the same, and we thus deduce that the identity
$F_\del^+(z)\equiv 1$ holds for $z\in I\cap \Sigma_0^\del$ in the
sense of a boundary value taken from the upper half-plane.  But by
analytic continuation it also holds in the whole region of interest:
$0<\Im(z)<\epsilon$ and $\Re(z)\in I$.  This completes the proof that
$\mat{X}(z)$ has no jump discontinuity along the vertical segments
between the blue and yellow regions illustrated in the upper
half-plane in Figure~\ref{fig:SigmaSD}.

Now let $\mat{X}^{\nab,-}(z)$ denote the matrix $\mat{X}(z)$ defined
by \eqref{eq:XD} with $\mat{D}(z)$ given by \eqref{eq:Xbluelower}, and
let $\mat{X}^{\del,-}(z)$ denote the matrix $\mat{X}(z)$ defined by
\eqref{eq:XD} with $\mat{D}(z)$ given by \eqref{eq:Xyellowlower}.
We will now show that $\mat{X}^{\nab,-}(z)$ and $\mat{X}^{\del,-}(z)$
are the same analytic function in the common region
$-\epsilon<\Im(z)<0$ and $\Re(z)\in I$, where $I$ is a transition
band.  As before, by direct calculation we have
\begin{equation}
e^{(N(d_N-c)-\kappa)g(z)\sigma_3}e^{-(N\ell_c+\gamma)\sigma_3/2}
\mat{X}^{\nab,-}(z)^{-1}\mat{X}^{\del,-}(z)e^{(N\ell_c+\gamma)\sigma_3/2}
e^{-(N(d_N-c)-\kappa)g(z)\sigma_3} =\mat{A}^-(z)\,,
\end{equation}
where
\begin{equation}
\begin{array}{rcl}
A_{11}^-(z)&:=&\displaystyle\frac{1+e^{iN\theta^0(z)}}{T_\nab(z)^{1/2}
T_\del(z)^{1/2}}\,,\\\\ A_{12}^-(z)&:=&\displaystyle
ie^{-\eta(z)}e^{-N(\ell_c-2(d_N-c)g(z)-i\theta_I^\del(z))}
\left[\frac{1+e^{iN\theta^0(z)}}{T_\nab(z)^{1/2}T_\del(z)^{1/2}}-
F^-_\del(z)^{-1}\right]\,,\\\\ A_{21}^-(z)&:=&\displaystyle
ie^{\eta(z)}e^{N(\ell_c-2(d_N-c)g(z)-i\theta_I^\nab(z))}
\left[F_\nab^-(z)^{-1}-\frac{1+e^{iN\theta^0(z)}}
{T_\nab(z)^{1/2}T_\del(z)^{1/2}}\right]\,,\\\\
A_{22}^-(z)&:=&\displaystyle
\frac{1+e^{iN\theta^0(z)}}{T_\nab(z)^{1/2}T_\del(z)^{1/2}}
e^{iN(\theta_I^\del(z)-\theta_I^\nab(z))} +
T_\nab(z)^{1/2}T_\del(z)^{1/2} -e^{iN\theta(z)}\left[F^-_\nab(z) +
F^-_\del(z)\right]\,,
\end{array}
\label{eq:A-}
\end{equation}
and where
\begin{equation}
\begin{array}{rcl}
F^-_\nab(z)&:=&\displaystyle
\frac{\displaystyle T_\del(z)^{1/2}\prod_{n\in\del}(z-x_{N,n})}
{\displaystyle T_\nab(z)^{1/2}\prod_{n\in\nab}(z-x_{N,n})}
e^{N(\ell_c-2(d_N-c)g(z)-i\theta_I^\nab(z)-V(z)-i\theta^0(z)/2)}\,,\\\\
F^-_\del(z)&:=&\displaystyle
\frac{\displaystyle T_\nab(z)^{1/2}\prod_{n\in\nab}(z-x_{N,n})}
{\displaystyle T_\del(z)^{1/2}\prod_{n\in\del}(z-x_{N,n})}
e^{-N(\ell_c-2(d_N-c)g(z)-i\theta_I^\del(z)-V(z)+i\theta^0(z)/2)}\,.
\end{array}
\label{eq:F-}
\end{equation}
Taking from Proposition~\ref{prop:tanal} the identity
\eqref{eq:Tsqrtidentity_down}, valid for $\Im(z)<0$, and using the identity
\eqref{eq:thetaInabdel}, these formulae simplify:
\begin{equation}
\begin{array}{rcl}
A_{11}^-(z)&=& 1\\\\ A_{12}^-(z)&=&\displaystyle
ie^{-\eta(z)}e^{-N(\ell_c-2(d_N-c)g(z)-i\theta_I^\del(z))}
\left[1-F^-_\del(z)^{-1}\right]\,,\\\\
A_{21}^-(z)&=&\displaystyle
ie^{\eta(z)}e^{N(\ell_c-2(d_N-c)g(z)-i\theta_I^\nab(z))}
\left[F^-_\nab(z)^{-1}-1\right]\,,\\\\
A_{22}^-(z)&=&\displaystyle
1+e^{iN\theta^0(z)}\left[2-F_\nab^-(z)-F_\del^-(z)\right]\,.
\end{array}
\end{equation}
Therefore again the problem reduces to showing that $F^-_\nab(z)\equiv
1$ and $F^-_\del(z)\equiv 1$.

Taking the boundary value of the function $F_\nab^-(z)$ from the lower
half-plane on the real interval $I\cap\Sigma_0^\nab$, we may
substitute for the ratio of products from
\eqref{eq:TAfuncdef} and use the identity
\eqref{eq:Tsqrtidentity_down} from Proposition~\ref{prop:tanal} along
with $\theta_I^\nab(z)\equiv \theta(z)$.  Since
\eqref{eq:rhodef}, \eqref{eq:gdef}, and \eqref{eq:theta} imply that
that for all real $z\in[a,b]$,
\begin{equation}
2(d_N-c)g_-(z)+i\theta(z) = -2c\int_a^b\log|z-x|\,d\mu_{\rm min}^c(x)+
2\int_{\Sigma_0^\del}\log|z-x|\rho^0(x)\,dx\,.
\end{equation}
where $g_-(z)$ indicates a boundary value taken from the lower
half-plane, the definitions \eqref{eq:fielddef} and
\eqref{eq:variationalderivative} along with the equilibrium condition
\eqref{eq:equilibrium} show that $F_\nab^-(z)\equiv 1$ in the sense of
a boundary value taken from the lower half-plane on
$I\cap\Sigma_0^\nab$.  But by analytic continuation, this identity
also holds throughout the region $-\epsilon<\Im(z)<0$ and $\Re(z)\in
I$.

To show that $F_\del^-(z)\equiv 1$ in the region $-\epsilon<\Im(z)<0$
and $\Re(z)\in I$, we repeat the above arguments but take the boundary
value from the lower-half plane in the interval $I\cap \Sigma_0^\del$,
where the identity \eqref{eq:TBfuncdef} may be used to eliminate the
ratio of products and where the identity $\theta_I^\del(z)\equiv
\theta(z)$ holds.  This completes the proof that $\mat{X}(z)$ has no
jump discontinuity along the vertical segments between the yellow and
blue regions illustrated in the lower half-plane in
Figure~\ref{fig:SigmaSD}.
\end{proof}

\begin{remark}
Part of the significance of Proposition~\ref{prop:disappear} is that
all essential dependence on the set $Y_\infty$, the choice of which
was somewhat arbitrary, has disappeared.  In particular, when we
approximate $\mat{X}(z)$ in the limit of large $N$, we will be able to
obtain error estimates that are of the same magnitude regardless of
the number of transition points, or indeed regardless of whether there
are any transition points at all.  This is an improvement over the
bounds stated in our announcement
\cite{BaikKMM03} which identified different estimates in two cases
(there called Case I and Case II) depending on whether any transition
points are present.
\end{remark}

\begin{figure}[h]
\begin{center}
\input{SigmaSD_components.pstex_t}
\end{center}
\caption{\em Components of the oriented contour $\Sigma_{\rm SD}$.}
\label{fig:SigmaSD_components}
\end{figure}
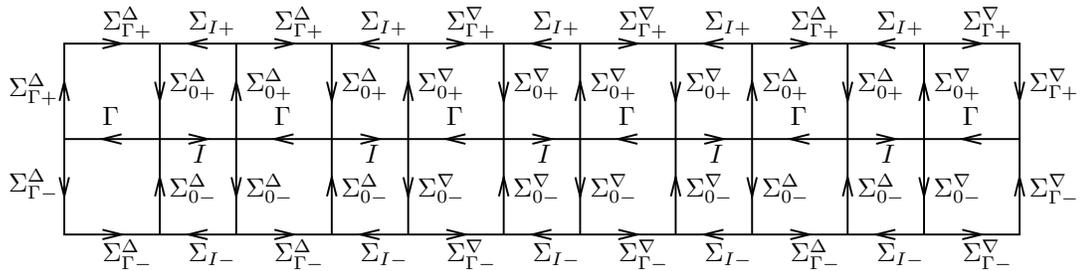

Having defined the matrix $\mat{X}(z)$ explicitly in terms of the
solution $\mat{P}(z;N,k)$ of Interpolation Problem~\ref{rhp:DOP} by
the formula \eqref{eq:XD} with $\mat{D}(z)$ given by
\eqref{eq:Dfirst}--\eqref{eq:Xyellowlower} 
allows us to replace that problem with an equivalent problem for the
new unknown $\mat{X}(z)$.  This is advantageous because the problem
whose solution is $\mat{X}(z)$ is more amenable to analysis.  In order
to correctly pose the problem, we must introduce some additional
notation for particular segments of $\Sigma_{\rm SD}$.  Vertical
segments of $\Sigma_{\rm SD}$ that are connected to band endpoints
will be denoted by $\Sigma_{0\pm}^\nab$ or $\Sigma_{0\pm}^\del$
\label{symbol:Sigma0pmnabdel}
depending on whether the endpoint lies in $\Sigma_0^\nab$ or
$\Sigma_0^\del$; the additional subscript indicates whether the
segment lies in the upper ($+$) or lower ($-$) half-plane.  Horizontal
segments lying above (below) bands will be denoted by $\Sigma_{I+}$
($\Sigma_{I-}$)\label{symbol:SigmaIpm}.  Horizontal segments lying above (below) voids will
be denoted by $\Sigma_{\Gamma+}^\nab$ 
($\Sigma_{\Gamma-}^\nab$)\label{symbol:SigmaGammapm}.
Horizontal segments lying above (below) saturated regions will be
denoted by $\Sigma_{\Gamma+}^\del$ ($\Sigma_{\Gamma-}^\del$).
Finally, each vertical segment passing through an endpoints $a$ or $b$
will be denoted by the same symbol as the component of $\Sigma_{\rm
SD}$ to which it is joined at $|\Im(z)|=\epsilon$.  See
Figure~\ref{fig:SigmaSD_components}.

The problem equivalent to Interpolation Problem~\ref{rhp:DOP} is the
subject of the following proposition.
\begin{prop}
The matrix $\mat{X}(z)$ defined by \eqref{eq:XD} and
\eqref{eq:Dfirst}--\eqref{eq:Xyellowlower} is the unique solution of
the following Riemann-Hilbert problem.
\end{prop}

\begin{rhp}
Find a $2\times 2$ matrix $\mat{X}(z)$ with the following properties:
\begin{enumerate}
\item
{\bf Analyticity}: $\mat{X}(z)$ is an analytic function of $z$ for
$z\in\cx\setminus \Sigma_{\rm SD}$.
\item
{\bf Normalization}: As $z\rightarrow\infty$,
\begin{equation}
\mat{X}(z)={\mathbb I} + O\left(\frac{1}{z}\right)\,.
\end{equation}
\item
{\bf Jump Conditions}: $\mat{X}(z)$ takes uniformly continuous
boundary values on $\Sigma_{\rm SD}$ from each connected component of
$\cx\setminus\Sigma_{\rm SD}$.  For each non-self-intersection point
$z\in\Sigma_{\rm SD}$ we denote by $\mat{X}_+(z)$ ($\mat{X}_-(z)$) the
limit of $\mat{X}(w)$ as $w\rightarrow z$ from the left (right).
Letting $g_+(z)+g_-(z)$ for real $z$ denote the sum of boundary values
taken by $g(z)$ from the upper and lower half-planes, the boundary
values taken on $\Sigma_{\rm SD}$ by $\mat{X}(z)$ satisfy the
following conditions.  For $z$ in a void $\Gamma\subset\Sigma_0^\nab$,
\begin{equation}
\mat{X}_+(z)=
\mat{X}_-(z)\left(\begin{array}{cc}
e^{-iN\theta_\Gamma}e^{-i\phi_\Gamma} &
iT_\nab(z)e^{\gamma-\eta(z)+\kappa(g_+(z)+g_-(z))}e^{-N\xi_\Gamma(z)}\\\\
0 & e^{iN\theta_\Gamma}e^{i\phi_\Gamma}
\end{array}\right)\,.
\end{equation}
For $z$ in a saturated region $\Gamma\subset\Sigma_0^\del$,
\begin{equation}
\mat{X}_+(z)=
\mat{X}_-(z)\left(\begin{array}{cc}
e^{-iN\theta_\Gamma}e^{-i\phi_\Gamma} & 0 \\\\
iT_\del(z)e^{\eta(z)-\gamma-\kappa(g_+(z)+g_-(z))}e^{-N\xi_\Gamma(z)}
& e^{iN\theta_\Gamma}e^{i\phi_\Gamma}
\end{array}
\right)\,.
\end{equation}
For $z$ in any band $I$,
\begin{equation}
\mat{X}_+(z)=
\mat{X}_-(z)\left(\begin{array}{cc}
0 & -ie^{\gamma-\eta(z)+\kappa(g_+(z)+g_-(z))} \\\\
-ie^{\eta(z)-\gamma-\kappa(g_+(z)+g_-(z))} & 0
\end{array}\right)\,.
\end{equation}
For $z$ in any vertical segment $\Sigma_{0\pm}^\nab$ meeting the real
axis at an endpoint $z_0$ of a band $I$,
\begin{equation}
\mat{X}_+(z)=
\mat{X}_-(z)\left(\begin{array}{cc}
T_\nab(z)^{\pm 1/2} & 0 \\\\\displaystyle
-iT_\nab(z)^{-1/2}e^{\eta(z)-\gamma-2\kappa g(z)}e^{\pm iN\theta(z_0)}
\exp\left(\pm 2\pi iNc\int_{z_0}^z\psi_I(s)\,ds\right) &
T_\nab(z)^{\mp 1/2}
\end{array}\right)\,.
\end{equation}
For $z$ in any vertical segment $\Sigma_{0\pm}^\del$ meeting the real
axis at an endpoint $z_0$ of a band $I$,
\begin{equation}
\mat{X}_+(z)=\mat{X}_-(z)
\left(\begin{array}{cc}
T_\del(z)^{\mp 1/2} & \displaystyle
-iT_\del(z)^{-1/2}e^{\gamma-\eta(z)+2\kappa g(z)} e^{\mp
iN\theta(z_0)}\exp\left(\pm 2\pi
iNc\int_{z_0}^z\overline{\psi}_I(z)\,ds
\right) \\\\0 & T_\del(z)^{\pm 1/2}
\end{array}\right)\,.
\end{equation}
For $z$ in any segment $\Sigma_{\Gamma\pm}^\nab$ parallel to a void
$\Gamma\subset\Sigma_0^\nab$ or with $\Re(z)=a$ or $\Re(z)=b$,
\begin{equation}
\mat{X}_+(z)=\mat{X}_-(z)
\left(\begin{array}{cc}
1 & iY(z)e^{\gamma-\eta(z)+2\kappa g(z)} e^{\mp iN\theta_\Gamma}e^{\mp
iN\theta^0(z)}e^{-N\xi_\Gamma(z)}
\\\\
0 & 1
\end{array}\right)\,.
\end{equation}
For $z$ in any segment $\Sigma_{\Gamma\pm}^\del$ parallel to a
saturated region $\Gamma\subset\Sigma_0^\del$ or with $\Re(z)=a$ or
$\Re(z)=b$,
\begin{equation}
\mat{X}_+(z)=\mat{X}_-(z)\left(\begin{array}{cc}
1 & 0 \\\\ iY(z)^{-1}e^{\eta(z)-\gamma-2\kappa g(z)}e^{\pm
iN\theta_\Gamma} e^{\mp iN\theta^0(z)}e^{-N\xi_\Gamma(z)} & 1
\end{array}\right)\,.
\end{equation}
To express as concisely as possible the relationship between the
boundary values taken by $\mat{X}(z)$ on segments $\Sigma_{I\pm}$
parallel to a band $I$ it is convenient to choose some fixed $y\in I$
and then define $y_N$ for each $N\in\nat$ by the rule
\begin{equation}
N\int_a^{y_N}\rho^0(x)\,dx = \left\lceil
N\int_a^y\rho^0(x)\,dx\right\rceil\,,
\end{equation}
which may be compared with (\ref{eq:yquantize}).  Thus, if $I$ is a
transition band we may take $y_N$ to be the transition point
$y_{k,N}\in Y_N$ contained therein.  Otherwise we may think of $y_N$
as a ``virtual transition point''.  With the sequence
$\{y_N\}_{N=0}^\infty$ so determined, we have that for $z$ in any
segment $\Sigma_{I\pm}$ parallel to any band $I$,
\begin{equation}
\mat{X}_+(z)=\mat{X}_-(z)\left(\begin{array}{cc}
T_\del(z)^{-1/2} &v_{12}^\pm(z)
\\\\
v_{21}^\pm(z) & T_\nab(z)^{-1/2}
\end{array}\right)^{\pm 1}\,,
\end{equation}
where
\begin{equation}
\begin{array}{rcl}
v_{12}^\pm(z)&:=&\displaystyle
\mp i T_\del(z)^{-1/2}e^{\gamma-\eta(z)+2\kappa g(z)}e^{\mp iN\theta(y_N)}
\exp\left(\pm 2\pi iNc\int_{y_N}^z\overline{\psi}_I(s)\,ds\right)\,,
\\\\
v_{21}^\pm(z)&:=&\displaystyle
\mp iT_\nab(z)^{-1/2}e^{\eta(z)-\gamma-2\kappa g(z)}e^{\pm iN\theta(y_N)}
\exp\left(\pm 2\pi iNc\int_{y_N}^z\psi_I(s)\,ds\right)\,.
\end{array}
\end{equation}
\end{enumerate}
\label{rhp:X}
\end{rhp}

\begin{proof}
The domain of analyticity of $\mat{X}(z)$ is clear from the nature of
the definition \eqref{eq:XD} with
\eqref{eq:Dfirst}--\eqref{eq:Xyellowlower}, and from 
Proposition~\ref{prop:disappear}.  The normalization condition follows
from the corresponding normalization of $\mat{P}(z;N,k)$ and the from
\eqref{eq:rhonorm}.  The continuity of the boundary values is obvious
everywhere except on the real axis, but here the poles in
$\mat{P}(z;N,k)$ are cancelled by corresponding zeros in the boundary
values of $T_\nab(z)^{1/2}$ and $T_\del(z)^{1/2}$.  Finally, the jump
conditions are a direct consequence of the continuity of
$\mat{P}(z;N,k)$ and the known discontinuities of $\mat{D}(z)$.

This shows that $\mat{X}(z)$ defined by \eqref{eq:XD} with
\eqref{eq:Dfirst}--\eqref{eq:Xyellowlower} indeed satisfies all of the
conditions of Riemann-Hilbert Problem~\ref{rhp:X}.  The uniqueness of
the solution follows from Liouville's Theorem because the matrix ratio
of any two solutions is necessarily an entire function of $z$ that
tends to the identity matrix as $z\rightarrow\infty$.
\end{proof}

\section{Asymptotic Analysis}
\label{sec:asymptotics}
In this section we provide all the tools for a complete asymptotic
analysis of discrete orthogonal polynomials with a large class of
(generally nonclassical) weights, in the joint limit of large degree
and a large number of nodes.  These results will then be used in
\S~\ref{sec:asymptoticspi} to establish precise convergence theorems about the
discrete orthogonal polynomials and in
\S~\ref{sec:universality} to prove a number of universality results concerning
statistics of related discrete orthogonal polynomial ensembles.

\subsection{Construction of a Parametrix for $\mat{X}(z)$.}
\subsubsection{Outer asymptotics.}
Our immediate goal is to use the deformations we have carried out to
construct a model for the matrix $\mat{X}(z)$ that we expect to be
asymptotically accurate pointwise in $z$ as $N\rightarrow\infty$.  The
proof of validity will be given in \S~\ref{sec:error}.

The basic observation at this point, which we will justify more
precisely in \S~\ref{sec:error}, is that the jump matrix relating
$\mat{X}_+(z)$ and $\mat{X}_-(z)$ in Riemann-Hilbert
Problem~\ref{rhp:X} is closely approximated by the identity matrix in
the limit $N\rightarrow\infty$ for $z\in\Sigma_{\rm SD}\setminus
[a,b]$.  Moreover, the jump matrix in any gap $\Gamma\subset [a,b]$ is
closely approximated in the same limit by a constant matrix
$e^{-iN\theta_\Gamma\sigma_3}e^{-i\phi_\Gamma\sigma_3}$.  Neglecting
the errors on an ad-hoc basis leads to a model Riemann-Hilbert
problem. \label{symbol:matrixdotX}
\begin{rhp}
Let $\{\Gamma_j=(\beta_{j-1},\alpha_j)\,,\text{for $j=1,\dots,G$}\}$
denote the set of interior gaps in $(a,b)$, and let 
%the corresponding real constant values of
%$\theta(z)$ and $\phi(z)$ be $\theta_{\Gamma_i}$ and
%$\phi_{\Gamma_i}$.  Let 
the bands be denoted by $\{I_j=(\alpha_j,\beta_j)\,,\text{for
$j=0,\dots,G$}\}$.  Let $\Sigma_{\rm model}$ denote the interval $[\alpha_0,\beta_G]$, oriented from left to right.  Find a
$2\times 2$ matrix $\dot{\mat{X}}(z)$ with the following properties:
\begin{enumerate}
\item
{\bf Analyticity}: $\dot{\mat{X}}(z)$ is an analytic function of $z$
for $z\in\cx\setminus \Sigma_{\rm model}$.
\item
{\bf Normalization}: As $z\rightarrow\infty$,
\begin{equation}
\dot{\mat{X}}(z)={\mathbb I}+O\left(\frac{1}{z}\right)\,.
\end{equation}
\item
{\bf Jump Conditions}: $\dot{\mat{X}}(z)$ takes continuous boundary
values on $\Sigma_{\rm model}$ except at the endpoints of the bands,
where inverse fourth-root singularities are admitted.  For $z\in\Sigma_{\rm model}$, let $\dot{\mat{X}}_+(z)$ ($\dot{\mat{X}}_-(z)$) denote the boundary value taken by $\dot{\mat{X}}(z)$ on the left (right) of $\Sigma_{\rm model}$ according to its orientation.  For $z$ in the gap
$\Gamma_j$, the boundary values satisfy
\begin{equation}
\dot{\mat{X}}_+(z)=\dot{\mat{X}}_-(z)\left(\begin{array}{cc}
e^{iN\theta_{\Gamma_j}}e^{i\phi_{\Gamma_j}} & 0 \\ \\0 &
e^{-iN\theta_{\Gamma_j}}e^{-i\phi_{\Gamma_j}}
\end{array}\right)\,,
\label{eq:dotXgapjump}
\end{equation}
where the constant $\theta_{\Gamma_j}$ is defined by
\eqref{eq:thetaGammajvoid} or \eqref{eq:thetaGammajsat} depending on
whether $\Gamma_j$ is a void or a saturated region, and
$\phi_{\Gamma_j}$ is defined by \eqref{eq:phiGamma}, while for $z$ in
any band $I_j$, the boundary values satisfy
\begin{equation}
\dot{\mat{X}}_+(z)=\dot{\mat{X}}_-(z)\left(\begin{array}{cc}
0 & -ie^{\gamma-\eta(z)+\kappa(g_+(z)+g_-(z))} \\
\\-ie^{\eta(z)-\gamma-\kappa(g_+(z)+g_-(z))} & 0\end{array}\right)\,.
\label{eq:dotXbandjump}
\end{equation}
Here the expression $g_+(z)+g_-(z)$ refers to the sum of the boundary values
taken for $z\in I_j\subset\mathbb{R}$ from the upper and lower half-planes.
\end{enumerate}
\label{rhp:theta}
\end{rhp}

\begin{figure}[h]
\begin{center}
\input{SigmaModel.pstex_t}
\end{center}
\caption{\em The contour $\Sigma_{\rm model}$ corresponding to the
hypothetical equilibrium measure illustrated in Figure~\ref{fig:intervals} shown
against the dashed background of $\Sigma_{\rm SD}$.  Note that by contrast with $\Sigma_{\rm SD}$, the gap
intervals $\Gamma_j$ are now oriented from left to right.  Thus the
boundary value $\dot{\mat{X}}_+(z)$ ($\dot{\mat{X}}_-(z)$) 
refers to a limit from the upper (lower) half-plane.}
\label{fig:SigmaModel}
\end{figure}
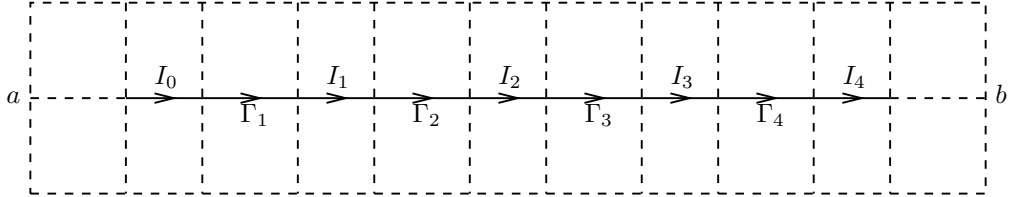

The contour $\Sigma_{\rm model}$ \label{symbol:Sigmamodel}
corresponding to the hypothetical
situation first illustrated in Figure~\ref{fig:intervals} is shown in
Figure~\ref{fig:SigmaModel}.
Problems of this sort are solved in terms of Riemann theta functions
of genus $G$, where $G+1$ is the number of bands $I_0,\dots,I_G$ (see,
for example, \cite{DeiftKMVZ99}).  Our subsequent analysis and error
estimates will not rely heavily on the specific formulae for the
solution, although as is clear from \S~\ref{sec:actualtheorems} these
details do emerge in the leading-order asymptotics justified by our
analysis.  For completeness, the solution of Riemann-Hilbert
Problem~\ref{rhp:theta} is explained in
Appendix~\ref{sec:thetasolve}.
%Note however, that it is in our solution of Riemann-Hilbert
%Problem~\ref{rhp:theta} that we determine the value of the constant
%$\gamma$ introduced with the transformation (\ref{eq:SfromR}); see
%(\ref{eq:gammadef}).

The essential facts we will require later are the following.
\begin{prop}
Riemann-Hilbert Problem~\ref{rhp:theta} has a unique solution
$\dot{\mat{X}}(z)$ that is uniformly bounded with bound independent of
$N$ in any neighborhood that does not contain any of the endpoints of
the bands $I_0,\dots,I_G$.  Although the numbers $\phi_{\Gamma_j}$
depend on the choice of transition points in the set $Y_N$, the
combination $\dot{\mat{X}}(z)e^{\kappa g(z)\sigma_3}$ is independent
of any particular choice of transition points.  Also,
$\det(\dot{\mat{X}}(z))=1$.
\label{prop:dotXbound}
\end{prop}
\begin{proof}
A solution is developed in detail in Appendix~\ref{sec:thetasolve} and
uniqueness can be established by an argument based on Liouville's
Theorem.  A similar argument proves that $\det(\dot{\mat{X}}(z))=1$.
The uniform boundedness of $\dot{\mat{X}}(z)$ away from the band
endpoints and the invariance of the combination
$\dot{\mat{X}}(z)e^{\kappa g(z)\sigma_3}$ are consequences of the
solution formulae given in Appendix~\ref{sec:thetasolve}; a discussion
of these features can be found there.
\end{proof}
The boundary values taken by the solution of Riemann-Hilbert
Problem~\ref{rhp:theta} have the following useful properties.
\begin{prop}
For $z$ in any interior gap (void or saturated region)
$\Gamma_j=(\beta_{j-1},\alpha_j)\subset\Sigma_{\rm model}$, we have the identity
\begin{equation}
\dot{\mat{X}}_+(z)e^{\kappa g_+(z)\sigma_3}=
\left(\dot{\mat{X}}_+(z)e^{\kappa g_+(z)\sigma_3}\right)^*
\left(\begin{array}{cc}e^{iN\theta_{\Gamma_j}} & 0 \\\\
0 & e^{-iN\theta_{\Gamma_j}}\end{array}\right)\,,
\end{equation}
where the star denotes componentwise complex conjugation.  
Similarly, for real $z<\alpha_0$, we have
\begin{equation}
\dot{\mat{X}}(z)e^{\kappa g_+(z)\sigma_3}=
\left(\dot{\mat{X}}(z)e^{\kappa g_+(z)\sigma_3}\right)^*\left(\begin{array}{cc}
e^{-2\pi iNc} & 0 \\\\ 0 & e^{2\pi iNc}\end{array}\right)\,,
\end{equation}
and for real $z>\beta_G$, we have
\begin{equation}
\dot{\mat{X}}(z)e^{\kappa g(z)\sigma_3}=
\left(\dot{\mat{X}}(z)e^{\kappa g(z)\sigma_3}\right)^*\,.
\end{equation}
Moreover, the product $p(z):=\dot{X}_{11}(z)\dot{X}_{12}(z)$ extends
to $\mathbb{C}\setminus ([\alpha_0,\beta_0]\cup\dots\cup
[\alpha_G,\beta_G])$ as a real-analytic function satisfying $p(z)<0$
for all real $z<\alpha_0$ and $p(z)>0$ for all
real $z>\beta_G$.  For all $j=1,\dots,G$, there
is a real number $z_j\in [\beta_{j-1},\alpha_j]$ such that $p(z)>0$
for $\beta_{j-1}< z < z_j$ and $p(z)<0$ for $z_j<z <\alpha_j$.  
If in fact $z_j\in (\beta_{j-1},\alpha_j)$, then $z_j$ is a simple
zero of $p(z)$.
%In particular, $p(z_j)=0$.  
The zeros $z_j$ depend on the 
parameter $\kappa$ in a quasiperiodic fashion with $G$ frequencies
that depend on the parameters $c\in (0,1)$ and $N$, 
the function $\eta(z)$, and the equilibrium measure.  Generically, $z_j\in(\beta_{j-1},\alpha_j)$, and the
situation in which $z_j=\beta_{j-1}$ or $z_j=\alpha_j$ for
some $j$ should be regarded as exceptional.  In the generic case, the boundary values
$\dot{X}_{11+}(z)$ and $\dot{X}_{12+}(z)$ are analytic at $z=z_j$ and thus 
either $\dot{X}_{11+}(z)$ has a
simple zero only at $z=z_j$ and $\dot{X}_{12+}(z)$ is bounded away
from zero in $\Gamma_j$, or $\dot{X}_{12+}(z)$ has a simple zero only
at $z=z_j$ and $\dot{X}_{11+}(z)$ is bounded away from zero in
$\Gamma_j$.

For $z$ in any band $I_j$, the identity
\begin{equation}
\dot{\mat{X}}_+(z)e^{\kappa g_+(z)\sigma_3}=
\left(\dot{\mat{X}}_+(z)e^{\kappa g_+(z)\sigma_3}\right)^*
\left(\begin{array}{cc} 0 & -ie^{\gamma-\eta(z)}\\\\
-ie^{\eta(z)-\gamma} & 0\end{array}\right)
\label{eq:dotXbandsymmetry}
\end{equation}
holds.  Furthermore, for $z\in I_j$ the elements of
$\dot{\mat{X}}_+(z)$ are strictly nonzero.
\label{prop:dotXsymmetry}
\end{prop}
\begin{proof}
The matrix $\mat{M}(z):=\dot{\mat{X}}(z)e^{\kappa g(z)\sigma_3}$ and
the corresponding matrix $\mat{N}(z):=\mat{M}(z^*)^*$ are both
analytic for $z\in\mathbb{C}\setminus (-\infty,\beta_G]$, where
$\beta_G$ is the rightmost band endpoint.  As $z\rightarrow\infty$, we
have $\mat{M}(z)e^{-\kappa\log(z)\sigma_3}=\mathbb{I}+O(1/z)$ and also
$\mat{N}(z)e^{-\kappa\log(z)\sigma_3}=\mathbb{I}+O(1/z)$.
Furthermore, it is easily checked that at each point $z\in
(-\infty,\beta_G]$, we have
$\mat{M}_-(z)^{-1}\mat{M}_+(z)=\mat{N}_-(z)^{-1}\mat{N}_+(z)$.  This
means that both matrices satisfy the same Riemann-Hilbert problem.
Uniqueness of solutions for this problem follows as usual from
Liouville's Theorem.  Thus, $\mat{M}(z)=\mat{N}(z)=\mat{M}(z^*)^*$.
The claimed relations follow from the jump relations for
$\dot{\mat{X}}(z)$ since for each real $z$,
$\mat{M}_-(z)=\mat{M}_+(z)^*$.

Suppose that at some point $z$ in a band $I_j$ we have $M_{11+}(z)=0$.
Then it follows from (\ref{eq:dotXbandsymmetry}) that $M_{12+}(z)=0$
also.  But this implies that $\det(\mat{M}(z))=0$ which contradicts
the fact that (see Proposition~\ref{prop:dotXbound})
$\det(\mat{M}(z))=1$.  In a similar way, one sees that any other
matrix element of $\mat{M}_\pm(z)$ having a zero in $I_j$ leads to a
contradiction.

The fact that the product $p(z)=\dot{X}_{11}(z)\dot{X}_{12}(z)$
extends to the complement of the bands $I_j$ as an analytic function
follows from the jump condition (\ref{eq:dotXgapjump}) and the
analyticity of $\dot{\mat{X}}(z)$ for
$z\in\mathbb{C}\setminus\Sigma_{\rm model}$.  The sign of $p(z)$ is
discussed in detail in Appendix~\ref{sec:thetasolve}.
\end{proof}

By using the explicit formulae given in Appendix~\ref{sec:thetasolve},
one can obtain the identities $W(z)\equiv \dot{X}_{11}(z)e^{\kappa g(z)}$
and $Z(z)\equiv\dot{X}_{12}(z)e^{-\kappa g(z)}$, where $W(z)$ and $Z(z)$
are the functions defined in \eqref{eq:Wdefine} and \eqref{eq:Zdefine} respectively.  

\subsubsection{Inner asymptotics near band edges.}
\label{sec:Airy}
In any neighborhood of a point in the interior of either
$\Sigma_0^\nab$ or $\Sigma_0^\del$ that marks the boundary between a
band and a gap, the pointwise asymptotics used to arrive at the jump
conditions for the matrix $\dot{\mat{X}}(z)$ starting from those for
the matrix $\mat{X}(z)$ are not uniformly valid.  It is therefore
necessary to construct a local approximation to $\mat{X}(z)$ near such
points using different techniques.  We refer to these boundary points
separating bands from gaps as {\em band edges}.  We want to stress
that band edges are to be distinguished from transition points making
up the set $Y_N$ defined in \S~\ref{sec:choiceofdel}.  Our method will be to
define in a disc of fixed size near each band edge a matrix that
exactly satisfies the jump conditions of $\mat{X}(z)$ and that matches
well onto the outer asymptotics given by $\dot{\mat{X}}(z)$ at the
boundary of the disc.

The distinguishing characteristic of a band edge $z=z_0$ is that in
the adjacent gap $\Gamma$ the function $\rho(z)$ is identically zero
since the equilibrium measure $\mu_{\rm min}^c$ realizes the lower constraint
for $z\in\Gamma$ if $\Gamma\subset\Sigma_0^\nab$ or the upper
constraint for $z\in\Gamma$ if $\Gamma\subset\Sigma_0^\del$, and
meanwhile in the adjacent band $\rho(z)$ is a nonzero analytic
function that vanishes at the band edge.  The nature of the vanishing
of $\rho(z)$ at the band edge must be understood before a local
approximation can be constructed.  Consider $\delta
E_c/\delta\mu-\ell_c$ where the variational derivative is evaluated on
the equilibrium measure.  In the band this quantity is identically zero
according to the equilibrium condition (\ref{eq:equilibrium}).  On the
other hand if $\rho^0(\cdot)$ and $V(\cdot)$ are analytic functions
then the function $\Psi(z)$ \label{symbol:Psi} defined for $z\in (a,b)$ by
\begin{equation}
\Psi(z):=V(z)+\int_{\Sigma_0^\nab}\log|z-x|\rho^0(x)\,dx -\int_{\Sigma_0^\del}\log|z-x|\rho^0(x)\,dx
\label{eq:Psidef}
\end{equation}
extends analytically into the upper half-plane (it is analytic in a
neighborhood of $z_0$ as long as $z_0$ is in the interior of either
$\Sigma_0^\nab$ or $\Sigma_0^\del$).   Since
\begin{equation}
\Psi(z)+2(d_N-c)\int_a^b\log|z-x|\rho(x)\,dx = \frac{\delta E_c}{\delta \mu}(z)\,,\hspace{0.2 in}\mbox{for $z\in (a,b)$,}
\label{eq:Psivariationalderivative}
\end{equation}
where the variational
derivative is evaluated on the equilibrium measure, 
we have
\begin{equation}
0 \equiv \Psi(z)+2(d_N-c)\int_a^b\log(z-x)\rho(x)\,dx - \ell_c -2\pi i
(d_N-c)\int_{z_0}^b\rho(x)\,dx + 2\pi i
(d_N-c)\int_{z_0}^z\rho(x)\,dx\,,
\label{eq:equivUHP}
\end{equation}
for $z$ near $z_0$ with $\Im(z)>0$.
Only the last integral involves contour integration off of the real
axis, and the integrand denotes the analytic function $\rho(\cdot)$ of
the band.  At the same time, the quantity $\delta
E_c/\delta\mu-\ell_c$ extends into the upper half-plane from the gap
$\Gamma$ as
\begin{equation}
\frac{\delta E_c}{\delta\mu}-\ell_c\Bigg|_{z\in\Gamma}=
\Psi(z)+2(d_N-c)\int_a^b\log(z-x)\rho(x)\,dx - \ell_c
-2\pi i (d_N-c)\int_{z_0}^b\rho(x)\,dx
\end{equation}
since $\rho(\cdot)\equiv 0$ for $z\in\Gamma$.  We therefore deduce
that
\begin{equation}
\frac{\delta E_c}{\delta\mu}-\ell_c\Bigg|_{z\in\Gamma}=
-2\pi i(d_N-c)\left[\int_{z_0}^z\rho(x)\,dx\right]_+
\label{eq:gapplus}
\end{equation}
where on the right-hand side the integrand is the continuation of the
analytic function $\rho(\cdot)$ defined in the adjacent band, and the
subscript denotes the boundary value taken on the gap $\Gamma$ from
the upper half-plane.  Using virtually the same arguments but
continuing all quantities into the lower half-plane, we find that
\begin{equation}
\frac{\delta E_c}{\delta\mu}-\ell_c\Bigg|_{z\in\Gamma}=
2\pi i(d_N-c)\left[\int_{z_0}^z\rho(x)\,dx\right]_-\,.
\label{eq:gapminus}
\end{equation}
Combining (\ref{eq:gapplus}) and (\ref{eq:gapminus}) reveals the
identity
\begin{equation}
\left[\int_{z_0}^z\rho(x)\,dx\right]_++\left[\int_{z_0}^z\rho(x)\,dx\right]_-=0
\label{eq:rhoplusminus}
\end{equation}
which holds for all $z$ in the gap when the integrand $\rho(\cdot)$ is
analytically extended about $z_0$ from the band.  Differentiating this
identity with respect to $z$, we discover that $\rho(z)^2$ extends
from the band to a complex annulus surrounding $z_0$ as a
single-valued analytic function that vanishes as $z\rightarrow z_0$
within the band (at least).  Moreover, it follows from
(\ref{eq:equivUHP}) that $\rho(z)^2$ is analytic at $z_0$ as well, and
so is necessarily of the form $\rho(z)^2= (z-z_0)^pe^{f(z)}$ where
$p=1,2,3,\dots$ and $f(z)$ is analytic at $z_0$.  

Clearly, only odd values of the positive integer $p$ are consistent
with (\ref{eq:rhoplusminus}).  However, even more is true.  If the
band edge point satisfies $z_0\in\Sigma_0^\nab$, then the combination
$(c-d_N)\rho(x)$ can be seen by \eqref{eq:rhodef} to be strictly
positive for $x$ in the band adjacent to $z_0$, and furthermore the
adjacent gap is a {\em void}, and thus from \eqref{eq:voidinequality}
we see that the common left-hand side of \eqref{eq:gapplus} and
\eqref{eq:gapminus} is strictly positive for $z\in\Gamma$.  
Similarly if the band edge
point satisfies $z_0\in\Sigma_0^\del$, then the combination
$(c-d_N)\rho(x)$ is strictly negative for $x$ in the band adjacent to
$z_0$, and the adjacent gap is a {\em saturated region} so that
\eqref{eq:saturatedregioninequality} makes the common left-hand side
of \eqref{eq:gapplus} and \eqref{eq:gapminus} strictly negative for
$z\in\Gamma$.  In both cases, we can easily see that the equations
\eqref{eq:gapplus} and \eqref{eq:gapminus} will only be consistent with
the assumption that $\rho(z)^2 = (z-z_0)^pe^{f(z)}$ for analytic $f(z)$
and $p=1,3,5,7,\dots$ if we discard the values $p=3,7,11,\dots$.  

Therefore, using only the assumption that $\rho^0(\cdot)$ and
$V(\cdot)$ are analytic functions, we have shown that at each band
edge $z_0$ in the interior of $\Sigma_0^\nab$ or $\Sigma_0^\del$ the
positive analytic function $\rho(\cdot)$ vanishes like $(z-z_0)^{p/2}$
where $p$ is of the form $p=1+4m$ for $m=0,1,2,3,\dots$.  This is the
general character of the vanishing of $\rho(\cdot)$ at band edges when
$V(\cdot)$ and $\rho^0(\cdot)$ are analytic functions, and it is quite similar
to the characterization of the local behavior of the equilibrium measure
(without upper constraint) near band edges as explained in \cite{DeiftKM98}.

As mentioned in
\S~\ref{sec:C3} ({\em cf.} in particular \eqref{eq:rootlower}
and
\eqref{eq:rootupper}), we will for simplicity consider only the generic
situation when $p=1$ at all band edges.  There are four cases.  Let
$ h<1$ \label{symbol:hparameter} be an arbitrary fixed positive parameter.

\subsubsection*{Left band edge with $z_0=\alpha\in\Sigma_0^\nab$ (lower constraint).}
Let $\Gamma$ denote the void to the left of $\alpha$; then
$e^{iN\theta(\alpha)}=e^{iN\theta_\Gamma}$.  Let $I$ denote the band to the
right of $\alpha$.  Consider $D^{\nab,L}_\Gamma$ 
\label{symbol:DnabL} to be an open disc
centered at $z=\alpha$ of radius $h\epsilon$.  Note that for
$\epsilon$ sufficiently small this radius will be less than half the
distance to the nearest distinct band edge and $D^{\nab,L}_\Gamma$
will be disjoint from the endpoints $\{a,b\}$.  We divide
$D^{\nab,L}_\Gamma\setminus(D^{\nab,L}_\Gamma\cap \Sigma_{\rm SD})$
into \label{symbol:DnabLquadrants} open quadrants:
\begin{equation}
\begin{array}{rcl}
D^{\nab,L}_{\Gamma,I}&=&\displaystyle
D^{\nab,L}_\Gamma\cap\left\{z\,\Big|\,z\neq
\alpha,\,0<\arg(z-\alpha)<\frac{\pi}{2}\right\}\,,\\\\
D^{\nab,L}_{\Gamma,II}&=&\displaystyle
D^{\nab,L}_\Gamma\cap\left\{z\,\Big|\,z\neq
\alpha,\,\frac{\pi}{2}<\arg(z-\alpha)<\pi\right\}\,,\\\\
D^{\nab,L}_{\Gamma,III}&=&\displaystyle
D^{\nab,L}_\Gamma\cap\left\{z\,\Big|\,z\neq
\alpha,\,-\pi<\arg(z-\alpha)<-\frac{\pi}{2}\right\}\,,\\\\
D^{\nab,L}_{\Gamma,IV}&=&\displaystyle
D^{\nab,L}_\Gamma\cap\left\{z\,\Big|\,z\neq
\alpha,\,-\frac{\pi}{2}<\arg(z-\alpha)<0\right\}\,.
\end{array}
\end{equation}
Now we introduce a local change of variables in $D^{\nab,L}_\Gamma$.
We set
\begin{equation}
\mat{Z}^{\nab,L}_\Gamma(z):=
\left\{\begin{array}{ll}
\mat{X}(z)e^{(\gamma-\eta(z)+2\kappa g(z))\sigma_3/2}e^{-iN\theta_\Gamma\sigma_3/2}\,,&\hspace{0.2 in}\mbox{for $z\in D^{\nab,L}_{\Gamma,I}$},\\\\
\mat{X}(z)T_\nab(z)^{\sigma_3/2}e^{(\gamma-\eta(z)+2\kappa g(z))\sigma_3/2}e^{-iN\theta_\Gamma\sigma_3/2}\,,&\hspace{0.2 in}\mbox{for $z\in D^{\nab,L}_{\Gamma,II}$},\\\\
\mat{X}(z)T_\nab(z)^{\sigma_3/2}e^{(\gamma-\eta(z)+2\kappa g(z))\sigma_3/2}e^{iN\theta_\Gamma\sigma_3/2}\,,&\hspace{0.2 in}\mbox{for $z\in D^{\nab,L}_{\Gamma,III}$},\\\\
\mat{X}(z)e^{(\gamma-\eta(z)+2\kappa g(z))\sigma_3/2}e^{iN\theta_\Gamma\sigma_3/2}\,,&\hspace{0.2 in}\mbox{for $z\in D^{\nab,L}_{\Gamma,IV}$}.
\end{array}\right.
\label{eq:XZnabL}
\end{equation}
\label{symbol:ZnabL}

According to \eqref{eq:rootlower}, the equation $\zeta=\tau_\Gamma^{\nab,L}(z)$ defined by \eqref{eq:zetadefAL}
gives an invertible conformal mapping taking, for $\epsilon$ sufficiently
small, the fixed disc $D_\Gamma^{\nab,L}$ to a neighborhood of $\zeta=0$ in the $\zeta$-plane
that scales like $N^{2/3}$.
%change of variables
%\begin{equation}
%\zeta=\tau^{\nab,L}_\Gamma(z):=
%\left(2\pi Nc\int_{\alpha}^zd\mu_{\rm min}^c(x)\right)^{2/3}
%\label{eq:zetadefAL}
%\end{equation}
%defined for $z>\alpha$ in $D^{\nab,L}_\Gamma$ (in the adjacent band)
%extends to all of $D^{\nab,L}_\Gamma$ as an analytic transformation
%$\zeta\leftrightarrow z$.  
The transformation
$\tau^{\nab,L}_\Gamma(z)$ maps ${\mathbb R}\cap D^{\nab,L}_\Gamma$ to
$\mathbb R$, taking $z=\alpha$ to $\zeta=0$ and is
orientation-preserving since $d\tau^{\nab,L}_\Gamma/dz(\alpha)$ is
real and positive.  The segments $\arg(z-\alpha)=\pm\pi/2$ in
$D^{\nab,L}_\Gamma$ are mapped to arcs in the $\zeta$-plane that are
tangent to the imaginary axis at $\zeta=0$ and that converge to the
rays $\arg(\zeta)=\pm\pi/2$ as $N\rightarrow\infty$ uniformly for
$\zeta$ in compact sets.  The exact jump conditions satisfied by the
boundary values of $\mat{Z}^{\nab,L}_\Gamma(z)$ on $\Sigma_{\rm
SD}\cap D^{\nab,L}_\Gamma$ may be written in terms of the new
coordinate $\zeta$ as follows:
\begin{equation}
\begin{array}{rcll}
\mat{Z}^{\nab,L}_{\Gamma+}(z)&=&\displaystyle\mat{Z}^{\nab,L}_{\Gamma-}(z)
\left(\begin{array}{cc}
1 & ie^{-(-\zeta)^{3/2}}\\\\ 0 & 1\end{array}\right)\,,&\hspace{0.2
in}
\mbox{for $z\in\Gamma\cap D_\Gamma^{\nab,L}$,}\\\\
\mat{Z}^{\nab,L}_{\Gamma+}(z)&=&\displaystyle
\mat{Z}^{\nab,L}_{\Gamma-}(z)\left(\begin{array}{cc}
0 & -i \\\\ -i & 0\end{array}\right)\,,&\hspace{0.2 in}
\mbox{for $z\in I\cap D_\Gamma^{\nab,L}$,}\\\\
\mat{Z}^{\nab,L}_{\Gamma+}(z)&=&\displaystyle\mat{Z}^{\nab,L}_{\Gamma-}(z)
\left(\begin{array}{cc}
1 & 0\\\\ -ie^{i\zeta^{3/2}} & 1\end{array}\right)\,,&\hspace{0.2 in}
\mbox{for $z\in\Sigma_{0+}^\nab\cap D_\Gamma^{\nab,L}$,}\\\\
\mat{Z}^{\nab,L}_{\Gamma+}(z)&=&\displaystyle\mat{Z}^{\nab,L}_{\Gamma-}(z)
\left(\begin{array}{cc}
1 & 0\\\\ -ie^{-i\zeta^{3/2}} & 1\end{array}\right)\,,&\hspace{0.2 in}
\mbox{for $z\in\Sigma_{0-}^\nab\cap D_\Gamma^{\nab,L}$.}
\end{array}
\end{equation}
Here, the subscripts ``$+$'' and ``$-$'' refer to boundary values taken
on $\Sigma_{\rm SD}\cap D_\Gamma^{\nab,L}$ respectively from the left and right relative to the orientation of $\Sigma_{\rm SD}$.

At the same time, we can define a ``comparison matrix''
$\dot{\mat{Z}}_\Gamma^{\nab,L}(z)$ 
\label{symbol:dotZnabL} from $\dot{\mat{X}}(z)$ by the
relation
\begin{equation}
\dot{\mat{Z}}^{\nab,L}_\Gamma(z):=\dot{\mat{X}}(z)e^{(\gamma-\eta(z)+2\kappa g(z))\sigma_3/2}e^{-iN{\rm sgn}(\Im(z))\theta_\Gamma\sigma_3/2}\,,\hspace{0.2 in}\mbox{for $z\in D^{\nab,L}_{\Gamma}\setminus(
D^{\nab,L}_{\Gamma}\cap\Sigma_{\rm SD})$}.
\label{eq:XZdotnabL}
\end{equation}
Note the difference (a factor of $T_\nab(z)^{\sigma_3/2}$ in quadrants
$II$ and $III$) between the transformation (\ref{eq:XZdotnabL}) and
the transformation (\ref{eq:XZnabL}).  This matrix extends to an
analytic function in $D^{\nab,L}_\Gamma$ with the exception of $z\in
I\cap D_\Gamma^{\nab,L}$, where it satisfies
\begin{equation}
\dot{\mat{Z}}^{\nab,L}_{\Gamma+}(z)=\dot{\mat{Z}}^{\nab,L}_{\Gamma-}(z)
\left(\begin{array}{cc}
0 & -i \\\\ -i & 0\end{array}\right)\hspace{0.2 in}
\mbox{for $z\in I\cap D_\Gamma^{\nab,L}$.}
\end{equation}
Again, the subscripts indicate boundary values consistent with the orientation of $\Sigma_{\rm SD}$, with ``$+$'' indicating approach from the left and ``$-$'' indicating approach from the right.  
Because the matrix elements of $\dot{\mat{Z}}(z)$ blow up no worse
than $(z-\alpha)^{-1/4}$, it is easy to see that
$\dot{\mat{Z}}_\Gamma^{\nab,L}(z)$ can be represented in the form
\label{symbol:HGammanabL}
\begin{equation}
\dot{\mat{Z}}_\Gamma^{\nab,L}(z)=
\mat{H}_\Gamma^{\nab,L}(z)\cdot \frac{1}{\sqrt{2}}
(-\tau^{\nab,L}_\Gamma(z))^{\sigma_3/4}
\left(\begin{array}{cc}1 & 1\\\\-1 & 1\end{array}\right)\,,
\label{eq:HdefAL}
\end{equation}
where $\mat{H}_\Gamma^{\nab,L}(z)$ is analytic in $D_\Gamma^{\nab,L}$.
The relations (\ref{eq:HdefAL}) and (\ref{eq:XZdotnabL}) together with
(\ref{eq:zetadefAL}) serve as a definition of
$\mat{H}_\Gamma^{\nab,L}(z)$ in terms of the solution
$\dot{\mat{X}}(z)$ of Riemann-Hilbert Problem~\ref{rhp:theta}.

Since the image of the boundary of $D_\Gamma^{\nab,L}$ in the
$\zeta$-plane expands as $N\rightarrow\infty$ with $\epsilon$ held
fixed, and since on the boundary $\mat{Z}_\Gamma^{\nab,L}(z)$ and
$\dot{\mat{Z}}_\Gamma^{\nab,L}(z)$ should be comparable, we propose to
concretely determine an approximation of $\mat{Z}_\Gamma^{\nab,L}(z)$
for $z\in D_\Gamma^{\nab,L}$ by solving the following Riemann-Hilbert
problem.
\label{symbol:hatZnabL}
\begin{rhp}
Let $C_+$ be a contour connecting the origin to infinity lying
entirely within a symmetrical sector about the positive imaginary axis
of opening angle strictly less than $\pi/3$.  Let $C_-$ denote the
complex-conjugate of $C_+$.  Find a $2\times 2$ matrix
$\hat{\mat{Z}}^{\nab,L}(\zeta)$ with the following properties:
\begin{enumerate}
\item
{\bf Analyticity}: $\hat{\mat{Z}}^{\nab,L}(\zeta)$ is an analytic
function of $\zeta$ for $\zeta\in\mathbb{C}\setminus(\mathbb{R}\cup
C_+\cup C_-)$.
\item
{\bf Normalization}: As $\zeta\rightarrow\infty$,
\begin{equation}
\hat{\mat{Z}}^{\nab,L}(\zeta)\cdot \frac{1}{\sqrt{2}}\left(\begin{array}{cc}
1 & -1\\\\1 & 1\end{array}\right)(-\zeta)^{-\sigma_3/4} = {\mathbb I}
+ O\left(\frac{1}{\zeta}\right)\,,
\label{eq:Airynormalize}
\end{equation}
uniformly with respect to direction.
\item
{\bf Jump Conditions}: $\hat{\mat{Z}}^{\nab,L}(\zeta)$ takes
continuous boundary values from each sector of its analyticity.  The
boundary values satisfy
\begin{equation}
\begin{array}{rcll}
\hat{\mat{Z}}^{\nab,L}_+(\zeta)&=&\displaystyle\hat{\mat{Z}}^{\nab,L}_-(\zeta)
\left(\begin{array}{cc}
1 & ie^{-(-\zeta)^{3/2}}\\\\ 0 & 1\end{array}\right)\,,&\hspace{0.2
in}
\mbox{for $\zeta\in{\mathbb R}$ and $\zeta<0$,}\\\\
\hat{\mat{Z}}^{\nab,L}_+(\zeta)&=&\displaystyle\hat{\mat{Z}}^{\nab,L}_-(\zeta)
\left(\begin{array}{cc}
0 & -i \\\\ -i & 0\end{array}\right)\,,&\hspace{0.2 in}
\mbox{for $\zeta\in{\mathbb R}$ and $\zeta>0$,}\\\\
\hat{\mat{Z}}^{\nab,L}_+(\zeta)&=&\displaystyle\hat{\mat{Z}}^{\nab,L}_-(\zeta)
\left(\begin{array}{cc}
1 & 0\\\\ -ie^{i\zeta^{3/2}} & 1\end{array}\right)\,,&\hspace{0.2 in}
\mbox{for $\zeta\in C_+$,}\\\\
\hat{\mat{Z}}^{\nab,L}_+(\zeta)&=&\displaystyle\hat{\mat{Z}}^{\nab,L}_-(\zeta)
\left(\begin{array}{cc}
1 & 0\\\\ -ie^{-i\zeta^{3/2}} & 1\end{array}\right)\,,&\hspace{0.2 in}
\mbox{for $\zeta\in C_-$.}
\end{array}
\end{equation}
To determine the boundary values, the contours on the real
$\zeta$-axis are oriented away from the origin, and the contours $C_+$
and $C_-$ are oriented toward the origin.  As usual, ``$+$'' indicates approach from the left and ``$-$'' indicates approach from the right.
\end{enumerate}
\label{rhp:Airy}
\end{rhp}
Note that the asymptotic behavior of $\hat{\mat{Z}}^{\nab,L}(\zeta)$
is chosen to match the explicit terms in
$\dot{\mat{Z}}_\Gamma^{\nab,L}(z)$ with the exception of the
holomorphic prefactor $\mat{H}^{\nab,L}_\Gamma(z)$, the effect of
which will be included after solving for
$\hat{\mat{Z}}^{\nab,L}(\zeta)$.  The solution of Riemann-Hilbert
Problem~\ref{rhp:Airy} was first found in \cite{DeiftZ95}, and we
provide it in the notation of our problem for completeness.

\begin{prop}[Deift and Zhou]\label{propAirynabL}
The unique solution of Riemann-Hilbert Problem~\ref{rhp:Airy} is given
by the following explicit formulae.  Let
\begin{equation}
w:=\left(\frac{3}{4}\right)^{2/3}\zeta\,.
\end{equation}
For $\zeta$ between the positive real axis and the contour $C_+$:
\begin{equation}
\hat{\mat{Z}}^{\nab,L}(\zeta):=
\left(\begin{array}{cc} \displaystyle e^{\frac{2\pi i}{3}}\sqrt{2\pi}
\left(\frac{3}{4}\right)^{-\frac{1}{6}}
e^{\frac{2iw^{3/2}}{3}}Ai'
\left(e^{\frac{\pi i}{3}}w \right) &
\displaystyle e^{\frac{5\pi i}{6}}
\sqrt{2\pi}\left(\frac{3}{4}\right)^{-\frac{1}{6}}
e^{-\frac{2iw^{3/2}}{3}}Ai'
\left(e^{-\frac{\pi i}{3}}w\right)\\\\
\displaystyle e^{-\frac{2\pi i}{3}}
\sqrt{2\pi}\left(\frac{3}{4}\right)^{\frac{1}{6}}
e^{\frac{2iw^{3/2}}{3}}Ai\left(e^{\frac{\pi i}{3}} w
\right) & \displaystyle e^{\frac{\pi i}{6}}
\sqrt{2\pi}\left(\frac{3}{4}\right)^{\frac{1}{6}}
e^{-\frac{2iw^{3/2}}{3}}Ai\left(e^{-\frac{\pi i}{3}} w
\right)
\end{array}\right)\,.
\label{eq:Airyfirst}
\end{equation}
For $\zeta$ between the positive real axis and the contour $C_-$:
\begin{equation}
\hat{\mat{Z}}^{\nab,L}(\zeta):=
\left(\begin{array}{cc}
\displaystyle e^{-\frac{2\pi i}{3}}
\sqrt{2\pi}\left(\frac{3}{4}\right)^{-\frac{1}{6}}
e^{-\frac{2iw^{3/2}}{3}} Ai'\left(e^{-\frac{\pi i}{3}} w
\right) &
\displaystyle e^{-\frac{5\pi i}{6}}
\sqrt{2\pi}\left(\frac{3}{4}\right)^{-\frac{1}{6}}
e^{\frac{2iw^{3/2}}{3}} Ai'\left(e^{\frac{\pi i}{3}}w \right)
\\\\
\displaystyle e^{\frac{2\pi i}{3}}
\sqrt{2\pi}\left(\frac{3}{4}\right)^{\frac{1}{6}}
e^{-\frac{2iw^{3/2}}{3}} Ai\left(e^{-\frac{\pi i}{3}}w
\right) &
\displaystyle e^{-\frac{\pi i}{6}}
\sqrt{2\pi}\left(\frac{3}{4}\right)^{\frac{1}{6}}
e^{\frac{2iw^{3/2}}{3}} Ai\left(e^{\frac{\pi i}{3}}w
\right)
\end{array}\right)\,.
\end{equation}
For $\zeta$ between the contour $C_+$ and the negative real axis:
\begin{equation}
\hat{\mat{Z}}^{\nab,L}(\zeta):=
\left(\begin{array}{cc}
\displaystyle
-\sqrt{2\pi}\left(\frac{3}{4}\right)^{-\frac{1}{6}}
e^{\frac{2(-w)^{3/2}}{3}}Ai'\left(-w\right) &
\displaystyle e^{\frac{5\pi i}{6}}
\sqrt{2\pi}\left(\frac{3}{4}\right)^{-\frac{1}{6}}
e^{-\frac{2(-w)^{3/2}}{3}} Ai'\left(e^{-\frac{\pi i}{3}}w
\right)\\\\
\displaystyle
-\sqrt{2\pi}\left(\frac{3}{4}\right)^{\frac{1}{6}}
e^{\frac{2(-w)^{3/2}}{3}}Ai\left(-w\right) &
\displaystyle e^{\frac{\pi i}{6}}
\sqrt{2\pi}\left(\frac{3}{4}\right)^{\frac{1}{6}}
e^{-\frac{2(-w)^{3/2}}{3}} Ai\left(e^{-\frac{\pi i}{3}}w
\right)
\end{array}\right)\,.
\label{eq:Airythird}
\end{equation}
Finally, for $\zeta$ between the contour $C_-$ and the negative real
axis:
\begin{equation}
\hat{\mat{Z}}^{\nab,L}(\zeta):=
\left(\begin{array}{cc}
\displaystyle
-\sqrt{2\pi}\left(\frac{3}{4}\right)^{-\frac{1}{6}}
e^{\frac{2(-w)^{3/2}}{3}}Ai'\left(-w\right) &
\displaystyle
e^{-\frac{5\pi
i}{6}}\sqrt{2\pi}\left(\frac{3}{4}\right)^{-\frac{1}{6}}
e^{-\frac{2(-w)^{3/2}}{3}} Ai'\left(e^{\frac{\pi i}{3}}w \right)
\\\\
\displaystyle
-\sqrt{2\pi}\left(\frac{3}{4}\right)^{\frac{1}{6}}e^{\frac{2(-w)^{3/2}}{3}}
Ai\left(-w\right) &
\displaystyle e^{-\frac{\pi i}{6}}
\sqrt{2\pi}\left(\frac{3}{4}\right)^{\frac{1}{6}}
e^{-\frac{2(-w)^{3/2}}{3}} Ai\left(e^{\frac{\pi i}{3}}w \right)
\end{array}\right)\,.
\label{eq:Airylast}
\end{equation}
\end{prop}

\begin{proof}
The jump conditions are easily verified with the help of the identity
\begin{equation}
Ai(z) + e^{\frac{2\pi i}{3}}Ai(e^{\frac{2\pi i}{3}}z)+e^{-\frac{2\pi
i}{3}}Ai(e^{-\frac{2\pi i}{3}}z)=0\,.
\label{eq:AiryTrinity}
\end{equation}
The asymptotics are verified with the use of the steepest descent
asymptotic formulae
\begin{equation}
\begin{array}{rcl}
Ai(z)&=&\displaystyle
\frac{1}{2\sqrt{\pi}}z^{-1/4}e^{-2z^{3/2}/3}(1+O(z^{-3/2}))\\\\
Ai'(z)&=&\displaystyle
-\frac{1}{2\sqrt{\pi}}z^{1/4}e^{-2z^{3/2}/3}(1+O(z^{-3/2}))
\end{array}
\label{eq:AiryAsymp}
\end{equation}
both of which hold as $z\rightarrow\infty$ with $-\pi<\arg(z)<\pi$.
In fact, these calculations show that the $O(\zeta^{-1})$ error term
in the normalization condition (\ref{eq:Airynormalize}) is of a more
precise form, namely
\begin{equation}
\hat{\mat{Z}}^{\nab,L}(\zeta)\cdot \frac{1}{\sqrt{2}}\left(\begin{array}{cc}
1 & -1\\\\1 & 1\end{array}\right)(-\zeta)^{-\sigma_3/4} =
\left(\begin{array}{cc}
1+O(\zeta^{-3/2}) & O(\zeta^{-1})\\\\O(\zeta^{-2}) & 1+O(\zeta^{-3/2})
\end{array}\right)\,.
\label{eq:Airydecay}
\end{equation}
In this sense the decay rate to the identity matrix of $1/\zeta$ is
only sharp in one of the matrix elements, with the remaining matrix
elements exhibiting more rapid decay.  Uniqueness of the solution
follows from Liouville's Theorem.
\end{proof}

The contours $C_\pm$ in Riemann-Hilbert Problem~\ref{rhp:Airy} are
chosen so that in $\tau_\Gamma^{\nab,L}(D_\Gamma^{\nab,L})$ they agree
with the images under $\tau_\Gamma^{\nab,L}$ of the segments
$\Sigma_{0\pm}^\nab\cap D_\Gamma^{\nab,L}$.  Thus, the sectorial
condition on $C_\pm$ can be satisfied by taking the contour parameter
$\epsilon$ controlling the radius of $D_\Gamma^{\nab,L}$ to be
sufficiently small.  We now define a local parametrix for $\mat{X}(z)$
by the formula \label{symbol:XhatnabL}
\begin{equation}\label{eq;XhatnabL}
\hat{\mat{X}}_\Gamma^{\nab,L}(z):=
\left\{\begin{array}{ll}
\mat{H}_\Gamma^{\nab,L}(z)\hat{\mat{Z}}^{\nab,L}(\tau_\Gamma^{\nab,L}(z))
e^{(\eta(z)-\gamma-2\kappa
g(z))\sigma_3/2}e^{iN\theta_\Gamma\sigma_3/2}\,, &\hspace{0.2
in}\mbox{for $z\in D_{\Gamma,I}^{\nab,L}$,}\\\\
\mat{H}_\Gamma^{\nab,L}(z)\hat{\mat{Z}}^{\nab,L}(\tau_\Gamma^{\nab,L}(z))
T_\nab(z)^{-\sigma_3/2}e^{(\eta(z)-\gamma-2\kappa
g(z))\sigma_3/2}e^{iN\theta_\Gamma\sigma_3/2}\,,&\hspace{0.2
in}\mbox{for $z\in D_{\Gamma,II}^{\nab,L}$,}\\\\
\mat{H}_\Gamma^{\nab,L}(z)\hat{\mat{Z}}^{\nab,L}(\tau_\Gamma^{\nab,L}(z))
T_\nab(z)^{-\sigma_3/2}e^{(\eta(z)-\gamma-2\kappa
g(z))\sigma_3/2}e^{-iN\theta_\Gamma\sigma_3/2}\,,&\hspace{0.2
in}\mbox{for $z\in D_{\Gamma,III}^{\nab,L}$,}\\\\
\mat{H}_\Gamma^{\nab,L}(z)\hat{\mat{Z}}^{\nab,L}(\tau_\Gamma^{\nab,L}(z))
e^{(\eta(z)-\gamma-2\kappa
g(z))\sigma_3/2}e^{-iN\theta_\Gamma\sigma_3/2}\,,&\hspace{0.2
in}\mbox{for $z\in D_{\Gamma,IV}^{\nab,L}$.}
\end{array}\right.
\end{equation}
Note that in this formula, the transformation
$\tau_\Gamma^{\nab,L}(\cdot)$ and the matrix
$\mat{H}^{\nab,L}_\Gamma(z)$ will be different in neighborhoods
$D_\Gamma^{\nab,L}$ corresponding to different left band edges in
$\Sigma_0^\nab$, being defined locally by (\ref{eq:zetadefAL}),
(\ref{eq:XZdotnabL}), and (\ref{eq:HdefAL}).

\subsubsection*{Right band edge with $z_0=\beta\in\Sigma_0^\nab$ (lower constraint).}
With $\Gamma$ denoting the void to the right of the band edge $\beta$
and $I$ denoting the adjacent band on the left of $\beta$, we let
$D_\Gamma^{\nab,R}$ \label{symbol:DnabR}
be a disc centered at $z=\beta$ with radius
$h\epsilon$.  The four open quadrants of
$D_\Gamma^{\nab,R}\setminus(D_\Gamma^{\nab,R}\cap \Sigma_{\rm SD})$
are defined as \label{symbol:DnabRquadrants}
\begin{equation}
\begin{array}{rcl}
D^{\nab,R}_{\Gamma,I}&=&\displaystyle
D^{\nab,R}_\Gamma\cap\left\{z\,\Big|\,z\neq
\beta,\,0<\arg(z-\beta)<\frac{\pi}{2}\right\}\,,\\\\
D^{\nab,R}_{\Gamma,II}&=&\displaystyle
D^{\nab,R}_\Gamma\cap\left\{z\,\Big|\,z\neq
\beta,\,\frac{\pi}{2}<\arg(z-\beta)<\pi\right\}\,,\\\\
D^{\nab,R}_{\Gamma,III}&=&\displaystyle
D^{\nab,R}_\Gamma\cap\left\{z\,\Big|\,z\neq
\beta,\,-\pi<\arg(z-\beta)<-\frac{\pi}{2}\right\}\,,\\\\
D^{\nab,R}_{\Gamma,IV}&=&\displaystyle
D^{\nab,R}_\Gamma\cap\left\{z\,\Big|\,z\neq
\beta,\,-\frac{\pi}{2}<\arg(z-\beta)<0\right\}\,.
\end{array}
\end{equation}
We introduce the local change of dependent variable \label{symbol:ZnabR}
\begin{equation}
\mat{Z}_\Gamma^{\nab,R}(z):=\left\{\begin{array}{ll}
\mat{X}(z)T_\nab(z)^{\sigma_3/2}e^{(\gamma-\eta(z)+2\kappa g(z))\sigma_3/2}e^{-iN\theta_\Gamma\sigma_3/2}\,, &\hspace{0.2 in}
\mbox{for $z\in D_{\Gamma,I}^{\nab,R}$},\\\\
\mat{X}(z)e^{(\gamma-\eta(z)+2\kappa g(z))\sigma_3/2}e^{-iN\theta_\Gamma\sigma_3/2}\,,&\hspace{0.2 in}\mbox{for $z\in D_{\Gamma,II}^{\nab,R}$},\\\\
\mat{X}(z)e^{(\gamma-\eta(z)+2\kappa g(z))\sigma_3/2}e^{iN\theta_\Gamma\sigma_3/2}\,, &\hspace{0.2 in}\mbox{for $z\in
D_{\Gamma,III}^{\nab,R}$},\\\\
\mat{X}(z)T_\nab(z)^{\sigma_3/2}e^{(\gamma-\eta(z)+2\kappa g(z))\sigma_3/2}e^{iN\theta_\Gamma\sigma_3/2}\,, &\hspace{0.2 in}\mbox{for $z\in
D_{\Gamma,IV}^{\nab,R}$},
\end{array}\right.
\end{equation}
(recall that $e^{iN\theta_\Gamma}=e^{iN\theta(\beta)}$),  and 
the local conformal change of independent variable $\zeta=\tau_\Gamma^{\nab,R}(z)$
defined by \eqref{eq:zetadefAR}.
%
%introduce the transformation
%(change of independent variable)
%\begin{equation}
%\zeta=\tau_\Gamma^{\nab,R}(z):=
%\left(-2\pi Nc\int_{\beta}^zd\mu_{\rm min}^c(x)\right)^{2/3}
%\label{eq:zetadefAR}
%\end{equation}
%initially for $z<\beta$ in $D_\Gamma^{\nab,R}$ (that is, in the band
%adjacent to $\beta$). With the assumption of $p=1$ this extends to all
%of $D_\Gamma^{\nab,R}$ as an analytic change of coordinates.  
The
mapping is orientation-reversing, taking $z<\beta$ to $\zeta>0$ and
$z>\beta$ to $\zeta<0$.  By taking $\epsilon$ sufficiently small, 
the radius $h\epsilon$ of
$D_\Gamma^{\nab,R}$ will be small enough that the images under
$\tau_\Gamma^{\nab,R}$ of the segments $\arg(z-\beta)=\pm\pi/2$ in
$D_\Gamma^{\nab,R}$ lie within a symmetrical sector of the imaginary
$\zeta$-axis of opening angle strictly less than $\pi/3$.  The exact
jump conditions satisfied by $\mat{Z}_\Gamma^{\nab,R}(z)$ in
$D_\Gamma^{\nab,R}$ may be written in terms of $\zeta$ as
\begin{equation}
\begin{array}{rcll}
\mat{Z}^{\nab,R}_{\Gamma+}(z)&=&\displaystyle \mat{Z}^{\nab,R}_{\Gamma-}(z)
\left(\begin{array}{cc} 1 & ie^{-(-\zeta)^{3/2}}\\\\ 0 & 1\end{array}\right)\,,&
\hspace{0.2 in}\mbox{for $z\in\Gamma\cap D_\Gamma^{\nab,R}$,}\\\\
\mat{Z}^{\nab,R}_{\Gamma+}(z)&=&\displaystyle \mat{Z}^{\nab,R}_{\Gamma-}(z)
\left(\begin{array}{cc} 0 & -i\\\\ -i & 0\end{array}\right)\,,&
\hspace{0.2 in}\mbox{for $z\in I\cap D_\Gamma^{\nab,R}$,}\\\\
\mat{Z}^{\nab,R}_{\Gamma+}(z)&=&\displaystyle\mat{Z}^{\nab,R}_{\Gamma-}(z)
\left(\begin{array}{cc} 1 & 0\\\\ -ie^{i\zeta^{3/2}} & 1\end{array}\right)\,,&
\hspace{0.2 in}\mbox{for $z\in\Sigma_{0-}^\nab\cap D_\Gamma^{\nab,R}$,}\\\\
\mat{Z}^{\nab,R}_{\Gamma+}(z)&=&\displaystyle\mat{Z}^{\nab,R}_{\Gamma-}(z)
\left(\begin{array}{cc} 1 & 0\\\\ -ie^{-i\zeta^{3/2}} & 1\end{array}\right)\,,&
\hspace{0.2 in}\mbox{for $z\in\Sigma_{0+}^\nab\cap D_\Gamma^{\nab,R}$.}
\end{array}
\label{eq:nabRexact}
\end{equation}
The subscripts ``$+$'' and ``$-$'' indicate respectively boundary values taken from the left and right of $\Sigma_{\rm SD}$ with respect to its orientation.
The ``comparison matrix'' \label{symbol:dotZnabR}
\begin{equation}
\dot{\mat{Z}}_\Gamma^{\nab,R}(z):=\dot{\mat{X}}(z)e^{(\gamma-\eta(z)+2\kappa g(z))\sigma_3/2}e^{-iN{\rm
sgn}(\Im(z)) \theta_\Gamma\sigma_3/2}\,, \hspace{0.2 in}
\mbox{for $z\in D_\Gamma^{\nab,R}\setminus (D_\Gamma^{\nab,R}\cap \Sigma_{\rm SD})$}
\label{eq:ddotVdefAR}
\end{equation}
satisfies the same jump condition for $z\in I\cap D_\Gamma^{\nab,R}$
as $\mat{Z}_\Gamma^{\nab,R}(z)$, but is otherwise analytic in
$D_\Gamma^{\nab,R}$ and can be written in the form \label{symbol:HGammanabR}
\begin{equation}
\dot{\mat{Z}}_\Gamma^{\nab,R}(z):=
\mat{H}_\Gamma^{\nab,R}(z)\cdot \frac{1}{\sqrt{2}}
(-\tau_\Gamma^{\nab,R}(z))^{\sigma_3/4}\left(\begin{array}{cc} i &
-i\\\\-i & -i\end{array}
\right)\,,
\label{eq:HdefAR}
\end{equation}
where $\mat{H}_\Gamma^{\nab,R}(z)$ is a holomorphic factor for $z\in
D_\Gamma^{\nab,R}$.  To come up with a matrix satisfying the jump
conditions of $\mat{Z}_\Gamma^{\nab,R}(z)$ that is a good match to
$\dot{\mat{Z}}_\Gamma^{\nab,R}(z)$ on the boundary of $D_\Gamma^{\nab,R}$, we
consider the solution $\hat{\mat{Z}}^{\nab,L}(\zeta)$ of
Riemann-Hilbert Problem~\ref{rhp:Airy} with the contours $C_\pm$
chosen such that $C_\pm\cap \tau_{\Gamma}^{\nab,R}(D_\Gamma^{\nab,R})
= \tau_\Gamma^{\nab,R}(\Sigma_{0\mp}^\nab)$, and we set
\label{symbol:hatZnabR}
\begin{equation}
\hat{\mat{Z}}^{\nab,R}(\zeta):=\hat{\mat{Z}}^{\nab,L}(\zeta)\cdot i\sigma_3\,.
\label{eq:VhatRA}
\end{equation}
\begin{prop}
The matrix $\hat{\mat{Z}}^{\nab,R}(\zeta)$ defined by
(\ref{eq:VhatRA}) is an analytic function of $\zeta$ for
$\zeta\in\mathbb{C}\setminus(\mathbb{R}\cup C_+\cup C_-)$ that
satisfies the normalization condition
\begin{equation}
\hat{\mat{Z}}^{\nab,R}(\zeta)\cdot \frac{1}{\sqrt{2}}\left(\begin{array}{cc}
-i & i\\\\i & i\end{array}\right)(-\zeta)^{-\sigma_3/4} = {\mathbb I}
+ O\left(\frac{1}{\zeta}\right)\,,
\end{equation}
as $\zeta\rightarrow\infty$, uniformly with respect to direction.
Moreover, $\hat{\mat{Z}}^{\nab,R}(\zeta)$ takes continuous boundary
values from each sector of its analyticity that with
$\zeta=\tau_\Gamma^{\nab,R}(z)$ satisfy the exact same set of
relations (\ref{eq:nabRexact}) as $\mat{Z}^{\nab,R}_\Gamma(z)$.
\end{prop}
We may construct a local parametrix for $\mat{X}(z)$ in
$D_\Gamma^{\nab,R}$ as follows: \label{symbol:XhatnabR}
\begin{equation}\label{eq;XhatnabR}
\hat{\mat{X}}_\Gamma^{\nab,R}(z):=\left\{
\begin{array}{ll}
\mat{H}_\Gamma^{\nab,R}(z)\hat{\mat{Z}}^{\nab,R}(\tau_\Gamma^{\nab,R}(z))
T_\nab(z)^{-\sigma_3/2}e^{(\eta(z)-\gamma-2\kappa
g(z))\sigma_3/2}e^{iN\theta_\Gamma\sigma_3/2}\,,&\hspace{0.2
in}\mbox{for $z\in D_{\Gamma,I}^{\nab,R}$,}\\\\
\mat{H}_\Gamma^{\nab,R}(z)\hat{\mat{Z}}^{\nab,R}(\tau_\Gamma^{\nab,R}(z))
e^{(\eta(z)-\gamma-2\kappa
g(z))\sigma_3/2}e^{iN\theta_\Gamma\sigma_3/2}\,,&\hspace{0.2
in}\mbox{for $z\in D_{\Gamma,II}^{\nab,R}$,}\\\\
\mat{H}_\Gamma^{\nab,R}(z)\hat{\mat{Z}}^{\nab,R}(\tau_\Gamma^{\nab,R}(z))
e^{(\eta(z)-\gamma-2\kappa
g(z))\sigma_3/2}e^{-iN\theta_\Gamma\sigma_3/2}\,,&\hspace{0.2
in}\mbox{for $z\in D_{\Gamma,III}^{\nab,R}$,}\\\\
\mat{H}_\Gamma^{\nab,R}(z)\hat{\mat{Z}}^{\nab,R}(\tau_\Gamma^{\nab,R}(z))
T_\nab(z)^{-\sigma_3/2}e^{(\eta(z)-\gamma-2\kappa
g(z))\sigma_3/2}e^{-iN\theta_\Gamma\sigma_3/2}\,,&\hspace{0.2
in}\mbox{for $z\in D_{\Gamma,IV}^{\nab,R}$.}\end{array}\right.
\end{equation}
Again, the transformation $\tau_\Gamma^{\nab,R}(\cdot)$ and the matrix
$\mat{H}_\Gamma^{\nab,R}(z)$ will be different in neighborhoods
$D_\Gamma^{\nab,R}$ corresponding to different right band edges in
$\Sigma_0^\nab$.

\subsubsection*{Left band edge with $z_0=\alpha\in\Sigma_0^\del$ (upper constraint).}
Letting $\Gamma$ denote the saturated region to the left of $\alpha$,
$I$ denote the band to the right, and $D_\Gamma^{\del,L}$ 
\label{symbol:DdelL} denote a
disc centered at $z=\alpha$ with radius $h\epsilon$, we partition
the disc into quadrants: \label{symbol:DdelLquadrants}
\begin{equation}
\begin{array}{rcl}
D^{\del,L}_{\Gamma,I}&=&\displaystyle
D^{\del,L}_\Gamma\cap\left\{z\,\Big|\,z\neq
\alpha,\,0<\arg(z-\alpha)<\frac{\pi}{2}\right\}\,,\\\\
D^{\del,L}_{\Gamma,II}&=&\displaystyle
D^{\del,L}_\Gamma\cap\left\{z\,\Big|\,z\neq
\alpha,\,\frac{\pi}{2}<\arg(z-\alpha)<\pi\right\}\,,\\\\
D^{\del,L}_{\Gamma,III}&=&\displaystyle
D^{\del,L}_\Gamma\cap\left\{z\,\Big|\,z\neq
\alpha,\,-\pi<\arg(z-\alpha)<-\frac{\pi}{2}\right\}\,,\\\\
D^{\del,L}_{\Gamma,IV}&=&\displaystyle
D^{\del,L}_\Gamma\cap\left\{z\,\Big|\,z\neq
\alpha,\,-\frac{\pi}{2}<\arg(z-\alpha)<0\right\}\,.
\end{array}
\end{equation}
Next we set \label{symbol:ZdelL}
\begin{equation}
\mat{Z}_\Gamma^{\del,L}(z):=\left\{\begin{array}{ll}
\mat{X}(z)e^{(\gamma-\eta(z)+2\kappa g(z))\sigma_3/2}e^{-iN\theta_\Gamma\sigma_3/2}\,,&\hspace{0.2 in}\mbox{for $z\in D_{\Gamma,I}^{\del,L}$,}\\\\
\mat{X}(z)T_\del(z)^{-\sigma_3/2}e^{(\gamma-\eta(z)+2\kappa g(z))\sigma_3/2}e^{-iN\theta_\Gamma\sigma_3/2}\,,&\hspace{0.2 in}\mbox{for $z\in D_{\Gamma,II}^{\del,L}$,}\\\\
\mat{X}(z)T_\del(z)^{-\sigma_3/2}e^{(\gamma-\eta(z)+2\kappa g(z))\sigma_3/2}e^{iN\theta_\Gamma\sigma_3/2}\,, &\hspace{0.2 in}\mbox{for $z\in
D_{\Gamma,III}^{\del,L}$,}\\\\
\mat{X}(z)e^{(\gamma-\eta(z)+2\kappa g(z))\sigma_3/2}e^{iN\theta_\Gamma\sigma_3/2}\,,&\hspace{0.2 in}\mbox{for $z\in D_{\Gamma,IV}^{\del,L}$,}
\end{array}\right.
\end{equation}
where we recall that $e^{iN\theta_\Gamma}=e^{iN\theta(\alpha)}$, and consider the conformal mapping $\zeta=\tau_\Gamma^{\del,L}(z)$ defined
by \eqref{eq:zetadefBL}.
%\begin{equation}
%\zeta=\tau_\Gamma^{\del,L}(z):=\left(2\pi Nc\int_{\alpha}^z %\left[\frac{1}{c}\rho^0(x)-
%\frac{d\mu^c_{\rm min}}{dx}(x)\right]\,dx\right)^{2/3}
%\label{eq:zetadefBL}
%\end{equation}
%at first for $z>\alpha$ in $D_\Gamma^{\del,L}$.  Since the equilibrium %measure
%does not achieve the upper constraint in $I\cap D_\Gamma^{\del,L}$,
%such values of $z$ are mapped to positive values of $\zeta$.  Assuming
%$p=1$, we obtain by analytic continuation an orientation-preserving
%holomorphic change of coordinates $z\leftrightarrow\zeta$ defined for
%$z\in D_\Gamma^{\del,L}$.  
We choose the parameter $\epsilon$ controlling the radius of $D_\Gamma^{\del,L}$
to be sufficiently small that the images
$\tau_\Gamma^{\del,L}(\Sigma_{0\pm}^\del\cap D_\Gamma^{\del,L})$ lie
within a symmetrical sector of the imaginary $\zeta$-axis of opening
angle strictly less than $\pi/3$.  The exact jump conditions satisfied
by the matrix $\mat{Z}_\Gamma^{\del,L}(z)$ may be written in terms of
$\zeta$ in a simple way:
\begin{equation}
\begin{array}{rcll}
\mat{Z}^{\del,L}_{\Gamma+}(z)&=&\displaystyle\mat{Z}^{\del,L}_{\Gamma-}(z)
\left(\begin{array}{cc} 1 & 0 \\\\ ie^{-(-\zeta)^{3/2}} & 1\end{array}\right)\,,&
\hspace{0.2 in}\mbox{for $z\in\Gamma\cap D_\Gamma^{\del,L}$,}\\\\
\mat{Z}^{\del,L}_{\Gamma+}(z)&=&\displaystyle\mat{Z}^{\del,L}_{\Gamma-}(z)
\left(\begin{array}{cc} 0 & -i \\ \\-i & 0\end{array}\right)\,,&
\hspace{0.2 in}\mbox{for $z\in I\cap D_\Gamma^{\del,L}$,}\\\\
\mat{Z}^{\del,L}_{\Gamma+}(z)&=&\displaystyle\mat{Z}^{\del,L}_{\Gamma-}(z)
\left(\begin{array}{cc} 1 & -ie^{i\zeta^{3/2}} \\\\ 0 & 1\end{array}\right)\,,&
\hspace{0.2 in}\mbox{for $z\in\Sigma_{0+}^\del\cap D_\Gamma^{\del,L}$,}\\\\
\mat{Z}^{\del,L}_{\Gamma+}(z)&=&\displaystyle\mat{Z}^{\del,L}_{\Gamma-}(z)
\left(\begin{array}{cc} 1 & -ie^{-i\zeta^{3/2}} \\\\ 0 & 1\end{array}\right)\,,&
\hspace{0.2 in}\mbox{for $z\in\Sigma_{0-}^\del\cap D_\Gamma^{\del,L}$.}
\end{array}
\label{eq:delLexact}
\end{equation}
The subscripts ``$+$'' and ``$-$'' refer respectively to boundary values taken on the oriented contour $\Sigma_{\rm SD}$ from the left and right.  The ``comparison matrix'' defined by the formula \label{symbol:dotZdelL}
\begin{equation}
\dot{\mat{Z}}_\Gamma^{\del,L}(z):=\dot{\mat{X}}(z)e^{(\gamma-\eta(z)+2\kappa g(z))\sigma_3/2}e^{-iN{\rm sgn}(\Im(z))
\theta_\Gamma\sigma_3/2}
\label{eq:ddotVdefBL}
\end{equation}
satisfies the same jump condition for $z\in I\cap D_\Gamma^{\del,L}$
as $\mat{Z}_\Gamma^{\del,L}(z)$ but is otherwise analytic in
$D_\Gamma^{\del,L}$, and thus may be written in the form
\label{symbol:HGammadelL}
\begin{equation}
\dot{\mat{Z}}_\Gamma^{\del,L}(z)=
\mat{H}_\Gamma^{\del,L}(z)\cdot\frac{1}{\sqrt{2}}
(-\tau_\Gamma^{\del,L}(z))^{\sigma_3/4}\left(\begin{array}{cc} i & i
\\\\ i & -i
\end{array}\right)\,.
\label{eq:HdefBL}
\end{equation}
The quotient matrix $\mat{H}_\Gamma^{\del,L}(z)$ is holomorphic in
$D_\Gamma^{\del,L}$.  Finding a matrix with the same jump conditions
as $\mat{Z}_\Gamma^{\del,L}(z)$ and matching onto
$\dot{\mat{Z}}_\Gamma^{\del,L}(z)$ at the boundary of
$D_\Gamma^{\del,L}$ leads us to recall the matrix
$\hat{\mat{Z}}^{\nab,L}(\zeta)$ solving Riemann-Hilbert
Problem~\ref{rhp:Airy} with the contours $C_\pm$ taken to be such that
for each $N$, $C_\pm\cap
\tau_\Gamma^{\del,L}(D_\Gamma^{\del,L})=\tau_\Gamma^{\del,L}(\Sigma_{0\pm}^\del)$
and to set \label{symbol:hatZdelL}
\begin{equation}
\hat{\mat{Z}}^{\del,L}(\zeta):=\hat{\mat{Z}}^{\nab,L}(\zeta)\cdot i\sigma_1\,.
\label{eq:VhatLB}
\end{equation}
\begin{prop}
The matrix $\hat{\mat{Z}}^{\del,L}(\zeta)$ defined by
(\ref{eq:VhatLB}) is an analytic function of $\zeta$ for
$\zeta\in\mathbb{C}\setminus (\mathbb{R}\cup C_+\cup C_-)$ that
satisfies the normalization condition
\begin{equation}
\hat{\mat{Z}}^{\del,L}(\zeta)\cdot \frac{1}{\sqrt{2}}\left(\begin{array}{cc}
-i & -i\\\\-i & i\end{array}\right)(-\zeta)^{-\sigma_3/4} = {\mathbb
I} + O\left(\frac{1}{\zeta}\right)\,,
\end{equation}
as $\zeta\rightarrow\infty$, uniformly with respect to direction.
Moreover, $\hat{\mat{Z}}^{\del,L}(\zeta)$ takes continuous boundary
values from each sector of its analyticity that with
$\zeta=\tau_\Gamma^{\del,L}(z)$ satisfy the exact same set of
relations (\ref{eq:delLexact}) as $\mat{Z}_\Gamma^{\del,L}(z)$.
\end{prop}

We construct a local parametrix for $\mat{X}(z)$ with the formula
\label{symbol:XhatdelL}
\begin{equation}\label{eq;XhatdelL}
\hat{\mat{X}}_\Gamma^{\del,L}(z):=\left\{\begin{array}{ll}
\mat{H}^{\del,L}_\Gamma(z)\hat{Z}^{\del,L}(\tau_\Gamma^{\del,L}(z))
e^{(\eta(z)-\gamma-2\kappa
g(z))\sigma_3/2}e^{iN\theta_\Gamma\sigma_3/2}\,, &\hspace{0.2
in}\mbox{for $z\in D_{\Gamma,I}^{\del,L}$,}\\\\
\mat{H}^{\del,L}_\Gamma(z)\hat{Z}^{\del,L}(\tau_\Gamma^{\del,L}(z))
T_\del(z)^{\sigma_3/2} e^{(\eta(z)-\gamma-2\kappa
g(z))\sigma_3/2}e^{iN\theta_\Gamma\sigma_3/2}\,, &\hspace{0.2
in}\mbox{for $z\in D_{\Gamma,II}^{\del,L}$,}\\\\
\mat{H}^{\del,L}_\Gamma(z)\hat{Z}^{\del,L}(\tau_\Gamma^{\del,L}(z))
T_\del(z)^{\sigma_3/2} e^{(\eta(z)-\gamma-2\kappa
g(z))\sigma_3/2}e^{-iN\theta_\Gamma\sigma_3/2}\,, &\hspace{0.2
in}\mbox{for $z\in D_{\Gamma,III}^{\del,L}$,}\\\\
\mat{H}^{\del,L}_\Gamma(z)\hat{Z}^{\del,L}(\tau_\Gamma^{\del,L}(z))
e^{(\eta(z)-\gamma-2\kappa
g(z))\sigma_3/2}e^{-iN\theta_\Gamma\sigma_3/2}\,, &\hspace{0.2
in}\mbox{for $z\in D_{\Gamma,IV}^{\del,L}$.}
\end{array}\right.
\end{equation}
As before, the transformation $\tau_\Gamma^{\del,L}(\cdot)$ and the
matrix $\mat{H}_\Gamma^{\del,L}(z)$ will be different in different
neighborhoods $D_\Gamma^{\del,L}$ corresponding to different left band
edges in $\Sigma_0^\del$.

\subsubsection*{\bf Right band edge with $z_0=\beta\in\Sigma_0^\del$ (upper constraint).}
With $\Gamma$ denoting the saturated region to the right of $\beta$
and $I$ denoting the band to the left, we work in a disc
$D_\Gamma^{\del,R}$ \label{symbol:DdelR}
centered at $z=\beta$ with radius $h\epsilon$,
and partition the disc into quadrants: \label{symbol:DdelRquadrants}
\begin{equation}
\begin{array}{rcl}
D^{\del,R}_{\Gamma,I}&=&\displaystyle
D^{\del,R}_\Gamma\cap\left\{z\,\Big|\,z\neq
\beta,\,0<\arg(z-\beta)<\frac{\pi}{2}\right\}\,,\\\\
D^{\del,R}_{\Gamma,II}&=&\displaystyle
D^{\del,R}_\Gamma\cap\left\{z\,\Big|\,z\neq
\beta,\,\frac{\pi}{2}<\arg(z-\beta)<\pi\right\}\,,\\\\
D^{\del,R}_{\Gamma,III}&=&\displaystyle
D^{\del,R}_\Gamma\cap\left\{z\,\Big|\,z\neq
\beta,\,-\pi<\arg(z-\beta)<-\frac{\pi}{2}\right\}\,,\\\\
D^{\del,R}_{\Gamma,IV}&=&\displaystyle
D^{\del,R}_\Gamma\cap\left\{z\,\Big|\,z\neq
\beta,\,-\frac{\pi}{2}<\arg(z-\beta)<0\right\}\,.
\end{array}
\end{equation}
We then set \label{symbol:ZdelR}
\begin{equation}
\mat{Z}_\Gamma^{\del,R}(z):=\left\{\begin{array}{ll}
\mat{X}(z)T_\del(z)^{-\sigma_3/2}
e^{(\gamma-\eta(z)+2\kappa
g(z))\sigma_3/2}e^{-iN\theta_\Gamma\sigma_3/2}
\,,&\hspace{0.2 in}\mbox{for $z\in D_{\Gamma,I}^{\del,R}$,}\\\\
\mat{X}(z)
e^{(\gamma-\eta(z)+2\kappa
g(z))\sigma_3/2}e^{-iN\theta_\Gamma\sigma_3/2}
\,,&\hspace{0.2 in}\mbox{for $z\in D_{\Gamma,II}^{\del,R}$,}\\\\
\mat{X}(z)
e^{(\gamma-\eta(z)+2\kappa
g(z))\sigma_3/2}e^{iN\theta_\Gamma\sigma_3/2}
\,,&\hspace{0.2 in}\mbox{for $z\in D_{\Gamma,III}^{\del,R}$,}\\\\
\mat{X}(z)T_\del(z)^{-\sigma_3/2}
e^{(\gamma-\eta(z)+2\kappa
g(z))\sigma_3/2}e^{iN\theta_\Gamma\sigma_3/2}
\,,&\hspace{0.2 in}\mbox{for $z\in D_{\Gamma,IV}^{\del,R}$,}
\end{array}\right.
\end{equation}
where we recall that $e^{iN\theta_\Gamma}=e^{iN\theta(\beta)}$, and consider the conformal mapping $\zeta=\tau_\Gamma^{\del,R}(z)$ defined by \eqref{eq:zetadefBR}.
%relation
%\begin{equation}
%\zeta=\tau_\Gamma^{\del,R}(z):=
%\left(-2\pi N c\int_{\beta}^z\left[\frac{1}{c}\rho^0(x)-
%\frac{d\mu_{\rm min}^c}{dx}(x)\right]\,dx\right)^{2/3}
%\label{eq:zetadefBR}
%\end{equation}
%for $z\in D_\Gamma^{\del,R}$ and $z<\beta$.  For such $z$, we have
%$\zeta>0$ since the upper constraint is not achieved in $I\cap
%D_\Gamma^{\del,R}$.  When $p=1$ this extends to an invertible analytic
This is an orientation-reversing transformation of the neighborhood
$D_\Gamma^{\del,R}$ of $z=\beta$ in the $z$-plane to a neighborhood of
the origin in the $\zeta$-plane.  By making $\epsilon$ small enough, the radius of $D_\Gamma^{\del,R}$ will be so
small that the images
$\tau_\Gamma^{\del,R}(\Sigma_{0\pm}^\del\cap D_\Gamma^{\del,R})$ lie
within a symmetrical sector of the imaginary $\zeta$-axis of opening
angle strictly less than $\pi/3$.  The matrix
$\mat{Z}_\Gamma^{\del,R}(z)$ then satisfies exactly the following jump
conditions:
\begin{equation}
\begin{array}{rcll}
\mat{Z}_{\Gamma+}^{\del,R}(z)&=&\displaystyle\mat{Z}_{\Gamma-}^{\del,R}(z)
\left(\begin{array}{cc}
1 & 0 \\\\ ie^{-(-\zeta)^{3/2}} & 1
\end{array}\right)\,,&\hspace{0.2 in}\mbox{for $z\in\Gamma\cap D_\Gamma^{\del,R}$,}\\\\
\mat{Z}_{\Gamma+}^{\del,R}(z)&=&\displaystyle\mat{Z}_{\Gamma-}^{\del,R}(z)
\left(\begin{array}{cc}
0 & -i \\\\ -i & 0
\end{array}\right)\,,&\hspace{0.2 in}\mbox{for $z\in I\cap D_\Gamma^{\del,R}$,}\\\\
\mat{Z}_{\Gamma+}^{\del,R}(z)&=&\displaystyle\mat{Z}_{\Gamma-}^{\del,R}(z)
\left(\begin{array}{cc}
1 & -ie^{i\zeta^{3/2}} \\ \\0 & 1
\end{array}\right)\,,&\hspace{0.2 in}\mbox{for $z\in\Sigma_{0-}^\del\cap D_\Gamma^{\del,R}$,}\\\\
\mat{Z}_{\Gamma+}^{\del,R}(z)&=&\displaystyle\mat{Z}_{\Gamma-}^{\del,R}(z)
\left(\begin{array}{cc}
1 & -ie^{-i\zeta^{3/2}} \\\\ 0 & 1
\end{array}\right)\,,&\hspace{0.2 in}\mbox{for $z\in\Sigma_{0+}^\del\cap D_\Gamma^{\del,R}$.}
\end{array}
\label{eq:delRexact}
\end{equation}
The subscripts ``$+$'' and ``$-$'' respectively indicate boundary values taken on $\Sigma_{\rm SD}$ from the left and right.
The ``comparison matrix'' \label{symbol:dotZdelR}
\begin{equation}
\dot{\mat{Z}}_\Gamma^{\del,R}(z):=\dot{\mat{X}}(z)e^{(\gamma-\eta(z)+2\kappa g(z))\sigma_3/2}e^{-iN{\rm
sgn}(\Im(z))\theta_\Gamma
\sigma_3/2}
\label{eq:ddotVdefBR}
\end{equation}
satisfies the same jump condition for $z\in I\cap D_\Gamma^{\del,R}$
as does $\mat{Z}_\Gamma^{\del,R}(z)$ and is otherwise analytic in
$D_\Gamma^{\del,R}$; it may be written in the form
\label{symbol:HGammadelR}
\begin{equation}
\dot{\mat{Z}}_\Gamma^{\del,R}(z)=
\mat{H}_\Gamma^{\del,R}(z)\cdot\frac{1}{\sqrt{2}}
(-\tau_\Gamma^{\del,R}(z))^{\sigma_3/4}\left(\begin{array}{cc} -1 & -1 \\ \\1 & -1
\end{array}\right)\,.
\label{eq:HdefBR}
\end{equation}
The quotient $\mat{H}^{\del,R}_\Gamma(z)$ is holomorphic in
$D_\Gamma^{\del,R}$.  A matrix that satisfies the same jump conditions
as $\mat{Z}_\Gamma^{\del,R}(z)$ and matches well onto
$\dot{\mat{Z}}_\Gamma^{\del,R}(z)$ may be obtained by considering the
matrix $\hat{\mat{Z}}^{\nab,L}(\zeta)$ satisfying Riemann-Hilbert
Problem~\ref{rhp:Airy} with the contours $C_\pm$ chosen so that
$C_\pm\cap\tau_\Gamma^{\del,R}(D_\gamma^{\del,R})=\tau_\Gamma^{\del,R}(\Sigma_{0\mp}^\del)$,
and set \label{symbol:hatZdelR}
\begin{equation}
\hat{\mat{Z}}^{\del,R}(\zeta):=
\hat{\mat{Z}}^{\nab,L}(\zeta)\cdot \sigma_1\sigma_3\,.
\label{eq:VhatRB}
\end{equation}

\begin{prop}
The matrix $\hat{\mat{Z}}^{\del,R}(\zeta)$ defined by
(\ref{eq:VhatRB}) is an analytic function of $\zeta$ for
$\zeta\in\mathbb{C}\setminus(\mathbb{R}\cup C_+\cup C_-)$ that
satisfies the normalization condition
\begin{equation}
\hat{\mat{Z}}^{\del,R}(\zeta)\cdot \frac{1}{\sqrt{2}}\left(\begin{array}{cc}
1 & 1\\\\-1 & 1\end{array}\right)(-\zeta)^{-\sigma_3/4} = {\mathbb I}
+ O\left(\frac{1}{\zeta}\right)\,,
\end{equation}
as $\zeta\rightarrow\infty$, uniformly with respect to direction.
Moreover, $\hat{\mat{Z}}^{\del,R}(\zeta)$ takes continuous boundary
values from each sector of its analyticity that with
$\zeta=\tau_\Gamma^{\del,R}(z)$ satisfy the exact same set of
relations (\ref{eq:delRexact}) as $\mat{Z}_\Gamma^{\del,R}(z)$.
\end{prop}

We use $\hat{\mat{Z}}^{\del,R}(\zeta)$ to construct a local parametrix
for $\mat{X}(z)$ in $D_\Gamma^{\del,R}$ by the scheme:
\label{symbol:XhatdelR}
\begin{equation}\label{eq;XhatdelR}
\hat{\mat{X}}_\Gamma^{\del,R}(z):=\left\{\begin{array}{ll}
\mat{H}_\Gamma^{\del,R}(z)\hat{\mat{Z}}^{\del,R}(\tau_\Gamma^{\del,R}(z))
T_\del(z)^{\sigma_3/2}e^{(\eta(z)-\gamma-2\kappa
g(z))\sigma_3/2}e^{iN\theta_\Gamma\sigma_3/2}\,,&\hspace{0.2 in}
\mbox{for $z\in D_{\Gamma,I}^{\del,R}$,}\\\\
\mat{H}_\Gamma^{\del,R}(z)\hat{\mat{Z}}^{\del,R}(\tau_\Gamma^{\del,R}(z))
e^{(\eta(z)-\gamma-2\kappa
g(z))\sigma_3/2}e^{iN\theta_\Gamma\sigma_3/2}\,,&\hspace{0.2 in}
\mbox{for $z\in D_{\Gamma,II}^{\del,R}$,}\\\\
\mat{H}_\Gamma^{\del,R}(z)\hat{\mat{Z}}^{\del,R}(\tau_\Gamma^{\del,R}(z))
e^{(\eta(z)-\gamma-2\kappa
g(z))\sigma_3/2}e^{-iN\theta_\Gamma\sigma_3/2}\,,&\hspace{0.2 in}
\mbox{for $z\in D_{\Gamma,III}^{\del,R}$,}\\\\
\mat{H}_\Gamma^{\del,R}(z)\hat{\mat{Z}}^{\del,R}(\tau_\Gamma^{\del,R}(z))
T_\del(z)^{\sigma_3/2}e^{(\eta(z)-\gamma-2\kappa
g(z))\sigma_3/2}e^{-iN\theta_\Gamma\sigma_3/2}\,,&\hspace{0.2 in}
\mbox{for $z\in D_{\Gamma,IV}^{\del,R}$.}
\end{array}\right.
\end{equation}
Once again, the transformation $\tau_\Gamma^{\del,R}(\cdot)$ and the
matrix $\mat{H}_\Gamma^{\del,R}(z)$ will be different in different
neighborhoods $D_\Gamma^{\del,R}$ corresponding to different right
band edges in $\Sigma_0^\del$.

\subsubsection*{Common properties of the four local approximations.}
The important properties of the local approximations are summarized in
the following proposition.
\begin{prop}
Although originally defined in the four open quadrants within each
disc, each function $\hat{\mat{X}}_\Gamma^{\nab,L}(z)\mat{X}(z)^{-1}$,
$\hat{\mat{X}}_\Gamma^{\nab,R}(z)\mat{X}(z)^{-1}$,
$\hat{\mat{X}}_\Gamma^{\del,L}(z)\mat{X}(z)^{-1}$, and
$\hat{\mat{X}}_\Gamma^{\del,R}(z)\mat{X}(z)^{-1}$ has a continuous and
hence analytic extension to the full interior of the corresponding
disc.  For each sufficiently small $\epsilon>0$ there is a constant
$M_\epsilon>0$ such that on the boundary of each disc centered at a
band edge $z=z_0$ we have
\begin{equation}
\sup_{|z-z_0|=h\epsilon}\|\hat{\mat{X}}_\Gamma^{*,*}(z)\dot{\mat{X}}(z)^{-1}
-{\mathbb I}\|\le \frac{M_\epsilon}{N}
\label{eq:delDmatch}
\end{equation}
for sufficiently large $N$.  Here $\hat{\mat{X}}_\Gamma^{*,*}(z)$
refers to any of the four different types of local parametrix.
\label{prop:Dproperties}
\end{prop}

\begin{proof}
The analyticity of $\hat{\mat{X}}_\Gamma^{*,*}(z)\mat{X}(z)^{-1}$
throughout $D_\Gamma^{*,*}$ follows directly from the construction in
each case, in that there is no approximation of the jump matrix.

%To obtain the estimate \eqref{eq:delDmatch}, we
%begin by noting some properties of the holomorphic prefactors
%$\mat{H}_\Gamma^{*,*}(z)$.  The matrices $\mat{H}_\Gamma^{*,*}(z)$ all
%have determinant one. The (analytic) matrix elements of
%$\mat{H}_\Gamma^{\nab,L}(z)$ and $\mat{H}_\Gamma^{\del,R}(z)$ are real
%when $z$ is real, and those of $\mat{H}_\Gamma^{\nab,R}(z)$ and
%$\mat{H}_\Gamma^{\del,L}(z)$ are imaginary when $z$ is real.  The
%matrix elements of the product $\mat{H}_\Gamma^{*,*}(z)N^{\sigma_3/6}$
%are bounded in $D_\Gamma^{*,*}$ uniformly with respect to $N$.

%These facts follow from the formula for $\mat{H}_\Gamma^{*,*}(z)$:
%One then uses Proposition~\ref{prop:dotXbound} to calculate the
%determinant, and the reality conditions on the real axis follow, upon
%taking a boundary value from the upper half-plane on the real axis,
%from Proposition~\ref{prop:dotXsymmetry}.  The bounds on the matrix
%elements follow from Proposition~\ref{prop:dotXbound} once again and
%the expressions (\ref{eq:zetadefAL}), (\ref{eq:zetadefAR}),
%(\ref{eq:zetadefBL}), and (\ref{eq:zetadefBR}) for the mappings
%$\tau_\Gamma^{*,*}(z)$.

To prove \eqref{eq:delDmatch}, first note that since each band edge
point $z_0$ is bounded away from all transition points $y_{k,N}\in
Y_N$ and from the endpoints $\{a,b\}$, Proposition~\ref{prop:tasymp}
guarantees that for $|z-z_0|\le h\epsilon$,
\begin{equation}
\hat{\mat{X}}_\Gamma^{*,*}(z)=\mat{H}^{*,*}_\Gamma(z)
\hat{\mat{Z}}^{*,*}(\tau_\Gamma^{*,*}(z))e^{(\eta(z)-\gamma-2\kappa g(z))\sigma_3/2}e^{iN{\rm sgn}(\Im(z))\theta_\Gamma\sigma_3/2}
\left(\mathbb{I}+\mat{G}(z)\right)\,,
\end{equation}
where for some constant $J_\epsilon>0$,
\begin{equation}
\sup_{|z-z_0|<h\epsilon}\|\mat{G}(z)\|\le\frac{J_\epsilon}{N}\,.
\end{equation}
Since according to Proposition~\ref{prop:dotXbound},
$\dot{\mat{X}}(z)$ is uniformly bounded for $|z-z_0|=h\epsilon$ and
has determinant one, it follows that a related constant
$\tilde{J}_\epsilon>0$ exists such that a similar estimate holds:
\begin{equation}
\sup_{|z-z_0|=h\epsilon}\|\dot{\mat{X}}(z)\mat{G}(z)\dot{\mat{X}}(z)^{-1}\|
\le\frac{\tilde{J}_\epsilon}{N}
\end{equation}
for all sufficiently large $N$.  Next, we recall the formula for the
holomorphic prefactors $\mat{H}_\Gamma^{*,*}(z)$:
\begin{equation}
\mat{H}_\Gamma^{*,*}(z)=\dot{\mat{X}}(z)e^{(\gamma-\eta(z)+2\kappa g(z))\sigma_3/2}e^{-iN{\rm sgn}(\Im(z))
\theta_{\Gamma}\sigma_3/2}\mat{C}^{*,*}
[-\tau_\Gamma^{*,*}(z)]^{-\sigma_3/4}\,,
\label{eq:Hdefine}
\end{equation}
where the constant matrices $\mat{C}^{*,*}$ are given by
\begin{equation}
\begin{array}{rclrcl}
\mat{C}^{\nab,L}&:=&\displaystyle
\frac{1}{\sqrt{2}}\left(\begin{array}{cc}1&-1\\\\1&1\end{array}\right)\,, &
\mat{C}^{\nab,R}&:=&\displaystyle
\frac{1}{\sqrt{2}}\left(\begin{array}{cc}-i & i\\\\i & i\end{array}\right)\,,
\\\\
\mat{C}^{\del,L}&:=&\displaystyle
\frac{1}{\sqrt{2}}\left(\begin{array}{cc}-i & -i\\\\-i & i\end{array}\right)\,,
&\mat{C}^{\del,R}&:=&\displaystyle
\frac{1}{\sqrt{2}}\left(\begin{array}{cc} -1 & 1\\\\-1 & -1\end{array}\right)\,.
\end{array}
\end{equation}
Thus, we have
\begin{equation}
\hat{\mat{X}}_\Gamma^{*,*}(z)\dot{\mat{X}}(z)^{-1}=\mat{H}_\Gamma^{*,*}(z)
\hat{\mat{Z}}^{*,*}(\tau_\Gamma^{*,*}(z))\mat{C}^{*,*}[-\tau_\Gamma^{*,*}(z)]^{-\sigma_3/4}\mat{H}^{*,*}_\Gamma(z)^{-1}\left(\mathbb{I}+\dot{\mat{X}}(z)
\mat{G}(z)\dot{\mat{X}}(z)^{-1}\right)\,.
\end{equation}
Using (\ref{eq:Hdefine}) again this can be written as
\begin{equation}
\begin{array}{rcl}
\hat{\mat{X}}_\Gamma^{*,*}(z)
\dot{\mat{X}}(z)^{-1} &=&
\mat{W}^{*,*}(z)\\\\
&&\displaystyle \,\,\,\cdot\,\,\,
[-\tau_\Gamma^{*,*}(z)]^{-\sigma_3/4}\cdot
\hat{\mat{Z}}^{*,*}[\tau_\Gamma^{*,*}(z)]\mat{C}^{*,*}
[-\tau_\Gamma^{*,*}(z)]^{-\sigma_3/4}\cdot
[-\tau_\Gamma^{*,*}(z)]^{\sigma_3/4}\\\\ &&\displaystyle
\,\,\,\cdot\,\,\,
\mat{W}^{*,*}(z)^{-1}\left(\mathbb{I}+\dot{\mat{X}}(z)\mat{G}(z)
\dot{\mat{X}}(z)^{-1}\right)
\end{array}
\end{equation}
where
\begin{equation}
\mat{W}^{*,*}(z):=\dot{\mat{X}}(z)e^{(\gamma-\eta(z)+2\kappa g(z))\sigma_3/2}e^{-iN{\rm sgn}(\Im(z))\theta_\Gamma\sigma_3/2}
\mat{C}^{*,*}
\end{equation}
is a matrix that is, according to Proposition~\ref{prop:dotXbound},
uniformly bounded when $|z-z_0|=h\epsilon$.  But, from
\eqref{eq:Airydecay}, we get that for $|z-z_0|=h\epsilon$, which
corresponds to $\tau_\Gamma^{*,*}(z)$ of size $N^{2/3}$,
\begin{equation}
\begin{array}{l}
\displaystyle
(-\tau_\Gamma^{*,*}(z))^{-\sigma_3/4}\hat{\mat{Z}}^{*,*}(\tau_\Gamma^{*,*}(z))
\mat{C}^{*,*}(-\tau_\Gamma^{*,*}(z))^{-\sigma_3/4}
(-\tau_\Gamma^{*,*}(z))^{\sigma_3/4} \\\\
\displaystyle\hspace{0.4 in}\begin{array}{cl}
=&\displaystyle (-\tau_\Gamma^{*,*}(z))^{-\sigma_3/4}\left[{\mathbb
I}+
\left(\begin{array}{cc} O(\tau_\Gamma^{*,*}(z)^{-3/2}) &
O(\tau_\Gamma^{*,*}(z)^{-1}) \\ O(\tau_\Gamma^{*,*}(z)^{-2}) &
O(\tau_\Gamma^{*,*}(z)^{-3/2})
\end{array}\right)\right](-\tau_\Gamma^{*,*}(z))^{\sigma_3/4} \\\\
=&{\mathbb I}+O(\tau_\Gamma^{*,*}(z)^{-3/2})
\end{array}
\end{array}
\end{equation}
which is of order $1/N$ as desired when $|z-z_0|=h\epsilon$.  This
establishes (\ref{eq:delDmatch}).
\end{proof}

\subsubsection{Definition of the parametrix $\hat{\mat{X}}(z)$.}
The parametrix $\hat{\mat{X}}(z)$ \label{symbol:hatX}
is an explicit, global approximation
of $\mat{X}(z)$ the validity of which we will establish in
\S~\ref{sec:error}.  It is defined for
$z\in\mathbb{C}\setminus(\Sigma_{\rm SD}\cup \{\mbox{disc
boundaries}\})$ as follows.  About each left band edge $\alpha\in
(a,b)$ where the lower constraint becomes active in a void $\Gamma$ we
have placed a disc $D_\Gamma^{\nab,L}$.
For $z\in D_\Gamma^{\nab,L}\cap (\mathbb{C}\setminus\Sigma_{\rm
SD})$ we set
\begin{equation}\label{eq:XfromXdotnabL}
\hat{\mat{X}}(z):=\hat{\mat{X}}_\Gamma^{\nab,L}(z)\,.
\end{equation}
About each right band edge $\beta\in (a,b)$ where the lower constraint
becomes active in a void $\Gamma$ we have placed a disc
$D_\Gamma^{\nab,R}$. For $z\in
D_\Gamma^{\nab,R}\cap(\mathbb{C}\setminus\Sigma_{\rm SD})$ we set
\begin{equation}\label{eq:XfromXdotnabR}
\hat{\mat{X}}(z):=\hat{\mat{X}}_\Gamma^{\nab,R}(z)\,.
\end{equation}
About each left band edge $\alpha\in (a,b)$ where the upper constraint
becomes active in a saturated region $\Gamma$ we have placed a disc
$D_\Gamma^{\del,L}$. For $z\in D_\Gamma^{\del,L}\cap
(\mathbb{C}\setminus\Sigma_{\rm SD})$ we set
\begin{equation}\label{eq:XfromXdotdelL}
\hat{\mat{X}}(z):=\hat{\mat{X}}_\Gamma^{\del,L}(z)\,.
\end{equation}
About each right band edge $\beta\in (a,b)$ where the upper constraint
becomes active in a saturated region $\Gamma$ we have placed a disc
$D_\Gamma^{\del,R}$. For $z\in
D_\Gamma^{\del,R}\cap(\mathbb{C}\setminus\Sigma_{\rm SD})$ we set
\begin{equation}\label{eq:XfromXdotdelR}
\hat{\mat{X}}(z):=\hat{\mat{X}}_\Gamma^{\del,R}(z)\,.
\end{equation}
Finally, for all $z\in\mathbb{C}\setminus\Sigma_{\rm SD}$ lying
outside the closure of all discs, we set
\begin{equation}\label{eq:XfromXdotlast}
\hat{\mat{X}}(z):=\dot{\mat{X}}(z)\,.
\end{equation}

\subsection{Error estimation.}
\label{sec:error}
\begin{figure}[h]
\begin{center}
\input{SigmaE.pstex_t}
\end{center}
\caption{\em The contour $\Sigma_E$ lies in the region $a\le\Re(z)\le b$ and $|\Im(z)|\le\epsilon$.  The circles of radius $h\epsilon$ ($h<1$) are all oriented in the clockwise direction.}
\label{fig:SigmaE}
\end{figure}
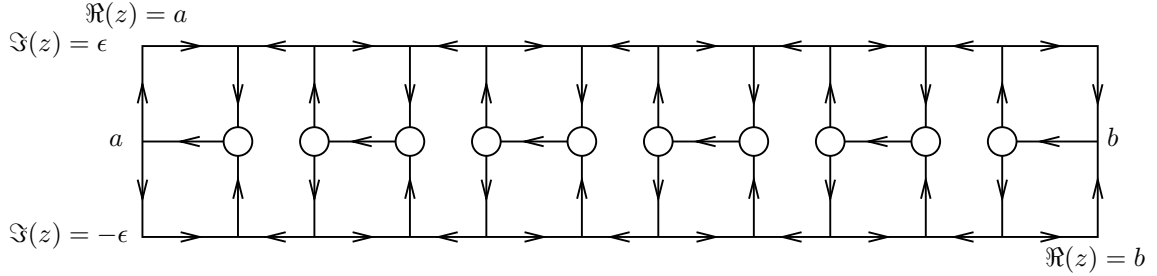

To compare the (unknown) solution $\mat{X}(z)$ of Riemann-Hilbert
Problem~\ref{rhp:X} to the explicit parametrix $\hat{\mat{X}}(z)$, we
consider the error matrix $\mat{E}(z)$ \label{symbol:matrixE} defined by
\begin{equation}\label{eq:XfromEandXhat}
\mat{E}(z):=\mat{X}(z)\hat{\mat{X}}(z)^{-1}\,.
\end{equation}
A direct calculation shows that this matrix has a continuous (and thus
analytic) extension to each band $I$ and also to the interior of each
disc $D_\Gamma^{*,*}$.  In other words, $\mat{E}(z)$ is analytic for
$z\in\mathbb{C}\setminus\Sigma_E$ where $\Sigma_E$ \label{symbol:SigmaE}
is the contour
pictured in Figure~\ref{fig:SigmaE}.
We want to deduce, for sufficiently small positive $\epsilon$, an
estimate for $\mat{E}(z)-\mathbb{I}$ that is valid in a neighborhood
of an arbitrary point of $[a,b]$.  In order to do this, it is useful
to first introduce an intermediate matrix $\mat{F}(z)$ which will
differ from $\mat{E}(z)$ only near all gaps $\Gamma$ and near the
endpoints $a$ and $b$.

\begin{figure}[h]
\begin{center}
\input{SigmaF.pstex_t}
\end{center}
\caption{\em The contour $\Sigma_F$.  Dashed lines indicate contour segments of $\Sigma_E$ to which $\mat{F}(z)$ has a continuous and hence analytic extension.  As with $\Sigma_E$, the disc boundaries are oriented in the clockwise direction.
The circular boundaries of the half-discs $D_a$ and $D_b$ are also
oriented in the clockwise direction.  The lower boundaries of all
regions $L^\nab_\Gamma$ and $L^\del_\Gamma$ are oriented from right to
left.}
\label{fig:SigmaF}
\end{figure}
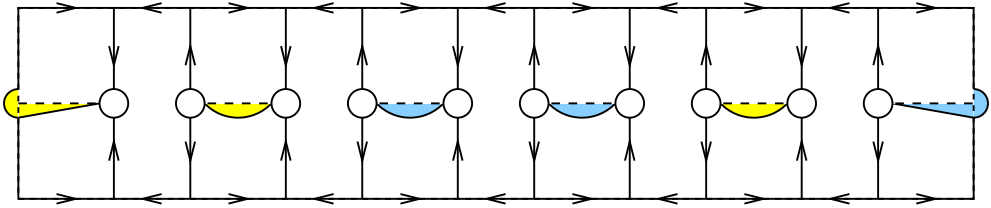

For each void interval $\Gamma$ that lies between two consecutive
bands, let $x$ and $y$ be the points where $\Gamma$ meets the
boundaries of the discs $D_\Gamma^{\nab,R}$ and $D_\Gamma^{\nab,L}$,
and let $L^\nab_\Gamma$ denote the open chord (that is, the part of a
disc bounded by a circular arc of the boundary and the straight line
segment joining the endpoints of the arc) determined by the points
$x$, $(x+y)/2-ih\epsilon$, and $y$.  If the lower constraint is
satisfied at $z=a$ and $\Gamma$ is the corresponding void interval
that meets the boundary of the disc $D_\Gamma^{\nab,L}$ at a point
$x$, then we let $L^\nab_\Gamma$ \label{symbol:LnabGamma} denote the
open triangle with vertices $a$, $x$, and $a-ih\epsilon$.  If the
lower constraint is satisfied at $z=b$ and $\Gamma$ is the
corresponding void interval that meets the boundary of the disc
$D_\Gamma^{\nab,R}$ at a point $x$, then we let $L^\nab_\Gamma$ denote
the open triangle with vertices $b$, $x$, and $b-ih\epsilon$.  The
various regions $L^\nab_\Gamma$ lie in the range $a<\Re(z)<b$ and are
illustrated with blue shading in Figure~\ref{fig:SigmaF}.  We make the
change of variables \label{symbol:matrixF}
\begin{equation}
\mat{F}(z)=\mat{E}(z)\dot{\mat{X}}(z)\left(\begin{array}{cc}
1 & -iT_\nab(z)e^{\gamma-\eta(z)+2\kappa
g(z)}e^{iN\theta_\Gamma}e^{-N\xi_\Gamma(z)}
\\\\0 & 1\end{array}\right)
\dot{\mat{X}}(z)^{-1}\,,\hspace{0.2 in}\mbox{for $z\in L_\Gamma^\nab$.}
\label{eq:EFfirst}
\end{equation}
For each saturated region $\Gamma$ that lies between two consecutive
bands, let $x$ and $y$ be the points where $\Gamma$ meets the
boundaries of the discs $D_\Gamma^{\del,R}$ and $D_\Gamma^{\del,L}$,
and let $L^\del_\Gamma$ \label{symbol:LdelGamma} 
denote the open chord determined by the points
$x$, $(x+y)/2-ih\epsilon$, and $y$.  If the upper constraint is
satisfied at $z=a$ and $\Gamma$ is the corresponding saturated region
that meets the boundary of the disc $D_\Gamma^{\del,L}$ at a point
$x$, then we let $L^\del_\Gamma$ denote the open triangle with
vertices $a$, $x$, and $a-ih\epsilon$.  If the upper constraint is
satisfied at $z=b$ and $\Gamma$ is the corresponding saturated region
that meets the boundary of the disc $D_\Gamma^{\del,R}$ at a point
$x$, then we let $L^\del_\Gamma$ denote the open triangle with
vertices $b$, $x$, and $b-ih\epsilon$.  The various regions
$L^\del_\Gamma$ lie in the range $a<\Re(z)<b$ and are illustrated with
yellow shading in Figure~\ref{fig:SigmaF}.  We make the change of
variables
\begin{equation}
\mat{F}(z)=\mat{E}(z)\dot{\mat{X}}(z)
\left(\begin{array}{cc}
1 & 0 \\\\ -iT_\del(z)e^{\eta(z)-\gamma-2\kappa
g(z)}e^{-iN\theta_\Gamma} e^{-N\xi_\Gamma(z)} & 1\end{array}\right)
\dot{\mat{X}}(z)^{-1}\,,\hspace{0.2 in}\mbox{for $z\in L^\del_\Gamma$.}
\end{equation}
Next, we define two open half-discs:
$D_a=\{z|\Re(z)<a,|z-a|<h\epsilon\}$ and
$D_b=\{z|\Re(z)>b,|z-b|<h\epsilon\}$.  In each of these half-discs
centered at an endpoint where the lower constraint is active
(indicated with blue shading in Figure~\ref{fig:SigmaF}) we set
\begin{equation}
\mat{F}(z)=\mat{E}(z)\dot{\mat{X}}(z)\left(\begin{array}{cc}
1 & \displaystyle
-ie^{\gamma-\eta(z)}\frac{\displaystyle\prod_{n\in\del}(z-x_{N,n})}{\displaystyle
\prod_{n\in\nab}(z-x_{N,n})}e^{N(\ell_c-V(z)-i\theta^0(z)/2)}e^{(k-\#\del)g(z)}
\\\\
0 & 1\end{array}\right)
\dot{\mat{X}}(z)^{-1}\,,
\label{eq:EFabvoid}
\end{equation}
and in each half-disc centered at an endpoint where the upper
constraint is active (indicated with yellow shading in
Figure~\ref{fig:SigmaF}) we set
\begin{equation}
\mat{F}(z)=\mat{E}(z)\dot{\mat{X}}(z)\left(
\begin{array}{cc}
1 & 0 \\\\\displaystyle
-ie^{\eta(z)-\gamma}\frac{\displaystyle\prod_{n\in\nab}(z-x_{N,n})}
{\displaystyle\prod_{n\in\del}(z-x_{N,n})}
e^{N(V(z)-\ell_c-i\theta^0(z)/2)}e^{(\#\del-k)g(z)} &
1\end{array}\right)
\dot{\mat{X}}(z)^{-1}\,.
\label{eq:EFabsat}
\end{equation}
It is important to observe that the matrix relating $\mat{F}(z)$ and
$\mat{E}(z)$ in (\ref{eq:EFabvoid}) and (\ref{eq:EFabsat}) is always
an analytic function of $z$ in the half-disc under consideration.
Indeed, the poles are all located in $[a,b]$, $e^{g(z)}$ is analytic
for $z\in\mathbb{C}\setminus[a,b]$, and $k-\#\del\in\mathbb{Z}$.  For
all remaining $z\in\mathbb{C}\setminus\Sigma_E$, we set
$\mat{F}(z)=\mat{E}(z)$.

\begin{Lemma}\label{lemEFab}
The matrix $\mat{F}(z)$ admits a continuous and hence analytic
extension to the upper boundaries of all regions $L_\Gamma^\nab$ and
$L_\Gamma^\del$, as well as to the segments $\Re(z)=a$ and $\Re(z)=b$
with $|\Im(z)|<h\epsilon$.
\end{Lemma}

\begin{proof}
This is rather straightforward to show once one makes the following
observations. First, in the quarter discs $D_a\cap \{z|\Im(z)>0\}$ and
$D_b\cap\{z|\Im(z)>0\}$ centered at endpoints where the lower
constraint holds, we have the identity
\begin{equation}
-ie^{\gamma-\eta(z)}\frac{\displaystyle\prod_{n\in\del}(z-x_{N,n})}
{\displaystyle\prod_{n\in\nab}(z-x_{N,n})}
e^{N(\ell_c-V(z)-i\theta^0(z)/2)}e^{(k-\#\del)g(z)} =
-iY(z)e^{\gamma-\eta(z)+2\kappa
g(z)}e^{-iN\theta_\Gamma}e^{-iN\theta^0(z)}e^{-N\xi_\Gamma(z)}\,.
\label{eq:uhplower}
\end{equation}
Here $\Gamma$ refers to the void that is adjacent to the endpoint.  If
the upper constraint is active, we have in the same region the
identity
\begin{equation}
-ie^{\eta(z)-\gamma}\frac{\displaystyle\prod_{n\in\nab}(z-x_{N,n})}
{\displaystyle\prod_{n\in\del}(z-x_{N,n})}e^{N(V(z)-\ell_c-i\theta^0(z)/2)}
e^{(\#\del-k)g(z)}=-iY(z)^{-1}e^{\eta(z)-\gamma-2\kappa g(z)}
e^{iN\theta_\Gamma}e^{-iN\theta^0(z)}e^{-N\xi_\Gamma(z)}\,.
\label{eq:uhpupper}
\end{equation}
On the other hand, in the quarter discs $D_a\cap\{z|\Im(z)<0\}$ and
$D_b\cap\{z|\Im(z)<0\}$ centered at endpoints where the lower
constraint holds, we have
\begin{equation}
-ie^{\gamma-\eta(z)}\frac{\displaystyle\prod_{n\in\del}(z-x_{N,n})}
{\displaystyle\prod_{n\in\nab}(z-x_{N,n})}
e^{N(\ell_c-V(z)-i\theta^0(z)/2)}e^{(k-\#\del)g(z)} =
-iY(z)e^{\gamma-\eta(z)+2\kappa
g(z)}e^{iN\theta_\Gamma}e^{-N\xi_\Gamma(z)}\,.
\label{eq:lhplower}
\end{equation}
If the upper constraint is active then in the same region
\begin{equation}
-ie^{\eta(z)-\gamma}\frac{\displaystyle\prod_{n\in\nab}(z-x_{N,n})}
{\displaystyle\prod_{n\in\del}(z-x_{N,n})}e^{N(V(z)-\ell_c-i\theta^0(z)/2)}
e^{(\#\del-k)g(z)}= -iY(z)^{-1}e^{\eta(z)-\gamma-2\kappa
g(z)}e^{-iN\theta_\Gamma}e^{-N\xi_\Gamma(z)}\,.
\label{eq:lhpupper}
\end{equation}
The claimed continuity of $\mat{F}(z)$ follows from these identities
upon using the definition $\mat{E}(z)=\mat{X}(z)\dot{\mat{X}}(z)^{-1}$
(since $\hat{\mat{X}}(z)=\dot{\mat{X}}(z)$ in for all relevant $z$ in
the current context), the jump conditions satisfied by $\mat{X}(z)$
and $\dot{\mat{X}}(z)$, and the relations
(\ref{eq:Tsqrtidentity_down}) connecting $T_\del(z)$, $T_\nab(z)$, and
$Y(z)$ for $\Im(z)<0$.
\end{proof}

The contour $\Sigma_F$ \label{symbol:SigmaF} 
where $\mat{F}(z)$ fails to be analytic is
shown in Figure~\ref{fig:SigmaF}.
In order to estimate $\mat{E}(z)^{-1}\mat{F}(z)-\mathbb{I}$, and
subsequently to estimate $\mat{F}(z)-\mathbb{I}$, we will now need to
recall the behavior of the functions $T_\nab(z)$ and
$T_\del(z)$ in the asymptotic limit $N\rightarrow\infty$.

It follows from Proposition~\ref{prop:tasymp} that in each region
$L_\Gamma^\nab$ (respectively $L_\Gamma^\del$), $T_\nab(z)$
(respectively $T_\del(z)$) is uniformly bounded.  Furthermore, in any
half disc $D_a$ or $D_b$ centered at an endpoint where the lower
constraint is active the function $Y(z)$ is uniformly bounded, and in
any half disc centered at an endpoint where the upper constraint is
active the function $Y(z)^{-1}$ is uniformly bounded.  Using the
identities (\ref{eq:uhplower})--(\ref{eq:lhpupper}) and the
variational inequalities (\ref{eq:voidinequality}) and
(\ref{eq:saturatedregioninequality}) which control
$\Re(\xi_\Gamma(z))$ in these regions (and noting that in particular
$\Re(\xi_\Gamma(a))>0$ and $\Re(\xi_\Gamma(b))>0$ by assumption ---
see \S~\ref{sec:C3}), we have the following result:
\begin{Lemma}
Let the contour parameter $\epsilon>0$ be sufficiently small.  Then,
there are constants $C_{1,\epsilon}>0$ and $C_{2,\epsilon}>0$ such
that for all sufficiently large $N$,
\begin{equation}
\sup_{z\in\mathbb{C}\setminus (\Sigma_E\cup\Sigma_F)}
\|\mat{E}(z)^{-1}\mat{F}(z)-\mathbb{I}\|\le C_{1,\epsilon}e^{-C_{2,\epsilon}N}\,.
\end{equation}
Here $\|\cdot\|$ denotes an arbitrary matrix norm.
\label{lemma:EFcompare}
\end{Lemma}

Being obtained from $\mat{X}(z)$ satisfying Riemann-Hilbert
Problem~\ref{rhp:X} by explicit transformations involving the
parametrix $\hat{\mat{X}}(z)$ as well as the explicit relations
(\ref{eq:EFfirst})--(\ref{eq:EFabsat}), the (unknown) matrix
$\mat{F}(z)$ is the solution of a Riemann-Hilbert problem as well:
\label{symbol:vF}
\begin{rhp}
Find a $2\times 2$ matrix $\mat{F}(z)$ with the following properties:
\begin{enumerate}
\item
{\bf Analyticity}: $\mat{F}(z)$ is an analytic function of $z$ for
$z\in\cx\setminus \Sigma_{F}$.
\item
{\bf Normalization}: As $z\rightarrow\infty$,
\begin{equation}
\mat{F}(z)={\mathbb I} + O\left(\frac{1}{z}\right)\,.
\end{equation}
\item
{\bf Jump Conditions}: $\mat{F}(z)$ takes uniformly continuous
boundary values on $\Sigma_{F}$ from each connected component of
$\cx\setminus\Sigma_{F}$.  For each non-self-intersection point
$z\in\Sigma_{F}$ we denote by $\mat{F}_+(z)$ ($\mat{F}_-(z)$) the
limit of $\mat{F}(w)$ as $w\rightarrow z$ from the left (right).  The
boundary values satisfy the jump condition
$\mat{F}_+(z)=\mat{F}_-(z)\mat{v}_{\mat{F}}(z)$, where for $z$ on the
lower boundary of a region $L^\nab_\Gamma$ below a void
$\Gamma\subset\Sigma_0^\nab$,
\begin{equation}
\mat{v}_{\mat{F}}(z)=
\dot{\mat{X}}(z)\left(\begin{array}{cc}
1 & iT_\nab(z)e^{\gamma-\eta(z)+2\kappa
g(z)}e^{iN\theta_\Gamma}e^{-N\xi_\Gamma(z)}\\\\ 0 & 1
\end{array}\right)\dot{\mat{X}}(z)^{-1}\,.
\end{equation}
For $z$ on the lower boundary of a region $L^\del_\Gamma$ below a
saturated region $\Gamma\subset\Sigma_0^\del$,
\begin{equation}
\mat{v}_{\mat{F}}(z)=\dot{\mat{X}}(z)\left(\begin{array}{cc}
1& 0 \\\\ iT_\del(z)e^{\eta(z)-\gamma-2\kappa
g(z)}e^{-iN\theta_\Gamma}e^{-N\xi_\Gamma(z)} & 1
\end{array}
\right)\dot{\mat{X}}(z)^{-1}\,.
\end{equation}
For $z$ in any vertical segment $\Sigma_{0\pm}^\nab\cap\Sigma_F$
meeting the boundary of a disc centered at an endpoint $z_0$ of a band
$I$,
\begin{equation}
\mat{v}_{\mat{F}}(z)=
\dot{\mat{X}}(z)\left(\begin{array}{cc}
T_\nab(z)^{\pm 1/2} & 0 \\\\\displaystyle
-iT_\nab(z)^{-1/2}e^{\eta(z)-\gamma-2\kappa g(z)}e^{\pm iN\theta(z_0)}
\exp\left(\pm 2\pi iNc\int_{z_0}^z\psi_I(s)\,ds\right) &
T_\nab(z)^{\mp 1/2}
\end{array}\right)\dot{\mat{X}}(z)^{-1}\,.
\label{eq:vF0pmnab}
\end{equation}
For $z$ in any vertical segment $\Sigma_{0\pm}^\del\cap\Sigma_F$
meeting the boundary of a disc centered at an endpoint $z_0$ of a band
$I$,
\begin{equation}
\mat{v}_{\mat{F}}(z)=\dot{\mat{X}}(z)
\left(\begin{array}{cc}
T_\del(z)^{\mp 1/2} & \displaystyle
-iT_\del(z)^{-1/2}e^{\gamma-\eta(z)+2\kappa g(z)} e^{\mp
iN\theta(z_0)}\exp\left(\pm 2\pi
iNc\int_{z_0}^z\overline{\psi}_I(z)\,ds
\right) \\\\0 & T_\del(z)^{\pm 1/2}
\end{array}\right)\dot{\mat{X}}(z)^{-1}\,.
\label{eq:vF0pmdel}
\end{equation}
For $z$ in any segment $\Sigma_{\Gamma\pm}^\nab\cap \Sigma_F$ parallel
to a void $\Gamma\subset\Sigma_0^\nab$ or with $\Re(z)=a$ or
$\Re(z)=b$,
\begin{equation}
\mat{v}_{\mat{F}}(z)=\dot{\mat{X}}(z)
\left(\begin{array}{cc}
1 & iY(z)e^{\gamma-\eta(z)+2\kappa g(z)}e^{\mp iN\theta_\Gamma}e^{\mp
iN\theta^0(z)}e^{-N\xi_\Gamma(z)}\\\\ 0 & 1
\end{array}\right)\dot{\mat{X}}(z)^{-1}\,,
\end{equation}
and for $z$ in the semicircular boundary of a half-disc $D_a$ or $D_b$
centered at an endpoint where the lower constraint is active,
\begin{equation}
\mat{v}_{\mat{F}}(z)=\dot{\mat{X}}(z)
\left(\begin{array}{cc}
1 &
\displaystyle ie^{\gamma-\eta(z)}\frac{\displaystyle\prod_{n\in\del}(z-x_{N,n})}
{\displaystyle\prod_{n\in\nab}(z-x_{N,n})}
e^{N(\ell_c-V(z)-i\theta^0(z)/2)}e^{(k-\#\del)g(z)}
\\\\
0 & 1
\end{array}\right)\dot{\mat{X}}(z)^{-1}\,.
\end{equation}
For $z$ in any segment $\Sigma_{\Gamma\pm}^\del$ parallel to a
saturated region $\Gamma\subset\Sigma_0^\del$ or with $\Re(z)=a$ or
$\Re(z)=b$,
\begin{equation}
\mat{v}_{\mat{F}}(z)=\dot{\mat{X}}(z)\left(\begin{array}{cc}
1 & 0 \\\\ iY(z)^{-1}e^{\eta(z)-\gamma-2\kappa g(z)}e^{\pm
iN\theta_\Gamma}e^{\mp iN\theta^0(z)}e^{-N\xi_\Gamma(z)} & 1
\end{array}\right)\dot{\mat{X}}(z)^{-1}\,,
\end{equation}
and for $z$ in the semicircular boundary of a half-disc $D_a$ or $D_b$
centered at an endpoint where the upper constraint is active,
\begin{equation}
\mat{v}_{\mat{F}}(z)=\dot{\mat{X}}(z)\left(\begin{array}{cc}
1 & 0 \\\\ \displaystyle
ie^{\eta(z)-\gamma}\frac{\displaystyle\prod_{n\in\nab}(z-x_{N,n})}
{\displaystyle\prod_{n\in\del}(z-x_{N,n})}
e^{N(V(z)-\ell_c-i\theta^0(z)/2)}e^{(\#\del-k)g(z)} & 1
\end{array}\right)\dot{\mat{X}}(z)^{-1}\,.
\end{equation}
With a sequence $\{y_N\}_{N=0}^\infty$ determined as in the
formulation of Riemann-Hilbert Problem~\ref{rhp:X}, we have that for
$z$ in any segment $\Sigma_{I\pm}$ parallel to any band $I$,
\begin{equation}
\mat{v}_{\mat{F}}(z)=\dot{\mat{X}}(z)\left(\begin{array}{cc}
T_\del(z)^{-1/2} &v_{12}^\pm(z)
\\\\
v_{21}^\pm(z) & T_\nab(z)^{-1/2}
\end{array}\right)^{\pm 1}\dot{\mat{X}}(z)^{-1}\,,
\label{eq:vFIpm}
\end{equation}
where
\begin{equation}
\begin{array}{rcl}
v_{12}^\pm(z)&:=&\displaystyle
\mp i T_\del(z)^{-1/2}e^{\gamma-\eta(z)+2\kappa g(z)}e^{\mp iN\theta(y_N)}
\exp\left(\pm 2\pi iNc\int_{y_N}^z\overline{\psi}_I(s)\,ds\right)\,,
\\\\
v_{21}^\pm(z)&:=&\displaystyle
\mp iT_\nab(z)^{-1/2}e^{\eta(z)-\gamma-2\kappa g(z)}e^{\pm iN\theta(y_N)}
\exp\left(\pm 2\pi iNc\int_{y_N}^z\psi_I(s)\,ds\right)\,.
\end{array}
\end{equation}
Finally, for $z$ in the clockwise-oriented boundary of any disc
$D_\Gamma^{\nab,L}$,
\begin{equation}
\mat{v}_{\mat{F}}(z)=\hat{\mat{X}}_\Gamma^{\nab,L}(z)\dot{\mat{X}}(z)^{-1}\,,
\end{equation}
for $z$ in the clockwise-oriented boundary of any disc
$D_\Gamma^{\nab,R}$,
\begin{equation}
\mat{v}_{\mat{F}}(z)=\hat{\mat{X}}_\Gamma^{\nab,R}(z)\dot{\mat{X}}(z)^{-1}\,,
\end{equation}
for $z$ in the clockwise-oriented boundary of any disc
$D_\Gamma^{\del,L}$,
\begin{equation}
\mat{v}_{\mat{F}}(z)=\hat{\mat{X}}_\Gamma^{\del,L}(z)\dot{\mat{X}}(z)^{-1}\,,
\end{equation}
and for $z$ in the clockwise-oriented boundary of any disc
$D_\Gamma^{\del,R}$,
\begin{equation}
\mat{v}_{\mat{F}}(z)=\hat{\mat{X}}_\Gamma^{\del,R}(z)\dot{\mat{X}}(z)^{-1}\,,
\end{equation}
\end{enumerate}
\label{rhp:F}
\end{rhp}

We have the following characterization of the jump matrix for
$\mat{F}(z)$.
\begin{Lemma}
Let the parameter $\epsilon>0$ of the contour $\Sigma_F$ be
sufficiently small.  Then there is a constant $C_\epsilon>0$ such that
\begin{equation}
\sup_{z\in\Sigma_F}\|\mat{v}_{\mat{F}}(z)-\mathbb{I}\|\le\frac{C_\epsilon}{N}
\end{equation}
holds for sufficiently large $N$.
\label{Lemma:vF}
\end{Lemma}

\begin{proof}
The estimates on the boundaries of the discs $D_\Gamma^{*,*}$ all
follow from Proposition~\ref{prop:Dproperties}.  For the remaining
parts of $\Sigma_F$, we note that by Proposition~\ref{prop:dotXbound}
$\dot{\mat{X}}(z)$ and $\dot{\mat{X}}(z)^{-1}$ are both uniformly
bounded for $z\in\Sigma_F$; thus it suffices to prove an estimate of
the same order for $\dot{\mat{X}}(z)^{-1}
\mat{v}_{\mat{F}}(z)\dot{\mat{X}}(z)$.  Now using Proposition~\ref{prop:tasymp}
one sees that the diagonal entries in (\ref{eq:vF0pmnab}),
(\ref{eq:vF0pmdel}), and (\ref{eq:vFIpm}) all differ from one by a
quantity of order $1/N$.  All off-diagonal matrix elements are
exponentially small as $N\rightarrow\infty$ for two different reasons.
First, we recall the variational inequalities
(\ref{eq:voidinequality}) and (\ref{eq:saturatedregioninequality})
that hold on the real axis in the voids and saturated regions
respectively; these control the off-diagonal entries of
$\dot{\mat{X}}(z)^{-1}
\mat{v}_{\mat{F}}(z)\dot{\mat{X}}(z)$ involving a factor $e^{-N\xi_\Gamma(z)}$ for $\epsilon$ sufficiently small that the inequality $\Re(\xi_\Gamma(z))>0$ holds on relevant portions of $\Sigma_F$ as it does in the gap $\Gamma\subset [a,b]$.  Second, we recall the inequalities
$\psi_I(x)>0$ and $\overline{\psi}_I(x)>0$ that hold in each band $I$
together with the presumed square-root vanishing of $\psi_I(x)$ at
band edges where the lower constraint becomes active and of
$\overline{\psi}_I(x)$ at band edges where the upper constraint
becomes active; these facts control the off-diagonal entries of
$\dot{\mat{X}}(z)^{-1}
\mat{v}_{\mat{F}}(z)\dot{\mat{X}}(z)$ in the segments $\Sigma_{0\pm}^\nab\cap\Sigma_F$, $\Sigma_{0\pm}^\del\cap\Sigma_F$, and $\Sigma_{I\pm}$.
\end{proof}

\begin{Lemma}
Let the contour parameter $\epsilon>0$ be sufficiently small.  Then
Riemann-Hilbert Problem~\ref{rhp:F} has a unique solution for
sufficiently large $N$, and the solution has the Cauchy integral
representation
\begin{equation}
\mat{F}(z)={\mathbb I} + \int_{\Sigma_{F}}(z-s)^{-1}\mat{m}(s)\,ds
\label{eq:FCauchy}
\end{equation}
where $\mat{m}(\cdot)$ is an arcwise-continuous matrix function in
$L^2(\Sigma_{F})$.  There is a constant $L_\epsilon>0$ such that
\begin{equation}
\|\mat{m}\|_2\le\frac{L_\epsilon}{N}
\label{eq:L2density}
\end{equation}
holds for all sufficiently large $N$.  Also, $\det(\mat{F}(z))=1$ for
all $z\in\cx\setminus\Sigma_{\mat{F}}$.
\label{Lemma:Fexist}
\end{Lemma}

\begin{proof}
This is essentially a consequence of the theory of matrix
Riemann-Hilbert problems with $L^2$ boundary values and uniformly
near-identity jump matrices (see, for example, \cite{Zhou89}).  The
key idea is that it is possible to convert the Riemann-Hilbert problem
into a system of singular integral equations of the form
$(1-B)\mat{u}={\mathbb I}$ where $B$ is a singular integral operator
acting on matrix functions $\mat{u}(z)$ defined on $\Sigma_{F}$; then the
desired density $\mat{m}(z)$ is proportional to both $\mat{u}(z)$ and
$\mat{v}_{\mat{F}}(z)-{\mathbb I}$.  The operator $B$ can be written
as a composition of multiplication by $\mat{v}_{\mat{F}}(z)-{\mathbb
I}$ and a singular integral operator with Cauchy kernel.  It is a deep
result of modern harmonic analysis
\cite{CoifmanMM82} that the Cauchy-kernel singular integral operators
are bounded in $L^2$ on contours that may be decomposed as finite
unions of graphs of Lipschitz functions (an appropriate Lipschitz
condition should be also satisfied at each self-intersection point).
The norm of $B$ in $L^2(\Sigma_{F})$ is proportional to the product of
$\|\mat{v}_{\mat{F}}(z)-\mathbb{I}\|_\infty$ which we know can be made
arbitrarily small according to Lemma~\ref{Lemma:vF}, and the $L^2$
norm of a Cauchy integral over $\Sigma_F$ which is finite if
$\epsilon$ is taken to be sufficiently small (this makes all
self-intersections of $\Sigma_F$ non-tangential).  Thus, for
sufficiently large $N$ we will have $\|B\|_2<1$ and the integral
equation for $\mat{u}(z)$ can be solved in $L^2(\Sigma_{F})$ by a
Neumann series:  $\mat{u}(z)=\mathbb{I}+B\mathbb{I}+B^2\mathbb{I} +\dots$.

We therefore have the invertibility of the operator $1-B$ for
sufficiently large $N$, with $\|(1-B)^{-1}\|_2$ bounded independently
of $N$, and thus the existence of $\mat{u}\in L^2(\Sigma_{F})$.  Moreover,
since the total length of $\Sigma_{F}$ is independent of $N$, we get
$\|\mat{u}\|_2=\|(1-B)^{-1}{\mathbb I}\|_2$ being bounded uniformly with
respect to $N$ as well.  This proves
\eqref{eq:L2density}, since $\mat{m}(z)$ is proportional to the
product of $\mat{u}(z)$ and $\mat{v}_{\mat{F}}(z)-{\mathbb I}$.

The fact that the boundary values taken by the solution $\mat{F}(z)$
supplied by the $L^{2}$ theory are in fact uniformly continuous along
the boundary of each connected component of $\cx\setminus\Sigma_{F}$
warrants some additional explanation.  Indeed, the $L^{2}$ theory only
guarantees a solution of the Riemann-Hilbert problem taking boundary
values in the $L^{2}$ sense.  However, since the jump matrix
$\mat{v}_{\mat{F}}(z)$ is analytic on each arc of $\Sigma_{F}$, it
follows that both $\mat{F}_{+}(z)$ and $\mat{F}_{-}(z)$ may be
continued analytically through to the opposite side of each arc, and
then from Morera's Theorem we deduce that not only is $\mat{F}(z)$
continuous up to the boundary, but also both boundary values are
analytic functions of $z$.  That uniform continuity extends even to
self-intersection points of $\Sigma_{F}$ can be shown using the
compatibility of the limiting values of $\mat{v}_{\mat{F}}(z)$ along
all arcs meeting at such a point; namely the cyclic product of the
limiting values is the identity matrix for all self-intersection
points.  Thus, the unique $L^{2}$ solution is in fact a classical
solution of the Riemann-Hilbert problem.
\end{proof}

Thus we arrive at the main result of this section.
\begin{prop}
Let the contour parameter $\epsilon$ be sufficiently small.  Then for
each closed set $K\subset\mathbb{C}\setminus\Sigma_F$, not necessarily
bounded, there is a constant $Q_{K,\epsilon}>0$ such that
\begin{equation}
\sup_{z\in K}\|\mat{E}(z)-\mathbb{I}\|\le\frac{Q_{K,\epsilon}}{N}
\end{equation}
holds for all sufficiently large $N$.  Recall that
$\mat{E}(z)=\mat{X}(z)\hat{\mat{X}}(z)^{-1}$.
\label{prop:E}
\end{prop}

\begin{proof}
From \eqref{eq:FCauchy} and \eqref{eq:L2density} we obtain the desired
estimate for the matrix $\mat{F}(z)$.  To complete the proof, we
recall Lemma~\ref{lemma:EFcompare}.
\end{proof}

\section{Discrete Orthogonal Polynomials:  
Proofs of Theorems Stated in \S\ref{sec:actualtheorems}}
\label{sec:asymptoticspi}

In this section, we start with the exact formula for $\mat{X}(z)$
valid in the entire complex $z$-plane:
\begin{equation}
\mat{X}(z)=\mat{E}(z)\hat{\mat{X}}(z)\,.
\label{eq:Rformula}
\end{equation}
This is written in terms of the explicit global parametrix and the
matrix $\mat{E}(z)$ which while not explicit is characterized by
Proposition~\ref{prop:E}.  We then work backwards to the matrix
$\mat{P}(z;N,k)$ and therefore obtain exact formulae for the monic
polynomials $\pi_{N,k}(z)$ valid in the whole complex plane as well as
the normalization constants $\gamma_{N,k}$ and recurrence coefficients
$a_{N,k}$ and $b_{N,k}$ in terms of the matrix elements of
$\hat{\mat{X}}(z)$ and $\mat{E}(z)$, and their asymptotics for large
$z$.  Then, under various conditions on $z$ we extract simple
asymptotic formulae by direct asymptotic expansion of the exact
formulae.  In particular, we will obtain Plancherel-Rotach type
asymptotics of the monic polynomials $\pi_{N,k}(z)$ for real $z$ in
the interval $[a,b]$ of accumulation of the discrete nodes of support
of the weights.

\subsection{Asymptotic analysis of $\mat{P}(z;N,k)$ for $z$ outside the
interval $[a,b]$.}
\label{sec:outside}

\subsubsection{Asymptotic behavior of $\pi_{N,k}(z)$ for $z$ outside the
interval $[a,b]$.  Proof of Theorem~\ref{theorem:outside}.}  Let
$K\subset\mathbb{C}\setminus [a,b]$ be a fixed closed set, not
necessarily bounded.  The parameter $\epsilon$ in the contour
$\Sigma_F$ may then be fixed at such a sufficiently small positive
value that $K\cap\Sigma_F=\emptyset$ and $K$ is contained in the
unbounded component of $\mathbb{C}\setminus\Sigma_F$.  For $z\in K$ we
thus have from \eqref{eq:XD}, \eqref{eq:Dfirst},
\eqref{eq:XfromXdotlast}, and \eqref{eq:Rformula}, we have
\begin{equation}
%\begin{array}{rcl}
\mat{P}(z;N,k) =
%&=&\displaystyle \mat{Q}(z;N,k)\prod_{j\in\del}(z-x_{N,j})^{\sigma_3}\\\\
%&=&\displaystyle
%\mat{R}(z)\prod_{j\in\del}(z-x_{N,j})^{\sigma_3}\\\\
%&=&\displaystyle e^{-(N\ell_c+\gamma)\sigma_3/2}\mat{S}(z)
%e^{(N\ell_c+\gamma)\sigma_3/2}e^{(k-\#\del)g(z)\sigma_3}
%\prod_{j\in\del}(z-x_{N,j})^{\sigma_3}\\\\
%&=& \displaystyle e^{-(N\ell_c+\gamma)\sigma_3/2}\mat{X}(z)
%e^{(N\ell_c+\gamma)\sigma_3/2}e^{(k-\#\del)g(z)\sigma_3}
%\prod_{j\in\del}(z-x_{N,j})^{\sigma_3}\\\\
%&=&
\displaystyle e^{-(N\ell_c+\gamma)\sigma_3/2}\mat{E}(z)\dot{\mat{X}}(z)
e^{(N\ell_c+\gamma)\sigma_3/2}e^{(k-\#\del)g(z)\sigma_3}
\prod_{n\in\del}(z-x_{N,n})^{\sigma_3}\,.
%\end{array}
\end{equation}
Since $z\in K$ is bounded away from $[a,b]$ we use the midpoint rule
to obtain
\begin{equation}
e^{(k-\#\del)g(z)}\prod_{n\in\del}(z-x_{N,n}) =
e^{(k-\#\del)g(z)}\exp\left(N\int_{\Sigma_0^\del}\log(z-x)\rho^0(x)\,dx
\right)\cdot\left(1+O\left(\frac{1}{N}\right)\right)\,,
\end{equation}
where the error term is uniform for $z\in K$.  Combining this result
with (\ref{eq:rhodef}) and (\ref{eq:gdef}), and recalling that
$k=cN+\kappa$ we get
\begin{equation}
e^{(k-\#\del)g(z)}\prod_{n\in\del}(z-x_{N,n})=e^{\kappa
g(z)}e^{NL_c(z)}\cdot\left( 1+O\left(\frac{1}{N}\right)\right)\,,
\end{equation}
where $L_c(z)$ is defined in \eqref{eq:Lcdef}.  Note that the product
$e^{\kappa g(z)}e^{NL_c(z)}$ is analytic for $z\in\mathbb{C}\setminus
[a,b]$.  In particular, this analysis leads to the following formula
\begin{equation}
P_{11}(z;N,k)=\left[E_{11}(z)\dot{X}_{11}(z)e^{\kappa g(z)} +
E_{12}(z)\dot{X}_{21}(z)e^{\kappa
g(z)}\right]e^{NL_c(z)}\cdot\left(1+O\left(\frac{1}{N}\right)
\right)\,.
\end{equation}
We use Proposition~\ref{prop:E} to estimate $\mat{E}(z)-\mathbb{I}$
and Proposition~\ref{prop:dotXbound} to characterize
$\dot{X}_{11}(z)$. The proof is complete upon noting that
$W(z)=\dot{X}_{11}(z)e^{\kappa g(z)}$, using the formulae for
$\dot{X}_{11}(z)$ obtained in Appendix~\ref{sec:thetasolve}, and
recalling from Proposition~\ref{prop:solnrhp} that
$P_{11}(z;N,k)=\pi_{N,k}(z)$.

\subsubsection{Asymptotic behavior of the leading coefficients 
$\gamma_{N,k}$ and of the recurrence coefficients $a_{N,k}$ and
$b_{N,k}$.  Proof of Theorem~\ref{theorem:cs}.}  Taking the set $K$ in
the proof of Theorem~\ref{theorem:outside} above to be unbounded
allows us to consider $z\rightarrow\infty$.  For arbitrary fixed $N$,
we have the expansion
\begin{equation}
e^{(k-\#\del)g(z)}\prod_{n\in\del}(z-x_{N,n}) = e^{\kappa
(g(z)-\log(z))}z^k\left(1 + \frac{H_{k,N}}{z} + O\left(\frac{1}{z^2}
\right)\right)\,,
\end{equation}
as $z\rightarrow\infty$, where
\begin{equation}
H_{k,N}:= N\int_{\Sigma_0^\del}x\rho^0(x)\,dx-\sum_{n\in\del}x_{N,n}-
Nc\int_a^bx\,d\mu_{\rm min}^c(x)\,.
\end{equation}
The matrices $\mat{E}(z)$ and
$\dot{\mat{X}}(z)e^{\kappa(g(z)-\log(z))\sigma_3}$ have asymptotic
expansions for large $z$ of the form
\begin{equation}
\begin{array}{rcl}
\mat{E}(z)&=&\displaystyle \mathbb{I}+\frac{1}{z}\mat{E}^{(1)} + \frac{1}{z^2}\mat{E}^{(2)} + O\left(\frac{1}{z^3}\right)\\\\
\dot{\mat{X}}(z)e^{\kappa(g(z)-\log(z))\sigma_3} &=&
\displaystyle\mathbb{I} + \frac{1}{z}\mat{B}^{(1)} + \frac{1}{z}\mat{B}^{(2)} +O\left(\frac{1}{z^3}\right)
\end{array}
\end{equation}
as $z\rightarrow\infty$.  In terms of these coefficients we thus have
for each fixed $N$ the expansions
\begin{equation}
z^{-k}P_{11}(z)=1+\frac{H_{k,N}+B^{(1)}_{11} + E^{(1)}_{11}}{z} +
O\left(\frac{1}{z^2}\right)\,,
\end{equation}
\begin{equation}
z^{-k}P_{21}(z)=\frac{e^{N\ell_c+\gamma}}{z}(B^{(1)}_{21}+E^{(1)}_{21})
+ O\left(\frac{1}{z^2}\right)\,,
\end{equation}
\begin{equation}
\begin{array}{rcl}
z^kP_{12}(z)&=&\displaystyle
\frac{e^{-N\ell_c-\gamma}}{z}(B^{(1)}_{12}+E^{(1)}_{12}) \\\\
&&\displaystyle\,\,\,+\,\,\, \frac{e^{-N\ell_c-\gamma}}{z^2}
(B^{(2)}_{12}+E^{(2)}_{12}+E_{11}^{(1)}B^{(1)}_{12}+E^{(1)}_{12}B^{(1)}_{22}-H_{k,N}B^{(1)}_{12}-H_{k,N}E^{(1)}_{12})\\\\
&&\displaystyle\,\,\,+\,\,\,O\left(\frac{1}{z^3}\right)
\,.
\end{array}
\end{equation}
Comparing with (\ref{eq:3termcoeffs}), we therefore have the following
exact formulae in which $H_{k,N}$ does not appear:
\begin{equation}
\begin{array}{rcl}
\gamma_{N,k}&=&\displaystyle\frac{e^{(N\ell_c+\gamma)/2}}{\sqrt{B^{(1)}_{12}+E^{(1)}_{12}}}\\\\
\gamma_{N,k-1}&=&e^{(N\ell_c+\gamma)/2}\sqrt{B^{(1)}_{21}+E^{(1)}_{21}}\,,\\\\
b_{N,k-1}&=&\sqrt{(B^{(1)}_{12}+E^{(1)}_{12})(B^{(1)}_{21}+E^{(1)}_{21})}\,,\\\\
a_{N,k}&=&\displaystyle
B^{(1)}_{11}+E^{(1)}_{11}+\frac{B^{(2)}_{12}+E^{(2)}_{12}+E^{(1)}_{11}B^{(1)}_{12}
+ E^{(1)}_{12}B^{(1)}_{22}}{B^{(1)}_{12}+E^{(1)}_{12}}\,.
\end{array}
\end{equation}
Now, for sufficiently large $z$, we have $\mat{E}(z)=\mat{F}(z)$, and
therefore Lemma~\ref{Lemma:Fexist} and in particular the Cauchy
integral representation (\ref{eq:FCauchy}) of $\mat{F}(z)$ implies
that the coefficients $E^{(1)}_{jk}$ and $E^{(2)}_{jk}$ are all of
order $1/N$ as $N\rightarrow\infty$.  Furthermore,
$\dot{\mat{X}}(z)e^{\kappa(g(z)-\log(z))}$ is a matrix that for some
fixed $R>0$ is analytic and uniformly bounded (independently of $N$)
for $|z|>R$, which implies that the coefficients $B_{jk}^{(1)}$ and
$B_{jk}^{(2)}$ remain bounded as $N\rightarrow\infty$.  In fact, for
sufficiently large $N$, $B_{12}^{(1)}$ and $B_{21}^{(1)}$ are bounded
away from zero, and thus
\begin{equation}
\begin{array}{rcl}
\gamma_{N,k}&=&\displaystyle \frac{e^{(N\ell_c+\gamma)/2}}{\sqrt{B_{12}^{(1)}}}\left(1+O\left(\frac{1}{N}\right)\right)\,,\\\\
\gamma_{N,k-1}&=&\displaystyle e^{(N\ell_c+\gamma)/2}\sqrt{B_{21}^{(1)}}\left(1+O\left(\frac{1}{N}\right)\right)\,,\\\\
b_{N,k-1}&=&\displaystyle
\sqrt{B_{12}^{(1)}B_{21}^{(1)}}\left(1+O\left(\frac{1}{N}\right)\right)\,,\\\\
a_{N,k}&=&\displaystyle B_{11}^{(1)} +
\frac{B_{12}^{(2)}}{B_{12}^{(1)}} + O\left(\frac{1}{N}\right)\,.
\end{array}
\end{equation}
Using the formulae obtained in Proposition~\ref{prop:dotXasymp}
established in Appendix~\ref{sec:thetasolve} then completes the proof.
It should be remarked that the quantities $B_{12}^{(1)}$ and $B_{21}^{(1)}$
are necessarily positive, since $\ell_c$ and $\gamma$ are real.

\subsection{Asymptotic behavior of $\pi_{N,k}(z)$ for $z$ near a void of
$[a,b]$.  Proof of Theorem~\ref{theorem:lower}.}
\label{sec:lower}
The variational inequality (\ref{eq:voidinequality}) holds strictly
throughout the closed interval $J\subset [a,b]$, and while it is
possible for either $a$ or $b$ to be an endpoint of $J$, neither
endpoint of $J$ may be a band edge.  We choose the contour parameter
$\epsilon$ to be sufficiently small that Proposition~\ref{prop:E}
controls $\mat{E}(z)-\mathbb{I}$, and then take $\delta$ to be small
enough that $K_J^\delta\cap
\Sigma_F=\emptyset$.
%Consider the compact set
%\begin{equation}
%K_J^\epsilon:=\bigcup_{w\in J}\{z\in\mathbb{C}\,\,\,\mbox{such %that}\,\,\,|z-w|\le h\epsilon/2\}\,.
%\label{eq:Kz1z2eps}
%\end{equation}
Then, for all $z\in K_J^\delta$, regardless of whether $\Im(z)$ is
positive or negative, or of whether $\Re(z)\in (a,b)$ or not, we have
the exact formula
\begin{equation}
\pi_{N,k}(z)=\left[E_{11}(z)\dot{X}_{11}(z) + E_{12}(z)\dot{X}_{21}(z)
\right]e^{(k-\#\del)g(z)}\prod_{n\in\del}(z-x_{N,n})\,.
\label{eq:piNkexactouter}
\end{equation}
This follows from \eqref{eq:XD}, \eqref{eq:Dfirst},
\eqref{eq:Dnexttovoid},
\eqref{eq:XfromXdotlast}, \eqref{eq:Rformula}, and 
Proposition~\ref{prop:solnrhp}.  It is not hard to verify that the
right-hand side extends analytically to the whole compact set
$K_J^\delta$.  Since each node $x_{N,n}$ with $n\in\del$ is bounded
away from $K_J^\delta$, we may approximate the product to within a
relative error of order $1/N$ uniform in $K_J^\delta$ to find
\begin{equation}
\pi_{N,k}(z)=\left[E_{11}(z)\dot{X}_{11}(z) + E_{12}(z)\dot{X}_{21}(z)
\right]e^{\kappa g(z)}e^{NL_c(z)}\left(1+O\left(\frac{1}{N}\right)\right)\,.
\end{equation}
Here we have used (\ref{eq:rhodef}) and (\ref{eq:gdef}) and
$k=Nc+\kappa$, and $L_c(z)$ is defined by (\ref{eq:Lcdef}).  Finally,
using Proposition~\ref{prop:E} to estimate $\mat{E}(z)-\mathbb{I}$,
and Proposition~\ref{prop:dotXbound} to uniformly bound
$\dot{X}_{11}(z)$ and $\dot{X}_{21}(z)$ for $z\in K_J^\delta$, we
arrive at
\begin{equation}
\pi_{N,k}(z)=e^{NL_c(z)}\left[\dot{X}_{11}(z)e^{\kappa g(z)} + O\left(\frac{1}{N}\right)\right]\,.
\label{eq:voidtheoremprerequisite}
\end{equation}
We recall the definition \eqref{eq:LcbarGammadef} of the analytic
function $\overline{L}_c^\Gamma(z)$, and note that $W(z)=\dot{X}_{11}(z)e^{\kappa g(z)}$.

The estimate (\ref{eq:voidtheoremestimate}) follows from
(\ref{eq:voidtheoremprerequisite}) because
$e^{N(L_c(z)-\overline{L}_c^\Gamma(z))}$ is uniformly bounded in
$K_J^\delta$.  Indeed, we have
\begin{equation}
L_c(z)-\overline{L}_c^\Gamma(z)=-\frac{i\theta_\Gamma}{2}\cdot
{\rm sgn}(\Im(z))\,.
%\left[\pi\int_{y<x\in\Sigma_0^\del}\rho^0(x)\,dx 
%-\frac{1}{2}\theta(x_\Gamma)\right]\,,
\label{eq:Lcdiff}
\end{equation}
%where $x_\Gamma$ is any point in the void $\Gamma$, and where $y$ is
%the nearest transition point to the right of $\Gamma$ (or
%$y=b$ if there is no transition point to the right of $\Gamma$).  By
%the quantization condition (\ref{eq:yquantize}) the integral is a
%rational number $M/N$ with denominator $N$.  
Thus the right-hand side of \eqref{eq:Lcdiff} is simply a different
imaginary constant in each half-plane.  This also establishes the
uniform boundedness of $A^\nab_\Gamma(z)$ when we use
Proposition~\ref{prop:dotXbound} to bound $\dot{X}_{11}(z)$.
%The invariance of $\dot{X}_{11}(z)e^{\kappa g(z)}$ is also guaranteed by 
%Proposition~\ref{prop:dotXbound}, and since $L_c(z)$ and 
%$\overline{L}_c^\Gamma(z)$ do not involve transition points, 
%this proves the invariance of $A^\nab_\Gamma(z)$.  
The analyticity of $A^\nab_\Gamma(z)$ in $K_J^\delta$ is a consequence
of the jump condition satisfied by $\dot{\mat{X}}(z)$ in the void
$\Gamma$; using $+$ ($-$) to denote boundary values taken on the real
axis from above (below), we have for real $z\in K_J^\delta$
\begin{equation}
\begin{array}{rcl}
A^\nab_{\Gamma+}(z)&=&e^{-iN\theta_\Gamma/2}\dot{X}_{11+}(z)e^{\kappa g_+(z)}\\\\
&=&e^{iN\theta_\Gamma/2}\dot{X}_{11-}(z)e^{i\phi_\Gamma}
e^{\kappa g_+(z)}\\\\&=&e^{iN\theta_\Gamma/2}\dot{X}_{11-}(z)e^{\kappa g_-(z)}\\\\&=&A^\nab_{\Gamma-}(z)\,,
\end{array}
\end{equation}
since by definition for $z$ in a void $\Gamma$, $i\phi_\Gamma=\kappa
g_-(z)-\kappa g_+(z)$.  Finally, the reality of $A^\nab_\Gamma(z)$
when $z$ is real and the information concerning its possible zero
follow from Proposition~\ref{prop:dotXsymmetry}.

\subsection{Asymptotic behavior of $\pi_{N,k}(z)$ for $z$ near a
saturated region of $[a,b]$.  }
\label{sec:saturated}
\subsubsection{Asymptotics valid away from hard edges.
Proof of Theorem~\ref{theorem:upper}.}
\label{sec:saturated:away}
Because the closed interval $J\subset\Gamma$ is bounded away from all of
the points $a,\alpha_0,\beta_0,\dots,\alpha_G,\beta_G,b$, we may fix the
parameter $\epsilon>0$ sufficiently small that Proposition~\ref{prop:E} 
controls $\mat{E}(z)-\mathbb{I}$,
and then select $\delta>0$ small enough that $K_J^\delta\cap\Sigma_F=\emptyset$
where the compact set $K_J^\delta$ is defined by \eqref{eq:fleshedJ}.
%
%Let $J\subset(a,b)$ be a closed interval throughout which the
%equilibrium measure achieves the upper constraint and the variational
%inequality (\ref{eq:saturatedregioninequality}) holds strictly.  As in
%\S~\ref{sec:lower}, this excludes the possibility of either endpoint
%of $J$ being a band edge.  Unlike in \S~\ref{sec:lower} however, we
%allow neither $a$ nor $b$ to be an endpoint of $J$; this particular
%restriction will be removed below in \S~\ref{sec:saturated:hardedges}.
%Consider the associated compact set $K_J^\epsilon$ defined by
%(\ref{eq:Kz1z2eps}), and fix the parameter $\epsilon>0$ sufficiently
%small that Proposition~\ref{prop:E} holds,
%$K_J^\epsilon\cap\Sigma_F=\emptyset$, and
%\begin{equation}
%\min J > a+h\epsilon/2\hspace{0.2 in} 
%\mbox{and}\hspace{0.2 in} \max J <b-h\epsilon/2\,.
%\label{eq:Jawayab}
%\end{equation}
%The condition (\ref{eq:Jawayab}) is necessary to ensure that we 
%may write the following exact formula for $\pi_{N,k}(z)$:
For $z\in K_J^\delta$ we thus have the following exact formula:
\begin{equation}
\begin{array}{l}
\displaystyle\pi_{N,k}(z)=
\left[E_{11}(z)\dot{X}_{11}(z)+E_{12}(z)\dot{X}_{21}(z)\right]
e^{(k-\#\del)g(z)} \prod_{n\in\del}(z-x_{N,n}) \\\\
\displaystyle +\,\,
\left[E_{11}(z)\dot{X}_{12}(z)+E_{12}(z)\dot{X}_{22}(z)\right]
i{\rm sgn}(\Im(z)) e^{\eta(z)-\gamma}e^{N(V(z)-\ell_c-i{\rm sgn}(\Im(z))\theta^0(z)/2)}  e^{(\#\del-k)g(z)}\prod_{n\in\nab}(z-x_{N,n})\,.
\end{array}
\end{equation}
This formula follows from \eqref{eq:XD}, 
\eqref{eq:Dnexttosat},
\eqref{eq:XfromXdotlast}, \eqref{eq:Rformula}, and 
Proposition~\ref{prop:solnrhp},  
and the right-hand side extends analytically to the whole set $K_J^\delta$.
Using the definition (\ref{eq:TBfuncdef}) of the function $T_\del(z)$,
and its characterization for nonreal $z$ in
Proposition~\ref{prop:tanal}, and recalling the definition
(\ref{eq:Lcdef}), we can rewrite this formula as
\begin{equation}
\begin{array}{l}
\displaystyle \pi_{N,k}(z)=\displaystyle\left[\exp\left(-N\int_{\Sigma_0^\nab}\log(z-x)\rho^0(x)\,dx\right)\prod_{n\in\nab}(z-x_{N,n})\right]
e^{NL_c(z)-iN{\rm sgn}(\Im(z))\theta^0(z)/2}\\\\
\displaystyle\hspace{0.2 in}\times\,\,\,\Bigg(
\left[E_{11}(z)\dot{X}_{11}(z)e^{\kappa g(z)}+E_{12}(z)\dot{X}_{21}(z)e^{\kappa g(z)}\right]T_\del(z)^{-1}2\cos\left(\frac{N\theta^0(z)}{2}\right)\\\\
\displaystyle\hspace{0.3 in}+\,\,\,
\left[E_{11}(z)\dot{X}_{12}(z)e^{-\kappa g(z)}+E_{12}(z)\dot{X}_{22}(z)e^{-\kappa g(z)}\right]i{\rm sgn}(\Im(z))e^{\eta(z)-\gamma}e^{-N\xi_\Gamma(z)}
e^{iN{\rm sgn}(\Im(z))(\theta_\Gamma-\theta^0(z)/2)}\Bigg)\,.
\end{array}
\label{eq:satfinalexact}
\end{equation}
Since $K_J^\delta$ is bounded away from any nodes $x_{N,n}$ for
$n\in\nab$, the product on the first line of (\ref{eq:satfinalexact})
may be approximated in terms of an exponential of an integral up to a
relative error of order $1/N$ uniform in $K_J^\delta$.  From
Proposition~\ref{prop:tasymp} it follows that $T_\del(z)^{-1}-1$ is
also of order $1/N$ uniformly in $K_J^\delta$.  Finally, using
Proposition~\ref{prop:E} to estimate $\mat{E}(z)-\mathbb{I}$ and
Proposition~\ref{prop:dotXbound} to bound $\dot{\mat{X}}(z)$ uniformly
in $K_J^\delta$, we see that
\begin{equation}
\begin{array}{rcl}
\pi_{N,k}(z)&=&\displaystyle e^{NL_c(z)-iN{\rm sgn}(\Im(z))\theta^0(z)/2}\\\\
&&\displaystyle \times\,\,\,\left(\left(\dot{X}_{11}(z)e^{\kappa g(z)}
+
O\left(\frac{1}{N}\right)\right)
2\cos\left(\frac{N\theta^0(z)}{2}\right)+O\left(\exp\left(-N\inf_{z\in
K_J^\delta}\Re(\xi_\Gamma(z))\right)\right)\right)\,.
\end{array}
\label{eq:sattheoremprerequisite}
\end{equation}
The exponential estimate holds because ${\rm
sgn}(\Re(i\theta^0(z)))={\rm sgn}(\Im(z))$ for all $z\in
K_J^\delta$.  Also, since $\Re(\xi_\Gamma(z))$ is strictly positive
for all $z\in K_J^\delta$ (this is equivalent to the inequality
(\ref{eq:saturatedregioninequality}) being strict in $J$ and
$\delta$ being sufficiently small), this term is exponentially small
as $N\rightarrow\infty$.  We note that $W(z)=\dot{X}_{11}(z)e^{\kappa g(z)}$.  

The estimates (\ref{eq:sattheoremestimates}) follow from
(\ref{eq:sattheoremprerequisite}) because
$e^{N(L_c(z)-\overline{L}_c^\Gamma(z)-i{\rm
sgn}(\Im(z))\theta^0(z)/2)}$ is uniformly bounded in $K_J^\delta$.
Indeed, we have
\begin{equation}
L_c(z)-\overline{L}_c^\Gamma(z)-\frac{i}{2}{\rm sgn}(\Im(z))\theta^0(z) =
-\frac{i\theta_\Gamma}{2}\cdot{\rm sgn}(\Im(z))\,.
%\left[-\pi\int_{y<x\in\Sigma_0^\nab}\rho^0(x)\,dx-
%\frac{1}{2}\theta(x_\Gamma)\right]\,,
\end{equation}
%where $x_\Gamma$ is any point in the saturated region $\Gamma$, and
%where $y$ is the nearest transition point to the right of 
%$\Gamma$ (or $y=b$ if there is no transition point to
%the right of $\Gamma$).  By the quantization condition
%(\ref{eq:yquantize}) the integral is a rational number $M/N$ with
%denominator $N$.  
The rest of the proof follows that of Theorem~\ref{theorem:lower}.

\subsubsection{Asymptotics uniformly valid near hard edges.  Proof of Theorem~\ref{theorem:hardedges}.}
\label{sec:saturated:hardedges}
We will analyze the case where the saturated region is $\Gamma=(a,\alpha_0)$
and $J=[a,t]$ with $t\in\Gamma$ in detail.  The analysis in a saturated
region near $z=b$ is similar.

The upper constraint is active throughout $J$, and the
variational inequality (\ref{eq:saturatedregioninequality}) holds
strictly for all $z\in J$.  
%Let $J$ be a closed
%interval throughout which (\ref{eq:saturatedregioninequality}) holds
%with $a\le\min J <a+h\epsilon/2$ and $\max J>a$, and consider the
%compact set $K_J^\epsilon$ defined by (\ref{eq:Kz1z2eps}).  
We take the fixed parameter $\epsilon$ to be sufficiently small that
Proposition~\ref{prop:E} controls $\mat{E}(z)-\mathbb{I}$, and then
choose $\delta>0$ small enough that $K_J^\delta\cap\Sigma_F=\emptyset$
where $K_J^\delta$ is defined by
\eqref{eq:fleshedJ}.
%$K_J^\epsilon\cap\Sigma_F=\emptyset$ and that Proposition~\ref{prop:E}
%holds.  
The set $K_J^\delta$ is the closure of the union of two open sets:
$K_{J,{\rm out}}^\delta$ consisting of the points in the interior of
$K_J^\delta$ with $\Re(z)<a$ and $K_{J,{\rm in}}^\delta$ consisting of
the points in the interior of $K_J^\delta$ with $\Re(z)>a$.

For $z\in K_{J,{\rm in}}^\delta$, the exact formula
(\ref{eq:satfinalexact}) for $\pi_{N,k}(z)$ is valid.  Since
$K_{J,{\rm in}}^\delta$ is not bounded away from $z=a$, we may no
longer neglect $T_\del(z)^{-1}-1$.  However, we may substitute from
Proposition~\ref{prop:tasymp} an asymptotic formula for
$T_\del(z)^{-1}$ that is uniformly valid in $K_{J,{\rm in}}^\delta$.
The remaining approximations we make for $z\in K_{J,{\rm
in}}^\delta$ are exactly the same as in the proof of
Theorem~\ref{theorem:upper}.

On the other hand, for $z\in K_{J,{\rm out}}^\delta$, the exact
formula (\ref{eq:piNkexactouter}) holds.  Using (\ref{eq:Yfuncdef}),
we may write this in the form
\begin{equation}
\begin{array}{rcl}
\pi_{N,k}(z)&=&\displaystyle
\left[E_{11}(z)\dot{X}_{11}(z)+E_{12}(z)\dot{X}_{21}(z)\right]
e^{(k-\#\del)g(z)}Y(z)\\\\
&&\displaystyle\times\,\,\,
\left[\exp\left(-N\int_{\Sigma_0^\nab}\log(z-x)\rho^0(x)\,dx\right)
\prod_{n\in\nab}(z-x_{N,n})\right]
\exp\left(N\int_{\Sigma_0^\del}\log(z-x)\rho^0(x)\,dx\right)\,.
\end{array}
\end{equation}
The terms in the large square brackets may be estimated using the
midpoint rule to approximate the integral in the exponent; these terms
are thus of the form $1+O(1/N)$ uniformly for $z\in K_{J,{\rm
out}}^\delta$.  From Proposition~\ref{prop:tasymp} we may substitute
an asymptotic formula for $Y(z)$ that is uniformly valid in $K_{J,{\rm
out}}^\delta$.  Using Proposition~\ref{prop:E} to estimate
$\mat{E}(z)-\mathbb{I}$ uniformly for $z\in K_{J,{\rm out}}^\delta$
and Proposition~\ref{prop:dotXbound} to uniformly bound
$\dot{X}_{11}(z)$ and $\dot{X}_{12}(z)$ in the same region, we obtain
an asymptotic expression for $\pi_{N,k}(z)$ that is valid in
$K_{J,{\rm out}}^\delta$.  To write this expression, we note that
$W(z)=\dot{X}_{11}(z)e^{\kappa g(z)}$.

The two asymptotic formulae so-obtained are uniformly valid right up
to the line $\Re(z)=a$ that divides $K_J^\delta$ into two parts.
Moreover, it is an exercise to check that the formulae agree for
$\Re(z)=a$.  In this way, we obtain a uniform approximation for
$\pi_{N,k}(z)$ for $z$ near $a$ that is an analytic function, and the
proof is complete.

\subsubsection{Asymptotics of zeros of $\pi_{N,k}(z)$ in saturated regions.
Proof of Theorem~\ref{theorem:exponential}.}
\label{sec:saturated:zeros}
The zeros of the cosine function in (\ref{eq:piNkupperasymp}) and
(\ref{eq:piNkhardedge}) are exactly the nodes of orthogonalization
making up the set $X_N$. We may thus expect that there should be a
zero of $\pi_{N,k}(z)$ very close to each node $x_{N,n}$ in a
saturated region.  To make this precise, we now study how the zeros of
the leading term in (\ref{eq:piNkupperasymp}) are perturbed by the
term $\delta_N(z)$.  Neither $\varepsilon_N(z)$ nor $\delta_N(z)$ in
(\ref{eq:piNkupperasymp}) are purely real for real $z$ (although the
imaginary part of $\varepsilon_N(z)$ is necessarily exponentially
small for real $z$ to balance with that of $\delta_N(z)$ since
$\pi_{N,k}(z)$ is a real polynomial).  However, from
(\ref{eq:satfinalexact}) we get the exact formula
\begin{equation}
\Re(\delta_N(z)) = \left(B_\Gamma^\del(z)\sin\left(\frac{N\theta^0(z)}{2}\right) + \sigma_N(z)\right)e^{\eta(z)-\gamma-N\xi_\Gamma(z)}
\end{equation}
where $\sigma_N(z)$ is uniformly of order $1/N$ for
$z\in\mathbb{R}\cap K_J^\delta$ with $K_J^\delta$ defined by
\eqref{eq:fleshedJ} for $\delta$ small enough, and
\begin{equation}
B_\Gamma^\del(z):=e^{N(L_c(z)-\overline{L}_c^\Gamma(z))}e^{-iN{\rm sgn}(\Im(z))\theta^0(z)/2}\dot{X}_{12}(z)e^{-\kappa g(z)}e^{iN{\rm sgn}(\Im(z))\theta_\Gamma}\,.
\label{eq:Bdelgamma}
\end{equation}
Note that (\ref{eq:Bdelgamma}) apparently defines $B_\Gamma^\del(z)$
for $\Im(z)\neq 0$, but it is easy to check that this definition
extends analytically to a real function for real $z$.

Now, if the saturated region is the interval $\Gamma=(a,\alpha_0)$,
then Proposition~\ref{prop:dotXsymmetry} guarantees that
$A_\Gamma^\del(z)$ and $B_\Gamma^\del(z)$ are bounded away from zero
and have opposite signs.  Since $\theta^0(z)$ is a strictly decreasing
function of $z$ for $z\in (a,b)$, it follows that there is a zero of
$\pi_{N,k}(z)$ exponentially close to but strictly greater than each
node $x_{N,n}$ in the interval $J$ (and no other zeros in $J$).
Similarly, if the saturated region is the interval
$\Gamma=(\beta_G,b)$, then Proposition~\ref{prop:dotXsymmetry}
guarantees that $A_\Gamma^\del(z)$ and $B_\Gamma^\del(z)$ are bounded
away from zero and have the same sign.  From this it follows that
there is a zero of $\pi_{N,k}(z)$ exponentially close to but strictly
less than each node $x_{N,n}$ in the interval $J$, and no other zeros
in $J$.  Note that with the use of the asymptotic formulae given in
Theorem~\ref{theorem:hardedges}, it follows that these conclusions
even hold true if the interval $J$ under consideration contains either
$z=a$ or $z=b$ as an endpoint.

If the saturated region is $\Gamma=\Gamma_j=(\beta_{j-1},\alpha_j)$
for some $j=1,\dots,G$, then Proposition~\ref{prop:dotXsymmetry}
implies that the product $A_\Gamma^\del(z)B_\Gamma^\del(z)$ vanishes
at exactly one point $z=z_j$ in $\Gamma=\Gamma_j$.  If $z_j<\min J$
then $A_\Gamma^\del(z)$ and $B_\Gamma^\del(z)$ are bounded away from
zero and have opposite signs for $z\in J$, and thus there is a zero of
$\pi_{N,k}(z)$ exponentially close to but strictly greater than each
node $x_{N,n}$ in $J$, and no other zeros in $J$.  If $z_j>\max J$
then $A_\Gamma^\del(z)$ and $B_\Gamma^\del(z)$ have the same sign and
thus there is a zero of $\pi_{N,k}(z)$ exponentially close to but
strictly less than each node $x_{N,n}$ in $J$, and no other zeros in
$J$.

Continuing with the case $\Gamma=\Gamma_j=(\beta_{j-1},\alpha_j)$,
suppose that $z_j$ lies in the interior of $J\subset\Gamma_j$.  If it
is $B_\Gamma^\del(z)$ that vanishes at $z=z_j$, then it is clear that
$\pi_{N,k}(z)$ has a zero exponentially close to each node $x_{N,n}$
in $J$ and no other zeros in $J$.  Moreover, in this case there is a
neighborhood of $z_j$ of length proportional to $1/N$ outside of which
$\Re(\delta_N(z))$ has the same sign as its leading term; thus with
the possible exception of a bounded number of nodes surrounding
$z=z_j$ the zeros exponentially localized near the nodes lying to the
left (right) of $z=z_j$ lie to the left (right) of the nearest node.
In fact, Proposition~\ref{prop:confine} guarantees that this situation
persists inward from the left and right to a single interval between
two consecutive nodes $[x_{N,m},x_{N,m+1}]$ that contains no zeros of
$\pi_{N,k}(z)$ at all, and such that there is a zero exponentially
close to but to the left of $x_{N,m}$ and another zero exponentially
close to but to the right of $x_{N,m+1}$.  Thus in this situation, the
interval $J$ contains precisely one less than the maximum possible
number of zeros of $\pi_{N,k}(z)$ since there is exactly one
consecutive pair of nodes that do not have any zero between them.

On the other hand if it is $A_\Gamma^\del(z)$ that vanishes at $z=z_j$
in the interior of $J$, then in addition to the zeros of the cosine
function, there is a single zero of
$A_\Gamma^\del(z)+\Re(\varepsilon_N(z))$, say $z=z_{j,N}$, that is
subjected to perturbation.  The zeros of the cosine lying to the left
(right) of $z=z_{j,N}$ are easily seen (using
Proposition~\ref{prop:dotXsymmetry} to analyze the relative signs of
$A_\Gamma^\del(z)$ and $B_\Gamma^\del(z)$) to move under perturbation
an exponentially small amount to the left (right).  The ``spurious''
zero $z_{j,N}$ is also perturbed an exponentially small amount, and it
is easy to see that the closer $z_{j,N}$ lies to a node in $X_N$, the
more it is repelled by the perturbation.  Even in the degenerate case
that $z_{j,N}$ coincides exactly with a node, it is easy to see that
the perturbation always serves to unfold the double zero into two real
zeros of $\pi_{N,k}(z)$ both exponentially close to the same node,
with one on either side.  Thus in this situation, the interval $J$
always contains precisely the maximum possible number of zeros of
$\pi_{N,k}(z)$ (one zero between each consecutive pair of nodes), all
exponentially localized to nodes in $X_N$ with the possible exception
of exactly one, which necessarily corresponds to the zero $z_{j,N}$ of
$A_\Gamma^\del(z)+\Re(\varepsilon_N(z))$.  This completes the proof.

\subsection{Asymptotic behavior of $\pi_{N,k}(z)$ for $z$ near a band.}
\label{sec:band}
\subsubsection{Proof of Theorem~\ref{theorem:band}.}
The closed interval $J$ is necessarily bounded away from the two
nearest band edge points $z=\alpha_j$ and $z=\beta_j$.  Therefore,
given $\epsilon>0$ sufficiently small that Proposition~\ref{prop:E}
controls $\mat{E}(z)-\mathbb{I}$, we may choose $\delta>0$ small
enough that the set $K_J^\delta$ defined by \eqref{eq:fleshedJ}
satisfies $K_J^\delta\cap\Sigma_F=\emptyset$.

Suppose first that the band $I$ containing $J$ is not a transition
band, but rather is completely contained in $\Sigma_0^\nab$. Then,
from Proposition~\ref{prop:solnrhp}, \eqref{eq:XD}, \eqref{eq:Xblueupper},
\eqref{eq:Xbluelower}, \eqref{eq:XfromXdotlast}, and \eqref{eq:Rformula}, we
have the following exact formula for $\pi_{N,k}(z)$ in $K_J^\delta$:
\begin{equation}
\begin{array}{rcl}
\pi_{N,k}(z)&=&\displaystyle T_\nab(z)^{-1/2}
\left[\exp\left(-N\int_{\Sigma_0^\del}\log(z-x)\rho^0(x)\,dx\right)
\prod_{n\in\del}(z-x_{N,n})\right]e^{N\overline{L}_c^I(z)}(-1)^{M_I^\nab}\\\\
&&\displaystyle\times\,\,\,\Bigg[\left(E_{11}(z)\dot{X}_{11}(z)+E_{12}(z)\dot{X}_{21}(z)\right)
e^{\kappa g(z)}e^{-iN{\rm sgn}(\Im(z))\theta_I^\nab(z)/2} \\\\
&&\displaystyle\hspace{0.2 in}+\,\,\, i{\rm
sgn}(\Im(z))e^{\eta(z)-\gamma}\left(E_{11}(z)\dot{X}_{12}(z)+E_{12}(z)\dot{X}_{22}(z)\right)
e^{-\kappa g(z)}e^{iN{\rm sgn}(\Im(z))\theta_I^\nab(z)/2}\Bigg]\,,
\end{array}
\label{eq:bandnabexactfirst}
\end{equation}
where $\overline{L}_c^I(z)$ is defined by \eqref{eq:LcbarI}, and
\begin{equation}
M_I^\nab:=N\int_{y<x\in\Sigma_0^\del}\rho^0(x)\,dx
\end{equation}
where $y$ is the nearest transition point to the right of $J\subset
I$.  It follows from (\ref{eq:yquantize}) that
$M_I^\nab\in\mathbb{Z}$.  The right-hand side of
\eqref{eq:bandnabexactfirst} extends analytically to the whole compact
set $K_J^\delta$.  The terms in square brackets on the first line of
(\ref{eq:bandnabexactfirst}) are seen to be $1+O(1/N)$ uniformly for
$z\in K_J^\delta$ by a midpoint rule approximation argument (since
$K_J^\delta$ is in this case bounded away from any component of
$\Sigma_0^\del$).  Similarly, $T_\nab(z)^{-1/2}=1+O(1/N)$ uniformly
for $z\in K_J^\delta$ by Proposition~\ref{prop:tasymp}.
Proposition~\ref{prop:E} and Proposition~\ref{prop:dotXbound} then
imply that the terms in parentheses on the second line of
(\ref{eq:bandnabexactfirst}) are simply $\dot{X}_{11}(z) + O(1/N)$ and
that the terms in parentheses on the third line of
(\ref{eq:bandnabexactfirst}) are just $\dot{X}_{12}(z)+O(1/N)$, with
all errors uniform in $K_J^\delta$.  Thus one obtains an asymptotic
formula for $\pi_{N,k}(z)$ valid uniformly in $K_J^\delta$.

Next, suppose that the band $I$ containing $J$ is not a transition
band, but is rather completely contained in $\Sigma_0^\del$. In this case, 
from Proposition~\ref{prop:solnrhp}, \eqref{eq:XD}, 
\eqref{eq:Xyellowupper}, \eqref{eq:Xyellowlower}, \eqref{eq:XfromXdotlast},
and \eqref{eq:Rformula}, we have the following exact formula for
$\pi_{N,k}(z)$ in $K_J^\delta$:
\begin{equation}
\begin{array}{rcl}
\pi_{N,k}(z)&=&\displaystyle T_\del(z)^{-1/2}\left[
\exp\left(-N\int_{\Sigma_0^\nab}\log(z-x)\rho^0(x)\,dx\right)
\prod_{n\in\nab}(z-x_{N,n})\right] e^{N\overline{L}_c^I(z)}(-1)^{M_I^\del}\\\\
&&\displaystyle\times\,\,\,\Bigg[\left(E_{11}(z)\dot{X}_{11}(z)+E_{12}(z)\dot{X}_{21}(z)\right)e^{\kappa g(z)}e^{iN{\rm sgn}(\Im(z))[\theta^0(z)-\theta_I^\del(z)]/2}\\\\
&&\displaystyle\hspace{0.2 in}+\,\,\, i{\rm sgn}(\Im(z))e^{\eta(z)-\gamma}
\left(E_{11}(z)\dot{X}_{12}(z)+E_{12}(z)\dot{X}_{22}(z)\right)e^{\kappa g(z)}e^{-iN{\rm sgn}(\Im(z))[\theta^0(z)-\theta_I^\del(z)]/2}\Bigg]\,,
\end{array}
\label{eq:banddelexactfirst}
\end{equation}
where
\begin{equation}
M_I^\del:=N\int_{y<x\in\Sigma_0^\nab}\rho^0(x)\,dx\in\mathbb{Z}
\end{equation}
and $y$ is the nearest transition point to the right of $J\subset I$.
Once again, the right-hand side may be considered as an analytic
function in the set $K_J^\delta$.  Since $K_J^\delta$ is bounded away
from $\Sigma_0^\nab$ in this case, the expression
(\ref{eq:banddelexactfirst}) may be approximated in virtually the same
way as (\ref{eq:bandnabexactfirst}) in order to obtain a uniformly
valid asymptotic formula for $\pi_{N,k}(z)$.

Finally, suppose that the band $I$ containing $J$ is a transition
band, in which we must place a transition point $y\in Y_N$.  Recall
that $J$ is bounded away from the endpoints $\alpha_j$ and $\beta_j$
of $I=I_j$.  Thus, without any loss of generality, we may choose the
transition point $y\in I\cap Y_N$ such that either $y<\min J$ or
$y>\max J$.  This means that either $J\subset \Sigma_0^\nab$ or $J\subset
\Sigma_0^\del$, and we may analyze either \eqref{eq:bandnabexactfirst}
or \eqref{eq:banddelexactfirst} respectively, exactly as we have done above.

We now wish to write
the two exact formulae (\ref{eq:bandnabexactfirst}) and
(\ref{eq:banddelexactfirst}) in such a form that it is clear that the
limit $N\rightarrow\infty$ yields an  asymptotic formula that is
independent of whether $J\subset\Sigma_0^\nab$ or $J\subset\Sigma_0^\del$.
In fact, a direct calculation using (\ref{eq:thetanabI}) and
(\ref{eq:thetadelI}) along with the quantization condition
(\ref{eq:yquantize}) shows that
\begin{equation}
(-1)^{M_I^\nab}e^{\mp iN\theta_I^\nab(z)/2} =
(-1)^{M_I^\del}e^{ \pm iN[\theta^0(z)-\theta_I^\del(z)]/2} =
\exp\left(\pm i\pi Nc\left[\mu_{\rm min}^c([x,b])-\int_x^z\psi_I(s)\,ds\right]\right)
\end{equation}
where $x$ is any point or endpoint of the band $I$.  

Therefore, in considering the limit $N\rightarrow\infty$, it remains
to recall that Proposition~\ref{prop:dotXsymmetry} implies that
$\dot{X}_{11+}(z)$ does not vanish at any point of $I$, and that
$W(z)=\dot{X}_{11}(z)e^{\kappa g(z)}$.  This completes the proof.

\subsubsection{Asymptotic behavior of the zeros.  
Proof of Theorem~\ref{theorem:bandzeros}.}
Theorem~\ref{theorem:bandzeros} is a consequence of the estimate
(\ref{eq:bandestimate}) established in Theorem~\ref{theorem:band}, the
strict inequalities $0<d\mu_{\rm min}^c/dx<\rho^0(x)/c$ holding by definition
for $J\subset I$ because $I$ is a band, and from the strict inequality 
$A_I(x)>0$ for $x\in J\subset\mathbb{R}$ stated in Theorem~\ref{theorem:band}.

\subsection{Asymptotic behavior of $\pi_{N,k}(z)$ for $z$ near a band edge.}
\subsubsection{Band/void edges.  
Proof of Theorem~\ref{theorem:Airynab}.}  First consider a left band
endpoint $z=\alpha$ between a band $I$ (on the right) and a void
$\Gamma$ (on the left).  We take the contour parameter $\epsilon$
sufficiently small that Proposition~\ref{prop:E} controls
$\mat{E}(z)-\mathbb{I}$, and then choose $r>0$ sufficiently small that
the disc $|z-\alpha|\le r$ is contained in the disc
$D_\Gamma^{\nab,L}$.  For such $z$, we thus have the following exact
formula for $\pi_{N,k}(z)$:
\begin{equation}
\begin{array}{rcl}
\pi_{N,k}(z)&=&
\displaystyle -\sqrt{2\pi}
e^{(\eta(z)-\gamma)/2}e^{N\overline{L}_c^I(z)}\\\\
&&\displaystyle\times\,\,\,T_\nab(z)^{-1/2}
\left[\exp\left(-N\int_{\Sigma_0^\del}\log(z-x)\rho^0(x)\,dx\right)
\prod_{n\in\del}(z-x_{N,n})\right]\\\\
%&&\displaystyle\times\,\,\,(-1)^{M_\Gamma^{\nab,L}}\\\\
&&\displaystyle\times\,\,\,\Bigg[\left(\frac{3}{4}\right)^{1/6}
\left(E_{11}(z)H_{\Gamma,12}^{\nab,L}(z)+E_{12}(z)H_{\Gamma,22}^{\nab,L}(z)
\right)
Ai\left(-\left(\frac{3}{4}\right)^{2/3}\tau_\Gamma^{\nab,L}(z)\right) \\\\
&&\displaystyle\hspace{0.1 in}+\,\,\,
\left(\frac{3}{4}\right)^{-1/6}\left(E_{11}(z)H_{\Gamma,11}^{\nab,L}(z)+
E_{12}(z)H_{\Gamma,21}^{\nab,L}(z)\right)
Ai'\left(-\left(\frac{3}{4}\right)^{2/3}
\tau_\Gamma^{\nab,L}(z)\right)\Bigg]\,.
\end{array}
\label{eq:AirynabLexact}
\end{equation}
%where
%\begin{equation}
%M_\Gamma^{\nab,L}:=N\int_{y<x\in\Sigma_0^\del}\rho^0(x)\,dx
%\end{equation}
%and $y$ is the nearest transition point to the right of the band edge
%under consideration.  The quantization condition (\ref{eq:yquantize})
%on the transition points implies that $M_\Gamma^{\nab,L}$ is an
%integer.  
This follows from \eqref{eq:XD}, \eqref{eq:Dnexttovoid}, \eqref{eq:Xblueupper},
\eqref{eq:Xbluelower}, \eqref{eq;XhatnabL}, \eqref{eq:XfromXdotnabL},
\eqref{eq:Rformula}, and Proposition~\ref{prop:solnrhp}.

Recall from \S~\ref{sec:Airy} that $\mat{H}_\Gamma^{\nab,L}(z)$ is
analytic throughout $D_\Gamma^{\nab,L}$.  From the definition of this
function in terms of $\tau_\Gamma^{\nab,L}(z)$ and $\dot{\mat{X}}(z)$
it follows that the first column (second column) of
$\mat{H}_{\Gamma}^{\nab,L}(z)$ is uniformly bounded in
$D_\Gamma^{\nab,L}$ by a quantity of order $N^{-1/6}$ (of order
$N^{1/6}$).  Also, Proposition~\ref{prop:dotXsymmetry} implies that
the matrix elements of $\mat{H}_\Gamma^{\nab,L}(z)$ are real for real
$z$.
%Finally, since the definition of $\mat{H}_{\Gamma}^{\nab,L}(z)$
%involves $e^{iN\theta_\Gamma/2}$, the combination
%$(-1)^{M_\Gamma^{\nab,L}}\mat{H}_\Gamma^{\nab,L}(z)$ turns out not to
%depend at all on the choice of transition points $Y_N$.  
We may now use an argument based on the midpoint rule for Riemann
sums to approximate the terms in square brackets on the second line of
(\ref{eq:AirynabLexact}), recall Proposition~\ref{prop:tasymp} to
handle $T_\del(z)$ and use Proposition~\ref{prop:E} to estimate
$\mat{E}(z)-\mathbb{I}$.

Finally, we may observe from \S~\ref{sec:Airy} the relations
\begin{equation}
\begin{array}{rcl}
\displaystyle H_{\Gamma,12}^{\nab,L}(z)&=&\displaystyle -H_\Gamma^-(z)
(-\tau_\Gamma^{\nab,L}(z))^{1/4}\,,\\\\
\displaystyle H_{\Gamma,11}^{\nab,L}(z)&=&\displaystyle H_\Gamma^+(z)
(-\tau_\Gamma^{\nab,L}(z))^{-1/4}\,,
\end{array}
\end{equation}
where we have used the identities $W(z)=\dot{X}_{11}(z)e^{\kappa
g(z)}$ and $Z(z)=\dot{X}_{12}(z)e^{-\kappa g(z)}$, and the functions
$H_\Gamma^\pm(z)$ are defined by \eqref{eq:HGamma}.  This completes
the proof of the asymptotic formula \eqref{eq:AirynabLuniform} and the
corresponding error estimates.

Since $\tau_\Gamma^{\nab,L}(z)$ is uniformly bounded independently of
$N$ for $z$ in shrinking neighborhoods of the band edge $z=\alpha$ with
radius of order $N^{-2/3}$, we immediately obtain 
the asymptotic formula \eqref{eq:AirynabLinner} and the corresponding
error estimate.

Next, consider a right band endpoint $z=\beta$ between a band $I$ (on
the left) and a void $\Gamma$ (on the right).  Again taking $\epsilon$
small enough that Proposition~\ref{prop:E} controls
$\mat{E}(z)-\mathbb{I}$, we take the parameter $r$ small enough that
the disc $|z-\beta|\le r$ is contained in the disc
$D_\Gamma^{\nab,R}$.  In this case, we have the exact formula:
\begin{equation}
\begin{array}{rcl}
\pi_{N,k}(z)&=&\displaystyle 
-i\sqrt{2\pi}e^{(\eta(z)-\gamma)/2}e^{N\overline{L}_c^I(z)}
\\\\
&&\displaystyle\times\,\,\,
T_\nab(z)^{-1/2}
\left[\exp\left(-N\int_{\Sigma_0^\del}\log(z-x)\rho^0(x)\,dx\right)
\prod_{n\in\del}(z-x_{N,n})\right]\\\\
&&\displaystyle\times\,\,\,\Bigg[\left(\frac{3}{4}\right)^{1/6}
\left(E_{11}(z)H_{\Gamma,12}^{\nab,R}(z)+E_{12}(z)H_{\Gamma,22}^{\nab,R}(z)
\right)Ai\left(-\left(\frac{3}{4}\right)^{2/3}\tau_\Gamma^{\nab,R}(z)\right) \\\\
&&\displaystyle\hspace{0.1 in}+\,\,\,
\left(\frac{3}{4}\right)^{-1/6}\left(E_{11}(z)H_{\Gamma,11}^{\nab,R}(z)+
E_{12}(z)H_{\Gamma,21}^{\nab,R}(z)\right)Ai'\left(-\left(\frac{3}{4}\right)^{2/3}\tau_\Gamma^{\nab,R}(z)\right)\Bigg]\,.
\end{array}
\label{eq:AirynabRexact}
\end{equation}
This follows from \eqref{eq:XD}, \eqref{eq:Dnexttovoid}, 
\eqref{eq:Xblueupper}, \eqref{eq:Xbluelower}, \eqref{eq;XhatnabR},
\eqref{eq:XfromXdotnabR}, \eqref{eq:Rformula}, and 
Proposition~\ref{prop:solnrhp}.
%where with $y$ denoting the nearest transition point to the 
%right of the band edge,
%\begin{equation}
%M_\Gamma^{\nab,R}:=N\int_{y<x\in\Sigma_0^\del}\rho^0(x)\,dx\in\mathbb{Z}\,.
%\end{equation}
Once again, we see that the second line in \eqref{eq:AirynabRexact}
may be replaced by $1+O(1/N)$ uniformly for $|z-\beta|\le r$. Since
the first column of $\mat{H}_\Gamma^{\nab,R}(z)$ is uniformly of order
$N^{-1/6}$ and the second column of $\mat{H}_\Gamma^{\nab,R}(z)$ is
uniformly of order $N^{1/6}$, and since we have the exact
representations (from \S~\ref{sec:Airy})
\begin{equation}
\begin{array}{rcl}
\displaystyle H_{\Gamma,12}^{\nab,R}(z)&= & \displaystyle
iH_\Gamma^+(z)(-\tau_\Gamma^{\nab,R}(z))^{1/4}\\\\
\displaystyle H_{\Gamma,11}^{\nab,R}(z)&=&\displaystyle
-iH_\Gamma^-(z)(-\tau_\Gamma^{\nab,R}(z))^{-1/4}\,,
\end{array}
\end{equation}
we immediately obtain the asymptotic formula
\eqref{eq:AirynabRuniform} and the corresponding error estimates 
with the use of Proposition~\ref{prop:E}.  The asymptotic formula
\eqref{eq:AirynabRinner} and its error estimate then follow exactly as
before, since $\tau_\Gamma^{\nab,R}(z)$ remains uniformly bounded as
$N\rightarrow\infty$ if $|z-\beta|\le rN^{-2/3}$.  Note that in this
case the matrix elements of $\mat{H}_\Gamma^{\nab,R}$ are imaginary
for real $z$.

\subsubsection{Band/saturated region edges.  
Proof of Theorem~\ref{theorem:Airydel}.}  First consider the
neighborhood of a left band edge $z=\alpha$ separating a band $I$ (for
$z>\alpha$) from a saturated region $\Gamma$ (for $z<\alpha$).  We
choose the contour parameter $\epsilon$ sufficiently small that
Proposition~\ref{prop:E} controls the matrix $\mat{E}(z)-\mathbb{I}$.
Then we choose $r>0$ small enough that the disc $|z-\alpha|\le r$ is
contained within the disc $D_\Gamma^{\del,L}$.  In this case, we have
the following exact formula for $\pi_{N,k}(z)$:
\begin{equation}
\begin{array}{rcl}
\pi_{N,k}(z)&=&\displaystyle i\sqrt{2\pi}e^{(\eta(z)-\gamma)/2}
e^{N\overline{L}^I_c(z)}\\\\
&&\displaystyle\times\,\,\,
T_\del(z)^{-1/2}\left[\exp\left(-N\int_{\Sigma_0^\nab}\log(z-x)\rho^0(x)\,dx
\right)\prod_{n\in\nab}(z-x_{N,n})\right]
\\\\
%&&\displaystyle\times\,\,\, 
%(-1)^{M_\Gamma^{\del,L}}\\\\
&&\displaystyle\times\,\,\,\Bigg[\left(\frac{3}{4}\right)^{1/6}
\left(E_{11}(z)H_{\Gamma,12}^{\del,L}(z)+E_{12}(z)H_{\Gamma,22}^{\del,L}(z)\right)F^L_A(z)\\\\
&&\displaystyle \hspace{0.1 in}+\,\,\,
\left(\frac{3}{4}\right)^{-1/6}\left(E_{11}(z)H_{\Gamma,11}^{\del,L}(z)+E_{12}(z)H_{\Gamma,21}^{\del,L}(z)\right)F^L_B(z)\Bigg]
\end{array}
\label{eq:AirydelLexact}
\end{equation}
where $F^L_A(z)$ and $F^L_B(z)$ are the combinations of trigonometric
functions and Airy functions and their derivatives defined by
\eqref{eq:FALBL}.
This formula follows from
\eqref{eq:XD}, \eqref{eq:Dnexttosat}, \eqref{eq:Xyellowupper},
\eqref{eq:Xyellowlower}, \eqref{eq;XhatdelL}, \eqref{eq:XfromXdotdelL},
\eqref{eq:Rformula}, and Proposition~\ref{prop:solnrhp}.
%\begin{equation}
%W_L(z):= \cos\left(\frac{N\theta^0(z)}{2}\right)
%Bi\left(-\left(\frac{3}{4}\right)^{2/3}\tau_\Gamma^{\del,L}(z)\right)-
%\sin\left(\frac{N\theta^0(z)}{2}\right)
%Ai\left(-\left(\frac{3}{4}\right)^{2/3}\tau_\Gamma^{\del,L}(z)\right)\,,
%\end{equation}
%and
%\begin{equation}
%X_L(z):= \cos\left(\frac{N\theta^0(z)}{2}\right)
%Bi'\left(-\left(\frac{3}{4}\right)^{2/3}\tau_\Gamma^{\del,L}(z)\right)-
%\sin\left(\frac{N\theta^0(z)}{2}\right)
%Ai'\left(-\left(\frac{3}{4}\right)^{2/3}\tau_\Gamma^{\del,L}(z)\right)\,.
%\end{equation}
%and with $y$ denoting the nearest transition point to the right 
%of the band edge,
%\begin{equation}
%M_{\Gamma}^{\del,L}:=N\int_{y<x\in\Sigma_0^\nab}\rho^0(x)\,dx\in\mathbb{Z}\,.
%\end{equation}
The terms on the second line of the right-hand side in
\eqref{eq:AirydelLexact} are $1+O(1/N)$ as
$N\rightarrow\infty$ uniformly for $|z-\alpha|\le r$, as can be seen
from a midpoint rule approximation of the integral, and by using
Proposition~\ref{prop:tasymp}.  Proposition~\ref{prop:E} is then used to
control $\mat{E}(z)-\mathbb{I}$. Noting that the second column of
$\mat{H}_\Gamma^{\del,L}(z)$ is uniformly of order $N^{1/6}$ while the
first column of $\mat{H}_\Gamma^{\del,L}(z)$ is uniformly of order
$N^{-1/6}$,  and moreover recalling from \S~\ref{sec:Airy} the explicit
formulae
\begin{equation}
\begin{array}{rcl}
\displaystyle H_{\Gamma,12}^{\del,L}(z)&=&\displaystyle -iH_\Gamma^-(z)
\left(-\tau_\Gamma^{\del,L}(z)\right)^{1/4}\,,\\\\
\displaystyle H_{\Gamma,11}^{\del,L}(z)&=&\displaystyle -iH_\Gamma^+(z)
\left(-\tau_\Gamma^{\del,L}(z)\right)^{-1/4}\,,
\end{array}
\end{equation}
which also rely on the identities $W(z)=\dot{X}_{11}(z)e^{\kappa g(z)}$
and $Z(z)=\dot{X}_{12}(z)e^{-\kappa g(z)}$, the asymptotic formula
\eqref{eq:AirydelLuniform} is obtained along with the corresponding error 
estimates.  The asymptotic formula \eqref{eq:AirydelLinner} then
follows along with its error estimate by noting that
$\tau_\Gamma^{\del,L}(z)$ remains uniformly bounded as
$N\rightarrow\infty$ if $|z-\alpha|\le rN^{-2/3}$.

Next consider the neighborhood of a right band edge $z=\beta$
separating a band $I$ (for $z<\beta$) from a saturated region $\Gamma$
(for $z>\beta$).  Again take the contour parameter $\epsilon$
sufficiently small that Proposition~\ref{prop:E} provides a uniform
estimate of $\mat{E}(z)-\mathbb{I}$ on appropriate closed sets, and
then choose $r>0$ small enough that the disc $|z-\beta|\le r$ ls
contained within the disc $D_\Gamma^{\del,R}$.  Then we have for $z$
with $|z-\beta|\le r$ the exact formula:
\begin{equation}
\begin{array}{rcl}
\pi_{N,k}(z)&=&\displaystyle 
\sqrt{2\pi}e^{(\eta(z)-\gamma)/2}e^{N\overline{L}^I_c(z)}\\\\
&&\displaystyle\times\,\,\,
T_\del(z)^{-1/2}\left[\exp\left(-N\int_{\Sigma_0^\nab}\log(z-x)\rho^0(x)\,dx
\right)\prod_{n\in\nab}(z-x_{N,n})\right]
\\\\
%&&\displaystyle\times\,\,\, 
%(-1)^{M_\Gamma^{\del,R}}\\\\
&&\displaystyle\times\,\,\,\Bigg[\left(\frac{3}{4}\right)^{1/6}
\left(E_{11}(z)H_{\Gamma,12}^{\del,R}(z)+E_{12}(z)H_{\Gamma,22}^{\del,R}(z)
\right)F_A^R(z)\\\\
&&\displaystyle \hspace{0.1 in}+\,\,\,
\left(\frac{3}{4}\right)^{-1/6}
\left(E_{11}(z)H_{\Gamma,11}^{\del,R}(z)+E_{12}(z)H_{\Gamma,21}^{\del,R}(z)
\right)F_B^R(z)\Bigg]
\end{array}
\label{eq:AirydelRexact}
\end{equation}
where $F_A^R(z)$ and $F_B^R(z)$ are the expressions defined by 
\eqref{eq:FARBR}.  
This formula follows from \eqref{eq:XD}, \eqref{eq:Dnexttosat},
\eqref{eq:Xyellowupper}, \eqref{eq:Xyellowlower}, \eqref{eq;XhatdelR},
\eqref{eq:XfromXdotdelR}, \eqref{eq:Rformula}, and 
Proposition~\ref{prop:solnrhp}.
%\begin{equation}
%W_R(z):= \cos\left(\frac{N\theta^0(z)}{2}\right)
%Bi\left(-\left(\frac{3}{4}\right)^{2/3}\tau_\Gamma^{\del,R}(z)\right)+
%\sin\left(\frac{N\theta^0(z)}{2}\right)
%Ai\left(-\left(\frac{3}{4}\right)^{2/3}\tau_\Gamma^{\del,R}(z)\right)\,,
%\end{equation}
%and
%\begin{equation}
%X_R(z):= \cos\left(\frac{N\theta^0(z)}{2}\right)
%Bi'\left(-\left(\frac{3}{4}\right)^{2/3}\tau_\Gamma^{\del,R}(z)\right)+
%\sin\left(\frac{N\theta^0(z)}{2}\right)
%Ai'\left(-\left(\frac{3}{4}\right)^{2/3}\tau_\Gamma^{\del,R}(z)\right)\,,
%\end{equation}
%and with $y$ denoting the nearest transition point to the right of the band 
%edge,
%\begin{equation}
%M_{\Gamma}^{\del,R}:=N\int_{y<x\in\Sigma_0^\nab}\rho^0(x)\,dx\in\mathbb{Z}\,.
%\end{equation}
%With similar approximations made in these formulae as in the case of
%band edges near the lower constraint, we have the following. 
%Restricting $z$ to shrinking neighborhoods of the band edge $z=z_0$
%with radius of order $N^{-2/3}$, we then find a simple asymptotic
%formula for $\pi_{N,k}(z)$ in terms of $Ai(\cdot)$ and $Bi(\cdot)$
%alone.
Once again, the terms on the second line of the right-hand side of
\eqref{eq:AirydelRexact} can be approximated uniformly for $|z-\beta|\le r$
as $1+O(1/N)$ as $N\rightarrow\infty$.  Proposition~\ref{prop:E} again
guarantees that uniformly for $|z-\beta|\le r$ we have
$\mat{E}(z)-\mathbb{I}=O(1/N)$, and then noting that
$\mat{H}_\Gamma^{\del,R}(z)N^{\sigma_3/6}$ remains uniformly bounded
as $N\rightarrow\infty$ and more specifically that
\begin{equation}
\begin{array}{rcl}
\displaystyle H_{\Gamma,12}^{\del,R}(z)&=&\displaystyle H_\Gamma^-(z)
\left(-\tau_\Gamma^{\del,R}(z)\right)^{1/4}\,,\\\\
\displaystyle H_{\Gamma,11}^{\del,R}(z)&=&\displaystyle -H_\Gamma^+(z)
\left(-\tau_\Gamma^{\del,R}(z)\right)^{-1/4}\,,
\end{array}
\end{equation}
we complete the proof of the asymptotic formula \eqref{eq:AirydelRuniform}
and its corresponding error estimates.  Since $\tau_\Gamma^{\del,R}(z)$ is
uniformly bounded as $N\rightarrow\infty$ with $|z-\beta|\le rN^{-2/3}$, 
we then obtain immediately the asymptotic formula \eqref{eq:AirydelRinner}
and its corresponding error estimate.

%The phenomenon of a band edge meeting a gap where the equilibrium
%measure achieves the upper constraint is a feature that is only
%present for orthogonal polynomials in the discrete weights setting;
%thus the resulting asymptotic formulae given here representing
%$\pi_{N,k}(z)$ in terms of both $Ai(\cdot)$ and $Bi(\cdot)$ describe
%phenomena characteristic of the discrete weights under study.

\section{Universality:  Proofs 
of Theorems Stated in \S~\ref{sec:universalityactualtheorems}}
\label{sec:universality}

\subsection{Relation between correlation functions of dual ensembles.}

Since the holes are also governed by a discrete orthogonal polynomial
ensemble, the correlation functions for holes are again represented as
determinants involving the reproducing kernel, this time
corresponding to the dual weights. It turns out that there is a simple
relation between the correlation functions for particles and those
for holes.

\subsubsection{Probabilistic approach.}

Let $\overline{R}^{(N,\bar{k})}_m$ \label{symbol:dualcorr}
be the $m$-point correlation
function of the dual orthogonal polynomial ensemble for the holes.
Hence $\overline{R}^{(N,\bar{k})}_m$ is defined as in
\eqref{eq:correldef} with the replacement of $p^{(N,k)}$ by
$\overline{p}^{(N,\bar{k})}$. Let $\overline{K}_{N,\bar{k}}$ 
\label{symbol:dualK} denote
the reproducing kernel of the dual ensemble.
Then \eqref{eq:determinantalformulacorrelations}
implies that
\begin{equation}
  \overline{R}_m^{(N,\bar{k})}(x_1,\dots, x_m)
= \det \bigl( \overline{K}_{N,\bar{k}} (x_i, x_j) \bigr)_{1\le
i,j\le m}\, ,
\end{equation}
for nodes $x_1, \dots, x_m$.
Now, given nodes $x_1,\dots, x_m$,
\begin{equation}
\begin{split}
  &\mathbb{P} ( \text{there are particles at each of the nodes $x_1,\dots,
  x_m$}) \\
  &\quad  = \mathbb{P} ( \text{there are no holes at any of the nodes
$x_1,\dots,
  x_m$}) \\
  &\quad = 1 - \sum_{i=1}^{m} \mathbb{P} (\text{there is a hole at the node
  $x_i$})
  + \sum_{1\le i<j\le m} \mathbb{P} (\text{there are holes at
  both of the nodes $x_i$ and $x_j$}) \\
  &\quad  \quad
  - \sum_{1\le i<j<k\le m} \mathbb{P} (\text{there are holes at each of the nodes
  $x_i, x_j, x_k$})  + \cdots.
\end{split}
\end{equation}
Thus from \eqref{eq:correldef},
\begin{equation}
\begin{split}
  R_m^{(N,k)} (x_1, \dots, x_m)
  = 1 - \sum_{i=1}^{m} \overline{R}_1^{(N,\bar{k})} (x_i)
  + \sum_{1\le i<j\le m} \overline{R}_2^{(N,\bar{k})} (x_i, x_j)
  - \sum_{1\le i<j<k\le m} \overline{R}_3^{(N,\bar{k})} (x_i, x_j, x_k)  + \cdots
\end{split}
\label{eq:correlationdual}
\end{equation}
Therefore the determinantal formula
(\ref{eq:determinantalformulacorrelations}) for the correlation
functions implies the following.

\begin{prop}\label{prop:dualkernel}
Let $K_{N,k}$ be the reproducing kernel \eqref{eq:reproducingdef}
for the discrete orthogonal polynomial ensemble, and let
$\overline{K}_{N,\bar{k}}$ be the reproducing kernel of the
corresponding dual orthogonal polynomial ensemble. Then with
$\bar{k}=N-k$,
\begin{equation}\label{eq:correldual1}
  \det \bigl( K_{N,k} (x_i, x_j) \bigr)_{1\le i,j\le m}
  = \det \bigl( \delta_{ij}- \overline{K}_{N,\bar{k}} (x_i, x_j) \bigr)_{1\le
i,j\le m}.
\end{equation}
\end{prop}

In particular, when $m=1$, this result implies that for a node $x\in X_N$,
\begin{equation}\label{eq:correldual2}
  K_{N,k}(x,x) = 1- \overline{K}_{N,\bar{k}} (x,x)
\end{equation}
and then when $m=2$ we further discover that for nodes $x\neq y$,
\begin{equation}
  K_{N,k}(x,y)^2 = \overline{K}_{N,\bar{k}} (x,y)^2\,.
\label{eq:correldual3squared}
\end{equation}

\subsubsection{Direct approach.}
It is possible to establish these same results, and also to refine
(\ref{eq:correldual3squared}) by determining the relative sign of
$K_{N,k}(x,y)$ and $\overline{K}_{N,\bar{k}}(x,y)$, by using
Proposition~\ref{prop:solnrhp} regarding the solution formula for
Interpolation Problem~\ref{rhp:DOP} and the dual relation
\begin{equation}
\overline{\mat{P}}(z;N,\bar{k})=\sigma_1\mat{P}(z;N,k)
\prod_{n=0}^{N-1}(z-x_{N,n})^{-\sigma_3}\sigma_1\,,
\hspace{0.2 in}\bar{k}=N-k\,.
\label{eq:Pdualexplicit}
\end{equation}
Here, $\mat{P}(z;N,k)$ is the solution of Interpolation
Problem~\ref{rhp:DOP} with weights $\{w_{N,j}\}$ on the nodes $X_N$,
and $\overline{\mat{P}}(z;N,\bar{k})$ is the solution of Interpolation
Problem~\ref{rhp:DOP} with the dual weights $\{\overline{w}_{N,j}\}$
defined by (\ref{eq:dualweightsdefine}) and with the exponent $k$ in
the normalization condition replaced by $\bar{k}$.
Note that \eqref{eq:Pdualexplicit} implies in particular that if $z$ and $w$ are {\em not} nodes ($z,w\not\in X_N$), then
\begin{equation}
\left[\overline{\mat{P}}(z;N,\bar{k})^{-1}\overline{\mat{P}}(w;N,\bar{k})\right]_{21} = \left[\mat{P}(z;N,k)^{-1}\mat{P}(w;N,k)\right]_{12}\prod_{n=0}^{N-1}(z-x_{N,n})(w-x_{N,n})\,.
\label{eq:Pdualexplicit21}
\end{equation}

Suppose first that $n\neq m$ are distinct indices.  Then
\begin{equation}
\begin{array}{rcl}
\displaystyle \overline{K}_{N,\bar{k}}(x_{N,m},x_{N,n})&=&\displaystyle
\frac{\sqrt{\overline{w}_{N,m}\overline{w}_{N,n}}}{x_{N,m}-x_{N,n}}
\left(\overline{\mat{P}}(x_{N,m};N,\bar{k})^{-1}
\overline{\mat{P}}(x_{N,n};N,\bar{k})\right)_{21}\\\\
&=&\displaystyle
\frac{\sqrt{\overline{w}_{N,m}\overline{w}_{N,n}}}{x_{N,m}-x_{N,n}}
\mathop{\lim_{w\rightarrow x_{N,m}}}_{z\rightarrow x_{N,n}}
\left(\overline{\mat{P}}(w;N,\bar{k})^{-1}
\overline{\mat{P}}(z;N,\bar{k})\right)_{21}\\\\
&=&\displaystyle
\frac{\sqrt{\overline{w}_{N,m}\overline{w}_{N,n}}}{x_{N,m}-x_{N,n}}\\\\
&&\displaystyle\,\,\,\cdot\,\,\,
\mathop{\lim_{w\rightarrow x_{N,m}}}_{z\rightarrow x_{N,n}}
\prod_{j=0}^{N-1}(w-x_{N,j})(z-x_{N,j})
\cdot
\left(\mat{P}(w;N,k)^{-1}\mat{P}(z;N,k)\right)_{12}\,.
\end{array}
\end{equation}
where in going from the second to the third line we have used
(\ref{eq:Pdualexplicit21}).  The limiting operation is necessary because
while
$\left(\overline{\mat{P}}(w;N,\bar{k})^{-1}\mat{P}(z;N,\bar{k})\right)_{21}$
is analytic in $w$ and $z$ near $w=x_{N,m}$ and $z=x_{N,n}$,
$\left(\mat{P}(w;N,k)^{-1}\mat{P}(z;N,k)\right)_{12}$ has
singularities at these points.  Next, using the definition
(\ref{eq:dualweightsdefine}) of the dual weights, we obtain
\begin{equation}
\begin{array}{rcl}
\displaystyle \overline{K}_{N,\bar{k}}(x_{N,m},x_{N,n})&=&\displaystyle
\frac{(-1)^{m+n}}{\sqrt{w_{N,m}w_{N,n}}}\cdot
\frac{\displaystyle\mathop{\lim_{w\rightarrow x_{N,m}}}_{z\rightarrow x_{N,n}}
\left[(w-x_{N,m})(z-x_{N,n})\left(\mat{P}(w;N,k)^{-1}\mat{P}(z;N,k)\right)_{12}\right]}{x_{N,m}-x_{N,n}}\\\\
&=&\displaystyle
\frac{(-1)^{m+n}}{\sqrt{w_{N,m}w_{N,n}}}\cdot\frac{\displaystyle
\left[\mathop{\rm Res}_{w=x_{N,m}}\mat{P}(w;N,k)^{-1}\mathop{\rm Res}_{z=x_{N,n}}\mat{P}(z;N,k)\right]_{12}}{x_{N,m}-x_{N,n}}\,,
\end{array}
\end{equation}
where we have used the fact that $\det \mat{P}(z;N,k)=1$ which implies that
$\mat{P}(z;N,k)^{-1}$ has simple poles at the nodes just like $\mat{P}(z;N,k)$ does.  Now again because $\det \mat{P}(z;N,k)=1$, we obtain from (\ref{eq:poles}) that
\begin{equation}
\mathop{\rm Res}_{w=x_{N,m}}\mat{P}(w;N,k)^{-1}=
\lim_{w\rightarrow x_{N,m}}\left(\begin{array}{cc}0 & -w_{N,m} \\ \\
0 & 0\end{array}\right)\mat{P}(w;N,k)^{-1}\,.
\end{equation}
Using this, together with (\ref{eq:poles}), we arrive at
\begin{equation}
\begin{array}{rcl}\displaystyle
\overline{K}_{N,\bar{k}}(x_{N,m},x_{N,n})&=&\displaystyle
\frac{(-1)^{m+n}}{\sqrt{w_{N,m}w_{N,n}}}\cdot \frac{\displaystyle
-w_{N,m}w_{N,n}\left(\mat{P}(x_{N,m};N,k)^{-1}\mat{P}(x_{N,n};N,k)\right)_{21}}{x_{N,m}-x_{N,n}}\\\\
&=&\displaystyle
(-1)^{m+n+1}\sqrt{w_{N,m}w_{N,n}}\frac{\displaystyle
\left(\mat{P}(x_{N,m};N,k)^{-1}\mat{P}(x_{N,n};N,k)\right)_{21}}{x_{N,m}-x_{N,n}}\\\\
&=& \displaystyle(-1)^{m+n+1}K_{N,k}(x_{N,m},x_{N,n})\,.
\end{array}
\end{equation}
Thus, we have proved the following,  a more specific version of (\ref{eq:correldual3squared}).
\begin{prop}
For distinct nodes $x=x_{N,m}$ and $y=x_{N,n}$ in $X_N$,
\begin{equation}
\overline{K}_{N,\bar{k}}(x,y)=(-1)^{m+n+1}K_{N,k}(x,y)
\end{equation}
where $\bar{k}=N-k$.
\label{prop:KKbaroffdiag}
\end{prop}

Now, we consider the reproducing kernel and its dual on the diagonal.  We begin with
\begin{equation}
\begin{array}{rcl}\displaystyle
\overline{K}_{N,\bar{k}}(x_{N,m},x_{N,m})&=&\displaystyle
\overline{w}_{N,m}\left[
\frac{d}{dz}\overline{\mat{P}}(z;N,\bar{k})^{-1}\Bigg|_{z=x_{N,m}}\overline{\mat{P}}(x_{N,m};N;\bar{k})\right]_{21}\\\\
&=&\displaystyle -\overline{w}_{N,m}\left[
\overline{\mat{P}}(x_{N,m};N,\bar{k})^{-1}\frac{d}{dz}\overline{\mat{P}}(z;N,\bar{k})\Bigg|_{z=x_{N,m}}\right]_{21}\,.
\end{array}
\end{equation}
But, using (\ref{eq:Pdualexplicit}), we see that
\begin{equation}
\begin{array}{rcl}\displaystyle
\overline{\mat{P}}(z;N,\bar{k})^{-1}\frac{d}{dz}\overline{\mat{P}}(z;N,\bar{k}) &= &\displaystyle
\sigma_1\prod_{j=0}^{N-1}(z-x_{N,j})^{\sigma_3}\left[\mat{P}(z;N,k)^{-1}\frac{d}{dz}\mat{P}(z;N,k)\right]\prod_{j=0}^{N-1}(z-x_{N,j})^{-\sigma_3}\sigma_1\\\\
&&\displaystyle\,\,\, +\,\,\,
\sigma_1\prod_{j=1}^{N-1}(z-x_{N,j})^{\sigma_3}\frac{d}{dz}\left[\prod_{j=0}^{N-1}(z-x_{N,j})^{-\sigma_3}\right]\sigma_1\,,
\end{array}
\end{equation}
and the second term is a diagonal matrix.  Consequently,
\begin{equation}
\overline{K}_{N,\bar{k}}(x_{N,m},x_{N,m}) = -\overline{w}_{N,m}
\lim_{z\rightarrow x_{N,m}}\left(\prod_{j=0}^{N-1}(z-x_{N,m})^2
\left[\mat{P}(z;N,k)^{-1}\frac{d}{dz}\mat{P}(z;N,k)\right]_{12}\right)\,.
\end{equation}
From Proposition~\ref{prop:solnrhp}, we then get
\begin{equation}
\begin{array}{l}\displaystyle
\left[\mat{P}(z;N,k)^{-1}\frac{d}{dz}\mat{P}(z;N,k)\right]_{12}\\\\
\displaystyle \hspace{0.5 in}
=\,\,\,\sum_{n=0}^{N-1}\sum_{j=0}^{N-1}\frac{P_{11}(x_{N,n};N,k)P_{21}(x_{N,j};N,k)-P_{11}(x_{N,j};N,k)P_{21}(x_{N,n};N,k)}{(z-x_{N,n})(z-x_{N,j})^2}w_{N,n}w_{N,j}\\\\
\displaystyle \hspace{0.5 in}
=\,\,\,\mathop{\sum\sum}_{n\neq j}
\frac{P_{11}(x_{N,n};N,k)P_{21}(x_{N,j};N,k)-P_{11}(x_{N,j};N,k)P_{21}(x_{N,n};N,k)}{(z-x_{N,n})(z-x_{N,j})^2}w_{N,n}w_{N,j}\,.
\end{array}
\end{equation}
Therefore,
\begin{equation}
\begin{array}{rcl}\displaystyle
\overline{K}_{N,\bar{k}}(x_{N,m},x_{N,m})&=&\displaystyle
-\overline{w}_{N,m}w_{N,m}\mathop{\prod_{j=0}}_{j\neq m}^{N-1}(x_{N,m}-x_{N,j})^2\\\\
&&\displaystyle\,\,\,\cdot
\mathop{\sum_{n=0}}_{n\neq m}^{N-1}\frac{P_{11}(x_{N,n};N,k)P_{21}(x_{N,m};N,k)-P_{11}(x_{N,m};N,k)P_{21}(x_{N,n};N,k)}{x_{N,m}-x_{N,n}}w_{N,n}\,,
\end{array}
\end{equation}
and using (\ref{eq:dualweightsdefine}), this becomes
\begin{equation}
\overline{K}_{N,\bar{k}}(x_{N,m},x_{N,m})=
-\mathop{\sum_{n=0}}_{n\neq m}^{N-1}\frac{P_{11}(x_{N,n};N,k)P_{21}(x_{N,m};N,k)-P_{11}(x_{N,m};N,k)P_{21}(x_{N,n};N,k)}{x_{N,m}-x_{N,n}}w_{N,n}\,.
\end{equation}
Now for $z\in\mathbb{C}\setminus X_N$, we have $\det\mat{P}(z;N,k)=1$, and taking the limit $z\rightarrow x_{N,m}$ with the use of the explicit formula for $\mat{P}(z;N,k)$ furnished by Proposition~\ref{prop:solnrhp} yields the
identity
\begin{equation}
\begin{array}{l}\displaystyle
w_{N,m}\left[P_{21}(x_{N,m};N,k)\frac{d}{dz}P_{11}(z;N,k)\Bigg|_{z=x_{N,m}}-P_{11}(x_{N,m};N,k)\frac{d}{dz}P_{21}(z;N,k)\Bigg|_{z=x_{N,m}}\right] \,\,\,+ \\\\
\displaystyle \hspace{0.4 in}\mathop{\sum_{n=0}}_{n\neq m}^{N-1}\frac{P_{11}(x_{N,m};N,k)P_{21}(x_{N,n};N,k)-P_{11}(x_{N,n};N,k)P_{21}(x_{N,m};N,k)}{x_{N,m}-x_{N,n}}w_{N,n} \,\,\,=\,\,\, 1\,.
\end{array}
\end{equation}
So, we have (again using $\det\mat{P}(z;N,k)=1$),
\begin{equation}
\overline{K}_{N,\bar{k}}(x_{N,m},x_{N,m})=1-w_{N,m}\left[
\frac{d}{dz}\mat{P}(z;N,k)^{-1}\Bigg|_{z=x_{N,m}}\mat{P}(x_{N,m};N,k)\right]_{21} = 1-K_{N,k}(x_{N,m},x_{N,m})
\end{equation}
which completes the direct proof of the following.
\begin{prop}
For any node $x\in X_N$,
\begin{equation}
\overline{K}_{N,\bar{k}}(x,x)=1-K_{N,k}(x,x)
\end{equation}
where $\bar{k}=N-k$.
\label{prop:KKbardiag}
\end{prop}

Combining Propositions~\ref{prop:KKbaroffdiag} and
\ref{prop:KKbardiag}, we therefore may write for any given set of
nodes $x_1,\dots,x_m$,
\begin{equation}
\left(\overline{K}_{N,\bar{k}}(x_i,x_j)\right)_{1\le i,j \le m} = \mat{D}
\left(\delta_{ij}-K_{N,k}(x_i,x_j)\right)_{1\le i,j\le m}\mat{D}\,,
\end{equation}
where $\mat{D}:={\rm diag}(1,-1,1,-1,\dots,(-1)^{m+1})$.
Taking determinants then yields another independent proof of Proposition~\ref{prop:dualkernel}.

\begin{remark}
  The dual ensemble is useful for several reasons.  Of course, the
  statistics of holes are often of independent interest.  But even if
  one is only interested in particle statistics, the dual ensemble is
  very helpful in the analysis of statistics near saturated regions of
  the node space $X_N$ where the upper constraint is active for the
  particle weights.  It follows from Proposition~\ref{prop:flip} that
  each saturated region for the particle weights with $k$ particles is
  a void for the (dual) hole weights with $\bar{k}=N-k$ holes.  In
  this way, each calculation valid for the particle ensemble near a
  void automatically translates via
  Proposition~\ref{prop:KKbaroffdiag} and
  Proposition~\ref{prop:KKbardiag} into a statement about particle
  statistics near saturated regions.
\end{remark}

\subsection{Exact formulae for $K_{N,k}(x,y)$.}
The following result will be used often below to obtain formulae for
$K_{N,k}(x,y)$ and $K_{N,k}(x,x)$ in various regions of $[a,b]$.
\begin{Lemma}
%\label{lemunvi}
Let $x$ be any node satisfying $x\in X_N\cap \Sigma_0^\nab$.  Then
\begin{equation}\label{eq;univ1}
\begin{array}{l}
\displaystyle
  \sqrt{w(x)} e^{(N\ell_c+\gamma)/2}e^{(k-\#\del)g_+(x)} \prod_{n\in\del} |x-x_{N,n}| \\\\
\displaystyle\hspace{0.4 in}  =\,\,\, \displaystyle
e^{(-\eta(x)+\gamma+2\kappa g_+(x))/2}e^{-iN\theta(x)/2} \frac{
\displaystyle e^{-\frac{1}{2}N \left[ \frac{\delta E_c}{\delta
\mu}(x)-\ell_c \right]}}{\sqrt{2\pi N\rho^0(x)}}
 T_{\nab}(x)^{1/2}
  \,.
\end{array}
\end{equation}
Here the variational derivative is evaluated on the equilibrium
measure $\mu_{\rm min}^c$, and $g_+(x)$ denotes the boundary value taken by $g(z)$ as $z\rightarrow x$ with $\Im(z)>0$.
\label{lem;commonfactor}
\end{Lemma}

\begin{proof}
  Let $x=x_{N,j} \in X_N\cap\Sigma^\nab_0$. Hence $j\in\nab$.
  Substituting for $w_N(\cdot)$ from \eqref{eq:weightrewrite} and
  \eqref{eq:weightform}, and using the fact that $x=x_{N,j}\in X_N$, we get
\begin{equation}
w(x) e^{N\ell_c+\gamma}e^{2(k-\#\del)g_+(x)} \prod_{n\in\del}
(x-x_{N,n})^2 = (-1)^{N-1-j} e^{-NV(x)-\eta(x))+ N\ell_c+\gamma +
2(k-\#\del)g_+(x)} \frac{\displaystyle \prod_{n\in\del}
  (x_{N,j}-x_{N,n})}{\displaystyle \prod_{\substack{n\in\nab \\ n\neq j} }
  (x_{N,j}-x_{N,n})} \,.
\end{equation}
But, using \eqref{eq:TAfuncdef}, we have
\begin{equation}
\begin{array}{l}
\displaystyle
\frac{\displaystyle\prod_{n\in\del}(x_{N,j}-x_{N,n})}
{\displaystyle\prod_{\substack{n\in\nab\\n\neq j}}(x_{N,j}-x_{N,n})}
=
\lim_{z\rightarrow x_{N,j}}(z-x_{N,j})
\frac{\displaystyle\prod_{n\in\del}(z-x_{N,n})}
{\displaystyle\prod_{n\in\nab}(z-x_{N,n})}\\\\
\displaystyle\hspace{0.2 in}=\,\,\, \lim_{z\rightarrow x_{N,j}}
\frac{z-x_{N,j}}{\displaystyle 2\cos\left(\frac{N\theta^0(z)}{2}\right)}
T_\nab(z)\exp\left(-N\left[\int_{\Sigma_0^\nab}\log|z-s|\rho^0(s)\,ds-
\int_{\Sigma_0^\del}\log|z-s|\rho^0(s)\,ds\right]\right)\,.
\end{array}
\end{equation}
The limit of the fraction can be taken using l'H\^opital's rule, and
the remaining factors are continuous for real $z$.  Thus, we arrive at
\begin{equation}
\frac{\displaystyle\prod_{n\in\del}(x_{N,j}-x_{N,n})}
{\displaystyle\prod_{\substack{n\in\nab\\n\neq
j}}(x_{N,j}-x_{N,n})} = \frac{\displaystyle T_\nab(x_{N,j}) \exp
\left(-N\left[\int_{\Sigma_0^\nab}\log|x_{N,j}-s|\rho^0(s)\,ds-
\int_{\Sigma_0^\del}\log|x_{N,j}-s|\rho^0(s)\,ds\right]\right)}
{\displaystyle 2\pi
N\rho^0(x_{N,j})\sin\left(\frac{N\theta^0(x_{N,j})}{2}\right)}\,.
\end{equation}
From the definition of $\rho^0$ and \eqref{eq:BS},
\begin{equation}
  \theta^0(x_{N,j}) = \pi\frac{2N-2j-1}{N}, \qquad
\sin\left(\frac{N\theta^0(x_{N,j})}{2}\right) = (-1)^{N-j-1}.
\end{equation}
Therefore, recalling the definition \eqref{eq:gdef} and
\eqref{eq:rhodef} of the complex phase function $g(z)$, the definition
\eqref{eq:variationalderivative} of the variational derivative of the
energy functional $E_c[\cdot]$ and the definition of $\theta(z)$
\eqref{eq:theta}, we obtain an identity that is the square of
\eqref{eq;univ1}.  By directly comparing the arguments of both sides of
\eqref{eq;univ1} one verifies that the square root has been taken
consistently.
\end{proof}

The following elementary lemma will be useful.
\begin{Lemma}
Let $f(x)$ and $M(x,y)$ be differentiable functions with $M(x,x)\equiv 0$.
Then
\begin{equation}
\frac{\partial}{\partial x}\left[f(x)f(y)M(x,y)\right]_{y=x}=f(x)^2\frac{\partial}{\partial x}M(x,y)\Bigg|_{y=x}\,.
\end{equation}
\label{lem:extractfactors}
\end{Lemma}
%Also the following simple facts of matrix algebra will be
%useful:
%\begin{equation}
%\label{eq;conjsigma31}
%  h^{\sigma_3} \begin{pmatrix} 1 & a\\0&1 \end{pmatrix} =
%  \begin{pmatrix} 1 & ah^2 \\0&1 \end{pmatrix}  h^{\sigma_3}\,,
%\qquad
%  h^{\sigma_3} \begin{pmatrix} 1 & 0\\a&1 \end{pmatrix} =
%  \begin{pmatrix} 1 & 0 \\ah^{-2}&1 \end{pmatrix} h^{\sigma_3}\,,
%\end{equation}
%\begin{equation}\label{eq;conjsigma32}
%  \biggl[ \begin{pmatrix} 1 & a\\0&1 \end{pmatrix} A \begin{pmatrix} 1 & b\\0&1
%  \end{pmatrix} \biggr]_{21} = A_{21}\,,
%\end{equation}
%\begin{equation}\label{eq;conjsigma325}
%  \biggl[ \begin{pmatrix} 1 & 0\\a&1 \end{pmatrix} A \begin{pmatrix} 1 & 0\\b&1
%  \end{pmatrix} \biggr]_{21} = A_{21} + aA_{11}+bA_{22}+abA_{12}\,,
%\end{equation}
%\begin{equation}\label{eq;conjsigma3255}
%  \biggl[ \begin{pmatrix} 1 & a\\0&1 \end{pmatrix} A \begin{pmatrix} 1 & b\\0&1
%  \end{pmatrix} \biggr]_{12} = A_{12} + aA_{22}+bA_{11}+abA_{21}\,,
%\end{equation}
%\begin{equation}\label{eq;conjsigma33}
%  \bigl( g^{\sigma_3} A h^{-\sigma_3} \bigr)_{21} =
%  (gh)^{-1}(A)_{21}\,,
%  \qquad
%  \bigl( g^{\sigma_3} A h^{-\sigma_3} \bigr)_{12} =
%  (gh)(A)_{12}\,.
%\end{equation}

We will now use these results to express $K_{N,k}(x,y)$ in terms of
the piecewise analytic global parametrix $\hat{\mat{X}}(z)$ and the
error matrix $\mat{E}(z)$, for $x$ and $y$ in different parts of the
interval $[a,b]$ of accumulation of the nodes.  The first result in
this direction is the following. \label{symbol:matrixB}\label{symbol:vandw}
\begin{prop}
Let $x$ and $y$ be distinct nodes in a band $I$, both lying in the same component of $\Sigma_0^\nab$ and lying outside all discs $D_\Gamma^{\nab,*}$.  Then
\begin{equation}
K_{N,k}(x,y)=\frac{1}{2\pi N\sqrt{\rho^0(x)\rho^0(y)}} \frac{\mat{v}^T
e^{iN\theta(x)\sigma_3/2}\mat{B}(x)^{-1}\mat{B}(y)e^{-iN\theta(y)\sigma_3/2}
\mat{w}}{x-y}
\label{eq:Kbandoffdiag}
\end{equation}
and
\begin{equation}
K_{N,k}(x,x)=\frac{1}{2\pi N\rho^0(x)}\left[2\pi Nc\frac{d\mu_{\rm min}^c}{dx}(x)- \mat{v}^Te^{iN\theta(x)\sigma_3/2}\mat{B}(x)^{-1}\mat{B}'(x)
e^{-iN\theta(x)\sigma_3/2}\mat{w}\right]\,,
\label{eq:Kbanddiag}
\end{equation}
where
\begin{equation}
\mat{v}:=\left(\begin{array}{c} -i \\ 1\end{array}\right)\,,\hspace{0.2 in}
\mat{w}:=\left(\begin{array}{c}
1 \\ i\end{array}\right)\,,
\end{equation}
(note that $\mat{v}^T\mat{w}=0$) and
\begin{equation}
\mat{B}(x):=
\mat{E}_+(x)\dot{\mat{X}}_+(x)e^{(\kappa g_+(x)+\gamma/2-\eta(x)/2)\sigma_3}\,,
\label{eq:Bdef}
\end{equation}
and the subscript ``$+$'' denotes the boundary value taken as $z\rightarrow x$ with $\Im(z)>0$.
\label{prop:Kbandexact}
\end{prop}

\begin{proof}
For distinct nodes $x$ and $y$, we begin with
\begin{equation}
K_{N,k}(x,y)=\sqrt{w(x)w(y)}\frac{\displaystyle\left[\mat{P}(x;N,k)^{-1}
\mat{P}(y;N,k)\right]_{21}}{x-y}
\end{equation}
and defining the quotient by l'H\^opital's rule,
\begin{equation}
K_{N,k}(x,x)=w(x)\frac{\partial}{\partial x}\left[\mat{P}(x;N,k)^{-1}
\mat{P}(y;N,k)\right]_{21}\Bigg|_{y=x}\,.
\end{equation}
Now, for any real $x\in\Sigma_0^\nab$, we
have from \eqref{eq:PtoQ} and \eqref{eq:QtoR1}
that
\begin{equation}
 \mat{P}(x;N,k)= \mat{R}_+(x) \left(\begin{array}{cc}
  1& \displaystyle ie^{- iN\theta^0(x)/2}e^{-NV_N(x)}\frac{\displaystyle
\prod_{n\in\del}(x-x_{N,n})}{\displaystyle\prod_{n\in\nab}
(x-x_{N,n})} \\\\ 0&1  \end{array}\right)
  \left[ \prod_{n\in\del}(x-x_{N,n}) \right]^{\sigma_3}\,,
\label{eq;Pnab}
\end{equation}
where $\mat{R}_+(x)$ denotes the boundary value taken from the upper
half-plane (from the left-hand side of the contour $\Sigma$; see
Figure~\ref{fig:Sigma}).  Thus,
\begin{equation}\label{eq:PinRSnab}
\begin{split}
  [\mat{P}(x;N,k)^{-1}\mat{P}(y;N,k)]_{21} &=
  [\mat{R}_+(x)^{-1}\mat{R}_+(y)]_{21}
  \prod_{n\in\del}(x-x_{N,n})(y-x_{N,n}) \\
  &= [\mat{S}_+(x)^{-1}\mat{S}_+(y)]_{21} e^{N\ell_c+\gamma}
  e^{(k-\#\del)(g_+(x)+g_+(y))}
  \prod_{n\in\del}(x-x_{N,n})(y-x_{N,n})
\end{split}
\end{equation}
where the second equality follows from \eqref{eq:SfromR}. When we further suppose that $x$ and $y$ lie within the same component of $\Sigma_0^\nab$ this formula may be rewritten as
\begin{equation}
\left[\mat{P}(x;N,k)^{-1}\mat{P}(y;N,k)\right]_{21}=
\left[\mat{S}_+(x)^{-1}\mat{S}_+(y)\right]_{21}
e^{N\ell_c+\gamma}e^{(k-\#\del)(g_+(x)+g_+(y))}
\prod_{n\in\del}|x-x_{N,n}||y-x_{N,n}|\,.
\label{eq:Pkernelsamecomponent}
\end{equation}
Letting $x$ and $y$ lie in a band $I\subset\Sigma_0^\nab$, we have
from \eqref{eq:X-nab-upper-def}, \eqref{eq:XfromEandXhat} and
\eqref{eq:XfromXdotlast} that for $z=x$ or $z=y$,
\begin{equation}
  \mat{S}_+(z)=
  \mat{E}(z)\dot{\mat{X}}_+(z)
  \left(\begin{array}{cc}
1 & 0 \\\\
ie^{\eta(z)-\gamma-2\kappa g_+(z)}e^{iN\theta(z)} &
1\end{array}\right) T_{\nab}(z)^{-\sigma_3/2} \,,
\label{eq:KbandexactS}
\end{equation}
and thus
\begin{equation}
\begin{array}{rcl}\displaystyle
\left[\mat{S}_+(x)^{-1}\mat{S}_+(y)\right]_{21}&=&\displaystyle
T_\nab(x)^{-1/2}e^{-(\kappa g_+(x)+\gamma/2-\eta(x)/2)}e^{iN\theta(x)/2}\\\\&&\displaystyle\,\,\,\cdot\,\,\,
T_\nab(y)^{-1/2}e^{-(\kappa g_+(y)+\gamma/2-\eta(y)/2)}e^{iN\theta(y)/2}\\\\
&&\displaystyle\,\,\,\cdot\,\,\,
\mat{v}^Te^{iN\theta(x)\sigma_3/2}\mat{B}(x)^{-1}\mat{B}(y)e^{-iN\theta(y)\sigma_3/2}\mat{w}\,,
\end{array}
\end{equation}
since $\mat{E}(z)$ is analytic in the band so that
$\mat{E}(z)=\mat{E}_+(z)$.  Now we substitute into
\eqref{eq:Pkernelsamecomponent}:
\begin{equation}
\left[\mat{P}(x;N,k)^{-1}\mat{P}(y;N,k)\right]_{21}=f(x)f(y)\mat{v}^T
e^{iN\theta(x)\sigma_3/2}\mat{B}(x)^{-1}\mat{B}(y)e^{-iN\theta(y)\sigma_3/2}\mat{w}
\label{eq:PxinvPy21}
\end{equation}
where
\begin{equation}
f(z):=T_\nab(z)^{-1/2}e^{-(\kappa g_+(z)-\eta(z)/2)}e^{N\ell_c/2}e^{(k-\#\del)g_+(z)}\prod_{n\in\del}|z-x_{N,n}|\,.
\end{equation}
Now \eqref{eq:PxinvPy21} holds for any $x$ and $y$ in the same band of $\Sigma_0^\nab$, and when we specialize to nodes $x,y\in X_N$, we obtain formulae for the reproducing kernel.  Therefore,
\begin{equation}
K_{N,k}(x,y)=\sqrt{w(x)w(y)}f(x)f(y)\frac{\mat{v}^Te^{iN\theta(x)\sigma_3/2}\mat{B}(x)^{-1}\mat{B}(y)e^{-iN\theta(y)\sigma_3/2}\mat{w}}{x-y}
\end{equation}
and using Lemma~\ref{lem:extractfactors},
\begin{equation}
K_{N,k}(x,x)=w(x)f(x)^2\frac{\partial}{\partial x}\left[\mat{v}^Te^{iN\theta(x)\sigma_3/2}\mat{B}(x)^{-1}\mat{B}(y)e^{-iN\theta(y)\sigma_3/2}\mat{w}\right]_{y=x}\,.
\label{eq:Kbanddiagalmost}
\end{equation}
Since $x\in X_N$ and $y\in X_N$, we may use Lemma~\ref{lem;commonfactor} along with the equilibrium condition \eqref{eq:equilibrium} that holds for $x$ and $y$ in a band $I$ to deduce
\begin{equation}
\sqrt{w(x)}f(x)=\frac{1}{\sqrt{2\pi N\rho^0(x)}}\hspace{0.2 in}
\text{and}\hspace{0.2 in}
\sqrt{w(y)}f(y)=\frac{1}{\sqrt{2\pi N\rho^0(y)}}\hspace{0.2 in}
\mbox{for $x$ and $y$ in $X_N\cap I\subset \Sigma_0^\nab$.}
\end{equation}
This proves (\ref{eq:Kbandoffdiag}).  To complete the proof of (\ref{eq:Kbanddiag}), we carry out the differentiation in (\ref{eq:Kbanddiagalmost}), noting that
by definition (see \eqref{eq:rhodef} and \eqref{eq:theta})
\begin{equation}
\theta'(x)=2\pi c \frac{d\mu_{\rm min}^c}{dx}(x)\,,\hspace{0.2 in}x\in\Sigma_0^\nab\,.
\end{equation}
%Hence using \eqref{eq;conjsigma325}
%\begin{equation}
%\begin{split}
%  [\mat{P}(z)^{-1}\mat{P}(w)]_{21} = &
%  \biggl\{ \bigl[(\mat{E}\dot{\mat{X}})(z)^{-%1}(\mat{E}\dot{\mat{X}})(w)\bigr]_{21}
%  -ie^{\eta(z)-\gamma-2\kappa g(z)}e^{iN\theta^\nab_I(z)}
%  \bigl[(\mat{E}\dot{\mat{X}})(z)^{-1}(\mat{E}\dot{\mat{X}})(w)\bigr]_{11} %\\
%  & +ie^{\eta(w)-\gamma-2\kappa g(w)}e^{iN\theta^\nab_I(w)}
%  \bigl[(\mat{E}\dot{\mat{X}})(z)^{-1}(\mat{E}\dot{\mat{X}})(w)\bigr]_{22} %\\
%  & + e^{\eta(z)-\gamma-2\kappa g(z)}e^{iN\theta^\nab_I(z)}
%  e^{\eta(w)-\gamma-2\kappa g(w)}e^{iN\theta^\nab_I(w)}
%  \bigl[(\mat{E}\dot{\mat{X}})(z)^{-1}(\mat{E}\dot{\mat{X}})(w)\bigr]_{12}
%  \biggr\} \\
%  & \times T_{\nab}(z)^{-1/2}T_{\nab}(w)^{-1/2}e^{N\ell_c+\gamma}
%  e^{(k-\#\del)(g(z)+g(w))}
%  \prod_{n\in\del}(z-x_{N,n})(w-x_{N,n}) .
%\end{split}
%\end{equation}
\end{proof}
\label{symbol:aandb}

\begin{prop}
Let $x$ and $y$ be distinct nodes in a void $\Gamma$  lying outside all
discs $D_\Gamma^{\nab,*}$.  Then
\begin{equation}
K_{N,k}(x,y)=\frac{T_\nab(x)^{1/2}T_\nab(y)^{1/2}e^{-\frac{1}{2}N\left[\frac{\delta E_c}{\delta\mu}(x)-\ell_c\right]}
e^{-\frac{1}{2}N\left[\frac{\delta E_c}{\delta\mu}(y)-\ell_c\right]}}
{2\pi N\sqrt{\rho^0(x)\rho^0(y)}}\cdot
\frac{\mat{a}^Te^{iN\theta_\Gamma\sigma_3/2}\mat{B}(x)^{-1}\mat{B}(y)e^{-iN\theta_\Gamma\sigma_3/2}\mat{b}}{x-y}
\label{eq:voidkernelexactoffdiag}
\end{equation}
and
\begin{equation}
K_{N,k}(x,x)=-\frac{T_\nab(x)e^{-N\left[\frac{\delta E_c}{\delta\mu}(x)-\ell_c\right]}}{2\pi N\rho^0(x)}\cdot
\mat{a}^Te^{iN\theta_\Gamma\sigma_3/2}\mat{B}(x)^{-1}\mat{B}'(x)
e^{-iN\theta_\Gamma\sigma_3/2}\mat{b}\,,
\label{eq:voidkernelexactdiag}
\end{equation}
where
\begin{equation}
\mat{a}:=\left(\begin{array}{c} 0 \\ 1\end{array}\right)\,,\hspace{0.2 in}
\mat{b}:=\left(\begin{array}{c} 1 \\ 0\end{array}\right)\,,
\end{equation}
(note that $\mat{a}^T\mat{b}=0$), $\mat{B}(x)$ is defined by (\ref{eq:Bdef}), and the variational derivative is evaluated on the
equilibrium measure $\mu_{\rm min}^c$.
\label{prop:voidkernelexact}
\end{prop}

\begin{proof}
  Note that the two points $x$ and $y$ lying in the same void interval
  necessarily belong to the same component of $\Sigma_0^\nab$.  The
  proof follows that of Proposition~\ref{prop:Kbandexact} with only a
  few modifications.  First, in place of \eqref{eq:KbandexactS} we
  have the simpler relation
\begin{equation}
\mat{S}_+(x)=\mat{E}_+(x)\dot{\mat{X}}_+(x)\,.
\end{equation}
Next, when we use Lemma~\ref{lem;commonfactor} we must retain the exponentials involving the variational derivative since in place of \eqref{eq:equilibrium} we have the variational inequality \eqref{eq:voidinequality} because $x$ and $y$ are in a void $\Gamma$.
Finally, we recall that the function $e^{iN\theta(x)}$ 
takes the constant value $e^{iN\theta_\Gamma}$
throughout $\Gamma$.
\end{proof}

Recall the definition of the mappings $\tau_\Gamma^{\nab, L}$ and
$\tau_\Gamma^{\nab, R}$ given in \eqref{eq:zetadefAL} and
\eqref{eq:zetadefAR} respectively. \label{symbol:AnabL}
\begin{prop}
\label{prop:Airykernelexact}
Let $x$ and $y$ be distinct nodes in a disc $D_\Gamma^{\nab,L}$.  Then
\begin{equation}
K_{N,k}(x,y)=\frac{1}{N^{2/3}\sqrt{\rho^0(x)\rho^0(y)}}\cdot
\frac{\mat{q}_\Gamma^{\nab,L}(x)^T\mat{A}_\Gamma^{\nab,L}(x)^{-1}\mat{A}_\Gamma^{\nab,L}(y)\mat{r}_\Gamma^{\nab,L}(y)}{x-y}
\label{eq:Airykernelexact}
\end{equation}
and
\begin{equation}
K_{N,k}(x,x)=\frac{1}{N^{2/3}\rho^0(x)}\left[-\mat{q}_\Gamma^{\nab,L}(x)^T
\mat{A}_\Gamma^{\nab,L}(x)^{-1}\frac{d\mat{A}_\Gamma^{\nab,L}}{dx}(x)
\mat{r}_\Gamma^{\nab,L}(x)-\mat{q}_\Gamma^{\nab,L}(x)^T
\frac{d\mat{r}_\Gamma^{\nab,L}}{dx}(x)\right]\,,
\label{eq:Airykernelexactdiag}
\end{equation}
where
\begin{equation}
\mat{q}_\Gamma^{\nab,L}(x):=\left(\begin{array}{c}\displaystyle -Ai\left(-\left(\frac{3}{4}\right)^{2/3}\tau_\Gamma^{\nab,L}(x)\right)\\\\
\displaystyle
N^{-1/3}Ai'\left(-\left(\frac{3}{4}\right)^{2/3}\tau_\Gamma^{\nab,L}(x)\right)
\end{array}\right)\,,\hspace{0.2 in}
\mat{r}_\Gamma^{\nab,L}(x):=\left(\begin{array}{c}\displaystyle
N^{-1/3}Ai'\left(-\left(\frac{3}{4}\right)^{2/3}\tau_\Gamma^{\nab,L}(x)\right)\\\\
\displaystyle
Ai\left(-\left(\frac{3}{4}\right)^{2/3}\tau_\Gamma^{\nab,L}(x)\right)
\end{array}\right)\,,
\end{equation}
(note that $\mat{q}_\Gamma^{\nab,L}(x)^T\mat{r}_\Gamma^{\nab,L}(x)
\equiv 0$), and
\begin{equation}
\mat{A}_\Gamma^{\nab,L}(x):=\mat{E}(x)
\mat{H}_\Gamma^{\nab,L}(x)N^{\sigma_3/6}\left(\frac{3}{4}\right)^{-\sigma_3/6}\,.
\label{eq:ALdef}
\end{equation}
Similarly, if $x$ and $y$ are distinct nodes in a disc
$D_\Gamma^{\nab,R}$, then
\begin{equation}
K_{N,k}(x,y)=-\frac{1}{N^{2/3}\sqrt{\rho^0(x)\rho^0(y)}}\cdot
\frac{\mat{q}_\Gamma^{\nab,R}(x)\mat{A}_\Gamma^{\nab,R}(x)^{-1}
\mat{A}_\Gamma^{\nab,R}(y)\mat{r}_\Gamma^{\nab,R}(y)}{x-y}
\label{eq:Airyexactright}
\end{equation}
and
\begin{equation}
K_{N,k}(x,x)=-\frac{1}{N^{2/3}\rho^0(x)}\left[
-\mat{q}_\Gamma^{\nab,R}(x)^T\mat{A}_\Gamma^{\nab,R}(x)^{-1}
\frac{d\mat{A}_\Gamma^{\nab,R}}{dx}(x)\mat{r}_\Gamma^{\nab,R}(x)-
\mat{q}_\Gamma^{\nab,R}(x)^T\frac{d\mat{r}_\Gamma^{\nab,R}}{dx}(x)
\right]\,,
\label{eq:Airyexactrightdiag}
\end{equation}
where
\begin{equation}
\mat{q}_\Gamma^{\nab,R}(x):=\left(\begin{array}{c}
\displaystyle -Ai\left(-\left(\frac{3}{4}\right)^{2/3}\tau_\Gamma^{\nab,R}(x)
\right)\\\\
\displaystyle N^{-1/3}Ai'\left(-\left(\frac{3}{4}\right)^{2/3}\tau_\Gamma^{\nab,R}(x)
\right)\end{array}\right)\,,\hspace{0.2 in}
\mat{r}_\Gamma^{\nab,R}(x):=\left(\begin{array}{c}
\displaystyle N^{-1/3}Ai'\left(-\left(\frac{3}{4}\right)^{2/3}\tau_\Gamma^{\nab,R}(x)
\right)\\\\
\displaystyle Ai\left(-\left(\frac{3}{4}\right)^{2/3}\tau_\Gamma^{\nab,R}(x)
\right)\end{array}\right)\,,
\end{equation}
(note again that
$\mat{q}_\Gamma^{\nab,R}(x)^T\mat{r}_\Gamma^{\nab,R}(x)\equiv 0$), and
\begin{equation}
\mat{A}_\Gamma^{\nab,R}(x):=\mat{E}(x)\mat{H}_\Gamma^{\nab,R}(x)N^{\sigma_3/6}
\left(\frac{3}{4}\right)^{-\sigma_3/6}\,.
\label{eq:ARdef}
\end{equation}
\end{prop}
\label{symbol:AnabR}

\begin{proof}
Using the fact that $x$ and $y$ necessarily lie in the same component of $\Sigma_0^\nab$, we still have the relation \eqref{eq:Pkernelsamecomponent},
where the subscript ``$+$'' indicates a boundary value taken from the upper half-plane.  Now the matrix $\mat{S}(z)$ is analytic for all $z\in D_\Gamma^{\nab,L}\cap \mathbb{C}_+$, so to obtain a formula for $\mat{S}(z)$ we may choose arbitrarily whether to consider $z$ in quadrant I or quadrant II of $D_\Gamma^{\nab,L}$ (the answer is necessarily the same).  For concreteness, we choose to evaluate $\mat{S}_+(x)$ by taking a limit from  $D_{\Gamma, II}^{\nab,L}$ (above the void $\Gamma$). In this region, $\mat{S}(z)\equiv\mat{X}(z)\equiv \mat{E}(z)\hat{\mat{X}}_\Gamma^{\nab,L}(z)$, so from \eqref{eq;XhatnabL} we then obtain that for $x$ and $y$ in $D_{\Gamma, II}^{\nab,L}$,
\begin{equation}
\begin{array}{rcl}\displaystyle
\left[\mat{S}(x)^{-1}\mat{S}(y)\right]_{21}&=&\displaystyle
T_\nab(x)^{-1/2}e^{(\eta(x)-\gamma-2\kappa g(x))/2}e^{iN\theta_\Gamma/2}\cdot
 T_\nab(y)^{-1/2}e^{(\eta(y)-\gamma-2\kappa g(y))/2}e^{iN\theta_\Gamma/2}
 \\\\&&\displaystyle\,\,\,
 \cdot\,\,\,
 N^{-1/3} \left(\frac{3}{4}\right)^{1/3} \left[\mat{G}(x)^{-1}
\mat{A}_\Gamma^{\nab,L}(x)^{-1}\mat{A}_\Gamma^{\nab,L}(y)
 \mat{G}(y)\right]_{21}\, ,
 \end{array}
 \end{equation}
where
\begin{equation}
 \mat{G}(z) :=
 \left(\frac{3}{4}\right)^{\sigma_3/6}
 N^{-\sigma_3/6}
 \hat{\mat{Z}}^{\nab,L}(\tau_\Gamma^{\nab,L}(z))
 N^{\sigma_3/6}\left(\frac{3}{4}\right)^{-\sigma_3/6}\,.
\end{equation}
 Using the explicit formula for
$\hat{\mat{Z}}^{\nab,L}(\zeta)$ furnished by \eqref{eq:Airythird} of
Proposition~\ref{propAirynabL}, we then obtain
\begin{equation}
\begin{array}{rcl}
\displaystyle \left[\mat{G}(x)^{-1}
\mat{A}_\Gamma^{\nab,L}(x)^{-1}\mat{A}_\Gamma^{\nab,L}(y)
\mat{G}(y)\right]_{21}&=&\displaystyle
2\pi \left(\frac{3}{4}\right)^{-1/3} N^{2/3}e^{(-\tau_\Gamma^{\nab,L}(x))^{3/2}/2} e^{(-\tau_\Gamma^{\nab,L}(y))^{3/2}/2} \\\\
&&\displaystyle\,\,\,\cdot\,\,\,\mat{q}_\Gamma^{\nab,L}(x)^T
\mat{A}_\Gamma^{\nab,L}(x)^{-1}\mat{A}_\Gamma^{\nab,L}(y)
\mat{r}_\Gamma^{\nab,L}(y)\,.
\end{array}
\end{equation}
Substituting into \eqref{eq:Pkernelsamecomponent} gives
\begin{equation}
\left[\mat{P}(x;N,k)^{-1}\mat{P}(y;N,k)\right]_{21}=N^{1/3}
b(x)b(y)\mat{q}_\Gamma^{\nab,L}(x)^T\mat{A}_\Gamma^{\nab,L}(x)^{-1}
\mat{A}_\Gamma^{\nab,L}(y)\mat{r}_\Gamma^{\nab,L}(y)
\end{equation}
where
\begin{equation}
b(z):=\sqrt{2\pi} e^{(-\tau_\Gamma^{\nab,L}(z))^{3/2}/2}T_\nab(z)^{-1/2}e^{-(\kappa g_+(z)-\eta(z))/2}e^{iN\theta_\Gamma/2}
e^{N\ell_c/2}e^{(k-\#\del)g_+(z)}\prod_{n\in\del}|z-x_{N,n}|\,.
\end{equation}
Here, by $(-\tau_\Gamma^{\nab,L}(z))^{3/2}$ we understand the boundary
value taken on $\mathbb{R}$ from $\Im(z)>0$, or equivalently
$\Im(\tau_\Gamma^{\nab,L}(z))>0$.  The subscript ``$+$'' on $g_+(z)$
denotes the same limit.  Using Lemma~\ref{lem;commonfactor}, we see
that for any node $z$ in $D_\Gamma^{\nab,L}$,
\begin{equation}
\sqrt{w(z)}b(z)=\frac{\displaystyle
e^{(-\tau_\Gamma^{\nab,L}(z))^{-3/2}/2}e^{-iN(\theta(z)-\theta_\Gamma)/2}e^{-\frac{1}{2}N\left[\frac{\delta E_c}{\delta\mu}(z)-\ell_c\right]}}{\sqrt{N\rho^0(z)}}\,,
\end{equation}
where the variational derivative is evaluated for $\mu=\mu_{\rm
min}^c$.  Now, if the node $z$ lies in the void $\Gamma$, then
$\theta(z)=\theta_\Gamma$ modulo $2\pi /N$, but from
(\ref{eq:gapplus}) and (\ref{eq:zetadefAL}), we see that
$(-\tau_\Gamma^{\nab,L}(z))^{-3/2} = N[\delta
E_c/\delta\mu(z)-\ell_c]$.  On the other hand, if the node $z$ lies in
the adjacent band, then from \eqref{eq:equilibrium} we have $\delta
E_c/\delta\mu(z)-\ell_c=0$, but again (\ref{eq:zetadefAL}) gives the
identity $(-\tau_\Gamma^{\nab,L}(z))^{-3/2} =
iN(\theta(z)-\theta_\Gamma)$ modulo $2\pi i$.  Thus, for all nodes
$z$ in $D_\Gamma^{\nab,L}$, we have
$\sqrt{w(z)}b(z)=1/\sqrt{N\rho^0(z)}$.  This proves
\eqref{eq:Airykernelexact}.  Using Lemma~\ref{lem:extractfactors} we
also obtain \eqref{eq:Airykernelexactdiag}.

The proofs of (\ref{eq:Airyexactright}) and
(\ref{eq:Airyexactrightdiag}) are analogous.  It is perhaps noteworthy
that the origin of the leading minus sign in these formulae is the
factor $i\sigma_3$ relating $\hat{\mat{Z}}^{\nab,L}(\zeta)$ and
$\hat{\mat{Z}}^{\nab,R}(\zeta)$ (see \eqref{eq:VhatRA}).
\end{proof}

\begin{remark}
  In each case we may verify after the fact that for a node $x\in X_N$,
\begin{equation}
K_{N,k}(x,x)=\mathop{\lim_{z,w\rightarrow x}}_{z,w\in\mathbb{C},\,\, z\neq w} K_{N,k}(z,w)\,,
\label{eq:amazinglimit}
\end{equation}
that is, in each region $K_{N,k}(x,y)$ may be viewed as an analytic
function of two complex variables sampled at the discrete nodes
$X_N\times X_N$.  This is not obvious from the definition.  Indeed,
the definition \eqref{eq:reproducingdef} of $K_{N,k}(x,y)$ can {\em a
priori} only be evaluated when $x$ and $y$ are both nodes due to the
factor $\sqrt{w(x)w(y)}$.  If the weights were given in the form
$w_{N,n}=w(x_{N,n})$ for some analytic function $w(x)$ there would be
a direct interpretation of the limit process \eqref{eq:amazinglimit}.
However, the weights under consideration (given by
\eqref{eq:weightform}) do not have the exact form of an analytic
function simply sampled at the nodes, due to the presence of a factor
involving an essentially discrete product over nodes.  Indeed, the
derivation of the exact formulae above for $K_{N,k}(x,y)$ both on and
off the diagonal made explicit use of the fact that $x$ and $y$ are
discrete nodes via Lemma~\ref{lem;commonfactor}.
\end{remark}

The following result can also be extracted from the proofs of
Propositions~\ref{prop:voidkernelexact} and
\ref{prop:Airykernelexact}.
\begin{prop}
\label{prop:Kinoutexact}
  Let $x$ and $y$ be nodes in the same component of $\Sigma_0^\nab$.
If $y$ lies in a disc
  $D_\Gamma^{\nab,L}$ and $x$ lies outside the disc but in the adjacent
  void $\Gamma$.  Then
\begin{equation}
K_{N,k}(x,y)=-\frac{T_\nab(x)^{1/2}e^{-\frac{1}{2}N\left[\frac{\delta E_c}{\delta\mu}(x)-\ell_c\right]}}{N^{5/6}\sqrt{2\pi\rho^0(x)\rho^0(y)}}\cdot
\frac{\mat{a}^Te^{iN\theta_\Gamma\sigma_3/2}\mat{B}(x)^{-1}
\mat{A}_\Gamma^{\nab,L}(y)\mat{r}_\Gamma^{\nab,L}(y)}{x-y}\,.
\end{equation}
Similarly, if $y$ lies in a disc $D_\Gamma^{\nab,R}$ and $x$ lies outside the disc but in the adjacent void $\Gamma$, then
\begin{equation}
K_{N,k}(x,y)=-i\frac{T_\nab(x)^{1/2}e^{-\frac{1}{2}N\left[\frac{\delta E_c}{\delta\mu}(x)-\ell_c\right]}}{N^{5/6}\sqrt{2\pi\rho^0(x)\rho^0(y)}}\cdot
\frac{\mat{a}^Te^{iN\theta_\Gamma\sigma_3/2}\mat{B}(x)^{-1}
\mat{A}_\Gamma^{\nab,R}(y)\mat{r}_\Gamma^{\nab,R}(y)}{x-y}\,.
\end{equation}
Here the notation on the right-hand side is the same as in Proposition~\ref{prop:voidkernelexact} and Proposition~\ref{prop:Airykernelexact}.
\end{prop}

\subsection{Asymptotic formulae for $K_{N,k}(x,y)$ and universality.}

\begin{Lemma}
  Fix a closed interval $F\subset[a,b]$ that contains none of the band
  endpoints $\alpha_0,\dots,\alpha_G$ and $\beta_0,\dots,\beta_G$.
  Without loss of generality, fix the contour parameter $\epsilon>0$
  sufficiently small that $F$ lies outside all discs $D_\Gamma^{*,*}$
  and that Proposition~\ref{prop:E} controls $\mat{E}(z)-\mathbb{I}$.
  Then there is a constant $C_F>0$ such that for all $N$ sufficiently
  large,
\begin{equation}
\sup_{x\in F} \|\mat{B}'(x)\| \le C_F\hspace{0.2 in}
\mbox{and}\hspace{0.2 in}
\sup_{x,y\in F} \frac{\|\mat{B}(x)^{-1}\mat{B}(y)-\mathbb{I}\|}{|x-y|}
\le C_F\,,
\end{equation}
where $\|\cdot\|$ denotes a matrix norm and $\mat{B}(x)$ is defined by
\eqref{eq:Bdef} for arbitrary $x\in [a,b]$ (note that $\mat{B}(x)$
depends on $\epsilon$ via $\mat{E}(x)$).
\label{lem;Mestimateband}
\end{Lemma}
% For any $z, w$ in an interval of either bands, voids or
%  saturated regions bounded away from the endpoints of the band/gap
%  edges by a fixed distance,
%\begin{equation}
%  \frac{
%\bigl(\mat{E}_+\dot{\mat{X}}_+\bigr)(z)^{-1}
%(\mat{E}_+\dot{\mat{X}}_+)(w) - I}{z-w}
%\end{equation}
%is bounded uniformly with respect to $N$.
%\label{lem;Mestimateband}
%\end{Lemma}

\begin{proof}
  The matrix $\mat{W}(z):=\dot{\mat{X}}(z)e^{(\kappa
    g(z)+\gamma/2-\eta(z)/2)\sigma_3}$ can be analytically continued
  through the interval $F$ from the upper half-plane by a jump
  relation of the form (see Riemann-Hilbert Problem~\ref{rhp:theta})
  $\mat{W}_+(z)=\mat{W}_-(z)\mat{v}$ where $\mat{v}$ is a constant
  matrix (with respect to $z$) whose entries are uniformly bounded as
  $N\rightarrow\infty$ (for $z$ in a void or saturated region
  $\Gamma_i$ we have $\mat{v}=e^{iN\theta_{\Gamma_i}\sigma_3}$ and for
  $z$ in a band $I$ we have $\mat{v}=-i\sigma_1$).  Since $\mat{W}(z)$
  is uniformly bounded for $z\in\mathbb{C}\setminus\Sigma_{\rm model}$
  bounded away from the band endpoints (from
  Proposition~\ref{prop:dotXbound} and \eqref{eq:etaNcontrol} as well
  as the assumption that $\kappa$ remains bounded as
  $N\rightarrow\infty$), it follows that the analytic continuation of
  $\mat{W}_+(z)$ from $F$ is uniformly bounded in a fixed complex
  neighborhood $G$ of $F$ as $N\rightarrow\infty$.  Cauchy's Theorem
  applied on a closed contour in $G$ encircling $F$ then shows that
  $\mat{W}_+(z)$ and all its derivatives remain uniformly bounded in
  $F$ as $N\rightarrow\infty$.

  The same is true of the matrix $\mat{E}(z)$.  Indeed, if $F$ is a
  subinterval of a band $I$, then $\mat{E}(z)$ is already analytic in
  a complex neighborhood $G$ of $F$, and is uniformly bounded in $G$
  as $N\rightarrow\infty$ according to Proposition~\ref{prop:E}.  The
  uniform boundedness of all derivatives of $\mat{E}_+(z)=\mat{E}(z)$
  for $z\in F$ then follows from Cauchy's Theorem.  On the other hand,
  if $F$ is a subinterval of a void or saturated region, then the
  analytic continuation of $\mat{E}_+(z)$ to the neighborhood $G$ is
  accomplished by the formula $\mat{E}_+(z)=\mat{F}(z)$ where
  $\mat{F}(z)$ is the solution of Riemann-Hilbert Problem~\ref{rhp:F}.
  Since $\mat{F}(z)$ is uniformly bounded in $G$, again Cauchy's
  Theorem implies that all derivatives of $\mat{E}_+(z)$ are uniformly
  bounded for $z\in F$.
  
  Combining these results using $\mat{B}(x)=\mat{E}_+(x)\mat{W}_+(x)$
  establishes that $\mat{B}'(x)$ remains uniformly bounded in $F$ as
  $N\rightarrow\infty$.  The boundedness of the difference quotient
  follows from this result and the uniform boundedness of $\mat{B}(x)$
  itself, since $\det(\mat{B}(x))=1$.
\end{proof}

\begin{Lemma}\label{lem;EdgeEHZ}
Fix a value of the contour parameter $\epsilon>0$ sufficiently small
that Proposition~\ref{prop:E} controls $\mat{E}(z)-\mathbb{I}$ on
appropriate closed sets.  Then for each disc $D_\Gamma^{\nab,L}$ there
is a constant $C_\Gamma^{\nab,L}>0$ and for each disc
$D_\Gamma^{\nab,R}$ there is a constant $C_\Gamma^{\nab,R}>0$ such
that for all $N$ sufficiently large,
\begin{equation}
\sup_{x\in D_\Gamma^{\nab,L}\cap\mathbb{R}}\left\|\frac{d\mat{A}^{\nab,L}_\Gamma}{dx}(x)\right\|\le C_\Gamma^{\nab,L}\hspace{0.2 in}
\text{and}\hspace{0.2 in}
\sup_{x,y\in D_\Gamma^{\nab,L}\cap\mathbb{R}}\frac{\|\mat{A}^{\nab,L}_\Gamma(x)^{-1}\mat{A}_\Gamma^{\nab,L}(y)-\mathbb{I}\|}{|x-y|}\le C_\Gamma^{\nab,L}
\end{equation}
and
\begin{equation}
\sup_{x\in D_\Gamma^{\nab,R}\cap\mathbb{R}}\left\|\frac{d\mat{A}^{\nab,R}_\Gamma}{dx}(x)\right\|\le C_\Gamma^{\nab,R}\hspace{0.2 in}
\text{and}\hspace{0.2 in}
\sup_{x,y\in D_\Gamma^{\nab,R}\cap\mathbb{R}}\frac{\|\mat{A}^{\nab,R}_\Gamma(x)^{-1}\mat{A}_\Gamma^{\nab,R}(y)-\mathbb{I}\|}{|x-y|}\le C_\Gamma^{\nab,R}\,,
\end{equation}
where $\|\cdot\|$ denotes a matrix norm, and $\mat{A}_\Gamma^{\nab,L}(x)$ is defined by \eqref{eq:ALdef} and $\mat{A}_\Gamma^{\nab,R}(x)$ is defined by \eqref{eq:ARdef}.

Also, for the same constants and for sufficiently large $N$,
\begin{equation}
\sup_{x\in D_\Gamma^{\nab,L}\cap\mathbb{R}}\|\mat{q}_\Gamma^{\nab,L}(x)\|\le
C_\Gamma^{\nab,L}\hspace{0.2 in}\text{and}\hspace{0.2 in}
\sup_{x\in D_\Gamma^{\nab,L}\cap\mathbb{R}}\|\mat{r}_\Gamma^{\nab,L}(x)\|\le
C_\Gamma^{\nab,L}
\end{equation}
and
\begin{equation}
\sup_{x\in D_\Gamma^{\nab,R}\cap\mathbb{R}}\|\mat{q}_\Gamma^{\nab,R}(x)\|\le
C_\Gamma^{\nab,R}\hspace{0.2 in}\text{and}\hspace{0.2 in}
\sup_{x\in D_\Gamma^{\nab,R}\cap\mathbb{R}}\|\mat{r}_\Gamma^{\nab,R}(x)\|\le
C_\Gamma^{\nab,R}\,.
\end{equation}

Finally, there is a constant $K>0$ such that for sufficiently large $N$,
\begin{equation}
\mathop{\sup_{x\in D_\Gamma^{\nab,L}\cap\mathbb{R}}}_{x<\alpha}
\|\mat{q}_\Gamma^{\nab,L}(x)\|\le\frac{C_\Gamma^{\nab,L}e^{-NK(\alpha-x)^{3/2}}}{N^{1/6}}\hspace{0.2 in}\text{and}\hspace{0.2 in}
\mathop{\sup_{x\in D_\Gamma^{\nab,L}\cap\mathbb{R}}}_{x<\alpha}
\|\mat{r}_\Gamma^{\nab,L}(x)\|\le\frac{C_\Gamma^{\nab,L}e^{-NK(\alpha-x)^{3/2}}}{N^{1/6}}\,,
\end{equation}
where $\alpha$ is the band edge point at the center of the disc $D_\Gamma^{\nab,L}$, and
\begin{equation}
\mathop{\sup_{x\in D_\Gamma^{\nab,R}\cap\mathbb{R}}}_{x>\beta}
\|\mat{q}_\Gamma^{\nab,R}(x)\|\le\frac{C_\Gamma^{\nab,R}e^{-NK(x-\beta)^{3/2}}}{N^{1/6}}\hspace{0.2 in}\text{and}\hspace{0.2 in}
\mathop{\sup_{x\in D_\Gamma^{\nab,R}\cap\mathbb{R}}}_{x>\beta}
\|\mat{r}_\Gamma^{\nab,R}(x)\|\le\frac{C_\Gamma^{\nab,R}e^{-NK(x-\beta)^{3/2}}}{N^{1/6}}\,,
\end{equation}
where $\beta$ is the band edge point at the center of the disc $D_\Gamma^{\nab,R}$.
\end{Lemma}

\begin{proof}
The statements concerning the matrices $\mat{A}_\Gamma^{\nab,L}(x)$
and $\mat{A}_\Gamma^{\nab,R}(x)$ are elementary consequences of two
facts.  First, from Proposition~\ref{prop:E}, we have that
$\mat{E}(z)$ is analytic and remains uniformly bounded as
$N\rightarrow\infty$ in each disc $D_\Gamma^{\nab,L}$ or
$D_\Gamma^{\nab,R}$.  Next (see \S~\ref{sec:Airy}) the product
$\mat{H}_\Gamma^{\nab,L}(z)N^{\sigma_3/6}$ is analytic in each disc
$D_\Gamma^{\nab,L}$ and remains uniformly bounded there as
$N\rightarrow\infty$, while the product
$\mat{H}_\Gamma^{\nab,R}(z)N^{\sigma_3/6}$ is analytic in each disc
$D_\Gamma^{\nab,R}$ and remains uniformly bounded there as
$N\rightarrow\infty$.  It follows from Cauchy's Theorem applied on the
boundary of each disc that all derivatives of
$\mat{A}_\Gamma^{\nab,L}(z)$ are uniformly bounded independent of $N$
in $D_\Gamma^{\nab,L}$, and the same holds for
$\mat{A}_\Gamma^{\nab,R}(z)$ in $D_\Gamma^{\nab,R}$.  The boundedness
of the difference quotients then follows since $\det
(\mat{A}_\Gamma^{\nab,L}(x))=1$ in $D_\Gamma^{\nab,L}$ and
$\det(\mat{A}_\Gamma^{\nab,R}(x))=1$ in $D_\Gamma^{\nab,R}$.

The statements concerning the vectors $\mat{q}_\Gamma^{\nab,L}(x)$, $\mat{r}_\Gamma^{\nab,L}(x)$, $\mat{q}_\Gamma^{\nab,R}(x)$, and $\mat{r}_\Gamma^{\nab,R}(x)$ are obtained from the asymptotic formulae \eqref{eq:AiryAsymp} and from the elementary estimates
holding for all $x>0$:
\begin{equation}
|Ai(x)|\le \frac{Ce^{-2x^{3/2}/3}}{(1+x)^{1/4}}\hspace{0.2 in}
\text{and}\hspace{0.2 in}
|Ai'(x)|\le C (1+x)^{1/4}e^{-2x^{3/2}/3}\,,
\end{equation}
where $C>0$ is some appropriate constant.  Then one uses the fact that
in each case the argument of the Airy functions is $N^{2/3}$ times an analytic function of $x$ that has a nonvanishing derivative and is independent of $N$.
\end{proof}

\subsubsection{Universal statistics for particles in a band.  Proof of Theorem~\ref{theorem;bulk} and Theorem~\ref{theorem;bulk2}.}
\label{sec:univband}

Consider a fixed closed interval $F$ in the interior of a band $I$.
We can easily establish the following asymptotic formulae uniformly valid
in $F$.
\begin{Lemma}
Let $F$ be a fixed closed interval in the interior of a band $I$.  Then there is a constant $C_F>0$ such that for all sufficiently large $N$,
\begin{equation}
\max_{x\in X_N\cap F}\left|K_{N,k}(x,x)-\frac{c}{\rho^0(x)}\frac{d\mu_{\rm min}^c}{dx}(x)\right|\le\frac{C_F}{N}\,,
\label{eq:Kbanddiagasymp}
\end{equation}
and
\begin{equation}
\max_{x,y\in X_N\cap F} \left|
K_{N,k}(x,y)-\frac{1}{N\pi\sqrt{\rho^0(x)\rho^0(y)}}\frac{\displaystyle
\sin\left(\frac{N}{2}(\theta(x)-\theta(y))\right)}{x-y}\right| \le
\frac{C_F}{N}\,  \label{eq:Kbandoffdiagasymp}
\end{equation}
where $\theta(z)$ is defined in \eqref{eq:theta}. Also, for some
other constant $C_F'>0$ and $N$ sufficiently large,
\begin{equation}
\max_{x,y\in X_N\cap F}\left|(x-y)K_{N,k}(x,y)\right|\le\frac{C_F'}{N}\,.
\label{eq:Kbandoffdiagasympcrude}
\end{equation}
\label{lem;Kbanduniform}
\end{Lemma}

\begin{proof}
First, note that without any loss of generality, we may suppose that $F$ lies
in $\Sigma_0^\nab$.  Indeed, if $I$ is a transition band this can be arranged
by judicious choice of the transition points $Y_\infty$. But even if $I$
lies between two saturated regions, a pair of artificial transition
points may be introduced in $I$, one on each side of $F$ in order to
``switch'' $F$ back into $\Sigma_0^\nab$.
Then, applying
Lemma~\ref{lem;Mestimateband} to the exact formula \eqref{eq:Kbanddiag} established in Proposition~\ref{prop:Kbandexact} then yields \eqref{eq:Kbanddiagasymp}.  Similarly, applying Lemma~\ref{lem;Mestimateband} to \eqref{eq:Kbandoffdiag}
and using the identity
\begin{equation}
\mat{v}^Te^{iN(\theta(x)-\theta(y))\sigma_3/2}  \mat{w}=2\sin\left(\frac{N}{2}(\theta(x)-\theta(y))\right)
\end{equation}
proves \eqref{eq:Kbandoffdiagasymp}.  Then \eqref{eq:Kbandoffdiagasympcrude} follows from \eqref{eq:Kbandoffdiagasymp} .
\end{proof}

The estimate \eqref{eq:Kbandoffdiagasympcrude} shows that the reproducing kernel is concentrated near the
diagonal, and a nonzero limit for $K_{N,k}(x,y)$ as
$N\rightarrow\infty$ may only be expected for nodes $x$ and $y$ in $F$ with $x-y$ of
size bounded by $1/N$.
To find the limit, we will now localize by considering a finite number
of nodes near a certain fixed $x\in F$ (being fixed as
$N\rightarrow\infty$, $x$ is not necessarily a node).
Since
\begin{equation}
\begin{array}{rcl}
\displaystyle
R^{(N,k)}_1(x)\cdot N\rho^0(x)&=&\displaystyle
\mathbb{E}(\text{number of particles per node near $x$})\cdot\frac{\text{number of nodes}}{\text{unit length}}\\\\
&=&\displaystyle\mathbb{E}\left(\frac{\text{number of particles}}{\text{unit length}}\right)\,,
\end{array}
\end{equation}
from \eqref{eq:determinantalformulacorrelations} and
\eqref{eq:Kbanddiagasymp} we see that the asymptotic mean spacing
between particles near $x\in F$ is $\delta(x)/N$ where $\delta(x)$
is defined by \eqref{eq:scaledparticlespacing}.
For $\xi_N$ and $\eta_N$ in some bounded set $D$, we thus consider
nodes $z$ and $w$ defined by
\begin{equation}\label{eq;zwscaled}
  z:= x+ \xi_N\frac{\delta(x)}{N}\,,\hspace{0.2 in}
 w:= x+ \eta_N\frac{\delta(x)}{N}\,.
\end{equation}
Note that the admissible values of $\xi_N$ and $\eta_N$ are finite in
number and are asymptotically equally spaced with spacing
$(\rho^0(x)\delta(x))^{-1}$ (because $z$ and $w$ are both nodes in $X_N$).
Since $F\subset I$ and thus neither constraint is active, we
have $0<(\rho^0(x)\delta(x))^{-1}<1$.

Now, from Taylor's Theorem and \eqref{eq;zwscaled} we have
\begin{equation}
\begin{array}{rcl}\displaystyle
\theta(z)-\theta(w)&=&\displaystyle
\theta\left(x+\xi_N\frac{\delta(x)}{N}\right)-\theta\left(x+\eta_N\frac{\delta(x)}{N}\right)\\\\
&=&\displaystyle
\theta\left(x+\eta_N\frac{\delta(x)}{N} + (\xi_N-\eta_N)\frac{\delta(x)}{N}\right)-\theta\left(x+\eta_N\frac{\delta(x)}{N}\right)\\\\
&=&\displaystyle \theta'\left(x+\eta_N\frac{\delta(x)}{N}\right)(\xi_N-\eta_N)\frac{\delta(x)}{N} + \frac{\theta''(\sigma)}{2}(\xi_N-\eta_N)^2\frac{\delta(x)^2}{N^2}\\\\
&=&\displaystyle \theta'(x)(\xi_N-\eta_N)\frac{\delta(x)}{N} + \theta''(\tau)
(\xi_N-\eta_N)\eta_N\frac{\delta(x)^2}{N^2}+\frac{\theta''(\sigma)}{2}
(\xi_N-\eta_N)^2\frac{\delta(x)^2}{N^2}\\\\
&=&\displaystyle \frac{2\pi}{N}(\xi_N-\eta_N) +\theta''(\tau)
(\xi_N-\eta_N)\eta_N\frac{\delta(x)^2}{N^2}+\frac{\theta''(\sigma)}{2}
(\xi_N-\eta_N)^2\frac{\delta(x)^2}{N^2}
\end{array}
\end{equation}
for some $\sigma$ and $\tau$ near $x\in F$.  In the last line we have used $\theta'(x)=2\pi cd\mu_{\rm min}^c/dx(x)$.
Therefore,
\begin{equation}
\frac{\displaystyle\sin\left(\frac{N}{2}(\theta(z)-\theta(w))\right)}{z-w}=\frac{N}{\delta(x)}\left[\frac{\sin(\pi(\xi_N-\eta_N))}{\xi_N-\eta_N} -
\cos(q)\left(\theta''(\tau)\eta_N + \frac{\theta''(\sigma)}{2}(\xi_N-\eta_N)\right)\frac{\delta(x)^2}{2N}\right]\,,
\end{equation}
for some $q\in\mathbb{R}$.  

Since the node density is analytic and positive, we have
\begin{equation}
\frac{1}{\sqrt{\rho^0(z)\rho^0(w)}}=\frac{1}{\rho^0(x)} + O\left(\frac{1}{N}\right)
\end{equation}
because we are assuming $\xi_N$ and $\eta_N$ to remain bounded as $N\rightarrow\infty$.  Combining these results with Lemma~\ref{lem;Kbanduniform} proves the following.
\begin{Lemma}
Fix $x$ in the interior of any band $I$, and consider $\xi_N$ and $\eta_N$ to lie in a fixed bounded discrete set $D$ such that $z$ and $w$ defined by \eqref{eq;zwscaled} lie in the set of nodes $X_N$.  Then there is
a constant $C_D(x)>0$ such that for all sufficiently large $N$,
\begin{equation}
\max_{\xi_N,\eta_N\in D}\left| K_{N,k}(z,w)-\frac{c}{\rho^0(x)}\frac{d\mu_{\rm min}^c}{dx}(x)
S(\xi_N,\eta_N)\right|\le \frac{C_D(x)}{N}\,.
\end{equation}
\label{lem;Kasympsinekernel}
\end{Lemma}

Applying Lemma~\ref{lem;Kbanduniform} and Lemma~\ref{lem;Kasympsinekernel} to the determinantal formula \eqref{eq:determinantalformulacorrelations}, we immediately obtain corresponding asymptotics for all multipoint correlation functions, which completes the proof of Theorem~\ref{theorem;bulk}.

We now give the proof of Theorem~\ref{theorem;bulk2}.
From Lemma~\ref{lem;Kasympsinekernel} and the asymptotic equal spacing of $\xi_N$ and $\eta_N$, it follows that if $x_{N,i}$ and $x_{N,j}$ are two nodes in $X_N$ such that $x_{N,i}\rightarrow x$ and $x_{N,j}\rightarrow x$ while $i-j$ remains fixed as $N\rightarrow\infty$, and $x$ is in the interior of a band $I$, then
\begin{equation}
K_{N,k}(x_{N,i},x_{N,j})={\cal S}_{ij}(x)+O\left(\frac{1}{N}\right)\,.
\label{eq:discretesinekernelequalspacing}
\end{equation}

Recall the formula \eqref{eq;Am2} for $A_m^{(N,k)}(B)$ and its
interpretation \eqref{eq;Am1} as a probability. The operator $K_{N,k}\bigl|_{B_N}$ acts on $\ell^2(B_N)$
with the kernel given by
\begin{equation}
  K_{N,k}\left( x+ (x_i-x), x+(x_j-x) \right)\,,
\end{equation}
where the $x_i$ are the nodes in $B_N$.  The first result
is that as $N\rightarrow\infty$,
\begin{equation}\label{eq;limittoS1}
  \det \left(1-tK_{N,k}\bigl|_{B_N} \right)
  = \det\left( 1 - t\mathcal{S}(x)\bigl|_{\mathbb{B}} \right) +
O\left(\frac{1}{N}\right)
\end{equation}
holds uniformly for $t$ in compact sets in $\cx$.
This follows from the analytic dependence of
determinants of matrices of fixed finite dimension on the matrix
elements, using Lemma~\ref{lem;Kasympsinekernel}, and then using \eqref{eq:discretesinekernelequalspacing}. The statement \eqref{eq;limittoS2} then follows
from the analyticity of the left-hand side of \eqref{eq;limittoS1} in $t$.

\subsubsection{Correlation functions for particles in voids.  Proof of Theorem~\ref{theorem:onepointvoid} and Theorem~\ref{theorem:normalvoid}.}
\label{sec:univvoid}
Let $F=[u,v]$ be a fixed closed interval in a void $\Gamma$ such that $u\not\in\{\beta_0,\dots,\beta_G\}$ and $v\not\in\{\alpha_0,\dots,\alpha_G\}$.  We admit the possibility that $u=a$ or $v=b$.  Applying Lemma~\ref{lem;Mestimateband} to the exact formulae
\eqref{eq:voidkernelexactoffdiag} and \eqref{eq:voidkernelexactdiag} of Proposition~\ref{prop:voidkernelexact}, and taking into account the variational inequality \eqref{eq:voidinequality}, we arrive at the following.

\begin{Lemma}
\label{lem;Kvoiduniform}
Let $F$ be a fixed closed interval in a void $\Gamma$ that is bounded away from all bands.  Then there is a constant $C_F>0$ such that for all $N$ sufficiently large,
\begin{equation}
\max_{x,y\in X_N\cap F} \left|K_{N,k}(x,y)e^{\frac{1}{2}N\left[\frac{\delta E_c}{\delta\mu}(x)-\ell_c\right]}e^{\frac{1}{2}N\left[\frac{\delta E_c}{\delta\mu}(x)-\ell_c\right]}\right|\le \frac{C_F}{N}\,,
\label{eq:voiduniformestimate}
\end{equation}
where the variational derivatives are evaluated on the equilibrium measure.
\end{Lemma}
Applying this result to the formula \eqref{eq:determinantalformulacorrelations} for the correlation functions, we complete the proof of Theorem~\ref{theorem:onepointvoid}.

Now we prove Theorem~\ref{theorem:normalvoid}.  From \eqref{eq:voiduniformestimate} we see that the reproducing
kernel $K_{N,k}(x,y)$ is uniformly exponentially small for nodes
in $F$.  In particular, and unlike in the bands, the kernel is not
concentrated near the diagonal.  In fact, concentration of
$K_{N,k}(x,y)$ for $x$ and $y$ in a set $F$ bounded away from the
bands requires the existence of a local minimum of $\delta
E_c/\delta\mu-\ell_c$ at some point $x\in F$.  Indeed, suppose
first that the local minimum occurs at some point in the {\em
interior} of $F$ and is genuine, so that the expansion \eqref{eq:variationalTaylor} holds
with $W>0$ and $H > 0$ as $z\rightarrow x$ (by assumption the variational derivative is an analytic function of $z$ in the void $\Gamma$).  The proper scaling is evidently then $z-x = O(N^{-1/2})$.  Thus, we consider nodes $z$ and $w$ near $x$ of the form
\begin{equation}
\label{eq:zwvoidscaled}
z=x+\frac{\xi_N}{H\sqrt{N}}\,,\hspace{0.2 in}w=x+\frac{\eta_N}{H\sqrt{N}}\,,
\end{equation}
for $\xi_N$ and $\eta_N$ in some bounded discrete set $D$ such
that $z$ and $w$ are in the set of nodes $X_N$. We then have from
Proposition \ref{prop:voidkernelexact} that
\begin{equation}
K_{N,k}(z,w)=\frac{e^{-NW}}{N}\cdot\left[
q_N(x)+
O\left(\frac{1}{\sqrt{N}}\right)\right]\cdot e^{-(\xi_N^2 + \eta_N^2)/2} \,,
\end{equation}
where the error is uniform for $\xi_N$ and $\eta_N$ in $D$, and
\begin{equation}
q_N(x):=-\frac{T_\nab(x)}{2\pi\rho^0(x)}\cdot \mat{a}^Te^{iN\theta_\Gamma\sigma_3/2}\mat{B}(x)^{-1}\mat{B}'(x)e^{-iN\theta_\Gamma\sigma_3/2}\mat{b}
\label{eq:qN}
\end{equation}
is uniformly bounded as $N\rightarrow\infty$.
Now the asymptotic spacing between points in the discrete set $D$ is $H/(\sqrt{N}\rho^0(x))$ which goes to zero as $N\rightarrow\infty$.  Thus, for any fixed interval $[A,B]\subset\mathbb{R}$,
\begin{equation}
\begin{array}{rcl}\displaystyle
E_{\rm int}([A,B];x,H,N)&
%:=&\displaystyle
%\mathbb{E}(\text{number of particles at nodes $z$ of the form %(\ref{eq:zwvoidscaled}) with $A\le \xi_N \le B$})\\\\ &
= &\displaystyle\mathop{\sum_{\xi_N\in D}}_{A\le \xi_N\le B}R_1^{(N,k)}\left(x+\frac{\xi_N}{H\sqrt{N}}\right)\\\\
&=&\displaystyle
\mathop{\sum_{\xi_N\in D}}_{A\le \xi_N\le B}K_{N,k}\left(x+\frac{\xi_N}{H\sqrt{N}},x+\frac{\xi_N}{H\sqrt{N}}\right)\\\\
&=&\displaystyle
\frac{e^{-NW}}{N}\mathop{\sum_{\xi_N\in D}}_{A\le \xi_N\le B}\left[q_N(x)+O\left(\frac{1}{\sqrt{N}}\right)\right]\cdot
e^{-\xi_N^2}\\\\
&=&\displaystyle
\frac{e^{-NW}\rho^0(x)}{H\sqrt{N}}\mathop{\sum_{\xi_N\in D}}_{A\le \xi_N\le B}\left[q_N(x)+O\left(\frac{1}{\sqrt{N}}\right)\right]\cdot
e^{-\xi_N^2}\frac{H}{\sqrt{N}\rho^0(x)}\\\\
&=&\displaystyle
\frac{e^{-NW}\rho^0(x)}{H\sqrt{N}}\left[q_N(x)\int_A^Be^{-\xi^2}\,d\xi + O\left(\frac{1}{\sqrt{N}}\right)\right]\,.
\end{array}
\end{equation}

The statement \eqref{eq:normalvoid} will be established if we can bound
$q_N(x)$ away from zero as $N\rightarrow\infty$.   Now, for $x$ in a void $\Gamma$, and for $N$ sufficiently large that $\mat{E}(z)-\mathbb{I}$ is sufficiently small, $q_N(x)$ will be bounded away from zero if the Wronskian
\begin{equation}
W[\dot{X}_{11+}e^{\kappa g_+},\dot{X}_{21+}e^{\kappa g_+}](x):=
\dot{X}_{12+}(x)e^{\kappa g_+(x)}\frac{d}{dx}\dot{X}_{11+}(x)e^{\kappa g_+(x)}-
\dot{X}_{11+}(x)e^{\kappa g_+(x)}\frac{d}{dx}\dot{X}_{21+}(x)e^{\kappa g_+(x)}\,,
\end{equation}
where the subscript ``$+$'' indicates a boundary value taken from the
upper half-plane, is bounded away from zero.  The Wronskian is not
identically zero in any subinterval of $\Gamma$, for the following
reasons.  If $W[\dot{X}_{11+}e^{\kappa g_+},\dot{X}_{21+}e^{\kappa
g_+}](x)$ were identically zero as a function of $x$, then there would
necessarily be an interval in which $\dot{X}_{11+}(x)$ and
$\dot{X}_{21+}(x)$ are proportional by a constant multiplier.
Analytically extending this proportionality toward $z=\infty$, we see
from the normalization condition
($\dot{\mat{X}}(z)\rightarrow\mathbb{I}$ as $z\rightarrow\infty$) that
we would have to have $\dot{X}_{21}(z)\equiv 0$.  The jump condition
for $\dot{\mat{X}}(z)$ in any band then forces $\dot{X}_{22}(z)\equiv
0$ in addition, contradicting the fact that
$\det(\dot{\mat{X}}(z))=1$.

Since $\dot{\mat{X}}(z)$ takes analytic boundary values in the void $\Gamma$, it follows that $W[\dot{X}_{11+}e^{\kappa g_+},\dot{X}_{21+}e^{\kappa g_+}](x)$ has only isolated zeros in $\Gamma$
for each value of $N$.  From the exact solution formulae given in Appendix~\ref{sec:thetasolve}, the number of zeros in $\Gamma$ is finite and remains uniformly bounded as $N\rightarrow\infty$, and the zeros move quasiperiodically as $N$ varies.   For a given $x\in\Gamma$, either we have
$W[\dot{X}_{11+}e^{\kappa g_+},\dot{X}_{21+}e^{\kappa g_+}](x)= 0$ for all $N\in\mathbb{Z}$, or for each sufficiently small $\epsilon>0$ we may extract a subsequence of $N$ values for which
$|W[\dot{X}_{11+}e^{\kappa g_+},\dot{X}_{21+}e^{\kappa g_+}](x)|\ge \epsilon$.  The first situation may only occur for a finite number of $x\in\Gamma$.  This completes the proof of Theorem~\ref{theorem:normalvoid}.

\subsubsection{Correlation functions for particles in saturated regions.  Proof of Theorem~\ref{theorem:onepointsaturatedregion} and Theorem~\ref{theorem:normalsaturatedregion}.}
Let $F$ be a fixed closed interval in a saturated region $\Gamma$ that is bounded away from the bands, but which may have either $a$ or $b$ as an endpoint if the upper constraint is active there.  We may exploit the dual ensemble (for the holes) to analyze the particle statistics in $F$.  According to Proposition~\ref{prop:flip}, the equilibrium measures for the particle ensemble with $k$ particles and for the dual hole ensemble with $\bar{k}=N-k$ holes are explicitly related, and $F$ lies in a void for the hole ensemble.  Consequently, the results of our analysis for $x$ and $y$ in a void hold true for the dual kernel $\overline{K}_{N,\bar{k}}(x,y)$ and the corresponding hole correlation functions.  To recover results for the kernel $K_{N,k}(x,y)$ and the corresponding particle correlation functions in $F$, we simply apply
Propositions~\ref{prop:KKbaroffdiag} and \ref{prop:KKbardiag}.
This proves Theorem~\ref{theorem:onepointsaturatedregion}.

Combining the above duality arguments with the proof of Theorem~\ref{theorem:normalvoid} proves Theorem~\ref{theorem:normalsaturatedregion}.  
 
\subsubsection{Universal statistics for particles near band edges.
Proof of Theorem~\ref{theorem:Airyunivcorrvoid}, Theorem~\ref{theorem:Airyunivcorrsat}, Theorem~\ref{thm:Airydet}, and
Theorem~\ref{thm:Airydetforhole}.}
\label{sec:univedgevoid}
Near the edge of a band, the equilibrium measure vanishes, and
hence in the band the scaling $x-y=O(1/N)$ is not
correct as the one-point function vanishes in the limit
$N\to\infty$. Also in the void or the saturated region near the
band edge, the positive constant \eqref{eq:voidKFconst} or
\eqref{eq:satKFconst} is no longer bounded away from zero. Therefore we
need to introduce a different scaling near a band edge to find the correct
scaling limit.

Under the generic simplifying assumptions listed in \S~\ref{sec:C3}, the density $d\mu^c_{\rm min}/dx$ of the equilibrium measure
vanishes like a square root at each band edge adjacent to a void, and
at a band edge adjacent to a saturated region, the ``dual equilibrium
measure'' $\rho^0(x)/c- d\mu^c_{\rm min}/dx(x)$ vanishes like a square
root. In this case, it turns out that the proper scaling is to
consider nodes $x$ satisfying
\begin{equation}
  x-\alpha =O(N^{-2/3}) \hspace{0.2 in}\text{or}\hspace{0.2 in}
  x-\beta=O(N^{-2/3})
\end{equation}
depending on whether we consider a left band edge $z_0=\alpha$ or
a right band edge $z_0=\beta$.
Below, we will show that the limiting correlation
function under the above scaling is given by the Airy kernel as in
the so-called edge scaling limit of the Gaussian unitary ensemble of random matrix theory, and also in the context of ensembles of
more general Hermitian matrices of invariant measure (see e.g.
\cite{TracyW94} and \cite{BleherI99}).

We begin with the following lemma, which is the analogue of
Lemma~\ref{lem;Kbanduniform}.
\begin{Lemma}
\label{lem;Kedgeuniform}
For each disc $D_\Gamma^{\nab,L}$ there is
a constant $C_\Gamma^{\nab,L}>0$ such that for all sufficiently large $N$,
\begin{equation}
K_{N,k}(x,x)=-\frac{t'(x)}{N^{1/3}\rho^0(x)}
\left[Ai'\left(N^{2/3}t(x)\right)^2-Ai\left(N^{2/3}t(x)\right)Ai''\left(N^{2/3}t(x)\right)\right] + \varepsilon^{(1)}_N(x)\,,
\end{equation}
and
\begin{equation}
K_{N,k}(x,y)=-\frac{1}{N\sqrt{\rho^0(x)\rho^0(y)}}\cdot\frac{\displaystyle Ai\left(N^{2/3}t(x)\right)Ai'\left(N^{2/3}t(y)\right)-Ai'\left(N^{2/3}t(x)\right)Ai\left(N^{2/3}t(y)\right)}{x-y}+\varepsilon^{(2)}_N(x,y)\,,
\end{equation}
where $t(x):=-(3/4)^{2/3}N^{-2/3}\tau_\Gamma^{\nab,L}(x)$
is a real-analytic function in $D_\Gamma^{\nab,L}$ that is independent of $N$ and strictly decreasing along the real axis, and
\begin{equation}
\max_{x\in X_N\cap D_\Gamma^{\nab,L}}\left|\varepsilon^{(1)}_N(x)\right|\le\frac{C_\Gamma^{\nab,L}}{N^{2/3}}\hspace{0.2 in}\text{and}\hspace{0.2 in}
\max_{x,y\in X_N\cap D_\Gamma^{\nab,L}}\left|\varepsilon^{(2)}_N(x,y)\right|\le\frac{C_\Gamma^{\nab,L}}{N^{2/3}}\,,
\end{equation}
and also for some constant $K>0$ we have the one-sided estimates
\begin{equation}
\begin{array}{rcl}\displaystyle
\mathop{\max_{x\in X_N\cap D_\Gamma^{\nab,L}}}_{x<\alpha}
\left|\varepsilon_N^{(1)}(x)\right|&\le &\displaystyle
\frac{C_\Gamma^{\nab,L}e^{-2NK(\alpha-x)^{3/2}}}{N}\,,\\\\\displaystyle\mathop{\max_{x,y\in X_N\cap D_\Gamma^{\nab,L}}}_{x,y<\alpha}
\left|\varepsilon_N^{(2)}(x,y)\right|&\le &\displaystyle
\frac{C_\Gamma^{\nab,L}e^{-NK(\alpha-x)^{3/2}}e^{-NK(\alpha-y)^{3/2}}}{N}\,,
\end{array}
\end{equation}
where $\alpha$ is the band edge point at the center of the disc $D_\Gamma^{\nab,L}$.
Similarly, for each disc $D_\Gamma^{\nab,R}$ there exists a constant $C_\Gamma^{\nab,R}>0$ such that for all sufficiently large $N$,
\begin{equation}
K_{N,k}(x,x)=\frac{t'(x)}{N^{1/3}\rho^0(x)}\left[Ai'\left(N^{2/3}t(x)\right)^2-Ai\left(N^{2/3}t(x)\right)Ai''\left(N^{2/3}t(x)\right)\right] +\varepsilon^{(1)}_N(x)\,,
\end{equation}
and
\begin{equation}
K_{N,k}(x,y)=\frac{1}{N\sqrt{\rho^0(x)\rho^0(y)}}\cdot\frac{\displaystyle Ai\left(N^{2/3}t(x)\right)Ai'\left(N^{2/3}t(y)\right)-Ai'\left(N^{2/3}t(x)\right)Ai\left(N^{2/3}t(y)\right)}{x-y} + \varepsilon^{(2)}_N(x,y)\,,
\end{equation}
where now $t(x):=-(3/4)^{2/3}N^{-2/3}\tau_\Gamma^{\nab,R}(x)$
is a real-analytic function in $D_\Gamma^{\nab,R}$ that is independent of $N$ and strictly increasing along the real axis, and
\begin{equation}
\max_{x\in X_N\cap D_\Gamma^{\nab,R}}\left|\varepsilon^{(1)}_N(x)\right|\le\frac{C_\Gamma^{\nab,R}}{N^{2/3}}\,,\hspace{0.2 in}\text{and}\hspace{0.2 in}
\max_{x,y\in X_N\cap D_\Gamma^{\nab,R}}\left|\varepsilon^{(2)}_N(x,y)\right|\le\frac{C_\Gamma^{\nab,R}}{N^{2/3}}\,,
\end{equation}
and also for some constant $K>0$ we have the one-sided estimates
\begin{equation}
\begin{array}{rcl}\displaystyle
\mathop{\max_{x\in X_N\cap D_\Gamma^{\nab,R}}}_{x>\beta}
\left|\varepsilon_N^{(1)}(x)\right|&\le &\displaystyle
\frac{C_\Gamma^{\nab,R}e^{-2NK(x-\beta)^{3/2}}}{N}\,,\\\\\displaystyle
\mathop{\max_{x,y\in X_N\cap D_\Gamma^{\nab,R}}}_{x,y>\beta}
\left|\varepsilon_N^{(2)}(x,y)\right|&\le &\displaystyle\frac{C_\Gamma^{\nab,R}e^{-NK(x-\beta)^{3/2}}e^{-NK(y-\beta)^{3/2}}}{N}\,,
\end{array}
\end{equation}
where $\beta$ is the band edge point at the center of the disc $D_\Gamma^{\nab,R}$.
\end{Lemma}

\begin{proof}
This follows from Proposition~\ref{prop:Airykernelexact} and Lemma~\ref{lem;EdgeEHZ}.
\end{proof}

Now we localize near the diagonal by considering $\xi_N$ and $\eta_N$ to lie in a fixed bounded set such that
\begin{equation}
x=z_0+\frac{\xi_N}{t'(z_0)N^{2/3}}\hspace{0.2 in}\text{and}\hspace{0.2 in}
y=z_0+\frac{\eta_N}{t'(z_0)N^{2/3}}
\label{eq:Airylocalize}
\end{equation}
are nodes.  Here $z_0=\alpha$ or $z_0=\beta$ is the band edge, which is independent of $N$.  Because $t(z_0)=0$, the Airy kernel $A(\xi_N,\eta_N)$
defined by \eqref{eq:Airykernel}
will appear in the asymptotics with $\xi_N$ and $\eta_N$ considered bounded.  Although for each $N$ the possible values of $\xi_N$ and $\eta_N$ are discrete, their spacing tends to zero like $N^{-1/3}$, and
in this sense the Airy kernel, unlike the discrete sine kernel,  may be thought of as a continuous function of two independent variables.

Now in a disc $D_\Gamma^{\nab,L}$ centered at a left band edge $z_0=\alpha$, a direct calculation using \eqref{eq:zetadefAL} shows that $t'(\alpha)=-\left(\pi c B^L_\alpha\right)^{2/3}$, where $B^L_\alpha$ is defined
in \eqref{eq:Bleft}.
Similarly, in a disc $D_\Gamma^{\nab,R}$ centered at a right band edge $z_0=\beta$, one may use \eqref{eq:zetadefAR} to see that $t'(\beta)=\left(\pi c B^R_\beta\right)^{2/3}$ where $B^R_\beta$ is
defined in \eqref{eq:Bright}.
With the help of Lemma~\ref{lem;Kedgeuniform} we may prove the following result.
\begin{Lemma}
For each fixed $M>0$ and each left band edge $\alpha$ there is a constant
$C_\alpha(M)>0$ such that for sufficiently large $N$,
\begin{equation}
\mathop{\max_{x,y\in X_N}}_{\alpha-MN^{-1/2}<x,y<\alpha+MN^{-2/3}}\left|K_{N,k}(x,y)-\frac{\displaystyle\left(\pi c B^L_\alpha\right)^{2/3}}{N^{1/3}\rho^0(\alpha)}A(\xi_N,\eta_N)\right|\le
\frac{C_\alpha(M)}{N^{2/3}}\,,
\end{equation}
where $B^L_\alpha$ is defined via a limit from the adjacent band from (\ref{eq:Bleft}) and
$\xi_N$ and $\eta_N$ are defined in terms of $x$ and $y$ using \eqref{eq:Airylocalize} and $t'(\alpha)=-\left(\pi c B^L_\alpha\right)^{2/3}$.
Similarly, for each fixed $M>0$ and each right band edge $\beta$ there is a constant $C_\beta(M)>0$ such that for sufficiently large $N$,
\begin{equation}
\mathop{\max_{x,y\in X_N}}_{\beta-MN^{-2/3}<x,y<\beta+MN^{-1/2}}\left|K_{N,k}(x,y)-\frac{\displaystyle\left(\pi c B^R_\beta\right)^{2/3}}{N^{1/3}\rho^0(\beta)}A(\xi_N,\eta_N)\right|\le
\frac{C_\beta(M)}{N^{2/3}}\,,
\end{equation}
where $B^R_\beta$ is defined via a limit from the adjacent band from (\ref{eq:Bright}) and $\xi_N$ and $\eta_N$ are defined in terms of $x$ and $y$ using \eqref{eq:Airylocalize} and $t'(\beta)=\left(\pi c B^R_\beta\right)^{2/3}$.
\label{lem;Kedgeinner}
\end{Lemma}

\begin{proof}
We show how the the computation works for $x$ and $y$ near a left endpoint $\alpha$.  The calculation near $\beta$ is similar.  From Lemma~\ref{lem;Kedgeuniform} we have
\begin{equation}
K_{N,k}(x,y)=-\frac{1}{N^{1/3}\sqrt{\rho^0(x)\rho^0(y)}}\cdot
\frac{t(x)-t(y)}{x-y}\cdot A(N^{2/3}t(x),N^{2/3}t(y))+\varepsilon_N^{(2)}(x,y)\,.
\end{equation}
With $\alpha-MN^{-1/2}<x,y<\alpha+MN^{-2/3}$, we have
\begin{equation}
\frac{1}{\sqrt{\rho^0(x)\rho^0(y)}}\cdot\frac{t(x)-t(y)}{x-y}=\frac{t'(\alpha)}{\rho^0(\alpha)}+O\left(\frac{1}{N^{1/2}}\right)\,.
\end{equation}
Since all partial derivatives of the Airy kernel $A(\xi_N,\eta_N)$ tend rapidly to zero as $\xi_N$ and $\eta_N$ tend to $+\infty$ (while under our assumptions on $x$ and $y$, $\xi_N$ and $\eta_N$ are bounded below by a fixed constant, they may grow in the positive direction like $N^{1/6}$), we then obtain with
$\alpha-MN^{-1/2}<x,y<\alpha+MN^{-2/3}$,
\begin{equation}
A(N^{2/3}t(x),N^{2/3}t(y))=A(\xi_N,\eta_N)+O\left(\frac{1}{N^{1/3}}\right)\,.
\end{equation}
Combining these estimates with the uniform estimate $\varepsilon_N(x,y)=O(N^{-2/3})$ furnished by Lemma~\ref{lem;Kedgeuniform} gives the desired result.
\end{proof}

Applying this result to the determinantal formula \eqref{eq:determinantalformulacorrelations} for the correlation functions
completes the proof of Theorem~\ref{theorem:Airyunivcorrvoid}.

Theorem~\ref{theorem:Airyunivcorrsat} follows from Theorem~\ref{theorem:Airyunivcorrvoid} with the use of
the relation
\eqref{eq:correlationdual} connecting the correlation functions of the
particle and hole (dual) ensembles.  One also uses Proposition~\ref{prop:flip} to change the square-root behavior of
the equilibrium measure density near the upper constraint into square-root vanishing for the equilibrium measure density corresponding to the dual ensemble.

Now we turn our attention to the proof of Theorem~\ref{thm:Airydet}.
Here we present in detail a proof of \eqref{eq:Airydetlimit}. The proof of
\eqref{eq:Airydetlimit2} is analogous and is left to the reader. The starting point is the
fact that, according to \eqref{eq;Am2}, for any real $s$,
\begin{equation}\label{eq:largest}
\begin{array}{l}\displaystyle
\mathbb{P}\left((x_{\rm min}-\alpha)\cdot(\pi NcB^L_\alpha)^{2/3}\ge -s\right) \\\\
\displaystyle\hspace{0.4 in}=\,\,\,  \mathbb{P} \left( x_{\rm min} \ge \alpha-\frac{s}{(\pi NcB^L_\alpha)^{2/3}}\right)\\\\
  \displaystyle \hspace{0.4 in}=\,\,\,\mathbb{P}\left(\text{there are no particles
  at any nodes $x_{N,j}$ satisfying $\displaystyle x_{N,j}<\alpha-\frac{s}{(\pi NcB^L_\alpha)^{2/3}}$}\right) \\\\
  \displaystyle \hspace{0.4 in}=\,\,\,\det (\mathbb{I}-K_{N,k}|_{L_s}),
\end{array}
\end{equation}
where $L_s:=\{y\in X_N \,\,\,\text{such that}\,\,\, y<\alpha-s/(\pi NcB^L_\alpha)^{2/3}\}$ is the (finite, for each $N$) set of nodes that lie strictly to the left of $\alpha-s/(\pi NcB^L_\alpha)^{2/3}$. Since the right-hand side of \eqref{eq:largest}
is the determinant of a finite matrix that we would like to compare with
a Fredholm determinant, we will
first define an
integral operator $\tilde{\cal A}_N|_{[s,\infty)}$ acting on $L^2[s,\infty)$
with a kernel $\tilde{A}_N(\xi,\eta)$ such that the Fredholm determinant $\det(1-\tilde{\cal A}_N|_{[s,\infty)})$ has precisely the same
value for each $N$ as the matrix determinant in \eqref{eq:largest}.  Moreover, it will be obvious from the construction that the kernel $\tilde{A}_N(\xi,\eta)$ will approximate the Airy kernel $A(\xi,\eta)$ at least pointwise.

Let $M(s)$ denote the index of the rightmost node $x_{N,M(s)}$ lying
strictly to the left of $\alpha-s/(\pi NcB^L_\alpha)^{2/3}$, and let
$x_{N,-1}<a$ be defined by
\begin{equation}
\int_{x_{N,-1}}^a\rho^0(x)\,dx=\frac{1}{2N}\,.
\end{equation}
We define a kernel $\tilde{A}_N(\xi,\eta)$ on $[s,\infty)\times [s,\infty)$
by setting
\begin{equation}
\begin{array}{l}\displaystyle
\tilde{A}_N(\xi,\eta):= \frac{N^{1/3}}{(\pi c B^L_\alpha)^{2/3}}\sqrt{\rho^0\left(\alpha-\frac{\xi}{(\pi NcB^L_\alpha)^{2/3}}\right)
\rho^0\left(\alpha-\frac{\eta}{(\pi NcB^L_\alpha)^{2/3}}\right)}K_{N,k}(x_{N,i},x_{N,j})
\\\\
\hspace{0.4 in}\displaystyle
\begin{array}{lrcccl}
\text{if} & (\pi NcB^L_\alpha)^{2/3}(\alpha-x_{N,i})&\le &\xi & < &(\pi NcB^L_\alpha)^{2/3}(\alpha-x_{N,i-1})\\
\text{and} & (\pi NcB^L_\alpha)^{2/3}(\alpha-x_{N,j}) &\le &\eta  &< &(\pi NcB^L_\alpha)^{2/3}(\alpha-x_{N,j-1})\,,
\end{array}
\end{array}
\end{equation}
for all pairs of integers $i$ and $j$ satisfying $0\le i,j \le M(s)$, and $\tilde{A}_N(\xi,\eta):=0$ for all
other $\xi\in [s,\infty)$ and $\eta\in [s,\infty)$.

By a direct computation, we have for each positive integer $p$,
\begin{equation}
\begin{array}{l}\displaystyle
\int_s^\infty\dots\int_s^\infty \det(\tilde{A}_N(\xi_m,\xi_n))_{1\le m,n \le  p} \,d\xi_1\dots d\xi_p \\\\
\displaystyle\hspace{0.4 in}=\,\,\,
\sum_{i_1=0}^{M(s)}\dots\sum_{i_p=0}^{M(s)}\det(K_{N,k}(x_{N,i_m},x_{N,i_n}))_{1\le m,n\le p}\left\{N\int_{x_{N,i_1-1}}^{x_{N,i_1}}\rho^0(x_1)\,dx_1
\cdots N\int_{x_{N,i_p-1}}^{x_{N,i_p}}\rho^0(x_p)\,dx_p\right\}\\\\
\hspace{0.4 in}\displaystyle =\,\,\,
\sum_{i_1=0}^{M(s)}\dots\sum_{i_p=0}^{M(s)}\det(K_{N,k}(x_{N,i_m},x_{N,i_n}))_{1\le m,n\le p}\,.
\end{array}
\end{equation}
This calculation uses the quantization rule \eqref{eq:BS} that defines
the positions of the nodes in terms of the function $\rho^0(x)$.  The
infinite series formula for the Fredholm determinant then implies that
$\det(1-\tilde{\cal
  A}_N|_{[s,\infty)})=\det(\mathbb{I}-K_{N,k}|_{L_s})$.

Therefore to prove \eqref{eq:Airydetlimit}, we need to show that
as $N\to\infty$,
\begin{equation}
  \tilde{\cal A}_N|_{[s,\infty)} \to {\cal A}|_{[s,\infty)} \qquad \text{in trace norm.}
\end{equation}
Since $\tilde{\cal A}_N|_{[s,\infty)}$ and ${\cal A}|_{[s,\infty)}$ are both
positive trace class operators, the following two conditions
\cite{SimonTr} imply convergence in trace norm:
\begin{itemize}
\item[(a)]
  ${\rm tr}\, \tilde{\cal A}_N|_{[s,\infty)} \to {\rm tr}\, {\cal A}|_{[s,\infty)}$
\item[(b)]
$\tilde{\cal A}_N|_{[s,\infty)} \to {\cal A}|_{[s,\infty)}, \qquad
  \text{in the weak-$*$ topology.}$
\end{itemize}

For the purpose of establishing these two conditions, the following
properties of the kernels $\tilde{A}_N(\xi,\eta)$ and $A(\xi,\eta)$
are essential ingredients.
\begin{Lemma}
\label{lem;tildeA}
For each fixed $\xi$ and $\eta$,
\begin{equation}
\lim_{N\rightarrow\infty}\tilde{A}_N(\xi,\eta)=A(\xi,\eta)\,.
\end{equation}
Also, there are positive constants $C$ and $D$ such that the estimate
\begin{equation}
\left|\tilde{A}_{N}(\xi,\eta)\right|\le Ce^{-D(|\xi|^{3/2}+|\eta|^{3/2})}\,,
\end{equation}
holds for all $\xi>s$ and $\eta>s$ (the constants $C$ and $D$ depend
on $s$ but not on $N$).  An estimate of the same form holds with
$\tilde{A}_N(\xi,\eta)$ replaced by $A(\xi,\eta)$.
\end{Lemma}
\begin{proof}
  The pointwise convergence follows from Lemma~\ref{lem;Kedgeinner},
  since $\xi$ and $\eta$ fixed corresponds to $x-\alpha$ and
  $y-\alpha$ of order $N^{-2/3}$.  To obtain the claimed estimates,
  one uses Lemma~\ref{lem;Kedgeuniform} when $\xi$ and $\eta$ are of
  order $N^{2/3}$ such that the corresponding values of $x$ and $y$
  are in the disc $D_\Gamma^{\nab,L}$ surrounding the band edge
  $\alpha$.  When $\xi$ and $\eta$ are such that the corresponding $x$
  and $y$ values are both outside the disc, one uses
  Lemma~\ref{lem;Kvoiduniform} to obtain an exponential estimate in
  terms of the variables $\xi$ and $\eta$.  Finally, when $x$ is in
  the disc and $y$ is outside the disc (or vice-versa), we may use the
  exact representation given by Proposition~\ref{prop:Kinoutexact} and
  similar calculations.
\end{proof}

We first prove (a). From our definition of $\tilde{\cal A}_N|_{[s,\infty)}$,
\begin{equation}\label{eq:trofU}
\begin{array}{rcl}\displaystyle
  {\rm tr}\, \tilde{\cal A}_N|_{[s,\infty)} &=&\displaystyle \int_{[s,\infty)}\tilde{A}_N(\xi,\xi)\,d\xi\\\\
  &=&\displaystyle
  \int_{[s,N^{1/6})}\tilde{A}_N(\xi,\xi)\,d\xi +
  \int_{[N^{1/6},\infty)}
  \tilde{A}_N(\xi,\xi)\,d\xi\,.
  \end{array}
  \end{equation}
Applying Lemma~\ref{lem;tildeA} we see that the second integral is
exponentially small as $N\rightarrow\infty$. On the other hand,
from the definition of $\tilde{A}_N$, the first integral satisfies
\begin{equation}
  \sum_{i=M(N^{1/6})+1}^{M(s)} K_{N,k}(x_{N,i}, x_{N,i}) \le
  \int_s^{N^{1/6}}\tilde{A}_N(\xi,\xi)\,d\xi
  \le
  \sum_{i=M(N^{1/6})}^{M(s)+1} K_{N,k}(x_{N,i}, x_{N,i})\,.
\end{equation}
Using Lemma~\ref{lem;Kedgeinner} and the fact that each of the above sums
consists of $O(N^{1/2})$ terms, we find
\begin{equation}
\sum_{i=M(N^{1/6})+1}^{M(s)}
A(\xi^{(i)}_N,\xi^{(i)}_N)\,\Delta\xi + O(N^{-1/6})\le
\int_s^{N^{1/6}}\tilde{A}_N(\xi,\xi)\,d\xi \le
\sum_{i=M(N^{1/6})}^{M(s)+1}
A(\xi^{(i)}_N,\xi^{(i)}_N)\,\Delta\xi + O(N^{-1/6})\,,
\end{equation}
where 
\begin{equation}
\xi_N^{(i)}=(\pi Nc B^L_\alpha)^{2/3}(\alpha-x_{N,i})\hspace{0.2 in}\text{and}
\hspace{0.2 in} \Delta\xi:=\frac{(\pi c B^L_\alpha)^{2/3}}{N^{1/3}\rho^0(\alpha)}\,.
\end{equation}
Given the asymptotic equal spacing of $\Delta\xi$ between consecutive
points $\xi_N^{(i)}$ in the limit $N\rightarrow\infty$, we see that both sums
above are in fact Riemann sums:
\begin{equation}
\lim_{N\rightarrow\infty}\sum_{i=M(N^{1/6})+1}^{M(s)}A(\xi_N^{(i)},\xi_N^{(i)})\,\Delta\xi = 
\lim_{N\rightarrow\infty}\sum_{i=M(N^{1/6})}^{M(s)+1}A(\xi_N^{(i)},\xi_N^{(i)})\,\Delta\xi 
=\int_s^\infty A(\xi,\xi)\,d\xi\,,
\end{equation}
which proves (a).

In order to check the condition (b), we need to show that for any
$f, g\in L^2[s,\infty)$
\begin{equation}\label{eq:Airycheck2}
  \int_s^\infty \int_s^\infty f(\xi)^* \tilde{A}_N(\xi,\eta)  g(\eta) \,d\xi\, d\eta
  \to \int_s^\infty \int_s^\infty f(\xi)^* A(\xi,\eta)g(\eta) \,d\xi\,
  d\eta\,,
\end{equation}
as $N\rightarrow\infty$, where the asterisk denotes complex conjugation.
But from Lemma~\ref{lem;tildeA}, 
\begin{equation}
f(\xi)^*\tilde{A}_N(\xi,\eta)g(\eta)\to f(\xi)^*A(\xi,\eta)g(\eta)
\end{equation}
as $N\rightarrow\infty$ for almost every $\xi$ and $\eta$, and also
$|f(\xi)^*\tilde{A}_N(\xi,\eta)g(\eta)|\le C|f(\xi)|e^{-D|\xi|^{3/2}}
|g(\eta)|e^{-D|\eta|^{3/2}}$, a bound that is independent of $N$.  By
Cauchy-Schwarz,
\begin{equation}
\int_s^\infty \int_s^\infty C|f(\xi)|e^{-D|\xi|^{3/2}}
|g(\eta)|e^{-D|\eta|^{3/2}}\,d\xi\,d\eta
\le C\|f\|_2\|g\|_2\int_s^\infty e^{-2D|x|^{3/2}}\,dx<\infty
\end{equation}
so the desired result follows from the Lebesgue Dominated Convergence
Theorem.
Hence both conditions (a) and (b) hold and this completes the
proof of \eqref{eq:Airydetlimit}.

The proof of Theorem~\ref{thm:Airydetforhole} follows from that of Theorem~\ref{thm:Airydet} by duality.

\appendix

\section{The Explicit Solution of Riemann-Hilbert Problem~\ref{rhp:theta}}
\label{sec:thetasolve}
\subsection{Obtaining piecewise constant jump matrices:  the transformation
$\dot{\mat{X}}(z)\to\mat{Y}^\sharp(z)$.}
The first step in solving Riemann-Hilbert Problem~\ref{rhp:theta}
is to introduce a change of variables leading to a piecewise-constant
jump matrix.
%Let $\alpha_j$ and $\beta_j$ denote respectively the left and right %endpoints of the band interval $I_j\subset\Sigma_{\rm model}$, for %$j=0,1,\dots,G$.  
Suppose that $h(z)$ is a function analytic for $z\in\mathbb{C}\setminus (-\infty,\beta_G]$ and consider the change of variables \label{symbol:matrixY}
\begin{equation}
\mat{Y}(z):=\dot{\mat{X}}(z)e^{(\kappa g(z)-h(z))\sigma_3}\,.
\label{eq:dotXtoY}
\end{equation}
Then, since by definition $-i\phi_{\Gamma_j}=\kappa g_+(z)-\kappa g_-(z)$ when
$z$ is in any gap $\Gamma_j$, and for $z\in (-\infty,\beta_G)$ setting
$h_\pm(z):=\lim_{\epsilon\downarrow 0}h(z\pm i\epsilon)$, we have the jump condition
\begin{equation}
\mat{Y}_+(z)=\mat{Y}_-(z)\left(\begin{array}{cc}
e^{iN\theta_{\Gamma_j}-h_+(z)+h_-(z)} & 0 \\\\ 0 & e^{-iN\theta_{\Gamma_j}+h_+(z)-h_-(z)}\end{array}\right)
\end{equation}
for $z\in\Gamma_j$ for $j=1,\dots,G$, and
\begin{equation}
\mat{Y}_+(z)=\mat{Y}_-(z)\left(\begin{array}{cc} 0 &
-ie^{\gamma-\eta(z)+h_+(z)+h_-(z)} \\\\
-ie^{\eta(z)-\gamma-h_+(z)-h_-(z)} & 0\end{array}\right)
\end{equation}
for $z$ in any band $I_j$ for $j=0,1,\dots,G$.
Finally, since for all real $z<\alpha_0$ we have $g_+(z)-g_-(z)=2\pi i$, we have introduced a new discontinuity into $\mat{Y}(z)$ by the change of
variables (\ref{eq:dotXtoY}):
\begin{equation}
\mat{Y}_+(z)=\mat{Y}_-(z)\left(\begin{array}{cc} e^{2\pi i\kappa-h_+(z)+h_-(z)} & 0 \\\\ 0 & e^{-2\pi i\kappa +h_+(z)-h_-(z)}\end{array}\right)
\end{equation}
for $z\in (-\infty,\alpha_0)$.

In order to arrive at a problem with piecewise-constant jump matrices that is still normalized to the identity matrix as $z\rightarrow\infty$, we thus insist that $h(z)$ be the solution of the following scalar Riemann-Hilbert problem:
\begin{rhp}
Find a scalar function $h(z)$ with the following properties:
\begin{enumerate}
\item {\bf Analyticity}: $h(z)$ is an analytic function of $z$ for
$z\in\mathbb{C}\setminus (-\infty,\beta_G]$.
\item {\bf Normalization}:  As $z\rightarrow\infty$,
\begin{equation}
h(z)=\kappa\log(z)+O\left(\frac{1}{z}\right)\,.
\label{eq:hnorm}
\end{equation}
\item {\bf Jump Conditions}:  $h(z)$ takes piecewise-continuous boundary values
on $(-\infty,\beta_G]$ with jump discontinuities only allowed at the band endpoints.  For real $z$, let $h_\pm(z):=\lim_{\epsilon\downarrow 0}h(z\pm i\epsilon)$.  For $z$ in the gap 
$\Gamma_j=(\beta_{j-1},\alpha_j)$, $j=1,\dots,G$, the boundary values
satisfy
\begin{equation}
h_+(z)-h_-(z)=ic_j
\label{eq:hgapjump}
\end{equation}
where $c_1,\dots,c_G$ are some real constants.  For $z$ in any  band $I_j=(\alpha_j,\beta_j)$, $j=0,\dots,G$,
the boundary values satisfy
\begin{equation}
h_+(z)+h_-(z)=\eta(z)-\gamma\,,
\label{eq:hbandjump}
\end{equation}
where $\gamma$ is a real constant (the same constant for all bands).  Finally, for real $z<\alpha_0$,
\begin{equation}
h_+(z)-h_-(z)=2 \pi i\kappa\,.
\label{eq:hleftjump}
\end{equation}
\end{enumerate}
The determination of the constants $c_1,\dots,c_G$ and the constant $\gamma$ is part of the problem.
\label{rhp:h}
\end{rhp}
%This Riemann-Hilbert problem only has a solution if the constant $\gamma$,
%that was introduced for later convenience in the transformation %(\ref{eq:SfromR}), has a particular value (see (\ref{eq:gammadef}) below).
%When $\gamma$ has the correct value, t

To solve Riemann-Hilbert Problem~\ref{rhp:h} it is easiest to first solve for $h'(z)$.  Evidently the function $h'(z)$ should be analytic for $z\in\mathbb{C}\setminus\cup_j I_j$; in each band $I_j$ the boundary values should satisfy $h'_+(z)+h'_-(z)=\eta'(z)$.  As $z\rightarrow\infty$, we require the normalization condition $h'(z)=\kappa/z + O(z^{-2})$.
In order to obtain a formula for $h'(z)$,  
recall the 
analytic function $R(z)$ defined for $z\in\mathbb{C}\setminus\cup_kI_k$
by \eqref{eq:Rofzdefine}, and
%to satisfy
%\begin{equation}
%R(z)^2=\prod_{k=0}^G(z-\alpha_k)(z-\beta_k)\,,\hspace{0.2 in}
%\mbox{and}\hspace{0.2 in}R(z)\sim z^{G+1}\hspace{0.1 in}\mbox{as %$z\rightarrow\infty$.}
%\end{equation}
set 
\begin{equation}
h'(z)=\frac{k(z)}{R(z)}
\end{equation}
to introduce a new unknown function $k(z)$.
Evidently, $k(z)$ must be analytic in $\mathbb{C}\setminus\cup_jI_k$ and its boundary values $k_\pm(z):=\lim_{\epsilon\downarrow 0}k(z\pm i\epsilon)$ for $z$ in a band necessarily satisfy 
\begin{equation}
k_+(z)-k_-(z)=\eta'(z)R_+(z)\,, \hspace{0.2 in}
\text{for $z$ in any band $I_j$\,,}
\end{equation}
since
$R_+(z)+R_-(z)=0$ holds for $z$ in the bands with $R_\pm(z):=\lim_{\epsilon\downarrow 0}R(z\pm i\epsilon)$.  Taking into account the required asymptotic behavior of $k(z)$ for large $z$ implied by \eqref{eq:hnorm}, we solve for $k(z)$ in terms of a Cauchy integral:
\begin{equation}
k(z)=\frac{1}{2\pi i}\int_{\cup_jI_j}\frac{\eta'(x)R_+(x)}{x-z}\,dx +
\kappa z^G + \sum_{p=0}^{G-1} f_pz^p
\end{equation}
This is the general solution for $k(z)$, and this establishes the formula
\eqref{eq:hprime} for $h'(z)$.  The constants $f_0,\dots,f_{G-1}$ would seem at this point to be arbitrary; the next step is therefore to explain how they are determined.

%We have found the formula
%\begin{equation}
%h'(z)=\frac{1}{2\pi iR(z)}\int_{\cup_jI_j}\frac{\eta'(x)R_+(x)}{x-z}\,dx +
%\frac{1}{R(z)}\left[\kappa z^G + \sum_{p=0}^{G-1}f_pz^p\right]\,.
%\label{eq:hprime}
%\end{equation}
Note that the inverse square-root singularities present in $h'(z)$ at the
band endpoints are integrable, so $h(z)$ will indeed have piecewise-continuous boundary values as required.  To complete the solution of Riemann-Hilbert Problem~\ref{rhp:h} we must ensure that the identity $h_+'(z)+h_-'(z)=\eta'(z)$ holding in each distinct band $I_j$ actually implies that $h_+(z)+h_-(z)=\eta(z)-\gamma$ holds with the same integration constant $-\gamma$ in each band.  We thus require that the conditions
\eqref{eq:intszero} all hold.
%for $k=1,\dots,G$,
%\begin{equation}
%\int_{\Gamma_k}h'(z)\,dz = 0\,.
%\label{eq:intszero}
%\end{equation}
Substituting into \eqref{eq:intszero} from \eqref{eq:hprime}, we obtain a square linear system of equations on the unknowns $f_1,\dots,f_G$:
\begin{equation}
\sum_{m=0}^{G-1}f_m\int_{\Gamma_l}\frac{z^m\,dz}{R(z)} = -\int_{\Gamma_l}
\left[\frac{1}{2\pi i}\int_{\cup_jI_j}\frac{\eta'(x)R_+(x)}{x-z}\,dx+
\kappa z^G\right]\frac{dz}{R(z)}\,,\hspace{0.2 in}\text{for $l=1,\dots,G$\,.}
\label{eq:linsysf}
\end{equation}
The linear system (\ref{eq:linsysf}) is invertible.  The determinant of
the coefficient matrix is easily seen by multilinearity to be
\begin{equation}
\det\left(\left\{\int_{\Gamma_l}\frac{s^{m-1}\,ds}{R(s)}\right\}_{1\le l,m\le G}\right) =
\int_{\Gamma_1}\cdots\int_{\Gamma_G}D_G(s_1,\dots,s_G)\frac{ds_G}{R(s_G)}\cdots\frac{ds_1}{R(s_1)}\,,
\end{equation}
where $D_G(s_1,\dots,s_G)$ is the Vandermonde determinant $\det(\{s_l^{m-1}\}_{1\le l,m\le G})$.  Since the gaps $\Gamma_1,\dots,\Gamma_G$ are separated from each other by the bands (that is, $\Gamma_j=(\beta_{j-1},\alpha_j)$ and $\alpha_j<\beta_j$ for all $j$), the strict inequalities  $s_1<s_2<\cdots <s_G$ hold throughout the range of integration, which implies that $D_G(s_1,\dots,s_G)$ is of one sign.  Similarly, the product
$R(s_1)R(s_2)\cdots R(s_G)$ is also of one sign.  This proves that the determinant of (\ref{eq:linsysf}) is nonzero.  Note that the constants $f_0,\dots,f_{G-1}$ solving (\ref{eq:linsysf}) are all real because $R_+(x)$ is purely imaginary in the bands, and $R(z)$ is purely real in the gaps.

With the real constants $f_0,\dots,f_{G-1}$ determined in this way, and taking into account the normalization condition \eqref{eq:hnorm} on $h(z)$ as $z\rightarrow\infty$, we see that $h(z)$ must be given in terms of 
$h'(z)$ by the integral formula \eqref{eq:hformula}.
%have an unambiguous formula for $h(z)$ itself:
%\begin{equation}
%h(z)=\kappa\log(z)+\int_z^\infty\left[\frac{\kappa}{s}-h'(s)\right]\,ds
%\label{eq:hformula}
%\end{equation}
%where the path of integration lies in $\mathbb{C}\setminus (-%\infty,\beta_G]$, 
Note that in \eqref{eq:hformula}, the point at infinity can be approached in any direction since $\kappa/z-h'(z)=O(z^{-2})$ as $z\rightarrow\infty$.
In particular, if we consider $z<\alpha_0$ and take paths of integration
with $\arg(s-\alpha_0)=\pm\pi$ to compute $h_\pm(z)$, then it is easy to see that the condition (\ref{eq:hleftjump}) is satisfied for $z<\alpha_0$.
The formula (\ref{eq:hformula}) clearly satisfies relations of the form
(\ref{eq:hgapjump}) for $j=1,\dots,G$; moreover the constants $c_j$ defined by \eqref{eq:csdef}, with $h(z)$ given by (\ref{eq:hformula}), are all real.  Furthermore, we obtain
the formula \eqref{eq:gammadefine} for the integration constant $\gamma$.  %(\ref{eq:hformula}) determines the integration constant $\gamma$:
%\begin{equation}
%\gamma:= \eta(\beta_G)-2h(\beta_G)\,.
%\label{eq:gammadef}
%\end{equation}
%Thus, in order for Riemann-Hilbert Problem~\ref{rhp:h} to be solvable, we %need to have introduced this particular value of $\gamma$ (which is easily %seen to depend 
Clearly, $\gamma$ only depends on the function $\eta(z)$ and the configuration of endpoints $\alpha_0<\beta_0<\alpha_1<\beta_1<\dots<\alpha_G<\beta_G$.
%with the change of variables (\ref{eq:SfromR})  in %\S~\ref{sec:useofequilibrium}.

Given a configuration of endpoints, the solution $h(z)$ of 
Riemann-Hilbert Problem~\ref{rhp:h} depends additionally on the data $(\kappa,\eta(\cdot))$.  A key property of the function $h(z)$, easily verified by superposition, is the following.
\begin{prop}
Fix a configuration of endpoints.  Let $h_0(z)$ be the solution of Riemann-Hilbert Problem~\ref{rhp:h} corresponding to the data $(0,\eta(\cdot))$ with the real constants in (\ref{eq:hgapjump}) denoted by $c_j^{(0)}$ and with the integration constant in (\ref{eq:hbandjump}) denoted by $\gamma^{(0)}$.  Let $h_1(z)$ be the solution of Riemann-Hilbert Problem~\ref{rhp:h} corresponding to the data $(1,0)$ with the real constants in (\ref{eq:hgapjump}) denoted by $\omega_j$ and with the integration constant in (\ref{eq:hbandjump}) denoted by $\gamma^{(1)}$.  Finally, let $h(z)$ be the solution of Riemann-Hilbert Problem~\ref{rhp:h} corresponding to general data $(\kappa,\eta(\cdot))$, with constants $c_j$ and $\gamma$.  Then,
\begin{equation}
\begin{array}{rcl}
h(z)&=&h_0(z)+\kappa h_1(z)\\\\
c_j&=&c_j^{(0)}+\omega_j\kappa\,,\hspace{0.2 in}\mbox{for $j=1,\dots,G$}\\\\
\gamma&=&\gamma^{(0)}+\gamma^{(1)}\kappa\,.
\end{array}
\end{equation}
The quantities $\omega_j$, which will have the interpretation of frequencies, are independent of $\kappa$ and $\eta(\cdot)$, depending only on the value of the parameter $c\in(0,1)$, the functions $V(\cdot)$ and $\rho^0(\cdot)$, and the corresponding equilibrium measure.  The quantities $c_j^{(0)}$ are similar, but depend additionally on the function $\eta(\cdot)$.
\end{prop}

The function $\mat{Y}(z)$ related to $\dot{\mat{X}}(z)$ by (\ref{eq:dotXtoY}) is analytic for $z\in\mathbb{C}\setminus\Sigma_{\rm model}$ on which it takes
boundary values that are continuous except at the band edges where inverse fourth-root singularities may exist.  The boundary values $\mat{Y}_\pm(z):=\lim_{\epsilon\downarrow 0} \mat{Y}(z\pm i\epsilon)$ for $z\in\Sigma_{\rm model}$ satisfy the jump relations
\begin{equation}
\mat{Y}_+(z)=\mat{Y}_-(z)\left(\begin{array}{cc} e^{iN\theta_{\Gamma_j}}e^{-i(c_j^{(0)}+\omega_j\kappa)} & 0 \\\\
0 & e^{-iN\theta_{\Gamma_j}}e^{i(c_j^{(0)}+\omega_j\kappa)}\end{array}\right)
\end{equation}
for $z\in\Gamma_j$, and
\begin{equation}
\mat{Y}_+(z)=\mat{Y}_-(z)\left(\begin{array}{cc} 0&-i\\\\-i & 0\end{array}
\right)
\end{equation}
for $z\in I_j$.  From (\ref{eq:hnorm}) and the asymptotic relation $g(z)=\log(z)+O(1/z)$ as $z\rightarrow\infty$ which follows from (\ref{eq:rhonorm}), we see that $\mat{Y}(z)=\mathbb{I}+O(z^{-1})$ as $z\rightarrow\infty$.

To find $\mat{Y}(z)$, and thus to explain the solution of Riemann-Hilbert Problem~\ref{rhp:theta}, first consider the matrix $\mat{Y}^\sharp(z)$
\label{symbol:matrixYsharp}
related to $\mat{Y}(z)$ by
\begin{equation}
\mat{Y}^\sharp(z):=\left\{\begin{array}{ll}
\mat{Y}(z)\,, &\hspace{0.2 in}\Im(z)>0\,,\\\\
\mat{Y}(z)\left(\begin{array}{cc} 0 & -i \\\\-i & 0\end{array}\right)\,,
&\hspace{0.2 in}\Im(z)<0\,.
\end{array}\right.
\end{equation}
This matrix only tends to the identity as $z\rightarrow\infty$ with $\Im(z)>0$.  However, the advantage is that now the jump will be
characterized everywhere by piecewise-constant off-diagonal matrices.
Namely, if we introduce the notation $\Gamma_0$ \label{symbol:Gammazero}
for ${\mathbb
R}\setminus\Sigma_{\rm model}=(-\infty,\alpha_0)\cup (\beta_G,\infty)$, then
$\mat{Y}^\sharp(z)$ is continuous and thus analytic for $z$ in any
of the bands $I_j$.  On the other hand, letting $\mat{Y}^\sharp_\pm(z):=\lim_{\epsilon\downarrow 0}\mat{Y}^\sharp(z\pm i\epsilon)$ for real $z$, we see that 
for $z$ in the interval $\Gamma_j$,
\begin{equation}
\mat{Y}^\sharp_+(z)=\mat{Y}^\sharp_-(z)\left(\begin{array}
{cc} 0 & ie^{-iN\theta_{\Gamma_j}}e^{i(c_j^{(0)}+\omega_j\kappa)} \\ \\ie^{iN\theta_{\Gamma_j}}e^{-i(c_j^{(0)}+\omega_j\kappa)} & 0
\end{array}\right)\,,
\end{equation}
for $j=1,2,\dots,G$, and for $z\in\Gamma_0$,
\begin{equation}
\mat{Y}^\sharp_+(z)=\mat{Y}^\sharp_-(z)\left(\begin{array}{cc}
0 & i \\\\ i & 0\end{array}\right)\,.
\end{equation}
Thus, $\mat{Y}^\sharp(z)$ is analytic for $z\in\mathbb{C}\setminus\Sigma_{\rm model}'$, \label{symbol:Sigmamodelprime} where
$\Sigma_{\rm
model}':=\Gamma_0\cup\dots\cup\Gamma_G$.  

\subsection{Construction of $\mat{Y}^\sharp(z)$ by means of hyperelliptic function theory.}
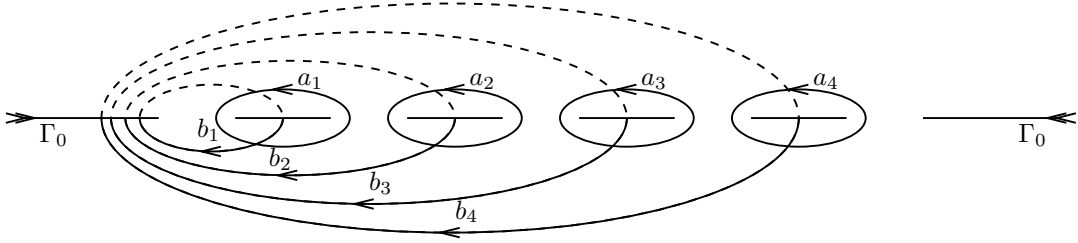
\begin{figure}[h]
\begin{center}
\input{SigmaModelHomology.pstex_t}
\end{center}
\caption{\em The contour $\Sigma_{\rm model}'$ and the homology basis
$a_1,\dots,a_G$ and $b_1,\dots,b_G$ in the associated two-sheeted
Riemann surface $S$.  The sheet on which the cycles are shown with solid curves is that on which $y^S(P)=y(z(P))$ (on the other sheet $y^S(P)=-y(z(P))$).}
\label{fig:homology}
\end{figure}

We will first develop the solution assuming
that $G>0$.  Along with the contour $\Sigma_{\rm model}'$, we associate the
hyperelliptic Riemann surface $S$ 
\label{symbol:RiemannSurface} whose model is two copies of the complex
plane cut and identified along $\Sigma_{\rm model}'$.  Such a surface comes equipped with a function $z:S\to \mathbb{C}$, $P\mapsto z(P)$ 
\label{symbol:z(P)} 
that realizes the identification of each sheet of $S$ with the complex plane.  Each point
$z\in\mathbb{C}$ with the exception of the endpoints $\alpha_0,\dots,\beta_G$ has two preimages on $S$.  Let $y(z)$ be the function analytic for $z\in\mathbb{C}\setminus\Sigma_{\rm model}'$ that satisfies
\begin{equation}
y(z)^2 = -\prod_{k=0}^G (z-\alpha_k)(z-\beta_k)\,,\hspace{0.2 in}
\text{and}\hspace{0.2 in}\text{$y(z)\sim iz^{G+1}$ as
$z\rightarrow\infty$ with 
$\Im(z)>0$\,.}
\end{equation}
This function may be analytically continued to all of the Riemann
surface $S$ with the exception of the two preimages of $z=\infty$ as a
function \label{symbol:yS} $y^S(P)$, for $P\in S$.  We may distinguish the two sheets of
$S$ according to whether $y^S(P)=y(z(P))$ or $y^S(P)=-y(z(P))$, as
long as $z(P)\in\cx\setminus\Sigma_{\rm model}'$.  The polynomial
relation that realizes $S$ as an algebraic curve is then
\begin{equation}
y^S(P)^2 = \prod_{j=0}^G(z(P)-\alpha_j)(z(P)-\beta_j)\,, \hspace{0.2 in}
P\in S\,.
\end{equation}
We next specify on $S$ a basis of homology cycles: \label{symbol:homology}
closed contours $a_k$
encircling $\Gamma_k$ on the sheet where $y^S(P)=y(z(P))$ for $k=1,\dots,G$ in the counterclockwise
direction and conjugate contours $b_k$ oriented in the clockwise
direction and chosen exactly so that $b_k$ intersects only the closed
cycle $a_k$, and precisely once, from the left of $a_k$.  The homology
basis is illustrated in Figure~\ref{fig:homology}.

Also, let a vector of holomorphic differentials $\mat{m}^S(P)\in\cx^G$ be defined for $P\in S$ \label{symbol:holodiffs}
to have components 
\begin{equation}
m^S_p(P):=\frac{z(P)^{p-1}}{y^S(P)}\,dz(P)\,,\hspace{0.2 in}\text{for $p=1,\dots,G$\,.}
\end{equation}
A corresponding vector $\mat{m}(z)\in\cx^G$ may be defined for $z\in\cx\setminus
\Sigma_{\rm model}'$ to have components $m_p(z):=z^{p-1}/y(z)$ for
$p=1,\dots,G$.   We define
a $G\times G$ constant matrix $\mat{A}$ of coefficients so that
\begin{equation}
\oint_{a_j}\mat{A}\mat{m}^S(P)=2\pi
i\mat{e}^{(j)} \,,\hspace{0.2 in}\text{for $k=1,\dots,G$\,,}
\end{equation}
where $\mat{e}^{(j)}$ are the standard unit vectors in
$\cx^G$.  These equations determining the matrix $\mat{A}$ may be written
in the equivalent form \eqref{eq:Adetermine}, which makes the integration
concrete and also makes it clear that the elements of the matrix
$\mat{A}$ real.  We use the notation $\mat{a}^{(1)},\dots,\mat{a}^{(G)}$ to denote (in order) the columns of $\mat{A}$.  The vector ${\mat A}{\mat m}^S(P)$ is the
vector of
{\em normalized holomorphic differentials} on $S$, the normalization being relative to the cycles $a_1,\dots,a_G$.  With $\mat{A}$ so determined,
we construct vectors $\mat{b}^{(j)}\in\cx^G$ by defining
\begin{equation}
\mat{b}^{(j)}:=\oint_{b_j}\mat{A}\mat{m}^S(P)\,,
\label{eq:BdetermineS}
\end{equation}
and we denote by $\mat{B}$ the matrix whose columns are in order
$\mat{b}^{(1)}\dots \mat{b}^{(G)}$.  The definition \eqref{eq:BdetermineS} may be writen in
the equivalent form \eqref{eq:Bdetermine}, which makes the integration concrete.  The matrix $\mat{B}$ is real, symmetric,
and negative definite, and thus $\mat{B}$ defines for $\mat{w}\in\cx^G$ a Riemann theta
function $\Theta(\mat{w})$ by the Fourier series \eqref{eq:thetafunction}.  %From $\mat{B}$ we may define a Riemann theta function
%for $\mat{w}\in\cx^G$ by the Fourier series
%\begin{equation}
%\Theta(\mat{w}):=\sum_{\mat{n}\in\Z^G}
%t_{\mat{n}}e^{\mat{n}^T\mat{w}}\,,\hspace{0.2 in}
%\mbox{with Fourier coefficients}\hspace{0.2 in}
%t_{\mat{n}}:=\exp\left(\frac{1}{2}\mat{n}^T\mat{B}\mat{n}\right)\,.
%\label{eq:thetaFourier}
%\end{equation}

Given a base point $P_0\in S$, the Abel-Jacobi mapping $\mat{w}^S(P)$
is defined by \label{symbol:AbelmapS}
\begin{equation}
\mat{w}^S(P):=\int_{P_0}^P\mat{A}\mat{m}^S
\end{equation}
and since the path of integration on $S$ is not specified, the mapping
is made well-defined by taking the range to be the {\em Jacobian variety} ${\rm Jac}(S) =\cx^G/\Lambda$ where $\Lambda$ is the integer lattice with basis
vectors $2\pi i\mat{e}^{(1)},\dots,2\pi i\mat{e}^{(G)}$ and $\mat{b}^{(1)},
\dots,\mat{b}^{(G)}$.  The definition \eqref{eq:Abelmap} of $\mat{w}(z)$
is a concrete version of the Abel-Jacobi mapping with base point $z=\alpha_0$.
%With the help of the normalized holomorphic differentials, we now
%introduce what is essentially the Abel-Jacobi mapping.  For
%$z\in\cx\setminus{\mathbb R}$, let
%\begin{equation}
%\mat{w}(z):=\int_{\alpha_0}^z\mat{A}\mat{m}(z)\,dz
%\end{equation}
%where for $\Im(z)>0$ the path of integration is restricted to the
%upper half-plane but is otherwise arbitrary, and similarly for
%$\Im(z)<0$ the path is in the lower half-plane.
Since $\mat{m}(z)$ behaves like
\begin{equation}
\mat{m}(z)=-i\frac{{\rm sgn}(\Im(z))}{z^2}\mat{a}^{(G)} + O\left(\frac{1}{z^3}\right)
\hspace{0.2 in}\mbox{as $z\rightarrow\infty$,}
\end{equation}
we see that $\mat{m}(z)$ is integrable at infinity in the two half-planes.
The asymptotic behavior of $\mat{w}(z)$ may be easily computed:
\begin{equation}
\mat{w}(z)=\left\{\begin{array}{ll}
\displaystyle \mat{w}_+(\infty) + i\mat{a}^{(G)}\frac{1}{z} +
O\left(\frac{1}{z^2}\right)\,, &\hspace{0.2 in}\text{as $z\rightarrow\infty$ with $\Im(z)>0$\,,}\\\\
\displaystyle\mat{w}_-(\infty)-i\mat{a}^{(G)}\frac{1}{z}+
O\left(\frac{1}{z^2}\right)\,, &\hspace{0.2 in}\text{as $z\rightarrow\infty$ with $\Im(z)<0$\,,}\end{array}\right.
\end{equation}
where the special values $\mat{w}_\pm(\infty)$ are defined by \eqref{eq:wpminfty}.

For $z\in{\mathbb
R}$, we denote the boundary values taken by $\mat{m}(z)$ and $\mat{w}(z)$ on ${\mathbb R}$ from the 
half-planes 
$\cx_\pm$ by $\mat{m}_\pm(z)$ and $\mat{w}_\pm(z)$.  The boundary values $\mat{w}_\pm(z)$ are continuous
functions with the following expressions:
\begin{equation}
\mat{w}_\pm(z)=\int_{\alpha_0}^z\mat{A}\mat{m}_\pm(x)\,dx\hspace{0.2 in}
\mbox{for $z\in\Gamma_0$,}
\label{eq:AbelpminGamma0}
\end{equation}
\begin{equation}
\mat{w}_\pm(z)=-\frac{1}{2}\mat{b}^{(j)}\mp \sum_{k=1}^{j-1}\pi i\mat{e}^{(k)}
+\int_{\beta_{j-1}}^z\mat{A}\mat{m}_\pm(x)\,dx\hspace{0.2 in}\mbox{for $z\in\Gamma_j$, $j=1,\dots,G$,}
\end{equation}
\begin{equation}
\mat{w}_\pm(z)=\int_{\alpha_0}^z\mat{A}\mat{m}(x)\,dx\,,\hspace{0.2 in}
\mbox{for $z\in I_0$,}
\end{equation}
and
\begin{equation}
\mat{w}_\pm(z)=-\frac{1}{2}\mat{b}^{(j)}\mp\sum_{k=1}^j \pi i\mat{e}^{(k)}
+\int_{\alpha_j}^z\mat{A}\mat{m}(x)\,dx\,, \hspace{0.2 in}
\mbox{for $z\in I_j$, $j=1,\dots,G$.}
\end{equation}
If $z$ lies in the right half of $\Gamma_0$, we interpret the integral
in (\ref{eq:AbelpminGamma0}) as lying always on $\Gamma_0$ and passing
through the point at infinity.  Since $\mat{m}(z)$ is analytic for
$z\in\cx\setminus\Sigma_{\rm model}'$ and when $z\in\Sigma_{\rm
model}'$ we have $\mat{m}_+(z)+\mat{m}_-(z)=0$, the boundary values of
$\mat{w}(z)$ on the real axis are related as follows:
\begin{equation}
\mat{w}_+(z)=-\mat{w}_-(z) - \mat{b}^{(j)}\,,\hspace{0.2 in}\mbox{for $z\in \Gamma_j$, $j=1,\dots,G$,}
\label{eq:AbelInGaps}
\end{equation}
and
\begin{equation}
\mat{w}_+(z)=-\mat{w}_-(z)\,,\hspace{0.2 in}\mbox{for $z\in\Gamma_0$,}
\label{eq:AbelInGap0}
\end{equation}
and
\begin{equation}
\mat{w}_+(z)=\mat{w}_-(z)-\sum_{k=1}^j2\pi i\mat{e}^{(k)}\,,\hspace{0.2 in}
\mbox{for $z\in I_j$, $j=0,\dots,G$.}
\label{eq:AbelInBands}
\end{equation}
By the $2\pi i$-periodicity of $\Theta(\mat{w})$ in each coordinate
direction of $\cx^G$, we see from (\ref{eq:AbelInBands}) that for any
vector $\mat{q}\in \cx^G$, the function \label{symbol:shiftedtheta}
\begin{equation}
f(z;\mat{q}):=\Theta(\mat{w}(z)-\mat{q})
\end{equation}
is analytic in $\cx\setminus\Sigma_{\rm model}'$, and in fact takes
continuous boundary values on $\Sigma_{\rm model}'$.  Moreover, using
the facts
\begin{equation}
\Theta(-\mat{w})=\Theta(\mat{w})\,, \hspace{0.2 in}\mbox{and}\hspace{0.2 in}
\Theta(\mat{w}\pm\mat{b}^{(j)})=e^{-B_{jj}/2}e^{\pm w_j}\Theta(\mat{w})\,,
\end{equation}
holding for all $\mat{w}\in\mathbb{C}^G$ and easily derived directly from the Fourier series
(\ref{eq:thetafunction}), we find from (\ref{eq:AbelInGap0}) that
if for real $z$ we define $f_\pm(z;\mat{q}):=\lim_{\epsilon\downarrow 0}
f(z\pm i\epsilon;\mat{q})$, then
\begin{equation}
f_+(z;\mat{q})=f_-(z;-\mat{q})\hspace{0.2 in}\mbox{for $z\in\Gamma_0$,}
\label{eq:fpm0}
\end{equation}
and from (\ref{eq:AbelInGaps}) that for $j=1,\dots,G$,
\begin{equation}
\begin{array}{rcl}
f_+(z;\mat{q})&=&\displaystyle
e^{-B_{jj}/2-w_{j-}(z)}e^{q_j}f_-(z;-\mat{q})\\\\
&=&\displaystyle
e^{B_{jj}/2+w_{j+}(z)}e^{q_j}f_-(z;-\mat{q})
\end{array}
\label{eq:fpmk}
\end{equation}
when $z\in\Gamma_j$.  Now with the vector $\mat{r}$ defined componentwise by \eqref{eq:vectorr} and with the frequency vector $\mat{\Omega}$ having  components $\omega_1,\dots,\omega_G$,  consider the quotient functions
\label{symbol:gzq}
\begin{equation}
g^\pm(z;\mat{q}):=\frac{f(z;\pm\mat{q}\pm i\mat{r}\mp i\kappa\mat{\Omega})}{f(z;\pm\mat{q})}=
\frac{\Theta(\mat{w}(z)\mp\mat{q}\mp i\mat{r}\pm i\kappa\mat{\Omega})}{\Theta(\mat{w}(z)\mp\mat{q})}\,.
\end{equation}
As long as the denominator does not vanish identically, it will have
at most $G$ zeros on $\cx\setminus\Sigma_{\rm model}'$ (more precisely, replacing $\mat{w}(z)$ by $\mat{w}^S(P)$, the resulting function of $P$ will have exactly $G$
zeros on the Riemann surface $S$, counting multiplicity, and these may occur on either of the two sheets).
These quotient functions $g^\pm(z;\mat{q})$, when $\mat{q}$ is such that
the denominator is not identically zero, are meromorphic functions for
$z\in\cx\setminus\Sigma_{\rm model}'$.  From (\ref{eq:fpm0}) and
(\ref{eq:fpmk}) we then find that
\begin{equation}
g^\pm_+(z;\mat{q})=e^{\pm iN\theta_{\Gamma_j}}e^{\mp i(c_j^{(0)}+\omega_j\kappa)}g^\mp_-(z;\mat{q})
\hspace{0.2 in}
\mbox{for $z\in\Gamma_j$, $j=1,\dots,G$,}
\end{equation}
and for $z\in\Gamma_0$,
\begin{equation}
g^\pm_+(z;\mat{q})=g^\mp_-(z;\mat{q})\,.
\end{equation}
The subscripts again denote boundary values taken as real $z$ is approached from the upper and lower half-planes, just as for $f(z;\mat{q})$.
The
functions $g^\pm(z;\mat{q})$ also take finite values as $z\rightarrow\infty$
separately in each half-plane, and in particular we have the asymptotic formula:
\begin{equation}
\frac{g^\pm(z;\mat{q})}{g^\pm_+(\infty;\mat{q})}=1 + \frac{i}{z}\mat{a}^{(G)}\cdot\left[\frac{\nabla\Theta(\mat{w}_+(\infty)\mp\mat{q}\mp i\mat{r}\pm i\kappa\mat{\Omega})}{\Theta(\mat{w}_+(\infty)\mp\mat{q}\mp i\mat{r}\pm i\kappa\mat{\Omega})}-\frac{\nabla\Theta(\mat{w}_+(\infty)\mp\mat{q})}{\Theta(\mat{w}_+(\infty)\mp\mat{q})}\right] +
O\left(\frac{1}{z^2}\right)\,,
\end{equation}
as $z\rightarrow\infty$ with $\Im(z)>0$, where
$\nabla$ denotes the gradient vector in $\mathbb{C}^G$, and thus $\mat{a}^{(G)}\cdot\nabla$ is a derivative in the direction
of $\mat{a}^{(G)}$.

Comparing the desired jump relations satisfied by $\mat{Y}^\sharp(z)$ in
the gaps with the jump relations satisfied by the functions $g^\pm(z;\mat{q})$, we are led to the strategy of constructing the matrix $\mat{Y}^\sharp(z)$ from
the quotient functions $g^\pm(z;\mat{q})$ by choosing $\mat{q}$ appropriately.
The functions $g^\pm(z;\mat{q})$ have poles in $\cx\setminus
\Sigma_{\rm model}'$ corresponding to the zeros of the denominators, however they are also
typically finite at the endpoints $\alpha_0,\dots,\beta_G$.
Since we can admit mild singularities in $\mat{Y}^\sharp(z)$ at the endpoints, we may introduce additional functional factors with such singularities that also have
zeros that
cancel any poles in $g^\pm(z;\mat{q})$.  Thus, we may seek the matrix elements of $\mat{Y}^\sharp(z)$ in the form of products of  $g^\pm(z;\mat{q})$ with these functional factors
and choosing the vectors $\mat{q}$ appropriately.

To introduce the correct functional factors, recall the function $\lambda(z)$ defined for $z\in\cx\setminus\Sigma_{\rm model}'$ by
\eqref{eq:lambda}
%Let $\lambda(z)$ be defined for $z\in\cx\setminus\Sigma_{\rm model}'$ by
%\begin{equation}
%\lambda(z)^4=\prod_{k=0}^G\frac{z-\alpha_k}{z-\beta_k}\,,\hspace{0.2 in}
%\mbox{and $\lambda(z)\rightarrow 1$ as $z\rightarrow\infty$ with %$\Im(z)>0$.}
%\end{equation}
and the corresponding functions $u(z)$ and $v(z)$ defined in the same
domain by \eqref{eq:uv}.
Noting that $e^{i\pi/4}\lambda(z)$ is a real-analytic function that is positive 
for $z\in {\mathbb R}\setminus\Sigma_{\rm model}'$, we obtain the identity
\begin{equation}
v(z)=-u(z^*)^*\,.  
\end{equation}
%Now set
%\begin{equation}
%u(z):=\frac{1}{2}\left[\lambda(z)+\frac{1}{\lambda(z)}\right]
%\hspace{0.2 in}\mbox{and}
%\hspace{0.2 in}
%v(z):=\frac{1}{2i}\left[\lambda(z)-\frac{1}{\lambda(z)}\right]=-u(z^*)^*\,,
%\end{equation}
%for $z\in\cx\setminus\Sigma_{\rm model}'$.  
Both functions $u(z)$ and $v(z)$ are
analytic throughout their domain of definition.  
The boundary values $u_\pm(z):=\lim_{\epsilon\downarrow 0}u(z\pm i\epsilon)$ and $v_\pm(z):=\lim_{\epsilon\downarrow 0}v(z\pm i\epsilon)$ taken on $\Sigma_{\rm model}'$ have mild singularities at the
endpoints, but are otherwise continuous and satisfy
\begin{equation}
\text{$u_+(z)=-v_-(z)$ and $v_+(z)=u_-(z)$ for $z\in\Sigma_{\rm model}'$\,.}
\end{equation}
Also, we clearly have $u(z)\rightarrow 1$ and $v(z)\rightarrow 0$ as
$z\rightarrow\infty$ with $\Im(z)>0$; more precisely,
\begin{equation}
u(z)=1 + O\left(\frac{1}{z^2}\right)\hspace{0.1 in}\mbox{and}
\hspace{0.1 in}
v(z)=\frac{1}{4iz}\sum_{k=0}^G(\beta_k-\alpha_k) + \frac{1}{8iz^2}
\sum_{k=0}^G(\beta_k^2-\alpha_k^2) + O\left(\frac{1}{z^3}\right)
\end{equation}
as $z\rightarrow \infty$ with $\Im(z)>0$.

To locate the zeros of $u(z)$ and $v(z)$, note that
\begin{equation}
u(z)u(z^*)^*=-u(z)v(z)=\frac{i}{4\lambda(z)^2}\left[\lambda(z)^4-1\right]
\label{eq:uustar}
\end{equation}
which proves that any zeros of $u(z)$ must be real.  But the
right-hand side does not vanish for $z\in{\mathbb
C}\setminus\Sigma_{\rm model}'$, and therefore strictly speaking
$u(z)$ is nonzero in its domain of definition.  However, the
right-hand side of (\ref{eq:uustar}) has exactly one simple zero
$z=x_j$ in the interior of $\Gamma_j$ for each $j=1,\dots,G$, and no
other zeros.  These are precisely the $G$ roots of the polynomial
equation \eqref{eq:polynomialequation}.  Therefore, the boundary values $u_\pm(z)$ can have
zeros.  Parallel arguments apply to $v(z)$.  Since $\lambda_+(z):=\lim_{\epsilon\downarrow 0}\lambda(z+i\epsilon)>0$ for
$z\in\Sigma_{\rm model}'$, we deduce finally that
\begin{equation}
\text{$u_-(x_j)=0$ and $v_+(x_j)=0$ for $j=1,\dots,G$\,,}
\end{equation}
where $x_1,\dots,x_G$ are the roots of \eqref{eq:polynomialequation}
and $\beta_{j-1}<x_j<\alpha_j$.

To build matrix elements of $\mat{Y}^\sharp(z)$ out of products of $g^\pm( z;\mat{q})$ with
$u(z)$ or $v(z)$, the vector $\mat{q}$ should be chosen to align the poles
of the functions $g^\pm(z;\mat{q})$ with the zeros of the boundary values of $u(z)$ or $v(z)$.  An important observation at this point is that the aggregates of points
$D_\pm:=\{x_j\pm i0\}$ form {\em nonspecial divisors}, meaning that if the
expression $(t_{G-1}z^{G-1}+t_{G-2}z^{G-2}+\dots + t_0)/y(z)$ is made
to vanish at all of the $G$ points in either $D_+$ or $D_-$ by an appropriate choice of the coefficients $t_j$, then it vanishes identically.
This implies that the function $f(z;\mat{q})$ will have exactly the
same zeros as $u(z)$, with the same multiplicity, if one takes $\mat{q}=\mat{q}_u$, with $\mat{q}_u$ defined by \eqref{eq:quv} in terms of the Abel-Jacobi mapping evaluated on the divisor $D_-$ and the vector $\mat{k}$ of Riemann constants defined by \eqref{eq:RiemannK}.
%\begin{equation}
%\mat{q}=\mat{q}_u:=\sum_{j=1}^G\mat{w}_-(x_j) + \mat{k}
%\end{equation}
Similarly, the function $f(z;\mat{q}_v)$ has exactly the same zeros as
$v(z)$ with the same multiplicity when $\mat{q}_v$ is defined by \eqref{eq:quv} in terms of
the Abel-Jacobi mapping evaluated on the divisor $D_+$ and the vector
$\mat{k}$.
%and will have exactly the same zeros as $v(z)$ if
%\begin{equation}
%\mat{q}=\mat{q}_v:=\sum_{j=1}^G\mat{w}_+(x_j) + \mat{k}
%\end{equation}
%where the vector $\mat{k}$ of Riemann constants can be defined as
%\begin{equation}
%\mat{k}:=\left\{\begin{array}{ll}
%\displaystyle \pi i\sum_{j\scriptstyle\mbox{ odd}}\mat{e}^{(j)}+
%\frac{1}{2}\sum_{j=1}^G\mat{b}^{(j)} \,, & G\mbox{ odd}\\\\
%\displaystyle \pi i\sum_{j\scriptstyle\mbox{ even}}\mat{e}^{(j)}+
%\frac{1}{2}\sum_{j=1}^G\mat{b}^{(j)} \,, & G\mbox{ even}\,.
%\end{array}\right.
%\label{eq:RiemannK}
%\end{equation}
Note that as a consequence of (\ref{eq:AbelInGaps}) and
(\ref{eq:RiemannK}), $\mat{q}_u+\mat{q}_v=0$ modulo $2\pi i{\mathbb Z}^G$.
Also, $\mat{q}_u-\mat{q}_v^*=0$ modulo $2\pi i{\mathbb Z}^G$.
We therefore can see that all four functions
$u(z)g^\pm(z;\pm \mat{q}_u)$ and
$v(z)g^\pm(z;\pm\mat{q}_v)=v(z)g^\pm(z;\mp\mat{q}_u)$ are analytic for
$z\in\cx\setminus\Sigma_{\rm model}'$ and take continuous boundary
values on $\Sigma_{\rm model}'$ with the exception of the band endpoints
$\alpha_0,\dots,\beta_G$ where they all have negative
one-fourth power singularities.  With the help of these functions, we
may now assemble the solution of Riemann-Hilbert Problem~\ref{rhp:theta}.
First we write down a formula for $\mat{Y}^\sharp(z)$ by setting
\begin{equation}
\mat{Y}^\sharp(z):=\left(\begin{array}{cc}
\displaystyle u(z)\frac{g^+(z;\mat{q}_u)}{g^+_+(\infty;\mat{q}_u)} &
\displaystyle iv(z)\frac{g^-(z;-\mat{q}_v)}{g^-_-(\infty;-\mat{q}_v)}\\\\
\displaystyle iv(z)\frac{g^+(z;\mat{q}_v)}{g^+_-(\infty;\mat{q}_v)} &
\displaystyle u(z)\frac{g^-(z;-\mat{q}_u)}{g^-_+(\infty;-\mat{q}_u)}
\end{array}\right)\hspace{0.2 in}\text{for $G>0$.}
\end{equation}
If $G=0$, then the Riemann theta functions are not necessary, and we
have simply
\begin{equation}
\mat{Y}^\sharp(z):=\left(\begin{array}{cc}
u(z) & iv(z)\\iv(z) & u(z)\end{array}\right)\,,\hspace{0.2 in}\text{for $G=0$.}
\end{equation}

\subsection{The matrix $\dot{\mat{X}}(z)$ and its properties.}
In both cases, $G=0$ and $G>0$, going back to the solution $\dot{\mat{X}}(z)$ of Riemann-Hilbert Problem~\ref{rhp:theta}
requires multiplying $\mat{Y}^\sharp(z)$ on the right by $i\sigma_1$ 
for $\Im(z)<0$ to recover $\mat{Y}(z)$, followed by multiplication on the right by $e^{(h(z)-\kappa g(z))\sigma_3}$ to obtain $\dot{\mat{X}}(z)$.  We have proved
the following.
\begin{prop}
The unique solution of Riemann-Hilbert Problem~\ref{rhp:theta} is given
by the following explicit formulae:
\begin{equation}
\dot{\mat{X}}(z):=\left\{
\begin{array}{ll}\displaystyle
\left(\begin{array}{cc}
\displaystyle u(z)\frac{g^+(z;\mat{q}_u)}{g^+_+(\infty;\mat{q}_u)}e^{h(z)} &
\displaystyle iv(z)\frac{g^-(z;-\mat{q}_v)}{g^-_-(\infty;-\mat{q}_v)}e^{-h(z)}\\\\
\displaystyle iv(z)\frac{g^+(z;\mat{q}_v)}{g^+_-(\infty;\mat{q}_v)}e^{h(z)} &
\displaystyle u(z)\frac{g^-(z;-\mat{q}_u)}{g^-_+(\infty;-\mat{q}_u)}e^{-h(z)}
\end{array}\right)e^{-\kappa g(z)\sigma_3}\,,&\text{$\Im(z)>0$ and $G>0$}\,,
\\\\
\left(\begin{array}{cc}
\displaystyle -v(z)\frac{g^-(z;-\mat{q}_v)}{g^-_-(\infty;-\mat{q}_v)}e^{h(z)}
&
\displaystyle iu(z)\frac{g^+(z;\mat{q}_u)}{g^+_+(\infty;\mat{q}_u)}e^{-h(z)} \\\\
\displaystyle iu(z)\frac{g^-(z;-\mat{q}_u)}{g^-_+(\infty;-\mat{q}_u)}e^{h(z)} &
\displaystyle -v(z)\frac{g^+(z;\mat{q}_v)}{g^+_-(\infty;\mat{q}_v)}e^{-h(z)}
\end{array}\right)e^{-\kappa g(z)\sigma_3}\,, & \text{$\Im(z)<0$ and $G>0$}\,.
\end{array}\right.
\end{equation}
\begin{equation}
\dot{\mat{X}}(z):=\left\{
\begin{array}{ll}\displaystyle
\left(\begin{array}{cc}
\displaystyle u(z)e^{h(z)}&
\displaystyle iv(z)e^{-h(z)}\\\\
\displaystyle iv(z)e^{h(z)} &
\displaystyle u(z)e^{-h(z)}
\end{array}\right)e^{-\kappa g(z)\sigma_3}\,,&\text{$\Im(z)>0$ and $G=0$}\,,
\\\\
\left(\begin{array}{cc}
\displaystyle -v(z)e^{h(z)}
&
\displaystyle iu(z)e^{-h(z)} \\\\
\displaystyle iu(z)e^{h(z)} &
\displaystyle -v(z)e^{-h(z)}
\end{array}\right)e^{-\kappa g(z)\sigma_3}\,, & \text{$\Im(z)<0$ and $G=0$}\,.
\end{array}\right.
\end{equation}
\end{prop}
In verifying the solution, it is useful to observe in addition that
$g^\pm_+(\infty;\mat{q})=g^\mp_-(\infty;\mat{q})$ for any
$\mat{q}\in\cx^G$.

\begin{prop}
Let the coefficients $B_{jk}^{(1)}$ and $B_{jk}^{(2)}$ be defined in terms of the elements of the matrix $\dot{\mat{X}}(z)e^{\kappa(g(z)-\log(z))\sigma_3}$ as follows:
\begin{equation}
\dot{\mat{X}}(z)e^{\kappa(g(z)-\log(z))\sigma_3}=\mathbb{I} + \frac{1}{z}\mat{B}^{(1)} + \frac{1}{z^2}\mat{B}^{(2)} + O\left(\frac{1}{z^3}\right)
\end{equation}
as $z\rightarrow\infty$.  Then
\begin{equation}
B_{12}^{(1)}=\left\{\begin{array}{ll}
\displaystyle\frac{1}{4}(\beta_0-\alpha_0)\,, &\hspace{0.2 in}\text{for $G=0$}\\\\
\displaystyle \frac{\Theta(\mat{w}_+(\infty)-\mat{q}_v+ i\mat{r}-i\kappa\mat{\Omega})\Theta(\mat{w}_-(\infty)-\mat{q}_v)}
{\Theta(\mat{w}_-(\infty)-\mat{q}_v+ i\mat{r}-i\kappa\mat{\Omega})\Theta(\mat{w}_+(\infty)-\mat{q}_v)}\frac{1}{4}\sum_{j=0}^G(\beta_j-\alpha_j)\,,&\hspace{0.2 in}
\text{for $G>0$.}
\end{array}\right.
\end{equation}
\begin{equation}
B_{21}^{(1)}=\left\{\begin{array}{ll}
\displaystyle\frac{1}{4}(\beta_0-\alpha_0)\,, &\hspace{0.2 in}\text{for
$G=0$}\\\\
\displaystyle
\frac{\Theta(\mat{w}_+(\infty)-\mat{q}_v-i\mat{r}+i\kappa\mat{\Omega})\Theta(\mat{w}_-(\infty)-\mat{q}_v)}
{\Theta(\mat{w}_-(\infty)-\mat{q}_v-i\mat{r}+i\kappa\mat{\Omega})\Theta(\mat{w}_+(\infty)-\mat{q}_v)}\frac{1}{4}\sum_{j=0}^G(\beta_j-\alpha_j)\,, &\hspace{0.2 in}\text{for $G>0$.}
\end{array}\right.
\end{equation}
Also, if $G=0$ then
\begin{equation}
B_{11}^{(1)}+\frac{B_{12}^{(2)}}{B_{12}^{(1)}}=\frac{1}{2}(\beta_0+\alpha_0)
\,,
\end{equation}
and if $G>0$ then
\begin{equation}
\begin{array}{rcl}\displaystyle
B_{11}^{(1)}+\frac{B_{12}^{(2)}}{B_{12}^{(1)}}&=&\displaystyle
\frac{1}{2}\frac{\displaystyle\sum_{j=0}^G(\beta_j^2-\alpha_j^2)}{\displaystyle\sum_{j=0}^G(\beta_j-\alpha_j)}\\\\
&&\displaystyle\,\,\,+\,\,\,
\frac{i\mat{a}^{(G)}\cdot\nabla\Theta(\mat{w}_+(\infty)+\mat{q}_v-i\mat{r}+i\kappa\mat{\Omega})}
{\Theta(\mat{w}_+(\infty)+\mat{q}_v-i\mat{r}+i\kappa\mat{\Omega})}-
\frac{i\mat{a}^{(G)}\cdot\nabla\Theta(\mat{w}_+(\infty)+\mat{q}_v)}
{\Theta(\mat{w}_+(\infty)+\mat{q}_v)}\\\\
&&\displaystyle\,\,\,+\,\,\,
\frac{i\mat{a}^{(G)}\cdot\nabla\Theta(\mat{w}_+(\infty)-\mat{q}_v+i\mat{r}-i\kappa\mat{\Omega})}
{\Theta(\mat{w}_+(\infty)-\mat{q}_v+i\mat{r}-i\kappa\mat{\Omega})}-
\frac{i\mat{a}^{(G)}\cdot\nabla\Theta(\mat{w}_+(\infty)-\mat{q}_v)}
{\Theta(\mat{w}_+(\infty)-\mat{q}_v)}\,.
\end{array}
\end{equation}
\label{prop:dotXasymp}
\end{prop}
\begin{proof}
This follows directly from the explicit formulae for $\dot{\mat{X}}(z)$ and
the fact that $\mat{q}_u+\mat{q}_v=0$ modulo $2\pi i\mathbb{Z}^G$.  It is also perhaps useful to point out that it is never necessary to use an explicit expression for the leading coefficient of $h(z)-\kappa\log(z)$ as $z\rightarrow\infty$; although this coefficient appears in $B_{12}^{(2)}$ and also in $B_{11}^{(1)}$ it cancels out of the particular combination $B_{11}^{(1)}+B_{12}^{(2)}/B_{12}^{(1)}$.
\end{proof}

\subsubsection{Completion of the proof of Proposition~\ref{prop:dotXbound}.}
The uniform boundedness of $\dot{\mat{X}}(z)$ for $z$ bounded away from any band endpoints $\alpha_0,\dots,\beta_G$ follows from the corresponding property of the functions $u(z)$ and $v(z)$, and the manner in which the
large parameter $N$ enters into the argument of the Riemann theta functions as a
real phase $\mat{r}-\kappa\mat{\Omega}$ that is independent of $z$ (recall that $\Theta(\mat{w})$ is periodic with period $2\pi$ in each imaginary coordinate direction in $\cx^G$).  To see the independence of the combination $\dot{\mat{X}}(z)e^{\kappa g(z)\sigma_3}$ from the arbitrary locations of any transition points in $Y_N$, observe that the combination $\mat{Y}(z)e^{-h(z)\sigma_3}$ can only involve the function $g(z)$ through
the endpoints of the bands, which clearly do not depend on any arbitrary choice of transition points in transition bands.
%, and through the constant values $\theta_{\Gamma_j}$ of the function %$\theta(z)$ in the voids and saturated regions.  But the constants %$\theta_\Gamma$ only appear in terms of the form $e^{\pm %iN\theta_{\Gamma_j}}$, and using (\ref{eq:yquantize}) it follows that these %terms are independent of the choice of transition points.

\subsubsection{Completion of the proof of Proposition~\ref{prop:dotXsymmetry}.}
Note that for real $z$, the product $p(z):=\dot{X}_{11}(z)\dot{X}_{12}(z)$
can equivalently be written in terms of the elements of the matrix $\mat{Y}^\sharp(z)$ as $p(z)=Y^\sharp_{11+}(z)Y^\sharp_{12+}(z)$.  Now, $\mat{Y}^\sharp(z)$ satisfies the symmetry
\begin{equation}
\mat{Y}^\sharp(z)=\mat{Y}^\sharp(z^*)^*\left(\begin{array}{cc} 0 & -i\\\\
-i & 0\end{array}\right)\,.
\label{eq:Ysharpsymmetry}
\end{equation}
Indeed, the left and right-hand sides of (\ref{eq:Ysharpsymmetry}) both have the same asymptotic behavior as $z\rightarrow\infty$ regardless of whether
$\Im(z)>0$ or $\Im(z)<0$, and satisfy the same jump conditions for $z\in\Gamma_j$, $j=0,\dots G$.  In other words, both sides of (\ref{eq:Ysharpsymmetry}) solve the same Riemann-Hilbert problem.  A uniqueness argument based on Liouville's Theorem thus proves (\ref{eq:Ysharpsymmetry}).  Using (\ref{eq:Ysharpsymmetry}) we may also write $p(z)=-iY^\sharp_{11+}(z)Y^\sharp_{11-}(z)^*$ when $z$ is real.

Let us now consider the zeros of the function $Y^\sharp_{11}(z)$.
Recall that $u(z)$ is nonzero for $z\in\mathbb{C}\setminus\Sigma_{\rm model}'$,
and for $z\in\Sigma_{\rm model}'$ we have that $u_+(z)$ is bounded away from zero while $u_-(z)$ vanishes only at a single point $x_j$ in each interior gap $\Gamma_j=(\beta_{j-1},\alpha_j)$ for $j=1,\dots,G$.
This shows that if $G=0$ then (since there are no interior gaps) $p(z)$ is strictly nonzero for $z\in\Sigma_{\rm model}'$.  If $G>0$, then the zeros of $u_-(z)$ on $\Sigma_{\rm model}'$ are cancelled (by construction) by corresponding zeros of the entire function $\Theta(\mat{w}(z)-\mat{q}_u)$
in the denominator of $Y^\sharp_{11}(z)$.  Thus for $G>0$ the zeros of $u(z)$ are precisely the zeros of the numerator $\Theta(\mat{w}(z)-\mat{q}_u-i\mat{r}+i\kappa\mat{\Omega})$.  Recall that $\mat{r}$ and $\kappa\mat{\Omega}$ are real.  By Jacobi inversion theory, it can thus be shown that $\Theta(\mat{w}(z)-\mat{q}_u-i\mat{r}+i\kappa\mat{\Omega})$ has exactly one zero on either the upper or lower edge of each cut $\Gamma_j$ for $j=1,\dots,G$ (in a nongeneric situation the zero may lie at one or the other endpoint of $\Gamma_j$).  As the parameter $\kappa$ is varied continuously with $\mat{q}_u$, $\mat{r}$, and $\mat{\Omega}$ held fixed, the $G$ zeros of $\Theta(\mat{w}(z)-\mat{q}_u-i\mat{r}+i\kappa\mat{\Omega})$ oscillate about the cuts (moving one way along the upper edge and the other way along the lower edge) in a quasiperiodic fashion.  The actual behavior of the zeros is more complicated since the parameter $\kappa$ can in fact only be incremented by integers (recall that the polynomial degree $k$ is $cN+\kappa$); thus if the frequency vector $\mat{\Omega}$ is a rational multiple of a lattice vector in $2\pi \mathbb{Z}^G$ then the motion of the zeros will be periodic rather than quasiperiodic.  In an extremely nongeneric situation the zeros of $\Theta(\mat{w}(z)-\mat{q}_u-i\mat{r}+i\kappa\mat{\Omega})$ may be located at the endpoints of the interior gaps $\Gamma_j$ for all admissible $\kappa$.
Therefore,  if $\dot{X}_{11\pm}(z)$ is bounded away from zero for $z\in\Gamma_j$, $j=1,\dots,G$, then $\dot{X}_{11\mp}(z)$ necessarily has a simple zero in the interior of $\Gamma_j$.  This shows that the product $p(z)$ vanishes at exactly one point $z=z_j$ in the interval $[\beta_{j-1},\alpha_j]$ for $j=1,\dots,G$ and $G>0$.

Now, for real $z$ we write
\begin{equation}
p(z)=-iu_+(z)u_-(z)^*\left[\frac{Y^\sharp_{11+}(z)}{u_+(z)}\right]
\left[\frac{Y^\sharp_{11-}(z)}{u_-(z)}\right]^*\,.
\end{equation}
It follows from (\ref{eq:uustar}) that for $z$ real, $-iu_+(z)u_-(z)^*$ is a real function that satisfies
\begin{equation}
\begin{array}{rcll}
-iu_+(z)u_-(z)^*&<&0\,,&\hspace{0.2 in}\mbox{for $z<\alpha_0$}\\\\
-iu_+(z)u_-(z)^*&>&0\,,&\hspace{0.2 in}\mbox{for $z>\beta_G$}\\\\
-iu_+(z)u_-(z)^*&\rightarrow &-\infty\,,&\hspace{0.2 in}
\mbox{as $z\uparrow\alpha_j$ for $j=0,\dots,G$}\\\\
-iu_+(z)u_-(z)^*&\rightarrow &+\infty\,,&\hspace{0.2 in}
\mbox{as $z\downarrow\beta_j$ for $j=0,\dots,G$.}
\end{array}
\end{equation}
At the same time, we have the existence of the following finite limits
\begin{equation}
\begin{array}{rcl}
A_j&:=&\displaystyle\lim_{z\uparrow\alpha_j}\left[\frac{Y^\sharp_{11+}(z)}{u_+(z)}\right]
\left[\frac{Y^\sharp_{11-}(z)}{u_-(z)}\right]^*\\\\
B_j&:=&\displaystyle\lim_{z\downarrow\beta_j}\left[\frac{Y^\sharp_{11+}(z)}{u_+(z)}\right]
\left[\frac{Y^\sharp_{11-}(z)}{u_-(z)}\right]^*\,.
\end{array}
\end{equation}
We clearly have $A_j\ge 0$ and $B_j\ge 0$ for all $j=0,\dots,G$; the limits are strictly positive unless one of the zeros $z_j$ occurs at an endpoint.  This proves that $p(z)>0$ for $z<z_j$ in $\Gamma_j$ and that $p(z)<0$ for $z>z_j$ in $\Gamma_j$, for $j=1,\dots,G$, while $p(z)<0$ for $z<\alpha_0$ and $p(z)>0$ for $z>\beta_G$.

\section{Construction of the Hahn Equilibrium Measure:  Proof of Theorem~\ref{theorem:hahn}}
\label{sec:hahnderivation}
\subsection{General strategy.  The one-band ansatz.}
The main idea is to begin with an ansatz that there is only one band,
a subinterval of $(0,1)$ of the form $(\alpha,\beta)$, where $\alpha$ and $\beta$ are to
be determined. Then using the ansatz we derive formulae for $\alpha$ and
$\beta$, the ``candidate'' equilibrium measure, and the corresponding
Lagrange multiplier. Of course, one must then check that the measure
produced by the ``one-band'' ansatz is consistent with the variational
problem.

The associated field \eqref{eq:fielddef} is
\begin{equation}
\begin{array}{rcl}
\varphi(x) &:=&\displaystyle
V^{\rm Hahn}(x;A,B)+\int_0^1\log|x-y|\rho^0(y)\,dy\\\\
& =& -(A+x)\log(A+x)-(B+1-x)\log(B+1-x) +x\log(x) +(1-x)\log(1-x)
\\\\
&&\displaystyle\hspace{0.3 in} +\,\,\, A\log(A)+(B+1)\log(B+1) -1\,,
\end{array}
\end{equation}
and hence
\begin{equation}
  \varphi'(x)= -\log(A+x)+\log(1+B-x) + \log(x)-\log(1-x)\,.
\end{equation}
In the presumed band, the candidate equilibrium measure (we will refer to
its density as $\psi(x)$)
satisfies the equilibrium condition \eqref{eq:equilibrium}.  Differentiating
this equation with respect to $x$, one finds that
\begin{equation}\label{eq;bandeq}
  \mbox{P. V.} \int_{0}^1 \frac{\psi(y)}{x-y}\, dy =
\frac1{2c}\varphi'(x)
\end{equation}
holds identically for $x$ in the (as yet unknown) band $\alpha<x<\beta$.
Introducing the Cauchy transform of $\psi(x)$,
\begin{equation}\label{}
  F(z) := \int_{0}^1 \frac{\psi(y)}{z-y}\, dy\,,
\hspace{0.2 in}  \mbox{for $z\in \cx\setminus[0,1]$}
\end{equation}
elementary properties of Cauchy integrals imply that if $F_+(x)$ denotes
the boundary value taken on $[0,1]$ from above, and $F_-(x)$ denotes the
corresponding boundary value taken from below, then
\begin{equation}\label{eq;CTproperties}
\begin{array}{rcl}
  F_+(x)- F_-(x) &=& - 2\pi i \psi(x)\,, \\\\
\displaystyle  \frac12\left( F_+(x)+F_-(x) \right) &=&
\displaystyle  \mbox{P. V.} \int_{0}^1 \frac{\psi(y)}{x-y}\, dy\,.
\end{array}
\end{equation}
Also as $z\to\infty$, the condition that $\psi(x)$ should be the
density of a probability measure implies that
\begin{equation}
  F(z) = \frac1{z} + O\left(\frac{1}{z^{2}}\right)\,,
\hspace{0.2 in}\mbox{as $z\rightarrow\infty$}\,.
\end{equation}
According to the one-band ansatz, the remaining intervals $(0,\alpha)$ and
$(\beta,1)$ are either voids or saturated regions.  So at this point, the
one-band ansatz bifurcates into four distinct cases that must be
investigated: these are the four configurations void-band-void,
saturated-band-void, saturated-band-saturated, and
void-band-saturated.  We will work out many of the details in the
void-band-void case, and then show how the analysis changes in the other
three configurations.

\subsection{The void-band-void configuration.}
In both voids $(0,\alpha)$ and $(\beta,1)$, we have the lower constraint in
force: $\psi(x)\equiv 0$. Hence from \eqref{eq;bandeq} and
\eqref{eq;CTproperties}, $F$ necessarily solves the following
scalar Riemann-Hilbert problem: $F(z)$ is analytic in
$\cx\setminus [\alpha,\beta]$ and satisfies the jump condition
\begin{equation}
  F_+(x)+F_-(x) = \frac1{c} \varphi(x)\,, \hspace{0.2 in}
\text{for $\alpha<x<\beta$}\,,
\end{equation}
and as $z\to\infty$,
\begin{equation}\label{eq;Fasym}
  F(z) = \frac1{z} + O\left(\frac{1}{z^2}\right)\,.
\end{equation}
The solution to this Riemann-Hilbert problem is given by the explicit formula
\begin{equation}\label{}
  F(z) = \frac{R(z)}{2\pi ic} \int_\alpha^\beta
  \frac{\varphi'(y)}{R_+(y)(y-z)}\,dy\,,\hspace{0.2 in}
  \text{for $z\in\mathbb{C}\setminus [\alpha,\beta]$\,,}
\label{eq;FVBVformula}
\end{equation}
where the subscript ``$+$'' indicates a boundary value taken from the upper half-plane, 
$R(z)^2 = (z-\alpha)(z-\beta)$, and
the square root $R(z)$ is defined to be analytic in
$\cx\setminus [\alpha,\beta]$ with the
condition that $R(z) \sim z$ as $z\to\infty$. The asymptotic
condition \eqref{eq;Fasym} on $F$ now implies, by explicit asymptotic
expansion of the formula \eqref{eq;FVBVformula}, the following two conditions
on the endpoints $\alpha$ and $\beta$:
\begin{equation}
\begin{array}{rcl}
\displaystyle   -\frac1{2\pi i}\int_\alpha^\beta \frac{\varphi'(y)}{R_+(y)}\, dy &=& 0
\\\\
\displaystyle   -\frac1{2\pi i}\int_\alpha^\beta \frac{y\varphi'(y)}{R_+(y)}\, dy &=&
c\,.
\end{array}
\end{equation}
In a standard application of formulae involving Cauchy integrals,
one can evaluate the above integrals and find that the endpoint equations
are equivalent to
\begin{eqnarray}
\label{eq;SDeq1}
  \cosh^{-1} \left(\frac{A+s}{d}\right) -
  \cosh^{-1}\left(\frac{B+1-s}{d}\right) - \cosh^{-1}\left(\frac{s}{d}\right) +
  \cosh^{-1} \left(\frac{1-s}{d}\right) &=& 0 \\
\label{eq;SDeq2}
  \sqrt{(A+s)^2-d^2} - \sqrt{s^2-d^2} +\sqrt{(1+B-s)^2-d^2} - \sqrt{(1-s)^2-d^2}
  &=& 2c+A+B \,,
\end{eqnarray}
where $s$ and $d$ are defined by
\begin{equation}
  s:= \frac{\beta+\alpha}{2}\,, \qquad d:=\frac{\beta-\alpha}{2}\,.
\end{equation}
Now using the addition formula $\cosh^{-1}(a)\pm
\cosh^{-1}(b)=\cosh^{-1}(ab\pm\sqrt{(a^2-1)(b^2-1)})$ twice, the
equation \eqref{eq;SDeq1} becomes
\begin{equation}\label{eq;SDeq3}
  (2+A+B)s-\sqrt{(A+s)^2-d^2}\sqrt{s^2-d^2} = B+1 -\sqrt{(B+1-s)^2-d^2}
\sqrt{(1-s)^2-d^2}\,.
\end{equation}
Thus $\alpha$ and $\beta$ are necessarily solutions of the system of equations
\eqref{eq;SDeq2} and \eqref{eq;SDeq3}. Now we take the square of
both sides of \eqref{eq;SDeq2} and add two times
\eqref{eq;SDeq3}, and we find
\begin{equation}\label{}
  \sqrt{(A+s)^2-d^2} - \sqrt{s^2-d^2} =
  \frac{(2c+A+B)^2+A^2-B^2}{2(2c+A+B)}\,.
\end{equation}
To simplify upcoming formulae, we set
\begin{eqnarray}
\label{eq;Wdef}
  W &=& \sqrt{(A+s)^2-d^2} \\
  X &=& \sqrt{s^2-d^2} \\
  Y &=& \sqrt{(B+1-s)^2-d^2} \\
\label{eq;Zdef}
  Z &=& \sqrt{(1-s)^2-d^2}\,.
\end{eqnarray}
With this notation, the equations \eqref{eq;SDeq3} and
\eqref{eq;SDeq2} for $s$ and $d$ (and hence for $\alpha$ and $\beta$)
become
\begin{eqnarray}
\label{eq;SDeq4}
  (2+A+B)s-WX &=& B+1 -YZ \\
\label{eq;SDeq5}
  Y-Z &=& 2c+A+B - (W-X)\,.
\end{eqnarray}
By taking the square of both sides of \eqref{eq;SDeq4} and adding
two times \eqref{eq;SDeq5}, we find
\begin{equation}\label{eq;WXeq1}
  W-X = \frac{(2c+A+B)^2+A^2-B^2}{2(2c+A+B)} =: K\,.
\end{equation}
On the other hand, by the definition of $W$ and $X$,
$(W-X)(W+X)=A^2+2As$, and thus \eqref{eq;WXeq1} implies that $W+Y=
(A^2+2As)/K$. Hence
\begin{eqnarray}
\label{eq;Wsol1}
  2W &=& K + \frac{A^2+2As}{K} \\
\label{eq;Xsol1}
  2X &=& -K + \frac{A^2+2As}{K}\,.
\end{eqnarray}
Also \eqref{eq;SDeq5} implies that $Y-Z= 2c+A+B - K$, and from the
definitions of $Y$ and $Z$, we get $(Y-Z)(Y+Z)=B^2+2B(1-s)$. Therefore we
find
\begin{eqnarray}
\label{eq;Ysol1}
  2Y &=& 2c+A+B-K + \frac{B^2+2B(1-s)}{2c+A+B-K} \\
\label{eq;Zsol1}
  2Z &=& -(2c+A+B-K) + \frac{B^2+2B(1-s)}{2c+A+B-K}\,.
\end{eqnarray}
Substituting \eqref{eq;Wsol1}, \eqref{eq;Xsol1}, \eqref{eq;Ysol1}
and \eqref{eq;Zsol1} into \eqref{eq;SDeq4}, we find a quadratic
equation in $s$:
\begin{equation}\label{eq;theSeq}
  (2+A+B)s - \frac14 \left\{ -K^2+\left( \frac{A^2+2As}{K} \right)^2
  \right\} = B+1 + \frac14 \left\{
  (2c+A+B-K)^2 - \left( \frac{B^2+2B(1-s)}{2c+A+B-K}
  \right)^2 \right\}
\end{equation}
where $K$ is defined in \eqref{eq;WXeq1}. The solutions to this quadratic
equation are
\begin{equation}
  s=s_1 := \frac{A(A+B)+(A+B)(B-A+2)c+(B-A+2)c^2}{(A+B+2c)^2}
\label{eq;SisS1}
\end{equation}
and
\begin{equation}
  s=s_2 := \frac{A(A+B)(1+B)+(A+B)(A+B+2)c+(A+B+2)c^2}{A^2-B^2}\,.
\end{equation}
Since $0\le \alpha+\beta \le 2$, we need to check which of these two roots
actually lie in $[0,1]$. For $s_2$, one sees that if $A< B$, then $s_2<0$,
and if
$A>B$, then by looking at the terms not involving $c$,
\begin{equation}\label{}
  s_2 -1 \ge \frac{A(A+B)(1+B)}{A^2-B^2} -1 =
  \frac{A^2B+AB^2+AB+B^2}{A^2-B^2}>0\,.
\end{equation}
On the other hand, the numerator of $s_1$ can be written as
\begin{equation}
  A(A+B)(1-c) + c \left\{ (A+B)(B+2)-Ac \right\} + (B+2)c^2\,.
\end{equation}
Since $0<c<1$, each term is positive, and thus $s_1>0$. Analogously, the
numerator of $1-s_1$ can be written as
\begin{equation}
  (A+B)B(1-c) + c\bigl\{ A(A+2) +B(A+2-c) \bigr\} +(A+2)c^2\,,
\end{equation}
and each term is positive as $c\in (0,1)$, which implies that
$s_1<1$. Thus $s_1\in (0,1)$, and we have found the root we need.

Substituting \eqref{eq;SisS1} into \eqref{eq;Xsol1} and
\eqref{eq;Zsol1}, we find
\begin{equation}
  X=X_0\,, \qquad Z=Z_0\,
\end{equation}
where
\begin{eqnarray}
  X_0 &:=& \frac{-c^2-(A+B)c+A}{(A+B+2c)^2}  \\
  Z_0 &:=& \frac{-c^2-(A+B)c+B}{(A+B+2c)^2}\,.
\end{eqnarray}
Note that
\begin{eqnarray}
  X_0 >0\,, &\qquad& \text{for $0< c<c_A$} \\
  X_0 <0\,, &\qquad&  \text{for $c>c_A$} \\
  Z_0 >0\,, &\qquad& \text{for $0< c<c_B$} \\
  Z_0 <0\,, &\qquad&  \text{for $c>c_B$} \,.
\end{eqnarray}
Thus the conditions $X, Z>0$ yield conditions on $c$ which are:
\begin{equation}
  0< c < \min(c_A, c_B)\,.
\end{equation}
Only for $c$ satisfying these inequalities can $\alpha$ and $\beta$ be found
from the equations $\beta+\alpha=2s$ and $\alpha\beta= s^2-d^2=Y$. The quadratic
equation for $\alpha$ and $\beta$ is exactly \eqref{eq;qeqLR}, and the explicit
solutions are given by \eqref{eq;L} and \eqref{eq;R}.

With the endpoints determined, the candidate density for the
equilibrium measure in the interesting region $\alpha<x<\beta$ can be obtained
by evaluating $F(z)$ and using
\begin{equation}
  \psi(x) = -\frac{1}{2\pi i} \left( F_+(x) - F_-(x) \right)\,.
\end{equation}
First we evaluate $F(z)$.  For $z\in\mathbb{C}\setminus [\alpha,\beta]$, we have
\begin{equation}
\begin{split}
  F(z) &= \frac{R(z)}{2\pi ic} \int_\alpha^\beta
  \frac{\varphi'(y)}{R_+(y)(y-z)}\,dy \\\\
  &= \frac{R(z)}{4\pi ic} \int_{\Gamma_0}
  \frac{\varphi'(y)}{R(y)(y-z)}\,dy
\end{split}
\end{equation}
where the closed contour $\Gamma_0$ encloses the interval $[\alpha,\beta]$ once
in the clockwise direction, and the inside of $\Gamma_0$ does not
include any points $y$ in the set $\{z\}\cup(-\infty,
0]\cup(1,\infty]$.  Noting that $\varphi'(z)$ is analytic in
$\mathbb{C}\setminus ([-A, 0]\cup [1,1+B])$, we deform the contour of
integration so that the integral over $\Gamma_0$ becomes the integral
over the union of the intervals $[-A,0]$ and $[1,1+B]$. Being careful
with branches of the various multivalued functions involved, we find that
\begin{equation}
\begin{split}
  F(z) = \frac1{2c}\varphi'(z) + \frac{R(z)}{2c}
\left( - \int_{-A}^0 \frac{dy}{R(y)(y-z)}
- \int_1^{1+B} \frac{dy}{R(y)(y-z)}
\right)\,.
\end{split}
\end{equation}
Here when $z$ is in either of the intervals $(-A,0)$ or $(1,1+B)$
where $F(z)$ is supposed to be analytic, the integral is interpreted
as the principal value.  This integral is equal to
\begin{equation}
    F(z) = \frac1{2c}\varphi'(z) - \frac{R(z)}{2c}
\left(  \int_0^A \frac{ds}{\sqrt{(s+\alpha)(s+\beta)}(s+z)}
+ \int_1^{1+B} \frac{ds}{\sqrt{(s-\alpha)(s-\beta)}(s-z)}
\right)\,.
\end{equation}
Now we use the following formula (see, for example, \cite{AbramowitzS65}),
\begin{equation}
 \int \frac{ds}{\sqrt{(s+a)(s+b)}(s+z)}
= \frac2{\sqrt{(z-a)(z-b)}} \left[
\log \left( \sqrt{\frac{s+b}{z-b}} + \sqrt{\frac{s+a}{z-a}} \right)
- \frac12 \log(s+z) \right]
\end{equation}
and evaluate the two integrals exactly.
The result of this calculation is that
for $z\in\mathbb{C}\setminus [\alpha,\beta]$,
\begin{equation}
\begin{split}
  F(z) &= \frac1{c}
\log \left( \sqrt{\frac{1+B-\beta}{z-\beta}} + \sqrt{\frac{1+B-\alpha}{z-\alpha}} \right)
- \frac1{c}
\log \left( \sqrt{\frac{1-\beta}{z-\beta}} + \sqrt{\frac{1-\alpha}{z-\alpha}} \right) \\
&\qquad
- \frac1{c}
\log \left( \sqrt{\frac{A+\beta}{z-\beta}} + \sqrt{\frac{A+\alpha}{z-\alpha}} \right)
+\frac1{c}
\log \left( \sqrt{\frac{\beta}{z-\beta}} + \sqrt{\frac{\alpha}{z-\alpha}} \right)
\,,
\end{split}
\end{equation}
where all the square root functions $\sqrt{w}$ are defined to be
analytic in $w\in\mathbb{C}\setminus (-\infty, 0]$ with the condition
that $\sqrt{w} >0$ for $w>0$, and the logarithm $\log(w)$ is defined
to be analytic in $w\in\mathbb{C}\setminus (-\infty, 0]$ with the
condition that $\log(w) >0$ for $w>1$.  Now using
$\log(a+ib)-\log(a-ib) = 2i \arctan( b/a)$, we obtain, for $x\in
(\alpha,\beta)$, the formula cited in Theorem~\ref{theorem:hahn} for the
equilibrium measure (as mentioned above, the candidate $\psi(x)$ we
have just constructed turns out to be the actual density of the
equilibrium measure).

The Lagrange multiplier $\ell_c$ can then be obtained
from the variational condition \eqref{eq:equilibrium},
which we may evaluate for any $x\in [\alpha,\beta]$.
Therefore, since $\psi(x)$ is only supported in $[\alpha,\beta]$ in the void-band-void
case under consideration,
\begin{equation}
  \ell_c
= -2c \int_\alpha^\beta \log(\beta-s) \psi(s)ds + \varphi(\beta)\,.
\end{equation}
Here we have arbitrarily picked $x=\beta$.  Now using the exact formula
for $\psi(x)$ and the identity
\begin{equation}
  \frac1{\beta-\alpha} \int_\alpha^\beta \log(\beta-s)
\arctan\left( k\sqrt{\frac{\beta-s}{s-\alpha}} \right)\, ds
= \frac{\pi k}{2(1+k)} \left(\log(\beta-\alpha)-1 \right)
+ \frac{\pi k}{1-k^2} \left[ \frac1{k}\log(1+k) -\log 2 \right]\,,
\end{equation}
we obtain the corresponding formula for the multiplier.

\subsection{The saturated-band-void configuration.}
Since $\psi(x)\equiv 1/c$ in the saturated region supposed to
be the interval $(0,\alpha)$, from
\eqref{eq;bandeq} and
\eqref{eq;CTproperties}, the Cauchy transform of the candidate
density $\psi(x)$, $F(z)$, is necessarily the solution of the following
scalar Riemann-Hilbert problem: $F(z)$ is analytic in
$\mathbb{C}\setminus [0,\beta]$ and satisfies the jump conditions
\begin{equation}
  F_+(x)-F_-(x) = -\frac{2\pi i}{c}
\end{equation}
for $0<x<\alpha$, and
\begin{equation}
  F_+(x)+F_-(x) = \frac1{c} \varphi(x)
\end{equation}
for $\alpha<x<\beta$,
and as $z\to\infty$,
\begin{equation}
  F(z) = \frac1{z} + O\left(\frac{1}{z^2}\right)\,.
\end{equation}
The solution is
\begin{equation}
  F(z) = \frac{R(z)}{2\pi ic} \left(
-\int_0^\alpha \frac{2\pi i}{R(y)(y-z)}\, dy
+ \int_\alpha^\beta \frac{\varphi'(y)}{R_+(y)(y-z)}\, dy \right)\,,
\end{equation}
where $R(z)$ denotes the same square root function as before.
The equations for the endpoints $\alpha$ and $\beta$ now include additional terms:
\begin{equation}
\begin{array}{rcl}
\displaystyle  \int_0^\alpha \frac1{R(y)}\, dy
- \frac1{2\pi i} \int_\alpha^\beta \frac{\varphi'(y)}{R_+(y)}\, dy &=& 0 \\\\\displaystyle
  \int_0^\alpha \frac{y}{R(y)}\, dy
- \frac1{2\pi i} \int_\alpha^\beta \frac{y\varphi'(y)}{R_+(y)}\, dy &=& c\,.
\end{array}
\end{equation}
Evaluating these integrals, these equations are equivalent to
({\em cf.} \eqref{eq;SDeq4} and \eqref{eq;SDeq5})
\begin{eqnarray}
  (A+B+2)s+WX &=& B+1-YZ \\
\label{eq;SDeqanz2}
  Y-Z &=& 2c+A+B -(W+X)
\end{eqnarray}
where $W,X,Y,Z$ are exactly as defined in \eqref{eq;Wdef}--\eqref{eq;Zdef}.
Similar reasoning then yields
\begin{eqnarray}
  2W &=& K + \frac{A^2+2As}{K} \\
  2X &=& K - \frac{A^2+2As}{K} \\
  2Y &=& 2c+A+B-K + \frac{B^2+2B(1-s)}{2c+A+B-K} \\
  2Z &=& -(2c+A+B-K) + \frac{B^2+2B(1-s)}{2c+A+B-K}.
\end{eqnarray}
Substituting these formulae into \eqref{eq;SDeqanz2}, we obtain {\em
exactly} the same equation as in the previous case, namely \eqref{eq;theSeq},
and we also find that
we must take the root $s=s_1$.  From this, we get
\begin{equation}
  X= -X_0, \qquad Z=Z_0\,.
\end{equation}
Thus $X, Z>0$ imply the following conditions on $c$:
\begin{equation}
  c_A< c<c_B\,.
\end{equation}
Hence it is necessary in the saturated-band-void case that $A<B$ (see
\eqref{eq:cAcB}).  Under these conditions, one finds that the
solutions $\alpha$ and $\beta$ are again given by the exactly the same formulae
as in the previous case --- only the conditions on $c$, $A$, and $B$
are different.  The candidate equilibrium measure in the band $(\alpha,\beta)$
and the corresponding Lagrange multiplier may now be found as before
by evaluating the integrals in the explicit formula for $F(z)$ and
taking boundary values on $(\alpha,\beta)$.

\subsection{The void-band-saturated configuration.}
The appropriate scalar Riemann-Hilbert problem for the Cauchy
transform of the candidate density is the following: $F(z)$ is
analytic in $\mathbb{C}\setminus [\alpha,1]$ and satisfies the jump
conditions
\begin{equation}
  F_+(x)-F_-(x) = -\frac{2\pi i}{c}
\end{equation}
for $\beta<x<1$, and
\begin{equation}
  F_+(x)+F_-(x) = \frac1{c} \varphi(x)
\end{equation}
for $\alpha<x<\beta$,
and as $z\to\infty$,
\begin{equation}
  F(z) = \frac1{z} + O\left(\frac{1}{z^2}\right)\,.
\end{equation}
The solution is
\begin{equation}
  F(z) = \frac{R(z)}{2\pi ic} \left(
\int_\alpha^\beta \frac{\varphi'(y)}{R_+(y)(y-z)}\, dy-\int_\beta^1 \frac{2\pi i}{R(y)(y-z)}\, dy
 \right)\,.
\end{equation}
By taking moments of $F(z)$ for large $z$ and analyzing the resulting
equations we find again the same quadratic equation for $s$, but this time we
get
\begin{equation}
  X=X_0, \qquad Z=-Z_0\,.
\end{equation}
These imply the following conditions on $c$:
\begin{equation}
  c_B<c<c_A\,,
\end{equation}
and therefore this configuration is only possible if $B<A$. Again one
then finds that $\alpha$ and $\beta$ are given by the same formulae as before, and
by evaluating $F(z)$ one can calculate the candidate density for $\alpha<x<\beta$
and the Lagrange multiplier $\ell_c$.

\subsection{The saturated-band-saturated configuration.}
The scalar Riemann-Hilbert problem for the Cauchy transform of
$\psi(x)$ is:  $F(z)$ is analytic in
$\mathbb{C}\setminus [0,1]$ and satisfies the jump conditions
\begin{equation}
  F_+(x)-F_-(x) = -\frac{2\pi i}{c}
\end{equation}
for $0<x<\alpha$ and $\beta<x<1$, and
\begin{equation}
  F_+(x)+F_-(x) = \frac1{c} \varphi(x)
\end{equation}
for $\alpha<x<\beta$,
and as $z\to\infty$,
\begin{equation}
  F(z) = \frac1{z} + O\left(\frac{1}{z^2}\right)\,.
\end{equation}
The solution is
\begin{equation}
  F(z) = \frac{R(z)}{2\pi ic} \left(
-\int_0^\alpha \frac{2\pi i}{R(y)(y-z)}\,dy + \int_\alpha^\beta \frac{\varphi'(y)}{R_+(y)(y-z)}\, dy-
\int_\beta^1 \frac{2\pi i}{R(y)(y-z)}\, dy
 \right)\,.
\end{equation}
By similar analysis, we arrive again at the same quadratic
equation for $s$, and find
\begin{equation}
  X=-X_0, \qquad Z=-Z_0\,.
\end{equation}
These imply the following conditions on $c$:
\begin{equation}
  \max(c_A, c_B) <c <1\,,
\end{equation}
and again the endpoints $\alpha$ and $\beta$ have the same expressions as
before.  Evaluating the integrals in $F(z)$ and taking boundary values
then gives the candidate density in $(\alpha,\beta)$, and the Lagrange multiplier
may then be found by direct integration.

\clearpage
\section{List of Important Symbols}
%\subsection*{Symbols introduced in \S~1.}
\tablehead{Symbol & Meaning & Page reference\\
\hline \\
}
\begin{supertabular}{lll}
$N$ & Number of nodes & \pageref{symbol:number}\\\\
$X_N$ & Set of nodes & \pageref{symbol:nodeset}\\\\
$x_{N,n}$ & Node & \pageref{symbol:node}\\\\
$w_{N,n}$ & Weight at the node $x_{N,n}$ & \pageref{symbol:weight}\\\\
$w(x)$ & Weight function defined for $x\in X_N$ &
\pageref{symbol:weightfunc}\\\\
$p_{N,k}(z)$ & Discrete orthonormal polynomial & \pageref{symbol:DOP}\\\\
$c_{N,k}^{(m)}$ & Coefficient of $z^m$ in $p_{N,k}(z)$ & \pageref{symbol:coeffs}\\\\
$\gamma_{N,k}$ & Leading coefficient of $p_{N,k}(z)$ & \pageref{symbol:leadingcoeff}\\\\
$\pi_{N,k}(z)$ & Monic discrete orthogonal polynomial & \pageref{symbol:MDOP}\\\\
$\rho^0(x)$ & Node density function & \pageref{symbol:rho0}\\\\
$[a,b]$ & Interval containing nodes & \pageref{symbol:nodeinterval}\\\\
$V_N(x)$ & Exponent of weights & \pageref{symbol:VN}\\\\
$V(x)$ & Fixed component of $V_N(x)$ & \pageref{symbol:V}\\\\
$\eta(x)$ & Correction to $NV(x)$ & \pageref{symbol:eta}\\\\
$k$ & Degree of polynomial, number of particles & \pageref{symbol:k}\\\\
$c$ & Asymptotic ratio of $k/N$ & \pageref{symbol:c}\\\\
$\kappa$ & Correction to $Nc$ & \pageref{symbol:kappa}\\\\
$\mat{P}(z;N,k)$ & Solution of Interpolation Problem~\ref{rhp:DOP} &
\pageref{symbol:P}\\\\
$a_{N,k}$ & Diagonal recurrence coefficients & \pageref{symbol:aNk}\\\\
$b_{N,k}$ & Off-diagonal recurrence coefficients & \pageref{symbol:bNk}\\\\
%\end{tabular}
%\subsection*{Symbols introduced in \S~2.}
%\begin{tabular}{lll}
%Symbol & Meaning & Page of definition or first appearance\\
%\hline \\
$\del $ & Subset of node indices where triangularity is reversed & 
\pageref{symbol:del}\\\\
$\mathbb{Z}_N$ & $\{0,1,2,\dots,N-1\}$ & \pageref{symbol:ZN}\\\\
$\#\del$ & Number of elements in $\del$ & \pageref{symbol:Numdel}\\\\
$\mat{Q}(z;N,k)$ & $\mat{P}(z;N,k)$ with residue triangularity modified &
\pageref{symbol:Q}\\\\
$\sigma_3$ & Pauli matrix & \pageref{symbol:sigma3}\\\\
$\nab$ & Complementary set to $\del$ in $\mathbb{Z}_N$ & \pageref{symbol:nab}
\\\\
$\overline{\mat{P}}(z;N,\overline{k})$ & Dual of $\mat{P}(z;N,k)$ & 
\pageref{symbol:Pbar}\\\\
$\bar{k}$ & $N-k$, number of holes & \pageref{symbol:bark}\\\\
$\sigma_1$ & Pauli matrix & \pageref{symbol:sigma1}\\\\
$\overline{w}_{N,n}$ & Dual weight at the node $x_{N,n}$ & 
\pageref{symbol:dualweights}\\\\
$\overline{\pi}_{N,\bar{k}}(z)$ & Dual of $\pi_{N,k}(z)$ &
\pageref{symbol:pidual}\\\\
$\overline{\gamma}_{N,\bar{k}-1}$ & Dual of $\gamma_{N,k}$ & 
\pageref{symbol:gammadual}\\\\
$\varphi(x)$ & External field & \pageref{symbol:varphi}\\\\
$E_c[\mu]$ & Energy functional & \pageref{symbol:Ec}\\\\
$\mu_{\rm min}^c$ & Equilibrium measure & \pageref{symbol:minimizer}\\\\
$F_c[\mu]$ & Modified energy functional & \pageref{symbol:Fc}\\\\
$\ell_c$ & Lagrange multiplier (Robin constant) & \pageref{symbol:ellc}\\\\
$\underline{\mathcal F}$ & Set where lower constraint holds & \pageref{symbol:underF}\\\\
$\overline{\mathcal F}$ & Set where upper constraint holds &
\pageref{symbol:overF}\\\\
$G$ & Genus of $S$ & \pageref{symbol:G}\\\\
$\alpha_0,\dots,\alpha_G$ & Left endpoints of bands & \pageref{symbol:alphas}\\\\
$\beta_0,\dots,\beta_G$ & Right endpoints of bands & \pageref{symbol:betas}\\\\
$I_0,\dots,I_G$ & Bands & \pageref{symbol:bands}\\\\
$\Gamma_1,\dots,\Gamma_G$ & Interior gaps (voids and saturated regions) & \pageref{symbol:gaps}\\\\
$\displaystyle \frac{\delta E_c}{\delta\mu}(x)$ &
Variational derivative of $E_c$ &\pageref{symbol:varderiv}\\\\
$L_c(z)$ & Complex logarithmic potential of $\mu_{\rm min}^c$ &
\pageref{symbol:Lc(z)}\\\\
$\overline{L}_c^\Gamma(z)$ & Continuation from $\Gamma$ of logarithmic potential of $\mu_{\rm min}^c$ & \pageref{symbol:LcbarGamma}\\\\
$\overline{L}_c^I(z)$ & Continuation from $I$ of logarithmic potential of $\mu_{\rm min}^c$ &\pageref{symbol:LcbarI}\\\\
$\xi_\Gamma(x)$ & Analytic function defined in gap $\Gamma$ &
\pageref{symbol:xiGamma}\\\\
$\psi_I(x)$ and $\overline{\psi}_I(x)$ & Analytic functions defined in band $I$ & \pageref{symbol:psiI}\\\\
$\tau_\Gamma^{\nab,L}(z)$ & Conformal mapping near band/void edge $\alpha$ &
\pageref{symbol:taunabL}\\\\
$\tau_\Gamma^{\nab,R}(z)$ & Conformal mapping near band/void edge $\beta$ &
\pageref{symbol:taunabR}\\\\
$\tau_\Gamma^{\del,L}(z)$ & Conformal mapping near band/saturated region edge 
$\alpha$ &\pageref{symbol:taudelL}\\\\
$\tau_\Gamma^{\del,R}(z)$ & Conformal mapping near band/saturated region edge
$\beta$ &\pageref{symbol:taudelR}\\\\
$\theta_{\Gamma_1},\dots,\theta_{\Gamma_G}$ & Constants defined in interior gaps $\Gamma_1,\dots,\Gamma_G$ &\pageref{symbol:thetaconst}\\\\
$\theta_{(a,\alpha_0)}$ and $\theta_{(\beta_G,b)}$ &
Constants defined in exterior gaps $(a,\alpha_0)$ and $(\beta_G,b)$&\pageref{symbol:thetaGammaleftright}\\\\
$\overline{V}_N(x)$ & Dual of $V_N(x)$ & \pageref{symbol:VbarN}\\\\
$\bar{\mu}_{\rm min}^{1-c}$ & Dual of $\mu_{\rm min}^c$ & 
\pageref{symbol:mubar}\\\\
$R(z)$ & Branch of square root of $(z-\alpha_0)\cdots(z-\beta_G)$ & \pageref{symbol:sqrt}\\\\
$h(z)$ & Solution of Riemann-Hilbert Problem~\ref{rhp:h} & \pageref{symbol:h(z)}\\\\
$\gamma $ & Correction to $N\ell_c$ & \pageref{symbol:gamma}\\\\
$c_j^{(0)}$ & $\kappa$-independent part of $h_+(z)-h_-(z)$ for $z\in\Gamma_j$ &
\pageref{symbol:cjzero}\\\\
$\omega_j$ & Coefficient of $\kappa$ in $h_+(z)-h_-(z)$ for $z\in\Gamma_j$ &
\pageref{symbol:omegaj}\\\\
$\mat{r}$ & Phase vector with components $N\theta_{\Gamma_j}-c_j^{(0)}$ &
\pageref{symbol:vectorr}\\\\
$\mat{\Omega}$ & Frequency vector with components $\omega_j$ &
\pageref{symbol:vectorOmega}\\\\
$y(z)$ & Branch of square root of $(z-\alpha_0)\cdots (z-\beta_G)$ &
\pageref{symbol:y(z)}\\\\
$m_1(z)\,dz,\dots,m_G(z)\,dz$ & Branches of holomorphic differentials $m_1^S(P),\dots,m_G^S(P)$ &
\pageref{symbol:vecm(z)}\\\\
$\mat{A}$ & Matrix with columns $\mat{a}^{(1)}$,\dots,$\mat{a}^{(G)}$ &
\pageref{symbol:Amatrix}\\\\
$\mat{B}$ & Riemann matrix with columns $\mat{b}^{(1)}$,\dots,$\mat{b}^{(G)}$ &
\pageref{symbol:Bmatrix}\\\\
$\mat{k}$ & Vector of Riemann constants & \pageref{symbol:Riemannk}\\\\
$\Theta(\mat{w})$ & Riemann theta function & \pageref{symbol:thetafunction}
\\\\
$\mat{w}(z)$ & Branch of Abel-Jacobi mapping & \pageref{symbol:Abelmap}\\\\
$\mat{w}_\pm(\infty)$ & Limiting values of $\mat{w}(z)$ as $z\rightarrow\infty$ &\pageref{symbol:wpminfty}\\\\
$u(z)$ and $v(z)$ & Factors in solution of Riemann-Hilbert Problem~\ref{rhp:theta} &
\pageref{symbol:u(z)}\\\\
$\mat{q}_u$  and $\mat{q}_v$ & Vectors in the Jacobian giving zeros of $u(z)$ and $v(z)$&
\pageref{symbol:qu}\\\\
$W(z)$ & Alternate notation for $\dot{X}_{11}(z)e^{\kappa g(z)}$ &
\pageref{symbol:W(z)}\\\\
$Z(z)$ & Alternate notation for $\dot{X}_{12}(z)e^{-\kappa g(z)}$ &
\pageref{symbol:Z(z)}\\\\
$H_\Gamma^\pm(z)$ & Factors in first row of 
$\mat{H}_\Gamma^{\nab,L}(z)$, $\mat{H}_\Gamma^{\nab,R}(z)$, 
$\mat{H}_\Gamma^{\del,L}(z)$, and $\mat{H}_\Gamma^{\del,R}(z)$ & 
\pageref{symbol:HGammapm}\\\\
$K_J^\delta$ & Compact complex neighborhood of a closed interval $J$ &
\pageref{symbol:KJdelta}\\\\
$\theta^0(z)$ & Phase variable related to $\rho^0(z)$ &
\pageref{symbol:theta0}\\\\
$w_{N,n}^{\rm Kraw}(p,q)$ & Krawtchouk weights & \pageref{symbol:wKraw}\\\\
$V_N^{\rm Kraw}(x;l)$ & Exponent of Krawtchouk weights & 
\pageref{symbol:VKraw}\\\\
$w_{N,n}(b,c,d)$ & Weight degenerating to $w^{\rm Hahn}_{N,n}(P,Q)$ and
$w^{\rm Assoc}_{N,n}(P,Q)$ & \pageref{symbol:Hahnraw}\\\\
$w^{\rm Hahn}_{N,n}(P,Q)$ & Hahn weights & \pageref{symbol:wHahn}\\\\
$w^{\rm Assoc}_{N,n}(P,Q)$ & Associated Hahn weights &
\pageref{symbol:wAssoc}\\\\
$V_N^{\rm Hahn}(x;P,Q)$ &Exponent of $w_{N,n}^{\rm Hahn}(P,Q)$ &
\pageref{symbol:VNHahn}\\\\
$V^{\rm Hahn}(x;A,B)$ & Fixed component of $V_N^{\rm Hahn}(x;NA+1,NB+1)$ &
\pageref{symbol:VHahn}\\\\
$\eta^{\rm Hahn}(x;P,Q)$ & Correction to $NV^{\rm Hahn}(x;A,B)$ &
\pageref{symbol:etaHahn}\\\\
$V_N^{\rm Assoc}(x;P,Q)$ & Exponent of $w_{N,n}^{\rm Assoc}(P,Q)$ &
\pageref{symbol:VNAssoc}\\\\
$V^{\rm Assoc}(x;A,B)$ & Fixed component of $V_N^{\rm Assoc}(x;NA+1,NB+1)$ &
\pageref{symbol:VAssoc}\\\\
$c_A$ and $c_B$ & Critical values of $c$ for the Hahn equilibrium measure &
\pageref{symbol:cAcB}\\\\
$p^{(N,k)}(x_1,\dots,x_k)$ & Joint probability distribution of $k$ particles &
\pageref{symbol:jointprobparticles}\\\\
$\mathbb{P}({\rm event})$ & Probability of an event &
\pageref{symbol:Prob}\\\\
$Z_{N,k}$ & Normalization constant for 
$p^{(N,k)}(x_1,\dots,x_k)$ & \pageref{symbol:partitionfunction}\\\\
$R_m^{(N,k)}(x_1,\dots,x_m)$ & $m$-point correlation function of
$k$-particle ensemble&
\pageref{symbol:RmNk}\\\\
%$R_1^{(N,k)}(x)$ & One point function (density of states) &
%\pageref{symbol:onepoint}\\\\
$\mathbb{E}(X)$ & Expected value of a random variable $X$ & 
\pageref{symbol:Expected}\\\\
$K_{N,k}(x,y)$ & Reproducing (Christoffel-Darboux) kernel &
\pageref{symbol:reproducing}\\\\
$A_m^{(N,k)}(B)$ & Local particle occupation probability & 
\pageref{symbol:AmNkB}\\\\
$\overline{p}^{(N,\bar{k})}(y_1,\dots,y_{\bar{k}})$ &
Joint probability distribution of $\bar{k}$ holes &
\pageref{symbol:jointprobholes}\\\\
$\overline{Z}_{N,\bar{k}}$ & Normalization constant 
for $\overline{p}^{(N,\bar{k})}(y_1,\dots,y_{\bar{k}})$ &
\pageref{symbol:Zbar}\\\\
$\mathfrak{a}$, $\mathfrak{b}$, and $\mathfrak{c}$ &
Dimensions of the $\mathfrak{abc}$-hexagon &\pageref{symbol:Hexabc}\\\\
$P_1$, $P_2$, $P_3$, $P_4$, $P_5$, and $P_6$ &
Vertices of the $\mathfrak{abc}$-hexagon &\pageref{symbol:vertices}\\\\
$\mathcal{L}$ & Hexagonal lattice within the $\mathfrak{abc}$-hexagon &
\pageref{symbol:hexlattice}\\\\
$\mathfrak{A}$, $\mathfrak{B}$, and $\mathfrak{C}$ & Rescaled $\mathfrak{a}$, $\mathfrak{b}$,
and $\mathfrak{c}$ &\pageref{symbol:alphabetagamma}\\\\
$\mathcal{L}_m$ & $m^{\rm th}$ vertical sublattice of $\mathcal{L}$ &
\pageref{symbol:hexlatticem}\\\\
$N(\mathfrak{a},\mathfrak{b},\mathfrak{c},m)$
& Number of points in $\mathcal{L}_m$ &\pageref{symbol:cardLm}\\\\
$\mathfrak{a}_m$ and $\mathfrak{b}_m$ & $|m-\mathfrak{a}|$ and 
$|m-\mathfrak{b}|$ & \pageref{symbol:ambm}\\\\
$Q_m$ & Lowest lattice point in $\mathcal{L}_m$ & \pageref{symbol:Qm}\\\\
$L_m$ & Number of holes in $\mathcal{L}_m$ & \pageref{symbol:holesinLm}\\\\
$\tilde{P}_m(x_1,\dots,x_\mathfrak{c})$ &
Probability of finding particles at $x_1,\dots,x_\mathfrak{c}$ in 
$\mathcal{L}_m$ & \pageref{symbol:particleprob}\\\\
$P_m(\xi_1,\dots,\xi_{L_m})$ & 
Probability of finding holes at $\xi_1,\dots,\xi_{L_m}$ in 
$\mathcal{L}_m$ & \pageref{symbol:holeprob}\\\\
$\tau$ & Rescaled location of $\mathcal{L}_m$ in the $\mathfrak{abc}$-hexagon &
\pageref{symbol:tauparameter}\\\\
$S(\xi,\eta)$ & Discrete sine kernel & \pageref{symbol:DSK}\\\\
$\mathcal{S}_{ij}(x)$ & Node index form of $S(\xi,\eta)$ &
\pageref{symbol:DSKindex}\\\\
$E_{\rm int}([A,B];x,H,N)$ &
Expected number of particles near $x$ & \pageref{symbol:Eint}\\\\
$M_{\rm int}([A,B];x,H,N)$ & Number of nodes near $x$ &
\pageref{symbol:Mint}\\\\
$A(\xi,\eta)$ & Airy kernel & \pageref{symbol:AiryK}\\\\
$x_{\rm min}$ and $x_{\rm max}$ & Nodes occupied by leftmost and
rightmost particles & \pageref{symbol:xminmax}\\\\
$\mathcal{A}|_{[s,\infty)}$ & Operator acting with kernel $A(\xi,\eta)$ on
$L^2[s,\infty)$ & \pageref{symbol:calA}\\\\
$h_{\rm min}$ and $h_{\rm max}$ & Nodes occupied by leftmost and
rightmost holes & \pageref{symbol:hminmax}\\\\
$Y_\infty=\{y_1,\dots,y_M\}$ & Limiting transition points &
\pageref{symbol:yk}\\\\
$Y_N=\{y_{1,N},\dots,y_{M,N}\}$ & Transition points &
\pageref{symbol:ykN}\\\\
$\Sigma_0^\nab$ and $\Sigma_0^\del$ & Complementary systems of subintervals
of $(a,b)$ & \pageref{symbol:Sigmanaughts}\\\\
$d_N$ & $\#\del/N$ & \pageref{symbol:dN}\\\\
$\epsilon$ & Contour parameter & \pageref{symbol:epsilon}\\\\
$\Sigma$ & Contour of discontinuity of $\mat{R}(z)$ &\pageref{symbol:Sigma}\\\\
$\Omega_\pm^\nab$ and $\Omega_\pm^\del$ & 
Compact regions of $\mathbb{C}\setminus\Sigma$ &\pageref{symbol:Omegapmnabdel}
\\\\
$\mat{R}(z)$ & Matrix unknown obtained from $\mat{Q}(z;N,k)$ &
\pageref{symbol:matrixR}\\\\
$\rho(x)$ & Density for $g(z)$ & \pageref{symbol:rho}\\\\
$g(z)$ & Complex logarithmic potential of $\rho(x)$ & \pageref{symbol:gdef}
\\\\
$\mat{S}(z)$ & Matrix unknown obtained from $\mat{R}(z)$ &
\pageref{symbol:matrixS}\\\\
$\theta(z)$ & Phase variable related to $\rho(x)$ & \pageref{symbol:thetaofz}
\\\\
$\phi(z)$ & Correction to $N\theta(z)$ & \pageref{symbol:phiofz}\\\\
$T_\nab(z)$ & Analytic function measuring discreteness in $\Sigma_0^\nab$ &
\pageref{symbol:Tnab}\\\\
$T_\del(z)$ & Analytic function measuring discreteness in $\Sigma_0^\del$ &
\pageref{symbol:Tdel}\\\\
$\phi_\Gamma$ & Constant value of $\phi(z)$ in gap $\Gamma$ & 
\pageref{symbol:phiGamma}\\\\
$\mat{L}_\pm(z)$ & Lower-triangular factors in jump for $\mat{S}(z)$ in bands
&\pageref{symbol:Lpm}\\\\
$\mat{J}(z)$ & Off-diagonal factor in jump for $\mat{S}(z)$ in bands &
\pageref{symbol:Jmatrix}\\\\
$\mat{U}_\pm(z)$ & Upper-triangular factors in jump for $\mat{S}(z)$ in bands
&\pageref{symbol:Upm}\\\\
$Y(z)$ & Scalar function related to $T_\nab(z)$ and $T_\del(z)$ &
\pageref{symbol:Y(z)}\\\\
%$D(z,x)$ & Difference quotient of $\theta^0(x)$ &
%\pageref{symbol:differencequotient}\\\\
%$m(x)$ & Mass integral associated with the density $\rho^0(x)$ &
%\pageref{symbol:m(x)}\\\\
%$s_{N,n}$ & Image of $x_{N,n}$ with equal spacing & %\pageref{symbol:straightened}\\\\
%$\tilde{T}_\nab(z)$ & $T_\nab(z)$ with equally spaced nodes &
%\pageref{symbol:TildeTnab}\\\\
%$\tilde{T}_\del(z)$ & $T_\del(z)$ with equally spaced nodes &
%\pageref{symbol:TildeTdel}\\\\
%$\tilde{T}_\nab^J(z)$ & Dominant contribution to $\tilde{T}_\nab(z)$ as %$N\rightarrow\infty$ &
%\pageref{symbol:TildeTnabJ}\\\\
$\theta_I^\nab(z)$ and $\theta_I^\del(z)$ & Analytic continuation of 
$\theta(z)$ from $I\cap\Sigma_0^\nab$ and $I\cap\Sigma_0^\del$  & \pageref{symbol:thetanabI}\\\\
$\mat{X}(z)$ & Solution of Riemann-Hilbert Problem~\ref{rhp:X} &
\pageref{symbol:matrixX}\\\\
$\Sigma_{\rm SD}$ & Contour of discontinuity of $\mat{X}(z)$ &
\pageref{symbol:SigmaSD}\\\\
$\mat{D}(z)$ & Matrix factor relating $\mat{X}(z)$ and $\mat{P}(z;N,k)$ &
\pageref{symbol:matrixD}\\\\
$\Sigma_{0\pm}^\nab$ and $\Sigma_{0\pm}^\del$ &
Vertical segments of $\Sigma_{\rm SD}$ connected to band endpoints &
\pageref{symbol:Sigma0pmnabdel}\\\\
$\Sigma_{I\pm}$ & Horizontal segments of $\Sigma_{\rm SD}$ parallel to a 
band $I$ &\pageref{symbol:SigmaIpm}\\\\
$\Sigma_{\Gamma\pm}$ & Horizontal segments of $\Sigma_{\rm SD}$ parallel to
a gap $\Gamma$ & \pageref{symbol:SigmaGammapm}\\\\
$\dot{\mat{X}}(z)$ & Solution of Riemann-Hilbert Problem~\ref{rhp:theta} &
\pageref{symbol:matrixdotX}\\\\
$\Sigma_{\rm model}$ & Contour of discontinuity of $\dot{\mat{X}}(z)$ &
\pageref{symbol:Sigmamodel}\\\\
$\Psi(z)$ & $\delta E_c/\delta\mu(z)-(d_N-c)(g_+(z)+g_-(z))$ &
\pageref{symbol:Psi}\\\\
$h$ & Additional contour parameter (with $\epsilon$) & 
\pageref{symbol:hparameter}\\\\
$D_\Gamma^{\nab,L}$ & Disc  centered at
band/void edge $z=\alpha$ & \pageref{symbol:DnabL}\\\\
$D_{\Gamma,I}^{\nab,L}$, $D_{\Gamma,II}^{\nab,L}$, $D_{\Gamma,III}^{\nab,L}$,
and $D_{\Gamma,IV}^{\nab,L}$ & Quadrants of $D_{\Gamma}^{\nab,L}$ &
\pageref{symbol:DnabLquadrants}\\\\
$\mat{Z}_\Gamma^{\nab,L}$ & Matrix proportional to $\mat{X}(z)$ in 
$D_\Gamma^{\nab,L}$ & \pageref{symbol:ZnabL} \\\\
$\dot{\mat{Z}}_\Gamma^{\nab,L}$ & Matrix proportional to $\dot{\mat{X}}(z)$
in $D_\Gamma^{\nab,L}$ & \pageref{symbol:dotZnabL}\\\\
$\mat{H}_\Gamma^{\nab,L}(z)$ & Holomorphic prefactor in
$\dot{\mat{Z}}_\Gamma^{\nab,L}(z)$ &\pageref{symbol:HGammanabL}\\\\
$\hat{\mat{Z}}^{\nab,L}(\zeta)$ & Explicit model for 
$\mat{Z}_\Gamma^{\nab,L}(z)$ & \pageref{symbol:hatZnabL}\\\\
$\hat{\mat{X}}_\Gamma^{\nab,L}(z)$ & Local parametrix for $\mat{X}(z)$
in $D_\Gamma^{\nab,L}$ & \pageref{symbol:XhatnabL}\\\\
$D_\Gamma^{\nab,R}$ & Disc  centered at
 band/void edge $z=\beta$ & \pageref{symbol:DnabR}\\\\
$D_{\Gamma,I}^{\nab,R}$, $D_{\Gamma,II}^{\nab,R}$, $D_{\Gamma,III}^{\nab,R}$,
and $D_{\Gamma,IV}^{\nab,R}$ & Quadrants of $D_{\Gamma}^{\nab,R}$ &
\pageref{symbol:DnabRquadrants}\\\\
$\mat{Z}_\Gamma^{\nab,R}$ & Matrix proportional to $\mat{X}(z)$ in 
$D_\Gamma^{\nab,R}$ & \pageref{symbol:ZnabR} \\\\
$\dot{\mat{Z}}_\Gamma^{\nab,R}$ & Matrix proportional to $\dot{\mat{X}}(z)$
in $D_\Gamma^{\nab,R}$ & \pageref{symbol:dotZnabR}\\\\
$\mat{H}_\Gamma^{\nab,R}(z)$ & Holomorphic prefactor in
$\dot{\mat{Z}}_\Gamma^{\nab,R}(z)$ &\pageref{symbol:HGammanabR}\\\\
$\hat{\mat{Z}}^{\nab,R}(\zeta)$ & Explicit model for 
$\mat{Z}_\Gamma^{\nab,R}(z)$ & \pageref{symbol:hatZnabR}\\\\
$\hat{\mat{X}}_\Gamma^{\nab,R}(z)$ & Local parametrix for $\mat{X}(z)$
in $D_\Gamma^{\nab,R}$ & \pageref{symbol:XhatnabR}\\\\
$D_\Gamma^{\del,L}$ & Disc centered at
band/saturated region edge $z=\alpha$ & \pageref{symbol:DdelL}\\\\
$D_{\Gamma,I}^{\del,L}$, $D_{\Gamma,II}^{\del,L}$, $D_{\Gamma,III}^{\del,L}$,
and $D_{\Gamma,IV}^{\del,L}$ & Quadrants of $D_{\Gamma}^{\del,L}$ &
\pageref{symbol:DdelLquadrants}\\\\
$\mat{Z}_\Gamma^{\del,L}$ & Matrix proportional to $\mat{X}(z)$ in 
$D_\Gamma^{\del,L}$ & \pageref{symbol:ZdelL} \\\\
$\dot{\mat{Z}}_\Gamma^{\del,L}$ & Matrix proportional to $\dot{\mat{X}}(z)$
in $D_\Gamma^{\del,L}$ & \pageref{symbol:dotZdelL}\\\\
$\mat{H}_\Gamma^{\del,L}(z)$ & Holomorphic prefactor in
$\dot{\mat{Z}}_\Gamma^{\del,L}(z)$ &\pageref{symbol:HGammadelL}\\\\
$\hat{\mat{Z}}^{\del,L}(\zeta)$ & Explicit model for 
$\mat{Z}_\Gamma^{\del,L}(z)$ & \pageref{symbol:hatZdelL}\\\\
$\hat{\mat{X}}_\Gamma^{\del,L}(z)$ & Local parametrix for $\mat{X}(z)$
in $D_\Gamma^{\del,L}$ & \pageref{symbol:XhatdelL}\\\\
$D_\Gamma^{\del,R}$ & Disc centered at
band/saturated region edge $z=\beta$ & \pageref{symbol:DdelR}\\\\
$D_{\Gamma,I}^{\del,R}$, $D_{\Gamma,II}^{\del,R}$, $D_{\Gamma,III}^{\del,R}$,
and $D_{\Gamma,IV}^{\del,R}$ & Quadrants of $D_{\Gamma}^{\del,R}$ &
\pageref{symbol:DdelRquadrants}\\\\
$\mat{Z}_\Gamma^{\del,R}$ & Matrix proportional to $\mat{X}(z)$ in 
$D_\Gamma^{\del,R}$ & \pageref{symbol:ZdelR} \\\\
$\dot{\mat{Z}}_\Gamma^{\del,R}$ & Matrix proportional to $\dot{\mat{X}}(z)$
in $D_\Gamma^{\del,R}$ & \pageref{symbol:dotZdelR}\\\\
$\mat{H}_\Gamma^{\del,R}(z)$ & Holomorphic prefactor in
$\dot{\mat{Z}}_\Gamma^{\del,R}(z)$ &\pageref{symbol:HGammadelR}\\\\
$\hat{\mat{Z}}^{\del,R}(\zeta)$ & Explicit model for 
$\mat{Z}_\Gamma^{\del,R}(z)$ & \pageref{symbol:hatZdelR}\\\\
$\hat{\mat{X}}_\Gamma^{\del,R}(z)$ & Local parametrix for $\mat{X}(z)$
in $D_\Gamma^{\del,R}$ & \pageref{symbol:XhatdelR}\\\\
$\hat{\mat{X}}(z)$ & Parametrix (global) for $\mat{X}(z)$ &
\pageref{symbol:hatX}\\\\
$\mat{E}(z)$ & Error matrix $\mat{X}(z)\hat{\mat{X}}(z)^{-1}$ &
\pageref{symbol:matrixE}\\\\
$\Sigma_E$ & Contour of discontinuity of $\mat{E}(z)$ &
\pageref{symbol:SigmaE}\\\\
$L_\Gamma^\nab$ & Region of deformation below a void $\Gamma$ &
\pageref{symbol:LnabGamma}\\\\
$\mat{F}(z)$ & Solution of Riemann-Hilbert Problem~\ref{rhp:F}& 
\pageref{symbol:matrixF}\\\\
$L_\Gamma^\del$ & Region of deformation below a saturated region $\Gamma$ &
\pageref{symbol:LdelGamma}\\\\
$\Sigma_F$ & Contour of discontinuity of $\mat{F}(z)$ &
\pageref{symbol:SigmaF}\\\\
$\mat{v}_{\mat{F}}(z)$ & Jump matrix for $\mat{F}(z)$ on $\Sigma_F$ &
\pageref{symbol:vF}\\\\
$\overline{R}^{(N,\bar{k})}_m(x_1,\dots,x_m)$ & $m$-point correlation
function of $\bar{k}$-hole ensemble & \pageref{symbol:dualcorr}\\\\
$\overline{K}_{N,\bar{k}}(x,y)$ &
Dual of $K_{N,k}(x,y)$ & \pageref{symbol:dualK}\\\\
$\mat{B}(x)$ & Matrix factor in exact formula for $K_{N,k}(x,y)$ &
\pageref{symbol:matrixB}\\\\
$\mat{v}$ and $\mat{w}$ & Vector factors in exact formula for $K_{N,k}(x,y)$
&\pageref{symbol:vandw}\\\\
$\mat{a}$ and $\mat{b}$ & Vector factors in exact formula for $K_{N,k}(x,y)$
&\pageref{symbol:aandb}\\\\
$\mat{A}_\Gamma^{\nab,L}(x)$, $\mat{q}_\Gamma^{\nab,L}(x)$, and
$\mat{r}_\Gamma^{\nab,L}(x)$ & Factors in exact formula for $K_{N,k}(x,y)$
&\pageref{symbol:AnabL}\\\\
$\mat{A}_\Gamma^{\nab,R}(x)$, $\mat{q}_\Gamma^{\nab,R}(x)$, and
$\mat{r}_\Gamma^{\nab,R}(x)$ & Factors in exact formula for $K_{N,k}(x,y)$
&\pageref{symbol:AnabL}\\\\
$\mat{Y}(z)$ & Matrix constructed from $\dot{\mat{X}}(z)$ and $h(z)$ &
\pageref{symbol:matrixY}\\\\
$\mat{Y}^\sharp(z)$ & Matrix directly related to $\mat{Y}(z)$ &
\pageref{symbol:matrixYsharp}\\\\
$\Gamma_0$ & $(-\infty,\alpha_0)\cup (\beta_G,\infty)$ & 
\pageref{symbol:Gammazero}\\\\
$\Sigma_{\rm model}'$ & Contour of discontinuity of $\mat{Y}^\sharp(z)$ &
\pageref{symbol:Sigmamodelprime}\\\\
$S$ & Hyperelliptic Riemann surface & \pageref{symbol:RiemannSurface}\\\\
$z(P)$ & Hyperelliptic sheet projection function & \pageref{symbol:z(P)}\\\\
$y^S(P)$ & Analytic continuation of $y(z(P))$ to $S$ & \pageref{symbol:yS}\\\\
$a_1,\dots,a_G$ and $b_1,\dots,b_G$ & Homology basis on $S$ &
\pageref{symbol:homology}\\\\
$m_1^S(P),\dots,m_G^S(P)$ & Holomorphic differentials (unnormalized) on $S$ &
\pageref{symbol:holodiffs}\\\\
$\mat{w}^S(P)$ & Abel-Jacobi mapping on $S$ & \pageref{symbol:AbelmapS}\\\\
$f(z;\mat{q})$ & Shifted Riemann theta function & \pageref{symbol:shiftedtheta}
\\\\
$g^\pm(z;\mat{q})$ & Ratios of shifted Riemann theta functions &
\pageref{symbol:gzq}
\\
\end{supertabular}

\end{document}

%% file: defects.pstex_t
\begin{picture}(0,0)%
\epsfig{file=defects.pstex}%
\end{picture}%
\setlength{\unitlength}{1973sp}%
\begingroup\makeatletter\ifx\SetFigFont\undefined%
\gdef\SetFigFont#1#2#3#4#5{%
  \reset@font\fontsize{#1}{#2pt}%
  \fontfamily{#3}\fontseries{#4}\fontshape{#5}%
  \selectfont}%
\fi\endgroup%
\begin{picture}(12316,5516)(-157,-5694)
\end{picture}

%% file: defectmotion.pstex_t
\begin{picture}(0,0)%
\epsfig{file=defectmotion.pstex}%
\end{picture}%
\setlength{\unitlength}{1973sp}%
\begingroup\makeatletter\ifx\SetFigFont\undefined%
\gdef\SetFigFont#1#2#3#4#5{%
  \reset@font\fontsize{#1}{#2pt}%
  \fontfamily{#3}\fontseries{#4}\fontshape{#5}%
  \selectfont}%
\fi\endgroup%
\begin{picture}(8716,4622)(-157,-4472)
\end{picture}

%% file: tilingbase.pstex_t
\begin{picture}(0,0)%
\epsfig{file=tilingbase.pstex}%
\end{picture}%
\setlength{\unitlength}{1973sp}%
\begingroup\makeatletter\ifx\SetFigFont\undefined%
\gdef\SetFigFont#1#2#3#4#5{%
  \reset@font\fontsize{#1}{#2pt}%
  \fontfamily{#3}\fontseries{#4}\fontshape{#5}%
  \selectfont}%
\fi\endgroup%
\begin{picture}(10074,9879)(64,-9328)
\put(4726,-1036){\makebox(0,0)[lb]{\smash{\SetFigFont{12}{14.4}{\rmdefault}{\mddefault}{\updefault}{\color[rgb]{0,0,0}$P_5$}%
}}}
\put(7276,-2686){\makebox(0,0)[lb]{\smash{\SetFigFont{12}{14.4}{\rmdefault}{\mddefault}{\updefault}{\color[rgb]{0,0,0}$P_4$}%
}}}
\put(676,-3061){\makebox(0,0)[lb]{\smash{\SetFigFont{12}{14.4}{\rmdefault}{\mddefault}{\updefault}{\color[rgb]{0,0,0}$P_6$}%
}}}
\put(4276,-9211){\makebox(0,0)[lb]{\smash{\SetFigFont{12}{14.4}{\rmdefault}{\mddefault}{\updefault}{\color[rgb]{0,0,0}$P_2$}%
}}}
\put(2326,-8686){\makebox(0,0)[lb]{\smash{\SetFigFont{11}{13.2}{\rmdefault}{\mddefault}{\updefault}{\color[rgb]{0,0,0}$\mathfrak{b}$}%
}}}
\put(5776,-8461){\makebox(0,0)[lb]{\smash{\SetFigFont{11}{13.2}{\rmdefault}{\mddefault}{\updefault}{\color[rgb]{0,0,0}$\mathfrak{a}$}%
}}}
\put(7351,-4786){\makebox(0,0)[lb]{\smash{\SetFigFont{11}{13.2}{\rmdefault}{\mddefault}{\updefault}{\color[rgb]{0,0,0}$\mathfrak{c}$}%
}}}
\put(6001,-1561){\makebox(0,0)[lb]{\smash{\SetFigFont{11}{13.2}{\rmdefault}{\mddefault}{\updefault}{\color[rgb]{0,0,0}$\mathfrak{b}$}%
}}}
\put(2626,-1636){\makebox(0,0)[lb]{\smash{\SetFigFont{11}{13.2}{\rmdefault}{\mddefault}{\updefault}{\color[rgb]{0,0,0}$\mathfrak{a}$}%
}}}
\put(601,-5236){\makebox(0,0)[lb]{\smash{\SetFigFont{11}{13.2}{\rmdefault}{\mddefault}{\updefault}{\color[rgb]{0,0,0}$\mathfrak{c}$}%
}}}
\put(301,-8011){\makebox(0,0)[lb]{\smash{\SetFigFont{12}{14.4}{\rmdefault}{\mddefault}{\updefault}{\color[rgb]{0,0,0}$P_1$}%
}}}
\put(7351,-7261){\makebox(0,0)[lb]{\smash{\SetFigFont{12}{14.4}{\rmdefault}{\mddefault}{\updefault}{\color[rgb]{0,0,0}$P_3$}%
}}}
\end{picture}

%% file: tilingrhombi.pstex_t
\begin{picture}(0,0)%
\includegraphics{tilingrhombi.pstex}%
\end{picture}%
\setlength{\unitlength}{1973sp}%
\begingroup\makeatletter\ifx\SetFigFont\undefined%
\gdef\SetFigFont#1#2#3#4#5{%
  \reset@font\fontsize{#1}{#2pt}%
  \fontfamily{#3}\fontseries{#4}\fontshape{#5}%
  \selectfont}%
\fi\endgroup%
\begin{picture}(6333,1861)(826,-4889)
\put(826,-4711){\makebox(0,0)[lb]{\smash{\SetFigFont{12}{14.4}{\rmdefault}{\mddefault}{\updefault}{\color[rgb]{0,0,0}Type I}%
}}}
\put(3301,-4786){\makebox(0,0)[lb]{\smash{\SetFigFont{12}{14.4}{\rmdefault}{\mddefault}{\updefault}{\color[rgb]{0,0,0}Type II}%
}}}
\put(6076,-4711){\makebox(0,0)[lb]{\smash{\SetFigFont{12}{14.4}{\rmdefault}{\mddefault}{\updefault}{\color[rgb]{0,0,0}Type III}%
}}}
\end{picture}

%% file: tilingandLm.pstex_t
\begin{picture}(0,0)%
\includegraphics{tilingandLm.pstex}%
\end{picture}%
\setlength{\unitlength}{1973sp}%
\begingroup\makeatletter\ifx\SetFigFont\undefined%
\gdef\SetFigFont#1#2#3#4#5{%
  \reset@font\fontsize{#1}{#2pt}%
  \fontfamily{#3}\fontseries{#4}\fontshape{#5}%
  \selectfont}%
\fi\endgroup%
\begin{picture}(10449,9795)(139,-9244)
\put(2176,-8686){\makebox(0,0)[lb]{\smash{\SetFigFont{12}{14.4}{\rmdefault}{\mddefault}{\updefault}{\color[rgb]{0,0,0}$Q_m$}%
}}}
\end{picture}

%% file: minimizer.pstex_t
\begin{picture}(0,0)%
\epsfig{file=minimizer.pstex}%
\end{picture}%
\setlength{\unitlength}{2960sp}%
\begingroup\makeatletter\ifx\SetFigFont\undefined%
\gdef\SetFigFont#1#2#3#4#5{%
  \reset@font\fontsize{#1}{#2pt}%
  \fontfamily{#3}\fontseries{#4}\fontshape{#5}%
  \selectfont}%
\fi\endgroup%
\begin{picture}(8400,6000)(601,-5761)
\put(5326,-4786){\makebox(0,0)[lb]{\smash{\SetFigFont{9}{10.8}{\rmdefault}{\mddefault}{\updefault}{\color[rgb]{0,0,0}$\psi_{I_2}(x)$}%
}}}
\put(2476,-4786){\makebox(0,0)[lb]{\smash{\SetFigFont{9}{10.8}{\rmdefault}{\mddefault}{\updefault}{\color[rgb]{0,0,0}$\Sigma_0^\del$}%
}}}
\put(6976,-4786){\makebox(0,0)[lb]{\smash{\SetFigFont{9}{10.8}{\rmdefault}{\mddefault}{\updefault}{\color[rgb]{0,0,0}$\Sigma_0^\del$}%
}}}
\put(4726,-511){\makebox(0,0)[lb]{\smash{\SetFigFont{9}{10.8}{\rmdefault}{\mddefault}{\updefault}{\color[rgb]{0,0,0}$\Sigma_0^\nab$}%
}}}
\put(8176,-511){\makebox(0,0)[lb]{\smash{\SetFigFont{9}{10.8}{\rmdefault}{\mddefault}{\updefault}{\color[rgb]{0,0,0}$\Sigma_0^\nab$}%
}}}
\put(1426,-5386){\makebox(0,0)[lb]{\smash{\SetFigFont{9}{10.8}{\rmdefault}{\mddefault}{\updefault}{\color[rgb]{0,0,0}$x=a$}%
}}}
\put(8476,-5311){\makebox(0,0)[lb]{\smash{\SetFigFont{9}{10.8}{\rmdefault}{\mddefault}{\updefault}{\color[rgb]{0,0,0}$x=b$}%
}}}
\put(4651,-1186){\makebox(0,0)[lb]{\smash{\SetFigFont{9}{10.8}{\rmdefault}{\mddefault}{\updefault}{\color[rgb]{0,0,0}$\rho^0(x)/c$}%
}}}
\put(4051,-3736){\makebox(0,0)[lb]{\smash{\SetFigFont{9}{10.8}{\rmdefault}{\mddefault}{\updefault}{\color[rgb]{0,0,0}$\psi_{I_1}(x)$}%
}}}
\put(6601,-2086){\makebox(0,0)[lb]{\smash{\SetFigFont{9}{10.8}{\rmdefault}{\mddefault}{\updefault}{\color[rgb]{0,0,0}$\psi_{I_3}(x)$}%
}}}
\put(3976,-1786){\makebox(0,0)[lb]{\smash{\SetFigFont{9}{10.8}{\rmdefault}{\mddefault}{\updefault}{\color[rgb]{0,0,0}$x=y_1$}%
}}}
\put(6526,-4111){\makebox(0,0)[lb]{\smash{\SetFigFont{9}{10.8}{\rmdefault}{\mddefault}{\updefault}{\color[rgb]{0,0,0}$x=y_2$}%
}}}
\put(7801,-1336){\makebox(0,0)[lb]{\smash{\SetFigFont{9}{10.8}{\rmdefault}{\mddefault}{\updefault}{\color[rgb]{0,0,0}$x=y_3$}%
}}}
\put(2851,-2086){\makebox(0,0)[lb]{\smash{\SetFigFont{9}{10.8}{\rmdefault}{\mddefault}{\updefault}{\color[rgb]{0,0,0}$\psi_{I_0}(x)$}%
}}}
\put(7951,-4036){\makebox(0,0)[lb]{\smash{\SetFigFont{9}{10.8}{\rmdefault}{\mddefault}{\updefault}{\color[rgb]{0,0,0}$\psi_{I_4}(x)$}%
}}}
\end{picture}

%% file: Sigma.pstex_t
\begin{picture}(0,0)%
\epsfig{file=Sigma.pstex}%
\end{picture}%
\setlength{\unitlength}{1579sp}%
\begingroup\makeatletter\ifx\SetFigFont\undefined%
\gdef\SetFigFont#1#2#3#4#5{%
  \reset@font\fontsize{#1}{#2pt}%
  \fontfamily{#3}\fontseries{#4}\fontshape{#5}%
  \selectfont}%
\fi\endgroup%
\begin{picture}(17175,4039)(-1424,-7169)
\put(14401,-4411){\makebox(0,0)[lb]{\smash{\SetFigFont{5}{6.0}{\rmdefault}{\mddefault}{\updefault}{\color[rgb]{0,0,0}{\normalsize $\Omega_+^\nab$}}%
}}}
\put(2551,-4411){\makebox(0,0)[lb]{\smash{\SetFigFont{5}{6.0}{\rmdefault}{\mddefault}{\updefault}{\color[rgb]{0,0,0}{\normalsize $\Omega_+^\del$}}%
}}}
\put( 76,-5236){\makebox(0,0)[lb]{\smash{\SetFigFont{5}{6.0}{\rmdefault}{\mddefault}{\updefault}{\color[rgb]{0,0,0}{\normalsize $a$}}%
}}}
\put(2551,-5986){\makebox(0,0)[lb]{\smash{\SetFigFont{5}{6.0}{\rmdefault}{\mddefault}{\updefault}{\color[rgb]{0,0,0}{\normalsize $\Omega_-^\del$}}%
}}}
\put(15751,-5236){\makebox(0,0)[lb]{\smash{\SetFigFont{5}{6.0}{\rmdefault}{\mddefault}{\updefault}{\color[rgb]{0,0,0}{\normalsize $b$}}%
}}}
\put(7951,-5911){\makebox(0,0)[lb]{\smash{\SetFigFont{5}{6.0}{\rmdefault}{\mddefault}{\updefault}{\color[rgb]{0,0,0}{\normalsize $\Omega_-^\nab$}}%
}}}
\put(7951,-4486){\makebox(0,0)[lb]{\smash{\SetFigFont{5}{6.0}{\rmdefault}{\mddefault}{\updefault}{\color[rgb]{0,0,0}{\normalsize $\Omega_+^\nab$}}%
}}}
\put(11851,-4411){\makebox(0,0)[lb]{\smash{\SetFigFont{5}{6.0}{\rmdefault}{\mddefault}{\updefault}{\color[rgb]{0,0,0}{\normalsize $\Omega_+^\del$}}%
}}}
\put(11851,-5911){\makebox(0,0)[lb]{\smash{\SetFigFont{5}{6.0}{\rmdefault}{\mddefault}{\updefault}{\color[rgb]{0,0,0}{\normalsize $\Omega_-^\del$}}%
}}}
\put(14401,-5986){\makebox(0,0)[lb]{\smash{\SetFigFont{5}{6.0}{\rmdefault}{\mddefault}{\updefault}{\color[rgb]{0,0,0}{\normalsize $\Omega_-^\nab$}}%
}}}
\put(-224,-3286){\makebox(0,0)[lb]{\smash{\SetFigFont{5}{6.0}{\rmdefault}{\mddefault}{\updefault}{\color[rgb]{0,0,0}{\normalsize $\Re(z)=a$}}%
}}}
\put(4426,-3286){\makebox(0,0)[lb]{\smash{\SetFigFont{5}{6.0}{\rmdefault}{\mddefault}{\updefault}{\color[rgb]{0,0,0}{\normalsize $\Re(z)=y_{1,N}$}}%
}}}
\put(9751,-7111){\makebox(0,0)[lb]{\smash{\SetFigFont{5}{6.0}{\rmdefault}{\mddefault}{\updefault}{\color[rgb]{0,0,0}{\normalsize $\Re(z)=y_{2,N}$}}%
}}}
\put(12601,-3286){\makebox(0,0)[lb]{\smash{\SetFigFont{5}{6.0}{\rmdefault}{\mddefault}{\updefault}{\color[rgb]{0,0,0}{\normalsize $\Re(z)=y_{3,N}$}}%
}}}
\put(14626,-7111){\makebox(0,0)[lb]{\smash{\SetFigFont{5}{6.0}{\rmdefault}{\mddefault}{\updefault}{\color[rgb]{0,0,0}{\normalsize $\Re(z)=b$}}%
}}}
\put(-1424,-6736){\makebox(0,0)[lb]{\smash{\SetFigFont{5}{6.0}{\rmdefault}{\mddefault}{\updefault}{\color[rgb]{0,0,0}{\normalsize $\Im(z)=-\epsilon$}}%
}}}
\put(-1424,-3736){\makebox(0,0)[lb]{\smash{\SetFigFont{5}{6.0}{\rmdefault}{\mddefault}{\updefault}{\color[rgb]{0,0,0}{\normalsize $\Im(z)=\epsilon$}}%
}}}
\end{picture}

%% file: SigmaSD.pstex_t
\begin{picture}(0,0)%
\epsfig{file=SigmaSD.pstex}%
\end{picture}%
\setlength{\unitlength}{1579sp}%
\begingroup\makeatletter\ifx\SetFigFont\undefined%
\gdef\SetFigFont#1#2#3#4#5{%
  \reset@font\fontsize{#1}{#2pt}%
  \fontfamily{#3}\fontseries{#4}\fontshape{#5}%
  \selectfont}%
\fi\endgroup%
\begin{picture}(17250,4039)(-1499,-7169)
\put( 76,-5236){\makebox(0,0)[lb]{\smash{\SetFigFont{5}{6.0}{\rmdefault}{\mddefault}{\updefault}{\color[rgb]{0,0,0}{\normalsize $a$}}%
}}}
\put(15751,-5236){\makebox(0,0)[lb]{\smash{\SetFigFont{5}{6.0}{\rmdefault}{\mddefault}{\updefault}{\color[rgb]{0,0,0}{\normalsize $b$}}%
}}}
\put(14776,-7111){\makebox(0,0)[lb]{\smash{\SetFigFont{5}{6.0}{\rmdefault}{\mddefault}{\updefault}{\color[rgb]{0,0,0}{\normalsize $\Re(z)=b$}}%
}}}
\put(-299,-3286){\makebox(0,0)[lb]{\smash{\SetFigFont{5}{6.0}{\rmdefault}{\mddefault}{\updefault}{\color[rgb]{0,0,0}{\normalsize $\Re(z)=a$}}%
}}}
\put(4501,-3361){\makebox(0,0)[lb]{\smash{\SetFigFont{5}{6.0}{\rmdefault}{\mddefault}{\updefault}{\color[rgb]{0,0,0}{\normalsize $\Re(z)=y_{1,N}$}}%
}}}
\put(9751,-7111){\makebox(0,0)[lb]{\smash{\SetFigFont{5}{6.0}{\rmdefault}{\mddefault}{\updefault}{\color[rgb]{0,0,0}{\normalsize $\Re(z)=y_{2,N}$}}%
}}}
\put(12526,-3286){\makebox(0,0)[lb]{\smash{\SetFigFont{5}{6.0}{\rmdefault}{\mddefault}{\updefault}{\color[rgb]{0,0,0}{\normalsize $\Re(z)=y_{3,N}$}}%
}}}
\put(-1499,-6736){\makebox(0,0)[lb]{\smash{\SetFigFont{5}{6.0}{\rmdefault}{\mddefault}{\updefault}{\color[rgb]{0,0,0}{\normalsize $\Im(z)=-\epsilon$}}%
}}}
\put(-1499,-3736){\makebox(0,0)[lb]{\smash{\SetFigFont{5}{6.0}{\rmdefault}{\mddefault}{\updefault}{\color[rgb]{0,0,0}{\normalsize $\Im(z)=\epsilon$}}%
}}}
\end{picture}

%% file: SigmaSD_components.pstex_t
\begin{picture}(0,0)%
\epsfig{file=SigmaSD_components.pstex}%
\end{picture}%
\setlength{\unitlength}{1579sp}%
\begingroup\makeatletter\ifx\SetFigFont\undefined%
\gdef\SetFigFont#1#2#3#4#5{%
  \reset@font\fontsize{#1}{#2pt}%
  \fontfamily{#3}\fontseries{#4}\fontshape{#5}%
  \selectfont}%
\fi\endgroup%
\begin{picture}(16025,3864)(-274,-7144)
\put(13426,-5536){\makebox(0,0)[lb]{\smash{\SetFigFont{5}{6.0}{\rmdefault}{\mddefault}{\updefault}{\color[rgb]{0,0,0}{\normalsize $I$}}%
}}}
\put(1201,-3436){\makebox(0,0)[lb]{\smash{\SetFigFont{5}{6.0}{\rmdefault}{\mddefault}{\updefault}{\color[rgb]{0,0,0}{\normalsize $\Sigma_{\Gamma+}^\del$}}%
}}}
\put(3901,-3436){\makebox(0,0)[lb]{\smash{\SetFigFont{5}{6.0}{\rmdefault}{\mddefault}{\updefault}{\color[rgb]{0,0,0}{\normalsize $\Sigma_{\Gamma+}^\del$}}%
}}}
\put(12001,-3436){\makebox(0,0)[lb]{\smash{\SetFigFont{5}{6.0}{\rmdefault}{\mddefault}{\updefault}{\color[rgb]{0,0,0}{\normalsize $\Sigma_{\Gamma+}^\del$}}%
}}}
\put(15751,-4411){\makebox(0,0)[lb]{\smash{\SetFigFont{5}{6.0}{\rmdefault}{\mddefault}{\updefault}{\color[rgb]{0,0,0}{\normalsize $\Sigma_{\Gamma+}^\nab$}}%
}}}
\put(14701,-3436){\makebox(0,0)[lb]{\smash{\SetFigFont{5}{6.0}{\rmdefault}{\mddefault}{\updefault}{\color[rgb]{0,0,0}{\normalsize $\Sigma_{\Gamma+}^\nab$}}%
}}}
\put(9301,-3436){\makebox(0,0)[lb]{\smash{\SetFigFont{5}{6.0}{\rmdefault}{\mddefault}{\updefault}{\color[rgb]{0,0,0}{\normalsize $\Sigma_{\Gamma+}^\nab$}}%
}}}
\put(6601,-3436){\makebox(0,0)[lb]{\smash{\SetFigFont{5}{6.0}{\rmdefault}{\mddefault}{\updefault}{\color[rgb]{0,0,0}{\normalsize $\Sigma_{\Gamma+}^\nab$}}%
}}}
\put(1201,-7086){\makebox(0,0)[lb]{\smash{\SetFigFont{5}{6.0}{\rmdefault}{\mddefault}{\updefault}{\color[rgb]{0,0,0}{\normalsize $\Sigma_{\Gamma-}^\del$}}%
}}}
\put(3901,-7086){\makebox(0,0)[lb]{\smash{\SetFigFont{5}{6.0}{\rmdefault}{\mddefault}{\updefault}{\color[rgb]{0,0,0}{\normalsize $\Sigma_{\Gamma-}^\del$}}%
}}}
\put(12001,-7086){\makebox(0,0)[lb]{\smash{\SetFigFont{5}{6.0}{\rmdefault}{\mddefault}{\updefault}{\color[rgb]{0,0,0}{\normalsize $\Sigma_{\Gamma-}^\del$}}%
}}}
\put(15751,-6061){\makebox(0,0)[lb]{\smash{\SetFigFont{5}{6.0}{\rmdefault}{\mddefault}{\updefault}{\color[rgb]{0,0,0}{\normalsize $\Sigma_{\Gamma-}^\nab$}}%
}}}
\put(14701,-7086){\makebox(0,0)[lb]{\smash{\SetFigFont{5}{6.0}{\rmdefault}{\mddefault}{\updefault}{\color[rgb]{0,0,0}{\normalsize $\Sigma_{\Gamma-}^\nab$}}%
}}}
\put(9301,-7086){\makebox(0,0)[lb]{\smash{\SetFigFont{5}{6.0}{\rmdefault}{\mddefault}{\updefault}{\color[rgb]{0,0,0}{\normalsize $\Sigma_{\Gamma-}^\nab$}}%
}}}
\put(6601,-7086){\makebox(0,0)[lb]{\smash{\SetFigFont{5}{6.0}{\rmdefault}{\mddefault}{\updefault}{\color[rgb]{0,0,0}{\normalsize $\Sigma_{\Gamma-}^\nab$}}%
}}}
\put(2551,-3436){\makebox(0,0)[lb]{\smash{\SetFigFont{5}{6.0}{\rmdefault}{\mddefault}{\updefault}{\color[rgb]{0,0,0}{\normalsize $\Sigma_{I+}$}}%
}}}
\put(5251,-3436){\makebox(0,0)[lb]{\smash{\SetFigFont{5}{6.0}{\rmdefault}{\mddefault}{\updefault}{\color[rgb]{0,0,0}{\normalsize $\Sigma_{I+}$}}%
}}}
\put(7951,-3436){\makebox(0,0)[lb]{\smash{\SetFigFont{5}{6.0}{\rmdefault}{\mddefault}{\updefault}{\color[rgb]{0,0,0}{\normalsize $\Sigma_{I+}$}}%
}}}
\put(10651,-3436){\makebox(0,0)[lb]{\smash{\SetFigFont{5}{6.0}{\rmdefault}{\mddefault}{\updefault}{\color[rgb]{0,0,0}{\normalsize $\Sigma_{I+}$}}%
}}}
\put(13351,-3436){\makebox(0,0)[lb]{\smash{\SetFigFont{5}{6.0}{\rmdefault}{\mddefault}{\updefault}{\color[rgb]{0,0,0}{\normalsize $\Sigma_{I+}$}}%
}}}
\put(13351,-7086){\makebox(0,0)[lb]{\smash{\SetFigFont{5}{6.0}{\rmdefault}{\mddefault}{\updefault}{\color[rgb]{0,0,0}{\normalsize $\Sigma_{I-}$}}%
}}}
\put(10651,-7086){\makebox(0,0)[lb]{\smash{\SetFigFont{5}{6.0}{\rmdefault}{\mddefault}{\updefault}{\color[rgb]{0,0,0}{\normalsize $\Sigma_{I-}$}}%
}}}
\put(7951,-7086){\makebox(0,0)[lb]{\smash{\SetFigFont{5}{6.0}{\rmdefault}{\mddefault}{\updefault}{\color[rgb]{0,0,0}{\normalsize $\Sigma_{I-}$}}%
}}}
\put(5251,-7086){\makebox(0,0)[lb]{\smash{\SetFigFont{5}{6.0}{\rmdefault}{\mddefault}{\updefault}{\color[rgb]{0,0,0}{\normalsize $\Sigma_{I-}$}}%
}}}
\put(2551,-7086){\makebox(0,0)[lb]{\smash{\SetFigFont{5}{6.0}{\rmdefault}{\mddefault}{\updefault}{\color[rgb]{0,0,0}{\normalsize $\Sigma_{I-}$}}%
}}}
\put(2251,-4411){\makebox(0,0)[lb]{\smash{\SetFigFont{5}{6.0}{\rmdefault}{\mddefault}{\updefault}{\color[rgb]{0,0,0}{\normalsize $\Sigma_{0+}^\del$}}%
}}}
\put(3451,-4411){\makebox(0,0)[lb]{\smash{\SetFigFont{5}{6.0}{\rmdefault}{\mddefault}{\updefault}{\color[rgb]{0,0,0}{\normalsize $\Sigma_{0+}^\del$}}%
}}}
\put(4951,-4411){\makebox(0,0)[lb]{\smash{\SetFigFont{5}{6.0}{\rmdefault}{\mddefault}{\updefault}{\color[rgb]{0,0,0}{\normalsize $\Sigma_{0+}^\del$}}%
}}}
\put(11551,-4411){\makebox(0,0)[lb]{\smash{\SetFigFont{5}{6.0}{\rmdefault}{\mddefault}{\updefault}{\color[rgb]{0,0,0}{\normalsize $\Sigma_{0+}^\del$}}%
}}}
\put(13051,-4411){\makebox(0,0)[lb]{\smash{\SetFigFont{5}{6.0}{\rmdefault}{\mddefault}{\updefault}{\color[rgb]{0,0,0}{\normalsize $\Sigma_{0+}^\del$}}%
}}}
\put(14251,-4411){\makebox(0,0)[lb]{\smash{\SetFigFont{5}{6.0}{\rmdefault}{\mddefault}{\updefault}{\color[rgb]{0,0,0}{\normalsize $\Sigma_{0+}^\nab$}}%
}}}
\put(10351,-4411){\makebox(0,0)[lb]{\smash{\SetFigFont{5}{6.0}{\rmdefault}{\mddefault}{\updefault}{\color[rgb]{0,0,0}{\normalsize $\Sigma_{0+}^\nab$}}%
}}}
\put(8851,-4411){\makebox(0,0)[lb]{\smash{\SetFigFont{5}{6.0}{\rmdefault}{\mddefault}{\updefault}{\color[rgb]{0,0,0}{\normalsize $\Sigma_{0+}^\nab$}}%
}}}
\put(7651,-4411){\makebox(0,0)[lb]{\smash{\SetFigFont{5}{6.0}{\rmdefault}{\mddefault}{\updefault}{\color[rgb]{0,0,0}{\normalsize $\Sigma_{0+}^\nab$}}%
}}}
\put(6151,-4411){\makebox(0,0)[lb]{\smash{\SetFigFont{5}{6.0}{\rmdefault}{\mddefault}{\updefault}{\color[rgb]{0,0,0}{\normalsize $\Sigma_{0+}^\nab$}}%
}}}
\put(2251,-6061){\makebox(0,0)[lb]{\smash{\SetFigFont{5}{6.0}{\rmdefault}{\mddefault}{\updefault}{\color[rgb]{0,0,0}{\normalsize $\Sigma_{0-}^\del$}}%
}}}
\put(3451,-6061){\makebox(0,0)[lb]{\smash{\SetFigFont{5}{6.0}{\rmdefault}{\mddefault}{\updefault}{\color[rgb]{0,0,0}{\normalsize $\Sigma_{0-}^\del$}}%
}}}
\put(4951,-6061){\makebox(0,0)[lb]{\smash{\SetFigFont{5}{6.0}{\rmdefault}{\mddefault}{\updefault}{\color[rgb]{0,0,0}{\normalsize $\Sigma_{0-}^\del$}}%
}}}
\put(11551,-6061){\makebox(0,0)[lb]{\smash{\SetFigFont{5}{6.0}{\rmdefault}{\mddefault}{\updefault}{\color[rgb]{0,0,0}{\normalsize $\Sigma_{0-}^\del$}}%
}}}
\put(13051,-6061){\makebox(0,0)[lb]{\smash{\SetFigFont{5}{6.0}{\rmdefault}{\mddefault}{\updefault}{\color[rgb]{0,0,0}{\normalsize $\Sigma_{0-}^\del$}}%
}}}
\put(10351,-6061){\makebox(0,0)[lb]{\smash{\SetFigFont{5}{6.0}{\rmdefault}{\mddefault}{\updefault}{\color[rgb]{0,0,0}{\normalsize $\Sigma_{0-}^\nab$}}%
}}}
\put(8851,-6061){\makebox(0,0)[lb]{\smash{\SetFigFont{5}{6.0}{\rmdefault}{\mddefault}{\updefault}{\color[rgb]{0,0,0}{\normalsize $\Sigma_{0-}^\nab$}}%
}}}
\put(7651,-6061){\makebox(0,0)[lb]{\smash{\SetFigFont{5}{6.0}{\rmdefault}{\mddefault}{\updefault}{\color[rgb]{0,0,0}{\normalsize $\Sigma_{0-}^\nab$}}%
}}}
\put(6151,-6061){\makebox(0,0)[lb]{\smash{\SetFigFont{5}{6.0}{\rmdefault}{\mddefault}{\updefault}{\color[rgb]{0,0,0}{\normalsize $\Sigma_{0-}^\nab$}}%
}}}
\put(1201,-4936){\makebox(0,0)[lb]{\smash{\SetFigFont{5}{6.0}{\rmdefault}{\mddefault}{\updefault}{\color[rgb]{0,0,0}{\normalsize $\Gamma$}}%
}}}
\put(3901,-4936){\makebox(0,0)[lb]{\smash{\SetFigFont{5}{6.0}{\rmdefault}{\mddefault}{\updefault}{\color[rgb]{0,0,0}{\normalsize $\Gamma$}}%
}}}
\put(6601,-4936){\makebox(0,0)[lb]{\smash{\SetFigFont{5}{6.0}{\rmdefault}{\mddefault}{\updefault}{\color[rgb]{0,0,0}{\normalsize $\Gamma$}}%
}}}
\put(9301,-4936){\makebox(0,0)[lb]{\smash{\SetFigFont{5}{6.0}{\rmdefault}{\mddefault}{\updefault}{\color[rgb]{0,0,0}{\normalsize $\Gamma$}}%
}}}
\put(12001,-4936){\makebox(0,0)[lb]{\smash{\SetFigFont{5}{6.0}{\rmdefault}{\mddefault}{\updefault}{\color[rgb]{0,0,0}{\normalsize $\Gamma$}}%
}}}
\put(14701,-4936){\makebox(0,0)[lb]{\smash{\SetFigFont{5}{6.0}{\rmdefault}{\mddefault}{\updefault}{\color[rgb]{0,0,0}{\normalsize $\Gamma$}}%
}}}
\put(14251,-6061){\makebox(0,0)[lb]{\smash{\SetFigFont{5}{6.0}{\rmdefault}{\mddefault}{\updefault}{\color[rgb]{0,0,0}{\normalsize $\Sigma_{0-}^\nab$}}%
}}}
\put(-274,-4486){\makebox(0,0)[lb]{\smash{\SetFigFont{5}{6.0}{\rmdefault}{\mddefault}{\updefault}{\color[rgb]{0,0,0}{\normalsize $\Sigma_{\Gamma+}^\del$}}%
}}}
\put(-274,-5986){\makebox(0,0)[lb]{\smash{\SetFigFont{5}{6.0}{\rmdefault}{\mddefault}{\updefault}{\color[rgb]{0,0,0}{\normalsize $\Sigma_{\Gamma-}^\del$}}%
}}}
\put(2626,-5536){\makebox(0,0)[lb]{\smash{\SetFigFont{5}{6.0}{\rmdefault}{\mddefault}{\updefault}{\color[rgb]{0,0,0}{\normalsize $I$}}%
}}}
\put(5326,-5536){\makebox(0,0)[lb]{\smash{\SetFigFont{5}{6.0}{\rmdefault}{\mddefault}{\updefault}{\color[rgb]{0,0,0}{\normalsize $I$}}%
}}}
\put(8026,-5536){\makebox(0,0)[lb]{\smash{\SetFigFont{5}{6.0}{\rmdefault}{\mddefault}{\updefault}{\color[rgb]{0,0,0}{\normalsize $I$}}%
}}}
\put(10726,-5536){\makebox(0,0)[lb]{\smash{\SetFigFont{5}{6.0}{\rmdefault}{\mddefault}{\updefault}{\color[rgb]{0,0,0}{\normalsize $I$}}%
}}}
\end{picture}

%% file: SigmaModel.pstex_t
\begin{picture}(0,0)%
\epsfig{file=SigmaModel.pstex}%
\end{picture}%
\setlength{\unitlength}{1579sp}%
\begingroup\makeatletter\ifx\SetFigFont\undefined%
\gdef\SetFigFont#1#2#3#4#5{%
  \reset@font\fontsize{#1}{#2pt}%
  \fontfamily{#3}\fontseries{#4}\fontshape{#5}%
  \selectfont}%
\fi\endgroup%
\begin{picture}(15525,3066)(226,-6694)
\put(15751,-5236){\makebox(0,0)[lb]{\smash{\SetFigFont{5}{6.0}{\rmdefault}{\mddefault}{\updefault}{\color[rgb]{0,0,0}{\normalsize $b$}}%
}}}
\put(2551,-4936){\makebox(0,0)[lb]{\smash{\SetFigFont{5}{6.0}{\rmdefault}{\mddefault}{\updefault}{\color[rgb]{0,0,0}{\normalsize $I_0$}}%
}}}
\put(5251,-4936){\makebox(0,0)[lb]{\smash{\SetFigFont{5}{6.0}{\rmdefault}{\mddefault}{\updefault}{\color[rgb]{0,0,0}{\normalsize $I_1$}}%
}}}
\put(7951,-4936){\makebox(0,0)[lb]{\smash{\SetFigFont{5}{6.0}{\rmdefault}{\mddefault}{\updefault}{\color[rgb]{0,0,0}{\normalsize $I_2$}}%
}}}
\put(3901,-5536){\makebox(0,0)[lb]{\smash{\SetFigFont{5}{6.0}{\rmdefault}{\mddefault}{\updefault}{\color[rgb]{0,0,0}{\normalsize $\Gamma_1$}}%
}}}
\put(6601,-5536){\makebox(0,0)[lb]{\smash{\SetFigFont{5}{6.0}{\rmdefault}{\mddefault}{\updefault}{\color[rgb]{0,0,0}{\normalsize $\Gamma_2$}}%
}}}
\put(9301,-5536){\makebox(0,0)[lb]{\smash{\SetFigFont{5}{6.0}{\rmdefault}{\mddefault}{\updefault}{\color[rgb]{0,0,0}{\normalsize $\Gamma_3$}}%
}}}
\put(12001,-5536){\makebox(0,0)[lb]{\smash{\SetFigFont{5}{6.0}{\rmdefault}{\mddefault}{\updefault}{\color[rgb]{0,0,0}{\normalsize $\Gamma_4$}}%
}}}
\put(10651,-4936){\makebox(0,0)[lb]{\smash{\SetFigFont{5}{6.0}{\rmdefault}{\mddefault}{\updefault}{\color[rgb]{0,0,0}{\normalsize $I_3$}}%
}}}
\put(13351,-4936){\makebox(0,0)[lb]{\smash{\SetFigFont{5}{6.0}{\rmdefault}{\mddefault}{\updefault}{\color[rgb]{0,0,0}{\normalsize $I_4$}}%
}}}
\put(226,-5236){\makebox(0,0)[lb]{\smash{\SetFigFont{5}{6.0}{\rmdefault}{\mddefault}{\updefault}{\color[rgb]{0,0,0}{\normalsize $a$}}%
}}}
\end{picture}

%% file: SigmaE.pstex_t
\begin{picture}(0,0)%
\epsfig{file=SigmaE.pstex}%
\end{picture}%
\setlength{\unitlength}{1579sp}%
\begingroup\makeatletter\ifx\SetFigFont\undefined%
\gdef\SetFigFont#1#2#3#4#5{%
  \reset@font\fontsize{#1}{#2pt}%
  \fontfamily{#3}\fontseries{#4}\fontshape{#5}%
  \selectfont}%
\fi\endgroup%
\begin{picture}(17250,4039)(-1499,-7169)
\put(15751,-5236){\makebox(0,0)[lb]{\smash{\SetFigFont{5}{6.0}{\rmdefault}{\mddefault}{\updefault}{\color[rgb]{0,0,0}{\normalsize $b$}}%
}}}
\put( 76,-5236){\makebox(0,0)[lb]{\smash{\SetFigFont{5}{6.0}{\rmdefault}{\mddefault}{\updefault}{\color[rgb]{0,0,0}{\normalsize $a$}}%
}}}
\put(-1499,-3736){\makebox(0,0)[lb]{\smash{\SetFigFont{5}{6.0}{\rmdefault}{\mddefault}{\updefault}{\color[rgb]{0,0,0}{\normalsize $\Im(z)=\epsilon$}}%
}}}
\put(14776,-7111){\makebox(0,0)[lb]{\smash{\SetFigFont{5}{6.0}{\rmdefault}{\mddefault}{\updefault}{\color[rgb]{0,0,0}{\normalsize $\Re(z)=b$}}%
}}}
\put(-299,-3286){\makebox(0,0)[lb]{\smash{\SetFigFont{5}{6.0}{\rmdefault}{\mddefault}{\updefault}{\color[rgb]{0,0,0}{\normalsize $\Re(z)=a$}}%
}}}
\put(-1499,-6736){\makebox(0,0)[lb]{\smash{\SetFigFont{5}{6.0}{\rmdefault}{\mddefault}{\updefault}{\color[rgb]{0,0,0}{\normalsize $\Im(z)=-\epsilon$}}%
}}}
\end{picture}

%% file: SigmaF.pstex_t
\begin{picture}(0,0)%
\epsfig{file=SigmaF.pstex}%
\end{picture}%
\setlength{\unitlength}{1579sp}%
\begingroup\makeatletter\ifx\SetFigFont\undefined%
\gdef\SetFigFont#1#2#3#4#5{%
  \reset@font\fontsize{#1}{#2pt}%
  \fontfamily{#3}\fontseries{#4}\fontshape{#5}%
  \selectfont}%
\fi\endgroup%
\begin{picture}(15498,3216)(352,-6769)
\end{picture}

%% file: SigmaModelHomology.pstex_t
\begin{picture}(0,0)%
\epsfig{file=SigmaModelHomology.pstex}%
\end{picture}%
\setlength{\unitlength}{1579sp}%
\begingroup\makeatletter\ifx\SetFigFont\undefined%
\gdef\SetFigFont#1#2#3#4#5{%
  \reset@font\fontsize{#1}{#2pt}%
  \fontfamily{#3}\fontseries{#4}\fontshape{#5}%
  \selectfont}%
\fi\endgroup%
\begin{picture}(16866,3870)(-332,-7069)
\put(4276,-4636){\makebox(0,0)[lb]{\smash{\SetFigFont{5}{6.0}{\rmdefault}{\mddefault}{\updefault}{\color[rgb]{0,0,0}{\normalsize $a_1$}}%
}}}
\put(6976,-4636){\makebox(0,0)[lb]{\smash{\SetFigFont{5}{6.0}{\rmdefault}{\mddefault}{\updefault}{\color[rgb]{0,0,0}{\normalsize $a_2$}}%
}}}
\put(9676,-4636){\makebox(0,0)[lb]{\smash{\SetFigFont{5}{6.0}{\rmdefault}{\mddefault}{\updefault}{\color[rgb]{0,0,0}{\normalsize $a_3$}}%
}}}
\put(12376,-4636){\makebox(0,0)[lb]{\smash{\SetFigFont{5}{6.0}{\rmdefault}{\mddefault}{\updefault}{\color[rgb]{0,0,0}{\normalsize $a_4$}}%
}}}
\put(5401,-6286){\makebox(0,0)[lb]{\smash{\SetFigFont{5}{6.0}{\rmdefault}{\mddefault}{\updefault}{\color[rgb]{0,0,0}{\normalsize $b_3$}}%
}}}
\put(2701,-5461){\makebox(0,0)[lb]{\smash{\SetFigFont{5}{6.0}{\rmdefault}{\mddefault}{\updefault}{\color[rgb]{0,0,0}{\normalsize $b_1$}}%
}}}
\put(15601,-5536){\makebox(0,0)[lb]{\smash{\SetFigFont{5}{6.0}{\rmdefault}{\mddefault}{\updefault}{\color[rgb]{0,0,0}{\normalsize $\Gamma_0$}}%
}}}
\put(226,-5536){\makebox(0,0)[lb]{\smash{\SetFigFont{5}{6.0}{\rmdefault}{\mddefault}{\updefault}{\color[rgb]{0,0,0}{\normalsize $\Gamma_0$}}%
}}}
\put(3826,-5911){\makebox(0,0)[lb]{\smash{\SetFigFont{5}{6.0}{\rmdefault}{\mddefault}{\updefault}{\color[rgb]{0,0,0}{\normalsize $b_2$}}%
}}}
\put(6751,-6736){\makebox(0,0)[lb]{\smash{\SetFigFont{5}{6.0}{\rmdefault}{\mddefault}{\updefault}{\color[rgb]{0,0,0}{\normalsize $b_4$}}%
}}}
\end{picture}